\documentclass[12pt]{amsbook}
\usepackage{amssymb}
\usepackage[all]{xy}

\usepackage[colorlinks=true,linkcolor=magenta,citecolor=magenta]{hyperref}
\makeindex

\newtheorem{theorem}{Theorem}[chapter]
\newtheorem{answer}[theorem]{Answer}
\newtheorem{definition}[theorem]{Definition}
\newtheorem{exercise}[theorem]{Exercise}
\newtheorem{fact}[theorem]{Fact}
\newtheorem{proposition}[theorem]{Proposition}
\newtheorem{question}[theorem]{Question}

\begin{document}

\title{Affine noncommutative geometry}

\author{Teo Banica}
\address{Department of Mathematics, University of Cergy-Pontoise, F-95000 Cergy-Pontoise, France. {\tt teo.banica@gmail.com}}

\subjclass[2010]{46L05}
\keywords{Affine geometry, Quantum manifold}

\begin{abstract}
This is an introduction to noncommutative geometry, from an affine viewpoint, that is, by using coordinates. The spaces $\mathbb R^N,\mathbb C^N$ have no free analogues in the operator algebra sense, but the corresponding unit spheres $S^{N-1}_\mathbb R,S^{N-1}_\mathbb C$ do have free analogues  $S^{N-1}_{\mathbb R,+},S^{N-1}_{\mathbb C,+}$. There are many examples of real algebraic submanifolds $X\subset S^{N-1}_{\mathbb R,+},S^{N-1}_{\mathbb C,+}$, some of which are of Riemannian flavor, coming with a Haar integration functional $\int:C(X)\to\mathbb C$, that we will study here. We will mostly focus on free geometry, but we will discuss as well some related geometries, called easy, completing the picture formed by the 4 main geometries, namely real/complex, classical/free.
\end{abstract}

\maketitle

\chapter*{Preface}

Classical geometry has its origins in classical mechanics, with some of its most fundamental objects, such as the conics, coming from the trajectories of planets and other celestial objects around the Sun. Similarly, quantum mechanics has inspired several theories of quantum geometry, more commonly called ``noncommutative geometry''.

\bigskip

The idea of noncommutative geometry goes back to Heisenberg. Back in the 1920s, the main problem in physics was that of understanding the mechanics of the hydrogen atom, and it was known since Bohr that the Maxwell equations do not work. Heisenberg came with a clever idea for solving the problem, namely looking for some sort of ``quantum trajectory'' for the electron, instead of a classical, honest trajectory, and with his mathematics involving the algebra $M_\infty(\mathbb C)$ of the complex infinite matrices.

\bigskip

A few years later Schr\"odinger came with something better, namely a PDE for the wave function of the electron. This improved Heisenberg's findings, with $M_\infty(\mathbb C)$ being now understood to correspond to the algebra of operators on the Hilbert space $H=L^2(\mathbb R^3)$ of such wave functions. Later, Pauli and others added a copy of $\mathbb C^2$ to this space, as to account for the electron spin, and in this form, that of the late 1920s, quantum mechanics was powerful enough for solving many questions, such as the structure of all atoms. Making Bohr's dream, who was the initiator of the whole program, come true.

\bigskip

Einstein disagreed with all this, saying that such things, probability, noncommutative geometry, you name it, while certainly great, should be regarded as being temporary, and so homework for us for going towards determinism, meaning true, honest geometry.

\bigskip

Generally speaking, modern physics is about making Einstein's dream come true. Quantum mechanics has evolved several times since the 1920s, with the fine structure of hydrogen, then with quantum electrodynamics, then with the discovery, at smaller scales, of quantum chromodynamics, then with the Standard Model, then with all sorts of efforts in general quantum field theory, and then with string theory, which is something geometric. Slowly but surely, we are going towards Einstein's determinism.

\bigskip

This being said, I don't know about you, but personally I'm still waiting for nuclear-powered cars, first for the unlimited horsepower, and then for not having to refuel. And also, why not for nuclear watches too, because every time the battery of my Casio gives up, a whole pain with replacing it. And shall we trust modern physics with coming up soon with concrete answers, to these very concrete life questions that we have.

\bigskip

Shall we perhaps downgrade a bit our dreams in physics? Noncommutative geometry, in its modern formulations, is a bit about this. Forget about Einstein's determinism, or rather leave that for later, and more modestly, try instead to have some sort of noncommutative geometry theory working, improving what Heisenberg was saying, and of course, with the whole thing being as close as possible to the modern advances in physics.

\bigskip

The credit for such ideas goes to Connes, who created in the 80s a noncommutative geometry theory which is definitely simple, beautiful, and modern too. Connes looked at the noncommutative manifolds $X$, with this meaning that $A=C(X)$ is an operator algebra, $A\subset B(H)$, which are smooth and Riemannian, in a certain technical sense, with remarkable results in connection with physics, obtained all over the 90s and 00s.

\bigskip

Our aim here is to talk about noncommutative geometry too, from a point of view very close to the one of Connes, but with slightly different motivations in mind. We will keep from Connes his two main principles, namely that the noncommutative manifolds $X$ should appear from operator algebras, and also, that they should be Riemannian. However, based on our belief that at very small scales, smaller than those of the Standard Model, there is no room for smoothness, we will drop the assumption that $X$ should be smooth, and so Riemannian for us will rather mean that $X$ is real algebraic, a bit a la Nash, and coming with an integration functional $\int:C(X)\to\mathbb C$, that we will be eager to compute explicitly, using techniques of Jones, Voiculescu and Woronowicz. 

\bigskip

This book will be purely mathematical. The applications to physics, involving some more mathematics, such as PDE over our free manifolds, will be hopefully discussed in a series of further books. As for nuclear cars and watches, I am currently trying to build some in my garage, for my personal usage, but things difficult here. More later.

\bigskip

This book is partly based on a number of recent papers on quantum groups and noncommutative geometry, and I am particularly grateful to Julien Bichon, for his heavy involvement in the subject. Many thanks go as well to my cats. Their timeless views and opinions, on everyone and everything, have always been of great help.

\bigskip

\

{\em Cergy, October 2026}

\smallskip

{\em Teo Banica}

\baselineskip=15.95pt
\tableofcontents
\baselineskip=14pt

\part{Affine geometry}

\ \vskip50mm

\begin{center}
{\em Hot chili peppers in the blistering sun

Dust on my face and my cape

Me and Magdalena on the run

I think this time we shall escape}
\end{center}

\chapter{Spheres and tori}

\section*{1a. Classical geometries}

We would like to develop some noncommutative geometry theory, that can be of help in quantum mechanics, a bit like classical geometry is of help in classical mechanics. So, this will be a book about mathematical physics. With mathematical physics meaning, as usual, mathematics developed with physics motivations in mind.

\bigskip

Before anything, let us recommend some reading. Physics and quantum mechanics in particular can be learned from the books of Feynman \cite{fey}, or Griffiths \cite{gr1}, \cite{gr2}, or Weinberg \cite{we1}, \cite{we2}. For quantum mechanics, a look at the old books of Dirac \cite{dir}, von Neumann \cite{von} and Weyl \cite{wey} can be instructive too. Also, never forget that physics is a whole, and do not hesitate to complete your electrodynamics and quantum mechanics knowledge with some solid classical mechanics, and some thermodynamics too.

\bigskip

Back to our goals, noncommutative geometry, a look at all this physics does not help that  much. There are certainly a few things to be learned, as for instance the fact that noncommutative geometry should have something to do with the linear operators $T:H\to H$ over a complex Hilbert space $H$. But passed that, we are a bit in the dark. As an example, the Hilbert space $H$ used to be something quite abstract, $H=l^2(\mathbb N)$, for Heisenberg, then something more concrete, $H=L^2(\mathbb R^3)$, for Schr\"odinger, and then Pauli and others added a copy of $K=\mathbb C^2$, in order to take into account the spin of the electron. And this was only what happened in the 1920s, and there is no telling of what happened afterwards, up to the present days. It is probably safe to say that no one really knows what $H$ is. And even worse, no one really knows if there is one such $H$ at all.

\bigskip

The same story goes with the linear operators $T:H\to H$. These used to be densely defined and unbounded, during the good old days of Heisenberg, Schr\"odinger and Dirac. But then, with quantum mechanics evolving into quantum electrodynamics, then into the Standard Model, and then into more technical versions of quantum field theory, such operators often became everywhere defined, and bounded, $T\in B(H)$. And there is even worse, because at the truly advanced level you often manage to find a way to deal with usual complex matrices, $T\in M_N(\mathbb C)$, or sometimes with random matrices.

\bigskip

Summarizing, all this looks complicated, and my proposal would be to leave physics for later, although please have a look at the above-mentioned physics books, because what we will be doing here, and I insist, is mathematical physics, and not pure mathematics, and you won't understand otherwise, and have some inspiration from pure mathematics instead. We have for instance the following question, that we can try to solve:

\begin{question}
What are the noncommutative analogues of the geometry of $\mathbb R^N$, and of the geometry of $\mathbb C^N$?
\end{question}

Here by ``noncommutative'' we mean with the standard coordinates not commuting, $x_ix_j\neq x_jx_i$. But this is something a bit vague, because shall we look here for some weakenings of the commutation relations $x_ix_j=x_jx_i$, and as we will soon see, there are plently of interesting choices here, or shall we just look, for simplifying and to start with, for ``free'' geometries, where $x_i,x_j$ are not subject to any kind of relation.

\bigskip 

As for ``geometry'', things are quite vague here as well. We have indeed algebraic geometry, differential geometry, symplectic geometry, Riemannian geometry, and many more. Also, when talking $\mathbb C^N$, our manifolds can be taken to be real, or complex. And also, all these geometries usually come in two flavors, affine and projective.

\bigskip

In short, we have an idea with our Question 1.1, but everything is still too vague. So here we are back to physics, and quantum mechanics, for inspiration, and based on some knowledge there, let us formulate the following fact:

\begin{fact}
The spaces $\mathbb R^N$ and $\mathbb C^N$ have no interesting free analogues, and this due to the fact that these spaces are not compact.
\end{fact}

Obviously, this is something subjective. My point comes from the fact that, while you can certainly talk about the free algebra generated by variables $x_1,\ldots,x_N$, with some care of course in relation with conjugation and the complex structure, there is no way of putting a reasonable norm on this algebra, due to the fact that $x_1,\ldots,x_N$ are unbounded. And so this algebra is just some abstract, pure mathematical beast, totally unrelated to analysis, physics, probability, quantum mechanics, you name it.

\bigskip

In short, you'll have to trust me here. And we will talk about this more in detail later in this chapter, when learning about operator algebras. By the way, let me point out that Question 1.1 was something controversial too, because we assumed somehow by definition that the interesting fields are $F=\mathbb R,\mathbb C$. Which is far from being something accepted, with many mathematical physicists, myself included, agreeing that, at a more advanced level, other fields than $F=\mathbb R,\mathbb C$ are interesting too. But this is another story.

\bigskip

Going ahead with our study, Fact 1.2 suggests replacing $F^N=\mathbb R^N,\mathbb C^N$ with a suitably chosen compact manifold $X\subset F^N$. And here, we have several choices. The first thought goes to the unit sphere $S\subset F^N$. But then, why not looking instead at the unit torus $T\subset F^N$, which in addition is a group. And then, talking now groups, why not being a bit advanced, and looking instead at the unitary group $U\subset\mathcal L(F^N)$. And then why not, for being even more advanced, looking at the reflection group $K\subset\mathcal L(F^N)$.

\bigskip

These choices are all reasonable, with the mathematics of each of $S,T,U,K$ being well-known to be interesting, and describing a bit of the mathematics of $F^N$ itself. So, why not putting all these objects $S,T,U,K$ together, as to have a complete, robust object replacing $F^N$. We are led in this way to the following answer to Question 1.1:

\begin{answer}
An affine noncommutative geometry should come from a diagram
$$
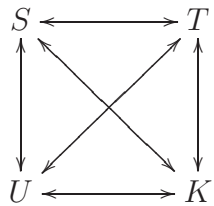
\xymatrix@R=50pt@C=50pt{
S\ar[r]\ar[d]\ar[dr]&T\ar[l]\ar[d]\ar[dl]\\
U\ar[u]\ar[ur]\ar[r]&K\ar[l]\ar[ul]\ar[u]
}$$
consisting of a sphere $S$, a torus $T$, a unitary group $U$, and a reflection group $K$.
\end{answer}

So, this will be our starting point, for the considerations in this book, and our guiding philosophy, in what follows, until the end. However, obviously, this is something quite advanced, and before starting our study, a few comments, and some advice:

\bigskip

(1) As already said, and above everything, this is something advanced. The above answer emerged in the late 10s, based on substantial work, a few dozen research papers, written all over the 90s, 00s, 10s. So, modesty, and it will take us a few 100 pages in this book, or perhaps the whole book, to understand what Answer 1.3 really says.

\bigskip

(2) Also importantly, Answer 1.3 talks about ``noncommutative geometries'' in general, and not about manifolds constructed inside such geometries. We will first try to understand these noncommutative geometries themselves, via some axiomatization and classification work. And the study of the corresponding manifolds will come after.

\bigskip

(3) Looking now at Answer 1.3 as it is, that sounds like some kind of pure mathematics, for the most involving algebraic geometry and Lie groups. Although not really necessary for reading this book, some knowledge here would be welcome. You can learn these for instance from Shafarevich \cite{sha} and Fulton-Harris \cite{fha}, respectively.

\bigskip

(4) In fact, we will see in a moment that what is really needed for understanding Answer 1.3 is rather basic analysis. So, getting now to the true prerequisites for the present book, these will be Rudin \cite{rud}. With perhaps a bit of familiarity with algebra too, say from Lang \cite{lng}, and a bit of functional analysis too, say from Lax \cite{lax}.
 
\bigskip

Getting started now, let us first discuss the case of the usual geometry, in $\mathbb R^N$. We must construct here the corresponding quadruplet $(S,T,U,K)$, as per the requirements of Answer 1.3, and the definitions here, all very familiar, are as follows:

\index{real torus}
\index{hyperoctahedral group}

\begin{definition}
The real sphere, torus, unitary group and reflection group are:
\begin{eqnarray*}
S^{N-1}_\mathbb R&=&\left\{x\in\mathbb R^N\Big|\sum_ix_i^2=1\right\}\\
T_N&=&\left\{x\in\mathbb R^N\Big|x_i=\pm\frac{1}{\sqrt{N}}\right\}\\
O_N&=&\left\{U\in M_N(\mathbb  R)\Big|U^t=U^{-1}\right\}\\
H_N&=&\left\{U\in M_N(-1,0,1)\Big|U^t=U^{-1}\right\}
\end{eqnarray*}
These are the usual sphere, cube, orthogonal group, and hyperoctahedral group.
\end{definition}

To be more precise here, all the objects on the right are certainly familiar, but the notations and terminology for them are perhaps not, and here are the details:

\bigskip

(1) The sphere $S^{N-1}_\mathbb R$ is the unit sphere of $\mathbb R^N$ as we know it, and we will often say sphere instead of unit sphere. As for the superscript $N-1$, which is very standard, this stands for the real dimension as a manifold, which is indeed $N-1$. 

\bigskip

(2) Regarding $T_N$, this is the standard cube in $\mathbb R^N$, with the $1/\sqrt{N}$ normalization being there in order to have an embedding $T_N\subset S^{N-1}_\mathbb R$, which will be useful for us. We also call $T_N$ torus, standing for ``discrete torus'', and more on this later.

\bigskip

(3) Regarding now $O_N$, this is the orthogonal group as we know it. We also call it unitary group of $\mathbb R^N$, because, a bit as for the cube/torus before, we are using here in this book a hybrid real/complex terminology for everything. More on this later.

\bigskip

(4) Finally, regarding $H_N$, this is the hyperoctahedral group, which is by definition the symmetry group of the hypercube in $\mathbb R^N$, which means our cube/torus $T_N$. The formula for $H_N$ in the statement is something elementary, coming from definitions.

\bigskip

Regarding now the correspondences between our objects, there are many ways of establishing them, depending on knowledge and taste. Assuming that you followed my advice, and that you are a bit familiar with basic algebraic geometry and Lie theory, you would probably say that our objects $(S,T,U,K)$ are trivially in correspondence with each other, QED. On the opposite, assuming that you are not familiar with all this, and that your mathematical background is basically the officially needed Rudin \cite{rud}, along with a bit of Lang \cite{lng} and Lax \cite{lax}, you might have a bit of troubles with all this.

\bigskip

Well, good news, establishing the correspondences for $\mathbb R^N$ is actually not crucial for us, at this point. And with this coming from the fact that, no matter what things can be said about $\mathbb R^N$, we will be doing noncommutative geometry in this book, and in the noncommutative setting things are far more rigid, and so the correspondences between $(S,T,U,K)$, even for $\mathbb R^N$, are to be discussed later, once we know what we're doing.

\bigskip

So, in the hope that you got my point. We are just having some preliminary fun, and we need to know, for getting started, that we are on the right track, and that our objects $(S,T,U,K)$ from Definition 1.4 are indeed in correspondence, be that a bit informal. So here is the statement, formulated informally, and coming with an informal proof, and with this, being informal, being, and I insist, the right thing to do, right now:

\begin{theorem}
We have a full set of correspondences, as follows,
$$\xymatrix@R=50pt@C=50pt{
S^{N-1}_\mathbb R\ar[r]\ar[d]\ar[dr]&T_N\ar[l]\ar[d]\ar[dl]\\
O_N\ar[u]\ar[ur]\ar[r]&H_N\ar[l]\ar[ul]\ar[u]
}$$
obtained via various results from basic geometry and group theory.
\end{theorem}

\begin{proof}
As already mentioned, there are several possible solutions to the problem, and all this is not crucial for us. Here is a way of constructing these correspondences:

\medskip

(1) $S^{N-1}_\mathbb R\leftrightarrow T_N$. Here $T_N$ comes from $S^{N-1}_\mathbb R$ via $|x_1|=\ldots=|x_N|$, while $S^{N-1}_\mathbb R$ appears from $T_N\subset\mathbb R^N$ by ``deleting'' this relation, while still keeping $\sum_ix_i^2=1$. 

\medskip

(2) $S^{N-1}_\mathbb R\leftrightarrow O_N$. This comes from the fact that $O_N$ is the isometry group of $S^{N-1}_\mathbb R$, and that, conversely, $S^{N-1}_\mathbb R$ appears as $\{Ux|U\in O_N\}$, where $x=(1,0,\ldots,0)$.

\medskip

(3) $S^{N-1}_\mathbb R\leftrightarrow H_N$. This is something trickier, but the passage can definitely be obtained, for instance via $T_N$, by using the constructions in (1) above and (5) below.

\medskip

(4) $T_N\leftrightarrow O_N$. Here $T_N\simeq\mathbb Z_2^N$ is a maximal torus of $O_N$, and the group $O_N$ itself can be reconstructed from this maximal torus, by using various methods.

\medskip

(5) $T_N\leftrightarrow H_N$. Here, similarly, $T_N\simeq\mathbb Z_2^N$ is a maximal torus of $H_N$, and the group $H_N$ itself can be reconstructed from this torus as a wreath product, $H_N=T_N\wr S_N$.

\medskip

(6) $O_N\leftrightarrow H_N$. This is once again something trickier, but the passage can definitely be obtained, for instance via $T_N$, by using the constructions in (4) and (5) above.
\end{proof}

The above result is of course something quite non-trivial, and having it understood properly would take some time. However, as already said, we will technically not need all this. Our purpose for the moment is just to explain our $(S,T,U,K)$ philosophy.

\bigskip

As a second basic example of geometry, we have the usual geometry of $\mathbb C^N$. Here the corresponding quadruplet $(S,T,U,K)$ can be constructed as follows:

\index{complex torus}
\index{complex reflection group}

\begin{definition}
The complex sphere, torus, unitary group and reflection group are:
\begin{eqnarray*}
S^{N-1}_\mathbb C&=&\left\{x\in\mathbb C^N\Big|\sum_i|x_i|^2=1\right\}\\
\mathbb T_N&=&\left\{x\in\mathbb C^N\Big|\ |x_i|=\frac{1}{\sqrt{N}}\right\}\\
U_N&=&\left\{U\in M_N(\mathbb  C)\Big|U^*=U^{-1}\right\}\\
K_N&=&\left\{U\in M_N(\mathbb T\cup\{0\})\Big|U^*=U^{-1}\right\}
\end{eqnarray*}
These are the usual complex sphere, torus, unitary group, and complex reflection group.
\end{definition}

As before, the superscript $N-1$ for the sphere does not fit with the rest, but is quite standard, somewhat coming from dimension considerations. We will use it as such. Also, the $1/\sqrt{N}$ factor is there in order to have an embedding $\mathbb T_N\subset S^{N-1}_\mathbb C$.

\bigskip

Also as before, in what regards the correspondences between our objects, there are many ways of establishing them, will all this being not crucial for us. In analogy with Theorem 1.5, let us formulate a second informal statement, as follows:

\begin{theorem}
We have a full set of correspondences, as follows,
$$\xymatrix@R=50pt@C=50pt{
S^{N-1}_\mathbb C\ar[r]\ar[d]\ar[dr]&\mathbb T_N\ar[l]\ar[d]\ar[dl]\\
U_N\ar[u]\ar[ur]\ar[r]&K_N\ar[l]\ar[ul]\ar[u]
}$$
obtained via various results from basic geometry and group theory.
\end{theorem}

\begin{proof}
We follow the proof in the real case, by making adjustments where needed, and with of course the reiterated comment that all this is not crucial for us:

\medskip

(1) $S^{N-1}_\mathbb C\leftrightarrow\mathbb T_N$. Same proof as before, using $|x_1|=\ldots=|x_N|$. 

\medskip

(2) $S^{N-1}_\mathbb C\leftrightarrow U_N$. Here ``isometry'' must be taken in an affine complex sense.

\medskip

(3) $S^{N-1}_\mathbb C\leftrightarrow K_N$. Trickier as before, best viewed by passing via $\mathbb T_N$.

\medskip

(4) $\mathbb T_N\leftrightarrow U_N$. Coming from the fact that $\mathbb T_N\simeq\mathbb T^N$ is a maximal torus of $U_N$.

\medskip

(5) $\mathbb T_N\leftrightarrow K_N$. Once again, maximal torus argument, and $K_N=\mathbb T_N\wr S_N$.

\medskip

(6) $U_N\leftrightarrow K_N$. Trickier as before, best viewed by passing via $\mathbb T_N$.
\end{proof}

As a conclusion, our $(S,T,U,K)$ philosophy seems to work, in the sense that these 4 objects, and the relations between them, encode interesting facts about $\mathbb R^N,\mathbb C^N$. Our plan in what follows will be that of leaving aside the complete understanding of what has been said above, and going directly for the noncommutative case. We will see that in the noncommutative setting things are more rigid, and therefore, simpler. And then, with this simpler theory axiomatized, we will come back of course to $\mathbb R^N,\mathbb C^N$, with full details about the correspondences between $(S,T,U,K)$ here, don't worry about that.

\section*{1b. Quantum spaces}

In order to talk about noncommutative geometry, the idea will be that of defining our quantum spaces $X$ as being abstract manifolds, whose coordinates $x_1,\ldots,x_N$ do not necessarily commute. Thus, we are in need of some good algebraic geometry correspondence, between such abstract spaces $X$, and the corresponding algebras of coordinates $A$. Following Heisenberg, von Neumann and many others, we will use here the correspondence $A=C(X)$ coming from operator algebra theory. Let us start with:

\index{Hilbert space}

\begin{definition}
A Hilbert space is a complex vector space $H$, given with a scalar product $<x,y>$, satisfying the following conditions:
\begin{enumerate}
\item $<x,y>$ is linear in $x$, and antilinear in $y$.

\item $\overline{<x,y>}=<y,x>$, for any $x,y$.

\item $<x,x>>0$, for any $x\neq0$.

\item $H$ is complete with respect to the norm $||x||=\sqrt{<x,x>}$.
\end{enumerate}
\end{definition}

Observe that we are using here mathematicians' convention for linearity, as opposed to Dirac's convention in \cite{dir}, used by physicists. Ironically, this change came from Dirac himself, who advised his students and mankind to ``shut up and compute'', in anything related to quantum mechanics. Many mathematicians, including myself, followed his advice, shut up and computed, and concluded that linearity at left is better.

\bigskip

Back to mathematics now, in the above definition, the fact that $||x||=\sqrt{<x,x>}$ is indeed a norm comes from the Cauchy-Schwarz inequality, $|<x,y>|\leq||x||\cdot||y||$, which itself comes from the fact that the following degree 2 polynomial, with $t\in\mathbb R$ and $w\in\mathbb T$, being positive, its discriminant must be negative:
$$f(t)=||x+wty||^2$$

In finite dimensions, any algebraic basis $\{f_1,\ldots,f_N\}$ can be turned into an orthonormal basis $\{e_1,\ldots,e_N\}$, by using the Gram-Schmidt procedure. Thus, we have $H\simeq\mathbb C^N$, with this latter space being endowed with its usual scalar product, namely:
$$<x,y>=\sum_ix_i\bar{y}_i$$

The same happens in infinite dimensions, once again by Gram-Schmidt, coupled if needed with the Zorn lemma, in case our space is really very big. In other words, any Hilbert space has an orthonormal basis $\{e_i\}_{i\in I}$, and we have:
$$H\simeq l^2(I)$$

\index{separable Hilbert space}

Of particular interest is the ``separable'' case, where $I$ is countable. According to the above, there is up to isomorphism only one Hilbert space here, namely:
$$H=l^2(\mathbb N)$$

All this is, however, quite tricky, and can be a bit misleading. Consider for instance the space $H=L^2[0,1]$ of square-summable functions $f:[0,1]\to\mathbb C$, with:
$$<f,g>=\int_0^1f(x)\overline{g(x)}dx$$

This space is of course separable, because we can use the basis $f_n=x^n$ with $n\in\mathbb N$, orthogonalized by Gram-Schmidt. However, the orthogonalization procedure is something non-trivial, and so the isomorphism $H\simeq l^2(\mathbb N)$ that we obtain is something non-trivial as well. Doing some computations here is actually an excellent exercise.

\bigskip

In what follows we will be interested in the linear operators $T:H\to H$ which are bounded. Regarding such operators, we have the following result:

\index{linear operator}
\index{bounded operator}
\index{Banach algebra}
\index{adjoint operator}
\index{norm of operators}

\begin{theorem}
Given a Hilbert space $H$, the linear operators $T:H\to H$ which are bounded, in the sense that 
$$||T||=\sup_{||x||\leq1}||Tx||$$
is finite, form a complex algebra $B(H)$, having the following properties:
\begin{enumerate}
\item $B(H)$ is complete with respect to $||.||$, so we have a Banach algebra. 

\item $B(H)$ has an involution $T\to T^*$, given by $<Tx,y>=<x,T^*y>$.
\end{enumerate}
In addition, the norm and involution are related by the formula $||TT^*||=||T||^2$.
\end{theorem}

\begin{proof}
The fact that we have indeed an algebra follows from:
$$||S+T||\leq||S||+||T||$$
$$||\lambda T||=|\lambda|\cdot||T||$$
$$||ST||\leq||S||\cdot||T||$$

(1) Assuming that $\{T_n\}\subset B(H)$ is Cauchy then $\{T_nx\}$ is Cauchy for any $x\in H$, so we can define indeed the limit $T=\lim_{n\to\infty}T_n$ by setting:
$$Tx=\lim_{n\to\infty}T_nx$$

(2) Here the existence of $T^*$ comes from the fact that $\varphi(x)=<Tx,y>$ being a linear form $H\to\mathbb C$, we must have $\varphi(x)=<x,T^*y>$, for a certain vector $T^*y\in H$. Moreover, since this vector is unique, $T^*$ is unique too, and we have as well:
$$(S+T)^*=S^*+T^*\quad,\quad 
(\lambda T)^*=\bar{\lambda}T^*$$
$$(ST)^*=T^*S^*\quad,\quad 
(T^*)^*=T$$

Observe also that we have indeed $T^*\in B(H)$, because:
\begin{eqnarray*}
||T||
&=&\sup_{||x||=1}\sup_{||y||=1}<Tx,y>\\
&=&\sup_{||y||=1}\sup_{||x||=1}<x,T^*y>\\
&=&||T^*||
\end{eqnarray*}

Regarding now the last assertion, we have the following estimate:
$$||TT^*||
\leq||T||\cdot||T^*||
=||T||^2$$

On the other hand, we have as well the following estimate:
\begin{eqnarray*}
||T||^2
&=&\sup_{||x||=1}|<Tx,Tx>|\\
&=&\sup_{||x||=1}|<x,T^*Tx>|\\
&\leq&||T^*T||
\end{eqnarray*}

By replacing $T\to T^*$ we obtain from this $||T||^2\leq||TT^*||$, and we are done.
\end{proof}

Observe that when $H$ comes with an orthonormal basis $\{e_i\}_{i\in I}$, the linear map $T\to M$ given by $M_{ij}=<Te_j,e_i>$ produces an embedding as follows:
$$B(H)\subset M_I(\mathbb C)$$

Moreover, in this picture the operation $T\to T^*$ takes a very simple form, namely:
$$(M^*)_{ij}=\overline{M}_{ji}$$

However, as explained before Theorem 1.9, it is better in general not to use bases, and this because very simple spaces like $L^2[0,1]$ do not have simple bases.

\bigskip

The conditions found in Theorem 1.9 suggest the following definition:

\index{operator algebra}
\index{Banach algebra}

\begin{definition}
A $C^*$-algebra is a complex algebra $A$, having:
\begin{enumerate}
\item A norm $a\to||a||$, making it a Banach algebra.

\item An involution $a\to a^*$, satisfying $||aa^*||=||a||^2$.
\end{enumerate}
\end{definition}

Generally speaking, the elements $a\in A$ are best thought of as being some kind of ``generalized operators'', on some Hilbert space which is not present.  By using this idea, one can emulate spectral theory in this setting, as follows:

\index{spectrum}
\index{spectral radius}
\index{self-adjoint operator}
\index{normal operator}
\index{rational calculus}

\begin{proposition}
Given $a\in A$, define its spectrum as being the set
$$\sigma(a)=\left\{\lambda\in\mathbb C\Big|a-\lambda\not\in A^{-1}\right\}$$
and its spectral radius $\rho(a)$ as the radius of the smallest centered disk containing $\sigma(a)$.
\begin{enumerate}
\item The spectrum of a norm one element is in the unit disk.

\item The spectrum of a unitary element $(a^*=a^{-1}$) is on the unit circle. 

\item The spectrum of a self-adjoint element ($a=a^*$) consists of real numbers. 

\item The spectral radius of a normal element ($aa^*=a^*a$) is equal to its norm.
\end{enumerate}
\end{proposition}

\begin{proof}
Our first claim is that for any polynomial $f\in\mathbb C[X]$, and more generally for any rational function $f\in\mathbb C(X)$ having poles outside $\sigma(a)$, we have:
$$\sigma(f(a))=f(\sigma(a))$$

This indeed something well-known for the usual matrices. In the general case, assume first that we have a polynomial, $f\in\mathbb C[X]$. If we pick an arbitrary number $\lambda\in\mathbb C$, and write $f(X)-\lambda=c(X-r_1)\ldots(X-r_k)$, we have then, as desired:
\begin{eqnarray*}
\lambda\notin\sigma(f(a))
&\iff&f(a)-\lambda\in A^{-1}\\
&\iff&c(a-r_1)\ldots(a-r_k)\in A^{-1}\\
&\iff&a-r_1,\ldots,a-r_k\in A^{-1}\\
&\iff&r_1,\ldots,r_k\notin\sigma(a)\\
&\iff&\lambda\notin f(\sigma(a))
\end{eqnarray*}

Assume now that we are in the general case, $f\in\mathbb C(X)$. We pick $\lambda\in\mathbb C$, we write $f=P/Q$, and we set $F=P-\lambda Q$. By using the above finding, we obtain, as desired:
\begin{eqnarray*}
\lambda\in\sigma(f(a))
&\iff&F(a)\notin A^{-1}\\
&\iff&0\in\sigma(F(a))\\
&\iff&0\in F(\sigma(a))\\
&\iff&\exists\mu\in\sigma(a),F(\mu)=0\\
&\iff&\lambda\in f(\sigma(a))
\end{eqnarray*}

Regarding now the assertions in the statement, these basically follows from this:

\medskip

(1) This comes from the following formula, valid when $||a||<1$:
$$\frac{1}{1-a}=1+a+a^2+\ldots$$

(2) Assuming $a^*=a^{-1}$, we have the following norm computations:
$$||a||=\sqrt{||aa^*||}=\sqrt{1}=1$$
$$||a^{-1}||=||a^*||=||a||=1$$

If we denote by $D$ the unit disk, we obtain from this, by using (1):
$$||a||=1\implies\sigma(a)\subset D$$
$$||a^{-1}||=1\implies\sigma(a^{-1})\subset D$$

On the other hand, by using the rational function $f(z)=z^{-1}$, we have:
$$\sigma(a^{-1})\subset D\implies \sigma(a)\subset D^{-1}$$

Now by putting everything together we obtain, as desired:
$$\sigma(a)\subset D\cap D^{-1}=\mathbb T$$

(3) This follows from (2). Indeed, for $t>>0$ we have:
$$\left(\frac{a+it}{a-it}\right)^*
=\frac{a-it}{a+it}
=\left(\frac{a+it}{a-it}\right)^{-1}$$

Thus the element $f(a)$ is a unitary, and by using (2) its spectrum is contained in $\mathbb T$. We conclude from this that we have:
$$f(\sigma(a))=\sigma(f(a))\subset\mathbb T$$

But this shows that we have $\sigma(a)\subset f^{-1}(\mathbb T)=\mathbb R$, as desired.

\medskip

(4) We already know that we have $\rho(a)\leq ||a||$, for any $a\in A$. For the reverse inequality, when $a$ is normal, we fix a number $\rho>\rho(a)$. We have then:
\begin{eqnarray*}
\int_{|z|=\rho}\frac{z^n}{z -a}\,dz
&=&\int_{|z|=\rho}\sum_{k=0}^\infty z^{n-k-1}a^k\,dz\\
&=&\sum_{k=0}^\infty\left(\int_{|z|=\rho}z^{n-k-1}dz\right)a^k\\
&=&a^{n-1}
\end{eqnarray*}

By applying the norm and taking $n$-th roots we obtain from this formula, modulo some elementary manipulations, the following estimate: 
$$\rho\geq\lim_{n\to\infty}||a^n||^{1/n}$$

Now recall that $\rho$ was by definition an arbitrary number satisfying $\rho>\rho(a)$. Thus, we have obtained the following estimate, valid for any $a\in A$:
$$\rho(a)\geq\lim_{n\to\infty}||a^n||^{1/n}$$

In order to finish, we must prove that when $a$ is normal, this estimate implies the missing estimate, namely $\rho(a)\geq||a||$. We can proceed in two steps, as follows:

\medskip

\underline{Step 1}. In the case $a=a^*$ we have $||a^n||=||a||^n$ for any exponent of the form $n=2^k$, by using the $C^*$-algebra condition $||aa^*||=||a||^2$, and by taking $n$-th roots we get:
$$\rho(a)\geq||a||$$

Thus, we are done with the self-adjoint case, with the result $\rho(a)=||a||$.

\medskip

\underline{Step 2}. In the general normal case $aa^*=a^*a$ we have $a^n(a^n)^*=(aa^*)^n$, and by using this, along with the result from Step 1, applied to $aa^*$, we obtain:
\begin{eqnarray*}
\rho(a)
&\geq&\lim_{n\to\infty}||a^n||^{1/n}
=\sqrt{\lim_{n\to\infty}||a^n(a^n)^*||^{1/n}}\\
&=&\sqrt{\lim_{n\to\infty}||(aa^*)^n||^{1/n}}
=\sqrt{\rho(aa^*)}\\
&=&\sqrt{||a||^2}
=||a||
\end{eqnarray*}

Thus, we are led to the conclusion in the statement.
\end{proof}

We can now formulate a key theorem, as follows:

\index{Gelfand theorem}
\index{commutative algebra}

\begin{theorem}[Gelfand]
If $X$ is a compact space, the algebra $C(X)$ of continuous functions $f:X\to\mathbb C$ is a commutative $C^*$-algebra, with structure as follows:
\begin{enumerate}
\item The norm is the usual sup norm, $||f||=\sup_{x\in X}|f(x)|$.

\item The involution is the usual involution, $f^*(x)=\overline{f(x)}$.
\end{enumerate}
Conversely, any commutative $C^*$-algebra is of the form $C(X)$, with its ``spectrum'' $X=Spec(A)$ appearing as the space of characters $\chi :A\to\mathbb C$.
\end{theorem}

\begin{proof}
In what regards the first assertion, everything here is trivial. Conversely, given a commutative $C^*$-algebra $A$, we can define $X$ to be the set of characters $\chi :A\to\mathbb C$, with the topology making continuous all the evaluation maps $ev_a:\chi\to\chi(a)$. Then $X$ is a compact space, and $a\to ev_a$ is a morphism of algebras:
$$ev:A\to C(X)$$

Our first claim is that $ev$ is involutive. Indeed, we can use the following formula:
$$a=\frac{a+a^*}{2}-i\cdot\frac{i(a-a^*)}{2}$$

Thus it is enough to prove the equality $ev_{a^*}=ev_a^*$ for self-adjoint elements $a$. But this is the same as proving that $a=a^*$ implies that $ev_a$ is a real function, which is in turn true, because $ev_a(\chi)=\chi(a)$ is an element of $\sigma(a)$, contained in $\mathbb R$. Thus,  claim proved. Finally, since $A$ is commutative, each element is normal, so $ev$ is isometric:
$$||ev_a||=\rho(a)=||a||$$

It remains to prove that $ev$ is surjective. But this follows from the Stone-Weierstrass theorem, because $ev(A)$ is a closed subalgebra of $C(X)$, which separates the points.
\end{proof}

The Gelfand theorem suggests formulating the following definition:

\index{compact quantum space}
\index{quantum space}

\begin{definition}
Given a $C^*$-algebra $A$, not necessarily commutative, we write
$$A=C(X)$$
and call the abstract object $X$ a ``compact quantum space''.
\end{definition}

In other words, the category of compact quantum spaces is by definition the category of $C^*$-algebras, with the arrows reversed. We will be back to this, with examples, and with some technical comments as well, including a modification, the idea being that the above definition is in fact quite naive, and needs a fix. More on this later. 

\bigskip

Let us discuss now the other basic result regarding the $C^*$-algebras, namely the GNS representation theorem. We will need some more spectral theory, as follows:

\index{positive operator}
\index{square root of operators}

\begin{proposition}
For a normal element $a\in A$, the following are equivalent:
\begin{enumerate}
\item $a$ is positive, in the sense that $\sigma(a)\subset[0,\infty)$.

\item $a=b^2$, for some $b\in A$ satisfying $b=b^*$.

\item $a=cc^*$, for some $c\in A$.
\end{enumerate}
\end{proposition}

\begin{proof}
This is something very standard, as follows:

\medskip

$(1)\implies(2)$ Since our element $a$ is normal the algebra $<a>$ that is generates is commutative, and by using the Gelfand theorem, we can set $b=\sqrt{a}$. 

\medskip

$(2)\implies(3)$ This is trivial, because we can set $c=b$. 

\medskip

$(3)\implies(1)$ We proceed by contradiction. By multiplying $c$ by a suitable element of $<cc^*>$, we are led to the existence of an element $d\neq0$ satisfying $-dd^*\geq0$. By writing now $d=x+iy$ with $x=x^*,y=y^*$ we have:
$$dd^*+d^*d=2(x^2+y^2)\geq0$$

Thus $d^*d\geq0$. But this contradicts the elementary fact that $\sigma(dd^*),\sigma(d^*d)$ must coincide outside $\{0\}$, which can be established by explicit inversion.
\end{proof}

Here is now the GNS representation theorem, along with the idea of the proof:

\index{GNS theorem}
\index{GNS construction}

\begin{theorem}[GNS theorem]
Let $A$ be a $C^*$-algebra.
\begin{enumerate}
\item $A$ appears as a closed $*$-subalgebra $A\subset B(H)$, for some Hilbert space $H$. 

\item When $A$ is separable (usually the case), $H$ can be chosen to be separable.

\item When $A$ is finite dimensional, $H$ can be chosen to be finite dimensional. 
\end{enumerate}
\end{theorem}

\begin{proof}
Let us first discuss the commutative case, $A=C(X)$. Our claim here is that if we pick a probability measure on $X$, we have an embedding as follows:
$$C(X)\subset B(L^2(X))\quad,\quad
f\to(g\to fg)$$

Indeed, given a function $f\in C(X)$, consider the operator $T_f(g)=fg$, acting on $H=L^2(X)$. Observe that $T_f$ is indeed well-defined, and bounded as well, because:
$$||fg||_2
=\sqrt{\int_X|f(x)|^2|g(x)|^2dx}
\leq||f||_\infty||g||_2$$

Thus, $f\to T_f$ provides us with a $C^*$-algebra embedding $C(X)\subset B(H)$, as claimed. In general now, assuming that a linear form $\varphi:A\to\mathbb C$ has some suitable positivity properties, making it analogous to the integration functionals $\int_X:A\to\mathbb C$ from the commutative case, we can define a scalar product on $A$, by the following formula:
$$<a,b>=\varphi(b^*a)$$

By completing we obtain a Hilbert space $H$, and we have an embedding as follows:
$$A\subset B(H)\quad,\quad 
a\to(b\to ab)$$

Thus we obtain the assertion (1), and a careful examination of the construction $A\to H$, outlined above, shows that the assertions (2,3) are in fact proved as well.
\end{proof}

So long for operator theory and operator algebras. Obviously, some non-trivial things going on here, and although the above is basically all that we need, in what follows, more familiarity with all this would be desirable. The learning here starts with Rudin \cite{rud}, perfectly mastered, and then with some basic functional analysis, basic operator theory, and basic operator algebras, say from Lax \cite{lax}. For more, a good, useful, and especially modern book is Blackadar \cite{bla}. And for even more, go with Connes \cite{co1}.

\bigskip

Be said in passing, speaking Connes, with our $(S,T,U,K)$ philosophy we are already a bit away from what he does, because in his vision, the noncommutative Riemannian manifolds $X$ do not need coordinates, while in our vision, based on Nash \cite{nas}, they do. But this is only a slight difference, everything here being heavily inspired by \cite{co1}.

\section*{1c. Free spheres}

With the above formalism is hand, we can go ahead, and construct two free quadruplets $(S,T,U,K)$, in analogy with those corresponding to the classical real and complex geometries. Let us begin with the spheres. Following \cite{ba1}, \cite{bgo}, we have:

\index{free sphere}
\index{free real sphere}
\index{free complex sphere}

\begin{definition}
We have free real and complex spheres, defined via
$$C(S^{N-1}_{\mathbb R,+})=C^*\left(x_1,\ldots,x_N\Big|x_i=x_i^*,\sum_ix_i^2=1\right)$$
$$C(S^{N-1}_{\mathbb C,+})=C^*\left(x_1,\ldots,x_N\Big|\sum_ix_ix_i^*=\sum_ix_i^*x_i=1\right)$$
where the symbol $C^*$ stands for universal enveloping $C^*$-algebra.
\end{definition}

All this deserves some explanations. Given an integer $N\in\mathbb N$, consider the free complex algebra on $2N$ variables, denoted $x_1,\ldots,x_N$ and $x_1^*,\ldots,x_N^*$:
$$A=\Big<x_1,\ldots,x_N,x_1^*,\ldots,x_N^*\Big>$$

This algebra has an involution $*:A\to A$, given by $x_i\leftrightarrow x_i^*$. Now let us consider the following $*$-algebra quotients of our $*$-algebra $A$:
\begin{eqnarray*}
A_R&=&A\Big/\Big<x_i=x_i^*,\sum_ix_i^2=1\Big>\\
A_C&=&A\Big/\Big<\sum_ix_ix_i^*=\sum_ix_i^*x_i=1\Big>
\end{eqnarray*}

Since the first relations imply the second ones, we have quotient maps as follows:
$$A\to A_C\to A_R$$

Our claim now is both $A_C,A_R$ admit enveloping $C^*$-algebras, in the sense that the biggest $C^*$-norms on these $*$-algebras are bounded. We only have to check this for the bigger algebra $A_C$. But here, our claim follows from the following estimate:
$$||x_i||^2
=||x_ix_i^*||
\leq\Big|\Big|\sum_ix_ix_i^*\Big|\Big|
=1$$

Summarizing, our claim is proved, so we can define $C(S^{N-1}_{\mathbb R,+}),C(S^{N-1}_{\mathbb C,+})$ as being the enveloping $C^*$-algebras of $A_R,A_C$, and so Definition 1.16 makes sense. 

\bigskip

In order to formulate now some results, let us introduce as well:

\index{classical version}
\index{liberation}

\begin{definition}
Given a compact quantum space $X$, its classical version is the usual compact space $X_{class}\subset X$ obtained by dividing $C(X)$ by its commutator ideal:
$$C(X_{class})=C(X)/I\quad,\quad 
I=<[a,b]>$$
In this situation, we also say that $X$ appears as a ``liberation'' of $X$.
\end{definition}

In other words, the space $X_{class}$ appears as the Gelfand spectrum of the commutative $C^*$-algebra $C(X)/I$. Observe in particular that $X_{class}$ is indeed a classical space. As a first result now, regarding the above free spheres, we have:

\begin{theorem}
We have embeddings of compact quantum spaces, as follows,
$$\xymatrix@R=15mm@C=15mm{
S^{N-1}_{\mathbb R,+}\ar[r]&S^{N-1}_{\mathbb C,+}\\
S^{N-1}_\mathbb R\ar[r]\ar[u]&S^{N-1}_\mathbb C\ar[u]
}$$
and the spaces on top appear as liberations of the spaces on the bottom.
\end{theorem}

\begin{proof}
The first assertion, regarding the inclusions, comes from the fact that at the level of the associated $C^*$-algebras, we have surjective maps, as follows:
$$\xymatrix@R=15mm@C=15mm{
C(S^{N-1}_{\mathbb R,+})\ar[d]&C(S^{N-1}_{\mathbb C,+})\ar[d]\ar[l]\\
C(S^{N-1}_\mathbb R)&C(S^{N-1}_\mathbb C)\ar[l]
}$$

For the second assertion, we must establish the following isomorphisms, where the symbol $C^*_{comm}$ stands for ``universal commutative $C^*$-algebra generated by'':
$$C(S^{N-1}_\mathbb R)=C^*_{comm}\left(x_1,\ldots,x_N\Big|x_i=x_i^*,\sum_ix_i^2=1\right)$$
$$C(S^{N-1}_\mathbb C)=C^*_{comm}\left(x_1,\ldots,x_N\Big|\sum_ix_ix_i^*=\sum_ix_i^*x_i=1\right)$$

As a  first observation, it is enough to establish the second isomorphism, because the first one will follow from it, simply by dividing by the relations $x_i=x_i^*$. So, consider the second universal commutative $C^*$-algebra $A$ constructed above. Since the standard coordinates on $S^{N-1}_\mathbb C$ satisfy the defining relations for $A$, we have a quotient map of as follows, mapping standard coordinates to standard coordinates:
$$A\to C(S^{N-1}_\mathbb C)$$

Conversely, let us write $A=C(S)$, by using the Gelfand theorem. The variables $x_1,\ldots,x_N$ become in this way true coordinates, providing us with an embedding $S\subset\mathbb C^N$. Also, the quadratic relations become $\sum_i|x_i|^2=1$, so we have $S\subset S^{N-1}_\mathbb C$. Thus, we have a quotient map $C(S^{N-1}_\mathbb C)\to A$, as desired, and this gives all the results.
\end{proof}

Summarizing, we are done with the spheres. Before getting into tori, let us talk about algebraic manifolds. By using the free spheres constructed above, we can formulate:

\index{algebraic manifold}
\index{real algebraic manifold}

\begin{definition}
A real algebraic manifold $X\subset S^{N-1}_{\mathbb C,+}$ is a closed quantum subspace defined, at the level of the corresponding $C^*$-algebra, by a formula of type
$$C(X)=C(S^{N-1}_{\mathbb C,+})\Big/\Big<f_i(x_1,\ldots,x_N)=0\Big>$$
for certain family of noncommutative polynomials, as follows:
$$f_i\in\mathbb C<x_1,\ldots,x_N>$$
We denote by $\mathcal C(X)$ the $*$-subalgebra of $C(X)$ generated by the coordinates $x_1,\ldots,x_N$. 
\end{definition}

As a basic example here, we have the free real sphere $S^{N-1}_{\mathbb R,+}$. The classical spheres $S^{N-1}_\mathbb C,S^{N-1}_\mathbb R$, and their real submanifolds, are covered as well by this formalism. At the level of the general theory, we have the following version of the Gelfand theorem:

\index{liberation}

\begin{theorem}
If $X\subset S^{N-1}_{\mathbb C,+}$ is an algebraic manifold, as above, we have
$$X_{class}=\left\{x\in S^{N-1}_\mathbb C\Big|f_i(x_1,\ldots,x_N)=0\right\}$$
and $X$ appears as a liberation of $X_{class}$.
\end{theorem}

\begin{proof}
This is something that we already met, in the context of the free spheres. In general, the proof is similar, by using the Gelfand theorem. Indeed, if we denote by $X_{class}'$ the manifold constructed in the statement, then we have a quotient map of $C^*$-algebras as follows, mapping standard coordinates to standard coordinates:
$$C(X_{class})\to C(X_{class}')$$

Conversely now, from $X\subset S^{N-1}_{\mathbb C,+}$ we obtain $X_{class}\subset S^{N-1}_\mathbb C$. Now since the relations defining $X_{class}'$ are satisfied by $X_{class}$, we obtain an inclusion $X_{class}\subset X_{class}'$. Thus, at the level of algebras of continuous functions, we have a quotient map of $C^*$-algebras as follows, mapping standard coordinates to standard coordinates:
$$C(X_{class}')\to C(X_{class})$$

Thus, we have constructed a pair of inverse morphisms, and we are done.
\end{proof}

Finally, once again at the level of the general theory, we have:

\index{equality of manifolds}
\index{isomorphism of manifolds}

\begin{definition}
We agree to identify two real algebraic submanifolds $X,Y\subset S^{N-1}_{\mathbb C,+}$ when we have a $*$-algebra isomorphism between $*$-algebras of coordinates
$$f:\mathcal C(Y)\to\mathcal C(X)$$
mapping standard coordinates to standard coordinates.
\end{definition}

We will see later the reasons for making this convention, coming from amenability.

\section*{1d. Free tori}

Let us go back now to our general $(S,T,U,K)$ program. Now that we are done with the free spheres, we can introduce as well free tori, as follows:

\index{free torus}
\index{free real torus}
\index{free complex torus}

\begin{definition}
We have free real and complex tori, defined via
$$C(T_N^+)=C^*\left(x_1,\ldots,x_N\Big|x_i=x_i^*,x_i^2=\frac{1}{N}\right)$$
$$C(\mathbb T_N^+)=C^*\left(x_1,\ldots,x_N\Big|x_ix_i^*=x_i^*x_i=\frac{1}{N}\right)$$
with the symbol $C^*$ standing as usual for universal enveloping $C^*$-algebra.
\end{definition}

The fact that these tori are indeed well-defined comes from the fact that they are algebraic manifolds, in the sense of Definition 1.19. In fact, we have:

\begin{proposition}
We have inclusions of algebraic manifolds, as follows:
$$\xymatrix@R=14mm@C=14mm{
S^{N-1}_{\mathbb R,+}\ar[r]&S^{N-1}_{\mathbb C,+}\\
T_N^+\ar[r]\ar[u]&\mathbb T_N^+\ar[u]
}$$
In addition, this is an intersection diagram, in the sense that $T_N^+=\mathbb T_N^+\cap S^{N-1}_{\mathbb R,+}$.
\end{proposition}

\begin{proof}
All this is clear indeed, by using the equivalence relation in Definition 1.21, in order to get rid of functional analytic issues at the $C^*$-algebra level.
\end{proof}

In analogy with Theorem 1.18, we have the following result:

\begin{theorem}
We have inclusions of algebraic manifolds, as follows,
$$\xymatrix@R=15mm@C=15mm{
T_N^+\ar[r]&\mathbb T_N^+\\
T_N\ar[r]\ar[u]&\mathbb T_N\ar[u]
}$$
and the manifolds on top appear as liberations of those of the bottom.
\end{theorem}

\begin{proof}
This follows exactly as Theorem 1.18, and the best here is in fact to invoke Theorem 1.20, which is there precisely for dealing with such situations.
\end{proof}

Summarizing, we have free spheres and tori, having quite similar properties. Let us further study the tori. Up to a rescaling, these are given by algebras generated by unitaries, so studying the algebras generated by unitaries will be our next task. The point is that we have many such algebras, coming from the following construction:

\index{group algebra}
\index{full group algebra}
\index{Pontrjagin dual}

\begin{theorem}
Let $\Gamma$ be a discrete group, and consider the complex group algebra $\mathbb C[\Gamma]$, with involution given by the fact that all group elements are unitaries, $g^*=g^{-1}$.
\begin{enumerate}
\item The maximal $C^*$-seminorm on $\mathbb C[\Gamma]$ is a $C^*$-norm, and the closure of $\mathbb C[\Gamma]$ with respect to this norm is a $C^*$-algebra, denoted $C^*(\Gamma)$.

\item When $\Gamma$ is abelian, we have an isomorphism $C^*(\Gamma)\simeq C(G)$, where $G=\widehat{\Gamma}$ is its Pontrjagin dual, formed by the characters $\chi:\Gamma\to\mathbb T$.
\end{enumerate}
\end{theorem}

\begin{proof}
All this is very standard, the idea being as follows:

\medskip

(1) In order to prove the result, we must find a $*$-algebra embedding $\mathbb C[\Gamma]\subset B(H)$, with $H$ being a Hilbert space. For this purpose, consider the space $H=l^2(\Gamma)$, having $\{h\}_{h\in\Gamma}$ as orthonormal basis. Our claim is that we have an embedding, as follows:
$$\pi:\mathbb C[\Gamma]\subset B(H)\quad,\quad 
\pi(g)(h)=gh$$

Indeed, since $\pi(g)$ maps the basis $\{h\}_{h\in\Gamma}$ into itself, this operator is well-defined, bounded, and is an isometry. It is also clear from the formula $\pi(g)(h)=gh$ that $g\to\pi(g)$ is a morphism of algebras, and since this morphism maps the unitaries $g\in\Gamma$ into isometries, this is a morphism of $*$-algebras. Finally, the faithfulness of $\pi$ is clear.

\medskip

(2) Since $\Gamma$ is abelian, the corresponding group algebra $A=C^*(\Gamma)$ is commutative. Thus, we can apply the Gelfand theorem, and we obtain $A=C(X)$, with:
$$X=Spec(A)$$

But the spectrum $X=Spec(A)$, consisting of the characters $\chi:C^*(\Gamma)\to\mathbb C$, can be identified with the Pontrjagin dual $G=\widehat{\Gamma}$, and this gives the result.
\end{proof}

The above result suggests the following definition:

\begin{definition}
Given a discrete group $\Gamma$, the compact quantum space $G$ given by
$$C(G)=C^*(\Gamma)$$
is called abstract dual of $\Gamma$, and is denoted $G=\widehat{\Gamma}$.
\end{definition}

\index{group dual}

This is in fact something which is not very satisfactory, in general, due to amenability issues. However, in the case of the finitely generated discrete groups $\Gamma=<g_1,\ldots,g_N>$, which is the one that we are interested in here, the corresponding duals appear as algebraic submanifolds $\widehat{\Gamma}\subset S^{N-1}_{\mathbb C,+}$, and the notion of equivalence from Definition 1.21 is precisely the one that we need, identifying full and reduced group algebras.

\bigskip

We can now refine our findings about tori, as follows:

\begin{theorem}
The basic tori are all group duals, as follows,
$$\xymatrix@R=16.2mm@C=16.2mm{
T_N^+\ar[r]&\mathbb T_N^+\\
T_N\ar[r]\ar[u]&\mathbb T_N\ar[u]
}
\qquad
\xymatrix@R=8mm@C=15mm{\\ =}
\qquad
\xymatrix@R=15mm@C=15mm{
\widehat{\mathbb Z_2^{*N}}\ar[r]&\widehat{F_N}\\
\mathbb Z_2^N\ar[r]\ar[u]&\mathbb T^N\ar[u]
}$$
where $F_N$ is the free group on $N$ generators, and $*$ is a group-theoretical free product.
\end{theorem}

\begin{proof}
The basic tori appear indeed as group duals, and together with the Fourier transform identifications from Theorem 1.25 (2), this gives the result.
\end{proof}

Let us try now to understand the correspondence between the spheres $S$ and tori $T$. We first have the following result, summarizing our knowledge so far:

\begin{theorem}
The four main quantum spheres produce the main quantum tori
$$\xymatrix@R=15mm@C=15mm{
S^{N-1}_{\mathbb R,+}\ar[r]&S^{N-1}_{\mathbb C,+}\\
S^{N-1}_\mathbb R\ar[r]\ar[u]&S^{N-1}_\mathbb C\ar[u]
}\qquad
\xymatrix@R=8mm@C=15mm{\\ \to}
\qquad
\xymatrix@R=16mm@C=16mm{
T_N^+\ar[r]&\mathbb T_N^+\\
T_N\ar[r]\ar[u]&\mathbb T_N\ar[u]
}$$
via the formula $T=S\cap\mathbb T_N^+$, with the intersection being taken inside $S^{N-1}_{\mathbb C,+}$.
\end{theorem}

\begin{proof}
This comes from the above results, the situation being as follows:

\medskip

(1) Free complex case. Here the formula in the statement reads $\mathbb T_N^+=S^{N-1}_{\mathbb C,+}\cap\mathbb T_N^+$. But this is something trivial, because we have $\mathbb T_N^+\subset S^{N-1}_{\mathbb C,+}$.

\medskip

(2) Free real case. Here the formula in the statement reads $T_N^+=S^{N-1}_{\mathbb R,+}\cap\mathbb T_N^+$. But this is something that we already know, from Proposition 1.23.

\medskip

(3) Classical complex case. Here the formula in the statement reads $\mathbb T_N=S^{N-1}_\mathbb C\cap\mathbb T_N^+$. But this is clear as well, the classical version of $\mathbb T_N^+$ being $\mathbb T_N$.

\medskip

(4) Classical real case. Here the formula in the statement reads $T_N=S^{N-1}_\mathbb R\cap\mathbb T_N^+$. But this follows by intersecting the formulae from the proof of (2) and (3).
\end{proof}

Importantly, the correspondence $S\to T$ found above is not the only one. In order to discuss this, let us start with a general result, as follows:

\index{toral isometry}

\begin{theorem}
Given an algebraic manifold $X\subset S^{N-1}_{\mathbb C,+}$, the category of tori $T\subset\mathbb T_N^+$ acting affinely on $X$, in the sense that we have a morphism of algebras as follows,
$$\Phi:C(X)\to C(X)\otimes C(T)\quad,\quad x_i\to x_i\otimes g_i$$
has a universal object, denoted $T^+(X)$, and called toral isometry group of $X$.
\end{theorem}

\begin{proof}
This is something a bit advanced, and we will talk more about affine actions, with full details, in chapter 3 below. This being said, our theorem as stated formally makes sense, so let us prove it. Assume that $X\subset S^{N-1}_{\mathbb C,+}$ comes as follows:
$$C(X)=C(S^{N-1}_{\mathbb C,+})\Big/\Big<f_\alpha(x_1,\ldots,x_N)=0\Big>$$

Consider now the following variables:
$$X_i=x_i\otimes g_i\in C(X)\otimes C(\mathbb T_N^+)$$

Our claim is that the torus $T=T^+(X)$ in the statement appears as follows:
$$C(T)=C(\mathbb T_N^+)\Big/\Big<f_\alpha(X_1,\ldots,X_N)=0\Big>$$

In order to prove this claim, we have to clarify how the relations $f_\alpha(X_1,\ldots,X_N)=0$ are interpreted inside $C(\mathbb T_N^+)$, and then show that $T$ is indeed a toral subgroup. So, pick one of the defining polynomials, $f=f_\alpha$, and write it as follows:
$$f(x_1,\ldots,x_N)=\sum_r\sum_{i_1^r\ldots i_{s_r}^r}\lambda_r\cdot x_{i_1^r}\ldots x_{i_{s_r}^r}$$

With $X_i=x_i\otimes g_i$ as above, we have the following formula:
$$f(X_1,\ldots,X_N)=\sum_r\sum_{i_1^r\ldots i_{s_r}^r}\lambda_rx_{i_1^r}\ldots x_{i_{s_r}^r}\otimes g_{i_1^r}\ldots g_{i_{s_r}^r}$$

Since the variables on the right span a certain finite dimensional space, the relations $f(X_1,\ldots,X_N)=0$ correspond to certain relations between the variables $g_i$. Thus, we have indeed a subspace $T\subset\mathbb T_N^+$, with a universal map, as follows:
$$\Phi:C(X)\to C(X)\otimes C(T)$$

In order to show now that $T$ is a group dual, consider the following elements:
$$g_i'=g_i\otimes g_i\quad,\quad X_i'=x_i\otimes g_i'$$

Then from $f(X_1,\ldots,X_N)=0$ we deduce that, with $\Delta(g)=g\otimes g$, we have:
$$f(X_1',\ldots,X_N')
=(id\otimes\Delta)f(X_1,\ldots,X_N)
=0$$

Thus we can map $g_i\to g_i'$, and it follows that $T$ is a group dual, as desired.
\end{proof}

We can now formulate a second result relating spheres and tori, as follows:

\begin{theorem}
The four main quantum spheres produce via
$$T=T^+(S)$$
the corresponding four main quantum tori.
\end{theorem}

\begin{proof}
This is something elementary, which can be established as follows:

\medskip

(1) Free complex case. Here is there is nothing to be proved, because we obviously have an action $\mathbb T_N^+\curvearrowright S^{N-1}_{\mathbb C,+}$, and this action can only be universal.

\medskip

(2) Free real case. Here the situation is similar, because we have an obvious action $T_N^+\curvearrowright S^{N-1}_{\mathbb R,+}$, and it is clear that this action can only be universal.

\medskip

(3) Classical complex case. Once again, we have a similar situation here, with the obvious action, namely $\mathbb T_N\curvearrowright S^{N-1}_\mathbb C$, being easily seen to be universal.

\medskip

(4) Classical real case. Here the obvious action, namely $T_N\curvearrowright S^{N-1}_\mathbb R$, is universal as well, the reasons for this coming from (2) and (3) above.
\end{proof}

As a conclusion now, following \cite{bb2} and related papers, we can formulate:

\index{baby noncommutative geometry}

\begin{definition}
A ``baby noncommutative geometry'' consists of a quantum sphere $S$ and a quantum torus $T$, which are by definition algebraic manifolds as follows,
$$S^{N-1}_\mathbb R\subset S\subset S^{N-1}_{\mathbb C+}$$
$$T_N\subset T\subset\mathbb T_N^+$$
which must be subject to the following compatibility conditions,
$$T=S\cap\mathbb T_N^+=T^+(S)$$
with the intersection being taken inside $S^{N-1}_{\mathbb C,+}$, and $T^+$ being the toral isometry group.
\end{definition}

With this notion in hand, our main results so far can be summarized as follows:

\begin{theorem}
We have $4$ baby noncommutative geometries, as follows,
$$\xymatrix@R=50pt@C=50pt{
\mathbb R^N_+\ar[r]&\mathbb C^N_+\\
\mathbb R^N\ar[u]\ar[r]&\mathbb C^N\ar[u]
}$$ 
with each symbol $\mathbb K^N_\times$ standing for the corresponding pair $(S,T)$.
\end{theorem}

\begin{proof}
This follows indeed from Theorem 1.28 and Theorem 1.30.
\end{proof}

In what follows we will extend our baby theory, with pairs of type $(U,K)$, consisting of unitary and reflection groups. This will lead to a theory which is more advanced.

\section*{1e. Exercises}

Generally speaking, at this point, more reading of mathematics and physics would be welcome. Here is however one concrete exercise, that you should definitely try:

\begin{exercise}
Establish correspondences as follows,
$$\xymatrix@R=50pt@C=50pt{
S^{N-1}_\mathbb R\ar[r]\ar[d]\ar[dr]&T_N\ar[l]\ar[d]\ar[dl]\\
O_N\ar[u]\ar[ur]\ar[r]&H_N\ar[l]\ar[ul]\ar[u]
}$$
by using results from basic geometry and group theory.
\end{exercise}

We have already talked about this, in the beginning of this chapter, and the problem now is to have the thing done. For a bonus point, do the complex case too.

\chapter{Quantum groups}

\section*{2a. Quantum groups}

We have seen so far that the pairs sphere/torus $(S,T)$ corresponding to the real and complex geometries, of $\mathbb R^N,\mathbb C^N$, have some natural free analogues. Our objective now will be that of adding to the picture a pair of quantum groups $(U,K)$, as to reach to a quadruplet of objects $(S,T,U,K)$, with relations between them, as follows:
$$
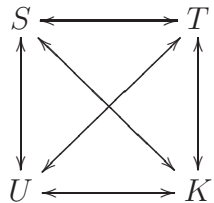
\xymatrix@R=50pt@C=50pt{
S\ar[r]\ar[d]\ar[dr]&T\ar[l]\ar[d]\ar[dl]\\
U\ar[u]\ar[ur]\ar[r]&K\ar[l]\ar[ul]\ar[u]
}$$

Before starting, some philosophical comments. You might argue that the pairs $(S,T)$ that we have look just fine, so why embarking into quantum groups, and complicating our theory with objects $(U,K)$. This is a reasonable criticism, and in answer:

\bigskip

(1) First of all, there is no sphere $S$ without corresponding rotation group $U$. With this meaning that, no matter what you want to do with $S$, of reasonably advanced type, like integrating over it, or looking at its Laplacian, and so on, you will certainly run into $U$. And for similar reasons, a bit more complicated, there is no $T$ without $K$ either.

\bigskip

(2) This being said, you will say, why not further studying $S,T$, say from a differential geometry viewpoint, and leaving $U,K$ for later. Well, this does not work. The problem is that $S,T$, at least in the free case, that we are very interested in here, while having a Laplacian, do not have a Dirac operator in the sense of Connes \cite{co1}. 

\bigskip

(3) In short, such ideas will not work, and we are led into $U,K$. By the way, meditating a bit about noncommutative differential geometry, at this point, is something recommended. And that you will have to do by yourself, the no-go results here being folklore. The needed read here is Connes \cite{co1}, with Blackadar \cite{bla} helping.

\bigskip

So, quantum groups. We will spend quite some time in introducing them, and working out their properties, and with this long series of things to be learned being good news, because the more theory we have about quantum groups, the more techniques we will have for dealing with $(U,K)$, and so with the whole quadruplets $(S,T,U,K)$.

\bigskip

The formalism that we need, coming from Woronowicz \cite{wo1}, is as follows:

\index{Woronowicz algebra}
\index{comultiplication}
\index{counit}
\index{antipode}

\begin{definition}
A Woronowicz algebra is a $C^*$-algebra $A$, given with a unitary matrix $u\in M_N(A)$ whose coefficients generate $A$, such that the formulae
$$\Delta(u_{ij})=\sum_ku_{ik}\otimes u_{kj}$$
$$\varepsilon(u_{ij})=\delta_{ij}$$
$$S(u_{ij})=u_{ji}^*$$
define morphisms of $C^*$-algebras as follows,
$$\Delta:A\to A\otimes A$$
$$\varepsilon:A\to\mathbb C$$
$$S:A\to A^{opp}$$
called comultiplication, counit and antipode.
\end{definition}

Obviously, this is something tricky, and we will see details in a moment, the idea being that these are the axioms which best fit with what we want to do, in this book. Let us also mention, technically, that $\otimes$ in the above can be any $C^*$-algebra tensor product, and with the precise choice being irrelevant, but more on this later. Also, $A^{opp}$ is the opposite algebra, with multiplication $a\cdot b=ba$, and more on this later too.

\bigskip

We say that $A$ is cocommutative when $\Sigma\Delta=\Delta$, where $\Sigma(a\otimes b)=b\otimes a$ is the flip. With this convention, we have the following key result, from Woronowicz \cite{wo1}:

\index{cocommutative algebra}
\index{commutative algebra}

\begin{theorem}
The following are Woronowicz algebras:
\begin{enumerate}
\item $C(G)$, with $G\subset U_N$ compact Lie group. Here the structural maps are:
$$\Delta(\varphi)=[(g,h)\to \varphi(gh)]$$
$$\varepsilon(\varphi)=\varphi(1)$$
$$S(\varphi)=[g\to\varphi(g^{-1})]$$

\item $C^*(\Gamma)$, with $F_N\to\Gamma$ finitely generated group. Here the structural maps are:
$$\Delta(g)=g\otimes g$$
$$\varepsilon(g)=1$$
$$S(g)=g^{-1}$$

\end{enumerate}
Moreover, we obtain in this way all the commutative/cocommutative algebras.
\end{theorem}

\begin{proof}
In both cases, we have to exhibit a certain matrix $u$:

\medskip

(1) Here we can use the matrix $u=(u_{ij})$ formed by matrix coordinates of $G$:
$$g=\begin{pmatrix}
u_{11}(g)&\ldots&u_{1N}(g)\\
\vdots&&\vdots\\
u_{N1}(g)&\ldots&u_{NN}(g)
\end{pmatrix}$$

(2) Here we can use the diagonal matrix formed by generators of $\Gamma$:
$$u=\begin{pmatrix}
g_1&&0\\
&\ddots&\\
0&&g_N
\end{pmatrix}$$

Finally, the last assertion follows from the Gelfand theorem, in the commutative case. In the cocommutative case, this is something more technical, to be discussed later.
\end{proof}

In general now, the structural maps $\Delta,\varepsilon,S$ have the following properties:

\index{square of antipode}

\begin{proposition}
Let $(A,u)$ be a Woronowicz algebra.
\begin{enumerate} 
\item $\Delta,\varepsilon$ satisfy the usual axioms for a comultiplication and a counit, namely:
$$(\Delta\otimes id)\Delta=(id\otimes \Delta)\Delta$$
$$(\varepsilon\otimes id)\Delta=(id\otimes\varepsilon)\Delta=id$$

\item $S$ satisfies the antipode axiom, on the $*$-subalgebra generated by entries of $u$: 
$$m(S\otimes id)\Delta=m(id\otimes S)\Delta=\varepsilon(.)1$$

\item In addition, the square of the antipode is the identity, $S^2=id$.
\end{enumerate}
\end{proposition}

\begin{proof}
Observe first that the result holds in the case where $A$ is commutative. Indeed, by using Theorem 2.2 (1) we can write:
$$\Delta=m^t$$
$$\varepsilon=u^t$$
$$S=i^t$$

The 3 conditions in the statement come then by transposition from the basic 3 group theory conditions satisfied by $m,u,i$, namely:
$$m(m\times id)=m(id\times m)$$
$$m(id\times u)=m(u\times id)=id$$
$$m(id\times i)\delta=m(i\times id)\delta=1$$

Here $\delta(g)=(g,g)$. Observe also that the last condition, $S^2=id$, is satisfied as well, coming from the identity $i^2=id$, which is a consequence of the group axioms.

\medskip

Observe also that the result holds as well in the case where $A$ is cocommutative, by using Theorem 2.2 (1). Indeed, the 3 formulae in the statement are all trivial, and the condition $S^2=id$ follows once again from the group theory formula $(g^{-1})^{-1}=g$.

\medskip

In the general case now, the proof goes as follows:

\medskip

(1) We have the following computation:
$$(\Delta\otimes id)\Delta(u_{ij})
=\sum_l\Delta(u_{il})\otimes u_{lj}
=\sum_{kl}u_{ik}\otimes u_{kl}\otimes u_{lj}$$

We have as well the following computation, which gives the first formula:
$$(id\otimes\Delta)\Delta(u_{ij})
=\sum_ku_{ik}\otimes\Delta(u_{kj})
=\sum_{kl}u_{ik}\otimes u_{kl}\otimes u_{lj}$$

On the other hand, we have the following computation:
$$(id\otimes\varepsilon)\Delta(u_{ij})
=\sum_ku_{ik}\otimes\varepsilon(u_{kj})
=u_{ij}$$

We have as well the following computation, which gives the second formula:
$$(\varepsilon\otimes id)\Delta(u_{ij})
=\sum_k\varepsilon(u_{ik})\otimes u_{kj}
=u_{ij}$$

(2) By using the fact that the matrix $u=(u_{ij})$ is unitary, we obtain:
\begin{eqnarray*}
m(id\otimes S)\Delta(u_{ij})
&=&\sum_ku_{ik}S(u_{kj})\\
&=&\sum_ku_{ik}u_{jk}^*\\
&=&(uu^*)_{ij}\\
&=&\delta_{ij}
\end{eqnarray*}

We have as well the following computation, which gives the result:
\begin{eqnarray*}
m(S\otimes id)\Delta(u_{ij})
&=&\sum_kS(u_{ik})u_{kj}\\
&=&\sum_ku_{ki}^*u_{kj}\\
&=&(u^*u)_{ij}\\
&=&\delta_{ij}
\end{eqnarray*}

(3) Finally, the formula $S^2=id$ holds as well on the generators, and we are done.
\end{proof}

Let us record as well the following technical result:

\index{biunitary}

\begin{proposition}
Given a Woronowicz algebra $(A,u)$, we have $u^t=\bar{u}^{-1}$, so $u$ is biunitary, in the sense that it is unitary, with unitary transpose.
\end{proposition}

\begin{proof}
We have the following computation, based on the fact that $u$ is unitary:
\begin{eqnarray*}
(uu^*)_{ij}=\delta_{ij}
&\implies&\sum_kS(u_{ik}u_{jk}^*)=\delta_{ij}\\
&\implies&\sum_ku_{kj}u_{ki}^*=\delta_{ij}\\
&\implies&(u^t\bar{u})_{ji}=\delta_{ij}
\end{eqnarray*}

Similarly, we have the following computation, once agan using the unitarity of $u$:
\begin{eqnarray*}
(u^*u)_{ij}=\delta_{ij}
&\implies&\sum_kS(u_{ki}^*u_{kj})=\delta_{ij}\\
&\implies&\sum_ku_{jk}^*u_{ik}=\delta_{ij}\\
&\implies&(\bar{u}u^t)_{ji}=\delta_{ij}
\end{eqnarray*}

Thus, we are led to the conclusion in the statement.
\end{proof}

Summarizing, the Woronowicz algebras appear to have nice properties. In view of Theorem 2.2 and of Proposition 2.3, we can formulate the following definition:

\index{quantum group}
\index{compact quantum group}
\index{discrete quantum group}
\index{Pontrjagin dual}
\index{group algebra}

\begin{definition}
Given a Woronowicz algebra $A$, we formally write
$$A=C(G)=C^*(\Gamma)$$
and call $G$ compact quantum group, and $\Gamma$ discrete quantum group.
\end{definition}

When $A$ is commutative and cocommutative, $G$ and $\Gamma$ are usual abelian groups, dual to each other. In general, we still agree to write $G=\widehat{\Gamma},\Gamma=\widehat{G}$, but in a formal sense. As a final piece of general theory now, let us complement Definition 2.1 with:

\begin{definition}
Given two Woronowicz algebras $(A,u)$ and $(B,v)$, we write 
$$A\simeq B$$
and identify the corresponding quantum groups, when we have an isomorphism 
$$<u_{ij}>\simeq<v_{ij}>$$
of $*$-algebras, mapping standard coordinates to standard coordinates.
\end{definition}

With this convention, which is in tune with our conventions for algebraic manifolds from chapter 1, and more on this later, any compact or discrete quantum group corresponds to a unique Woronowicz algebra, up to equivalence. Also, we can see now why in Definition 2.1 the choice of the exact $C^*$-algebra tensor product $\otimes$ is irrelevant. Indeed, no matter what tensor product $\otimes$ we use there, we end up with the same Woronowicz algebra, and the same compact and discrete quantum groups, up to equivalence. 

\bigskip

In practice, we will use in what follows the simplest such tensor product $\otimes$, which is the maximal one, obtained as completion of the usual algebraic tensor product with respect to the biggest $C^*$-norm. With the remark that this product is something rather abstract, and so can be treated, in practice, as a usual algebraic tensor product.

\bigskip

Moving ahead now, let us call corepresentation of $A$ any unitary matrix $v\in M_n(\mathcal A)$, where $\mathcal A=<u_{ij}>$, satisfying the same conditions are those satisfied by $u$, namely:
$$\Delta(v_{ij})=\sum_kv_{ik}\otimes v_{kj}$$
$$\varepsilon(v_{ij})=\delta_{ij}$$
$$S(v_{ij})=v_{ji}^*$$

\index{representation}
\index{corepresentation}

These corepresentations can be thought of as corresponding to the finite dimensional unitary smooth representations of the underlying compact quantum group $G$. Following Woronowicz \cite{wo1}, we have the following key result:

\index{Haar integration}
\index{Ces\`aro limit}

\begin{theorem}
Any Woronowicz algebra $A=C(G)$ has a Haar integration functional, 
$$\left(\int_G\otimes id\right)\Delta=\left(id\otimes\int_G\right)\Delta=\int_G(.)1$$
which can be constructed by starting with any faithful positive form $\varphi\in A^*$, and setting
$$\int_G=\lim_{n\to\infty}\frac{1}{n}\sum_{k=1}^n\varphi^{*k}$$
where $\phi*\psi=(\phi\otimes\psi)\Delta$. Moreover, for any corepresentation $v\in M_n(\mathbb C)\otimes A$ we have
$$\left(id\otimes\int_G\right)v=P$$
where $P$ is the orthogonal projection onto $Fix(v)=\big\{\xi\in\mathbb C^n\big|v\xi=\xi\big\}$.
\end{theorem}

\begin{proof}
Following \cite{wo1}, this can be done in 3 steps, as follows:

\medskip

(1) Given $\varphi\in A^*$, our claim is that the following limit converges, for any $a\in A$:
$$\int_\varphi a=\lim_{n\to\infty}\frac{1}{n}\sum_{k=1}^n\varphi^{*k}(a)$$

Indeed, we can assume, by linearity, that $a$ is the coefficient of a corepresentation:
$$a=(\tau\otimes id)v$$

But in this case, an elementary computation shows that we have the following formula, where $P_\varphi$ is the orthogonal projection onto the $1$-eigenspace of $(id\otimes\varphi)v$:
$$\left(id\otimes\int_\varphi\right)v=P_\varphi$$

(2) Since $v\xi=\xi$ implies $[(id\otimes\varphi)v]\xi=\xi$, we have $P_\varphi\geq P$, where $P$ is the orthogonal projection onto the following fixed point space:
$$Fix(v)=\left\{\xi\in\mathbb C^n\Big|v\xi=\xi\right\}$$

The point now is that when $\varphi\in A^*$ is faithful, by using a standard positivity trick, one can prove that we have $P_\varphi=P$. Assume indeed $P_\varphi\xi=\xi$, and let us set:
$$a=\sum_i\left(\sum_jv_{ij}\xi_j-\xi_i\right)\left(\sum_kv_{ik}\xi_k-\xi_i\right)^*$$

We must prove that we have $a=0$. Since $v$ is biunitary, we have:
\begin{eqnarray*}
a
&=&\sum_i\left(\sum_j\left(v_{ij}\xi_j-\frac{1}{N}\xi_i\right)\right)\left(\sum_k\left(v_{ik}^*\bar{\xi}_k-\frac{1}{N}\bar{\xi}_i\right)\right)\\
&=&\sum_{ijk}v_{ij}v_{ik}^*\xi_j\bar{\xi}_k-\frac{1}{N}v_{ij}\xi_j\bar{\xi}_i-\frac{1}{N}v_{ik}^*\xi_i\bar{\xi}_k+\frac{1}{N^2}\xi_i\bar{\xi}_i\\
&=&\sum_j|\xi_j|^2-\sum_{ij}v_{ij}\xi_j\bar{\xi}_i-\sum_{ik}v_{ik}^*\xi_i\bar{\xi}_k+\sum_i|\xi_i|^2\\
&=&||\xi||^2-<v\xi,\xi>-\overline{<v\xi,\xi>}+||\xi||^2\\
&=&2(||\xi||^2-Re(<v\xi,\xi>))
\end{eqnarray*}

By using now our assumption $P_\varphi\xi=\xi$, we obtain from this:
\begin{eqnarray*}
\varphi(a)
&=&2\varphi(||\xi||^2-Re(<v\xi,\xi>))\\
&=&2(||\xi||^2-Re(<P_\varphi\xi,\xi>))\\
&=&2(||\xi||^2-||\xi||^2)\\
&=&0
\end{eqnarray*}

Now since $\varphi$ is faithful, this gives $a=0$, and so $v\xi=\xi$. Thus $\int_\varphi$ is independent of $\varphi$, and is given on coefficients $a=(\tau\otimes id)v$ by the following formula:
$$\left(id\otimes\int_\varphi\right)v=P$$

(3) With the above formula in hand, the left and right invariance of $\int_G=\int_\varphi$ is clear on coefficients, and so in general, and this gives all the assertions. See \cite{wo1}. 
\end{proof}

Consider the dense $*$-subalgebra $\mathcal A\subset A$ generated by the coefficients of the fundamental corepresentation $u$, and endow it with the following scalar product: 
$$<a,b>=\int_Gab^*$$

Once again following Woronowicz \cite{wo1}, we have the following result:

\index{Peter-Weyl}
\index{characters}

\begin{theorem}
We have the following Peter-Weyl type results:
\begin{enumerate}
\item Any corepresentation decomposes as a sum of irreducible corepresentations.

\item Each irreducible corepresentation appears inside a certain $u^{\otimes k}$.

\item $\mathcal A=\bigoplus_{v\in Irr(A)}M_{\dim(v)}(\mathbb C)$, the summands being pairwise orthogonal.

\item The characters of irreducible corepresentations form an orthonormal system.
\end{enumerate}
\end{theorem}

\begin{proof}
All these results are from \cite{wo1}, the idea being as follows:

\medskip

(1) Given a corepresentation $v\in M_n(A)$, consider its interwiner algebra:
$$End(v)=\left\{T\in M_n(\mathbb C)\Big|Tv=vT\right\}$$

It is elementary to see that this is a finite dimensional $C^*$-algebra, and we conclude from this that we have a decomposition as follows:
$$End(v)=M_{n_1}(\mathbb C)\oplus\ldots\oplus M_{n_k}(\mathbb C)$$

To be more precise, such a decomposition appears by writing the unit of our algebra as a sum of minimal projections, as follows, and then working out the details:
$$1=p_1+\ldots+p_k$$

But this decomposition allows us to define subcorepresentations $v_i\subset v$, which are irreducible, so we obtain, as desired, a decomposition $v=v_1+\ldots+v_k$.

\medskip

(2) To any corepresentation $v\in M_n(A)$ we associate its space of coefficients, given by $C(v)=span(v_{ij})$. The construction $v\to C(v)$ is then functorial, in the sense that it maps subcorepresentations into subspaces. Observe also that we have:
$$\mathcal A=\sum_{k\in\mathbb N*\mathbb N}C(u^{\otimes k})$$

Now given an arbitrary corepresentation $v\in M_n(A)$, the corresponding coefficient space is a finite dimensional subspace $C(v)\subset\mathcal A$, and so we must have, for certain positive integers $k_1,\ldots,k_p$, an inclusion of vector spaces, as follows:
$$C(v)\subset C(u^{\otimes k_1}\oplus\ldots\oplus u^{\otimes k_p})$$

We deduce from this that we have an inclusion of corepresentations, as follows:
$$v\subset u^{\otimes k_1}\oplus\ldots\oplus u^{\otimes k_p}$$

Thus, by using (1), we are led to the conclusion in the statement.

\medskip

(3) By using (1) and (2), we obtain a linear space decomposition as follows:
$$\mathcal A
=\sum_{v\in Irr(A)}C(v)
=\sum_{v\in Irr(A)}M_{\dim(v)}(\mathbb C)$$

In order to conclude, it is enough to prove that for any two irreducible corepresentations $v,w\in Irr(A)$, the corresponding spaces of coefficients are orthogonal:
$$v\not\sim w\implies C(v)\perp C(w)$$ 

As a first observation, which follows from an elementary computation, for any two corepresentations $v,w$ we have a Frobenius type isomorphism, as follows:
$$Hom(v,w)\simeq Fix(\bar{v}\otimes w)$$

Now let us set $P_{ia,jb}=\int_Gv_{ij}w_{ab}^*$. According to Theorem 2.7, the matrix $P$ is the orthogonal projection onto the following vector space:
$$Fix(v\otimes\bar{w})
\simeq Hom(\bar{v},\bar{w})
=\{0\}$$

Thus we have $P=0$, and so $C(v)\perp C(w)$, which gives the result.

\medskip

(4) The algebra $\mathcal A_{central}$ contains indeed all the characters, because we have:
$$\Sigma\Delta(\chi_v)
=\sum_{ij}v_{ji}\otimes v_{ij}
=\Delta(\chi_v)$$

The fact that the characters span $\mathcal A_{central}$, and form an orthogonal basis of it, follow from (3). Finally, regarding the norm 1 assertion, consider the following integrals:
$$P_{ik,jl}=\int_Gv_{ij}v_{kl}^*$$

We know from Theorem 2.7 that these integrals form the orthogonal projection onto $Fix(v\otimes\bar{v})\simeq End(\bar{v})=\mathbb C1$. By using this fact, we obtain the following formula:
$$\int_G\chi_v\chi_v^*
=\sum_{ij}\int_Gv_{ii}v_{jj}^*
=\sum_i\frac{1}{N}
=1$$

Thus the characters have indeed norm 1, and we are done.
\end{proof}

We refer to Woronowicz \cite{wo1} for full details on all the above, and for some applications as well. Let us just record here the fact that in the cocommutative case, we obtain from (4) that the irreducible corepresentations must be all 1-dimensional, and so that we must have $A=C^*(\Gamma)$ for some discrete group $\Gamma$, as mentioned in Theorem 2.2. 

\bigskip

At a more technical level now, we have a number of more advanced results, from Woronowicz \cite{wo1}, \cite{wo2} and other papers, that must be known as well. We will present them quickly, and for details you check my book \cite{ba8}. First we have:

\index{amenability}
\index{Kesten amenability}

\begin{theorem}
Let $A_{full}$ be the enveloping $C^*$-algebra of $\mathcal A$, and let $A_{red}$ be the quotient of $A$ by the null ideal of the Haar integration. The following are then equivalent:
\begin{enumerate}
\item The Haar functional of $A_{full}$ is faithful.

\item The projection map $A_{full}\to A_{red}$ is an isomorphism.

\item The counit map $\varepsilon:A_{full}\to\mathbb C$ factorizes through $A_{red}$.

\item We have $N\in\sigma(Re(\chi_u))$, the spectrum being taken inside $A_{red}$.
\end{enumerate}
If this is the case, we say that the underlying discrete quantum group $\Gamma$ is amenable.
\end{theorem}

\begin{proof}
This is well-known in the group dual case, $A=C^*(\Gamma)$, with $\Gamma$ being a usual discrete group. In general, the result follows by adapting the group dual case proof:

\medskip

$(1)\iff(2)$ This simply follows from the fact that the GNS construction for the algebra $A_{full}$ with respect to the Haar functional produces the algebra $A_{red}$.

\medskip

$(2)\iff(3)$ Here $\implies$ is trivial, and conversely, a counit map $\varepsilon:A_{red}\to\mathbb C$ produces an isomorphism $A_{red}\to A_{full}$, via a formula of type $(\varepsilon\otimes id)\Phi$.

\medskip

$(3)\iff(4)$ Here $\implies$ is clear, coming from $\varepsilon(N-Re(\chi (u)))=0$, and the converse can be proved by doing some standard functional analysis.
\end{proof}

Yet another important result, also about the general Woronowicz algebras, and that we will be heavily using in what follows, is Tannakian duality, as follows:

\index{tensor category}
\index{Tannakian duality}

\begin{theorem}
The following operations are inverse to each other:
\begin{enumerate}
\item The construction $A\to C$, which associates to any Woronowicz algebra $A$ the tensor category formed by the intertwiner spaces $C_{kl}=Hom(u^{\otimes k},u^{\otimes l})$.

\item The construction $C\to A$, which associates to a tensor category $C$ the Woronowicz algebra $A$ presented by the relations $T\in Hom(u^{\otimes k},u^{\otimes l})$, with $T\in C_{kl}$.
\end{enumerate}
\end{theorem}

\begin{proof}
This is something quite deep, going back to Woronowicz \cite{wo2} in a slightly different form, and to Malacarne \cite{mal} in the simplified form presented above. The idea is that this can be proved by doing some abstract algebra, as follows:

\medskip

(1) We have indeed a construction $A\to C$ as above, whose output is a tensor $C^*$-subcategory with duals of the tensor $C^*$-category of Hilbert spaces.

\medskip

(2) We have as well a construction $C\to A$ as above, simply by dividing the free $*$-algebra on $N^2$ variables by the relations in the statement.

\medskip

Regarding now the bijection claim, some elementary algebra shows that $C=C_{A_C}$ implies $A=A_{C_A}$, and also that $C\subset C_{A_C}$ is automatic. Thus we are left with proving $C_{A_C}\subset C$. But this latter inclusion can be proved indeed, by doing some algebra, and using von Neumann's bicommutant theorem, in finite dimensions. See \cite{mal}. 
\end{proof}

As a concrete consequence of the above result, we have:

\index{algebraic manifold}

\begin{theorem}
We have an embedding as follows, using double indices,
$$G\subset S^{N^2-1}_{\mathbb C,+}\quad,\quad 
x_{ij}=\frac{u_{ij}}{\sqrt{N}}$$
making $G$ an algebraic submanifold of the free sphere.
\end{theorem}

\begin{proof}
The fact that we have an embedding as above follows from the fact that $u=(u_{ij})$ is biunitary, that we know from Proposition 2.4.  As for the algebricity claim, this follows from Theorem 2.10. Indeed, assuming that $A=C(G)$ is of the form $A=A_C$, it follows that $G$ is algebraic. But this is always the case, because we can take $C=C_A$.
\end{proof}

Observe that the embedding constructed above makes the link between our isomorphim conventions for quantum groups and for algebraic manifolds.

\section*{2b. Free rotations}

Let us get back now to our original objective, namely constructing pairs of quantum unitary and reflection groups $(O_N^+,H_N^+)$ and $(U_N^+,K_N^+)$, as to complete the pairs $(S^{N-1}_{\mathbb R,+},T_N^+)$ and $(S^{N-1}_{\mathbb C,+},\mathbb T_N^+)$ that we already have. Following Wang \cite{wa1}, we have:

\index{orthogonal quantum group}
\index{unitary quantum group}
\index{free quantum group}

\begin{theorem}
The following constructions produce compact quantum groups,
\begin{eqnarray*}
C(O_N^+)&=&C^*\left((u_{ij})_{i,j=1,\ldots,N}\Big|u=\bar{u},u^t=u^{-1}\right)\\
C(U_N^+)&=&C^*\left((u_{ij})_{i,j=1,\ldots,N}\Big|u^*=u^{-1},u^t=\bar{u}^{-1}\right)
\end{eqnarray*}
which appear respectively as liberations of the groups $O_N$ and $U_N$.
\end{theorem}

\begin{proof}
This first assertion follows from the elementary fact that if a matrix $u=(u_{ij})$ is orthogonal or biunitary, then so must be the following matrices:
$$u^\Delta_{ij}=\sum_ku_{ik}\otimes u_{kj}$$
$$u^\varepsilon_{ij}=\delta_{ij}$$
$$u^S_{ij}=u_{ji}^*$$

Indeed, the biunitarity of $u^\Delta$ can be checked by a direct computation. Regarding now the matrix $u^\varepsilon=1_N$, this is clearly biunitary. Also, regarding the matrix $u^S$, there is nothing to prove here either, because its unitarity its clear too. And finally, observe that if $u$ has self-adjoint entries, then so do the above matrices $u^\Delta,u^\varepsilon,u^S$.

\medskip

Thus our claim is proved, and we can define morphisms $\Delta,\varepsilon,S$ as in Definition 2.1, by using the universal properties of $C(O_N^+)$, $C(U_N^+)$. As for the second assertion, this follows exactly as for the free spheres, by adapting the sphere proof from chapter 1.
\end{proof}

The basic properties of $O_N^+,U_N^+$ can be summarized as follows:

\begin{theorem}
The quantum groups $O_N^+,U_N^+$ have the following properties:
\begin{enumerate}
\item The closed subgroups $G\subset U_N^+$ are exactly the $N\times N$ compact quantum groups. As for the closed subgroups $G\subset O_N^+$, these are those satisfying $u=\bar{u}$.

\item We have liberation embeddings $O_N\subset O_N^+$ and $U_N\subset U_N^+$, obtained by dividing the algebras $C(O_N^+),C(U_N^+)$ by their respective commutator ideals.

\item We have as well embeddings $\widehat{L}_N\subset O_N^+$ and $\widehat{F}_N\subset U_N^+$, where $L_N$ is the free product of $N$ copies of $\mathbb Z_2$, and where $F_N$ is the free group on $N$ generators.
\end{enumerate}
\end{theorem}

\begin{proof}
All these assertions are elementary, as follows:

\medskip

(1) This is clear from definitions, with the remark that, in the context of Definition 2.1, the formula $S(u_{ij})=u_{ji}^*$ shows that the matrix $\bar{u}$ must be unitary too. 

\medskip

(2) This follows from the Gelfand theorem. To be more precise, this shows that we have presentation results for $C(O_N),C(U_N)$, similar to those in Theorem 2.12, but with the commutativity between the standard coordinates and their adjoints added:
\begin{eqnarray*}
C(O_N)&=&C^*_{comm}\left((u_{ij})_{i,j=1,\ldots,N}\Big|u=\bar{u},u^t=u^{-1}\right)\\
C(U_N)&=&C^*_{comm}\left((u_{ij})_{i,j=1,\ldots,N}\Big|u^*=u^{-1},u^t=\bar{u}^{-1}\right)
\end{eqnarray*}

Thus, we are led to the conclusion in the statement.

\medskip

(3) This follows indeed from (1) and from Theorem 2.2, with the remark that with $u=diag(g_1,\ldots,g_N)$, the condition $u=\bar{u}$ is equivalent to $g_i^2=1$, for any $i$.
\end{proof}

The last assertion in Theorem 2.13 suggests the following construction:

\index{diagonal torus}
\index{group dual}

\begin{proposition}
Given a closed subgroup $G\subset U_N^+$, consider its ``diagonal torus'', which is the closed subgroup $T\subset G$ constructed as follows:
$$C(T)=C(G)\Big/\left<u_{ij}=0\Big|\forall i\neq j\right>$$
This torus is then a group dual, $T=\widehat{\Lambda}$, where $\Lambda=<g_1,\ldots,g_N>$ is the discrete group generated by the elements $g_i=u_{ii}$, which are unitaries inside $C(T)$.
\end{proposition}

\begin{proof}
Since $u$ is unitary, its diagonal entries $g_i=u_{ii}$ are unitaries inside $C(T)$. Moreover, from $\Delta(u_{ij})=\sum_ku_{ik}\otimes u_{kj}$ we obtain, when passing inside the quotient:
$$\Delta(g_i)=g_i\otimes g_i$$

It follows that we have $C(T)=C^*(\Lambda)$, modulo identifying as usual the $C^*$-completions of the various group algebras, and so that we have $T=\widehat{\Lambda}$, as claimed.
\end{proof}

With this notion in hand, Theorem 2.13 (3) reformulates as follows:

\begin{theorem}
The diagonal tori of the basic unitary groups are the basic tori:
$$\xymatrix@R=16.5mm@C=18mm{
O_N^+\ar[r]&U_N^+\\
O_N\ar[r]\ar[u]&U_N\ar[u]}
\qquad
\xymatrix@R=8mm@C=15mm{\\ \to}
\qquad
\xymatrix@R=16.5mm@C=18mm{
T_N^+\ar[r]&\mathbb T_N^+\\
T_N\ar[r]\ar[u]&\mathbb T_N\ar[u]}$$
In particular, the basic unitary groups are all distinct.
\end{theorem}

\begin{proof}
This is something clear and well-known in the classical case, and in the free case, this is a reformulation of Theorem 2.13 (3), which tells us that the diagonal tori of $O_N^+,U_N^+$, in the sense of Proposition 2.14, are the group duals $\widehat{L}_N,\widehat{F}_N$. 
\end{proof}

There is an obvious relation here with the considerations from chapter 1, that we will analyse later on. As a second result now regarding our free quantum groups, relating them this time to the free spheres constructed in chapter 1, we have:

\begin{proposition}
We have embeddings of algebraic manifolds as follows, obtained in double indices by rescaling the coordinates, $x_{ij}=u_{ij}/\sqrt{N}$:
$$\xymatrix@R=16.5mm@C=18mm{
O_N^+\ar[r]&U_N^+\\
O_N\ar[r]\ar[u]&U_N\ar[u]}
\qquad
\xymatrix@R=8mm@C=5mm{\\ \to}
\qquad
\xymatrix@R=15mm@C=14mm{
S^{N^2-1}_{\mathbb R,+}\ar[r]&S^{N^2-1}_{\mathbb C,+}\\
S^{N^2-1}_\mathbb R\ar[r]\ar[u]&S^{N^2-1}_\mathbb C\ar[u]
}$$
Moreover, the quantum groups appear from the quantum spheres via
$$G=S\cap U_N^+$$
with the intersection being computed inside the free sphere $S^{N^2-1}_{\mathbb C,+}$. 
\end{proposition}

\begin{proof}
As explained in Theorem 2.11, the biunitarity of the matrix $u=(u_{ij})$ gives an embedding of algebraic manifolds, as follows:
$$U_N^+\subset S^{N^2-1}_{\mathbb C,+}$$

Now since the relations defining $O_N,O_N^+,U_N\subset U_N^+$ are the same as those defining $S^{N^2-1}_\mathbb R,S^{N^2-1}_{\mathbb R,+},S^{N^2-1}_\mathbb C\subset S^{N^2-1}_{\mathbb C,+}$, this gives the result.
\end{proof}

\section*{2c. Free reflections}

Summarizing, in connection with our $(S,T,U,K)$ program, we have so far triples of type $(S,T,U)$, along with some correspondences between $S,T,U$. In order to introduce now the reflection groups $K$, things are more tricky, involving quantum permutation groups. Following Wang \cite{wa2}, these quantum groups are introduced as follows:

\index{quantum permutation}
\index{magic matrix}

\begin{theorem}
The following construction, where ``magic'' means formed of projections, which sum up to $1$ on each row and column,
$$C(S_N^+)=C^*\left((u_{ij})_{i,j=1,\ldots,N}\Big|u={\rm magic}\right)$$
produces a quantum group liberation of $S_N$. Moreover, the inclusion 
$$S_N\subset S_N^+$$
is an isomorphism at $N\leq3$, but not at $N\geq4$, where $S_N^+$ is not classical, nor finite.
\end{theorem}

\begin{proof}
We have several things to be proved, the idea being as follows:

\medskip

(1) The quantum group assertion follows by using the same arguments as those in the proof of Theorem 2.12. Consider indeed the following matrix: 
$$U_{ij}=\sum_ku_{ik}\otimes u_{kj}$$

As a first observation, the entries of this matrix are self-adjoint, $U_{ij}=U_{ij}^*$. In fact the entries $U_{ij}$ are orthogonal projections, because we have as well:
$$U_{ij}^2
=\sum_{kl}u_{ik}u_{il}\otimes u_{kj}u_{lj}
=\sum_ku_{ik}\otimes u_{kj}
=U_{ij}$$

In order to prove now that the matrix $U=(U_{ij})$ is magic, it remains to verify that the sums on the rows and columns are 1. For the rows, this can be checked as follows:
$$\sum_jU_{ij}
=\sum_{jk}u_{ik}\otimes u_{kj}
=\sum_ku_{ik}\otimes1
=1\otimes1$$

For the columns the computation is similar, as follows:
$$\sum_iU_{ij}
=\sum_{ik}u_{ik}\otimes u_{kj}
=\sum_k1\otimes u_{kj}
=1\otimes1$$

Thus the $U=(U_{ij})$ is magic, and so we can define a comultiplication map by using the universality property of $C(S_N^+)$, by setting $\Delta(u_{ij})=U_{ij}$. By using a similar reasoning, we can define as well a counit map by $\varepsilon(u_{ij})=\delta_{ij}$, and an antipode map by $S(u_{ij})=u_{ji}$. Thus the Woronowicz algebra axioms from Definition 2.1 are satisfied, and this finishes the proof of the first assertion, stating that $S_N^+$ is indeed a compact quantum group.

\medskip

(2) Observe now that we have an embedding of compact quantum groups $S_N\subset S_N^+$, obtained by using the standard coordinates of $S_N$, viewed as an algebraic group:
$$u_{ij}=\chi\left(\sigma\in S_N\Big|\sigma(j)=i\right)$$

By using the Gelfand theorem and working out the details, as we did with the free spheres are free unitary groups, the embedding $S_N\subset S_N^+$ is indeed a liberation. 

\medskip

(3) Finally, regarding the last assertion, the study here is as follows:

\medskip

\underline{Case $N=2$}. The result here is trivial, the $2\times2$ magic matrices being by definition as follows, with $p$ being a projection:
$$U=\begin{pmatrix}p&1-p\\1-p&p\end{pmatrix}$$

Indeed, this shows that the entries of a $2\times2$ magic matrix must pairwise commute, and so the algebra $C(S_2^+)$ follows to be commutative, which gives the result.

\medskip

\underline{Case $N=3$}. By using the same argument as in the $N=2$ case, and permuting rows and columns, it is enough to check that $u_{11},u_{22}$ commute. But this follows from:
\begin{eqnarray*}
u_{11}u_{22}
&=&u_{11}u_{22}(u_{11}+u_{12}+u_{13})\\
&=&u_{11}u_{22}u_{11}+u_{11}u_{22}u_{13}\\
&=&u_{11}u_{22}u_{11}+u_{11}(1-u_{21}-u_{23})u_{13}\\
&=&u_{11}u_{22}u_{11}
\end{eqnarray*}

Indeed, this gives $u_{22}u_{11}=u_{11}u_{22}u_{11}$, and then $u_{11}u_{22}=u_{22}u_{11}$, as desired.

\medskip

\underline{Case $N=4$}. In order to prove our various claims about $S_4^+$, consider the following matrix, with $p,q$ being projections, on some infinite dimensional Hilbert space:
$$U=\begin{pmatrix}
p&1-p&0&0\\
1-p&p&0&0\\
0&0&q&1-q\\
0&0&1-q&q
\end{pmatrix}$$ 

This matrix is magic, and if we choose $p,q$ as for the algebra $<p,q>$ to be not commutative, and infinite dimensional, we conclude that $C(S_4^+)$ is not commutative and infinite dimensional as well, and in particular is not isomorphic to $C(S_4)$.

\medskip

\underline{Case $N\geq5$}. Here we can use the standard embedding $S_4^+\subset S_N^+$, obtained at the level of the corresponding magic matrices in the following way:
$$u\to\begin{pmatrix}u&0\\ 0&1_{N-4}\end{pmatrix}$$

Indeed, with this embedding in hand, the fact that $S_4^+$ is a non-classical, infinite compact quantum group implies that $S_N^+$ with $N\geq5$ has these two properties as well.
\end{proof}

The above result came as a surprise at the time of \cite{wa2}, and there has been a lot of work since then, in order to understand what the quantum permutations really are, at $N\geq4$. We will be back to this, with further details, on several occasions. For the moment, let us just record the following alternative approach to $S_N^+$, also from Wang \cite{wa2}, which shows that we are not wrong with our formalism:

\begin{proposition}
The quantum group $S_N^+$ acts on the set $X=\{1,\ldots,N\}$, the corresponding coaction map $\Phi:C(X)\to C(X)\otimes C(S_N^+)$ being given by:
$$\Phi(e_i)=\sum_je_j\otimes u_{ji}$$
In fact, $S_N^+$ is the biggest compact quantum group acting on $X$, by leaving the counting measure invariant, in the sense that $(tr\otimes id)\Phi=tr(.)1$, where $tr(e_i)=\frac{1}{N},\forall i$.
\end{proposition}

\begin{proof}
Our claim is that given a compact matrix quantum group $G$, the following formula defines a morphism of algebras, which is a coaction map, leaving the trace invariant, precisely when the matrix $u=(u_{ij})$ is a magic corepresentation of $C(G)$: 
$$\Phi(e_i)=\sum_je_j\otimes u_{ji}$$

Indeed, let us first determine when $\Phi$ is multiplicative. We have:
$$\Phi(e_i)\Phi(e_k)
=\sum_{jl}e_je_l\otimes u_{ji}u_{lk}
=\sum_je_j\otimes u_{ji}u_{jk}$$

On the other hand, we have as well the following computation:
$$\Phi(e_ie_k)
=\delta_{ik}\Phi(e_i)
=\delta_{ik}\sum_je_j\otimes u_{ji}$$

We conclude that the multiplicativity of $\Phi$ is equivalent to the following conditions:
$$u_{ji}u_{jk}=\delta_{ik}u_{ji}\quad,\quad\forall i,j,k$$

Regarding now the unitality of $\Phi$, we have the following formula:
$$\Phi(1)
=\sum_i\Phi(e_i)
=\sum_{ij}e_j\otimes u_{ji}
=\sum_je_j\otimes\left(\sum_iu_{ji}\right)$$

Thus $\Phi$ is unital when $\sum_iu_{ji}=1$, $\forall j$. Finally, the fact that $\Phi$ is a $*$-morphism translates into $u_{ij}=u_{ij}^*$, $\forall i,j$. Summing up, in order for $\Phi(e_i)=\sum_je_j\otimes u_{ji}$ to be a morphism of $C^*$-algebras, the elements $u_{ij}$ must be projections, summing up to 1 on each row of $u$. Regarding now the preservation of the trace condition, observe that we have:
$$(tr\otimes id)\Phi(e_i)=\frac{1}{N}\sum_ju_{ji}$$

Thus the trace is preserved precisely when the elements $u_{ij}$ sum up to 1 on each of the columns of $u$. We conclude from this that $\Phi(e_i)=\sum_je_j\otimes u_{ji}$ is a morphism of $C^*$-algebras preserving the trace precisely when $u$ is magic, and since the coaction conditions on $\Phi$ are equivalent to the fact that $u$ must be a corepresentation, this finishes the proof of our claim. But this claim proves all the assertions in the statement.
\end{proof}

With the above results in hand, we can now introduce the quantum reflections:

\index{quantum reflection group}
\index{hyperoctahedral quantum group}
\index{free quantum group}

\begin{theorem}
The following constructions produce compact quantum groups,
\begin{eqnarray*}
C(H_N^+)&=&C^*\left((u_{ij})_{i,j=1,\ldots,N}\Big|u_{ij}=u_{ij}^*,\,(u_{ij}^2)={\rm magic}\right)\\
C(K_N^+)&=&C^*\left((u_{ij})_{i,j=1,\ldots,N}\Big|[u_{ij},u_{ij}^*]=0,\,(u_{ij}u_{ij}^*)={\rm magic}\right)
\end{eqnarray*}
which appear respectively as liberations of the reflection groups $H_N$ and $K_N$.
\end{theorem}

\begin{proof}
This can be proved in the usual way, with the first assertion coming from the fact that if $u$ satisfies the relations in the statement, then so do the matrices $u^\Delta,u^\varepsilon,u^S$, and with the second assertion coming as in the sphere case. See \cite{bb+}, \cite{bbc}.
\end{proof}

Summarizing, we are done with our construction task for the quadruplets $(S,T,U,K)$, in the free real and complex cases, and we can now formulate:

\begin{proposition}
We have a quadruplet as follows, called free real,
$$\xymatrix@R=50pt@C=50pt{
S^{N-1}_{\mathbb R,+}\ar@{-}[r]\ar@{-}[d]\ar@{-}[dr]&T_N^+\ar@{-}[l]\ar@{-}[d]\ar@{-}[dl]\\
O_N^+\ar@{-}[u]\ar@{-}[ur]\ar@{-}[r]&H_N^+\ar@{-}[l]\ar@{-}[ul]\ar@{-}[u]
}$$
and a quadruplet as follows, called free complex:
$$\xymatrix@R=50pt@C=50pt{
S^{N-1}_{\mathbb C,+}\ar@{-}[r]\ar@{-}[d]\ar@{-}[dr]&\mathbb T_N^+\ar@{-}[l]\ar@{-}[d]\ar@{-}[dl]\\
U_N^+\ar@{-}[u]\ar@{-}[ur]\ar@{-}[r]&K_N^+\ar@{-}[l]\ar@{-}[ul]\ar@{-}[u]
}$$
\end{proposition}

\begin{proof}
This is more of an empty statement, coming from the various constructions above, from chapter 1, and from the present chapter.
\end{proof}

Going ahead now, we must construct correspondences between our objects $(S,T,U,K)$, completing the work for the pairs $(S,T)$ started in chapter 1. This will take some time, and we will need some preliminaries. To start with, let us record the following result, which refines the various liberation statements formulated above:

\begin{theorem}
The quantum unitary and reflection groups are as follows,
$$\xymatrix@R=20pt@C=20pt{
&K_N^+\ar[rr]&&U_N^+\\
H_N^+\ar[rr]\ar[ur]&&O_N^+\ar[ur]\\
&K_N\ar[rr]\ar[uu]&&U_N\ar[uu]\\
H_N\ar[uu]\ar[ur]\ar[rr]&&O_N\ar[uu]\ar[ur]
}$$
and in this diagram, any face $P\subset Q,R\subset S$ has the property $P=Q\cap R$.
\end{theorem}

\begin{proof}
The fact that we have inclusions as in the statement follows from the definition of the various quantum groups involved. As for the various intersection claims, these follow as well from definitions. For some further details on all this, we refer to \cite{ba8}.
\end{proof}

\section*{2d. Diagrams, easiness}

In order to efficiently deal with the various quantum groups introduced above, we will need some specialized Tannakian duality results, in the spirit of the Brauer theorem \cite{bra}. Following \cite{bsp}, let us start with the following definition:

\index{partition}
\index{Kronecker symbol}

\begin{definition}
Associated to any partition $\pi\in P(k,l)$ between an upper row of $k$ points and a lower row of $l$ points is the linear map $T_\pi:(\mathbb C^N)^{\otimes k}\to(\mathbb C^N)^{\otimes l}$ given by 
$$T_\pi(e_{i_1}\otimes\ldots\otimes e_{i_k})=\sum_{j_1\ldots j_l}\delta_\pi\begin{pmatrix}i_1&\ldots&i_k\\ j_1&\ldots&j_l\end{pmatrix}e_{j_1}\otimes\ldots\otimes e_{j_l}$$
with the Kronecker type symbols $\delta_\pi\in\{0,1\}$ depending on whether the indices fit or not. 
\end{definition}

To be more precise, we agree to put the two multi-indices on the two rows of points, in the obvious way. The Kronecker symbols are then defined by $\delta_\pi=1$ when all the strings of $\pi$ join equal indices, and by $\delta_\pi=0$ otherwise. This construction is motivated by:

\begin{proposition}
The assignement $\pi\to T_\pi$ is categorical, in the sense that we have
$$T_\pi\otimes T_\sigma=T_{[\pi\sigma]}$$
$$T_\pi T_\sigma=N^{c(\pi,\sigma)}T_{[^\sigma_\pi]}$$
$$T_\pi^*=T_{\pi^*}$$
where $c(\pi,\sigma)$ are certain integers, coming from the erased components in the middle.
\end{proposition}

\begin{proof}
This comes from some routine computations, as follows:

\medskip

(1) The concatenation axiom follows from the following computation:
\begin{eqnarray*}
&&(T_\pi\otimes T_\sigma)(e_{i_1}\otimes\ldots\otimes e_{i_p}\otimes e_{k_1}\otimes\ldots\otimes e_{k_r})\\
&=&\sum_{j_1\ldots j_q}\sum_{l_1\ldots l_s}\delta_\pi\begin{pmatrix}i_1&\ldots&i_p\\j_1&\ldots&j_q\end{pmatrix}\delta_\sigma\begin{pmatrix}k_1&\ldots&k_r\\l_1&\ldots&l_s\end{pmatrix}e_{j_1}\otimes\ldots\otimes e_{j_q}\otimes e_{l_1}\otimes\ldots\otimes e_{l_s}\\
&=&\sum_{j_1\ldots j_q}\sum_{l_1\ldots l_s}\delta_{[\pi\sigma]}\begin{pmatrix}i_1&\ldots&i_p&k_1&\ldots&k_r\\j_1&\ldots&j_q&l_1&\ldots&l_s\end{pmatrix}e_{j_1}\otimes\ldots\otimes e_{j_q}\otimes e_{l_1}\otimes\ldots\otimes e_{l_s}\\
&=&T_{[\pi\sigma]}(e_{i_1}\otimes\ldots\otimes e_{i_p}\otimes e_{k_1}\otimes\ldots\otimes e_{k_r})
\end{eqnarray*}

(2) The composition axiom follows from the following computation:
\begin{eqnarray*}
&&T_\pi T_\sigma(e_{i_1}\otimes\ldots\otimes e_{i_p})\\
&=&\sum_{j_1\ldots j_q}\delta_\sigma\begin{pmatrix}i_1&\ldots&i_p\\j_1&\ldots&j_q\end{pmatrix}
\sum_{k_1\ldots k_r}\delta_\pi\begin{pmatrix}j_1&\ldots&j_q\\k_1&\ldots&k_r\end{pmatrix}e_{k_1}\otimes\ldots\otimes e_{k_r}\\
&=&\sum_{k_1\ldots k_r}N^{c(\pi,\sigma)}\delta_{[^\sigma_\pi]}\begin{pmatrix}i_1&\ldots&i_p\\k_1&\ldots&k_r\end{pmatrix}e_{k_1}\otimes\ldots\otimes e_{k_r}\\
&=&N^{c(\pi,\sigma)}T_{[^\sigma_\pi]}(e_{i_1}\otimes\ldots\otimes e_{i_p})
\end{eqnarray*}

(3) Finally, the involution axiom follows from the following computation:
\begin{eqnarray*}
&&T_\pi^*(e_{j_1}\otimes\ldots\otimes e_{j_q})\\
&=&\sum_{i_1\ldots i_p}<T_\pi^*(e_{j_1}\otimes\ldots\otimes e_{j_q}),e_{i_1}\otimes\ldots\otimes e_{i_p}>e_{i_1}\otimes\ldots\otimes e_{i_p}\\
&=&\sum_{i_1\ldots i_p}\delta_\pi\begin{pmatrix}i_1&\ldots&i_p\\ j_1&\ldots& j_q\end{pmatrix}e_{i_1}\otimes\ldots\otimes e_{i_p}\\
&=&T_{\pi^*}(e_{j_1}\otimes\ldots\otimes e_{j_q})
\end{eqnarray*}

Summarizing, our correspondence is indeed categorical. See \cite{bsp}.
\end{proof}

In analogy with the Tannakian categories, we have the following notion, from \cite{bsp}:

\index{category of partitions}

\begin{definition}
A collection of sets $D=\bigsqcup_{k,l}D(k,l)$ with $D(k,l)\subset P(k,l)$ is called a category of partitions when it has the following properties:
\begin{enumerate}
\item Stability under the horizontal concatenation, $(\pi,\sigma)\to[\pi\sigma]$.

\item Stability under vertical concatenation $(\pi,\sigma)\to[^\sigma_\pi]$, with matching middle symbols.

\item Stability under the upside-down turning $*$, with switching of colors, $\circ\leftrightarrow\bullet$.

\item Each set $P(k,k)$ contains the identity partition $||\ldots||$.

\item The sets $P(\emptyset,\circ\bullet)$ and $P(\emptyset,\bullet\circ)$ both contain the semicircle $\cap$.
\end{enumerate}
\end{definition} 

As a basic example, the set $D=P$ itself, formed by all partitions, is a category of partitions. The same goes for the category of pairings $P_2\subset P$. There are many other examples, and we will gradually explore them, in what follows.

\bigskip

Generally speaking, the axioms in Definition 2.24 can be thought of as being a ``delinearized version'' of the categorical conditions which are verified by the Tannakian categories. We have in fact the following result, going back to \cite{bsp}:

\index{Tannakian category}
\index{easy quantum group}
\index{easiness}

\begin{theorem}
Each category of partitions $D=(D(k,l))$ produces a family of compact quantum groups $G=(G_N)$, one for each $N\in\mathbb N$, via the formula
$$Hom(u^{\otimes k},u^{\otimes l})=span\left(T_\pi\Big|\pi\in D(k,l)\right)$$
which produces a Tannakian category, and therefore a closed subgroup $G_N\subset U_N^+$. The quantum groups which appear in this way are called ``easy''.
\end{theorem}

\begin{proof}
This follows indeed from Woronowicz's Tannakian duality, in its ``soft'' form from \cite{mal}, as explained in Theorem 2.10. Indeed, let us set:
$$C(k,l)=span\left(T_\pi\Big|\pi\in D(k,l)\right)$$

By using the axioms in Definition 2.24, and the categorical properties of the operation $\pi\to T_\pi$, from Proposition 2.23, we deduce that $C=(C(k,l))$ is a Tannakian category. Thus the Tannakian duality result applies, and gives the result.
\end{proof}

As a comment here, the word ``easy'' comes from what happens on the battleground, where we have many questions about quantum groups, and with Tannakian duality being our only serious tool. Thus, we can only call easy the quantum groups which are the simplest, with this meaning coming from partitions, from a Tannakian viewpoint.

\bigskip

Of course, you might find this terminology a bit strange, if you are new to the subject, but please believe me that everyone having worked on the subject struggled with easiness too. And, after fight comes some form of wisdom. Also, remember that as algebraic geometers, we are fans of Grothendieck, and his general idea of ``easy'' mathematics.

\bigskip

Be said in passing, in relation with this, modesty and everything, if you ever come across papers on easiness using alternative, complicated terms for easiness, better ignore them. Usually the more complicated the term used, the less funny the author.

\bigskip

Back to work now, we can formulate a general Brauer theorem, regarding the various quantum groups that we are interested in, as follows:

\index{Brauer theorem}
\index{orthogonal quantum group}
\index{unitary quantum group}
\index{quantum reflection group}
\index{hyperoctahedral quantum group}

\begin{theorem}
The basic quantum unitary and quantum reflection groups, namely
$$
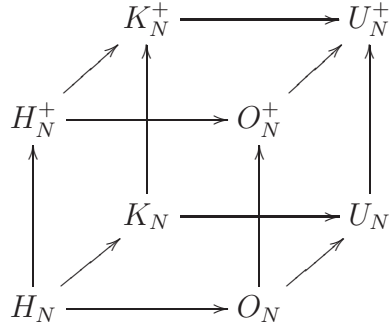
\xymatrix@R=18pt@C=18pt{
&K_N^+\ar[rr]&&U_N^+\\
H_N^+\ar[rr]\ar[ur]&&O_N^+\ar[ur]\\
&K_N\ar[rr]\ar[uu]&&U_N\ar[uu]\\
H_N\ar[uu]\ar[ur]\ar[rr]&&O_N\ar[uu]\ar[ur]
}$$
are all easy. The corresponding categories of partitions form an intersection diagram.
\end{theorem}

\begin{proof}
This is well-known, the categories being as follows, with $P_{even}$ being the category of partitions having even blocks, and with $\mathcal{P}_{even}(k,l)\subset P_{even}(k,l)$ consisting of the partitions satisfying $\#\circ=\#\bullet$ in each block, when flattening the partition:
$$
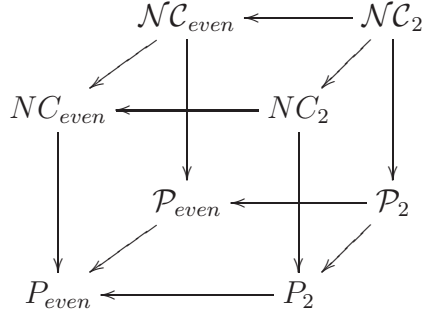
\xymatrix@R=18pt@C=2pt{
&\mathcal{NC}_{even}\ar[dl]\ar[dd]&&\mathcal {NC}_2\ar[dl]\ar[ll]\ar[dd]\\
NC_{even}\ar[dd]&&NC_2\ar[dd]\ar[ll]\\
&\mathcal P_{even}\ar[dl]&&\mathcal P_2\ar[dl]\ar[ll]\\
P_{even}&&P_2\ar[ll]
}$$

To be more precise, there is a long story with all this, with the results about $O_N,U_N$ going back to the 1937 paper of Brauer \cite{bra}, the results about $H_N,K_N$ being well-known too, for a long time, and with the quantum group results being more recent, from the 90s and 00s. We refer to \cite{ba8} for the whole story here, and in what concerns us, we can basically prove this, with the technology that we have, the idea being as follows:

\medskip

(1) The quantum group $U_N^+$ is defined via the following relations:
$$u^*=u^{-1}$$
$$u^t=\bar{u}^{-1}$$ 

But these relations tell us precisely that the following two operators must be in the associated Tannakian category $C$:
$$T_\pi\quad,\quad \pi={\ }^{\,\cap}_{\circ\bullet}$$
$$T_\pi\quad,\quad \pi={\ }^{\,\cap}_{\bullet\circ}$$

Thus the associated Tannakian category is $C=span(T_\pi|\pi\in D)$, with:
$$D
=<{\ }^{\,\cap}_{\circ\bullet}\,\,,{\ }^{\,\cap}_{\bullet\circ}>
={\mathcal NC}_2$$
  
Thus, we are led to the conclusion in the statement.

\medskip

(2) The quantum group $O_N^+\subset U_N^+$ is defined by imposing the following relations:
$$u_{ij}=\bar{u}_{ij}$$

But these relations tell us that the following operators must be in the associated Tannakian category $C$:
$$T_\pi\quad,\quad\pi=|^{\hskip-1.32mm\circ}_{\hskip-1.32mm\bullet}$$
$$T_\pi\quad,\quad \pi=|_{\hskip-1.32mm\circ}^{\hskip-1.32mm\bullet}$$

Thus the associated Tannakian category is $C=span(T_\pi|\pi\in D)$, with:
$$D
=<\mathcal{NC}_2,|^{\hskip-1.32mm\circ}_{\hskip-1.32mm\bullet},|_{\hskip-1.32mm\circ}^{\hskip-1.32mm\bullet}>
=NC_2$$
  
Thus, we are led to the conclusion in the statement.

\medskip

(3) The group $U_N\subset U_N^+$ is defined via the following relations:
$$[u_{ij},u_{kl}]=0$$
$$[u_{ij},\bar{u}_{kl}]=0$$

But these relations tell us that the following operators must be in the associated Tannakian category $C$:
$$T_\pi\quad,\quad \pi={\slash\hskip-2.1mm\backslash}^{\hskip-2.5mm\circ\circ}_{\hskip-2.5mm\circ\circ}$$
$$T_\pi\quad,\quad \pi={\slash\hskip-2.1mm\backslash}^{\hskip-2.5mm\circ\bullet}_{\hskip-2.5mm\bullet\circ}$$

Thus the associated Tannakian category is $C=span(T_\pi|\pi\in D)$, with:
$$D
=<\mathcal{NC}_2,{\slash\hskip-2.1mm\backslash}^{\hskip-2.5mm\circ\circ}_{\hskip-2.5mm\circ\circ},{\slash\hskip-2.1mm\backslash}^{\hskip-2.5mm\circ\bullet}_{\hskip-2.5mm\bullet\circ}>
=\mathcal P_2$$

Thus, we are led to the conclusion in the statement.

\medskip

(4) In order to deal now with $O_N$, we can simply use the following formula: 
$$O_N=O_N^+\cap U_N$$

At the categorical level, this tells us indeed that the associated Tannakian category is given by $C=span(T_\pi|\pi\in D)$, with:
$$D
=<NC_2,\mathcal P_2>
=P_2$$

Thus, we are led to the conclusion in the statement.

\medskip

(5) The proof for the reflection groups is similar, by first proving that $S_N$ is easy, corresponding to the category of all partitions $P$, and then by adding and suitably interpreting the reflection relations. We refer here to \cite{bb+}, \cite{bbc}, for details.

\medskip

(6) The proof for the quantum reflection groups is similar, by first proving that the quantum permutation group $S_N^+$ is easy, corresponding to the category of all noncrossing partitions $NC$, and then by adding and suitably interpreting the quantum reflection relations. As before, we refer here to \cite{bb+}, \cite{bbc}, for details.

\medskip

(7) As for the last assertion, which will be of use later on, this is something well-known and standard too. We refer here to \cite{bb+}, \cite{bbc}, and to \cite{ba8}, \cite{bsp} as well.
\end{proof}

Getting back now to our axiomatization questions, we must establish correspondences between our objects $(S,T,U,K)$, as a continuation of the work started in chapter 1, for the pairs $(S,T)$. Let us start by discussing the following correspondences:
$$U\to K\to T$$

We know from Theorem 2.15 that the correspondences $U\to T$ appear by taking the diagonal tori. In fact, the correspondences $K\to T$ appear by taking the diagonal tori as well, and the correspondences $U\to K$ are something elementary too, obtained by taking the ``reflection subgroup''. The complete statement here is as follows:

\index{reflection subgroup}
\index{diagonal torus}

\begin{theorem}
For the basic quadruplets $(S,T,U,K)$, the correspondences
$$\xymatrix@R=15mm@C=17mm{
O_N^+\ar[r]&U_N^+\\
O_N\ar[r]\ar[u]&U_N\ar[u]}
\quad
\xymatrix@R=8mm@C=15mm{\\ \to}
\quad
\xymatrix@R=15mm@C=17mm{
H_N^+\ar[r]&K_N^+\\
H_N\ar[r]\ar[u]&K_N\ar[u]}
\quad
\xymatrix@R=8mm@C=15mm{\\ \to}
\quad
\xymatrix@R=15mm@C=16mm{
T_N^+\ar[r]&\mathbb T_N^+\\
T_N\ar[r]\ar[u]&\mathbb T_N\ar[u]}$$
appear in the following way:
\begin{enumerate}
\item $U\to K$ appears by taking the reflection subgroup, $K=U\cap K_N^+$.

\item $U\to T$ appears by taking the diagonal torus, $T=U\cap\mathbb T_N^+$.

\item $K\to T$ appears as well by taking the diagonal torus, $T=K\cap\mathbb T_N^+$.
\end{enumerate}
\end{theorem}

\begin{proof}
This follows from the results that we already have, as follows:

\medskip

(1) This follows from Theorem 2.26, because the left face of the cube diagram there appears by intersecting the right face with the quantum group $K_N^+$.

\medskip

(2) This is something that we already know, from Theorem 2.15.

\medskip

(3) This follows exactly as in the unitary case, via the proof of Theorem 2.15.
\end{proof}

As a conclusion now, with respect to the ``baby theory'' developed in chapter 1, concerning the pairs $(S,T)$, we have some advances. First, we have completed the pairs $(S,T)$ there into quadruplets $(S,T,U,K)$. And second, we have established some correspondences between our objects, the situation here being as follows:
$$\xymatrix@R=50pt@C=50pt{
S\ar[r]&T\\
U\ar[r]\ar[ur]&K\ar[u]
}$$

There is still a long way to go, in order to establish a full set of correspondences, and to reach to an axiomatization, the idea being that the correspondences $S\leftrightarrow U$ can be established by using quantum isometries, and that the correspondences $T\to K\to U$ can be established by using advanced quantum group theory, and with all this heavily relying on the easiness theory developed above. We will discuss this in chapters 3-4 below.

\section*{2e. Exercises}

As before with the first chapter, the theory explained in the above will be of key importance in what follows, and as best homework on this, if interested, we can only recommend reading a good quantum group book. Here is however a first exercise:

\begin{exercise}
Develop a theory of finite quantum groups, in analogy with the usual group theory, by assuming that the Woronowicz algebras $A$ are finite dimensional.
\end{exercise}

This is a bit vague, and as bottom line, such a theory would need clear axioms, the Pontrjagin type duality worked out in detail, a few theorems, and examples.

\begin{exercise}
Work out what happens with the Haar integration and Peter-Weyl theory, theorems and proofs, in the classical group case, and in the group dual case.
\end{exercise}

This is something instructive, leading to a good functional analysis knowledge.

\begin{exercise}
Work out what happens with the amenability and Tannakian duality, theorems and proofs, in the classical group case, and in the group dual case.
\end{exercise}

As before with the previous exercise, this is something quite instructive.

\begin{exercise}
Prove that the dual of $S_5^+$ is not amenable.
\end{exercise}

As a hint here, try finding a subgroup which is not amenable. Also, as a bonus exercise, you can try as well to prove that the dual of $S_4^+$ is amenable, although this is quite tricky, and reputed to be undoable with bare hands. But who knows.

\chapter{Affine isometries}

\section*{3a. Quantum isometries}

We have seen so far that we have quadruplets $(S,T,U,K)$ consisting of a sphere $S$, a torus $T$, a unitary group $U$ and a reflection group $K$, corresponding to the four main geometries, namely real and complex, classical and free, which are as follows:
$$\xymatrix@R=50pt@C=50pt{
\mathbb R^N_+\ar[r]&\mathbb C^N_+\\
\mathbb R^N\ar[u]\ar[r]&\mathbb C^N\ar[u]
}$$

Here the upper symbols $\mathbb R^N_+,\mathbb C^N_+$ do not stand for the free analogues of $\mathbb R^N,\mathbb C^N$, which do not exist as such, but rather for the ``noncommutative geometry'' of these free analogues, which does exist, via the quadruplets $(S,T,U,K)$ that we constructed for them. As for the arrows, these stand for the obvious inclusions between the objects $S,T,U,K$. More on these notations in chapter 4, after axiomatizing everything.

\bigskip

We have now to work out the various correspondences between our objects $(S,T,U,K)$, as to reach to a full set of correspondences, in each of the above 4 cases. We know from chapters 1-2 that we already have 4 such correspondences. In this chapter we discuss 3 more correspondences, as to reach to a total of 7 correspondences, as follows:
$$\xymatrix@R=50pt@C=50pt{
S\ar[r]&T\\
U\ar[r]\ar[ur]&K\ar[u]
}\qquad
\xymatrix@R=23pt@C=50pt{\\ \to}
\qquad
\xymatrix@R=50pt@C=50pt{
S\ar[d]\ar[r]&T\ar[d]\\
U\ar[r]\ar[ur]\ar[u]&K\ar[u]
}$$

In order to connect the spheres and tori $(S,T)$ to the quantum groups $(U,K)$, the idea will be that of using quantum isometry groups. However, normally ``isometry'' comes from iso and metric, and so is something preserving the metric, and remember from the various discussions from chapters 1-2 that our objects $S,T,U,K$ are not exactly ``quantum metric spaces'' in some reasonable sense. Let us record this as a fact, to start with:

\begin{fact}
Our objects $S,T,U,K$ are not quantum metric spaces, in some reasonable sense, and so cannot have quantum isometry groups in a usual, iso\,+\,metric sense.
\end{fact}

Here the word ``reasonable'' can only suggest that we are into some controversies, and so are we, indeed. But dealing with this controversy is an easy task, because we can send packing any criticism with the following argument. The only reasonable notion of quantum metric is that of Connes \cite{co1}, based on a Dirac operator, and don't you dare to think otherwise, and since our manifolds $S,T,U,K$ do not have a Dirac operator a la Connes, they cannot be quantum metric spaces, in some reasonable sense. QED.

\bigskip

This being said, recall that $O_N,U_N$ are the isometry groups of $\mathbb R^N,\mathbb C^N$, with of course some care with respect to the complex structure when talking $U_N$. And now since $O_N^+,U_N^+$ are straightforward, very natural liberations of $O_N,U_N$, we can definitely think at $O_N^+,U_N^+$ as being ``quantum isometry groups'', in a somewhat abstract sense. So, regardless of Fact 3.1 says, we feel entitled to talk about quantum isometries, in our setting, be that in a bit abstract and unorthodox way, not exactly coming from iso\,+\,metric.

\bigskip

So, we have an idea here, in order to short-circuit Fact 3.1. And fortunately, in order to make now our point clear, and be able to rigorously talk about quantum isometries in our sense, pure mathematics comes to the rescue. In the classical case, we have indeed the following trivial speculation, taking us away from serious, metric geometry:

\index{isometry}
\index{affine isometry}

\begin{proposition}
Given a closed subset $X\subset S^{N-1}_\mathbb C$, the formula
$$G(X)=\left\{U\in U_N\Big|U(X)=X\right\}$$
defines a compact group of unitary matrices, or isometries, called affine isometry group of $X$. For the spheres $S^{N-1}_\mathbb R,S^{N-1}_\mathbb C$ we obtain in this way the groups $O_N,U_N$.
\end{proposition}

\begin{proof}
The fact that $G(X)$ as defined above is indeed a group is clear, its compactness is clear as well, and finally the last assertion is clear as well. In fact, all this works for any closed subset $X\subset\mathbb C^N$, but we are not interested here in such general spaces.
\end{proof}

Observe that in the case of the real and complex spheres, the affine isometry group $G(X)$ leaves invariant the Riemannian metric, because this metric is equivalent to the one inherited from $\mathbb C^N$, which is preserved by our isometries $U\in U_N$. Thus, we could have constructed as well $G(X)$ as being the group of metric isometries of $X$, with of course some extra care in relation with the complex structure, as for the complex sphere $X=S^{N-1}_\mathbb C$ to produce $G(X)=U_N$ instead of $G(X)=O_{2N}$. But, as already indicated in Fact 3.1, such things won't work for the free spheres, and so are to be avoided. 

\bigskip

The point now is that we have the following quantum analogue of Proposition 3.2, which is a perfect analogue, save for the fact that $X$ is now assumed to be algebraic, for some technical reasons, which allows us to talk about quantum isometry groups:

\index{quantum isometry}
\index{affine quantum isometry}
\index{quantum isometry group}

\begin{theorem}
Given an algebraic manifold $X\subset S^{N-1}_{\mathbb C,+}$, the category of closed subgroups $G\subset U_N^+$ acting affinely on $X$, in the sense that the formula
$$\Phi(x_i)=\sum_jx_j\otimes u_{ji}$$ 
defines a morphism of $C^*$-algebras, as follows, 
$$\Phi:C(X)\to C(X)\otimes C(G)$$
has a universal object, denoted $G^+(X)$, and called affine quantum isometry group of $X$. When $X$ is classical, $G^+(X)$ is a liberation of $G(X)$.
\end{theorem}

\begin{proof}
As it might be obvious from the above discussion, we are a bit into muddy waters here, with the result itself being some sort of matematical trick, in order to avoid serious geometry, which unfortunately does not exist in the free case. But, the statement makes sense as stated, so let us just prove it, and we will comment on it afterwards:

\medskip

(1) As a first observation, we have already met such a result, at the end of chapter 1, when talking about toral isometries. In general, the proof will be quite similar.

\medskip

(2) Another observation is that, in the case where $\Phi$ as above exists, this morphism is automatically a coaction, in the sense that it satisfies the following conditions:
$$(\Phi\otimes id)\Phi=(id\otimes\Delta)\Phi$$
$$(id\otimes\varepsilon)\Phi=id$$

(3) As a technical comment now, such coactions $\Phi$ can be thought of as coming from actions $G\curvearrowright X$ of the corresponding quantum group, written as follows:
$$X\times G\to X\quad,\quad (x,g)\to g(x)$$ 

So you might say why not doing it the other way around, by using morphisms of type $\Phi:C(X)\to C(G)\otimes C(X)$, corresponding to more familiar actions $(g,x)\to g(x)$. This is a good point, and in answer, at the basic level things are of course equivalent, and we refer here to \cite{bss} for more, but at the advanced level, and more specifically in the context of the actions from \cite{bbs}, discussed in chapter 16 below, things are definitely better written by using coactions $\Phi:C(X)\to C(X)\otimes C(G)$. So, in short, we had a left/right choice to be made, and based on some advanced considerations, we made the right choice.

\medskip

(4) In order to prove now the result, assume that $X\subset S^{N-1}_{\mathbb C,+}$ comes as follows:
$$C(X)=C(S^{N-1}_{\mathbb C,+})\Big/\Big<f_\alpha(x_1,\ldots,x_N)=0\Big>$$

Consider now the following variables:
$$X_i=\sum_jx_j\otimes u_{ji}\in C(X)\otimes C(U_N^+)$$

Our claim is that the quantum group in the statement $G=G^+(X)$ appears as:
$$C(G)=C(U_N^+)\Big/\Big<f_\alpha(X_1,\ldots,X_N)=0\Big>$$

(5) In order to prove this, we have to clarify how the relations $f_\alpha(X_1,\ldots,X_N)=0$ are interpreted inside $C(U_N^+)$, and then show that $G$ is indeed a quantum group. So, pick one of the defining polynomials, and write it as follows:
$$f_\alpha(x_1,\ldots,x_N)=\sum_r\sum_{i_1^r\ldots i_{s_r}^r}\lambda_r\cdot x_{i_1^r}\ldots x_{i_{s_r}^r}$$

With $X_i=\sum_jx_j\otimes u_{ji}$ as above, we have the following formula:
$$f_\alpha(X_1,\ldots,X_N)=\sum_r\sum_{i_1^r\ldots i_{s_r}^r}\lambda_r\sum_{j_1^r\ldots j_{s_r}^r}x_{j_1^r}\ldots x_{j_{s_r}^r}\otimes u_{j_1^ri_1^r}\ldots u_{j_{s_r}^ri_{s_r}^r}$$

Since the variables on the right span a certain finite dimensional space, the relations $f_\alpha(X_1,\ldots,X_N)=0$ correspond to certain relations between the variables $u_{ij}$. Thus, we have indeed a closed subspace $G\subset U_N^+$, with a universal map, as follows:
$$\Phi:C(X)\to C(X)\otimes C(G)$$

(6) In order to show now that $G$ is a quantum group, consider the following elements:
$$u_{ij}^\Delta=\sum_ku_{ik}\otimes u_{kj}\quad,\quad u_{ij}^\varepsilon=\delta_{ij}\quad,\quad u_{ij}^S=u_{ji}^*$$

Consider as well the following elements, with $\gamma\in\{\Delta,\varepsilon,S\}$:
$$X_i^\gamma=\sum_jx_j\otimes u_{ji}^\gamma$$

From the relations $f_\alpha(X_1,\ldots,X_N)=0$ we deduce that we have:
$$f_\alpha(X_1^\gamma,\ldots,X_N^\gamma)
=(id\otimes\gamma)f_\alpha(X_1,\ldots,X_N)
=0$$

Thus we can map $u_{ij}\to u_{ij}^\gamma$ for any $\gamma\in\{\Delta,\varepsilon,S\}$, and we are done.

\medskip

(7) Regarding now the last assertion, assume that we have $X\subset S^{N-1}_\mathbb C$, as in Proposition 3.2. In functional analytic terms, the definition of the group $G(X)$ there tells us that we must have a morphism $\Phi$ as in the statement. Thus we have $G(X)\subset G^+(X)$, and moreover, the classical version of $G^+(X)$ is the group $G(X)$, as desired.
\end{proof}

Before getting further, we should clarify the relation between Proposition 3.2, Theorem 3.3, and the ``toral isometry'' constructions from chapter 1. We have:

\begin{theorem}
Given an algebraic manifold $X\subset S^{N-1}_{\mathbb C,+}$, the category of closed subgroups $G\subset H$ acting affinely on $X$, with $H$ being one of the following quantum groups,
$$\xymatrix@R=18mm@C=22mm{
\mathbb T_N^+\ar[r]&K_N^+\ar[r]&U_N^+\\
\mathbb T_N\ar[r]\ar[u]&K_N\ar[r]\ar[u]&U_N\ar[u]
}$$
has a universal object, denoted respectively as follows,
$$\xymatrix@R=17mm@C=15mm{
T^+(X)\ar[r]&K^+(X)\ar[r]&G^+(X)\\
T(X)\ar[r]\ar[u]&K(X)\ar[r]\ar[u]&G(X)\ar[u]
}$$
which appears by intersecting $G^+(X)$ and $H$, inside $U_N^+$.
\end{theorem}

\begin{proof}
Here the assertion regarding $G^+(X)$ is something that we know, from Theorem 3.3, and all the other assertions follow from this, by intersecting with $H$.
\end{proof}

Summarizing, we have a reasonable notion of quantum isometry group, for the manifolds $X\subset S^{N-1}_{\mathbb C,+}$ that we are interested in, in this book, and we will heavily use this notion, in what follows. However, all this is tricky, and here is the story with this:

\bigskip

(1) Things go back to the paper of Goswami \cite{go1}, who proved there that any compact Riemannian manifold $X$ has a quantum isometry group $\mathcal G^+(X)$, liberating the usual isometry group $\mathcal G(X)$, and with isometry meaning here, of course, preserving the Riemannian metric. Moreover, Goswami proved in \cite{go1} that the same construction can be performed for a compact quantum Riemannian manifold in the sense of Connes \cite{co1}.

\bigskip

(2) All this is very interesting, and there has been a lot of work on the subject, by Goswami, his student Bhowmick, and their collaborators. Let us mention here the fundamental papers \cite{bg1}, \cite{bg2}, \cite{go2}, the rigidity theorem of Goswami in \cite{go3}, stating that $\mathcal G^+(X)=\mathcal G(X)$ when $X$ is connected, the work of Bhowmick et al. \cite{bdd}, \cite{bd+} on the isometries of the Chamseddine-Connes manifold \cite{cc1}, \cite{cc2}, and the book \cite{gbh}.

\bigskip

(3) In what regards the free case, however, this ``metric'' theory does not work, and we must basically trick as in Theorem 3.3. There are of course some versions of this, which might sound a bit more geometric, with some sort of fake spectral triple constructed in \cite{bgo}, then a beast called ``orthogonal filtration'' in \cite{bsk}, then a true and honest Laplacian constructed in \cite{dfw}. But, Fact 3.1 remains there as stated, end of the story.

\bigskip

(4) Importantly now, for $S^{N-1}_\mathbb R$ all sorts of quantum isometry groups that you can ever imagine all coincide, with $O_N$. And, as shown by computations in \cite{bgo} and related papers, the same kind of phenomenon holds for other types of spheres. And so, there is a bit of a non-problem with all this, and as long as you are interested in simple manifolds like spheres, like we do here, just use Theorem 3.3, no need to bother with more.

\bigskip

(5) Which of course does not mean that there are no interesting problems left. The main question here is that of axiomatizing and extending the Laplacian construction of Das-Franz-Wang in \cite{dfw}, then doing something hybrid between \cite{bsk}, \cite{go1}, meaning general construction of $\mathcal G^+(X)$, by using this Laplacian, and then finally comparing this quantum group with $G^+(X)$, with probably a first-class mathematical theorem at stake.

\bigskip

And this is all, for the moment. We will be back to this on several occasions, first at the end of the present chapter, with some details on the various generalities above, including the rigidity theorem of Goswami in \cite{go3}, which is the main result at the general level. And also, we will certainly talk about the work of Bhowmick et al. \cite{bdd}, \cite{bd+} on the isometries of the Chamseddine-Connes manifold \cite{cc1}, \cite{cc2}, later in this book. But for the moment, let us just enjoy Theorem 3.3 as it is, providing us with a simple definition for the quantum isometry groups, in our setting, and do some computations.

\section*{3b. Spheres and rotations}

In connection with our axiomatization questions for the quadruplets $(S,T,U,K)$, we can construct now the correspondences $S\to U$, in the following way:

\begin{theorem}
The quantum isometry groups of the basic spheres are
$$\xymatrix@R=15mm@C=14mm{
S^{N-1}_{\mathbb R,+}\ar[r]&S^{N-1}_{\mathbb C,+}\\
S^{N-1}_\mathbb R\ar[r]\ar[u]&S^{N-1}_\mathbb C\ar[u]
}
\qquad
\xymatrix@R=8mm@C=15mm{\\ \to}
\qquad
\xymatrix@R=16mm@C=18mm{
O_N^+\ar[r]&U_N^+\\
O_N\ar[r]\ar[u]&U_N\ar[u]}$$
modulo identifying, as usual, the various $C^*$-algebraic completions.
\end{theorem}

\begin{proof}
We have 4 results to be proved, and following \cite{bgo}, \cite{bg1} and related papers, where this result was established in its above form, we can proceed as follows:

\medskip

\underline{$S^{N-1}_{\mathbb C,+}$}. Let us first construct an action $U_N^+\curvearrowright S^{N-1}_{\mathbb C,+}$. We must prove here that the variables $X_i=\sum_jx_j\otimes u_{ji}$ satisfy the defining relations for $S^{N-1}_{\mathbb C,+}$, namely:
$$\sum_ix_ix_i^*=\sum_ix_i^*x_i=1$$

By using the biunitarity of $u$, we have the following computation:
\begin{eqnarray*}
\sum_iX_iX_i^*
&=&\sum_{ijk}x_jx_k^*\otimes u_{ji}u_{ki}^*\\
&=&\sum_jx_jx_j^*\otimes1\\
&=&1\otimes1
\end{eqnarray*}

Once again by using the biunitarity of $u$, we have as well:
\begin{eqnarray*}
\sum_iX_i^*X_i
&=&\sum_{ijk}x_j^*x_k\otimes u_{ji}^*u_{ki}\\
&=&\sum_jx_j^*x_j\otimes1\\
&=&1\otimes1
\end{eqnarray*}

Thus we have an action $U_N^+\curvearrowright S^{N-1}_{\mathbb C,+}$, which gives $G^+(S^{N-1}_{\mathbb C,+})=U_N^+$, as desired. 

\medskip

\underline{$S^{N-1}_{\mathbb R,+}$}. Let us first construct an action $O_N^+\curvearrowright S^{N-1}_{\mathbb R,+}$. We already know that the variables $X_i=\sum_jx_j\otimes u_{ji}$ satisfy the defining relations for $S^{N-1}_{\mathbb C,+}$, so we just have to check that these variables are self-adjoint. But this is clear from $u=\bar{u}$, as follows:
$$X_i^*
=\sum_jx_j^*\otimes u_{ji}^*
=\sum_jx_j\otimes u_{ji}
=X_i$$

Conversely, assume that we have an action $G\curvearrowright S^{N-1}_{\mathbb R,+}$, with $G\subset U_N^+$. The variables $X_i=\sum_jx_j\otimes u_{ji}$ must be then self-adjoint, and the above computation shows that we must have $u=\bar{u}$. Thus our quantum group must satisfy $G\subset O_N^+$, as desired.

\medskip

\underline{$S^{N-1}_\mathbb C$}. The fact that we have an action $U_N\curvearrowright S^{N-1}_\mathbb C$ is clear. Conversely, assume that we have an action $G\curvearrowright S^{N-1}_\mathbb C$, with $G\subset U_N^+$. We must prove that this implies $G\subset U_N$, and we will use a standard trick of Bhowmick-Goswami \cite{bg1}. We have:
$$\Phi(x_i)=\sum_jx_j\otimes u_{ji}$$

By multiplying this formula with itself we obtain:
$$\Phi(x_ix_k)=\sum_{jl}x_jx_l\otimes u_{ji}u_{lk}$$
$$\Phi(x_kx_i)=\sum_{jl}x_lx_j\otimes u_{lk}u_{ji}$$

Since the variables $x_i$ commute, these formulae can be written as:
$$\Phi(x_ix_k)=\sum_{j<l}x_jx_l\otimes(u_{ji}u_{lk}+u_{li}u_{jk})+\sum_jx_j^2\otimes u_{ji}u_{jk}$$
$$\Phi(x_ix_k)=\sum_{j<l}x_jx_l\otimes(u_{lk}u_{ji}+u_{jk}u_{li})+\sum_jx_j^2\otimes u_{jk}u_{ji}$$

Since the tensors at left are linearly independent, we must have:
$$u_{ji}u_{lk}+u_{li}u_{jk}=u_{lk}u_{ji}+u_{jk}u_{li}$$

By applying the antipode to this formula, then applying the involution, and then relabelling the indices, we succesively obtain:
$$u_{kl}^*u_{ij}^*+u_{kj}^*u_{il}^*=u_{ij}^*u_{kl}^*+u_{il}^*u_{kj}^*$$
$$u_{ij}u_{kl}+u_{il}u_{kj}=u_{kl}u_{ij}+u_{kj}u_{il}$$
$$u_{ji}u_{lk}+u_{jk}u_{li}=u_{lk}u_{ji}+u_{li}u_{jk}$$

Now by comparing with the original formula, we obtain from this:
$$u_{li}u_{jk}=u_{jk}u_{li}$$

In order to finish, it remains to prove that the coordinates $u_{ij}$ commute as well with their adjoints. For this purpose, we use a similar method. We have:
$$\Phi(x_ix_k^*)=\sum_{jl}x_jx_l^*\otimes u_{ji}u_{lk}^*$$
$$\Phi(x_k^*x_i)=\sum_{jl}x_l^*x_j\otimes u_{lk}^*u_{ji}$$

Since the variables on the left are equal, we deduce from this that we have:
$$\sum_{jl}x_jx_l^*\otimes u_{ji}u_{lk}^*=\sum_{jl}x_jx_l^*\otimes u_{lk}^*u_{ji}$$

Thus we have $u_{ji}u_{lk}^*=u_{lk}^*u_{ji}$, and so $G\subset U_N$, as claimed.

\medskip

\underline{$S^{N-1}_\mathbb R$}. The fact that we have an action $O_N\curvearrowright S^{N-1}_\mathbb R$ is clear. In what regards the converse, this follows by combining the results that we already have, as follows:
\begin{eqnarray*}
G\curvearrowright S^{N-1}_\mathbb R
&\implies&G\curvearrowright S^{N-1}_{\mathbb R,+},S^{N-1}_\mathbb C\\
&\implies&G\subset O_N^+,U_N\\
&\implies&G\subset O_N^+\cap U_N=O_N 
\end{eqnarray*}

Thus, we conclude that we have $G^+(S^{N-1}_\mathbb R)=O_N$, as desired.
\end{proof}

Let us discuss now the construction $U\to S$. In the classical case the situation is very simple, because the sphere $S=S^{N-1}$ appears by rotating the point $x=(1,0,\ldots,0)$ by the isometries in $U=U_N$. Moreover, the stabilizer of this action is the subgroup $U_{N-1}\subset U_N$ acting on the last $N-1$ coordinates, and so the sphere $S=S^{N-1}$ appears from the corresponding rotation group $U=U_N$ as an homogeneous space, as follows:
$$S^{N-1}=U_N/U_{N-1}$$

In functional analytic terms, all this becomes even simpler, the correspondence $U\to S$ being obtained, at the level of algebras of functions, as follows:
$$C(S^{N-1})\subset C(U_N)\quad,\quad 
x_i\to u_{1i}$$

In general now, the homogeneous space interpretation of $S$ as above fails, due to a number of subtle algebraic and analytic reasons, explained in \cite{bss} and related papers. However, we can have some theory going by using the functional analytic viewpoint, with an embedding $x_i\to u_{1i}$ as above. Let us start with the following observation:

\begin{proposition}
For the basic spheres, we have a diagram as follows,
$$\xymatrix@R=50pt@C=50pt{
C(S)\ar[r]^\Phi\ar[d]^\alpha&C(S)\otimes C(U)\ar[d]^{\alpha\otimes id}\\
C(U)\ar[r]^\Delta&C(U)\otimes C(U)
}$$
where the map on top is the affine coaction map,
$$\Phi(x_i)=\sum_jx_j\otimes u_{ji}$$
and the map on the left is given by $\alpha(x_i)=u_{1i}$.
\end{proposition}

\begin{proof}
The diagram in the statement commutes indeed on the standard coordinates, the corresponding arrows being as follows, on these coordinates:
$$\xymatrix@R=50pt@C=50pt{
x_i\ar[r]\ar[d]&\sum_jx_j\otimes u_{ji}\ar[d]\\
u_{1i}\ar[r]&\sum_ju_{1j}\otimes u_{ji}
}$$

Thus by linearity and multiplicativity, the whole the diagram commutes.
\end{proof}

As a consequence of Proposition 3.6, we have the following result:

\begin{theorem}
We have a quotient map and an inclusion as follows,
$$U\to S_U\subset S$$
with $S_U$ being the first row space of $U$, given by 
$$C(S_U)=<u_{1i}>\subset C(U)$$
at the level of the corresponding algebras of functions.
\end{theorem}

\begin{proof}
At the algebra level, we have an inclusion and a quotient map as follows:
$$C(S)\to C(S_U)\subset C(U)$$

Thus, we obtain the result, by transposing.
\end{proof}

We will prove in what follows that the inclusion $S_U\subset S$ constructed above is an isomorphism. This will produce the correspondence $U\to S$ that we are currently looking for. In order to do so, we will use the uniform integration over $S$, which can be introduced, in analogy with what happens in the classical case, in the following way:

\index{integration over spheres}

\begin{definition}
We endow each of the algebras $C(S)$ with its integration functional
$$\int_S:C(S)\to C(U)\to\mathbb C$$
obtained by composing the morphism of algebras given by 
$$x_i\to u_{1i}$$
with the Haar integration functional of the algebra $C(U)$.
\end{definition}

In order to efficiently integrate over the sphere $S$, and in the lack of some trick like spherical coordinates, we need to know how to efficiently integrate over the corresponding quantum isometry group $U$. There is a long story here, going back to the papers of Weingarten \cite{wei}, then Collins-\'Sniady \cite{csn} in the classical case, and to the more recent papers \cite{bb+}, \cite{bbc}, and then \cite{bsp}, in the quantum group case. Following \cite{bsp}, we have:

\index{Weingarten formula}
\index{Gram matrix}
\index{Weingarten matrix}

\begin{theorem}
Assuming that a compact quantum group $G\subset U_N^+$ is easy, coming from a category of partitions $D\subset P$, we have the Weingarten formula
$$\int_Gu_{i_1j_1}^{e_1}\ldots u_{i_kj_k}^{e_k}=\sum_{\pi,\sigma\in D(k)}\delta_\pi(i)\delta_\sigma(j)W_{kN}(\pi,\sigma)$$
for any indices $i_r,j_r\in\{1,\ldots,N\}$ and any exponents $e_r\in\{\emptyset,*\}$, where $\delta$ are the usual Kronecker type symbols, and where 
$$W_{kN}=G_{kN}^{-1}$$
is the inverse of the matrix $G_{kN}(\pi,\sigma)=N^{|\pi\vee\sigma|}$.
\end{theorem}

\begin{proof}
Let us arrange indeed all the integrals to be computed, at a fixed value of the exponent $k=(e_1\ldots e_k)$, into a single matrix, of size $N^k\times N^k$, as follows:
$$P_{i_1\ldots i_k,j_1\ldots j_k}=\int_Gu_{i_1j_1}^{e_1}\ldots u_{i_kj_k}^{e_k}$$

According to the construction of the Haar measure of Woronowicz \cite{wo1}, explained in chapter 2, this matrix $P$ is the orthogonal projection onto the following space:
$$Fix(u^{\otimes k})=span\left(\xi_\pi\Big|\pi\in D(k)\right)$$ 

In order to compute this projection, consider the following linear map:
$$E(x)=\sum_{\pi\in D(k)}<x,\xi_\pi>\xi_\pi$$

Consider as well the inverse $W$ of the restriction of $E$ to the following space:
$$span\left(T_\pi\Big|\pi\in D(k)\right)$$

By a standard linear algebra computation, it follows that we have:
$$P=WE$$

But the restriction of $E$ is the linear map corresponding to $G_{kN}$, so $W$ is the linear map corresponding to $W_{kN}$, and this gives the result. See \cite{bsp}.
\end{proof}

Following \cite{bgo}, we can now integrate over the spheres $S$, as follows:

\index{integration over spheres}
 
\begin{proposition}
The integration over the basic spheres is given by
$$\int_Sx_{i_1}^{e_1}\ldots x_{i_k}^{e_k}=\sum_\pi\sum_{\sigma\leq\ker i}W_{kN}(\pi,\sigma)$$
with $\pi,\sigma\in D(k)$, where $W_{kN}=G_{kN}^{-1}$ is the inverse of $G_{kN}(\pi,\sigma)=N^{|\pi\vee\sigma|}$. 
\end{proposition}

\begin{proof}
According to our conventions, the integration over $S$ is a particular case of the integration over $U$, via $x_i=u_{1i}$. By using the formula in Theorem 3.9, we obtain:
\begin{eqnarray*}
\int_Sx_{i_1}^{e_1}\ldots x_{i_k}^{e_k}
&=&\int_Uu_{1i_1}^{e_1}\ldots u_{1i_k}^{e_k}\\
&=&\sum_{\pi,\sigma\in D(k)}\delta_\pi(1)\delta_\sigma(i)W_{kN}(\pi,\sigma)\\
&=&\sum_{\pi,\sigma\in D(k)}\delta_\sigma(i)W_{kN}(\pi,\sigma)
\end{eqnarray*}

Thus, we are led to the formula in the statement.
\end{proof}

Again following \cite{ba1}, \cite{bgo}, we have the following key result:

\index{ergodicity}

\begin{theorem}
The integration functional of $S$ has the ergodicity property
$$\left(id\otimes\int_U\right)\Phi(x)=\int_Sx$$
where $\Phi:C(S)\to C(S)\otimes C(U)$ is the universal affine coaction map.
\end{theorem}

\begin{proof}
In the real case, $x_i=x_i^*$, it is enough to check the equality in the statement on an arbitrary product of coordinates, $x_{i_1}\ldots x_{i_k}$. The left term is as follows:
\begin{eqnarray*}
\left(id\otimes\int_U\right)\Phi(x_{i_1}\ldots x_{i_k})
&=&\sum_{j_1\ldots j_k}x_{j_1}\ldots x_{j_k}\int_Uu_{j_1i_1}\ldots u_{j_ki_k}\\
&=&\sum_{j_1\ldots j_k}\ \sum_{\pi,\sigma\in D(k)}\delta_\pi(j)\delta_\sigma(i)W_{kN}(\pi,\sigma)x_{j_1}\ldots x_{j_k}\\
&=&\sum_{\pi,\sigma\in D(k)}\delta_\sigma(i)W_{kN}(\pi,\sigma)\sum_{j_1\ldots j_k}\delta_\pi(j)x_{j_1}\ldots x_{j_k}
\end{eqnarray*}

Let us look now at the last sum on the right. The situation is as follows:

\medskip

(1) In the free case we have to sum quantities of type $x_{j_1}\ldots x_{j_k}$, over all choices of multi-indices $j=(j_1,\ldots,j_k)$ which fit into our given noncrossing pairing $\pi$, and just by using the condition $\sum_ix_i^2=1$, we conclude that the sum is 1. 

\medskip

(2) The same happens in the classical case. Indeed, our pairing $\pi$ can now be crossing, but we can use the commutation relations $x_ix_j=x_jx_i$, and the sum is again 1.

\medskip

Thus the sum on the right is 1, in all cases, and we obtain:
$$\left(id\otimes\int_U\right)\Phi(x_{i_1}\ldots x_{i_k})
=\sum_{\pi,\sigma\in D(k)}\delta_\sigma(i)W_{kN}(\pi,\sigma)$$

On the other hand, another application of the Weingarten formula gives:
\begin{eqnarray*}
\int_Sx_{i_1}\ldots x_{i_k}
&=&\int_Uu_{1i_1}\ldots u_{1i_k}\\
&=&\sum_{\pi,\sigma\in D(k)}\delta_\pi(1)\delta_\sigma(i)W_{kN}(\pi,\sigma)\\
&=&\sum_{\pi,\sigma\in D(k)}\delta_\sigma(i)W_{kN}(\pi,\sigma)
\end{eqnarray*}

Thus, we are done with the proof of the result, in the real case. In the complex case the proof is similar, by adding exponents everywhere. See \cite{ba1}, \cite{bgo}.
\end{proof} 

Still following \cite{ba1}, \cite{bgo}, we can now deduce a useful abstract characterization of the integration over the spheres, as follows:

\index{integration over spheres}

\begin{theorem}
There is a unique positive unital trace $tr:C(S)\to\mathbb C$ satisfying
$$(tr\otimes id)\Phi(x)=tr(x)1$$
where $\Phi$ is the coaction map of the corresponding quantum isometry group,
$$\Phi:C(S)\to C(S)\otimes C(U)$$
and this is the canonical integration, as constructed in Definition 3.8.
\end{theorem}

\begin{proof}
First of all, it follows from the Haar integral invariance condition for $U$ that the canonical integration has indeed the invariance property in the statement, namely:
$$(tr\otimes id)\Phi(x)=tr(x)1$$

In order to prove now the uniqueness, let $tr$ be as in the statement. We have:
\begin{eqnarray*}
tr\left(id\otimes\int_U\right)\Phi(x)
&=&\int_U(tr\otimes id)\Phi(x)\\
&=&\int_U(tr(x)1)\\
&=&tr(x)
\end{eqnarray*}

On the other hand, according to Theorem 3.11, we have as well:
$$
tr\left(id\otimes\int_U\right)\Phi(x)
=tr\left(\int_Sx\right)
=\int_Sx$$

We therefore conclude that $tr$ equals the standard integration, as claimed.
\end{proof}

Getting back now to our axiomatization questions, we have:

\begin{theorem}
The operation $S\to S_U$ produces a correspondence as follows,
$$\xymatrix@R=15mm@C=15mm{
S^{N-1}_{\mathbb R,+}\ar[r]&S^{N-1}_{\mathbb C,+}\\
S^{N-1}_\mathbb R\ar[r]\ar[u]&S^{N-1}_\mathbb C\ar[u]}
\qquad
\xymatrix@R=8mm@C=15mm{\\ \to}
\qquad
\xymatrix@R=17mm@C=16mm{
O_N^+\ar[r]&U_N^+\\
O_N\ar[r]\ar[u]&U_N\ar[u]
}$$
between basic unitary groups and the basic noncommutative spheres.
\end{theorem}

\begin{proof}
We use the ergodicity formula from Theorem 3.11, namely:
$$\left(id\otimes\int_U\right)\Phi=\int_S$$

We know that $\int_U$ is faithful on $\mathcal C(U)$, and that we have:
$$(id\otimes\varepsilon)\Phi=id$$

The coaction map $\Phi$ follows to be faithful as well. Thus for any $x\in\mathcal C(S)$ we have:
$$\int_Sxx^*=0\implies x=0$$

Thus $\int_S$ is faithful on $\mathcal C(S)$. But this shows that we have:
$$S=S_U$$

Thus, we are led to the conclusion in the statement.
\end{proof}

\section*{3c. Tori and reflections}

In relation with our initial goals for this chapter, we have satisfactory correspondences $S\leftrightarrow U$, and it remains to discuss the correspondence $T\to K$. Common sense suggests to get it via affine isometries as well, because in the classical case, we have:
$$K=G(T)$$

In the free case, however, things are quite tricky, with the naive formula $K=G^+(T)$ being wrong. In order to discuss this, and find the fix, we must compute the quantum isometry groups of the tori that we have. Quite surprisingly, this will lead us into subtle questions, in relation with $q=-1$ twists. To be more precise, we will need:

\index{twisting}
\index{twisted orthogonal group}
\index{twisted unitary group}

\begin{theorem}
The following constructions produce compact quantum groups,
\begin{eqnarray*}
C(\bar{O}_N)&=&C(O_N^+)\Big/\Big<u_{ij}u_{kl}=\pm u_{kl}u_{ij}\Big>\\
C(\bar{U}_N)&=&C(U_N^+)\Big/\Big<u_{ij}\dot{u}_{kl}=\pm\dot{u}_{kl}u_{ij}\Big>
\end{eqnarray*}
with the signs corresponding to anticommutation of different entries on same rows or same columns, and commutation otherwise, and where $\dot{u}$ stands for $u$ or for $\bar{u}$.
\end{theorem}

\begin{proof}
This is something well-known, coming from \cite{bbc} and subsequent papers, where these quantum groups were first introduced, the idea being as follows:

\medskip

(1) First of all, the operations $O_N\to\bar{O}_N$ and $U_N\to\bar{U}_N$ in the statement, obtained by replacing the commutation between the standard coordinates by some appropriate commutation/anticommutation, should be thought of as being $q=-1$ twists. 

\medskip

(2) However, this is not exactly the $q=-1$ twisting in the sense of Drinfeld \cite{dri} and Jimbo \cite{jim}, which produces non-semisimple objects, and so the result must be verified as such, independently of the $q=-1$ twisting literature related to \cite{dri}, \cite{jim}. 

\medskip

(3) But this is something elementary, which follows in the usual way, by considering the matrices $u^\Delta,u^\varepsilon,u^S$, defined by the same formulae as for $O_N^+,U_N^+$, and proving that these matrices satisfy the same relations as $u$. Indeed, let us first discuss the construction of the quantum group $\bar{O}_N$. We must prove that the algebra $C(\bar{O}_N)$ obtained from $C(O_N^+)$ via the relations in the statement has a comultiplication $\Delta$, a counit $\varepsilon$, and an antipode $S$. 
Regarding the construction of the comultiplication $\Delta$, let us set:
$$U_{ij}=\sum_ku_{ik}\otimes u_{kj}$$

For $j\neq k$ we have the following computation:
\begin{eqnarray*}
U_{ij}U_{ik}
&=&\sum_{s\neq t}u_{is}u_{it}\otimes u_{sj}u_{tk}+\sum_su_{is}u_{is}\otimes u_{sj}u_{sk}\\
&=&\sum_{s\neq t}-u_{it}u_{is}\otimes u_{tk}u_{sj}+\sum_su_{is}u_{is}\otimes(-u_{sk}u_{sj})\\
&=&-U_{ik}U_{ij}
\end{eqnarray*}

Also, for $i\neq k,j\neq l$ we have the following computation:
\begin{eqnarray*}
U_{ij}U_{kl}
&=&\sum_{s\neq t}u_{is}u_{kt}\otimes u_{sj}u_{tl}+\sum_su_{is}u_{ks}\otimes u_{sj}u_{sl}\\
&=&\sum_{s\neq t}u_{kt}u_{is}\otimes u_{tl}u_{sj}+\sum_s(-u_{ks}u_{is})\otimes(-u_{sl}u_{sj})\\
&=&U_{kl}U_{ij}
\end{eqnarray*}

Thus, we can define a comultiplication map for $C(\bar{O}_N)$, by setting:
$$\Delta(u_{ij})=U_{ij}$$

Regarding now the counit $\varepsilon$ and the antipode $S$, things are clear here, by using the same method, and with no computations needed, the formulae to be satisfied being trivially satisfied. We conclude that $\bar{O}_N$ is a compact quantum group.

\medskip

(4) The proof that the quantum space $\bar{U}_N$ in the statement is indeed a quantum group is similar, by adding $*$ exponents everywhere in the above computations.
\end{proof}

All the above might seem to be a bit ad-hoc, but there is way of doing the $q=-1$ twisting in a more conceptual way as well, by using representation theory and Tannakian duality. We will be back later to all this, in chapter 11 below, with full details.

\bigskip

Note in passing that all the above, while being something modest and strictly technical, needed in what follows, is also a polite way of saying that the Drinfeld-Jimbo construction \cite{dri}, \cite{jim} is somehow wrong at $q=-1$, and perhaps at other values of $q$ too. Which is of course a succulent topic, that we will keep for chapter 11. In the meantime, and as usual when it comes to controversies, we can only recommend some reading on all this. The paper of Drinfeld \cite{dri} is one of the best papers ever, and a must-read, in complement to the material from chapter 2. The original paper of Woronowicz \cite{wo1}, written as to cover the case $q>0$, and refurbished in our chapter 2 above as not to cover that $q>0$ case, due to our lack of trust in Drinfeld-Jimbo, is a must-read too. And for more on quantum groups, of all types, you have the books of Chari-Pressley \cite{cpr} and Majid \cite{maj}.

\bigskip

Speaking controversies, and being now a bit philosophers, we have been accumulating quite a few of them, throughout this book, and it is interesting to note that these are in fact all related. More precisely, the purely algebraic free versions of $\mathbb R^N,\mathbb C^N$, that we dismissed at the very beginning of this book, are interesting in connection with Drinfeld-Jimbo. Also, the Drinfeld-Jimbo construction, including the Woronowicz construction at $q>0$, is known to lead to smoothness in the sense of Connes. And so in short, by reading the present book, you not only learn about the fresh new skyscraper that we are attempting to build, but also about the ancient skyscraper nearby.

\bigskip

Now back to work, and to our axiomatization questions, we have:

\index{twisted orthogonal group}
\index{twisted unitary group}

\begin{theorem}
The quantum isometry groups of the basic tori are
$$\xymatrix@R=16mm@C=16mm{
T_N^+\ar[r]&\mathbb T_N^+\\
T_N\ar[r]\ar[u]&\mathbb T_N\ar[u]
}
\qquad
\xymatrix@R=8mm@C=15mm{\\ \to}
\qquad
\xymatrix@R=15mm@C=15.5mm{
H_N^+\ar[r]&K_N^+\\
\bar{O}_N\ar[r]\ar@{.}[u]&\bar{U}_N\ar@{.}[u]}$$
where $\bar{O}_N,\bar{U}_N$ are our standard $q=-1$ twists of $O_N,U_N$.
\end{theorem}

\begin{proof}
As a first observation, we have a mysterious lack of functoriality here, with the dotted lines standing for that, lack of inclusions there. But some quick thinking, based on our definition of the affine quantum isometry groups, tells us that there is no reason to have any kind of functoriality for such isometry groups, and so things fine. In practice now, there are 4 computations to be explained. In all cases we must find the conditions on a subgroup $G\subset U_N^+$ such that the following formula defines a coaction:
$$g_i\to\sum_jg_j\otimes u_{ji}$$

Since the coassociativity of such a map is automatic, we are left with checking that the map itself exists, and this is the same as checking that the following variables satisfy the same relations as the generators $g_i\in\Gamma$ of the discrete group $\Gamma=\widehat{T}$:
$$G_i=\sum_jg_j\otimes u_{ji}$$

(1) For $\Gamma=\mathbb Z_2^N$ the relations to be checked are as follows:
$$G_i=G_i^*\quad,\quad G_i^2=1\quad,\quad G_iG_j=G_jG_i$$

Regarding the first relation, namely $G_i=G_i^*$, by using $g_i=g_i^*$ this reads:
$$\sum_jg_j\otimes u_{ji}=\sum_jg_j\otimes u_{ji}^*$$

Now since the group generators $g_j$ are linearly independent, we obtain from this relation that we must have $u_{ij}=u_{ij}^*$ for any $i,j$. Thus, the condition on $G$ is:
$$G\subset O_N^+$$

We have the following formula, for the squares of our variables:
\begin{eqnarray*}
G_i^2
&=&\sum_{kl}g_kg_l\otimes u_{ki}u_{li}\\
&=&1+\sum_{k<l}g_kg_l\otimes(u_{ki}u_{li}+u_{li}u_{ki})
\end{eqnarray*}

We have as well the following formula, for the commutants:
\begin{eqnarray*}
\left[G_i,G_j\right]
&=&\sum_{kl}g_kg_l\otimes(u_{ki}u_{lj}-u_{kj}u_{li})\\
&=&\sum_{k<l}g_kg_l\otimes (u_{ki}u_{lj}-u_{kj}u_{li}+u_{li}u_{kj}-u_{lj}u_{ki})
\end{eqnarray*}

From the first relation we obtain $ab=-ba$ for $a\neq b$ on the same column of $u$, and by using the antipode, the same happens for rows. From the second relation we obtain:
$$[u_{ki},u_{lj}]=[u_{kj},u_{li}]\quad,\quad\forall k\neq l$$

Now by applying the antipode we obtain from this:
$$[u_{ik},u_{jl}]=[u_{jk},u_{il}]\quad,\quad\forall k\neq l$$

By relabelling, this gives the following formula:
$$[u_{ki},u_{lj}]=[u_{li},u_{kj}]\quad,\quad \forall i\neq j$$ 

Summing up, we are therefore led to the following conclusion:
$$[u_{ki},u_{lj}]=[u_{kj},u_{li}]=0\quad,\quad\forall i\neq j,k\neq l$$

Thus we must have $G\subset\bar{O}_N$, and this finishes the proof.

\medskip

(2) For $\Gamma=\mathbb Z_2^{*N}$ the relations to be checked are as follows:
$$G_i=G_i^*\quad,\quad G_i^2=1$$

As before, regarding the first relation, $G_i=G_i^*$, by using $g_i=g_i^*$ this reads:
$$\sum_jg_j\otimes u_{ji}=\sum_jg_j\otimes u_{ji}^*$$

Now since the group generators $g_j$ are linearly independent, we obtain from this relation that we must have $u_{ij}=u_{ij}^*$ for any $i,j$. Thus, the condition on $G$ is:
$$G\subset O_N^+$$

Also as before, in what regards the squares, we have:
\begin{eqnarray*}
G_i^2
&=&\sum_{kl}g_kg_l\otimes u_{ki}u_{li}\\
&=&1+\sum_{k\neq l}g_kg_l\otimes u_{ki}u_{li}
\end{eqnarray*}

Thus we obtain $G\subset H_N^+$, as claimed.

\medskip

(3) For $\Gamma=\mathbb Z^N$ the relations to be checked are as follows:
$$G_iG_i^*=G_i^*G_i=1\quad,\quad G_iG_j=G_jG_i$$

In what regards the first relation, we have the following formula:
\begin{eqnarray*}
G_iG_i^*
&=&\sum_{kl}g_kg_l^{-1}\otimes u_{ki}u_{li}^*\\
&=&1+\sum_{k\neq l}g_kg_l^{-1}\otimes u_{ki}u_{li}^*
\end{eqnarray*}

Also, we have the following formula:
\begin{eqnarray*}
G_i^*G_i
&=&\sum_{kl}g_k^{-1}g_l\otimes u_{ki}^*u_{li}\\
&=&1+\sum_{k\neq l}g_k^{-1}g_l\otimes u_{ki}^*u_{li}
\end{eqnarray*}

Finally, we have the following formula for the commutants:
\begin{eqnarray*}
\left[G_i,G_j\right]
&=&\sum_{kl}g_kg_l\otimes(u_{ki}u_{lj}-u_{kj}u_{li})\\
&=&\sum_{k<l}g_kg_l\otimes (u_{ki}u_{lj}-u_{kj}u_{li}+u_{li}u_{kj}-u_{lj}u_{ki})
\end{eqnarray*}

From the first relation we obtain $ab=-ba$ for $a\neq b$ on the same column of $u$, and by using the antipode, the same happens for rows. From the second relation we obtain:
$$[u_{ki},u_{lj}]=[u_{kj},u_{li}]\quad,\quad\forall k\neq l$$

By processing these formulae as before, in the proof of (1) above, we obtain from this that we must have $G\subset\bar{U}_N$, as claimed.

\medskip

(4) For $\Gamma=F_N$ the relations to be checked are as follows:
$$G_iG_i^*=G_i^*G_i=1$$

As before, in what regards the first relation, we have the following formula:
\begin{eqnarray*}
G_iG_i^*
&=&\sum_{kl}g_kg_l^{-1}\otimes u_{ki}u_{li}^*\\
&=&1+\sum_{k\neq l}g_kg_l^{-1}\otimes u_{ki}u_{li}^*
\end{eqnarray*}

Also as before, we have the following formula:
\begin{eqnarray*}
G_i^*G_i
&=&\sum_{kl}g_k^{-1}g_l\otimes u_{ki}^*u_{li}\\
&=&1+\sum_{k\neq l}g_k^{-1}g_l\otimes u_{ki}^*u_{li}
\end{eqnarray*}

By processing these formulae as before, in the proof of (2) above, we obtain from this that we must have $G\subset K_N^+$, as claimed.
\end{proof}

The above result is quite surprising, and does not fit with what happens in the classical case, where the classical isometry groups of the tori are the reflection groups. Thus, the above result is not exactly what we want. However, we can recycle it, as follows:

\begin{theorem}
The operation $T\to G^+(T)\cap K_N^+$ produces a correspondence
$$\xymatrix@R=15mm@C=15mm{
T_N^+\ar[r]&\mathbb T_N^+\\
T_N\ar[r]\ar[u]&\mathbb T_N\ar[u]
}
\qquad
\xymatrix@R=8mm@C=15mm{\\ \to}
\qquad
\xymatrix@R=15mm@C=15.5mm{
H_N^+\ar[r]&K_N^+\\
H_N\ar[r]\ar[u]&K_N\ar[u]}$$
between basic noncommutative tori, and basic quantum reflection groups.
\end{theorem}

\begin{proof}
The operation in the statement produces the following intersections:
$$\xymatrix@R=15mm@C=13mm{
H_N^+\ar[r]&K_N^+\\
\bar{O}_N\cap H_N^+\ar[r]\ar[u]&\bar{U}_N\cap K_N^+\ar[u]}$$

But a routine computation, coming from the fact that commutation + anticommutation means vanishing, gives the quantum groups in the statement. Indeed:

\medskip

(1) In what regards $\bar{U}_N\cap K_N^+$, here as explained above we can use the fact that commutation + anticommutation means vanishing, and we obtain, as desired:
\begin{eqnarray*}
\bar{U}_N\cap K_N^+
&=&(\bar{U}_N\cap K_N^+)_{class}\\
&=&\bar{U}_N\cap K_N\\
&=&K_N
\end{eqnarray*}

(2) In what regards $\bar{O}_N\cap H_N^+$, here we can proceed as follows:
\begin{eqnarray*}
\bar{O}_N\cap H_N^+
&=&\bar{O}_N\cap H_N^+\cap(\bar{U}_N\cap K_N^+)\\
&=&\bar{O}_N\cap H_N^+\cap K_N\\
&=&H_N
\end{eqnarray*}

Thus, we are led to the conclusion in the statement.
\end{proof}

As a conclusion to all this, we have now correspondences as follows:
$$
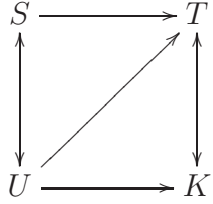
\xymatrix@R=50pt@C=50pt{
S\ar[d]\ar[r]&T\ar[d]\\
U\ar[r]\ar[ur]\ar[u]&K\ar[u]
}$$

Thus, in order to finish our axiomatization program for the abstract noncommutative geometries, we are left with constructing correspondences as follows:
$$
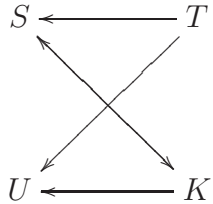
\xymatrix@R=50pt@C=50pt{
S\ar[dr]&T\ar[l]\ar[dl]\\
U&K\ar[ul]\ar[l]
}$$

We will be back to this in the next chapter, with the construction of some of these correspondences, and more specifically of those correspondences which are elementary to construct, and then with the axiomatization of the quadruplets of type $(S,T,U,K)$.

\section*{3d. Metric aspects}

Following now Goswami \cite{go1} and subsequent papers, let us comment on the ``metric'' aspects of our quantum isometry group construction. There are many things that can be said here, and the present section will be a modest introduction to all this. To start with, we have the following definition, which is something very standard in geometry:

\index{Riemannian manifold}
\index{Laplacian}
\index{Hodge Laplacian}

\begin{definition}
Given a compact Riemannian manifold $X$, we denote by $\Omega^1(X)$ the space of smooth $1$-forms on $X$, with scalar product given by
$$<\omega,\eta>=\int_X<\omega(x),\eta(x)>dx$$
and we construct the Hodge Laplacian $\Delta:L^2(X)\to L^2(X)$ by setting 
$$\Delta=d^*d$$
where $d:C^\infty(X)\to\Omega^1(X)$ is the usual differential map, and $d^*$ is its adjoint. 
\end{definition}

Observe the notational clash with the comultiplication for Woronowicz algebras, and with solving this clash being actually an open problem. Physicists like to use $\nabla^2$ for the Laplacian, but this is not very beautiful, as mathematicians we are just so used to $\Delta$. And with my hope here that, regardless of your main interests in mathematics, you teach from time to time PDE classes, as any serious mathematician should do.

\bigskip

Talking about notations, the problem comes from quantum groups, which were the last to come into play. Drinfeld, Jimbo and others used $(\Delta,\varepsilon,S)$, obviously algebra-inspired. Then Woronowicz came with $(\Phi,e,\kappa)$, for avoiding confusion with the Laplacian $\Delta$, with $\varepsilon$ from analysis, and with $S$ from Tomita-Takesaki theory. But then a bit later a younger, reckless generation came, including a former myself, thinking among others that Drinfeld-Jimbo is more interesting than Tomita-Takesaki, and reverting back to $(\Delta,\varepsilon,S)$. And here we are now, in more recent years, thinking about $\Delta$ and what to do with it.

\bigskip

But hey, to any problem there should be a solution. As already mentioned on several occasions, and more on this in a moment too, there are some deep problems in relation with the Laplacian in the noncommutative setting, namely axiomatization and general theory, including things like free harmonic functions, then heavy PDE theory to be developed, for the free manifolds, and then, as a culmination of all this, applications to physics, and more specifically, conjecturally, to QCD. And the one who will do all this will certainly have the knowledge and authority to decide what $\Delta$ should stand for.

\bigskip

Back to work now, and to Definition 3.17, we have the following standard result:

\begin{theorem}
The isometry group $\mathcal G(X)$ of a compact Riemannian manifold $X$ is the group of diffeomorphisms 
$$\varphi:X\to X$$
whose induced action on $C^\infty(X)$ commutes with the Hodge Laplacian $\Delta$. 
\end{theorem}

\begin{proof}
This is something well-known and standard, and for more on all this, basic Riemannian geometry and related topics, we refer to the book of do Carmo \cite{doc}, as well as to the book of Connes \cite{co1} and the paper of Goswami \cite{go1}.
\end{proof}

Based on the above, and following Goswami \cite{go1}, we can formulate:

\begin{definition}
The quantum isometry group $\mathcal G^+(X)$ of a compact Riemannian manifold $X$ is the biggest compact quantum group acting on $X$, via
$$\Phi:C(X)\to C(X)\otimes C(G)$$
with this coaction map commuting with the action of $\Delta$.
\end{definition}

This is something quite tricky. First, the coaction map $\Phi$ is by definition subject to the usual axioms for the algebraic coactions, namely:
$$(\Phi\otimes id)\Phi=(id\otimes\Delta)\Phi$$
$$(id\otimes\varepsilon)\Phi=id$$

In addition, $\Phi$ must be subject as well to the following smoothness assumption:
$$\Phi(C^\infty(X))\subset C^\infty(X)\otimes C(G)$$

As for the commutation condition with $\Delta$, this regards the canonical extension of the action to the space $L^2(X)$. And finally, and importantly, the above definition is something non-trivial, coming from a theorem, established in \cite{go1}, which states that a universal object $\mathcal G^+(X)$ as above exists indeed. For details here, we refer to \cite{go1}.

\bigskip

Regarding now the examples, we first have something that we know, as follows:

\begin{proposition}
The quantum isometry group of the $N$-simplex $X_N$ is
$$\mathcal G^+(X_N)=S_N^+$$
which is bigger than the usual isometry group $\mathcal G(X_N)=S_N$, at $N\geq4$.
\end{proposition}

\begin{proof}
Consider indeed the simplex $X_N\subset\mathbb R^N$, formed by definition by the standard basis $\{e_1,\ldots,e_N\}$ of $\mathbb R^N$. We know from chapter 2 that the symmetry and quantum symmetry groups of $X_N$, regarded as a set, are $S_N\subset S_N^+$. But this shows too that the classical and quantum isometry groups of $X_N$, regarded either as an algebraic manifold, as in Proposition 3.2 and Theorem 3.3, or as a Riemannian manifold, as in Theorem 3.18 and Definition 3.19, are $S_N\subset S_N^+$ as well. Finally, the fact that the inclusion $S_N\subset S_N^+$ is not an isomorphism at $N\geq4$ is something that we know too from chapter 2.
\end{proof}

It is possible to obtain more examples along the same lines, by looking at more general disconnected manifolds, and with the computation of $\mathcal G^+(X)$ for disconnected manifolds being actually a very interesting question. See \cite{go1}. In what regards the connected case, however, there have been a lot of computations by Bhowmick, Goswami and others, leading to the conjecture that we should have rigidity, $\mathcal G^+(X)=\mathcal G(X)$, in this case. And with this rigidity conjecture being now a theorem, due to Goswami \cite{go3}:

\index{Goswami theorem}
\index{eigenfunctions of the Laplacian}

\begin{theorem}
For a compact, connected Riemannian manifold $X$, the inclusion
$$\mathcal G(X)\subset\mathcal G^+(X)$$
is an isomorphism. That is, $X$ cannot have genuine quantum isometries.
\end{theorem}

\begin{proof}
There is a long story with this result, which solves an old conjecture, and whose proof is non-trivial, and for details, we refer to Goswami's paper \cite{go3}.
\end{proof}

In short, tough mathematics here, that we won't get into. This being said, in order to get a feeling for this, here is a particular case of Theorem 3.21, coming with proof:

\begin{proposition}
A compact connected Riemannian manifold $X$ cannot, in particular, have genuine group dual isometries.
\end{proposition}

\begin{proof}
Assume indeed that we have a group dual coaction, as follows:
$$\Phi:C(X)\to C(X)\otimes C^*(\Gamma)$$

Let $E=E_1\oplus E_2$ be the direct sum of two eigenspaces of the Laplacian $\Delta$. Pick a basis $\{x_i\}$ such that the corresponding corepresentation on $E$ becomes diagonal, in the sense that we have, for certain group elements $g_i\in\Gamma$:
$$\Phi(x_i)=x_i\otimes g_i$$

The formula $\Phi(x_ix_j)=\Phi(x_jx_i)$ reads then: 
$$x_ix_j\otimes g_ig_j=x_ix_j\otimes g_jg_i$$

Now since the eigenfunctions of $\Delta$ are well-known to form a domain, we obtain:
$$g_ig_j=g_jg_i$$

Similarly, $\Phi(x_i\bar{x}_j)=\Phi(\bar{x}_jx_i)$ gives $g_ig_j^{-1}=g_j^{-1}g_i$. Thus $\{g_i,g_i^{-1}\}$ pairwise commute, and with the eigenspace $E$ varying, this shows that $\Gamma$ must be abelian, as claimed. 
\end{proof}

The above is nice and fun, and there are probably some more things to be done here, along the lines of Theorem 3.4. However, as a word of warning, such ideas lead nowhere in the general context of Theorem 3.21. For the proof of that theorem, see \cite{go3}.

\bigskip

Getting back now to Goswami's foundational paper \cite{go1}, let us discuss the extension of the construction of $\mathcal G^+(X)$, to the case where $X$ is a noncommutative compact Riemannian manifold in the sense of Connes \cite{co1}. This is again heavy mathematics, with the full understanding of Connes' axiomatization in \cite{co1} requiring the reading of his subsequent ``reconstruction'' paper \cite{co2}, and also of his papers with Chamseddine \cite{cc1}, \cite{cc2} for examples and motivations, and with the work of Goswami \cite{go1} coming on top of that. So, let us be a bit informal here, and formulate things as follows:

\begin{theorem}
The theory of compact Riemannian manifolds $X$ can be extended into a theory of noncommutative compact Riemannian manifolds $X$, using spectral triples
$$X=(A,H,D)$$
in the sense of Connes. In this framework, we can talk about the corresponding quantum isometry groups $\mathcal G^+(X)$, constructed by using commutation with $D$.
\end{theorem}

\begin{proof}
This is something well-beyond the purposes of the present book, with the main references, including \cite{co1}, \cite{go1}, being those indicated above.
\end{proof}

Now back to our spheres, the free ones do not have spectral triples in the sense of Connes, but there are a few ways of talking about geometry and $\mathcal G^+(X)$, as follows:

\bigskip

(1) As explained in \cite{bgo}, it is possible to construct a Laplacian filtration for $S^{N-1}_{\mathbb R,+}$, meaning eigenspaces but no eigenvalues, and so no operator itself, as being the pullback of the Laplacian filtration for $S^{N-1}_\mathbb R$, via the embedding $S^{N-1}_\mathbb R\subset S^{N-1}_{\mathbb R,+}$. But that is enough in order to talk about $\mathcal G^+(S^{N-1}_{\mathbb R,+})$, with the result of course that this quantum group coincides with $G^+(S^{N-1}_{\mathbb R,+})=O_N^+$. For details here, and extensions, we refer to \cite{bgo}, \cite{bsk}.

\bigskip

(2) More recently, the paper of Das-Franz-Wang \cite{dfw} contains a proposal for the eigenvalues of the Laplacian of $S^{N-1}_{\mathbb R,+}$, motivated by their questions there, which is non-trivial, beautiful, and most likely correct, from a physical viewpoint. So, as a hot topic now, we have the question of extending the theory in \cite{dfw}, up to the limits of what can be done. And sky is the limit, when talking about what can be done with $\Delta$.

\bigskip

So, this is the situation, things doing well, and we refer to \cite{dfw} and related papers for more on all this. Let us mention, however, as a final comment on the subject, that something not to be ignored is Nash's theorem in \cite{nas}. That is one big result in mathematics, and extending it to the noncommutative setting is a key problem. And with this problem being not exactly ours, because are manifolds have coordinates, by definition. 

\section*{3e. Exercises} 

Things in this chapter have been a mix of basic theory and advanced mathematics, and  further meditating on all this, advanced aspects, after of course learning some more about them, is the exercise. Here are however a few concrete things, that you can try:

\begin{exercise}
Try to talk about orientation in free geometry. Also, try finding eigenvalues for the eigenspaces of the Laplacian on $S^{N-1}_{\mathbb R,+}$. Also, try constructing a spectral triple in the sense of Connes for $S^{N-1}_{\mathbb R,+}$, and explain what works, and what fails. And finally, try proving a Nash theorem for the spectral triples in the sense of Connes.
\end{exercise}

These are all difficult questions, and with the comment that browsing the internet won't help much, because the answers to most of these questions, and more specifically to those where the answer is negative, is rather folklore. And with this being such a pity, but this is the unfortunate etiquette in mathematics, and in mathematical physics too, and even in theoretical physics, don't write and publish negative results.

\chapter{Axiomatization}

\section*{4a. Basic quadruplets}

We finish here our axiomatization work. We recall that our goal is that of axiomatizing the quadruplets $(S,T,U,K)$ consisting of a quantum sphere, torus, unitary group and reflection group, with a full set of correspondences between them, as follows:
$$
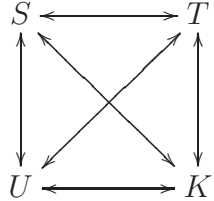
\xymatrix@R=50pt@C=50pt{
S\ar[r]\ar[d]\ar[dr]&T\ar[l]\ar[d]\ar[dl]\\
U\ar[u]\ar[ur]\ar[r]&K\ar[l]\ar[ul]\ar[u]
}$$

In order to discuss this, we first need precise definitions for all the objects involved. So, let us start with the following general definition:

\index{quantum sphere}
\index{quantum torus}
\index{quantum unitary group}
\index{quantum reflection group}

\begin{definition}
We call quantum sphere, torus, unitary group and reflection group the intermediate objects as follows,
$$S^{N-1}_\mathbb R\subset S\subset S^{N-1}_{\mathbb C,+}$$
$$T_N\subset T\subset\mathbb T_N^+$$
$$O_N\subset U\subset U_N^+$$
$$H_N\subset K\subset K_N^+$$
with $S$ being an algebraic manifold, and $T,U,K$ being compact quantum groups.
\end{definition}

Here, as usual, all the objects are taken up to the standard equivalence relation for quantum algebraic manifolds, discussed in chapter 1, as to avoid amenability issues.

\bigskip

As a first comment, what we are doing here is very straightforward, simply assuming that our objects $S,T,U,K$ lie somewhere between the minimal ones that we have, coming from $\mathbb R^N$, and the maximal possible ones, corresponding to the free analogue of $\mathbb C^N$. But this might actually seem a bit strange, because is it really a good idea to mix the real and complex cases. Good point, and we have several answers here, as follows:

\bigskip

(1) First of all, this is something technical, because we would like to deal at the same time with the real and complex cases, in order to simplify our axiomatization work. And if we want to distinguish between real and complex, we can always do that later. Although some hybrid things, for instance between $\mathbb R^N$ and $\mathbb C^N$, might be actually of interest.

\bigskip

(2) But, and here comes the point, due to some subtle reasons, we do not want in fact to distinguish between real and complex. We will see indeed later that we have an isomorphism between free projective quantum unitary groups $PO_N^+=PU_N^+$, which shows that in the free setting, the usual $\mathbb R/\mathbb C$ dichotomy tends to become ``blurred''.

\bigskip

(3) In short, our claim, which is quite bold, but is supported by some rock-solid and beautiful mathematical results, like $PO_N^+=PU_N^+$, is that the ground field for noncommutative geometry is not $F=\mathbb R$, nor $F=\mathbb C$, but rather some kind of mixture between $\mathbb R$ or $\mathbb C$. And with this to be taken, of course, in a philosophical sense.

\bigskip

(4) And with more supporting philosophy coming from physics. The field for classical mechanics is $F=\mathbb R$, while the field for quantum mechanics, or at least for QED, is rather $F=\mathbb C$. And the bet would be that the correct field for finer theories, bringing unification, should be something hybrid between $F=\mathbb R$ and $F=\mathbb C$. But more on this later.

\bigskip

Long story short, Definition 4.1 is not that bad as it is, and more on this later. At the level of the basic examples now, the situation is as follows:

\index{quadruplet}
\index{real classical geometry}
\index{complex classical grometry}
\index{free classical geometry}
\index{free complex geometry}

\begin{proposition}
We have ``basic'' quadruplets $(S,T,U,K)$ as follows:
\begin{enumerate}
\item A classical real and a classical complex quadruplet, as follows:
$$
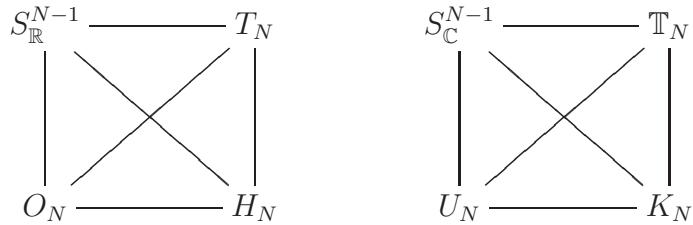
\xymatrix@R=50pt@C=50pt{
S^{N-1}_\mathbb R\ar@{-}[r]\ar@{-}[d]\ar@{-}[dr]&T_N\ar@{-}[l]\ar@{-}[d]\ar@{-}[dl]\\
O_N\ar@{-}[u]\ar@{-}[ur]\ar@{-}[r]&H_N\ar@{-}[l]\ar@{-}[ul]\ar@{-}[u]}
\qquad\qquad 
\xymatrix@R=50pt@C=50pt{
S^{N-1}_\mathbb C\ar@{-}[r]\ar@{-}[d]\ar@{-}[dr]&\mathbb T_N\ar@{-}[l]\ar@{-}[d]\ar@{-}[dl]\\
U_N\ar@{-}[u]\ar@{-}[ur]\ar@{-}[r]&K_N\ar@{-}[l]\ar@{-}[ul]\ar@{-}[u]}$$

\item A free real and a free complex quadruplet, as follows:
$$
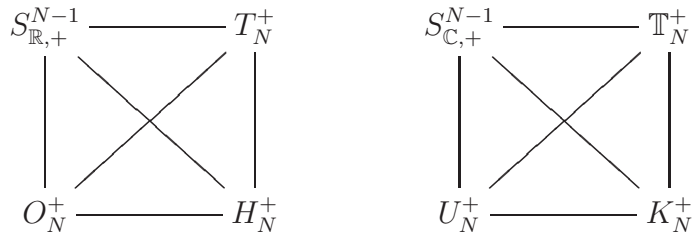
\xymatrix@R=50pt@C=50pt{
S^{N-1}_{\mathbb R,+}\ar@{-}[r]\ar@{-}[d]\ar@{-}[dr]&T_N^+\ar@{-}[l]\ar@{-}[d]\ar@{-}[dl]\\
O_N^+\ar@{-}[u]\ar@{-}[ur]\ar@{-}[r]&H_N^+\ar@{-}[l]\ar@{-}[ul]\ar@{-}[u]}
\qquad\qquad
\xymatrix@R=50pt@C=50pt{
S^{N-1}_{\mathbb C,+}\ar@{-}[r]\ar@{-}[d]\ar@{-}[dr]&\mathbb T_N^+\ar@{-}[l]\ar@{-}[d]\ar@{-}[dl]\\
U_N^+\ar@{-}[u]\ar@{-}[ur]\ar@{-}[r]&K_N^+\ar@{-}[l]\ar@{-}[ul]\ar@{-}[u]}$$
\end{enumerate}
\end{proposition}

\begin{proof}
This is more or less an empty statement, with the various objects appearing in the above diagrams being those constructed in chapters 1 and 2.
\end{proof}

Regarding now the correspondences between our objects $(S,T,U,K)$, we would like to have all 12 of them axiomatized. There is still quite some work to be done here, and in order to get started, let us begin with a summary of what we have so far:

\begin{theorem}
For the basic quadruplets, we have correspondences as follows,
$$\xymatrix@R=50pt@C=50pt{
S\ar[r]\ar[d]&T\ar[d]\\
U\ar[r]\ar[u]\ar[ur]&K\ar[u]
}$$
constructed via the following formulae:
\begin{enumerate}
\item $S=S_U$.

\item $T=S\cap\mathbb T_N^+=U\cap\mathbb T_N^+=K\cap\mathbb T_N^+$.

\item $U=G^+(S)$.

\item $K=U\cap K_N^+=K^+(T)$.
\end{enumerate}
\end{theorem}

\begin{proof}
This is a summary of what we already have, with the fact that the 7 correspondences in the statement work well for the 4 basic quadruplets, from Proposition 4.2, coming from the various results established in chapters 1-3:

\medskip

(1) The formula $S=S_U$ is from chapter 3, with the proof there being based on an ergodicity result, ultimately coming from easiness, and the Weingarten formula.

\medskip

(2) The formula $T=S\cap\mathbb T_N^+$ is from chapter 1, and this is something elementary, coming from definitions.

\medskip

(3) The formula $T=U\cap\mathbb T_N^+$ is from chapter 2, and this is once again something elementary, coming from definitions.

\medskip

(4) The formula $T=K\cap\mathbb T_N^+$ is once again from chapter 2, coming as before essentially from definitions.

\medskip

(5) The formula $U=G^+(S)$ is from chapter 3, with the proof being something quite standard, based on the tricks of Bhowmick-Goswami \cite{bgo}.

\medskip

(6) The formula $K=U\cap K_N^+$ is from chapter 2, and this is something elementary, coming from definitions.

\medskip

(7) The formula $K=K^+(T)$ is from chapter 3, and this is definitely something quite tricky, involving $q=-1$ twists.
\end{proof}

As a summary of the summary now, 7 correspondences done, 5 still to go.

\section*{4b. Easy geometries}

Our goal is that of having full correspondences between $(S,T,U,K)$. A key problem is that of passing from the discrete objects $(T,K)$ to the continuous ones $(S,U)$. We will solve this by doing some work at the quantum group level, in relation with $T,K,U$. To be more precise, we would like to have correspondences as follows:
$$T\to K\to U$$

This does not look very complicated, at the first glance, because you would say that $K,U$ can be reconstructed from the torus $T$ by some kind of Lie theory, and products, a bit as in the classical case. But, and here comes our point, in the free setting, and so in general too, we do not have Lie theory, so we must invent something else.

\bigskip

The answer will come from certain formulae involving the topological generation operation $<\,,>$ for the closed subgroups of $U_N^+$. So, let us begin by discussing this operation. This is closely related to the usual intersection operation $\cap$, again for the closed subgroups of $U_N^+$, and it is convenient to jointly talk about these two operations:

\index{intersection of quantum groups}
\index{generation operation}
\index{topological generation}

\begin{proposition}
The closed subgroups of $U_N^+$ are subject to operations as follows:
\begin{enumerate}
\item Intersection: $H\cap K$ is the biggest quantum subgroup of $H,K$.

\item Generation: $<H,K>$ is the smallest quantum group containing $H,K$.
\end{enumerate}
\end{proposition}

\begin{proof}
We must prove that the universal quantum groups in the statement exist indeed. For this purpose, let us pick writings as follows, with $I,J$ being Hopf ideals:
$$C(H)=C(U_N^+)/I\quad,\quad 
C(K)=C(U_N^+)/J$$

We can then construct our two universal quantum groups, as follows:
$$C(H\cap K)=C(U_N^+)/<I,J>$$
$$C(<H,K>)=C(U_N^+)/(I\cap J)$$

Thus, we are led to the conclusions in the statement.
\end{proof}

Let us develop now some basic theory for these operations, and for details in what follows, we refer to the book \cite{ba8}. In practice, $\cap$ can be computed by using:

\begin{proposition}
Assuming $H,K\subset G$, the intersection $H\cap K$ is given by
$$C(H\cap K)=C(G)/\{\mathcal R,\mathcal P\}$$
whenever we have formulae of type
$$C(H)=C(G)/\mathcal R\quad,\quad 
C(K)=C(G)/\mathcal P$$
with $\mathcal R,\mathcal P$ being sets of polynomial $*$-relations between the standard coordinates.
\end{proposition}

\begin{proof}
This follows from Proposition 4.4, and from the following trivial fact:
$$I=<\mathcal R>,J=<\mathcal P>
\ \implies\ <I,J>=<\mathcal R,\mathcal P>$$

Thus, we are led to the conclusion in the statement.
\end{proof}

In relation now with $<\,,>$, let us call Hopf image of a representation $C(G)\to A$ the smallest Hopf algebra quotient $C(L)$ producing a factorization as follows:
$$C(G)\to C(L)\to A$$

The fact that such a quotient exists indeed is routine, by dividing by a suitable ideal. This notion can be generalized to families of representations, and we have:

\index{Hopf image}

\begin{proposition}
Assuming $H,K\subset G$, the quantum group $<H,K>$ is such that
$$C(G)\to C(H\cap K)\to C(H),C(K)$$
is the joint Hopf image of the following quotient maps:
$$C(G)\to C(H),C(K)$$
\end{proposition}

\begin{proof}
In the particular case from the statement, the joint Hopf image appears as the smallest Hopf algebra quotient $C(L)$ producing factorizations as follows:
$$C(G)\to C(L)\to C(H),C(K)$$

Thus $L=<H,K>$, which leads to the conclusion in the statement.
\end{proof}

In the Tannakian setting now, we have the following result:

\begin{theorem}
The intersection and generation operations $\cap$ and $<\,,>$ can be constructed via the Tannakian correspondence $G\to C_G$, as follows:
\begin{enumerate}
\item Intersection: defined via $C_{G\cap H}=<C_G,C_H>$.

\item Generation: defined via $C_{<G,H>}=C_G\cap C_H$.
\end{enumerate}
\end{theorem}

\begin{proof}
This follows from Proposition 4.4, or rather from its proof, by taking $I,J$ to be the ideals coming from Tannakian duality, in its soft form, from chapter 2.
\end{proof}

In relation now with easiness, we first have the following result:

\begin{proposition}
Assuming that $H,K$ are easy, then so is $H\cap K$, and we have
$$D_{H\cap K}=<D_H,D_K>$$
at the level of the corresponding categories of partitions.
\end{proposition}

\begin{proof}
We have indeed the following computation:
\begin{eqnarray*}
C_{H\cap K}
&=&<C_H,C_K>\\
&=&<span(D_H),span(D_K)>\\
&=&span(<D_H,D_K>)
\end{eqnarray*}

Thus, by Tannakian duality we obtain the result.
\end{proof}

Regarding the generation operation, the situation is more complicated, as follows:

\begin{proposition}
Assuming that $H,K$ are easy, we have an inclusion 
$$<H,K>\subset\{H,K\}$$
coming from an inclusion of Tannakian categories as follows,
$$C_H\cap C_K\supset span(D_H\cap D_K)$$
where $\{H,K\}$ is the easy quantum group having as category of partitions $D_H\cap D_K$.
\end{proposition}

\begin{proof}
This follows from the properties of the generation operation, and from:
\begin{eqnarray*}
C_{<H,K>}
&=&C_H\cap C_K\\
&=&span(D_H)\cap span(D_K)\\
&\supset&span(D_H\cap D_K)
\end{eqnarray*}

Indeed, by Tannakian duality we obtain from this all the assertions.
\end{proof}

Summarizing, we have some problems here, and we must proceed as follows:

\index{easy generation}
\index{intersection of quantum groups}

\begin{theorem}
The intersection and easy generation operations $\cap$ and $\{\,,\}$ can be constructed via the Tannakian correspondence $G\to D_G$, as follows:
\begin{enumerate}
\item Intersection: defined via $D_{G\cap H}=<D_G,D_H>$.

\item Easy generation: defined via $D_{\{G,H\}}=D_G\cap D_H$.
\end{enumerate}
\end{theorem}

\begin{proof}
Here (1) is a result coming from Proposition 4.8, and (2) is more of an empty statement, related to the difficulties that we met in Proposition 4.9.
\end{proof}

With these generalities in hand, let us go back to our $T,U,K$ questions. Regarding $U,K$, we have the following summary of the results that we have so far, along with a few new things, in relation with the intersection and generation operations:

\begin{theorem}
The basic quantum unitary and reflection groups, namely
$$\xymatrix@R=18pt@C=18pt{
&K_N^+\ar[rr]&&U_N^+\\
H_N^+\ar[rr]\ar[ur]&&O_N^+\ar[ur]\\
&K_N\ar[rr]\ar[uu]&&U_N\ar[uu]\\
H_N\ar[uu]\ar[ur]\ar[rr]&&O_N\ar[uu]\ar[ur]
}$$
are all easy, and form an intersection/easy generation diagram, in the sense that any subsquare $P\subset Q,R\subset S$ of this diagram satisfies $Q\cap R=P$, $\{Q,R\}=S$.
\end{theorem}

\begin{proof}
We know from chapter 2 that the quantum unitary and reflection groups are all easy, the corresponding categories of partitions being as follows:
$$
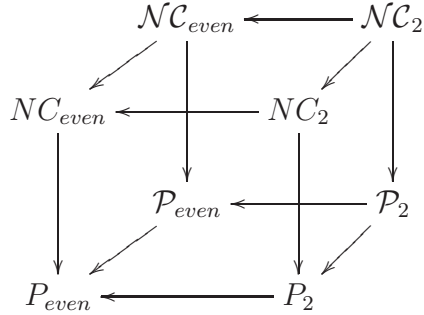
\xymatrix@R=18pt@C=2pt{
&\mathcal{NC}_{even}\ar[dl]\ar[dd]&&\mathcal {NC}_2\ar[dl]\ar[ll]\ar[dd]\\
NC_{even}\ar[dd]&&NC_2\ar[dd]\ar[ll]\\
&\mathcal P_{even}\ar[dl]&&\mathcal P_2\ar[dl]\ar[ll]\\
P_{even}&&P_2\ar[ll]
}$$

Now since these categories form an intersection and generation diagram, the quantum groups form an intersection and easy generation diagram, as claimed.
\end{proof}

Regarding now the tori $T$, the result is as follows:

\begin{theorem}
The diagonal tori of the basic unitary and reflection groups are
$$
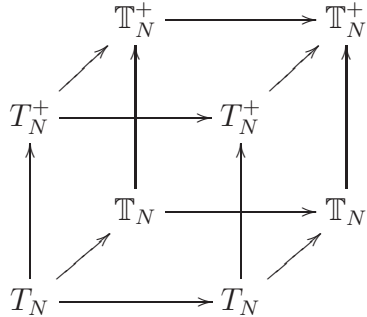
\xymatrix@R=17pt@C=17pt{
&\mathbb T_N^+\ar[rr]&&\mathbb T_N^+\\
T_N^+\ar[rr]\ar[ur]&&T_N^+\ar[ur]\\
&\mathbb T_N\ar[rr]\ar[uu]&&\mathbb T_N\ar[uu]\\
T_N\ar[uu]\ar[ur]\ar[rr]&&T_N\ar[uu]\ar[ur]
}$$
and these tori form an intersection/generation diagram, in the sense that any subsquare $P\subset Q,R\subset S$ of this diagram satisfies $Q\cap R=P$, $<Q,R>=S$.
\end{theorem}

\begin{proof}
The first assertion is something that we already know. As for the intersection and generation claim, this is something well-known, and elementary.
\end{proof}

As a first consequence of the above results, which is of interest for us, we have:

\index{reflection subgroup}
\index{liberation}

\begin{proposition}
The unitary quantum groups appear from their classical versions
$$
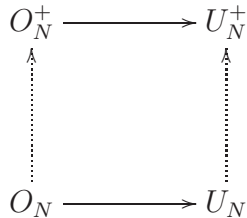
\xymatrix@R=50pt@C=50pt{
O_N^+\ar[r]&U_N^+\\
O_N\ar@.[u]\ar[r]&U_N\ar@.[u]
}$$
via $G=\{G_{class},K\}$, where $K\subset G$ is the quantum reflection subgroup, $K=G\cap K_N^+$.
\end{proposition}

\begin{proof}
We have two formulae to be established, the idea being as follows:

\medskip

(1) For the quantum group $O_N^+$ the classical version is $O_N$, and the corresponding reflection group is $H_N^+$, and from the fact that the front face of the quantum group diagram in Theorem 4.11 is an easy generation diagram we obtain, as desired:
$$O_N^+=\{O_N,H_N^+\}$$

(2) For the quantum group $U_N^+$ the classical version is $U_N$, and the corresponding reflection group is $K_N^+$, and from the fact that the rear face of the quantum group diagram in Theorem 4.11 is an easy generation diagram we obtain, as desired:
$$U_N^+=\{U_N,K_N^+\}$$

Thus, we are led to the conclusion in the statement.
\end{proof}

We can further reformulate the above result, in the following way:

\begin{proposition}
The unitary quantum groups appear from reflection subgroups
$$\xymatrix@R=52pt@C=52pt{
H_N^+\ar[r]&K_N^+\\
H_N\ar[u]\ar[r]&K_N\ar[u]
}\qquad
\xymatrix{\\ \to}
\qquad
\xymatrix@R=52pt@C=52pt{
O_N^+\ar[r]&U_N^+\\
O_N\ar[u]\ar[r]&U_N\ar[u]
}$$
via the easy generation formula $U=\{O_N,K\}$, computed inside $U_N^+$.
\end{proposition}

\begin{proof}
This is a reformulation of Proposition 4.13, as follows:

\medskip

(1) In the classical orthogonal case the formula is trivial, $O_N=\{O_N,H_N\}$.

\medskip

(2) In the free orthogonal case the formula etablished in Proposition 4.13 is precisely the one that we need, namely $O_N^+=\{O_N,H_N^+\}$.

\medskip

(3) In the classical unitary case, $U_N=\{O_N,K_N\}$ comes from the fact that the bottom face of the quantum group diagram in Theorem 4.11 is an easy generation diagram.

\medskip

(4) Finally, in the free unitary case, we have the following computation:
\begin{eqnarray*}
U_N^+
&=&\{U_N,K_N^+\}\\
&=&\{\{O_N,K_N\},K_N^+\}\\
&=&\{O_N,\{K_N,K_N^+\}\}\\
&=&\{O_N,K_N^+\}
\end{eqnarray*}

Thus, we are led to the conclusion in the statement.
\end{proof}

We can now update our main result, with 1 more correspondence, as follows:

\begin{theorem}
For the basic quadruplets, we have correspondences as follows,
$$\xymatrix@R=50pt@C=50pt{
S\ar[r]\ar[d]&T\ar[d]\\
U\ar[r]\ar[ur]\ar[u]&K\ar[u]\ar[l]
}$$
constructed via the following formulae:
\begin{enumerate}
\item $S=S_U$.

\item $T=S\cap\mathbb T_N^+=U\cap\mathbb T_N^+=K\cap\mathbb T_N^+$.

\item $U=G^+(S)=\{O_N,K\}$.

\item $K=U\cap K_N^+=K^+(T)$.
\end{enumerate}
\end{theorem}

\begin{proof}
This is an update of our main result so far, namely Theorem 4.3, by taking into account the findings from Proposition 4.14.
\end{proof}

Regarding the missing correspondences, namely $T\to S,U$ and $S\leftrightarrow K$, the situation here is more complicated, and we will discuss this later. We can however compose the correspondences that we have, and formulate, as a conclusion to what we did so far:

\index{quadruplet}
\index{easy geometry}

\begin{definition}
A quadruplet $(S,T,U,K)$ is said to produce an easy geometry when $U,K$ are easy, and one can pass from each object to all the other objects, as follows,
$$\begin{matrix}
S&=&S_{\{O_N,K^+(T)\}}&=&S_U&=&S_{\{O_N,K\}}\\
\\
S\cap\mathbb T_N^+&=&T&=&U\cap\mathbb T_N^+&=&K\cap\mathbb T_N^+\\
\\
G^+(S)&=&\{O_N,K^+(T)\}&=&U&=&\{O_N,K\}\\
\\
K^+(S)&=&K^+(T)&=&U\cap K_N^+&=&K
\end{matrix}$$
with the usual convention that all this is up to the equivalence relation.
\end{definition}

Observe that if we plug the data from any axiom line into the 3 other lines, we obtain axiomatizations in terms of one of $S,T,U,K$, that we can try to simplify afterwards. It is of course possible to axiomatize everything in terms of $ST,SU,SK,TU,TK,UK$ as well, and also in terms of $STU,STK,SUK,TUK$, and try to simplify afterwards.

\bigskip

In what follows we will not bother much with this, and use Definition 4.16 as it is. We will need that 12 correspondences, as results, and whether we call such results ``verifications of the axioms'' or ``basic properties of our geometry'' is irrelevant.

\bigskip

Regarding now the basic examples, we have here the following result:

\index{real classical geometry}
\index{complex classical grometry}
\index{free classical geometry}
\index{free complex geometry}

\begin{theorem}
We have $4$ basic easy geometries, denoted
$$\xymatrix@R=50pt@C=50pt{
\mathbb R^N_+\ar[r]&\mathbb C^N_+\\
\mathbb R^N\ar[u]\ar[r]&\mathbb C^N\ar[u]
}$$
which appear from quadruplets as above, as follows:
\begin{enumerate}
\item Classical real: produced by $(S^{N-1}_\mathbb R,T_N,O_N,H_N)$.

\item Classical complex: produced by $(S^{N-1}_\mathbb C,\mathbb T_N,U_N,K_N)$.

\item Free real: produced by $(S^{N-1}_{\mathbb R,+},T_N^+,O_N^+,H_N^+)$.

\item Free complex: produced by $(S^{N-1}_{\mathbb C,+},\mathbb T_N^+,U_N^+,K_N^+)$.
\end{enumerate}
\end{theorem}

\begin{proof}
This is something that we already know, which follows from Theorem 4.15, as explained in the discussion preceding Definition 4.16.
\end{proof}

It is possible to construct some further easy geometries in the above sense, and also to work out some classification results. To be more precise, the 4-diagram of geometries from Theorem 4.17 can be extended into a 9-diagram of geometries, as follows:
$$\xymatrix@R=40pt@C=40pt{
\mathbb R^N_+\ar[r]&\mathbb T\mathbb R^N_+\ar[r]&\mathbb C^N_+\\
\mathbb R^N_*\ar[u]\ar[r]&\mathbb T\mathbb R^N_*\ar[u]\ar[r]&\mathbb C^N_*\ar[u]\\
\mathbb R^N\ar[u]\ar[r]&\mathbb T\mathbb R^N\ar[u]\ar[r]&\mathbb C^N\ar[u]
}$$

Here the $*$ symbols stand for half-liberation, obtained by replacing the commutation relations $ab=ba$ with the half-commutation relations $abc=cba$. As for the products by $\mathbb T$, producing hybrid objects, between real and complex, these are quite standard, usual products by $\mathbb T$. Moreover, one can prove that, under some supplementary assumptions, these 9 easy geometries are the only ones. We will be back to this.

\section*{4c. Liberation theory}

Moving ahead now, if we want to improve the above, we have two problems which are still in need to be solved. First, we would like to understand the operation $K\to U$, without reference to easiness. And second, we would like to understand the operation $T\to U$. In short, we are back to the problem mentioned after Theorem 4.3, namely understanding the following operations, and this time without reference to easiness:
$$T\to K\to U$$

This is something quite subtle, which will take us into advanced quantum group theory. Let us start our discussion with the following definition:

\index{diagonal torus}
\index{reflection subgroup}
\index{soft liberation}
\index{hard liberation}
\index{liberation}

\begin{definition}
Consider a closed subgroup $G\subset U_N^+$, and let 
$$T\subset K\subset G$$
be its diagonal torus, and its reflection subgroup. The inclusion $G_{class}\subset G$ is called:
\begin{enumerate}
\item A soft liberation, when $G=<G_{class},K>$.

\item A hard liberation, when $G=<G_{class},T>$.
\end{enumerate}
\end{definition}

As a first remark, in relation with these notions, given a closed subgroup $G\subset U_N^+$ we have a diagram as follows, which is an intersection diagram:
$$\xymatrix@R=50pt@C=50pt{
T\ar[r]&K\ar[r]&G\\
T_{class}\ar[u]\ar[r]&K_{class}\ar[u]\ar[r]&G_{class}\ar[u]
}$$

With this picture in mind, the soft liberation condition states that the square on the right is a generation diagram. As for the hard liberation condition, which is something stronger, this states that the whole rectangle has the generation property.

\bigskip

Although many interesting subgroups $G\subset U_N^+$ appear as hard liberations, and we will see examples in a moment, we cannot expect this to happen in general. Indeed, a basic counterexample here is the quantum permutation group $G=S_N^+$, which is bigger than $G_{class}=S_N$ at $N\geq4$, but whose diagonal torus is trivial, $T=\{1\}$.

\bigskip

In order to comment on this, let us discuss now some weaker versions of the hard liberation property, involving spinned versions of the diagonal torus. We first have:

\index{spinned torus}

\begin{proposition}
Given a closed subgroup $G\subset U_N^+$ and a matrix $Q\in U_N$, we let $T_Q\subset G$ be the diagonal torus of $G$, with fundamental representation spinned by $Q$:
$$C(T_Q)=C(G)\Big/\left<(QuQ^*)_{ij}=0\Big|\forall i\neq j\right>$$
This torus is then a group dual, given by $T_Q=\widehat{\Lambda}_Q$, where $\Lambda_Q=<g_1,\ldots,g_N>$ is the discrete group generated by the elements 
$$g_i=(QuQ^*)_{ii}$$
which are unitaries inside $C(T_Q)$. We call these tori $T_Q$ the standard tori of $G$.
\end{proposition}

\begin{proof}
This follows from the general results for the diagonal torus from chapter 2, because, as said in the statement, $T_Q$ is by definition a certain diagonal torus. Equivalently, without using anything, since $v=QuQ^*$ is a unitary corepresentation, its diagonal entries $g_i=v_{ii}$, when regarded inside $C(T_Q)$, are unitaries, and satisfy:
$$\Delta(g_i)=g_i\otimes g_i$$

Thus $C(T_Q)$ is a group algebra, and more specifically we have $C(T_Q)=C^*(\Lambda_Q)$, where $\Lambda_Q=<g_1,\ldots,g_N>$ is the group in the statement, and this gives the result.
\end{proof}

The interest in the standard tori comes from the following result:

\begin{theorem}
Any torus $T\subset G$ appears as follows, for a certain $Q\in U_N$:
$$T\subset T_Q\subset G$$
In other words, any torus appears inside a standard torus.
\end{theorem}

\begin{proof}
Given a torus $T\subset G$, we have an inclusion as follows:
$$T\subset G\subset U_N^+$$

On the other hand, we know from chapter 2 that each torus $T\subset U_N^+$ has a fundamental corepresentation as follows, with $Q\in U_N$:
$$u=Q
\begin{pmatrix}
g_1\\
&\ddots\\
&&g_N
\end{pmatrix}
Q^*$$

But this shows that we have $T\subset T_Q$, and this gives the result.
\end{proof}

As an immediate consequence of the above result, we have:

\begin{proposition}
Let $G\subset U_N^+$ be a closed subgroup.
\begin{enumerate}
\item If $G$ is classical, its maximal tori $T\subset G$ are among the standard tori $T_Q$.

\item If $G$ is a group dual, or torus, $G$ itself is among its standard tori $T_Q$.
\end{enumerate}
\end{proposition}

\begin{proof}
Both these assertions follow from Theorem 4.20, and with the remark that we can proceed directly as well. Indeed, (1) follows by jointly diagonalizing the matrices $U\in T$, which produces a certain matrix $Q\in U_N$, as needed. As for (2), as mentioned in the proof of Theorem 4.20, this is something that we know from chapter 2. 
\end{proof}

Summarizing, associated to any closed subgroup $G\subset U_N^+$ is a whole family of tori, indexed by the unitaries $U\in U_N$, and this suggests the following definition:

\index{standard tori}
\index{skeleton}

\begin{definition}
Given a closed subgroup $G\subset U_N^+$, the collection of tori 
$$T=\left\{T_Q\subset G\big|Q\in U_N\right\}$$
which plays the role of a ``maximal torus'' for $G$, is called skeleton of $G$.
\end{definition}

Here the maximal torus claim comes from Proposition 4.21, and from the many more things that can be said, as theorems or conjectures, relating the combinatorics of $G$ to the combinatorics of $T$, and for a discussion here we refer to \cite{ba8}.

\bigskip

Getting back now to our generation questions, from Definition 4.18 and the comments afterwards, the notion of hard liberation involves the diagonal torus $T_1$, and in case that fails, as is for instance the case for $G=S_N^+$, the idea is to use more standard tori $T_Q$. And here we can use all such standard tori $T_Q$, or just a suitable selection of them.

\bigskip

Some further looking at the case $G=S_N^+$, which is quite technical and that we will not get into here in detail, suggests using the Fourier tori. Let us start with:

\begin{definition}
Associated to any finite abelian group $L$ is the Fourier transform
$$F_L:C(L)\to C^*(L)$$
which can be regarded as a usual unitary matrix, $F_L\in U_N$, where $N=|L|$.
\end{definition}

To be more precise, this is something that we basically know from chapter 1, when talking about finite abelian groups, and Pontrjagin duality for them. In practice, for the cyclic group $\mathbb Z_N$ the corresponding Fourier matrix is as follows, with $w=e^{2\pi i/N}$:
$$F_N=\frac{1}{\sqrt{N}}\,(w^{ij})_{ij}$$

In general, if can write $L=\mathbb Z_{N_1}\times\ldots\times\mathbb Z_{N_k}$, the corresponding Fourier matrix is:
$$F_L=F_{N_1}\otimes\ldots\otimes F_{N_k}$$

Consider now the set $\mathcal F_N\subset U_N$ formed by all these Fourier matrices, coming from the various abelian groups $L$ satisfying $|L|=N$:
$$\mathcal F_N=\left\{F_L\Big|\,|L|=N\right\}$$

With this notion in hand, we can now formulate, as a complement to Definition 4.18 above, the following technical definition:

\index{diagonal liberation}

\begin{definition}
A closed subgroup $G\subset U_N^+$ is called:
\begin{enumerate}
\item Generated by its tori, when $G=<(T_Q)_{Q\in U_N}>$.

\item Weakly generated by its tori, when $G=<G_{class},(T_Q)_{Q\in U_N}>$.

\item A Fourier liberation of $G_{class}$, when $G=<G_{class},(T_F)_{F\in\mathcal F_N}>$.
\end{enumerate}
\end{definition}

Obviously, this is something a bit tricky, and we are in fact now into wild quantum group territory. As a first observation, which is something clear, the relation of the above notions, GT, WG, FL with the notions SL, HL from Definition 4.18 is as follows:
$$
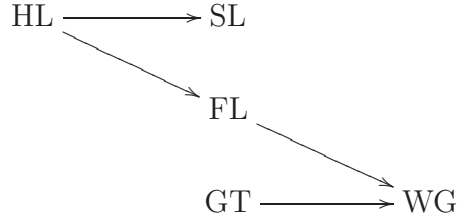
\xymatrix@R=20pt@C=50pt{
{\rm HL}\ar[r]\ar[dr]&{\rm SL}\\
&{\rm FL}\ar[dr]\\
&{\rm GT}\ar[r]&{\rm WG}
}$$

There are many things, either theorems or conjectures, that can be said about these properties, and with the study here being a subject of active research. It is known for instance that $S_N^+$ satisfies FL, and the conjecture is that any easy quantum group should satisfy FL, and also SL and GT. As for the general case, that of the arbitrary closed subgroups $G\subset U_N$, many of them are known to satisfy GT, but a safer conjecture would be that these satisfy WG. We refer to the book \cite{ba8} for a discussion of all this.

\bigskip

Long story short, the notions of soft and hard liberation from Definition 4.18 are just the tip of the iceberg, and the whole subject is quite technical. In relation now with the quantum groups that we are interested in, the result that we will need is as follows:

\index{soft liberation}
\index{hard liberation}

\begin{theorem}
The following happen:
\begin{enumerate}
\item $O_N^+,U_N^+$ appear as soft liberations of $O_N,U_N$.

\item $O_N^+,U_N^+$ appear as well as hard liberations of $O_N,U_N$.

\item $H_N^+,K_N^+$ appear as soft liberations of $H_N,K_N$.

\item $H_N^+,K_N^+$ do not appear as hard liberations of $H_N,K_N$.
\end{enumerate}
\end{theorem}

\begin{proof}
This result, that we will need in what follows, is something quite technical. In the lack of a simple and complete proof for all this, here is the idea:

\medskip

(1) This simply follows from (2) below. Normally there should be a simpler proof for this, by using Tannakian duality, but this is something which is not known yet.

\medskip

(2) A key result of Chirvasitu \cite{chi}, whose proof is quite technical, not to be explained here, states that we have the following generation formula, valid at any $N\geq3$:
$$O_N^+=<O_N,O_{N-1}^+>$$

With this in hand, the hard liberation formula $O_N^+=<O_N,T_N^+>$ can be proved by recurrence on $N$. Indeed, at $N=1$ there is nothing to prove, at $N=2$ this is something well-known, and elementary, as explained for instance in \cite{chi}, and in general, the recurrence step $N-1\to N$ can be established as follows:
\begin{eqnarray*}
O_N^+
&=&<O_N,O_{N-1}^+>\\
&=&<O_N,O_{N-1},T_{N-1}^+>\\
&=&<O_N,T_{N-1}^+>\\
&=&<O_N,T_N,T_{N-1}^+>\\
&=&<O_N,T_N^+>
\end{eqnarray*}

Regarding now $U_N^+=<U_N,\mathbb T_N^+>$, this follows from $O_N^+=<O_N,T_N^+>$, via standard lifting techniques. In order to discuss this, let us first examine the passage $O_N^+\to U_N^+$. We can construct a ``free complexification'' of $O_N^+$, as follows:
$$C(\widetilde{O_N^+})=<\tilde{v}_{ij}>\subset C(\mathbb T)*C(O_N)\quad,\quad \tilde{v}=zv$$

To be more precise, consider the free product on the right. If we denote by $z\in C(\mathbb T)$ the standard generator, which is the function $x\to x$, then $\tilde{v}=zv$ is a corepresentation of this free product, and so the algebra $<\tilde{v}_{ij}>$ is a Woronowicz algebra. The corresponding compact quantum group is denoted $\widetilde{O_N^+}$, and is called free complexification of $O_N^+$. Now, our claim is that the following embedding is an isomorphism:
$$\widetilde{O_N^+}\subset U_N^+$$

This claim can be proved in several ways. With the technology that we have so far, the simplest is to invoke easiness. Indeed, if we denote by $u$ the standard corepresentation of $U_N^+$, then the following inclusion of vector spaces follows to be an isomorphism, due to the fact that both spaces involved appear as the span of the same pairings:
$$Hom(u^{\otimes k},u^{\otimes l})\subset Hom((zv)^{\otimes k},(zv)^{\otimes l})$$

Thus, our embedding $\widetilde{O_N^+}\subset U_N^+$ has the property that it preserves the Hom spaces for the Peter-Weyl corepresentations, and from this, it follows from the Peter-Weyl theory, explained in chapter 2, that our embedding must be an isomorphism:
$$\widetilde{O_N^+}=U_N^+$$

Alternatively, a simple but advanced argument, based on free probability theory from \cite{vdn}, is that of saying that the characters of $\tilde{v}=zv$ and $u$ follow the same law, namely the Voiculescu circular law, and so the embedding $\widetilde{O_N^+}\subset U_N^+$ must be an isomorphism. In any case, either way we have proved our claim, and as a consequence of it, we have:
$$PO_N^+=PU_N^+$$

To be more precise, let us define the projective version of a closed subgroup $G\subset U_N^+$, with standard coordinates $u_{ij}$, to be the compact quantum group having as coordinates the variables $w_{ia,jb}=u_{ij}u_{ab}^*$. With this convention, we have then, as desired:
$$PU_N^+=P\widetilde{O_N^+}=PO_N^+$$

There are some alternative proofs as well of this fact, either by using the embedding $PO_N^+\subset PU_N^+$ and an easiness argument, a bit as before for $\widetilde{O_N^+}\subset U_N^+$, or, more quickly and conceptually, by using the embeddings $PO_N^+\subset PU_N^+\subset S_{M_N}^+$, with on the right a generalized quantum symmetry group, and a free probability argument, involving the Marchenko-Pastur law, which shows that our embeddings must be isomorphisms. We will be back to all this at various places, throughout this book, and in the meantime, for details we refer to \cite{ba8}. Now the point is that, as explained in \cite{chi}, from $PO_N^+=PU_N^+$ we obtain, via standard lifting arguments, that we have an isomorphism as follows:
$$U_N^+=<U_N,O_N^+>$$

By using this isomorphism and $O_N^+=<O_N,T_N^+>$, we obtain, as desired:
\begin{eqnarray*}
U_N^+
&=&<U_N,O_N^+>\\
&=&<U_N,O_N,T_N^+>\\
&=&<U_N,T_N^+>\\
&=&<U_N,\mathbb T_N^+>
\end{eqnarray*}

Still with me, I hope. There has been a lot of theory here, but as promised, we will be back to all this at various places, throughout this book.

\medskip

(3) This is something trivial, because $H_N^+,K_N^+$ equal their reflection subgroups.

\medskip

(4) This result, which is something quite surprising, is well-known, coming from the fact that the quantum group $H_N^{[\infty]}\subset H_N^+$ constructed by Raum-Weber in  \cite{rwe}, and its unitary counterpart $K_N^{[\infty]}\subset K_N^+$, have the same diagonal subgroups as $H_N^+,K_N^+$. Thus, the hard liberation procedure stops at $H_N^{[\infty]},K_N^{[\infty]}$, and cannot reach $H_N^+,K_N^+$. To be more precise, we can construct quantum groups $H_N^{[\infty]},K_N^{[\infty]}$ by using the relations $\alpha\beta\gamma=0$, for any $a\neq c$ on the same row or column of $u$, with the convention $\alpha=a,a^*$, and so on. These quantum groups appear as intermediate liberations, as follows: 
$$\xymatrix@R=15mm@C=17mm{
K_N\ar[r]&K_N^*\ar[r]&K_N^{[\infty]}\ar[r]&K_N^+\\
H_N\ar[r]\ar[u]&H_N^*\ar[r]\ar[u]&H_N^{[\infty]}\ar[r]\ar[u]&H_N^+\ar[u]}$$

Moreover, these quantum groups $H_N^{[\infty]},K_N^{[\infty]}$are easy, the corresponding categories $P_{even}^{[\infty]}\subset P_{even}$ and $\mathcal P_{even}^{[\infty]}\subset\mathcal P_{even}$ being generated by the following partition:
$$\eta=\ker\binom{i\ i\ j}{j\ i\ i}$$

In relation with our questions, since the relations $g_ig_ig_j=g_jg_ig_i$ are trivially satisfied for real reflections, the diagonal tori of these quantum groups coincide with those for $H_N^+,K_N^+$. Thus, the diagonal liberation procedure ``stops'' at $H_N^{[\infty]},K_N^{[\infty]}$.
\end{proof}

Now back to our axiomatization questions, as a first comment, in contrast to what happens in the classical case, where $K=<H_N,T>$, the correspondence $T\to K$ cannot be constructed via the hard generation formula $K=<H_N,T>$, due to Theorem 4.25 (4). Thus, our formula $K=K^+(T)$ is the only solution to the $T\to K$ probem.

\bigskip

As a second comment, all the above is interesting in connection with the cube formed by the quantum unitary and reflection groups. Let us recall from Theorem 4.11 that these quantum groups form an intersection and easy generation diagram, as follows:
$$
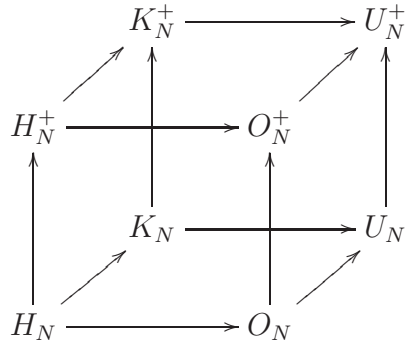
\xymatrix@R=20pt@C=20pt{
&K_N^+\ar[rr]&&U_N^+\\
H_N^+\ar[rr]\ar[ur]&&O_N^+\ar[ur]\\
&K_N\ar[rr]\ar[uu]&&U_N\ar[uu]\\
H_N\ar[uu]\ar[ur]\ar[rr]&&O_N\ar[uu]\ar[ur]
}$$

It is conjectured that this diagram should be a plain generation diagram, and the above results prove this conjecture for 5 of the faces. For the remaining face, namely the one on the left, the corresponding formula $K_N^+=<K_N,H_N^+>$ is not proved yet.

\bigskip

As yet another comment, the material in Theorem 4.25 is definitely waiting for more study. Indeed, we have the following Tannakian formulae:
$$C_{H\cap K}=<C_H,C_K>$$
$$C_{<H,K>}=C_H\cap C_K$$

Thus, from a Tannakian viewpoint, all the above results ultimately correspond to doing some combinatorics. To be more precise, the soft and hard generation properties in Definition 4.18 amount respectively in proving the following formulae:
$$C_G=<C_G,C_{U_N}>\cap\,C_K$$
$$C_G=<C_G,C_{U_N}>\cap\,C_T$$

In the easy case now, where $C_G=span(D)$, which is the case for the various quantum groups from Theorem 4.25, these two equalities reformulate as follows:
$$span(D)=span(D,\slash\hskip-2.1mm\backslash)\cap C_K$$
$$span(D)=span(D,\slash\hskip-2.1mm\backslash)\cap C_T$$

Thus, we are led into some combinatorics, which remains to be understood, in a direct way, without reference to algebra and recurrence methods. Open problem.

\bigskip

Getting back now to our axiomatization questions, we have the following refinement of Proposition 4.13, making no reference to easiness:

\index{hard generation}

\begin{proposition}
The unitary quantum groups appear from diagonal subgroups
$$\xymatrix@R=52pt@C=52pt{
T_N^+\ar[r]&\mathbb T_N^+\\
T_N\ar[u]\ar[r]&\mathbb T_N\ar[u]
}\qquad
\xymatrix{\\ \to}
\qquad
\xymatrix@R=52pt@C=52pt{
O_N^+\ar[r]&U_N^+\\
O_N\ar[u]\ar[r]&U_N\ar[u]
}$$
via the hard generation formula $U=<O_N,T>$, computed inside $U_N^+$.
\end{proposition}

\begin{proof}
This comes from the results in Theorem 4.25, as follows:

\medskip

(1) In the classical real case the condition to be verified is trivial, namely:
$$O_N=<O_N,T_N>$$

(2) In the free real case the condition to be verified is $O_N^+=<O_N,T_N^+>$. But this is exactly the hard liberation property of $O_N\subset O_N^+$, as explained in Theorem 4.25.

\medskip

(3) In the classical complex case the condition to be verified is as follows:
$$U_N=<O_N,\mathbb T_N>$$

But this is something well-known, coming for instance from the fact that the inclusion of compact Lie groups $\mathbb TO_N\subset U_N$ is maximal. For more details here, we refer to \cite{ba8}.

\medskip

(4) In the free complex case the condition is $U_N^+=<O_N,\mathbb T_N^+>$. But this comes from the hard liberation formula $U_N^+=<U_N,T_N^+>$ from Theorem 4.25, as follows:
\begin{eqnarray*}
U_N^+
&=&<U_N,\mathbb T_N^+>\\
&=&<O_N,\mathbb T_N,\mathbb T_N^+>\\
&=&<O_N,\mathbb T_N^+>
\end{eqnarray*}

Thus, we are led to the conclusions in the statement.
\end{proof}

Generally speaking, the same comments as those after Theorem 4.25 apply. In Tannakian formulation, the equalities to be proved are as follows:
$$C_U=span(P_2)\cap C_K$$
$$C_U=span(P_2)\cap C_T$$

Thus, we are led into some combinatorics, of basically the same type as the combinatorics needed for Theorem 4.25, which remains to be understood.

\bigskip

Alternatively, we have the question of understanding the formulae in Proposition 4.26 by using recurrence arguments, as in the proof of Theorem 4.25. But here, again, there is no simple proof known, at least so far. As before, we refer to \cite{ba8} for more on this.

\bigskip

We can now update our main result from the general, non-easy case, as follows:

\begin{theorem}
For the basic quadruplets, we have correspondences as follows,
$$\xymatrix@R=50pt@C=50pt{
S\ar[r]\ar[d]&T\ar[d]\ar[dl]\\
U\ar[r]\ar[ur]\ar[u]&K\ar[u]\ar[l]
}$$
constructed via the following formulae:
\begin{enumerate}
\item $S=S_U$.

\item $T=S\cap\mathbb T_N^+=U\cap\mathbb T_N^+=K\cap\mathbb T_N^+$.

\item $U=G^+(S)=<O_N,T>=<O_N,K>$.

\item $K=U\cap K_N^+=K^+(T)$.
\end{enumerate}
\end{theorem}

\begin{proof}
This is an update of our main result so far, namely Theorem 4.3, by taking into account the correspondences coming from Proposition 4.26.
\end{proof}

Summarizing, we have now a reasonable set of correspondences between our objects, which are constructed without any reference to easiness. Of course, all this was quite technical, and for further details we refer to \cite{ba8}. And with the remark that, in practice, all our geometries will be easy anyway, and so the easy geometry formalism developed in the previous section would do too, for most of what we want to do in this book.

\section*{4d. General axioms}

As already mentioned before, in chapter 1 and afterwards, in what regards the missing correspondences, $T\to S$ and $S\leftrightarrow K$, the situation here is quite complicated. In short, we have to give up now with our general principle of constructing all the correspondences independently of each other, and compose what we have. We are led to:

\index{quadruplet}

\begin{definition}
A quadruplet $(S,T,U,K)$ is said to produce a noncommutative geometry when one can pass from each object to all the other objects, as follows,
$$\begin{matrix}
S&=&S_{<O_N,T>}&=&S_U&=&S_{<O_N,K>}\\
\\
S\cap\mathbb T_N^+&=&T&=&U\cap\mathbb T_N^+&=&K\cap\mathbb T_N^+\\
\\
G^+(S)&=&<O_N,T>&=&U&=&<O_N,K>\\
\\
K^+(S)&=&K^+(T)&=&U\cap K_N^+&=&K
\end{matrix}$$
with the usual convention that all this is up to the equivalence relation.
\end{definition}

Observe the similarity with the axiomatics from the easy case. The same comments as those made after Definition 4.16, from the easy case, apply. 

\bigskip

To be more precise, if we plug the data from any axiom line into the 3 other lines, we obtain axiomatizations in terms of one of $S,T,U,K$, that we can try to simplify afterwards.  It is of course possible to axiomatize everything in terms of $ST,SU,SK,TU,TK,UK$ as well, and also in terms of $STU,STK,SUK,TUK$, and try to simplify afterwards. 

\bigskip

In what follows we will not bother much with this, and use Definition 4.28 as it is. We will need that 12 correspondences, as results, and whether we call such results ``verifications of the axioms'' or ``basic properties of our geometry'' is irrelevant.

\bigskip

Observe also that the above definition is independent from Definition 4.16, in the sense that an easy geometry in the sense of Definition 4.16 does not automatically satisfy the above axioms, or vice versa. However, we do not know counterexamples here.

\bigskip

As another technical comment, the previous work on this subject, that I did with Bichon in \cite{bb2}, was based on $(S,T,U)$ triples, but as explained there, this formalism, missing a lot of restrictions coming from $K$, is a bit too broad. As for the subsequent work, from a paper written by myself coming after \cite{bb2}, this was based on sextuplets $(S,\bar{S},T,U,\bar{U},K)$, with the bars  standing for twists, which is perhaps something quite natural, but which leads to too many correspondences between objects, namely 30. So, everything simpler or more complicated than Definition 4.28 has basically been tried, and as a conclusion to all this, Definition 4.28 is the good compromise for everything.

\bigskip

Regarding now the basic examples, these are of course the classical and free, real and complex geometries. To be more precise, we have the following result:

\index{real classical geometry}
\index{complex classical grometry}
\index{free classical geometry}
\index{free complex geometry}

\begin{theorem}
We have $4$ basic geometries, denoted
$$\xymatrix@R=50pt@C=50pt{
\mathbb R^N_+\ar[r]&\mathbb C^N_+\\
\mathbb R^N\ar[u]\ar[r]&\mathbb C^N\ar[u]
}$$
which appear from quadruplets as above, as follows:
\begin{enumerate}
\item Classical real: produced by $(S^{N-1}_\mathbb R,T_N,O_N,H_N)$.

\item Classical complex: produced by $(S^{N-1}_\mathbb C,\mathbb T_N,U_N,K_N)$.

\item Free real: produced by $(S^{N-1}_{\mathbb R,+},T_N^+,O_N^+,H_N^+)$.

\item Free complex: produced by $(S^{N-1}_{\mathbb C,+},\mathbb T_N^+,U_N^+,K_N^+)$.
\end{enumerate}
\end{theorem}

\begin{proof}
This is something that we already know, which follows from Theorem 4.27, as explained in the discussion preceding Definition 4.28.
\end{proof}

We will be back to more examples in chapters 9-12 below, and with some  classification results as well, the idea being that of looking for intermediate geometries on the horizontal, and on the vertical of the above diagram, and then combining these constructions. The conclusion there will be that the basic 4-diagram of geometries from Theorem 4.29 can be extended into a ``second level'' basic 9-diagram of geometries, as follows:
$$\xymatrix@R=40pt@C=40pt{
\mathbb R^N_+\ar[r]&\mathbb T\mathbb R^N_+\ar[r]&\mathbb C^N_+\\
\mathbb R^N_*\ar[u]\ar[r]&\mathbb T\mathbb R^N_*\ar[u]\ar[r]&\mathbb C^N_*\ar[u]\\
\mathbb R^N\ar[u]\ar[r]&\mathbb T\mathbb R^N\ar[u]\ar[r]&\mathbb C^N\ar[u]
}$$

But more on this later. Getting back now to abstract things, and to the axioms from Definition 4.28, let us recall that the correspondences there were partly obtained by composing. Here is an equivalent formulation of our axioms, which is more convenient, and that we will use in what follows, cutting some trivial redundancies:

\begin{theorem}
A quadruplet $(S,T,U,K)$ produces a noncommutative geometry when
$$\begin{matrix}
S&=&S_U\\
\\
S\cap\mathbb T_N^+&=&T&=&K\cap\mathbb T_N^+\\
\\
G^+(S)&=&<O_N,T>&=&U\\
\\
K^+(T)&=&U\cap K_N^+&=&K
\end{matrix}$$
with the usual convention that all this is up to the equivalence relation.
\end{theorem}

\begin{proof}
This follows indeed by examining the axioms in Definition 4.28, by cutting some trivial redundancies, and then by rescaling the whole table.
\end{proof}

We will use many times the above result, in what follows, so let us comment now, a bit informally, on the 7 axioms that we have, arranged in increasing order of complexity, based on the 4 computations that we have already:

\medskip

\begin{enumerate}
\item $T=S\cap\mathbb T_N^+$ is usually something quite trivial, and easy to check.

\medskip 

\item $T=K\cap\mathbb T_N^+$ is once again something quite trivial, and easy to check.

\medskip

\item $K=U\cap K_N^+$ is of the same nature, usually some trivial algebra.

\medskip

\item $U=G^+(S)$ is something more subtle, of algebraic geometric nature, which usually requires some tricks, in the spirit of Bhowmick-Goswami \cite{bg1}. These tricks  can actually get very complicated, and for many examples of quantum spheres $S$, the corresponding quantum isometry groups $G^+(S)$ are not known yet.

\medskip

\item $K=K^+(T)$ is something in the same spirit, but more complicated, with even the simplest possible non-trivial cases, namely the free real and complex ones, requiring subtle ingredients, such as a good knowledge of the $q=-1$ twisting.

\medskip

\item $S=S_U$ is something fairly heavy, requiring a good knowledge of the advanced representation theory and probability theory of compact quantum groups. Note that this is our only way here of getting to the sphere $S$.

\medskip

\item $U=<O_N,T >$ is something heavy too, requiring an excellent knowledge of the advanced representation theory of compact quantum groups. In fact, this is the key axiom, beating in complexity all the previous axioms, taken altogether.
\end{enumerate}

\medskip

Regarding now further work on these axioms, with new examples of geometries, and will classification results, we will discuss this later, in chapters 9-12 below. We will see there, among others, that under strong supplementary axioms, called ``purity'' and ``uniformity'', the 4 main geometries, from Theorem 4.29, are the only ones. 

\bigskip

In view of this, the question of developing the real and complex free geometries, which are the ``main'' non-classical geometries, appears. We will discuss this in chapters 5-8 below, with the construction of various ``free homogeneous spaces'', and we will come back to this later as well, in chapters 13-16 below, with more advanced results.

\bigskip

Finally, before getting into all this, some general comments on what we are doing, or rather on what our precise motivations are, are perhaps to be made at this point, now that we are already 100 pages into this book. Our motivations are, and no surprise here, unchanged since page 1, general quantum mechanics. However, this is related to some ongoing work, and in the lack of a precise idea, because this is how research usually goes, once you have a clear idea of what you're doing, you just do the computations and you're done, here are answers to some natural questions that you might have: 

\bigskip

(1) Shall we be worried by the fact that we are mixing $\mathbb R,\mathbb C$. Not at all, because as explained above, in the beginning of this chapter, there is some mathematical evidence supporting this, and more specifically results pointing towards the fact that ``free geometry is something hybrid between real and complex, and perhaps even scalarless''.

\bigskip

(2) What about arbitrary fields $F$ then, why not using them right from the beginning. Well, this is something rather for the future. It is not presently known what easiness over an arbitrary field $F$ means, although there has been some interesting work here by Bichon, on free quantum groups, and by Belinschi and others, on free probability.

\bigskip

(3) But is there anything from physics supporting all this. Sure yes, classical mechanics happens over $F=\mathbb R$, quantum mechanics happens over $F=\mathbb C$, and now go unify them. The unification should happen below the QCD level, and the conjecture would be that our $\mathbb R/\mathbb C$ or perhaps scalarless free geometry can be of help here.

\bigskip

(4) Sure yes, but is there anything concrete in this direction. Yes and no. There are countless things, for the most coming from the work of Connes, Jones, Voiculescu, pointing towards free geometry and its applications to quantum mechanics, and more specifically to a ``free mechanics'', below the QCD level. But nothing concrete, so far.

\bigskip

(5) Got it for physics, so now a math question, shall we be worried by the fact that our free spaces are all compact. This is one good question, and physically I would say that quarks as we know them are confined, with the strong force acting between them being short-range. And isn't all this suggesting the word ``compactness''.

\bigskip

(6) And a last question, mathematics or perhaps physics, what about smoothness, are we really sure about this. Well, think statistical mechanics, and call that mathematics or physics, as you wish. That teaches us that smoothness is a miracle. And we are not here for talking about miracles, but rather about basic things, aren't we.

\bigskip

And that is all, for the moment. We will of course regularly comment more on these topics, later in this book, once we will learn more things.

\section*{4e. Exercises} 

There are many possible exercises on the material above, which was quite varied, and also on the final axioms, that we have not explored yet. As a first exercise, we have:

\begin{exercise}
Try finding an easy geometry, or at least a candidate for an easy geometry, without full verification of the axioms, not among the $4$ main ones.
\end{exercise}

This question will be investigated later in this book, but thinking a bit at these questions, by yourself, will certainly not hurt, and will help understanding the material below. As a hint here, try finding something between $ab=ba$ and freeness.

\begin{exercise}
Establish the following isomorphism, as usual modulo equivalence
$$PO_N^+=PU_N^+$$
by using Tannakian duality and easiness, or by any other means.
\end{exercise}

This is something that we have talked about in the above, with comments on the multiple possible proofs of this. Pick one, and have it worked out, with full details.

\begin{exercise}
Establish the following isomorphism, as usual modulo equivalence
$$O_N^+=<O_N,O_{N-1}^+>$$
by using representation theory, or any other means.
\end{exercise}

This is something quite tricky, so in case you do not find the solution, you can look for it in the literature, and write down a brief account of what you found.

\begin{exercise}
Try axiomatizing the missing correspondences, namely
$$T\to S$$
$$S\leftrightarrow K$$
and in case you fail, explain at least what the difficulties are.
\end{exercise}

No comment here, as the author of the present book has failed so far in axiomatizing these correspondences, and does not want to talk about this.

\begin{exercise}
Axiomatize the abstract noncommutative geometries in our sense in terms of $S$, or of $T$, or of $U$, or of $K$, only.
\end{exercise}

This is something that we talked about in the above, that can be formally solved very quickly, simply by modifying the axioms, in the obvious way. The problem, however, is that of working out the simplifications that might appear in this way.

\part{Free manifolds}

\ \vskip50mm

\begin{center}
{\em In the town San Domingo

As we laughed and danced all night

To the thrub of flamingo guitars

Seemed a long long way from tomorrow's fight}
\end{center}

\chapter{Free integration}

\section*{5a. Weingarten formula}

We have seen so far that $\mathbb R^N,\mathbb C^N$ have no free analogues, in an analytic sense, but that the ``basic geometry'' of $\mathbb R^N,\mathbb C^N$, taken in a somewhat abstract sense, does have free analogues, that we can informally call ``basic geometry'' of $\mathbb R^N_+,\mathbb C^N_+$. Thus, we have 4 main geometries, classical/free, real/complex, forming a diagram as follows:
$$\xymatrix@R=50pt@C=50pt{
\mathbb R^N_+\ar[r]&\mathbb C^N_+\\
\mathbb R^N\ar[u]\ar[r]&\mathbb C^N\ar[u]
}$$

In this second part of the present book, we develop the geometry of $\mathbb R^N_+,\mathbb C^N_+$. To be more precise, each of these free geometries consists so far of 4 objects, namely a sphere $S$, a torus $T$, a unitary group $U$, and a reflection group $K$. We must on one hand study $S,T,U,K$, from a geometric perspective, and on the other hand construct other ``free manifolds'', as for instance suitable homogeneous spaces, and study them too.

\bigskip

Observe that all this is not exactly related to the axiomatization work from chapters 1-4. We will of course heavily use the various things that we learned there, but basically, what we want to do here is something new. We want to develop free geometry, real and complex, and our goals will be very explicit. As a basic question here, we have:

\begin{question}
What is a free manifold?
\end{question}

Unfortunately this is a difficult question, whose solution is not known yet. As an illustration, even in the quantum group case, after 30 long years of work on free quantum groups, it is still not known what a free quantum group exactly is. A quite reasonable definition seems to be $S_N^+\subset G\subset U_N^+$, but then comes the conjecture that such a quantum group must be easy, with this conjecture being something non-trivial.

\bigskip

In short, modesty. We have as starting point $S,T,U,K$, and this is not that bad, and we will slowly enlarge our menagery of free manifolds, not to the point of solving Question 5.1, but at least to the point of understanding what this question says.

\bigskip

For going ahead with more modesty, let us take $N=2$, in the real case. We know what the free circle $\bigcirc$ and the free square $\square$ are, and we also know what the symmetries of these free $\bigcirc$ and free $\square$ are, and the question is, shall we be awarded a PhD in noncommutative geometry for that. Ironically, probably yes, such basic things being not necessarily known by everyone. But leaving now aside academia and politics, the truth is that we are somewhere at the level of the ancient Greeks. Or even below, because the Greeks knew for instance what a conic is, along with many other things.

\bigskip

Which brings us into a second question, what kind of manifolds shall we look at, and what kind of geometry do we want to develop. A look at what we have, $S,T,U,K$, does not help much, because these are really very basic manifolds, having all geometric properties that you can ever dream of, and therefore belonging to all geometric theories that you can ever imagine. And so, we need some kind of plan here.

\bigskip

Looking at the story of classical geometry, you would say why not doing some algebraic geometry, reaching first to the level of the ancient Greeks, and then going up into more complicated things, towards analogues of what was doing the Italian school. But this is in fact completely unreasonable, because somewhere between ancient Greeks and more modern Italians we had Newton, who axiomatized classical mechanics by using geometry. And are we here for axiomatizing quantum mechanics by using free geometry, most likely leading to a Nobel Prize in physics, or shall we aim for something more modest.

\bigskip

Well, looks like we are completely lost. Again, we know what the free analogues of $\bigcirc$ and $\square$ are, and the question is, with this tremendous piece of knowledge, what's next. Fortunately there are already people who have thought about such things, in slightly different noncommutative geometry contexts, and we have, as a key piece of advice:

\begin{fact}[Connes principle]
Manifolds should be Riemannian.
\end{fact}

Here we are talking of course about noncommutative manifolds, because in what concerns the classical manifolds, this principle goes back to Riemann himself.

\bigskip

This being said, how can a free manifold be Riemannian, because we already know, as explained on several occasions in chapters 1-4, that such free manifolds are not smooth. And checking the mathematics and physics literature here, in look for new ideas, does not help much, because Riemannian geometry is always associated, in mathematics and physics, with all kinds of complicated differential geometry computations.

\bigskip

As a last-ditch attempt, let us ask instead an engineer. And the engineer actually answers something which is very interesting for us, namely:

\begin{fact}[Engineer's take]
Riemannian means that you can integrate over it.
\end{fact}

This is of course, technically speaking, not exactly correct. But hey, we are deep into the mud, and open to any piece of valuable advice. And valuable advice this is, because we know how to integrate over $S,T,U,K$, and so we should look for similar manifolds, having a sort of Haar measure, that we can compute via a Weingarten formula.

\bigskip

Be said in passing, if you're not familiar with engineers and engineering, what our engineer friend says in Fact 5.3 is in fact something quite subtle, with ``you can integrate'' rather meaning ``your computer can integrate''. Which is something which really fits with what we are doing, because the Weingarten formula can be implemented on a computer, and so is ready for ``quantum engineering'', whatever that might mean.

\bigskip

All this looks good, but as a last piece of philosophy now, aren't we going into some kind of extreme, with respect to smoothness, by following engineer's advice. Indeed, if we adopt this viewpoint, what shall we then think of the Riemannian manifolds which are smooth, but lack an efficient integration formula over them, such as a Weingarten formula? This is not en easy question, to put it this way, and in the lack of any academic willing to discuss such things, we will have to ask the cat. And cat says:

\begin{fact}[Cat's take]
Some manifolds are more Riemannian than other.
\end{fact}

Which sounds wise, there is now agreement between everyone involved so far. Be said in passing, what cat says agrees as well with Nash \cite{nas}, suggesting that the noncommutative Riemannian manifolds having coordinates are ``more Riemannian'' that those not having coordinates. And also with von Neumann \cite{von}, teaching us the noncommutative spaces having trace functionals $tr:L^\infty(X)\to\mathbb C$ are perhaps ``more Riemannian'' than those lacking this property. And also with Jones \cite{jo3}, and with Voiculescu \cite{vdn}, whose theories need trace functionals $tr:L^\infty(X)\to\mathbb C$, and with the ``quality'' of such a trace functional being directly related to the quality of your investigations. 

\bigskip

In short, problem solved. As a last thing, however, more on smoothness. Classical geometry teaches us that smoothness comes in several flavors, $C^1,C^2,\ldots,C^\infty$, with the mathematical reasons behind this being usually complicated, such as solutions of PDE, or singularities of algebraic manifolds, and so on, and with the physical reasons behind this being even more complicated, basically coming from statistical mechanics. So if there is one thing to be said, ``some manifolds are more smooth than other''. Obviously.

\bigskip

In what regards now our free manifolds, let us not forget that these have, as mentioned on several occasions in chapters 1-4, a Laplacian $\Delta$. However, no one knows how to construct this Laplacian in general, nor how to use it in relation to integration, nor how to use it in order to have some PDE running, on these manifolds. But one day, all this will be done, and our free manifolds will be entitled to be called ``half-smooth''. 

\bigskip

So, this will be our philosophy, for the next 100 pages to follow. As usual when regarding controversies, we can only recommend more reading on this, geometry at large, story of geometry, and Riemannian manifolds and related topics. Good references here are Shafarevich \cite{sha}, do Carmo \cite{doc} and Arnold \cite{ar2} for geometry, then von Neumann \cite{von} or Blackadar \cite{bla} for operator algebras, and then Connes \cite{co1} and Connes-Marcolli \cite{cma} for a mix of operator algebras and geometry, complemented perhaps with Gracia-Bond\'ia-V\'arilly-Figueroa \cite{gvf} and Landi \cite{lan}. And don't forget about Nash \cite{nas}.

\bigskip

Back to work now, our first task will be that of explaining how to integrate over $S,T,U,K$. In order to integrate over $U,K$, we can use the Weingarten formula \cite{csn}, \cite{wei}, whose general quantum group formulation, from \cite{bsp}, is as follows:

\index{Weingarten formula}
\index{Gram matrix}
\index{Weingarten matrix}

\begin{theorem}
Assuming that a closed subgroup $G\subset U_N^+$ is easy, coming from a category of partitions $D\subset P$, we have the Weingarten formula
$$\int_Gu_{i_1j_1}^{e_1}\ldots u_{i_kj_k}^{e_k}=\sum_{\pi,\sigma\in D(k)}\delta_\pi(i)\delta_\sigma(j)W_{kN}(\pi,\sigma)$$
where $\delta$ are Kronecker type symbols, and where the Weingarten matrix
$$W_{kN}=G_{kN}^{-1}$$
is the inverse of the Gram matrix $G_{kN}(\pi,\sigma)=N^{|\pi\vee\sigma|}$. This formula applies to all classical and free unitary and reflection groups $U,K$, which are all easy.
\end{theorem}

\begin{proof}
We know the Weingarten formula from chapter 3, the idea being that the integrals in the statement form the projection on the following space:
$$Fix(u^{\otimes k})=span\left(\xi_\pi\Big|\pi\in D(k)\right)$$ 

As for the easiness property of our various classical and free unitary and reflection groups $U,K$, this is something that we know too from chapter 3.
\end{proof}

Regarding now the integration over the tori $T$, this is something very simple, because we can use here the following fact, coming from the definition of group algebras:

\begin{theorem}
Given a finitely generated discrete group $\Gamma=<g_1,\ldots,g_N>$, the integrals over the corresponding torus $T=\widehat{\Gamma}$ are given by
$$\int_Tg_{i_1}^{e_1}\ldots g_{i_k}^{e_k}=\delta_{g_{i_1}^{e_1}\ldots g_{i_k}^{e_k},1}$$
for any indices $i_r\in\{1,\ldots,N\}$ and any exponents $e_r\in\{\emptyset,*\}$, with the Kronecker symbol on the right being a usual one, computed inside the group $\Gamma$.
\end{theorem}

\begin{proof}
This is something clear, coming from the fact that the Haar integration over the torus $T=\widehat{\Gamma}$ is given by the following formula:
$$\int_Tg=\delta_{g1}$$

Indeed, this formula defines a functional on the algebra $C(T)=C^*(\Gamma)$, which is obviously left and right invariant, and so is the Haar functional.
\end{proof}

Finally, regarding $S$, here the integrals appear as particular cases of the integrals over $U$, as explained in chapter 3, and we have a Weingarten formula, as follows:

\index{integration over spheres}
 
\begin{theorem}
The integration over a sphere $S$, which is such that $U=G^+(S)$ is easy, coming from a category of pairings $D$, is given by the Weingarten formula
$$\int_Sx_{i_1}^{e_1}\ldots x_{i_k}^{e_k}=\sum_\pi\sum_{\sigma\leq\ker i}W_{kN}(\pi,\sigma)$$
with $\pi,\sigma\in D(k)$, where $W_{kN}=G_{kN}^{-1}$ is the inverse of $G_{kN}(\pi,\sigma)=N^{|\pi\vee\sigma|}$. This formula applies to all classical and free spheres $S$, whose $U=G^+(S)$ are easy.
\end{theorem}

\begin{proof}
This is something that we know too from chapter 3, coming from the definition of the integration functional over $S$, as being the following composition:
$$\int_S:C(S)\to C(U)\to\mathbb C$$

Indeed, with this description of the integration functional in mind, we can compute this functional via the Weingarten formula for $U$, from Theorem 5.5, as follows:
\begin{eqnarray*}
\int_Sx_{i_1}^{e_1}\ldots x_{i_k}^{e_k}
&=&\int_Uu_{1i_1}^{e_1}\ldots u_{1i_k}^{e_k}\\
&=&\sum_{\pi,\sigma\in D(k)}\delta_\pi(1)\delta_\sigma(i)W_{kN}(\pi,\sigma)\\
&=&\sum_\pi\sum_{\sigma\leq\ker i}W_{kN}(\pi,\sigma)
\end{eqnarray*}

Thus, we are led to the formula in the statement.
\end{proof}

Summarizing, we know how to integrate over $S,T,U,K$, and with the remark, in relation with Fact 5.3, that our integration methods can be implemented on a computer, as engineers love it. Indeed, after implementing $D$, and then $G_{kN}(\pi,\sigma)=N^{|\pi\vee\sigma|}$, the big problem, which is that of inverting, $W_{kN}=G_{kN}^{-1}$, can be solved by the computer, and then you obtain all the integrals that you want just by summing. For some numerics here, you can check for instance the various papers citing Collins-\'Sniady \cite{csn}.

\section*{5b. Free probability}

We are not over with integration, because we have now to apply our various results above, to some suitable variables, and see what we get. The range of applications here is potentially infinite, and in the lack of a good high energy physics problem to be solved, and let us put that on our to-do list, we will do some pure mathematics.

\bigskip

The point indeed is that Voiculescu came in the 80s with a beautiful theory of free probability, explained in his book with Dykema and Nica \cite{vdn}, and we would like to know if our liberation considerations fit with this. More specifically, the correspondence between classical and free probability was axiomatized by Bercovici-Pata in \cite{bpa}, and we would like to know if our constructions $X\to X^+$ fit with this correspondence. 

\bigskip

And there is a long way to go here. First we must explain free probability, following Voiculescu-Dykema-Nica \cite{vdn}, so let start with the following standard definition:

\index{random variable}
\index{moments}
\index{colored integers}
\index{law}
\index{distribution}

\begin{definition}
Let $A$ be a $C^*$-algebra, given with a positive trace $tr$.
\begin{enumerate}
\item The elements $a\in A$ are called random variables.

\item The moments of such a variable are the numbers $M_k(a)=tr(a^k)$.

\item The law of such a variable is the functional $\mu_a:P\to tr(P(a))$.
\end{enumerate}
\end{definition}

Here $k=\circ\bullet\bullet\circ\ldots$ is as usual a colored integer, and the powers $a^k$ are defined by the usual formulae, namely $a^\emptyset=1,a^\circ=a,a^\bullet=a^*$ and multiplicativity. As for the polynomial $P$, this is by definition a noncommuting $*$-polynomial in one variable:
$$P\in\mathbb C<X,X^*>$$

Observe that the law is uniquely determined by the moments, because:
$$P(X)=\sum_k\lambda_kX^k\implies\mu_a(P)=\sum_k\lambda_kM_k(a)$$

In the self-adjoint case, $a=a^*$ the law is a usual probability measure, supported by the spectrum of $a$. This follows indeed from the Gelfand theorem, and the Riesz theorem. More generally, the same happens in the normal case, $aa^*=a^*a$, with the spectrum being now complex. However, in the non-normal case, $aa^*\neq a^*a$, such a probability measure describing the law $\mu_a$ does not exist, due to the following computation:
\begin{eqnarray*}
aa^*-a^*a\neq0
&\implies&(aa^*-a^*a)^2>0\\
&\implies&aa^*aa^*-aa^*a^*a-a^*aaa^*+a^*aa^*a>0\\
&\implies&tr(aa^*aa^*-aa^*a^*a-a^*aaa^*+a^*aa^*a)>0\\
&\implies&tr(aa^*aa^*+a^*aa^*a)>tr(aa^*a^*a+a^*aaa^*)\\
&\implies&tr(aa^*aa^*)>tr(aaa^*a^*)
\end{eqnarray*}

Indeed, assuming that $a$ has a probability measure as law, the above quantities would both appear by integrating $|z|^2$ with respect to this measure, which is contradictory.

\bigskip

Talking now probability, in a general sense, if there is one thing to be known here, this is the Central Limit Theorem (CLT). So, let us start with this:

\index{CLT}
\index{normal variable}
\index{Gaussian variable}
\index{Fourier transform}

\begin{theorem}[CLT]
Given real random variables $x_1,x_2,x_3,\ldots,$ which are i.i.d., centered, and with variance $t>0$, we have, with $n\to\infty$, in moments,
$$\frac{1}{\sqrt{n}}\sum_{i=1}^nx_i\sim g_t$$
where $g_t$ is the Gaussian law of parameter $t$, having as density $g_t=\frac{1}{\sqrt{2\pi t}}\,e^{-x^2/2t}dx$.
\end{theorem}

\begin{proof}
This is something standard, the proof being in three steps, as follows:

\medskip

(1) Linearization of the convolution. It well-known that the log of the Fourier transform $F_x(\xi)=E(e^{i \xi x})$ does the job, in the sense that if $x,y$ are independent, then:
$$F_{x+y}=F_xF_y$$

(2) Study of the limit. We have the following formula for a general Fourier transform $F_x(\xi)=E(e^{i \xi x})$, in terms of moments:
$$F_x(\xi)=\sum_{k=0}^\infty\frac{i^kM_k(x)}{k!}\,\xi^k$$

It follows that the Fourier transform of the variable in the statement is:
\begin{eqnarray*}
F(\xi)
&=&\left[F_x\left(\frac{\xi}{\sqrt{n}}\right)\right]^n\\
&=&\left[1-\frac{t\xi^2}{2n}+O(n^{-2})\right]^n\\
&\simeq&\left[1-\frac{t\xi^2}{2n}\right]^n\\
&\simeq&e^{-t\xi^2/2}
\end{eqnarray*}

(3) Gaussian laws. The Fourier transform of the Gaussian law is given by:
\begin{eqnarray*}
F_{g_t}(x)
&=&\frac{1}{\sqrt{2\pi t}}\int_\mathbb Re^{-y^2/2t+ixy}dy\\
&=&\frac{1}{\sqrt{2\pi t}}\int_\mathbb Re^{-(y/\sqrt{2t}-\sqrt{t/2}ix)^2-tx^2/2}dy\\
&=&\frac{1}{\sqrt{2\pi t}}\int_\mathbb Re^{-z^2-tx^2/2}\sqrt{2t}dz\\
&=&e^{-tx^2/2}
\end{eqnarray*}

Thus the variables on the left and on the right in the statement have the same Fourier transform, and so these variables follow the same law, as claimed.
\end{proof}

Following Voiculescu \cite{vdn}, in order to extend the CLT to the free setting, our starting point will be the following definition:

\index{freeness}

\begin{definition}
Given a pair $(A,tr)$, two subalgebras $B,C\subset A$ are called free when the following condition is satisfied, for any $b_i\in B$ and $c_i\in C$:
$$tr(b_i)=tr(c_i)=0\implies tr(b_1c_1b_2c_2\ldots)=0$$
Also, two noncommutative random variables $b,c\in A$ are called free when the $C^*$-algebras $B=<b>$, $C=<c>$ that they generate inside $A$ are free, in this sense.
\end{definition}

As a first observation, there is a similarity here with the classical notion of independence. Indeed, modulo some standard identifications, two subalgebras $B,C\subset L^\infty(X)$ are independent when the following condition is satisfied, for any $b\in B$ and $c\in C$: 
$$tr(bc)=tr(b)tr(c)$$

But this is equivalent to the following condition, which is similar to freeness:
$$tr(b)=tr(c)=0\implies tr(bc)=0$$

In short, freeness appears by definition as a kind of ``free analogue'' of independence. As a first result now regarding this notion, clarifying the basics, we have:

\begin{proposition}
Assuming that $B,C\subset A$ are free, the restriction of $tr$ to $<B,C>$ can be computed in terms of the restrictions of $tr$ to $B,C$. To be more precise,
$$tr(b_1c_1b_2c_2\ldots)=P\Big(\{tr(b_{i_1}b_{i_2}\ldots)\}_i,\{tr(c_{j_1}c_{j_2}\ldots)\}_j\Big)$$
where $P$ is certain polynomial in several variables, depending on the length of the word $b_1c_1b_2c_2\ldots$, and having as variables the traces of products of type
$$b_{i_1}b_{i_2}\ldots\quad,\quad 
c_{j_1}c_{j_2}\ldots$$
with the indices being chosen increasing, $i_1<i_2<\ldots$ and $j_1<j_2<\ldots$
\end{proposition}

\begin{proof}
We can start indeed our computation as follows:
\begin{eqnarray*}
tr(b_1c_1b_2c_2\ldots)
&=&tr\big[(b_1'+tr(b_1))(c_1'+tr(c_1))(b_2'+tr(b_2))(c_2'+tr(c_2))\ldots\ldots\big]\\
&=&tr(b_1'c_1'b_2'c_2'\ldots)+{\rm other\ terms}\\
&=&{\rm other\ terms}
\end{eqnarray*}

Observe that we have used here the freeness condition, in the following form:
$$tr(b_i')=tr(c_i')=0\implies tr(b_1'c_1'b_2'c_2'\ldots)=0$$

Thus, we are led into some sort of recurrence, as desired. For more on all this, including examples, we refer to the book of Voiculescu-Dykema-Nica \cite{vdn}.
\end{proof}

As a second result regarding the notion of freeness, which provides us with a useful class of examples, which can be used for various modeling purposes, we have:

\begin{proposition}
Given two algebras $(B,tr)$ and $(C,tr)$, the following hold:
\begin{enumerate}
\item $B,C$ are independent inside their tensor product $B\otimes C$, endowed with its canonical tensor product trace, given on basic tensors by $tr(b\otimes c)=tr(b)tr(c)$.

\item $B,C$ are free inside their free product $B*C$, endowed with its canonical free product trace, given by the formulae in Proposition 5.11.
\end{enumerate}
\end{proposition}

\begin{proof}
Both the assertions are clear from definitions, as follows:

\medskip

(1) This is clear, because we have by construction of the trace:
\begin{eqnarray*}
tr(bc)
&=&tr[(b\otimes1)(1\otimes c)]\\
&=&tr(b\otimes c)\\
&=&tr(b)tr(c)
\end{eqnarray*}

(2) This is clear again, the only point being that of showing that the notion of freeness, or the recurrence formulae in Proposition 5.11, can be used in order to construct a canonical free product trace, on the free product of the two algebras involved:
$$tr:B*C\to\mathbb C$$

But this can be done for instance by using a GNS construction. Indeed, by taking the free product of the GNS constructions for $(B,tr)$ and $(C,tr)$, we obtain a representation as follows, with the $*$ on the right being a free product of pointed Hilbert spaces:
$$B*C\to B(L^2(B)*L^2(C))$$

Now by composing with the linear form $T\to<T\xi,\xi>$, where $\xi=1_B=1_C$ is the common distinguished vector of $L^2(B)$ and $L^2(C)$, we obtain a linear form, as follows:
$$tr:B*C\to\mathbb C$$

It is routine then to check that $tr$ is indeed a trace, and this is the ``canonical free product trace'' from the statement. Then, an elementary computation shows that $B,C$ are indeed free inside $B*C$, with respect to this trace, as desired. 
\end{proof}

Finally, again following \cite{vdn}, we have the following more explicit modelling result:

\begin{theorem}
We have the following results, valid for group algebras:
\begin{enumerate}
\item $C^*(\Gamma),C^*(\Lambda)$ are independent inside $C^*(\Gamma\times\Lambda)$.

\item $C^*(\Gamma),C^*(\Lambda)$ are free inside $C^*(\Gamma*\Lambda)$.
\end{enumerate}
\end{theorem}

\begin{proof}
We can use here the general results in Proposition 5.12, along with the following two isomorphisms, which are both standard:
$$C^*(\Gamma\times\Lambda)=C^*(\Lambda)\otimes C^*(\Gamma)$$
$$C^*(\Gamma*\Lambda)=C^*(\Lambda)*C^*(\Gamma)$$

Alternatively, we can prove this directly, starting from definitions, by using the fact that each group algebra is spanned by the corresponding group elements.
\end{proof}

There are many interesting things that can be said about the analogy between independence and freeness. We have in particular the following remarkable result, due to Voiculescu \cite{vdn}, which is the key to everything in free probability:

\index{Cauchy transform}
\index{R-transform}

\begin{theorem}
Given a real probability measure $\mu$, consider its Cauchy transform
$$G_\mu(\xi)=\int_\mathbb R\frac{d\mu(t)}{\xi-t}$$
and define its $R$-transform as being the solution of the following equation:
$$G_\mu\left(R_\mu(\xi)+\frac{1}{\xi}\right)=\xi$$
The operation $\mu\to R_\mu$ linearizes then the free convolution operation.
\end{theorem}

\begin{proof}
In order to prove this, we need a good model for the free convolution. The best here is to use the semigroup algebra of the free semigroup on two generators:
$$A=C^*(\mathbb N*\mathbb N)$$

Indeed, we have some freeness in the semigroup setting, a bit in the same way as for the group algebras $C^*(\Gamma*\Lambda)$, from Theorem 5.13 (2), and in addition to this fact, and to what happens in the group algebra case, the following two key things happen:

\medskip

(1) The variables of type $S^*+f(S)$, with $S\in C^*(\mathbb N)$ being the shift, and with $f\in\mathbb C[X]$ being a polynomial, model in moments all the distributions $\mu:\mathbb C[X]\to\mathbb C$. This is indeed something elementary, which can be checked via a direct algebraic computation.

\medskip

(2) Given $f,g\in\mathbb C[X]$, the variables $S^*+f(S)$ and $T^*+g(T)$, where $S,T\in C^*(\mathbb N*\mathbb N)$ are the shifts corresponding to the generators of $\mathbb N*\mathbb N$, are free, and their sum has the same law as $S^*+(f+g)(S)$. This follows indeed by using a $45^\circ$ argument.

\medskip

With these results in hand, we can see  that the operation $\mu\to f$ linearizes the free convolution. We are therefore left with a computation inside $C^*(\mathbb N)$, whose conclusion is that $R_\mu=f$ can be recaptured from $\mu$ via the Cauchy transform $G_\mu$, as stated.
\end{proof}

We can now state and prove a free analogue of the CLT, from \cite{vdn}, as follows:

\index{FCLT}
\index{free CLT}

\begin{theorem}[FCLT]
Given self-adjoint variables $x_1,x_2,x_3,\ldots,$ which are f.i.d., centered, with variance $t>0$, we have, with $n\to\infty$, in moments,
$$\frac{1}{\sqrt{n}}\sum_{i=1}^nx_i\sim\gamma_t$$\
where $\gamma_t$ is the Wigner semicircle law of parameter $t$, given by $\gamma_t=\frac{1}{2\pi t}\sqrt{4t^2-x^2}dx$.
\end{theorem}

\begin{proof}
At $t=1$, the $R$-transform of the variable in the statement can be computed by using the linearization property with respect to the free convolution, and is:
$$R(\xi)
=nR_x\left(\frac{\xi}{\sqrt{n}}\right)
\simeq\xi$$

On the other hand, some elementary computations show that the Cauchy transform of the Wigner law $\gamma_1$ from the statement satisfies the following equation:
$$G_{\gamma_1}\left(\xi+\frac{1}{\xi}\right)=\xi$$

Thus we have $R_{\gamma_1}(\xi)=\xi$, which by the way follows as well from $S^*+S\sim\gamma_1$, and this gives the result. The passage to the general case, $t>0$, is routine.
\end{proof}

In the complex case now, we recall that the complex Gaussian law of parameter $t>0$ is defined as follows, with $a,b$ being independent, each following the law $g_t$:
$$G_t=law\left(\frac{1}{\sqrt{2}}(a+ib)\right)$$

With this convention, we have the following result:

\index{CCLT}

\begin{theorem}[CCLT]
Given complex random variables $x_1,x_2,x_3,\ldots,$ which are i.i.d., centered, and with variance $t>0$, we have, with $n\to\infty$, in moments,
$$\frac{1}{\sqrt{n}}\sum_{i=1}^nx_i\sim G_t$$
where $G_t$ is the complex Gaussian law of parameter $t$.
\end{theorem}

\begin{proof}
This follows indeed from the real CLT, established above, without new computations needed, just by taking real and imaginary parts.
\end{proof}

In the free case, the Voiculescu circular law of parameter $t>0$ is defined as follows, with $\alpha,\beta$ being independent, each following the law $\gamma_t$:
$$\Gamma_t=law\left(\frac{1}{\sqrt{2}}(\alpha+i\beta)\right)$$

With this convention, we have the following result, again from Voiculescu \cite{vdn}:

\index{FCCLT}
\index{Voiculescu circular law}
\index{circular variable}

\begin{theorem}[FCCLT]
Given noncommutative random variables $x_1,x_2,x_3,\ldots,$ which are f.i.d., centered, and with variance $t>0$, we have, with $n\to\infty$, in moments,
$$\frac{1}{\sqrt{n}}\sum_{i=1}^nx_i\sim\Gamma_t$$
where $\Gamma_t$ is the Voiculescu circular law of parameter $t$.
\end{theorem}

\begin{proof}
This follows indeed from the FCLT, by taking real and imaginary parts.
\end{proof}

With these ingredients in hand, let us go back now to our quantum groups. According to the Peter-Weyl theory from chapter 2, if there are some variables that we should look at, these are the characters. And here, for the unitary quantum groups, we have:

\index{main character}
\index{asymptotic character}

\begin{theorem}
With $N\to\infty$, the main characters 
$$\chi=\sum_{i=1}^Nu_{ii}$$
for the basic unitary quantum groups are as follows:
\begin{enumerate}
\item $O_N$: real Gaussian, following $g_1$.

\item $O_N^+$: semicircular, following $\gamma_1$.

\item $U_N$: complex Gaussian, following $G_1$.

\item $U_N^+$: circular, following $\Gamma_1$.
\end{enumerate}
\end{theorem}

\begin{proof}
We use the moment method, and combinatorics. For a closed subgroup $G_N\subset U_N^+$, we have, according to the Peter-Weyl type results of Woronowicz in \cite{wo1}:
$$\int_{G_N}\chi^k=\dim(Fix(u^{\otimes k}))$$

In the easy case now, where $G=(G_N)$ comes from a certain category of partitions $D$, the fixed point space on the right is spanned by the vectors $T_\pi$ with $\pi\in D(k)$. Now since by Lindst\"om \cite{lin} these vectors are linearly independent with $N\to\infty$, we have:
$$\lim_{N\to\infty}\int_{G_N}\chi^k=|D(k)|$$

Thus, we are led into some combinatorics, and the continuation is as follows:

\medskip

(1) For $O_N$ we have $D=P_2$, so we obtain as even asymptotic moments the numbers $|P_2(2k)|=k!!$, which are well-known to be the moments of the Gaussian law.

\medskip

(2) For $O_N^+$ we have $D=NC_2$, so we obtain as even asymptotic moments the Catalan numbers $|NC_2(2k)|=C_k$, which are the moments of the Wigner semicircle law.

\medskip

(3) For $U_N$ we have $D=\mathcal P_2$, and we can conclude as in the real case, involving $O_N$, by using this time moments with respect to colored integers, as in Definition 5.8.

\medskip

(4) For $U_N^+$ we have $D=\mathcal{NC}_2$, and again we can conclude as in the real case, involving $O_N^+$, by using moments with respect to colored integers, as in Definition 5.8.
\end{proof}

The above result is of course just the tip of the iceberg, and there are countless things that can be done, as a continuation of this. In what follows we will orient the discussion towards something rather theoretical, namely the Bercovici-Pata bijection \cite{bpa}.

\section*{5c. Truncated characters}

We have seen so far that for $O_N,O_N^+,U_N,U_N^+$, the asymptotic laws of the main characters are the laws $g_1,\gamma_1,G_1,\Gamma_1$ coming from the various classical and free CLT. This is certainly nice, but there is still one conceptual problem, coming from:

\begin{proposition}
The above convergences $law(\chi_u)\to g_1,\gamma_1,G_1,\Gamma_1$ are as follows: 
\begin{enumerate}
\item They are non-stationary in the classical case.

\item They are stationary in the free case, starting from $N=2$.
\end{enumerate}
\end{proposition}

\begin{proof}
This is something quite subtle, which can be proved as follows:

\medskip

(1) Here we can use an amenability argument, based on the Kesten criterion. Indeed, $O_N,U_N$ being coamenable, the upper bound of the support of the law of $Re(\chi_u)$ is precisely $N$, and we obtain from this that the law of $\chi_u$ itself depends on $N\in\mathbb N$.

\medskip

(2) Here the result follows from the well-known fact that the linear maps $T_\pi$ associated to the noncrossing pairings are linearly independent, at any $N\geq2$, which fact, which is non-trivial, follows itself either from the general theory developed by Jones in \cite{jo1}, in relation with the Temperley-Lieb algebra, or from Di Francesco \cite{dif}.
\end{proof}

Fortunately, the solution to the convergence question raised by Proposition 5.19 is quite simple. The idea will be that of improving our $g_1,\gamma_1,G_1,\Gamma_1$ results with certain $g_t,\gamma_t,G_t,\Gamma_t$ results, which will require $N\to\infty$ in both the classical and free cases, in order to hold at any $t$. Following \cite{bbc}, the definition that we will need is as follows:

\index{truncated character}

\begin{definition}
Given a Woronowicz algebra $(A,u)$, the variable
$$\chi_t=\sum_{i=1}^{[tN]}u_{ii}$$
is called truncation of the main character, with parameter $t\in(0,1]$.
\end{definition}

Our purpose in what follows will be that of proving that for $O_N,O_N^+,U_N,U_N^+$, the asymptotic laws of the truncated characters $\chi_t$ with $t\in(0,1]$ are the laws $g_t,\gamma_t,G_t,\Gamma_t$. This is something quite technical, but natural, motivated by the findings in Proposition 5.19, and also by a number of more advanced considerations, to become clear later on. So, let us do this. In order to study the truncated characters, we can use:

\begin{theorem}
The moments of the truncated characters are given by
$$\int_G(u_{11}+\ldots +u_{ss})^k=Tr(W_{kN}G_{ks})$$
and with $N\to\infty$ this quantity equals $(s/N)^k|D(k)|$.
\end{theorem}

\begin{proof}
The first assertion follows from the following computation:
\begin{eqnarray*}
\int_G(u_{11}+\ldots +u_{ss})^k
&=&\sum_{i_1=1}^{s}\ldots\sum_{i_k=1}^s\int u_{i_1i_1}\ldots u_{i_ki_k}\\
&=&\sum_{\pi,\sigma\in D(k)}W_{kN}(\pi,\sigma)\sum_{i_1=1}^{s}\ldots\sum_{i_k=1}^s\delta_\pi(i)\delta_\sigma(i)\\
&=&\sum_{\pi,\sigma\in D(k)}W_{kN}(\pi,\sigma)G_{ks}(\sigma,\pi)\\
&=&Tr(W_{kN}G_{ks})
\end{eqnarray*}

We have $G_{kN}(\pi,\sigma)=N^k$ for $\pi=\sigma$, and $G_{kN}(\pi,\sigma)\leq N^{k-1}$ for $\pi\neq\sigma$. Thus with $N\to\infty$ we have $G_{kN}\sim N^k1$, which gives:
\begin{eqnarray*}
\int_G(u_{11}+\ldots +u_{ss})^k
&=&Tr(G_{kN}^{-1}G_{ks})\\
&\sim&Tr((N^k1)^{-1} G_{ks})\\
&=&N^{-k}Tr(G_{ks})\\
&=&N^{-k}s^k|D(k)|
\end{eqnarray*}

Thus, we have obtained the formula in the statement.
\end{proof}

In order to process the above moment formula, we will need some more probability theory, both classical and free. Given a random variable $a$, we write:
$$\log F_a(-iz)=\sum_{n=1}^\infty k_n(a)\cdot\frac{z^n}{n!}\quad,\quad 
zR_a(z)=\sum_{n=1}^\infty\kappa_n(a)\cdot z^n$$

We call the above coefficients $k_n(a),\kappa_n(a)$ the cumulants, respectively free cumulants of $a$. With this notion in hand, we can define then more general quantities $k_\pi(a),\kappa_\pi(a)$, depending on arbitrary partitions $\pi\in P(k)$, which coincide with the above ones for the 1-block partitions, and then by multiplicativity over the blocks, and we have:

\index{cumulant}
\index{free cumulant}
\index{moment-cumulant formula}

\begin{theorem}
We have the classical and free moment-cumulant formulae
$$M_k(a)=\sum_{\pi\in P(k)}k_\pi(a)\quad,\quad 
M_k(a)=\sum_{\pi\in NC(k)}\kappa_\pi(a)$$
where $k_\pi(a),\kappa_\pi(a)$ are the generalized cumulants and free cumulants of $a$.
\end{theorem}

\begin{proof}
This is something very standard, due to Rota in the classical case, and to Speicher in the free case, obtained either by using the formulae of $F_a,R_a$, or by doing some direct combinatorics, based on the M\"obius inversion formula. See \cite{vdn}.
\end{proof}

We can now improve our results about characters, as follows:

\index{truncated character}
\index{asymptotic character}

\begin{theorem}
With $N\to\infty$, the laws of truncated characters are as follows:
\begin{enumerate}
\item For $O_N$ we obtain the Gaussian law $g_t$.

\item For $O_N^+$ we obtain the Wigner semicircle law $\gamma_t$.

\item For $U_N$ we obtain the complex Gaussian law $G_t$.

\item For $U_N^+$ we obtain the Voiculescu circular law $\Gamma_t$.
\end{enumerate}
\end{theorem}

\begin{proof}
With $s=[tN]$ and $N\to\infty$, the formula in Theorem 5.21 gives:
$$\lim_{N\to\infty}\int_{G_N}\chi_t^k=\sum_{\pi\in D(k)}t^{|\pi|}$$

By using now the formulae in Theorem 5.22, this gives the results.
\end{proof}

As an interesting consequence, related to \cite{bpa}, let us formulate as well:

\index{Bercovici-Pata bijection}

\begin{theorem}
The asymptotic laws of truncated characters for the operations 
$$O_N\to O_N^+\quad,\quad 
U_N\to U_N^+$$
are in Bercovici-Pata bijection, in the sense that the classical cumulants in the classical case equal the free cumulants in the free case.
\end{theorem}

\begin{proof}
This follows indeed from Theorem 5.23, and from the standard combinatorial interpretation of the Bercovici-Pata bijection \cite{bpa}, in terms of cumulants.
\end{proof}

Let us discuss now the integration over the spheres. Following \cite{bgo}, we have:

\begin{theorem}
With $N\to\infty$, the rescaled coordinates of the various spheres 
$$\sqrt{N}x_i\in C(S^{N-1}_{\times})$$
are as follows, with respect to the uniform integration:
\begin{enumerate}
\item $S^{N-1}_\mathbb R$: real Gaussian.

\item $S^{N-1}_{\mathbb R,+}$: semicircular.

\item $S^{N-1}_\mathbb C$: complex Gaussian.

\item $S^{N-1}_{\mathbb C,+}$: circular.
\end{enumerate}
\end{theorem}

\begin{proof}
This follows from Theorem 5.23, but we can use as well the Weingarten formula for the spheres, from Theorem 5.7. Indeed, we have the following estimate:
$$\int_{S^{N-1}_\times}x_{i_1}\ldots x_{i_k}\,dx\simeq N^{-k/2}\sum_{\sigma\in P_2^\times(k)}\delta_\sigma(i)$$

With this formula in hand, we can compute the asymptotic moments of each coordinate $x_i$. Indeed, by setting $i_1=\ldots=i_k=i$, all Kronecker symbols are 1, and we obtain:
$$\int_{S^{N-1}_\times}x_i^k\,dx\simeq N^{-k/2}|P_2^\times(k)|$$

But this gives the results, via the same combinatorics as before. See \cite{bgo}.
\end{proof}

\section*{5d. Poisson laws}

In order to discuss now the quantum reflection groups, we will need some more theory, namely Poisson limit theorems. In the classical case, we have the following result:

\index{Poisson limit theorem}
\index{PLT}
\index{Poisson law}

\begin{theorem}[PLT]
We have the following convergence, in moments,
$$\left(\left(1-\frac{t}{n}\right)\delta_0+\frac{t}{n}\,\delta_1\right)^{*n}\to p_t$$
the limiting measure being
$$p_t=\frac{1}{e^t}\sum_{k=0}^\infty\frac{t^k\delta_k}{k!}$$
which is the Poisson law of parameter $t>0$.
\end{theorem}

\begin{proof}
We recall that the Fourier transform is given by $F_f(x)=E(e^{ixf})$. We therefore obtain the following formula:
\begin{eqnarray*}
F_{p_t}(x)
&=&e^{-t}\sum_k\frac{t^k}{k!}F_{\delta_k}(x)\\
&=&e^{-t}\sum_k\frac{t^k}{k!}\,e^{ikx}\\
&=&e^{-t}\sum_k\frac{(e^{ix}t)^k}{k!}\\
&=&\exp(-t)\exp(e^{ix}t)\\
&=&\exp\left((e^{ix}-1)t\right)
\end{eqnarray*}

Let us denote by $\mu_n$ the measure under the convolution sign, namely:
$$\mu_n=\left(1-\frac{t}{n}\right)\delta_0+\frac{t}{n}\delta_1$$

We have the following computation: 
\begin{eqnarray*}
F_{\delta_r}(x)=e^{irx}
&\implies&F_{\mu_n}(x)=\left(1-\frac{t}{n}\right)+\frac{t}{n}e^{ix}\\
&\implies&F_{\mu_n^{*n}}(x)=\left(\left(1-\frac{t}{n}\right)+\frac{t}{n}e^{ix}\right)^n\\
&\implies&F_{\mu_n^{*n}}(x)=\left(1+\frac{(e^{ix}-1)t}{n}\right)^n\\
&\implies&F(x)=\exp\left((e^{ix}-1)t\right)
\end{eqnarray*}

Thus, we obtain the Fourier transform of $p_t$, as desired.
\end{proof}

In the free case, the result is as follows, with $\boxplus$ being the free convolution operation:

\index{free Poisson limit theorem}
\index{FPLT}
\index{free Poisson law}
\index{Marchenko-Pastur law}

\begin{theorem}[FPLT]
We have the following convergence, in moments,
$$\left(\left(1-\frac{t}{n}\right)\delta_0+\frac{t}{n}\,\delta_1\right)^{\boxplus n}\to\pi_t$$
the limiting measure being the Marchenko-Pastur law of parameter $t>0$, 
$$\pi_t=\max(1-t,0)\delta_0+\frac{\sqrt{4t-(x-1-t)^2}}{2\pi x}\,dx$$
also called free Poisson law of parameter $t>0$.
\end{theorem}

\begin{proof}
Consider the measure in the statement, under the convolution sign:
$$\mu=\left(1-\frac{t}{n}\right)\delta_0+\frac{t}{n}\delta_1$$

The Cauchy transform of this measure is elementary to compute, given by:
$$G_{\mu}(\xi)=\left(1-\frac{t}{n}\right)\frac{1}{\xi}+\frac{t}{n}\cdot\frac{1}{\xi-1}$$

In order to prove the result, we want to compute the following $R$-transform:
$$R
=R_{\mu^{\boxplus n}}(y)
=nR_\mu(y)$$

But the equation for this function $R$ is as follows:
$$\left(1-\frac{t}{n}\right)\frac{1}{y^{-1}+R/n}+\frac{t}{n}\cdot\frac{1}{y^{-1}+R/n-1}=y$$

By multiplying by $n/y$, this equation can be written as:
$$\frac{t+yR}{1+yR/n}=\frac{t}{1+yR/n-y}$$

With $n\to\infty$ we obtain $t+yR=t/(1-y)$, so we conclude that we have:
$$R=\frac{t}{1-y}$$

On the other hand the $R$-transform of $\pi_t$ is standard to compute, based on the degree 2 equation satisfied by the Cauchy transform, and we get $t/(1-y)$ as well, as desired.
\end{proof}

In order to get beyond this, let us introduce the following notions:

\index{compound Poisson law}
\index{free compound Poisson law}

\begin{definition}
Associated to any compactly supported positive measure $\rho$ on $\mathbb C$, not necessarily of mass $1$, are the probability measures
$$p_\rho=\lim_{n\to\infty}\left(\left(1-\frac{c}{n}\right)\delta_0+\frac{1}{n}\,\rho\right)^{*n}\quad,\quad 
\pi_\rho=\lim_{n\to\infty}\left(\left(1-\frac{c}{n}\right)\delta_0+\frac{1}{n}\,\rho\right)^{\boxplus n}$$
where $c=mass(\rho)$, called compound Poisson and compound free Poisson laws.
\end{definition}

In what follows we will be interested in the case where $\rho$ is discrete, as is for instance the case for $\rho=t\delta_1$ with $t>0$, which produces the Poisson and free Poisson laws. The following result allows one to detect compound Poisson/free Poisson laws:

\begin{theorem}
For a discrete measure, $\rho=\sum_{i=1}^sc_i\delta_{z_i}$ with $c_i>0$ and $z_i\in\mathbb R$,
$$F_{p_\rho}(y)=\exp\left(\sum_{i=1}^sc_i(e^{iyz_i}-1)\right)$$
$$R_{\pi_\rho}(y)=\sum_{i=1}^s\frac{c_iz_i}{1-yz_i}$$
where $F,R$ denote respectively the Fourier transform, and Voiculescu's $R$-transform.
\end{theorem}

\begin{proof}
Let $\mu_n$ be the measure in Definition 5.28, under the convolution signs:
$$\mu_n=\left(1-\frac{c}{n}\right)\delta_0+\frac{1}{n}\rho$$

In the classical case, we have the following computation:
\begin{eqnarray*}
&&F_{\mu_n}(y)=\left(1-\frac{c}{n}\right)+\frac{1}{n}\sum_{i=1}^sc_ie^{iyz_i}\\
&\implies&F_{\mu_n^{*n}}(y)=\left(\left(1-\frac{c}{n}\right)+\frac{1}{n}\sum_{i=1}^sc_ie^{iyz_i}\right)^n\\
&\implies&F_{p_\rho}(y)=\exp\left(\sum_{i=1}^sc_i(e^{iyz_i}-1)\right)
\end{eqnarray*}

In the free case now, we use a similar method. The Cauchy transform of $\mu_n$ is:
$$G_{\mu_n}(\xi)=\left(1-\frac{c}{n}\right)\frac{1}{\xi}+\frac{1}{n}\sum_{i=1}^s\frac{c_i}{\xi-z_i}$$

Consider now the $R$-transform of the measure $\mu_n^{\boxplus n}$, which is given by:
$$R_{\mu_n^{\boxplus n}}(y)=nR_{\mu_n}(y)$$

The above formula of $G_{\mu_n}$ shows that the equation for $R=R_{\mu_n^{\boxplus n}}$ is as follows:
\begin{eqnarray*}
&&\left(1-\frac{c}{n}\right)\frac{1}{y^{-1}+R/n}+\frac{1}{n}\sum_{i=1}^s\frac{c_i}{y^{-1}+R/n-z_i}=y\\
&\implies&\left(1-\frac{c}{n}\right)\frac{1}{1+yR/n}+\frac{1}{n}\sum_{i=1}^s\frac{c_i}{1+yR/n-yz_i}=1
\end{eqnarray*}

Now multiplying by $n$, rearranging the terms, and letting $n\to\infty$, we get:
\begin{eqnarray*}
&&\frac{c+yR}{1+yR/n}=\sum_{i=1}^s\frac{c_i}{1+yR/n-yz_i}\\
&\implies&c+yR_{\pi_\rho}(y)=\sum_{i=1}^s\frac{c_i}{1-yz_i}\\
&\implies&R_{\pi_\rho}(y)=\sum_{i=1}^s\frac{c_iz_i}{1-yz_i}
\end{eqnarray*}

This finishes the proof in the free case, and we are done.
\end{proof}

We also have the following technical result, providing a useful alternative to Definition 5.28, in order to detect the classical and free compound Poisson laws:

\begin{theorem}
For a discrete measure, written as $\rho=\sum_{i=1}^sc_i\delta_{z_i}$ with $c_i>0$ and $z_i\in\mathbb R$, we have the classical/free formulae
$$p_\rho/\pi_\rho=law\left(\sum_{i=1}^sz_i\alpha_i\right)$$
where the variables $\alpha_i$ are Poisson/free Poisson$(c_i)$, independent/free.
\end{theorem}

\begin{proof}
Let $\alpha$ be the sum of Poisson/free Poisson variables in the statement:
$$\alpha=\sum_{i=1}^sz_i\alpha_i$$

By using some well-known Fourier transform formulae, we have:
\begin{eqnarray*}
F_{\alpha_i}(y)=\exp(c_i(e^{iy}-1))
&\implies&F_{z_i\alpha_i}(y)=\exp(c_i(e^{iyz_i}-1))\\
&\implies&F_\alpha(y)=\exp\left(\sum_{i=1}^sc_i(e^{iyz_i}-1)\right)
\end{eqnarray*}

Also, by using some well-known $R$-transform formulae, we have:
\begin{eqnarray*}
R_{\alpha_i}(y)=\frac{c_i}{1-y}
&\implies&R_{z_i\alpha_i}(y)=\frac{c_iz_i}{1-yz_i}\\
&\implies&R_\alpha(y)=\sum_{i=1}^s\frac{c_iz_i}{1-yz_i}
\end{eqnarray*}

Thus we have indeed the same formulae as those which are needed.
\end{proof}

We refer to \cite{bpa}, \cite{vdn} for the general theory here, to \cite{ba8}, \cite{csn} for representation theory aspects, and to \cite{mpa}, \cite{wig} for random matrix aspects. In what follows we will only need the main examples of classical and free compound Poisson laws, which are the classical and free Bessel laws, constructed as follows:

\index{Bessel law}
\index{free Bessel law}

\begin{definition}
The Bessel and free Bessel laws are the compound Poisson laws
$$b^s_t=p_{t\varepsilon_s}\quad,\quad \beta^s_t=\pi_{t\varepsilon_s}$$
where $\varepsilon_s$ is the uniform measure on the $s$-th roots unity. In particular:
\begin{enumerate}
\item At $s=1$ we obtain the usual Poisson and free Poisson laws, $p_t,\pi_t$.

\item At $s=2$ we obtain the ``real'' Bessel and free Bessel laws, denoted $p^2_t,\pi^2_t$.

\item At $s=\infty$ we obtain the ``complex'' Bessel and free Bessel laws, denoted $P_t,\Pi_t$.
\end{enumerate}
\end{definition}

There is a lot of theory regarding these laws, involving classical and quantum reflection groups, subfactors and planar algebras, and free probability and random matrices. We refer here to \cite{bb+}, where these laws were introduced. Let us just record here:

\begin{theorem}
The moments of the various central limiting measures, namely
$$
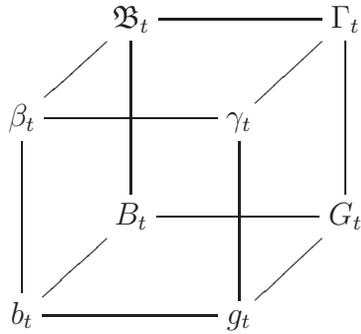
\xymatrix@R=20pt@C=22pt{
&\Pi_t\ar@{-}[rr]\ar@{-}[dd]&&\Gamma_t\ar@{-}[dd]\\
\pi^2_t\ar@{-}[rr]\ar@{-}[dd]\ar@{-}[ur]&&\gamma_t\ar@{-}[dd]\ar@{-}[ur]\\
&P_t\ar@{-}[rr]\ar@{-}[uu]&&G_t\ar@{-}[uu]\\
p^2_t\ar@{-}[uu]\ar@{-}[ur]\ar@{-}[rr]&&g_t\ar@{-}[uu]\ar@{-}[ur]
}$$
are always given by the same formula, involving partitions, namely
$$M_k=\sum_{\pi\in D(k)}t^{|\pi|}$$
with the sets of partitions $D(k)$ in question being respectively
$$
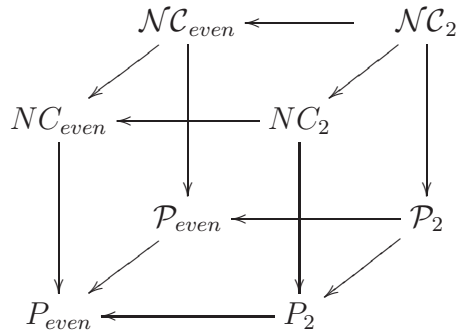
\xymatrix@R=20pt@C=2pt{
&\mathcal{NC}_{even}\ar[dl]\ar[dd]&&\ \ \ \mathcal{NC}_2\ \ \ \ar[ll]\ar[dd]\ar[dl]\\
NC_{even}\ar[dd]&&NC_2\ar[ll]\ar[dd]\\
&\mathcal P_{even}\ar[dl]&&\mathcal P_2\ar[ll]\ar[dl]\\
P_{even}&&P_2\ar[ll]
}$$
and with $|.|$ being the number of blocks. 
\end{theorem}

\begin{proof}
This follows indeed from our various moment results. See \cite{bb+}.
\end{proof}

Getting back now to our quantum reflection groups, we have:

\index{Bessel law}
\index{free Bessel law}
\index{Bercovici-Pata bijection}

\begin{theorem}
With $N\to\infty$, the laws of truncated characters are as follows:
\begin{enumerate}
\item For $H_N$ we obtain the Bessel law $p^2_t$.

\item For $H_N^+$ we obtain the free Bessel law $\pi^2_t$.

\item For $K_N$ we obtain the complex Bessel law $P_t$.

\item For $K_N^+$ we obtain the complex free Bessel law $\Pi_t$.
\end{enumerate}
Also, we have the Bercovici-Pata bijection for truncated characters.
\end{theorem}

\begin{proof}
At $t=1$ this follows by counting the partitions, a bit as in the continuous case, in the proof of Theorem 5.18. At $t\in(0,1)$ this is routine, by using the Weingarten formula, as in the continuous case, in the proof of Theorem 5.23. See \cite{bb+}.
\end{proof}

The results that we have so far, for the quantum unitary and refelection groups, are quite interesting, from a theoretical probability perspective, because we have:

\begin{theorem}
The laws of the truncated characters for the basic quantum groups,
$$
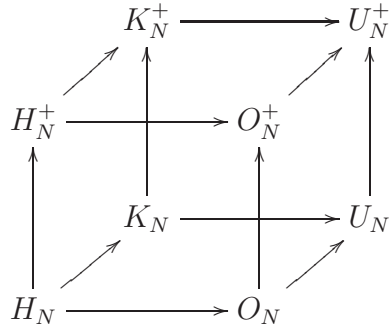
\xymatrix@R=18pt@C=18pt{
&K_N^+\ar[rr]&&U_N^+\\
H_N^+\ar[rr]\ar[ur]&&O_N^+\ar[ur]\\
&K_N\ar[rr]\ar[uu]&&U_N\ar[uu]\\
H_N\ar[uu]\ar[ur]\ar[rr]&&O_N\ar[uu]\ar[ur]
}$$
and the various classical and free central limiting measures, namely
$$
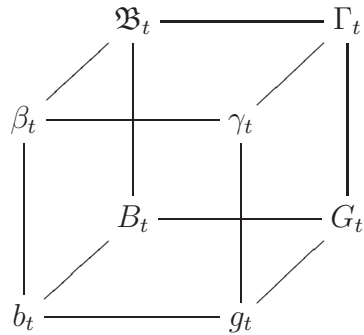
\xymatrix@R=20pt@C=22pt{
&\Pi_t\ar@{-}[rr]\ar@{-}[dd]&&\Gamma_t\ar@{-}[dd]\\
\pi^2_t\ar@{-}[rr]\ar@{-}[dd]\ar@{-}[ur]&&\gamma_t\ar@{-}[dd]\ar@{-}[ur]\\
&P_t\ar@{-}[rr]\ar@{-}[uu]&&G_t\ar@{-}[uu]\\
p^2_t\ar@{-}[uu]\ar@{-}[ur]\ar@{-}[rr]&&g_t\ar@{-}[uu]\ar@{-}[ur]
}$$
in the $N\to\infty$ limit.
\end{theorem}

\begin{proof}
This follows indeed by putting together the various results discussed above, concerning general free probability theory, and our computations here.
\end{proof}

Regarding now the tori, the situation here is more complicated, no longer involving the Bercovici-Pata bijection. Let us recall indeed that our tori and their duals are:
$$\xymatrix@R=50pt@C=50pt{
T_N^+\ar[r]&\mathbb T_N^+\\
T_N\ar[r]\ar[u]&\mathbb T_N\ar[u]
}\qquad
\xymatrix@R=25pt@C=50pt{\\ :\\}
\qquad
\xymatrix@R=50pt@C=50pt{
\mathbb Z_2^{*N}\ar[d]&F_N\ar[l]\ar[d]\\
\mathbb Z_2^N&\mathbb Z^N\ar[l]
}$$

We are interested in the computation of the laws of the associated truncated characters, which are the following variables, with $g_1,\ldots,g_N$ being the group generators: 
$$\chi_t=g_1+g_2+\ldots+g_{[tN]}$$

\index{Meixner law}
\index{free Meixner law}

By dilation we can assume $t=1$. For the complex tori, $\mathbb T_N\subset\mathbb T_N^+$, we are led into the computation of the Kesten measures for $F_N\to\mathbb Z^N$, and so into the Meixner/free Meixner correspondence. As for the real tori, $T_N\subset T_N^+$, here we are led into the computation of the Kesten measures for $\mathbb Z_2^{*N}\to\mathbb Z_2^N$, and so into a real version of this correspondence. These are both quite technical questions, that we will not get into, here.

\section*{5e. Exercises} 

As a first exercise in relation with the material in this chapter, we have:

\begin{exercise}
Work out the Weingarten formula for the classical spheres
$$S^{N-1}_\mathbb R\quad,\quad S^{N-1}_\mathbb C$$
in general, then at small $N\in\mathbb N$, and at big $N\in\mathbb N$, and find some applications of this.
\end{exercise}

Here the application part is a bit up to you. The classical spheres are very classical objects, and they appear in connection with many questions.

\begin{exercise}
Compute the laws of the truncated characters for the main tori
$$\xymatrix@R=50pt@C=50pt{
T_N^+\ar[r]&\mathbb T_N^+\\
T_N\ar[r]\ar[u]&\mathbb T_N\ar[u]
}$$
and then work out the asymptotics, with $N\to\infty$.
\end{exercise}

This is something that we briefly discussed at the end of this chapter, and there is definitely some interesting work to be done here.

\chapter{Basic manifolds}

\section*{6a. Partial isometries}

In this chapter and in the next two ones we keep building on the work started in the previous chapter, by systematically developing the real and complex free geometry. We will extend the family of objects $(S,T,U,K)$ that we have, first with some general homogeneous spaces, of ``quantum partial isometries'', and then with some generalizations of these spaces, that we will call ``affine homogeneous spaces''. We will also discuss, at the end of chapter 8, the axiomatization problem for the free manifolds.  

\bigskip

We will insist on probabilistic aspects, and this for two reasons. First, due to our belief, extensively explained in the beginning of the previous chapter, that given a manifold $X$, the thing to do with it is to compute its integration functional $tr:C(X)\to\mathbb C$, via a formula as explicit as possible, with the idea in mind that the various applications of $X$ to physics should involve precisely this integration functional $tr:C(X)\to\mathbb C$.

\bigskip

But then, there is a second reason as well, more technical and subtle. There are all sorts of ways of talking about ``liberation'', meaning operations of type $X\to X^+$, involving a group, a homogeneous space, or a more general manifold $X$. And the range of things that can be said here is endless, including the good, the bad, and the ugly:

\bigskip

(1) The ugly, to start with, means doing whatever quick algebra, without any idea in mind, such as erasing some commutation relations $ab=ba$, and then saying that you're done. With, as an illustrating example, talking about $\mathbb R^N\to\mathbb R^N_+$, simply by saying that $\mathbb R^N_+$ corresponds to the complex algebra $A=<x_1,\ldots,x_N>$ generated by $N$ free variables, which algebra $A$ has nothing to do with analysis and physics, in our opinion.

\bigskip

(2) The bad is something more subtle, meaning knowing what you're doing, but doing it badly. For instance the results $G(T_N)=H_N$ and $G^+(T_N)=\bar{O}_N$, from \cite{bbc}, suggest that $H_N\to\bar{O}_N$ is a liberation operation. Which might sound reasonable, algebrically speaking, but which analytically is something totally unplausible, because how on Earth could the free analogues of the Bessel laws for $H_N$ be the Gaussian laws for $\bar{O}_N$.

\bigskip

(3) The good, and no surprise here, is that of talking about $X\to X^+$, with a full knowledge of what this operation does, both algebrically and analytically, to the point of being 100\% sure that this is a ``true liberation''. And also, as per general mathematical physics requirements, with at least 1 motivation from physics in mind. As an example here, the liberation operation $H_N\to H_N^+$, also from \cite{bbc}, fulfills these requirements.

\bigskip

The thing now is that, in all the above, the delicate knowledge to be mastered is the analytic one. So, when could we say that $X\to X^+$ is a liberation, analytically speaking? And here the answer, coming from years of work and observations, including \cite{bbc}, is that ``this happens when $X\to X^+$ is compatible with the Bercovici-Pata bijection \cite{bpa}, or with some other well-established bijection from probability, such as the Meixner/free Meixner one''. Of course this is a bit vague, because our principle does not tell us at what exact variables to look at, and also there is usually a $N\to\infty$ limiting procedure appearing there, and so on. But, in practice, this remains an excellent principle, whose verification requires a good knowledge of the integration functional $tr:C(X)\to\mathbb C$.

\bigskip

In short, and getting back now to what was said before, we will be interested in what follows in homogeneous spaces $X$, and their integration functionals $tr:C(X)\to\mathbb C$. There has been quite some work on this subject, after the 2010 paper \cite{bgo} regarding the spheres, which launched everything, notably with the fundamental paper \cite{bsk}, and then \cite{ba5}, and then \cite{ba6}, that we will follow here, in this chapter and in the next two ones. 

\bigskip

In order to get started, with this program, we will first discuss, in this chapter, a class of homogeneous spaces which are of fairly general type, as follows:
$$X=(G_M\times G_N)\big/(G_L\times G_{M-L}\times G_{N-L})$$

These spaces cover indeed the quantum groups and the spheres. Also, they are quite concrete and useful objects, consisting of certain classes of ``partial isometries''. And also, importantly, in the discrete case, where $G=(G_N)$ is one of our easy quantum reflection groups, these spaces are very interesting, combinatorially. But more on this later.

\bigskip

We begin with a study in the classical case. Our starting point will be:

\index{partial isometry}

\begin{definition}
Associated to any integers $L\leq M,N$ are the spaces
$$O_{MN}^L=\left\{T:E\to F\ {\rm isometry}\Big|E\subset\mathbb R^N,F\subset\mathbb R^M,\dim_\mathbb RE=L\right\}$$
$$U_{MN}^L=\left\{T:E\to F\ {\rm isometry}\Big|E\subset\mathbb C^N,F\subset\mathbb C^M,\dim_\mathbb CE=L\right\}$$
where the notion of isometry is with respect to the usual real/complex scalar products.
\end{definition}

These spaces remind basic algebraic geometry, that you can learn by the way from Shafarevich \cite{sha}, or Harris \cite{har}, or Hartshorne \cite{hrt}, and more specifically Grassmannians, flag manifolds, Stiefel manifolds, and so on. However, all these latter manifolds are in fact projective, and our policy in our book will be to discuss projective geometry only at the end, in chapters 15-16 below. Thus, more on Grassmannians and related manifolds later, and for the moment we will stay affine, and use Definition 6.1 as it is.

\bigskip

As a first observation, in relation with our $(S,T,U,K)$ objects, it follows from definitions that at $L=M=N$ we obtain the orthogonal and unitary groups $O_N,U_N$:
$$O_{NN}^N=O_N\quad,\quad 
U_{NN}^N=U_N$$ 

Another interesting specialization is $L=M=1$. Here the elements of $O_{1N}^1$ are the isometries $T:E\to\mathbb R$, with $E\subset\mathbb R^N$ one-dimensional. But such an isometry is uniquely determined by $T^{-1}(1)\in\mathbb R^N$, which must belong to $S^{N-1}_\mathbb R$. Thus, we have $O_{1N}^1=S^{N-1}_\mathbb R$. Similarly, in the complex case we have $U_{1N}^1=S^{N-1}_\mathbb C$, and so our results here are:
$$O_{1N}^1=S^{N-1}_\mathbb R\quad,\quad 
U_{1N}^1=S^{N-1}_\mathbb C$$

Yet another interesting specialization is $L=N=1$. Here the elements of $O_{1N}^1$ are the isometries $T:\mathbb R\to F$, with $F\subset\mathbb R^M$ one-dimensional. But such an isometry is uniquely determined by $T(1)\in\mathbb R^M$, which must belong to $S^{M-1}_\mathbb R$. Thus, we have $O_{M1}^1=S^{M-1}_\mathbb R$. Similarly, in the complex case we have $U_{M1}^1=S^{M-1}_\mathbb C$, and so our results here are:
$$O_{M1}^1=S^{M-1}_\mathbb R\quad,\quad 
U_{M1}^1=S^{M-1}_\mathbb C$$

In general, the most convenient is to view the elements of $O_{MN}^L,U_{MN}^L$ as rectangular matrices, and to use matrix calculus for their study. We have indeed:

\begin{proposition}
We have identifications of compact spaces
$$O_{MN}^L\simeq\left\{U\in M_{M\times N}(\mathbb R)\Big|UU^t={\rm projection\ of\ trace}\ L\right\}$$
$$U_{MN}^L\simeq\left\{U\in M_{M\times N}(\mathbb C)\Big|UU^*={\rm projection\ of\ trace}\ L\right\}$$
with each partial isometry being identified with the corresponding rectangular matrix.
\end{proposition}

\begin{proof}
We can indeed identify the partial isometries $T:E\to F$ with their corresponding extensions $U:\mathbb R^N\to\mathbb R^M$, $U:\mathbb C^N\to\mathbb C^M$, obtained by setting:
$$U_{E^\perp}=0$$

Then, we can identify these latter linear maps $U$ with the corresponding rectangular matrices, and we are led to the conclusion in the statement.
\end{proof}

As an illustration, at $L=M=N$ we recover in this way the usual matrix description of $O_N,U_N$. Also, at $L=M=1$ we obtain the usual description of $S^{N-1}_\mathbb R,S^{N-1}_\mathbb C$, as row spaces over the corresponding groups $O_N,U_N$.  Finally, at $L=N=1$ we obtain the usual description of $S^{N-1}_\mathbb R,S^{N-1}_\mathbb C$, as column spaces over the corresponding groups $O_N,U_N$. 

\bigskip

Now back to the general case, observe that the isometries $T:E\to F$, or rather their extensions $U:\mathbb K^N\to\mathbb K^M$, with $\mathbb K=\mathbb R,\mathbb C$, obtained by setting $U_{E^\perp}=0$, can be composed with the isometries of $\mathbb K^M,\mathbb K^N$, according to the following scheme:
$$\xymatrix@R=18mm@C=18mm{
\mathbb K^N\ar[r]^{B^*}&\mathbb K^N\ar@.[r]^U&\mathbb K^M\ar[r]^A&\mathbb K^M\\
B(E)\ar@.[r]\ar[u]&E\ar[r]^T\ar[u]&F\ar@.[r]\ar[u]&A(F)\ar[u]
}$$

With the identifications in Proposition 6.2 made, the precise statement here is:

\begin{proposition}
We have an action map as follows, which is transitive,
$$O_M\times O_N\curvearrowright O_{MN}^L\quad,\quad 
(A,B)U=AUB^t$$
as well as an action map as follows, transitive as well,
$$U_M\times U_N\curvearrowright U_{MN}^L\quad,\quad 
(A,B)U=AUB^*$$
whose stabilizers are respectively the following groups,
$$O_L\times O_{M-L}\times O_{N-L}$$
$$U_L\times U_{M-L}\times U_{N-L}$$
computed by using the action on $U=(^1_0{\ }^0_0)$.
\end{proposition}

\begin{proof}
We have indeed action maps as in the statement, which are transitive. Let us compute now the stabilizer $G$ of the following point:
$$U=\begin{pmatrix}1&0\\0&0\end{pmatrix}$$

Since $(A,B)\in G$ satisfy $AU=UB$, their components must be of the following form:
$$A=\begin{pmatrix}x&*\\0&a\end{pmatrix}\quad,\quad 
B=\begin{pmatrix}x&0\\ *&b\end{pmatrix}$$

Now since $A,B$ are both unitaries, these matrices follow to be block-diagonal, and so:
$$G=\left\{(A,B)\Big|A=\begin{pmatrix}x&0\\0&a\end{pmatrix},B=\begin{pmatrix}x&0\\ 0&b\end{pmatrix}\right\}$$

The stabilizer of $U$ is then parametrized by triples $(x,a,b)$ belonging respectively to:$$O_L\times O_{M-L}\times O_{N-L}$$
$$U_L\times U_{M-L}\times U_{N-L}$$

Thus, we are led to the conclusion in the statement.
\end{proof}

Let us work out now the quotient space description of $O_{MN}^L,U_{MN}^L$. We have here:

\begin{theorem}
We have isomorphisms of homogeneous spaces as follows,
\begin{eqnarray*}
O_{MN}^L&=&(O_M\times O_N)/(O_L\times O_{M-L}\times O_{N-L})\\
U_{MN}^L&=&(U_M\times U_N)/(U_L\times U_{M-L}\times U_{N-L})
\end{eqnarray*}
with the quotient maps being given by $(A,B)\to AUB^*$, where $U=(^1_0{\ }^0_0)$.
\end{theorem}

\begin{proof}
This is just a reformulation of Proposition 6.3, by taking into account the fact that the fixed point used in the proof there was $U=(^1_0{\ }^0_0)$.
\end{proof}

Once again, the basic examples here come from the cases $L=M=N$ and $L=M=1$. At $L=M=N$ the quotient spaces at right are respectively:
$$O_N\quad,\quad U_N$$

At $L=M=1$ the quotient spaces at right are respectively:
$$O_N/O_{N-1}\quad,\quad U_N/U_{N-1}$$

In fact, in the general $L=M$ case we obtain the following spaces:
$$O_{MN}^M
=(O_M\times O_N)/(O_M\times O_{N-M})
=O_N/O_{N-M}$$
$$U_{MN}^M
=(U_M\times U_N)/(U_M\times U_{N-M})
=U_N/U_{N-M}$$

Similarly, the examples coming from the cases $L=M=N$ and $L=N=1$ are particular cases of the general $L=N$ case, where we obtain the following spaces:
$$O_{MN}^N
=(O_M\times O_N)/(O_M\times O_{M-N})
=O_N/O_{M-N}$$
$$U_{MN}^N
=(U_M\times U_N)/(U_M\times U_{M-N})
=U_N/U_{M-N}$$

Summarizing, in relation with our previous $(S,T,U,K)$ objects, we have here homogeneous spaces which unify the spheres with the unitary quantum groups.

\section*{6b. Free isometries}

We can now liberate the spaces $O_{MN}^L,U_{MN}^L$, as follows:

\index{free partial isometry}

\begin{definition}
Associated to any integers $L\leq M,N$ are the algebras
\begin{eqnarray*}
C(O_{MN}^{L+})&=&C^*\left((u_{ij})_{i=1,\ldots,M,j=1,\ldots,N}\Big|u=\bar{u},uu^t={\rm projection\ of\ trace}\ L\right)\\
C(U_{MN}^{L+})&=&C^*\left((u_{ij})_{i=1,\ldots,M,j=1,\ldots,N}\Big|uu^*,\bar{u}u^t={\rm projections\ of\ trace}\ L\right)
\end{eqnarray*}
with the trace being by definition the sum of the diagonal entries.
\end{definition}

Observe that the above universal algebras are indeed well-defined, as it was previously  the case for the free spheres, and this due to the trace conditions, which read: 
$$\sum_{ij}u_{ij}u_{ij}^*
=\sum_{ij}u_{ij}^*u_{ij}
=L$$

We have inclusions between the various spaces constructed so far, as follows:
$$\xymatrix@R=15mm@C=15mm{
O_{MN}^{L+}\ar[r]&U_{MN}^{L+}\\
O_{MN}^L\ar[r]\ar[u]&U_{MN}^L\ar[u]}$$

At the level of basic examples now, we first have the following result:

\begin{proposition}
At $L=M=1$ we obtain the following diagram,
$$\xymatrix@R=15mm@C=15mm{
S^{N-1}_{\mathbb R,+}\ar[r]&S^{N-1}_{\mathbb C,+}\\
S^{N-1}_\mathbb R\ar[r]\ar[u]&S^{N-1}_\mathbb C\ar[u]}$$
consisting of the basic real and complex spheres.
\end{proposition}

\begin{proof}
We recall that the various spheres involved are constructed as follows, with the symbol $\times$ standing for ``commutative'' and ``free'', respectively:
\begin{eqnarray*}
C(S^{N-1}_{\mathbb R,\times})&=&C^*_\times\left(z_1,\ldots,z_N\Big|z_i=z_i^*,\sum_iz_i^2=1\right)\\
C(S^{N-1}_{\mathbb C,\times})&=&C^*_\times\left(z_1,\ldots,z_N\Big|\sum_iz_iz_i^*=\sum_iz_i^*z_i=1\right)
\end{eqnarray*}

Now by comparing with the definition of $O_{1N}^{1\times},U_{1N}^{1\times}$, this proves our claim.
\end{proof}

Similarly, we have as well the following result:

\begin{proposition}
At $L=N=1$ we obtain the following diagram,
$$\xymatrix@R=15mm@C=15mm{
S^{M-1}_{\mathbb R,+}\ar[r]&S^{M-1}_{\mathbb C,+}\\
S^{M-1}_\mathbb R\ar[r]\ar[u]&S^{M-1}_\mathbb C\ar[u]}$$
consisting of the basic real and complex spheres.
\end{proposition}

\begin{proof}
This is similar to the proof of Proposition 6.6, coming from the definition of the various spheres involved, via some standard identifications.
\end{proof}

Finally, again regarding examples, we have as well the following result:

\begin{theorem}
At $L=M=N$ we obtain the following diagram,
$$\xymatrix@R=15mm@C=15mm{
O_N^+\ar[r]&U_N^+\\
O_N\ar[r]\ar[u]&U_N\ar[u]}$$
consisting of the groups $O_N,U_N$, and their liberations.
\end{theorem}

\begin{proof}
We recall that the various quantum groups in the statement are constructed as follows, with the symbol $\times$ standing once again for ``commutative'' and ``free'':
\begin{eqnarray*}
C(O_N^\times)&=&C^*_\times\left((u_{ij})_{i,j=1,\ldots,N}\Big|u=\bar{u},uu^t=u^tu=1\right)\\
C(U_N^\times)&=&C^*_\times\left((u_{ij})_{i,j=1,\ldots,N}\Big|uu^*=u^*u=1,\bar{u}u^t=u^t\bar{u}=1\right)
\end{eqnarray*}

On the other hand, according to Proposition 6.2 and to Definition 6.5, we have the following presentation results:
\begin{eqnarray*}
C(O_{NN}^{N\times})&=&C^*_\times\left((u_{ij})_{i,j=1,\ldots,N}\Big|u=\bar{u},uu^t={\rm projection\ of\ trace}\ N\right)\\
C(U_{NN}^{N\times})&=&C^*_\times\left((u_{ij})_{i,j=1,\ldots,N}\Big|uu^*,\bar{u}u^t={\rm projections\ of\ trace}\ N\right)
\end{eqnarray*}

We use now the standard fact that if $p=aa^*$ is a projection then $q=a^*a$ is a projection too. We use as well the following formulae:
$$Tr(uu^*)=Tr(u^t\bar{u})\quad,\quad 
Tr(\bar{u}u^t)=Tr(u^*u)$$

We therefore obtain the following formulae:
\begin{eqnarray*}
C(O_{NN}^{N\times})&=&C^*_\times\left((u_{ij})_{i,j=1,\ldots,N}\Big|u=\bar{u},\ uu^t,u^tu={\rm projections\ of\ trace}\ N\right)\\
C(U_{NN}^{N\times})&=&C^*_\times\left((u_{ij})_{i,j=1,\ldots,N}\Big|uu^*,u^*u,\bar{u}u^t,u^t\bar{u}={\rm projections\ of\ trace}\ N\right)
\end{eqnarray*}

Now observe that the conditions on the right are all of the form $(tr\otimes id)p=1$. To be more precise, $p$ must be as follows, for the above conditions:
$$p=uu^*,u^*u,\bar{u}u^t,u^t\bar{u}$$

We therefore obtain that, for any faithful state $\varphi$, we have:
$$(tr\otimes\varphi)(1-p)=0$$  

But this shows that the following projections must be all equal to the identity:
$$p=uu^*,u^*u,\bar{u}u^t,u^t\bar{u}$$

Thus, we are led to the conclusion in the statement.
\end{proof}

Regarding now the homogeneous space structure of $O_{MN}^{L\times},U_{MN}^{L\times}$, the situation here is more complicated in the free case than in the classical case, due to a number of reasons, of both algebraic and analytic nature. We first have the following result:

\begin{proposition}
The spaces $U_{MN}^{L\times}$ have the following properties:
\begin{enumerate}
\item We have an action $U_M^\times\times U_N^\times\curvearrowright U_{MN}^{L\times}$, given by:
$$u_{ij}\to\sum_{kl}u_{kl}\otimes a_{ki}\otimes b_{lj}^*$$

\item We have a map $U_M^\times\times U_N^\times\to U_{MN}^{L\times}$, given by: 
$$u_{ij}\to\sum_{r\leq L}a_{ri}\otimes b_{rj}^*$$
\end{enumerate}
Similar results hold for the spaces $O_{MN}^{L\times}$, with all the $*$ exponents removed.
\end{proposition}

\begin{proof}
In the classical case, consider the action and quotient maps:
$$U_M\times U_N\curvearrowright U_{MN}^L$$
$$U_M\times U_N\to U_{MN}^L$$

The transposes of these two maps are as follows, where $J=(^1_0{\ }^0_0)$:
\begin{eqnarray*}
\varphi&\to&((U,A,B)\to\varphi(AUB^*))\\
\varphi&\to&((A,B)\to\varphi(AJB^*))
\end{eqnarray*}

But with $\varphi=u_{ij}$ we obtain precisely the formulae in the statement. The proof in the orthogonal case is similar. Regarding now the free case, the proof goes as follows:

\medskip

(1) Assuming $uu^*u=u$, let us set:
$$U_{ij}=\sum_{kl}u_{kl}\otimes a_{ki}\otimes b_{lj}^*$$

We have then the following computation:
\begin{eqnarray*}
(UU^*U)_{ij}
&=&\sum_{pq}\sum_{klmnst}u_{kl}u_{mn}^*u_{st}\otimes a_{ki}a_{mq}^*a_{sq}\otimes b_{lp}^*b_{np}b_{tj}^*\\
&=&\sum_{klmt}u_{kl}u_{ml}^*u_{mt}\otimes a_{ki}\otimes b_{tj}^*\\
&=&\sum_{kt}u_{kt}\otimes a_{ki}\otimes b_{tj}^*\\
&=&U_{ij}
\end{eqnarray*}

Also, assuming that we have $\sum_{ij}u_{ij}u_{ij}^*=L$, we obtain:
\begin{eqnarray*}
\sum_{ij}U_{ij}U_{ij}^*
&=&\sum_{ij}\sum_{klst}u_{kl}u_{st}^*\otimes a_{ki}a_{si}^*\otimes b_{lj}^*b_{tj}\\
&=&\sum_{kl}u_{kl}u_{kl}^*\otimes1\otimes1\\
&=&L
\end{eqnarray*}

(2) Assuming $uu^*u=u$, let us set:
$$V_{ij}=\sum_{r\leq L}a_{ri}\otimes b_{rj}^*$$

We have then the following computation:
\begin{eqnarray*}
(VV^*V)_{ij}
&=&\sum_{pq}\sum_{x,y,z\leq L}a_{xi}a_{yq}^*a_{zq}\otimes b_{xp}^*b_{yp}b_{zj}^*\\
&=&\sum_{x\leq L}a_{xi}\otimes b_{xj}^*\\
&=&V_{ij}
\end{eqnarray*}

Also, assuming that we have $\sum_{ij}u_{ij}u_{ij}^*=L$, we obtain:
\begin{eqnarray*}
\sum_{ij}V_{ij}V_{ij}^*
&=&\sum_{ij}\sum_{r,s\leq L}a_{ri}a_{si}^*\otimes b_{rj}^*b_{sj}\\
&=&\sum_{l\leq L}1\\
&=&L
\end{eqnarray*}

By removing all the $*$ exponents, we obtain as well the orthogonal results.
\end{proof}

Let us examine now the relation between the above maps. In the classical case, given a quotient space $X=G/H$, the associated action and quotient maps are given by:
$$\begin{cases}
a:X\times G\to X&:\quad (Hg,h)\to Hgh\\
p:G\to X&:\quad g\to Hg
\end{cases}$$

Thus we have $a(p(g),h)=p(gh)$. In our context, a similar result holds: 

\begin{theorem}
With $G=G_M\times G_N$ and $X=G_{MN}^L$, where $G_N=O_N^\times,U_N^\times$, we have
$$\xymatrix@R=15mm@C=30mm{
G\times G\ar[r]^m\ar[d]_{p\times id}&G\ar[d]^p\\
X\times G\ar[r]^a&X
}$$
where $a,p$ are the action map and the map constructed in Proposition 6.9.
\end{theorem}

\begin{proof}
At the level of the associated algebras of functions, we must prove that the following diagram commutes, where $\Phi,\alpha$ are morphisms of algebras induced by $a,p$:
$$\xymatrix@R=15mm@C=25mm{
C(X)\ar[r]^\Phi\ar[d]_\alpha&C(X\times G)\ar[d]^{\alpha\otimes id}\\
C(G)\ar[r]^\Delta&C(G\times G)
}$$

When going right, and then down, the composition is as follows:
\begin{eqnarray*}
(\alpha\otimes id)\Phi(u_{ij})
&=&(\alpha\otimes id)\sum_{kl}u_{kl}\otimes a_{ki}\otimes b_{lj}^*\\
&=&\sum_{kl}\sum_{r\leq L}a_{rk}\otimes b_{rl}^*\otimes a_{ki}\otimes b_{lj}^*
\end{eqnarray*}

On the other hand, when going down, and then right, the composition is as follows, where $F_{23}$ is the flip between the second and the third components:
\begin{eqnarray*}
\Delta\pi(u_{ij})
&=&F_{23}(\Delta\otimes\Delta)\sum_{r\leq L}a_{ri}\otimes b_{rj}^*\\
&=&F_{23}\left(\sum_{r\leq L}\sum_{kl}a_{rk}\otimes a_{ki}\otimes b_{rl}^*\otimes b_{lj}^*\right)
\end{eqnarray*}

Thus the above diagram commutes indeed, and this gives the result.
\end{proof}

Summarizing, we have so far free analogues of the spaces of partial isometries $O_{MN}^L$ and $U_{MN}^L$, along with some information about their homogeneous space structure, which looks quite axiomatic, as formulated in Theorem 6.10. There are many things to be done, as a continuation of this, and we will do this slowly, our plan being as follows:

\bigskip

(1) In the remainder of this chapter we will discuss as well discrete versions of the above constructions, and then we will go into the thing to be done, namely study of the Haar functional, and verification of the Bercovici-Pata bijection. 

\bigskip

(2) And then, in chapters 7-8 below, we will discuss more abstract or more concrete versions of these constructions, following \cite{bss} and related papers, the idea being that both generalizing and particularizing are interesting topics to be discussed.

\bigskip

As a general comment now, I can feel that you are a bit puzzled by our strategy, because we are talking here about homogeneous spaces, without knowing what an homogeneous space is, in the quantum setting. To which I would answer, please relax, there is absolutely no hurry with that. We have our spaces, which is a good thing, our study is on the way, another good thing, and the discussion of the homogeneous space structure, which will be something abstract, inspired from what we found in Theorem 6.10, and which will bring 0 advances on our problems to be solved, will be surely done at some point, and more specifically in chapter 7 below, but absolutely no hurry with that.

\bigskip

Hope you got my point, the quantum homogeneous spaces are not the same thing as the classical homogeneous spaces, and their study is quite tricky, following a different path. That's how the quantum world is, sometimes similar to the classical one, but sometimes very different. In case you are not convinced, pick a sphere $S$, as those studied so far in this book, and try writing that as a quotient space, and deducing from this results which are better than those established so far in this book, about such spheres $S$.

\bigskip

That will not work. And you will join here cohorts of mathematicians, having tried to develop theories of quantum homogeneous spaces, in nice and gentle analogy with the theory of classical homogeneous spaces, with quite average results. In fact, it is the paper \cite{bgo}, dealing with noncommutative spheres $S$ in a radical new way, via basic algebra and probability, that launched the modern theory, that we are explaining here.

\section*{6c. Discrete extensions}

Let us discuss now some extensions of the above constructions. We will be mostly interested in the quantum reflection groups, but let us first discuss, with full details, the case of the quantum groups $S_N,S_N^+$. The starting point is the semigroup $\widetilde{S}_N$ of partial permutations. This is a quite familiar object in combinatorics, defined as follows:

\index{partial permutation}
\index{semigroup of partial permutations}

\begin{definition}
$\widetilde{S}_N$ is the semigroup of partial permutations of $\{1\,\ldots,N\}$,
$$\widetilde{S}_N=\left\{\sigma:X\simeq Y\Big|X,Y\subset\{1,\ldots,N\}\right\}$$
with the usual composition operation, $\sigma'\sigma:\sigma^{-1}(X'\cap Y)\to\sigma'(X'\cap Y)$.
\end{definition}

Observe that $\widetilde{S}_N$ is not simplifiable, because the null permutation $\emptyset\in\widetilde{S}_N$, having the empty set as domain/range, satisfies $\emptyset\sigma=\sigma\emptyset=\emptyset$, for any $\sigma\in\widetilde{S}_N$. Observe also that $\widetilde{S}_N$ has a ``subinverse'' map, sending $\sigma:X\to Y$ to its usual inverse $\sigma^{-1}:Y\simeq X$.

\bigskip

A first interesting result about this semigroup $\widetilde{S}_N$, which shows that we are dealing here with some non-trivial combinatorics, is as follows:

\begin{proposition}
The number of partial permutations is given by
$$|\widetilde{S}_N|=\sum_{k=0}^Nk!\binom{N}{k}^2$$
that is, $1,2,7,34,209,\ldots\,$, and we have the estimate
$$|\widetilde{S}_N|\simeq N!\sqrt{\frac{\exp(4\sqrt{N}-1)}{4\pi\sqrt{N}}}$$
in the $N\to\infty$ limit.
\end{proposition}

\begin{proof}
The first assertion is clear, because in order to construct a partial permutation $\sigma:X\to Y$ we must choose an integer $k=|X|=|Y|$, then we must pick two subsets $X,Y\subset\{1,\ldots,N\}$ having cardinality $k$, and there are $\binom{N}{k}$ choices for each, and finally we must construct a bijection $\sigma:X\to Y$, and there are $k!$ choices here. As for the estimate, which is non-trivial, this is however something standard, and well-known.
\end{proof}

Another result, which is trivial, but quite fundamental, is as follows:

\begin{proposition}
We have a semigroup embedding $u:\widetilde{S}_N\subset M_N(0,1)$, defined by 
$$u_{ij}(\sigma)=
\begin{cases}
1&{\rm if}\ \sigma(j)=i\\
0&{\rm otherwise}
\end{cases}$$
whose image are the matrices having at most one nonzero entry, on each row and column.
\end{proposition}

\begin{proof}
This is trivial from definitions, with $u:\widetilde{S}_N\subset M_N(0,1)$ extending the standard embedding $u:S_N\subset M_N(0,1)$, that we have been heavily using, so far.
\end{proof}

Let us discuss now the construction and main properties of the semigroup of quantum partial permutations $\widetilde{S}_N^+$, in analogy with the above. For this purpose, we use the above embedding $u:\widetilde{S}_N\subset M_N(0,1)$. Due to the formula $u_{ij}(\sigma)=\delta_{i\sigma(j)}$, the matrix $u=(u_{ij})$ is ``submagic'', in the sense that its entries are projections, which are pairwise orthogonal on each row and column. This suggests the following definition:

\index{submagic matrix}
\index{free partial permutation}
\index{quantum semigroup}
\index{subantipode map}

\begin{definition}
$C(\widetilde{S}_N^+)$ is the universal $C^*$-algebra generated by the entries of a $N\times N$ submagic matrix $u$, with comultiplication and counit maps given by
$$\Delta(u_{ij})=\sum_ku_{ik}\otimes u_{kj}$$
$$\varepsilon(u_{ij})=\delta_{ij}$$
where submagic means formed of projections, which are pairwise orthogonal on rows and columns. We call $\widetilde{S}_N^+$ semigroup of quantum partial permutations of $\{1,\ldots,N\}$.
\end{definition}

Here the fact that the morphisms of algebras $\Delta,\varepsilon$ as above exist indeed follows from the universality property of $C(\widetilde{S}_N^+)$, with the needed submagic checks being nearly identical to the magic checks for $C(S_N^+)$, from chapter 2. Observe also that the morphisms $\Delta,\varepsilon$ satisfy the usual axioms for a comultiplication and antipode, namely:
$$(\Delta\otimes id)\Delta=(id\otimes \Delta)\Delta$$
$$(\varepsilon\otimes id)\Delta=(id\otimes\varepsilon)\Delta=id$$

Thus, we have a bialgebra structure of $C(\widetilde{S}_N^+)$, which tells us that the underlying noncommutative space $\widetilde{S}_N^+$ is a compact quantum semigroup. This semigroup is of quite special type, because $C(\widetilde{S}_N^+)$ has as well a subantipode map, defined by:
$$S(u_{ij})=u_{ji}$$

To be more precise here, this map exists because the transpose of a submagic matrix is submagic too. As for the subantipode axiom satisfied by it, this is as follows, where $m^{(3)}$ is the triple multiplication, and $\Delta^{(2)}$ is the double comultiplication:
$$m^{(3)}(S\otimes id\otimes S)\Delta^{(2)}=S$$

Finally, observe that $\Delta,\varepsilon,S$ restrict to $C(\widetilde{S}_N)$, and correspond there, via Gelfand duality, to the usual multiplication, unit element, and subinversion map of $\widetilde{S}_N$.

\bigskip

As a conclusion to this discussion, the basic properties of the quantum semigroup $\widetilde{S}_N^+$ that we constructed in Definition 6.14 can be summarized as follows:

\begin{proposition}
We have maps as follows,
$$\begin{matrix}
C(\widetilde{S}_N^+)&\to&C(S_N^+)\\
\\
\downarrow&&\downarrow\\
\\
C(\widetilde{S}_N)&\to&C(S_N)
\end{matrix}
\quad \quad \quad:\quad \quad\quad
\begin{matrix}
\widetilde{S}_N^+&\supset&S_N^+\\
\\
\cup&&\cup\\
\\
\widetilde{S}_N&\supset&S_N
\end{matrix}$$
with the bialgebras at left corresponding to the quantum semigroups at right.
\end{proposition}

\begin{proof}
This is clear from the above discussion, and from the well-known fact that projections which sum up to $1$ are pairwise orthogonal.
\end{proof}

As a first example, we have $\widetilde{S}_1^+=\widetilde{S}_1$. At $N=2$ now, recall that the algebra generated by two free projections $p,q$ is isomorphic to the group algebra of $D_\infty=\mathbb Z_2*\mathbb Z_2$. We denote by $\varepsilon:C^*(D_\infty)\to\mathbb C1$ the counit map, given by the following formulae:
$$\varepsilon(1)=1$$
$$\varepsilon(\ldots pqpq\ldots)=0$$

With these conventions, we have the following result:

\begin{proposition}
We have an isomorphism
$$C(\widetilde{S}_2^+)\simeq\left\{(x,y)\in C^*(D_\infty)\oplus C^*(D_\infty)\Big|\varepsilon(x)=\varepsilon(y)\right\}$$
which is given by the formula
$$u=\begin{pmatrix}p\oplus 0&0\oplus r\\0\oplus s&q\oplus 0\end{pmatrix}$$ where $p,q$ and $r,s$ are the standard generators of the two copies of $C^*(D_\infty)$.
\end{proposition}

\begin{proof}
Consider an arbitrary $2\times 2$ matrix formed by projections:
$$u=\begin{pmatrix}P&R\\S&Q\end{pmatrix}$$

This matrix is submagic when the following conditions are satisfied:
$$PR=PS=QR=QS=0$$

But these conditions mean that $X=<P,Q>$ and $Y=<R,S>$ must commute, and must satisfy $xy=0$, for any $x\in X,y\in Y$. Thus, if we denote by $Z$ the universal non-unital algebra generated by two projections, we have an isomorphism as follows:
$$C(\widetilde{S}_2^+)\simeq\mathbb C1\oplus Z\oplus Z$$

Now since $C^*(D_\infty)=\mathbb C1\oplus Z$, we obtain an isomorphism as follows:
$$C(\widetilde{S}_2^+)\simeq\left\{(\lambda+a,\lambda+b)\Big|\lambda\in\mathbb C, a,b\in Z\right\}$$

Thus, we are led to the conclusion in the statement.
\end{proof}

Summarizing, the semigroups of partial permutations $\widetilde{S}_N$ have non-trivial liberations, a bit like the permutation groups $S_N$ used to have non-trivial liberations, and this starting from $N=2$ already. In order to reach now to homogeneous spaces, in the spirit of the partial isometry spaces discussed before, we can use the following simple observation:

\begin{proposition}
Any partial permutation $\sigma:X\simeq Y$ can be factorized as
$$\xymatrix@R=40pt@C=40pt
{X\ar[r]^{\sigma}\ar[d]_\gamma&Y\\\{1,\ldots,k\}\ar[r]_\beta&\{1,\ldots,k\}\ar[u]_\alpha}$$
with $\alpha,\beta,\gamma\in S_k$ being certain non-unique permutations, where $k=\kappa(\sigma)$.
\end{proposition}

\begin{proof}
Since we have $|X|=|Y|=k$, we can pick two bijections, as follows:
$$X\simeq\{1,\ldots,k\}\quad,\quad 
\{1,\ldots,k\}\simeq Y$$

We can complete then these bijections up to permutations $\gamma,\alpha\in S_N$. The remaining permutation $\beta\in S_k$ is then uniquely determined by $\sigma=\alpha\beta\gamma$, as desired.
\end{proof}

With a bit more work, this leads to homogeneous spaces, in the spirit of the partial isometry spaces discussed before. To be more precise, we have the following notion:

\begin{definition}
Associated to any partial permutation, written $\sigma:I\simeq J$ with $I\subset\{1,\ldots,N\}$ and $J\subset\{1,\ldots,M\}$, is the real/complex partial isometry
$$T_\sigma:span\left(e_i\Big|i\in I\right)\to span\left(e_j\Big|j\in J\right)$$
given on the standard basis elements by $T_\sigma(e_i)=e_{\sigma(i)}$.
\end{definition}

We denote by $S_{MN}^L$ the set of partial permutations $\sigma:I\simeq J$ as above, with range $I\subset\{1,\ldots,N\}$ and target $J\subset\{1,\ldots,M\}$, and with:
$$L=|I|=|J|$$

In analogy with the decomposition result $H_N^s=\mathbb Z_s\wr S_N$, we have:

\begin{proposition}
The space of partial permutations signed by elements of $\mathbb Z_s$,
$$H_{MN}^{sL}=\left\{T(e_i)=w_ie_{\sigma(i)}\Big|\sigma\in S_{MN}^L,w_i\in\mathbb Z_s\right\}$$
is isomorphic to the quotient space 
$$(H_M^s\times H_N^s)/(H_L^s\times H_{M-L}^s\times H_{N-L}^s)$$
and with this coming by looking at the stabilizer of $U=(^1_0{\ }^0_0)$.
\end{proposition}

\begin{proof}
This follows by adapting the computations in the proof of Proposition 6.3. Indeed, we have an action map as follows, which is transitive:
$$H_M^s\times H_N^s\to H_{MN}^{sL}\quad,\quad 
(A,B)U=AUB^*$$

Consider now the following point:
$$U=\begin{pmatrix}1&0\\0&0\end{pmatrix}$$

The stabilizer of this point is then the following group:
$$H_L^s\times H_{M-L}^s\times H_{N-L}^s$$

To be more precise, this group is embedded via:
$$(x,a,b)\to\left[\begin{pmatrix}x&0\\0&a\end{pmatrix},\begin{pmatrix}x&0\\0&b\end{pmatrix}\right]$$

But this gives the result.
\end{proof}

In the free case now, the idea is similar, by using inspiration from the construction of the quantum group $H_N^{s+}=\mathbb Z_s\wr_*S_N^+$ in \cite{bb+}. The result here is as follows:

\begin{proposition}
The compact quantum space $H_{MN}^{sL+}$ associated to the algebra
$$C(H_{MN}^{sL+})=C(U_{MN}^{L+})\Big/\left<u_{ij}u_{ij}^*=u_{ij}^*u_{ij}=p_{ij}={\rm projections},u_{ij}^s=p_{ij}\right>$$
has an action map, and is the target of a quotient map, as in Theorem 6.10.
\end{proposition}

\begin{proof}
We must show that if the variables $u_{ij}$ satisfy the relations in the statement, then these relations are satisfied as well for the following variables: 
$$U_{ij}=\sum_{kl}u_{kl}\otimes a_{ki}\otimes b_{lj}^*$$
$$V_{ij}=\sum_{r\leq L}a_{ri}\otimes b_{rj}^*$$

We use the fact that the standard coordinates $a_{ij},b_{ij}$ on the quantum groups $H_M^{s+},H_N^{s+}$ satisfy the following relations, for any $x\neq y$ on the same row or column of $a,b$:
$$xy=xy^*=0$$
 
We obtain, by using these relations:
\begin{eqnarray*}
U_{ij}U_{ij}^*
&=&\sum_{klmn}u_{kl}u_{mn}^*\otimes a_{ki}a_{mi}^*\otimes b_{lj}^*b_{mj}\\
&=&\sum_{kl}u_{kl}u_{kl}^*\otimes a_{ki}a_{ki}^*\otimes b_{lj}^*b_{lj}
\end{eqnarray*}

We have as well the following formula:
\begin{eqnarray*}
V_{ij}V_{ij}^*
&=&\sum_{r,t\leq L}a_{ri}a_{ti}^*\otimes b_{rj}^*b_{tj}\\
&=&\sum_{r\leq L}a_{ri}a_{ri}^*\otimes b_{rj}^*b_{rj}
\end{eqnarray*}

Consider now the following projections:
$$x_{ij}=a_{ij}a_{ij}^*\quad,\quad 
y_{ij}=b_{ij}b_{ij}^*\quad,\quad 
p_{ij}=u_{ij}u_{ij}^*$$

In terms of these projections, we have:
$$U_{ij}U_{ij}^*=\sum_{kl}p_{kl}\otimes x_{ki}\otimes y_{lj}$$
$$V_{ij}V_{ij}^*=\sum_{r\leq L}x_{ri}\otimes y_{rj}$$

By repeating the computation, we conclude that these elements are projections. Also, a similar computation shows that $U_{ij}^*U_{ij},V_{ij}^*V_{ij}$ are given by the same formulae.

Finally, once again by using the relations of type $xy=xy^*=0$, we have:
\begin{eqnarray*}
U_{ij}^s
&=&\sum_{k_rl_r}u_{k_1l_1}\ldots u_{k_sl_s}\otimes a_{k_1i}\ldots a_{k_si}\otimes b_{l_1j}^*\ldots b_{l_sj}^*\\
&=&\sum_{kl}u_{kl}^s\otimes a_{ki}^s\otimes(b_{lj}^*)^s
\end{eqnarray*}

We have as well the following formula:
\begin{eqnarray*}
V_{ij}^s
&=&\sum_{r_l\leq L}a_{r_1i}\ldots a_{r_si}\otimes b_{r_1j}^*\ldots b_{r_sj}^*\\
&=&\sum_{r\leq L}a_{ri}^s\otimes(b_{rj}^*)^s
\end{eqnarray*}

Thus the conditions of type $u_{ij}^s=p_{ij}$ are satisfied as well, and we are done.
\end{proof}

Let us discuss now the general case. We have the following result:

\begin{proposition}
The various spaces $G_{MN}^L$ constructed so far appear by imposing to the standard coordinates of $U_{MN}^{L+}$ the relations
$$\sum_{i_1\ldots i_s}\sum_{j_1\ldots j_s}\delta_\pi(i)\delta_\sigma(j)u_{i_1j_1}^{e_1}\ldots u_{i_sj_s}^{e_s}=L^{|\pi\vee\sigma|}$$
with $s=(e_1,\ldots,e_s)$ ranging over all the colored integers, and with $\pi,\sigma\in D(0,s)$.
\end{proposition}

\begin{proof}
According to the various constructions above, the relations defining $G_{MN}^L$ can be written as follows, with $\sigma$ ranging over a family of generators, with no upper legs, of the corresponding category of partitions $D$:
$$\sum_{j_1\ldots j_s}\delta_\sigma(j)u_{i_1j_1}^{e_1}\ldots u_{i_sj_s}^{e_s}=\delta_\sigma(i)$$

We therefore obtain the relations in the statement, as follows:
\begin{eqnarray*}
\sum_{i_1\ldots i_s}\sum_{j_1\ldots j_s}\delta_\pi(i)\delta_\sigma(j)u_{i_1j_1}^{e_1}\ldots u_{i_sj_s}^{e_s}
&=&\sum_{i_1\ldots i_s}\delta_\pi(i)\sum_{j_1\ldots j_s}\delta_\sigma(j)u_{i_1j_1}^{e_1}\ldots u_{i_sj_s}^{e_s}\\
&=&\sum_{i_1\ldots i_s}\delta_\pi(i)\delta_\sigma(i)\\
&=&L^{|\pi\vee\sigma|}
\end{eqnarray*}

As for the converse, this follows by using the relations in the statement, by keeping $\pi$ fixed, and by making $\sigma$ vary over all the partitions in the category.
\end{proof}

In the general case now, where $G=(G_N)$ is an arbitrary uniform easy quantum group, we can construct spaces $G_{MN}^L$ by using the above relations, and we have:

\begin{theorem}
The spaces $G_{MN}^L\subset U_{MN}^{L+}$ constructed by imposing the relations 
$$\sum_{i_1\ldots i_s}\sum_{j_1\ldots j_s}\delta_\pi(i)\delta_\sigma(j)u_{i_1j_1}^{e_1}\ldots u_{i_sj_s}^{e_s}=L^{|\pi\vee\sigma|}$$
with $\pi,\sigma$ ranging over all the partitions in the associated category, having no upper legs, are subject to an action map/quotient map diagram, as in Theorem 6.10.
\end{theorem}

\begin{proof}
We proceed as in the proof of Proposition 6.9. We must prove that, if the variables $u_{ij}$ satisfy the relations in the statement, then so do the following variables:
$$U_{ij}=\sum_{kl}u_{kl}\otimes a_{ki}\otimes b_{lj}^*$$
$$V_{ij}=\sum_{r\leq L}a_{ri}\otimes b_{rj}^*$$

Regarding the variables $U_{ij}$, the computation here goes as follows:
\begin{eqnarray*}
&&\sum_{i_1\ldots i_s}\sum_{j_1\ldots j_s}\delta_\pi(i)\delta_\sigma(j)U_{i_1j_1}^{e_1}\ldots U_{i_sj_s}^{e_s}\\
&=&\sum_{i_1\ldots i_s}\sum_{j_1\ldots j_s}\sum_{k_1\ldots k_s}\sum_{l_1\ldots l_s}u_{k_1l_1}^{e_1}\ldots u_{k_sl_s}^{e_s}\otimes \delta_\pi(i)\delta_\sigma(j)a_{k_1i_1}^{e_1}\ldots a_{k_si_s}^{e_s}\otimes(b_{l_sj_s}^{e_s}\ldots b_{l_1j_1}^{e_1})^*\\
&=&\sum_{k_1\ldots k_s}\sum_{l_1\ldots l_s}\delta_\pi(k)\delta_\sigma(l)u_{k_1l_1}^{e_1}\ldots u_{k_sl_s}^{e_s}\\
&=&L^{|\pi\vee\sigma|}
\end{eqnarray*}

For the variables $V_{ij}$ the proof is similar, as follows:
\begin{eqnarray*}
&&\sum_{i_1\ldots i_s}\sum_{j_1\ldots j_s}\delta_\pi(i)\delta_\sigma(j)V_{i_1j_1}^{e_1}\ldots V_{i_sj_s}^{e_s}\\
&=&\sum_{i_1\ldots i_s}\sum_{j_1\ldots j_s}\sum_{l_1,\ldots,l_s\leq L}\delta_\pi(i)\delta_\sigma(j)a_{l_1i_1}^{e_1}\ldots a_{l_si_s}^{e_s}\otimes(b_{l_sj_s}^{e_s}\ldots b_{l_1j_1}^{e_1})^*\\
&=&\sum_{l_1,\ldots,l_s\leq L}\delta_\pi(l)\delta_\sigma(l)\\
&=&L^{|\pi\vee\sigma|}
\end{eqnarray*}

Thus we have constructed an action map, and a quotient map, as in Proposition 6.9 above, and the commutation of the diagram in Theorem 6.10 is then trivial.
\end{proof}

\section*{6d. Integration theory}

Let us discuss now the integration over $G_{MN}^L$. We first have:

\begin{definition}
The integration functional of $G_{MN}^L$ is the composition
$$\int_{G_{MN}^L}:C(G_{MN}^L)\to C(G_M\times G_N)\to\mathbb C$$
of the representation $u_{ij}\to\sum_{r\leq L}a_{ri}\otimes b_{rj}^*$ with the Haar functional of $G_M\times G_N$.
\end{definition}

Observe that in the case $L=M=N$ we obtain the integration over $G_N$. Also, at $L=M=1$, or at $L=N=1$, we obtain the integration over the sphere. 

\bigskip

In the general case now, we first have the following result:

\begin{proposition}
The integration functional of $G_{MN}^L$ has the invariance property 
$$\left(\int_{G_{MN}^L}\!\otimes\ id\right)\Phi(x)=\int_{G_{MN}^L}x$$
with respect to the standard coaction map $\Phi(u_{ij})=\sum_{kl}u_{kl}\otimes a_{ki}\otimes b_{lj}^*$.
\end{proposition}

\begin{proof}
We restrict the attention to the orthogonal case, the proof in the unitary case being similar. We must check the following formula:
$$\left(\int_{G_{MN}^L}\!\otimes\ id\right)\Phi(u_{i_1j_1}\ldots u_{i_sj_s})=\int_{G_{MN}^L}u_{i_1j_1}\ldots u_{i_sj_s}$$

Let us compute the left term. This is given by:
\begin{eqnarray*}
X
&=&\left(\int_{G_{MN}^L}\!\otimes\ id\right)\sum_{k_xl_x}u_{k_1l_1}\ldots u_{k_sl_s}\otimes a_{k_1i_1}\ldots a_{k_si_s}\otimes b_{l_1j_1}^*\ldots b_{l_sj_s}^*\\
&=&\sum_{k_xl_x}\sum_{r_x\leq L}a_{k_1i_1}\ldots a_{k_si_s}\otimes b_{l_1j_1}^*\ldots b_{l_sj_s}^*\int_{G_M}a_{r_1k_1}\ldots a_{r_sk_s}\int_{G_N}b_{r_1l_1}^*\ldots b_{r_sl_s}^*\\
&=&\sum_{r_x\leq L}\sum_{k_x}a_{k_1i_1}\ldots a_{k_si_s}\int_{G_M}a_{r_1k_1}\ldots a_{r_sk_s}
\otimes\sum_{l_x}b_{l_1j_1}^*\ldots b_{l_sj_s}^*\int_{G_N}b_{r_1l_1}^*\ldots b_{r_sl_s}^*
\end{eqnarray*}

By using now the invariance property of the Haar functionals of $G_M,G_N$, we obtain:
\begin{eqnarray*}
X
&=&\sum_{r_x\leq L}\left(\int_{G_M}\!\otimes\ id\right)\Delta(a_{r_1i_1}\ldots a_{r_si_s})
\otimes\left(\int_{G_N}\!\otimes\ id\right)\Delta(b_{r_1j_1}^*\ldots b_{r_sj_s}^*)\\
&=&\sum_{r_x\leq L}\int_{G_M}a_{r_1i_1}\ldots a_{r_si_s}\int_{G_N}b_{r_1j_1}^*\ldots b_{r_sj_s}^*\\
&=&\left(\int_{G_M}\otimes\int_{G_N}\right)\sum_{r_x\leq L}a_{r_1i_1}\ldots a_{r_si_s}\otimes b_{r_1j_1}^*\ldots b_{r_sj_s}^*
\end{eqnarray*}

But this gives the formula in the statement, and we are done.
\end{proof}

We will prove now that the above functional is in fact the unique positive unital invariant trace on $C(G_{MN}^L)$. For this purpose, we will need the Weingarten formula:

\index{Weingarten formula}

\begin{theorem}
We have the Weingarten type formula
$$\int_{G_{MN}^L}u_{i_1j_1}\ldots u_{i_sj_s}=\sum_{\pi\sigma\tau\nu}L^{|\pi\vee\tau|}\delta_\sigma(i)\delta_\nu(j)W_{sM}(\pi,\sigma)W_{sN}(\tau,\nu)$$
where the matrices on the right are given by $W_{sM}=G_{sM}^{-1}$, with $G_{sM}(\pi,\sigma)=M^{|\pi\vee\sigma|}$.
\end{theorem}

\begin{proof}
We make use of the quantum group Weingarten formula, for which we refer to \cite{ba8}, \cite{bsp}. By using this formula for $G_M,G_N$, we obtain:
\begin{eqnarray*}
\int_{G_{MN}^L}u_{i_1j_1}\ldots u_{i_sj_s}
&=&\sum_{l_1\ldots l_s\leq L}\int_{G_M}a_{l_1i_1}\ldots a_{l_si_s}\int_{G_N}b_{l_1j_1}^*\ldots b_{l_sj_s}^*\\
&=&\sum_{l_1\ldots l_s\leq L}\sum_{\pi\sigma}\delta_\pi(l)\delta_\sigma(i)W_{sM}(\pi,\sigma)\sum_{\tau\nu}\delta_\tau(l)\delta_\nu(j)W_{sN}(\tau,\nu)\\
&=&\sum_{\pi\sigma\tau\nu}\left(\sum_{l_1\ldots l_s\leq L}\delta_\pi(l)\delta_\tau(l)\right)\delta_\sigma(i)\delta_\nu(j)W_{sM}(\pi,\sigma)W_{sN}(\tau,\nu)
\end{eqnarray*}

The coefficient being $L^{|\pi\vee\tau|}$, we obtain the formula in the statement.
\end{proof}

We can now derive an abstract characterization of the integration, as follows:

\begin{theorem}
The integration of $G_{MN}^L$ is the unique positive unital trace 
$$C(G_{MN}^L)\to\mathbb C$$
which is invariant under the action of the quantum group $G_M\times G_N$.
\end{theorem}

\begin{proof}
We use a standard method, from \cite{bgo}, \cite{bss}, the point being to show that we have the following ergodicity formula: 
$$\left(id\otimes\int_{G_M}\otimes\int_{G_N}\right)\Phi(x)=\int_{G_{MN}^L}x$$

We restrict the attention to the orthogonal case, the proof in the unitary case being similar. We must verify that the following holds:
$$\left(id\otimes\int_{G_M}\otimes\int_{G_N}\right)\Phi(u_{i_1j_1}\ldots u_{i_sj_s})=\int_{G_{MN}^L}u_{i_1j_1}\ldots u_{i_sj_s}$$

By using the Weingarten formula, the left term can be written as follows:
\begin{eqnarray*}
X
&=&\sum_{k_1\ldots k_s}\sum_{l_1\ldots l_s}u_{k_1l_1}\ldots u_{k_sl_s}\int_{G_M}a_{k_1i_1}\ldots a_{k_si_s}\int_{G_N}b_{l_1j_1}^*\ldots b_{l_sj_s}^*\\
&=&\sum_{k_1\ldots k_s}\sum_{l_1\ldots l_s}u_{k_1l_1}\ldots u_{k_sl_s}\sum_{\pi\sigma}\delta_\pi(k)\delta_\sigma(i)W_{sM}(\pi,\sigma)\sum_{\tau\nu}\delta_\tau(l)\delta_\nu(j)W_{sN}(\tau,\nu)\\
&=&\sum_{\pi\sigma\tau\nu}\delta_\sigma(i)\delta_\nu(j)W_{sM}(\pi,\sigma)W_{sN}(\tau,\nu)\sum_{k_1\ldots k_s}\sum_{l_1\ldots l_s}\delta_\pi(k)\delta_\tau(l)u_{k_1l_1}\ldots u_{k_sl_s}
\end{eqnarray*}

By using now the summation formula in Theorem 6.22, we obtain:
$$X=\sum_{\pi\sigma\tau\nu}L^{|\pi\vee\tau|}\delta_\sigma(i)\delta_\nu(j)W_{sM}(\pi,\sigma)W_{sN}(\tau,\nu)$$

Now by comparing with the Weingarten formula for $G_{MN}^L$, this proves our claim. Assume now that $\tau:C(G_{MN}^L)\to\mathbb C$ satisfies the invariance condition. We have:
\begin{eqnarray*}
\tau\left(id\otimes\int_{G_M}\otimes\int_{G_N}\right)\Phi(x)
&=&\left(\tau\otimes\int_{G_M}\otimes\int_{G_N}\right)\Phi(x)\\
&=&\left(\int_{G_M}\otimes\int_{G_N}\right)(\tau\otimes id)\Phi(x)\\
&=&\left(\int_{G_M}\otimes\int_{G_N}\right)(\tau(x)1)\\
&=&\tau(x)
\end{eqnarray*}

On the other hand, according to the formula established above, we have as well:
\begin{eqnarray*}
\tau\left(id\otimes\int_{G_M}\otimes\int_{G_N}\right)\Phi(x)
&=&\tau(tr(x)1)\\
&=&tr(x)
\end{eqnarray*}

Thus we obtain $\tau=tr$, and this finishes the proof.
\end{proof}

As a main application of the above results, we have:

\begin{proposition}
For a sum of coordinates of the following type,
$$\chi_E=\sum_{(ij)\in E}u_{ij}$$
with the coordinates not overlapping on rows and columns, we have
$$\int_{G_{MN}^L}\chi_E^s=\sum_{\pi\sigma\tau\nu}K^{|\pi\vee\tau|}L^{|\sigma\vee\nu|}W_{sM}(\pi,\sigma)W_{sN}(\tau,\nu)$$
where $K=|E|$ is the cardinality of the indexing set.
\end{proposition}

\begin{proof}
With $K=|E|$, we can write $E=\{(\alpha(i),\beta(i))\}$, for certain embeddings:
$$\alpha:\{1,\ldots,K\}\subset\{1,\ldots,M\}$$
$$\beta:\{1,\ldots,K\}\subset\{1,\ldots,N\}$$

In terms of these maps $\alpha,\beta$, the moment in the statement is given by:
$$M_s=\int_{G_{MN}^L}\left(\sum_{i\leq K}u_{\alpha(i)\beta(i)}\right)^s$$

By using the Weingarten formula, we can write this quantity as follows:
\begin{eqnarray*}
&&M_s\\
&=&\int_{G_{MN}^L}\sum_{i_1\ldots i_s\leq K}u_{\alpha(i_1)\beta(i_1)}\ldots u_{\alpha(i_s)\beta(i_s)}\\
&=&\sum_{i_1\ldots i_s\leq K}\sum_{\pi\sigma\tau\nu}L^{|\sigma\vee\nu|}\delta_\pi(\alpha(i_1),\ldots,\alpha(i_s))\delta_\tau(\beta(i_1),\ldots,\beta(i_s))W_{sM}(\pi,\sigma)W_{sN}(\tau,\nu)\\
&=&\sum_{\pi\sigma\tau\nu}\left(\sum_{i_1\ldots i_s\leq K}\delta_\pi(i)\delta_\tau(i)\right)L^{|\sigma\vee\nu|}W_{sM}(\pi,\sigma)W_{sN}(\tau,\nu)
\end{eqnarray*}

But, as explained before, in the proof of Theorem 6.25, the coefficient on the left in the last formula is $C=K^{|\pi\vee\tau|}$. We therefore obtain the formula in the statement.
\end{proof}

At a more concrete level now, we have the following conceptual result, making the link with the Bercovici-Pata bijection \cite{bpa}:

\index{Bercovici-Pata bijection}
\index{non-overlapping coordinates}

\begin{theorem}
In the context of the liberation operations 
$$O_{MN}^L\to O_{MN}^{L+}\quad,\quad 
U_{MN}^L\to U_{MN}^{L+}\quad,\quad 
H_{MN}^{sL}\to H_{MN}^{sL+}$$ 
the laws of the sums of non-overlapping coordinates,
$$\chi_E=\sum_{(ij)\in E}u_{ij}$$
are in Bercovici-Pata bijection, in the 
$$|E|=\kappa N,L=\lambda N,M=\mu N$$
regime and $N\to\infty$ limit.
\end{theorem}

\begin{proof}
We use formulae from \cite{bb+}, \cite{bbc}, \cite{bsp}. According to Proposition 6.27, in terms of $K=|E|$, the moments of the variables in the statement are given by:
$$M_s=\sum_{\pi\sigma\tau\nu}K^{|\pi\vee\tau|}L^{|\sigma\vee\nu|}W_{sM}(\pi,\sigma)W_{sN}(\tau,\nu)$$

We use now two standard facts, from \cite{bbc} and related papers, namely the fact that in the $N\to\infty$ limit the Weingarten matrix $W_{sN}$ is concentrated on the diagonal, and the fact that we have an inequality as follows, with equality precisely when $\pi=\sigma$:
$$|\pi\vee\sigma|\leq\frac{|\pi|+|\sigma|}{2}$$

Indeed, with these two ingredients in hand, we can now look in detail at what happens to our moment $M_s$ in the regime from the statement, namely:
$$K=\kappa N,L=\lambda N,M=\mu N,N\to\infty$$

In this regime, we obtain the following estimate:
\begin{eqnarray*}
M_s
&=&\sum_{\pi\sigma\tau\nu}K^{|\pi\vee\tau|}L^{|\sigma\vee\nu|}W_{sM}(\pi,\sigma)W_{sN}(\tau,\nu)\\
&\simeq&\sum_{\pi\tau}K^{|\pi\vee\tau|}L^{|\pi\vee\tau|}M^{-|\pi|}N^{-|\tau|}\\
&\simeq&\sum_\pi K^{|\pi|}L^{|\pi|}M^{-|\pi|}N^{-|\pi|}\\
&=&\sum_\pi\left(\frac{\kappa\lambda}{\mu}\right)^{|\pi|}
\end{eqnarray*}

In order to interpret this formula, we use general theory from \cite{bb+}, \cite{bbc}, \cite{bsp}:

\medskip

(1) For $G_N=O_N,\bar{O}_N/O_N^+$, the above variables $\chi_E$ follow to be asymptotically Gaussian/semicircular, of parameter $\frac{\kappa\lambda}{\mu}$, and hence in Bercovici-Pata bijection.

\medskip

(2) For $G_N=U_N,\bar{U}_N/U_N^+$ the situation is similar, with $\chi_E$ being asymptotically complex Gaussian/circular, of parameter $\frac{\kappa\lambda}{\mu}$, and in Bercovici-Pata bijection. 

\medskip

(3) Finally, for $G_N=H_N^s/H_N^{s+}$, the variables $\chi_E$ are asymptotically Bessel/free Bessel of parameter $\frac{\kappa\lambda}{\mu}$, and once again in Bercovici-Pata bijection.  
\end{proof}

The convergence in the above result is of course in moments, and we do not know whether some stronger convergence results can be formulated. Nor do we know whether one can use linear combinations of coordinates which are  more general than the sums $\chi_E$ that we consider. These are interesting questions, that we would like to raise here.

\bigskip

Also, there are several possible extensions of the above result, for instance by using twisting operations as well. We refer here to \cite{bb+}, \cite{bbc}, \cite{bss} and related papers.

\section*{6e. Exercises} 

Things have been quite advanced in this chapter, and our exercises will mostly focus on details and examples, in relation with the above. First, we have:

\begin{exercise}
Work out, with full details, the particular cases of the diagram
$$\xymatrix@R=15mm@C=15mm{
O_{MN}^{L+}\ar[r]&U_{MN}^{L+}\\
O_{MN}^L\ar[r]\ar[u]&U_{MN}^L\ar[u]}$$
at $L=M=1$, at $L=N=1$, at $L=M$, and at $L=N$.
\end{exercise} 

This is something that we talked about in the above, but without full details. The problem now is that of doing the complete work here, with full details.

\begin{exercise}
Work out, with full details, the particular cases of the spaces
$$G_{MN}^L\subset G_{MN}^{L+}$$
at $L=M$ and at $L=N$, for all the basic easy groups.
\end{exercise} 

As before, we should get row spaces here, and all this needs a complete proof.

\begin{exercise}
Develop a theory of semigroups and quantum semigroups
$$\xymatrix@R=15mm@C=15mm{
\widetilde{O}_N^+\ar[r]&\widetilde{U}_N^+\\
\widetilde{O}_N\ar[r]\ar[u]&\widetilde{U}_N\ar[u]}$$
and then a theory of semigroups and quantum semigroups
$$\xymatrix@R=15mm@C=15mm{
\widetilde{H}_N^+\ar[r]&\widetilde{K}_N^+\\
\widetilde{H}_N\ar[r]\ar[u]&\widetilde{K}_N\ar[u]}$$
and explain the relation with the spaces constructed in this chapter.
\end{exercise} 

Things are actually quite tricky here, in relation with the semigroup structure in the free case. Thus, the exercise asks to develop what can indeed be developed.

\chapter{Affine spaces}

\section*{7a. Quotient spaces}

In this chapter we eventually discuss some abstract aspects, regarding the homogeneous spaces, after about 150 pages of dealing with spheres $S$, and other examples. The reasons for this long delay come from the fact that the theory is quite tricky in the free setting, and so in the quantum setting in general. Any basic attempt of developing a nice, gentle theory in analogy with what is known about the classical homogeneous spaces fails, due to a number of subtle algebraic and analytic reasons, that you can only learn about after studying some examples. Which examples were duly studied in the preceding 150 pages, so we can now go ahead with abstractions, following \cite{ba5}, \cite{ba6}, \cite{bgo}, \cite{bss}.

\bigskip

You might of course smell some controversy in all this, and you are certainly right, because, no surprise, many people have tried, and this since the early 90s, to develop nice and gentle theories of quantum homogeneous spaces. However, from a modern perspective, the findings obtained in this way are rather no-go results. We refer to the papers \cite{ba5}, \cite{ba6}, \cite{bgo}, \cite{bss}, all written in the 10s, for a discussion here, and for references.

\bigskip

Finally, and again talking controversies, following our discussion from the beginning of chapter 5, where mathematician, engineer and cat were debating about noncommutative geometry, we will be obsessed as usual by computing the Haar integration $tr:C(X)\to\mathbb C$ on our homogeneous spaces $X$, and be rather weak on other geometric aspects.

\bigskip

Let us begin with some generalities regarding the quotient spaces, and more general homogeneous spaces. Regarding the quotients, we have the following construction:

\index{quotient space}
\index{homogeneous space}

\begin{proposition}
Given a quantum subgroup $H\subset G$, with associated quotient map $\rho:C(G)\to C(H)$, if we define the quotient space $X=G/H$ by setting
$$C(X)=\left\{f\in C(G)\Big|(\rho\otimes id)\Delta f=1\otimes f\right\}$$
then we have a coaction map as follows,
$$\Phi:C(X)\to C(X)\otimes C(G)$$
obtained as the restriction of the comultiplication of $C(G)$. In the classical case, we obtain in this way the usual quotient space $X=G/H$.
\end{proposition}

\begin{proof}
Observe that the linear subspace $C(X)\subset C(G)$ defined in the statement is indeed a subalgebra, because it is defined via a relation of type $\varphi(f)=\psi(f)$, with both $\varphi,\psi$ being morphisms of algebras. Observe also that in the classical case we obtain the algebra of continuous functions on the quotient space $X=G/H$, because:
\begin{eqnarray*}
(\rho\otimes id)\Delta f=1\otimes f
&\iff&(\rho\otimes id)\Delta f(h,g)=(1\otimes f)(h,g),\forall h\in H,\forall g\in G\\
&\iff&f(hg)=f(g),\forall h\in H,\forall g\in G\\
&\iff&f(hg)=f(kg),\forall h,k\in H,\forall g\in G
\end{eqnarray*}

Regarding now the construction of $\Phi$, observe that for $f\in C(X)$ we have: 
\begin{eqnarray*}
(\rho\otimes id\otimes id)(\Delta\otimes id)\Delta f
&=&(\rho\otimes id\otimes id)(id\otimes\Delta)\Delta f\\
&=&(id\otimes\Delta)(\rho\otimes id)\Delta f\\
&=&(id\otimes\Delta)(1\otimes f)\\
&=&1\otimes\Delta f
\end{eqnarray*}

Thus the condition $f\in C(X)$ implies $\Delta f\in C(X)\otimes C(G)$, and this gives the existence of $\Phi$. Finally, the other assertions are all clear.
\end{proof}

As an illustration, following \cite{bss}, in the group dual case we have:

\begin{proposition}
Assume that $G=\widehat{\Gamma}$ is a discrete group dual.
\begin{enumerate}
\item The quantum subgroups of $G$ are $H=\widehat{\Lambda}$, with $\Gamma\to\Lambda$ being a quotient group.

\item For such a quantum subgroup $\widehat{\Lambda}\subset\widehat{\Gamma}$, we have $\widehat{\Gamma}/\widehat{\Lambda}=\widehat{\Theta}$, where $\Theta=\ker(\Gamma\to\Lambda)$.
\end{enumerate}
\end{proposition}

\begin{proof}
This is well-known, the idea being as follows:

\medskip

(1) In one sense, this is clear. Conversely, since the algebra $C(G)=C^*(\Gamma)$ is cocommutative, so are all its quotients, and this gives the result.

\medskip

(2) Consider a quotient map $r:\Gamma\to\Lambda$, and denote by $\rho:C^*(\Gamma)\to C^*(\Lambda)$ its extension. Consider a group algebra element, written as follows:
$$f=\sum_{g\in\Gamma}\lambda_g\cdot g\in C^*(\Gamma)$$

We have then the following computation:
\begin{eqnarray*}
f\in C(\widehat{\Gamma}/\widehat{\Lambda})
&\iff&(\rho\otimes id)\Delta(f)=1\otimes f\\
&\iff&\sum_{g\in\Gamma}\lambda_g\cdot r(g)\otimes g=\sum_{g\in\Gamma}\lambda_g\cdot 1\otimes g\\
&\iff&\lambda_g\cdot r(g)=\lambda_g\cdot 1,\forall g\in\Gamma\\
&\iff&supp(f)\subset\ker(r)
\end{eqnarray*}

But this means that we have $\widehat{\Gamma}/\widehat{\Lambda}=\widehat{\Theta}$, with $\Theta=\ker(\Gamma\to\Lambda)$, as claimed.
\end{proof}

Given two compact quantum spaces $X,Y$, we say that $X$ is a quotient space of $Y$ when we have an embedding of $C^*$-algebras $\alpha:C(X)\subset C(Y)$. We have:

\begin{definition}
We call a quotient space $G\to X$ homogeneous when
$$\Delta(C(X))\subset C(X)\otimes C(G)$$
where $\Delta:C(G)\to C(G)\otimes C(G)$ is the comultiplication map.
\end{definition}

In other words, an homogeneous quotient space $G\to X$ is a quantum space coming from a subalgebra $C(X)\subset C(G)$, which is stable under the comultiplication. The relation with the quotient spaces from Proposition 7.1 is as follows:

\begin{theorem}
The following results hold:
\begin{enumerate}
\item The quotient spaces $X=G/H$ are homogeneous.

\item In the classical case, any homogeneous space is of type $G/H$.

\item In general, there are homogeneous spaces which are not of type $G/H$.
\end{enumerate}
\end{theorem}

\begin{proof}
Once again these results are well-known, the proof being as follows:

\medskip

(1) This is clear indeed from Proposition 7.1.

\medskip

(2) Consider a quotient map $p:G\to X$. The invariance condition in the statement tells us that we must have an action $G\curvearrowright X$, given by:
$$g(p(g'))=p(gg')$$

Thus, we have the following implication:
$$p(g')=p(g'')\implies p(gg')=p(gg''),\ \forall g\in G$$

Now observe that the following subset $H\subset G$ is a subgroup:
$$H=\left\{g\in G\Big|p(g)=p(1)\right\}$$

Indeed, $g,h\in H$ implies that we have:
$$p(gh)=p(g)=p(1)$$

Thus we have $gh\in H$, and the other axioms are satisfied as well. Our claim now is that we have an identification $X=G/H$, obtained as follows:
$$p(g)\to Hg$$

Indeed, the map $p(g)\to Hg$ is well-defined and bijective, because $p(g)=p(g')$ is equivalent to $p(g^{-1}g')=p(1)$, and so to $Hg=Hg'$, as desired. 

\medskip

(3) Given a discrete group $\Gamma$ and an arbitrary subgroup $\Theta\subset\Gamma$, the quotient space $\widehat{\Gamma}\to\widehat{\Theta}$ is homogeneous. Now by using Proposition 7.2, we can see that if $\Theta\subset\Gamma$ is not normal, the quotient space $\widehat{\Gamma}\to\widehat{\Theta}$ is not of the form $G/H$.
\end{proof}

With the above formalism in hand, let us try now to understand the general properties of the homogeneous spaces $G\to X$, in the sense of Theorem 7.4. We have:

\begin{proposition}
Assume that a quotient space $G\to X$ is homogeneous.
\begin{enumerate}
\item We have a coaction map as follows, obtained as restriction of $\Delta$:
$$\Phi:C(X)\to C(X)\otimes C(G)$$

\item We have the following implication:
$$\Phi(f)=f\otimes 1\implies f\in\mathbb C1$$

\item We have as well the following formula:
$$\left(id\otimes\int_G\right)\Phi f=\int_Gf$$

\item The restriction of $\int_G$ is the unique unital form satisfying $(\tau\otimes id)\Phi=\tau(.)1$.
\end{enumerate}
\end{proposition}

\begin{proof}
These results are all elementary, the proof being as follows:

\medskip

(1) This is clear from definitions, because $\Delta$ itself is a coaction.

\medskip

(2) Assume that $f\in C(G)$ satisfies $\Delta(f)=f\otimes 1$. By applying the counit we obtain:
$$(\varepsilon\otimes id)\Delta f=(\varepsilon\otimes id)(f\otimes 1)$$

We conclude from this that we have $f=\varepsilon(f)1$, as desired.

\medskip

(3) The formula in the statement, $(id\otimes\int_G)\Phi f=\int_Gf$, follows indeed from the left invariance property of the Haar functional of $C(G)$, namely:
$$\left(id\otimes\int_G\right)\Delta f=\int_Gf$$

(4) We use here the right invariance of the Haar functional of $C(G)$, namely:
$$\left(\int_G\otimes id\right)\Delta f=\int_Gf$$

Indeed, we obtain from this that $tr=(\int_G)_{|C(X)}$ is $G$-invariant, in the sense that:
$$(tr\otimes id)\Phi f=tr(f)1$$

Conversely, assuming that $\tau:C(X)\to\mathbb C$ satisfies $(\tau\otimes id)\Phi f=\tau(f)1$, we have:
\begin{eqnarray*}
\left(\tau\otimes\int_G\right)\Phi(f)
&=&\int_G(\tau\otimes id)\Phi(f)\\
&=&\int_G(\tau(f)1)\\
&=&\tau(f)
\end{eqnarray*}

On the other hand, we can compute the same quantity as follows:
\begin{eqnarray*}
\left(\tau\otimes\int_G\right)\Phi(f)
&=&\tau\left(id\otimes\int_G\right)\Phi(f)\\
&=&\tau(tr(f)1)\\
&=&tr(f)
\end{eqnarray*}

Thus we have $\tau(f)=tr(f)$ for any $f\in C(X)$, and this finishes the proof.
\end{proof}

Summarizing, we have a notion of noncommutative homogeneous space, which perfectly covers the classical case. In general, however, the group dual case shows that our formalism is more general than that of the quotient spaces $G/H$.

\section*{7b. Extended spaces}

We discuss now an extra issue, of analytic nature. The point indeed is that for one of the most basic examples of actions, namely $O_N^+\curvearrowright S^{N-1}_{\mathbb R,+}$, the associated morphism $\alpha:C(X)\to C(G)$ is not injective. The same is true for other basic actions, in the free setting. In order to include such examples, we must relax our axioms:

\index{homogeneous space}
\index{extended homogeneous space}

\begin{definition}
An extended homogeneous space over a compact quantum group $G$ consists of a morphism of $C^*$-algebras, and a coaction map, as follows,
$$\alpha:C(X)\to C(G)$$
$$\Phi:C(X)\to C(X)\otimes C(G)$$
such that the following diagram commutes
$$\xymatrix@R=16mm@C=20mm{
C(X)\ar[r]^\Phi\ar[d]_\alpha&C(X)\otimes C(G)\ar[d]^{\alpha\otimes id}\\
C(G)\ar[r]^\Delta&C(G)\otimes C(G)
}$$
and such that the following diagram commutes as well
$$\xymatrix@R=16mm@C=20mm{
C(X)\ar[r]^\Phi\ar[d]_\alpha&C(X)\otimes C(G)\ar[d]^{id\otimes\int}\\
C(G)\ar[r]^{\int(.)1}&C(X)
}$$
where $\int$ is the Haar integration over $G$. We write then $G\to X$.
\end{definition}

As a first observation, when the morphism $\alpha$ is injective we obtain an homogeneous space in the previous sense. The examples with $\alpha$ not injective, which motivate the above formalism, include the standard action $O_N^+\curvearrowright S^{N-1}_{\mathbb R,+}$, and the standard action $U_N^+\curvearrowright S^{N-1}_{\mathbb C,+}$. Following \cite{ba6}, here are a few general remarks on the above axioms:

\index{ergodicity}

\begin{proposition}
Assume that we have morphisms of $C^*$-algebras 
$$\alpha:C(X)\to C(G)$$
$$\Phi:C(X)\to C(X)\otimes C(G)$$ 
satisfying the coassociativity condition $(\alpha\otimes id)\Phi=\Delta\alpha$.
\begin{enumerate}
\item If $\alpha$ is injective on a dense $*$-subalgebra $A\subset C(X)$, and $\Phi(A)\subset A\otimes C(G)$, then $\Phi$ is automatically a coaction map, and is unique.

\item The ergodicity type condition $(id\otimes\int)\Phi=\int\alpha(.)1$ is equivalent to the existence of a linear form $\lambda:C(X)\to\mathbb C$ such that $(id\otimes\int)\Phi=\lambda(.)1$.
\end{enumerate}
\end{proposition}

\begin{proof}
This is something elementary, the idea being as follows:

\medskip

(1) Assuming that we have a dense $*$-subalgebra $A\subset C(X)$ as in the statement, satisying $\Phi(A)\subset A\otimes C(G)$, the restriction $\Phi_{|A}$ is given by:
$$\Phi_{|A}=(\alpha_{|A}\otimes id)^{-1}\Delta\alpha_{|A}$$

This restriction and is therefore coassociative, and unique. By continuity, the morphism $\Phi$ itself follows to be coassociative and unique, as desired.

\medskip

(2) Assuming $(id\otimes\int)\Phi=\lambda(.)1$, we have:
$$\left(\alpha\otimes\int\right)\Phi=\lambda(.)1$$

On the other hand, we have as well the following formula:
$$\left(\alpha\otimes\int\right)\Phi=\left(id\otimes\int\right)\Delta\alpha=\int\alpha(.)1$$

Thus we obtain $\lambda=\int\alpha$, as claimed.
\end{proof}

Given an extended homogeneous space $G\to X$ in our sense, with associated map $\alpha:C(X)\to C(G)$, we can consider the image of this latter map:
$$\alpha:C(X)\to C(Y)\subset C(G)$$

Equivalently, at the level of the associated noncommutative spaces, we can factorize the corresponding quotient map $G\to Y\subset X$. With these conventions, we have:

\begin{proposition}
Consider an extended homogeneous space $G\to X$.
\begin{enumerate}
\item $\Phi(f)=f\otimes 1\implies f\in\mathbb C1$.

\item $tr=\int\alpha$ is the unique unital $G$-invariant form on $C(X)$.

\item The image space obtained by factorizing, $G\to Y$, is homogeneous.
\end{enumerate}
\end{proposition}

\begin{proof}
We have several assertions to be proved, the idea being as follows:

\medskip

(1) This follows indeed from $(id\otimes\int)\Phi(f)=\int\alpha(f)1$, which gives $f=\int\alpha(f)1$.

\medskip

(2) The fact that $tr=\int\alpha$ is indeed $G$-invariant can be checked as follows:
\begin{eqnarray*}
(tr\otimes id)\Phi f
&=&(\smallint\alpha\otimes id)\Phi f\\
&=&(\smallint\otimes id)\Delta\alpha f\\
&=&\smallint\alpha(f)1\\
&=&tr(f)1
\end{eqnarray*}

As for the uniqueness assertion, this follows as before.

\medskip

(3) The condition $(\alpha\otimes id)\Phi=\Delta\alpha$, together with the fact that $i$ is injective, allows us to factorize $\Delta$ into a morphism $\Psi$, as follows:
$$\xymatrix@R=12mm@C=30mm{
C(X)\ar[r]^\Phi\ar[d]_\alpha&C(X)\otimes C(G)\ar[d]^{\alpha\otimes id}\\
C(Y)\ar@.[r]^\Psi\ar[d]_i&C(Y)\otimes C(G)\ar[d]^{i\otimes id}\\
C(G)\ar[r]^\Delta&C(G)\otimes C(G)
}$$

Thus the image space $G\to Y$ is indeed homogeneous, and we are done.
\end{proof}

Finally, still following \cite{ba6}, we have the following result:

\index{GNS construction}

\begin{theorem}
Let $G\to X$ be an extended homogeneous space, and construct quotients $X\to X'$, $G\to G'$ by performing the GNS construction with respect to $\int\alpha,\int$. Then $\alpha$ factorizes into an inclusion $\alpha':C(X')\to C(G')$, and we have an homogeneous space.
\end{theorem}

\begin{proof}
We factorize $G\to Y\subset X$ as above. By performing the GNS construction with respect to $\int i\alpha,\int i,\int$, we obtain a diagram as follows:
$$\xymatrix@R=12mm@C=30mm{
C(X)\ar[r]^p\ar[d]_\alpha&C(X')\ar[d]^{\alpha'}\ar[dr]^{tr'}\\
C(Y)\ar[r]^q\ar[d]_i&C(Y')\ar[d]^{i'}&\mathbb C\\
C(G)\ar[r]^r&C(G')\ar[ur]_{\int'}
}$$

Indeed, with $tr=\int\alpha$, the GNS quotient maps $p,q,r$ are defined respectively by:
\begin{eqnarray*}
\ker p&=&\left\{f\in C(X)\Big|tr(f^*f)=0\right\}\\
\ker q&=&\left\{f\in C(Y)\Big|\smallint(f^*f)=0\right\}\\
\ker r&=&\left\{f\in C(G)\Big|\smallint(f^*f)=0\right\}
\end{eqnarray*}

Next, we can define factorizations $i',\alpha'$ as above. Observe that $i'$ is injective, and that $\alpha'$ is surjective. Our claim now is that $\alpha'$ is injective as well. Indeed:
\begin{eqnarray*}
\alpha'p(f)=0
&\implies&q\alpha(f)=0\\
&\implies&\int\alpha(f^*f)=0\\
&\implies&tr(f^*f)=0\\
&\implies&p(f)=0
\end{eqnarray*}

We conclude that we have $X'=Y'$, and this gives the result.
\end{proof}

Summarizing, the basic homogeneous space theory from the classical case extends to the quantum group setting, with a few twists, both of algebraic and analytic nature. All the above was of course quite brief, and designed to best capture what happens in free geometry, but at the level of the general things that can be said about quantum homogeneous spaces, there is of course much more. We will be back to this.

\section*{7c. Affine spaces}

We discuss now an abstract extension of the constructions of manifolds that we have so far. The idea will be that of looking at certain classes of algebraic manifolds $X\subset S^{N-1}_{\mathbb C,+}$, which are homogeneous spaces, of a certain special type. Following \cite{ba6}, we have:

\index{homogeneous space}
\index{affine homogeneous space}
\index{ergodicity}

\begin{definition}
An affine homogeneous space over a closed subgroup $G\subset U_N^+$ is a closed subset $X\subset S^{N-1}_{\mathbb C,+}$, such that there exists an index set $I\subset\{1,\ldots,N\}$ such that
$$\alpha(x_i)=\frac{1}{\sqrt{|I|}}\sum_{j\in I}u_{ji}\quad,\quad 
\Phi(x_i)=\sum_jx_j\otimes u_{ji}$$
define morphisms of $C^*$-algebras, satisfying the following condition,
$$\left(id\otimes\int_G\right)\Phi=\int_G\alpha(.)1$$ 
called ergodicity condition for the action.
\end{definition}

Let us mention right away that this definition is something quite tricky, based on the explicit examples of homogeneous spaces that we have in mind, rather than on whatever abstract considerations, and that will take us some time to understand. 

\bigskip

To start with, as a basic example, $O_N^+\to S^{N-1}_{\mathbb R,+}$ is indeed affine in our sense, with $I=\{1\}$. The same goes for $U_N^+\to S^{N-1}_{\mathbb C,+}$, which is affine as well, also with $I=\{1\}$. 

\bigskip

Observe that the $1/\sqrt{|I|}$ constant appearing above is the correct one, because: 
\begin{eqnarray*}
\sum_i\left(\sum_{j\in I}u_{ji}\right)\left(\sum_{k\in I}u_{ki}\right)^*
&=&\sum_i\sum_{j,k\in I}u_{ji}u_{ki}^*\\
&=&\sum_{j,k\in I}(uu^*)_{jk}\\
&=&|I|
\end{eqnarray*}

As a first general result about such spaces, following \cite{ba6}, we have:

\begin{proposition}
Consider an affine homogeneous space $X$, as above.
\begin{enumerate}
\item The coaction condition $(\Phi\otimes id)\Phi=(id\otimes\Delta)\Phi$ is satisfied.

\item We have as well the formula $(\alpha\otimes id)\Phi=\Delta\alpha$.
\end{enumerate}
\end{proposition}

\begin{proof}
The coaction condition is clear. For the second formula, we first have:
\begin{eqnarray*}
(\alpha\otimes id)\Phi(x_i)
&=&\sum_k\alpha(x_k)\otimes u_{ki}\\
&=&\frac{1}{\sqrt{|I|}}\sum_k\sum_{j\in I}u_{jk}\otimes u_{ki}
\end{eqnarray*}

On the other hand, we have as well the following computation:
\begin{eqnarray*}
\Delta\alpha(x_i)
&=&\frac{1}{\sqrt{|I|}}\sum_{j\in I}\Delta(u_{ji})\\
&=&\frac{1}{\sqrt{|I|}}\sum_{j\in I}\sum_ku_{jk}\otimes u_{ki}
\end{eqnarray*}

Thus, by linearity, multiplicativity and continuity, we obtain the result.
\end{proof}

Summarizing, the terminology in Definition 7.10 is justified, in the sense that what we have there are indeed certain homogeneous spaces, of very special, ``affine" type. As a second result regarding such spaces, which closes the discussion in the case where $\alpha$ is injective, which is something that happens in many cases, we have:

\index{minimal homogeneous space}

\begin{theorem}
When $\alpha$ is injective we must have $X=X_{G,I}^{min}$, where:
$$C(X_{G,I}^{min})=\left<\frac{1}{\sqrt{|I|}}\sum_{j\in I}u_{ji}\Big|i=1,\ldots,N\right>\subset C(G)$$
Moreover, $X_{G,I}^{min}$ is affine homogeneous, for any $G\subset U_N^+$, and any $I\subset\{1,\ldots,N\}$.
\end{theorem}

\begin{proof}
The first assertion is clear from definitions. Regarding now the second assertion, consider the variables in the statement:
$$X_i=\frac{1}{\sqrt{|I|}}\sum_{j\in I}u_{ji}\in C(G)$$

In order to prove that we have $X_{G,I}^{min}\subset S^{N-1}_{\mathbb C,+}$, observe first that we have:
\begin{eqnarray*}
\sum_iX_iX_i^*
&=&\frac{1}{|I|}\sum_i\sum_{j,k\in I}u_{ji}u_{ki}^*\\
&=&\frac{1}{|I|}\sum_{j,k\in I}(uu^*)_{jk}\\
&=&1
\end{eqnarray*}

On the other hand, we have as well the following computation:
\begin{eqnarray*}
\sum_iX_i^*X_i
&=&\frac{1}{|I|}\sum_i\sum_{j,k\in I}u_{ji}^*u_{ki}\\
&=&\frac{1}{|I|}\sum_{j,k\in I}(\bar{u}u^t)_{jk}\\
&=&1
\end{eqnarray*}

Thus $X_{G,I}^{min}\subset S^{N-1}_{\mathbb C,+}$. Finally, observe that we have:
\begin{eqnarray*}
\Delta(X_i)
&=&\frac{1}{\sqrt{|I|}}\sum_{j\in I}\sum_ku_{jk}\otimes u_{ki}\\
&=&\sum_kX_k\otimes u_{ki}
\end{eqnarray*}

Thus we have indeed a coaction map, given by $\Phi=\Delta$. As for the ergodicity condition, namely $(id\otimes\int_G)\Delta=\int_G(.)1$, this holds as well, by definition of the integration functional $\int_G$. Thus, our axioms for affine homogeneous spaces are indeed satisfied.
\end{proof}

Our purpose now will be to show that the affine homogeneous spaces appear as follows, a bit in the same way as the discrete group algebras:
$$X_{G,I}^{min}\subset X\subset X_{G,I}^{max}$$

We make the standard convention that all the tensor exponents $k$ are ``colored integers'', that is, $k=e_1\ldots e_k$ with $e_i\in\{\circ,\bullet\}$, with $\circ$ corresponding to the usual variables, and with $\bullet$ corresponding to their adjoints. With this convention, we have:

\index{ergodicity}

\begin{proposition}
The ergodicity condition, namely
$$\left(id\otimes\int_G\right)\Phi=\int_G\alpha(.)1$$
is equivalent to the condition
$$(Px^{\otimes k})_{i_1\ldots i_k}=\frac{1}{\sqrt{|I|^k}}\sum_{j_1\ldots j_k\in I}P_{i_1\ldots i_k,j_1\ldots j_k}\quad,\quad\forall k,\forall i_1,\ldots,i_k$$
where $P$ is the matrix formed by the Peter-Weyl integrals of exponent $k$,
$$P_{i_1\ldots i_k,j_1\ldots j_k}=\int_Gu_{j_1i_1}^{e_1}\ldots u_{j_ki_k}^{e_k}$$
and where $(x^{\otimes k})_{i_1\ldots i_k}=x_{i_1}^{e_1}\ldots x_{i_k}^{e_k}$.
\end{proposition}

\begin{proof}
We have the following computation:
\begin{eqnarray*}
\left(id\otimes\int_G\right)\Phi(x_{i_1}^{e_1}\ldots x_{i_k}^{e_k})
&=&\sum_{j_1\ldots j_k}x_{j_1}^{e_1}\ldots x_{j_k}^{e_k}\int_Gu_{j_1i_1}^{e_1}\ldots u_{j_ki_k}^{e_k}\\
&=&\sum_{j_1\ldots j_k}P_{i_1\ldots i_k,j_1\ldots j_k}(x^{\otimes k})_{j_1\ldots j_k}\\
&=&(Px^{\otimes k})_{i_1\ldots i_k}
\end{eqnarray*}

On the other hand, we have as well the following computation:
\begin{eqnarray*}
\int_G\alpha(x_{i_1}^{e_1}\ldots x_{i_k}^{e_k})
&=&\frac{1}{\sqrt{|I|^k}}\sum_{j_1\ldots j_k\in I}\int_Gu_{j_1i_1}^{e_1}\ldots u_{j_ki_k}^{e_k}\\
&=&\frac{1}{\sqrt{|I|^k}}\sum_{j_1\ldots j_k\in I}P_{i_1\ldots i_k,j_1\ldots j_k}
\end{eqnarray*}

But this gives the formula in the statement, and we are done.
\end{proof}

As a consequence, we have the following result:

\index{maximal homogeneous space}

\begin{theorem}
We must have $X\subset X_{G,I}^{max}$, as subsets of $S^{N-1}_{\mathbb C,+}$, where:
$$C(X_{G,I}^{max})=C(S^{N-1}_{\mathbb C,+})\Big/\left<(Px^{\otimes k})_{i_1\ldots i_k}=\frac{1}{\sqrt{|I|^k}}\sum_{j_1\ldots j_k\in I}P_{i_1\ldots i_k,j_1\ldots j_k}\big|\forall k,\forall i_1,\ldots i_k\right>$$
Moreover, $X_{G,I}^{max}$ is affine homogeneous, for any $G\subset U_N^+$, and any $I\subset\{1,\ldots,N\}$.
\end{theorem}

\begin{proof}
Let us first prove that we have an action $G\curvearrowright X_{G,I}^{max}$. We must show here that the variables $X_i=\sum_jx_j\otimes u_{ji}$ satisfy the defining relations for $X_{G,I}^{max}$. We have:
\begin{eqnarray*}
(PX^{\otimes k})_{i_1\ldots i_k}
&=&\sum_{l_1\ldots l_k}P_{i_1\ldots i_k,l_1\ldots l_k}(X^{\otimes k})_{l_1\ldots l_k}\\
&=&\sum_{l_1\ldots l_k}P_{i_1\ldots i_k,l_1\ldots l_k}\sum_{j_1\ldots j_k}x_{j_1}^{e_1}\ldots x_{j_k}^{e_k}\otimes u_{j_1l_1}^{e_1}\ldots u_{j_kl_k}^{e_k}\\
&=&\sum_{j_1\ldots j_k}x_{j_1}^{e_1}\ldots x_{j_k}^{e_k}\otimes(u^{\otimes k}P^t)_{j_1\ldots j_k,i_1\ldots i_k}
\end{eqnarray*}

Since by Peter-Weyl the transpose of $P_{i_1\ldots i_k,j_1\ldots j_k}=\int_Gu_{j_1i_1}^{e_1}\ldots u_{j_ki_k}^{e_k}$ is the orthogonal projection onto $Fix(u^{\otimes k})$, we have $u^{\otimes k}P^t=P^t$. We therefore obtain:
\begin{eqnarray*}
(PX^{\otimes k})_{i_1\ldots i_k}
&=&\sum_{j_1\ldots j_k}P_{i_1\ldots i_k,j_1\ldots j_k}x_{j_1}^{e_1}\ldots x_{j_k}^{e_k}\\
&=&(Px^{\otimes k})_{i_1\ldots i_k}\\
&=&\frac{1}{\sqrt{|I|^k}}\sum_{j_1\ldots j_k\in I}P_{i_1\ldots i_k,j_1\ldots j_k}
\end{eqnarray*}

Thus we have an action $G\curvearrowright X_{G,I}^{max}$, and since this action is ergodic by Proposition 7.13, we have an affine homogeneous space, as claimed. 
\end{proof}

We can now merge the results that we have, and we obtain, following \cite{ba6}:

\begin{theorem}
Given a closed quantum subgroup $G\subset U_N^+$, and a set $I\subset\{1,\ldots,N\}$, if we consider the following $C^*$-subalgebra and the following quotient $C^*$-algebra,
\begin{eqnarray*}
C(X_{G,I}^{min})&=&\left<\frac{1}{\sqrt{|I|}}\sum_{j\in I}u_{ji}\Big|i=1,\ldots,N\right>\subset C(G)\\
C(X_{G,I}^{max})&=&C(S^{N-1}_{\mathbb C,+})\Big/\left<(Px^{\otimes k})_{i_1\ldots i_k}=\frac{1}{\sqrt{|I|^k}}\sum_{j_1\ldots j_k\in I}P_{i_1\ldots i_k,j_1\ldots j_k}\Big|\forall k,\forall i_1,\ldots i_k\right>
\end{eqnarray*}
then we have maps as follows,
$$G\to X_{G,I}^{min}\subset X_{G,I}^{max}\subset S^{N-1}_{\mathbb C,+}$$
the space $G\to X_{G,I}^{max}$ is affine homogeneous, and any affine homogeneous space $G\to X$ appears as an intermediate space $X_{G,I}^{min}\subset X\subset X_{G,I}^{max}$.
\end{theorem}

\begin{proof}
This follows indeed from the various results that we have, namely Theorem 7.12 and Theorem 7.14, regarding the minimal and maximal constructions.
\end{proof}

Summarizing, the situation with our affine homogeneous spaces is, from a point of view of abstract functional analysis, a bit similar to that of the full and reduced group algebras, with intermediate objects between them. We will be back to this, later on.

\bigskip

At the level of the general theory, based on Definition 7.10, we will need one more general result from \cite{ba6}, namely an extension of the Weingarten integration formula \cite{bbc}, \cite{csn}, \cite{wei}, to the affine homogeneous space setting, as follows:

\index{Weingarten formula}

\begin{theorem}
Assuming that $G\to X$ is an affine homogeneous space, with index set $I\subset\{1,\ldots,N\}$, the Haar integration functional $\int_X=\int_G\alpha$ is given by
$$\int_Xx_{i_1}^{e_1}\ldots x_{i_k}^{e_k}=\sum_{\pi,\sigma\in D}K_I(\pi)\overline{(\xi_\sigma)}_{i_1\ldots i_k}W_{kN}(\pi,\sigma)$$
where $\{\xi_\pi|\pi\in D\}$ is a basis of $Fix(u^{\otimes k})$, the quantities $K_I(\pi)$ are given by
$$K_I(\pi)=\frac{1}{\sqrt{|I|^k}}\sum_{j_1\ldots j_k\in I}(\xi_\pi)_{j_1\ldots j_k}$$
and $W_{kN}=G_{kN}^{-1}$ with $G_{kN}(\pi,\sigma)=<\xi_\pi,\xi_\sigma>$ is the associated Weingarten matrix.
\end{theorem}

\begin{proof}
By using the Weingarten formula for the quantum group $G$, in its abstract form, coming from Peter-Weyl theory, as discussed in chapter 2, we have:
\begin{eqnarray*}
\int_Xx_{i_1}^{e_1}\ldots x_{i_k}^{e_k}
&=&\frac{1}{\sqrt{|I|^k}}\sum_{j_1\ldots j_k\in I}\int_Gu_{j_1i_1}^{e_1}\ldots u_{j_ki_k}^{e_k}\\
&=&\frac{1}{\sqrt{|I|^k}}\sum_{j_1\ldots j_k\in I}\sum_{\pi,\sigma\in D}(\xi_\pi)_{j_1\ldots j_k}\overline{(\xi_\sigma)}_{i_1\ldots i_k}W_{kN}(\pi,\sigma)
\end{eqnarray*}

But this gives the formula in the statement, and we are done.
\end{proof}

Let us go back now to the ``minimal vs maximal'' discussion, in analogy with the group algebras. Again by following \cite{ba6}, here is a natural example of an intermediate space $X_{G,I}^{min}\subset X\subset X_{G,I}^{max}$, which will be of interest for us, in what follows:

\begin{theorem}
Given a closed quantum subgroup $G\subset U_N^+$, and a set $I\subset\{1,\ldots,N\}$, if we consider the following quotient algebra
$$C(X_{G,I}^{med})=C(S^{N-1}_{\mathbb C,+})\Big/\left<\sum_{j_1\ldots j_k}\xi_{j_1\ldots j_k}x_{j_1}^{e_1}\ldots x_{j_k}^{e_k}=\frac{1}{\sqrt{|I|^k}}\sum_{j_1\ldots j_k\in I}\xi_{j_1\ldots j_k}\Big|\forall k,\forall\xi\in Fix(u^{\otimes k})\right>$$
we obtain in this way an affine homogeneous space $G\to X_{G,I}$.
\end{theorem}

\begin{proof}
We know from Theorem 7.14 that $X_{G,I}^{max}\subset S^{N-1}_{\mathbb C,+}$ is constructed by imposing to the standard coordinates the conditions $Px^{\otimes k}=P^I$, where:
$$P_{i_1\ldots i_k,j_1\ldots j_k}=\int_Gu_{j_1i_1}^{e_1}\ldots u_{j_ki_k}^{e_k}$$
$$P^I_{i_1\ldots i_k}=\frac{1}{\sqrt{|I|^k}}\sum_{j_1\ldots j_k\in I}P_{i_1\ldots i_k,j_1\ldots j_k}$$

According to the Weingarten integration formula for $G$, we have:
\begin{eqnarray*}
(Px^{\otimes k})_{i_1\ldots i_k}&=&\sum_{j_1\ldots j_k}\sum_{\pi,\sigma\in D}(\xi_\pi)_{j_1\ldots j_k}\overline{(\xi_\sigma)}_{i_1\ldots i_k}W_{kN}(\pi,\sigma)x_{j_1}^{e_1}\ldots x_{j_k}^{e_k}\\
P^I_{i_1\ldots i_k}&=&\frac{1}{\sqrt{|I|^k}}\sum_{j_1\ldots j_k\in I}\sum_{\pi,\sigma\in D}(\xi_\pi)_{j_1\ldots j_k}\overline{(\xi_\sigma)}_{i_1\ldots i_k}W_{kN}(\pi,\sigma)
\end{eqnarray*}

Thus $X_{G,I}^{med}\subset X_{G,I}^{max}$, and the other assertions are standard as well.
\end{proof}

We can now put everything together, as follows:

\index{full version}
\index{reduced version}
\index{GNS construction}

\begin{theorem}
Given a closed subgroup $G\subset U_N^+$, and a subset $I\subset\{1,\ldots,N\}$, the affine homogeneous spaces over $G$, with index set $I$, have the following properties:
\begin{enumerate}
\item These are exactly the intermediate subspaces $X_{G,I}^{min}\subset X\subset X_{G,I}^{max}$ on which $G$ acts affinely, with the action being ergodic.

\item For the minimal and maximal spaces $X_{G,I}^{min}$ and $X_{G,I}^{max}$, as well as for the intermediate space $X_{G,I}^{med}$ constructed above, these conditions are satisfied.

\item By performing the GNS construction with respect to the Haar integration functional $\int_X=\int_G\alpha$ we obtain the minimal space $X_{G,I}^{min}$.
\end{enumerate}
We agree to identify all these spaces, via the GNS construction, and denote them $X_{G,I}$.
\end{theorem}

\begin{proof}
This follows indeed by combining the various results and observations formulated above. Once again, for full details on all these facts, we refer to \cite{ba6}.
\end{proof}

All this might seem of course a bit technical, but this is what comes out, as abstract general theory, from the various examples of homogeneous spaces studied so far in this book. With the remark that our formalism is quite advanced, in the sense that it is not very clear that these basic examples are indeed affine homogeneous spaces in our sense. But hey, that's how mathematics goes, sometimes a new definition takes some time to be understood. We will discuss all this, examples, in the remainder of this chapter.

\section*{7d. Basic examples}

Let us first discuss, again by following \cite{ba6} and related papers, some basic examples of affine homogeneous spaces, namely those coming from the classical groups, and those coming from the group duals. We will need the following technical result: 

\begin{proposition}
Assuming that a closed subset $X\subset S^{N-1}_{\mathbb C,+}$ is affine homogeneous over a classical group, $G\subset U_N$, then $X$ itself must be classical, $X\subset S^{N-1}_\mathbb C$.
\end{proposition}

\begin{proof}
We use the well-known fact that, since the standard coordinates $u_{ij}\in C(G)$ commute, the corepresentation $u^{\circ\circ\bullet\bullet}=u^{\otimes 2}\otimes\bar{u}^{\otimes 2}$ has the following fixed vector:
$$\xi=\sum_{ij}e_i\otimes e_j\otimes e_i\otimes e_j$$

With $k=\circ\circ\bullet\,\bullet$ and with this vector $\xi$, the ergodicity formula reads:
\begin{eqnarray*}
\sum_{ij}x_ix_jx_i^*x_j^*
&=&\frac{1}{\sqrt{|I|^4}}\sum_{i,j\in I}1\\
&=&1
\end{eqnarray*}

By using this formula, along with $\sum_ix_ix_i^*=\sum_ix_i^*x_i=1$, we obtain:
\begin{eqnarray*}
&&\sum_{ij}(x_ix_j-x_jx_i)(x_j^*x_i^*-x_i^*x_j^*)\\
&=&\sum_{ij}x_ix_jx_j^*x_i^*-x_ix_jx_i^*x_j^*-x_jx_ix_j^*x_i^*+x_jx_ix_i^*x_j^*\\
&=&1-1-1+1\\
&=&0
\end{eqnarray*}

We conclude that for any $i,j$ the following commutator vanishes:
$$[x_i,x_j]=0$$

By using now this commutation relation, plus once again the relations defining the free sphere $S^{N-1}_{\mathbb C,+}$, we have as well the following computation:
\begin{eqnarray*}
&&\sum_{ij}(x_ix_j^*-x_j^*x_i)(x_jx_i^*-x_i^*x_j)\\
&=&\sum_{ij}x_ix_j^*x_jx_i^*-x_ix_j^*x_i^*x_j-x_j^*x_ix_jx_i^*+x_j^*x_ix_i^*x_j\\
&=&\sum_{ij}x_ix_j^*x_jx_i^*-x_ix_i^*x_j^*x_j-x_j^*x_jx_ix_i^*+x_j^*x_ix_i^*x_j\\
&=&1-1-1+1\\
&=&0
\end{eqnarray*}

Thus we have $[x_i,x_j^*]=0$ as well, and so $X\subset S^{N-1}_\mathbb C$, as claimed. 
\end{proof}

We can now formulate the result in the classical case, as follows:

\index{bistochastic group}

\begin{theorem}
In the classical case, $G\subset U_N$, there is only one affine homogeneous space, for each index set $I=\{1,\ldots,N\}$, namely the quotient space 
$$X=G/(G\cap C_N^I)$$
where $C_N^I\subset U_N$ is the group of unitaries fixing the following vector,
$$\xi_I=\frac{1}{\sqrt{|I|}}(\delta_{i\in I})_i$$
which generalizes the complex bistochastic group, $C_N\subset U_N$.
\end{theorem}

\begin{proof}
Consider an affine homogeneous space $G\to X$. We already know from Proposition 7.19 that $X$ is classical. We will first prove that we have $X=X_{G,I}^{min}$, and then we will prove that $X_{G,I}^{min}$ equals the quotient space in the statement.

\medskip

(1) We use the well-known fact that the functional $E=(id\otimes\int_G)\Phi$ is the projection onto the fixed point algebra of the action, given by:
$$C(X)^\Phi=\left\{f\in C(X)\Big|\Phi(f)=f\otimes1\right\}$$

Thus our ergodicity condition, namely $E=\int_G\alpha(.)1$, shows that we must have: 
$$C(X)^\Phi=\mathbb C1$$

But in the classical case the condition $\Phi(f)=f\otimes 1$ reformulates as:
$$f(gx)=f(x)\quad,\quad\forall g\in G,x\in X$$

Thus, we recover in this way the usual ergodicity condition, stating that whenever a function $f\in C(X)$ is constant on the orbits of the action, it must be constant. Now observe that for an affine action, the orbits are closed. Thus an affine action which is ergodic must be transitive, and we deduce from this that we have:
$$X=X_{G,I}^{min}$$

(2) We know that the inclusion $C(X)\subset C(G)$ comes via:
$$x_i=\frac{1}{\sqrt{|I|}}\sum_{j\in I}u_{ji}$$

Thus, the quotient map $p:G\to X\subset S^{N-1}_\mathbb C$ is given by the following formula:
$$p(g)=\left(\frac{1}{\sqrt{|I|}}\sum_{j\in I}g_{ji}\right)_i$$

In particular, the image of the unit matrix $1\in G$ is the following vector:
\begin{eqnarray*}
p(1)
&=&\left(\frac{1}{\sqrt{|I|}}\sum_{j\in I}\delta_{ji}\right)_i\\
&=&\left(\frac{1}{\sqrt{|I|}}\delta_{i\in I}\right)_i\\
&=&\xi_I
\end{eqnarray*}

But this gives the quotient space result in the statement.

\medskip

(3) Finally, regarding the last assertion, stating that our group $C_N^I\subset U_N$ generalizes the complex bishochastic group $C_N\subset U_N$, this is more of a comment, coming from definitions. Indeed, $C_N$ consists by definition of the unitary matrices $g\in U_N$ which are bistochastic, meaning having the same sums on rows and columns. But this bistochasticity condition is equivalent to the following condition, with $\xi$ being the all-1 vector:
$$g\xi=\xi$$

Thus, our group $C_N^I\subset U_N$ generalizes indeed the group $C_N\subset U_N$, as claimed.
\end{proof}

Again by following \cite{ba6}, let us discuss now the group dual case. For simplifying, we will discuss the case of the ``diagonal'' embeddings only. Given a finitely generated discrete group $\Gamma=<g_1,\ldots,g_N>$, we can consider the following ``diagonal'' embedding: 
$$\widehat{\Gamma}\subset U_N^+\quad,\quad u_{ij}=\delta_{ij}g_i$$

With this convention, we have the following result:

\begin{theorem}
In the group dual case, $G=\widehat{\Gamma}$ with $\Gamma=<g_1,\ldots,g_N>$, we have
$$X=\widehat{\Gamma}_I\quad:\quad 
\Gamma_I=<g_i|i\in I>\subset\Gamma$$
for any affine homogeneous space $X$, when identifying full and reduced group algebras.
\end{theorem}

\begin{proof}
Assume indeed that we have an affine homogeneous space $G\to X$. In terms of the rescaled coordinates $h_i=\sqrt{|I|}x_i$, our axioms for $\alpha,\Phi$ read:
$$\alpha(h_i)=\delta_{i\in I}g_i$$
$$\Phi(h_i)=h_i\otimes g_i$$

As for the ergodicity condition, this translates as follows:
\begin{eqnarray*}
&&\left(id\otimes\int_G\right)\Phi(h_{i_1}^{e_1}\ldots h_{i_p}^{e_p})=\int_G\alpha(h_{i_1}^{e_p}\ldots h_{i_p}^{e_p})\\
&\iff&\left(id\otimes\int_G\right)(h_{i_1}^{e_1}\ldots h_{i_p}^{e_p}\otimes g_{i_1}^{e_1}\ldots g_{i_p}^{e_p})=\int_G\delta_{i_1\in I}\ldots\delta_{i_p\in I}g_{i_1}^{e_1}\ldots g_{i_p}^{e_p}\\
&\iff&\delta_{g_{i_1}^{e_1}\ldots g_{i_p}^{e_p},1}h_{i_1}^{e_1}\ldots h_{i_p}^{e_p}=\delta_{g_{i_1}^{e_1}\ldots g_{i_p}^{e_p},1}\delta_{i_1\in I}\ldots\delta_{i_p\in I}\\
&\iff&\left[g_{i_1}^{e_1}\ldots g_{i_p}^{e_p}=1\implies h_{i_1}^{e_1}\ldots h_{i_p}^{e_p}=\delta_{i_1\in I}\ldots\delta_{i_p\in I}\right]
\end{eqnarray*}

Now observe that from $g_ig_i^*=g_i^*g_i=1$ we obtain in this way: 
$$h_ih_i^*=h_i^*h_i=\delta_{i\in I}$$

Thus the elements $h_i$ vanish for $i\notin I$, and are unitaries for $i\in I$. We conclude that we have $X=\widehat{\Lambda}$, where $\Lambda=<h_i|i\in I>$ is the group generated by these unitaries. In order to finish now the proof, our claim is that for indices $i_x\in I$ we have:
$$g_{i_1}^{e_1}\ldots g_{i_p}^{e_p}=1\iff h_{i_1}^{e_1}\ldots h_{i_p}^{e_p}=1$$

Indeed, $\implies$ comes from the ergodicity condition, as processed above, and $\Longleftarrow$ comes from the existence of the morphism $\alpha$, which is given by $\alpha(h_i)=g_i$, for $i\in I$.
\end{proof}

Let us go back now to the general case, and discuss a number of further axiomatization issues, based on the examples that we have. We will need the following result:

\index{complex bistochastic group}

\begin{proposition}
The closed subspace $C_N^{I+}\subset U_N^+$ defined via
$$C(C_N^{I+})=C(U_N^+)\Big/\left<u\xi_I=\xi_I\right>$$
where $\xi_I=\frac{1}{\sqrt{|I|}}(\delta_{i\in I})_i$, is a compact quantum group.
\end{proposition}

\begin{proof}
We must check Woronowicz's axioms, and the proof goes as follows:

\medskip

(1) Let us set $U_{ij}=\sum_ku_{ik}\otimes u_{kj}$. We have then:
\begin{eqnarray*}
(U\xi_I)_i
&=&\frac{1}{\sqrt{|I|}}\sum_{j\in I}U_{ij}\\
&=&\frac{1}{\sqrt{|I|}}\sum_{j\in I}\sum_ku_{ik}\otimes u_{kj}\\
&=&\sum_ku_{ik}\otimes(u\xi_I)_k
\end{eqnarray*}

Since the vector $\xi_I$ is by definition fixed by $u$, we obtain:
\begin{eqnarray*}
(U\xi_I)_i
&=&\sum_ku_{ik}\otimes(\xi_I)_k\\
&=&\frac{1}{\sqrt{|I|}}\sum_{k\in I}u_{ik}\otimes1\\
&=&(u\xi_I)_i\otimes1\\
&=&(\xi_I)_i\otimes1
\end{eqnarray*}

Thus we can define indeed a comultiplication map, by $\Delta(u_{ij})=U_{ij}$.

\medskip

(2) In order to construct the counit map, $\varepsilon(u_{ij})=\delta_{ij}$, we must prove that the identity matrix $1=(\delta_{ij})_{ij}$ satisfies $1\xi_I=\xi_I$. But this is clear.

\medskip

(3) In order to construct the antipode, $S(u_{ij})=u_{ji}^*$, we must prove that the adjoint matrix $u^*=(u_{ji}^*)_{ij}$ satisfies $u^*\xi_I=\xi_I$. But this is clear from $u\xi_I=\xi_I$.
\end{proof}

Based on the computations that we have so far, we can formulate:

\begin{theorem}
Given a closed quantum subgroup $G\subset U_N^+$ and a set $I\subset\{1,\ldots,N\}$, we have a quotient map and an inclusion map as follows:
$$G/(G\cap C_N^{I+})\to X_{G,I}^{min}\subset X_{G,I}^{max}$$
These maps are both isomorphisms in the classical case. In general, they are both proper.
\end{theorem}

\begin{proof}
Consider the quantum group $H=G\cap C_N^{I+}$, which is by definition such that at the level of the corresponding algebras, we have:
$$C(H)=C(G)\Big/\left<u\xi_I=\xi_I\right>$$

In order to construct a quotient map $G/H\to X_{G,I}^{min}$, we must check that the defining relations for $C(G/H)$ hold for the standard generators $x_i\in C(X_{G,I}^{min})$. But if we denote by $\rho:C(G)\to C(H)$ the quotient map, then we have, as desired:
\begin{eqnarray*}
(id\otimes\rho)\Delta x_i
&=&(id\otimes\rho)\left(\frac{1}{\sqrt{|I|}}\sum_{j\in I}\sum_ku_{ik}\otimes u_{kj}\right)\\
&=&\sum_ku_{ik}\otimes(\xi_I)_k\\
&=&x_i\otimes1
\end{eqnarray*}

In the classical case, Theorem 7.20 shows that both the maps in the statement are isomorphisms. For the group duals, however, these maps are not isomorphisms, in general. This follows indeed from Theorem 7.21, and from the general theory in \cite{bss}.
\end{proof}

We discuss now a number of further examples. We will need:

\index{transposed quantum group}

\begin{proposition}
Given a compact matrix quantum group $G=(G,u)$, the pair 
$$G^t=(G,u^t)$$
where $(u^t)_{ij}=u_{ji}$, is a compact matrix quantum group as well.
\end{proposition}

\begin{proof}
The construction of the comultiplication is as follows, where $\Sigma$ is the flip:
\begin{eqnarray*}
\Delta^t[(u^t)_{ij}]=\sum_k(u^t)_{ik}\otimes(u^t)_{kj}
&\iff&\Delta^t(u_{ji})=\sum_ku_{ki}\otimes u_{jk}\\
&\iff&\Delta^t=\Sigma\Delta
\end{eqnarray*}

As for the corresponding counit and antipode, these can be simply taken to be $(\varepsilon,S)$, and the axioms of Woronowicz are then satisfied.
\end{proof}

We will need as well the following result, which is standard too:

\index{product of quantum groups}

\begin{proposition}
Given closed subgroups $G\subset U_N^+$ and $H\subset U_M^+$, with fundamental corepresentations $u=(u_{ij})$ and $v=(v_{ab})$, their product is a closed subgroup 
$$G\times H\subset U_{NM}^+$$
with fundamental corepresentation $w_{ia,jb}=u_{ij}\otimes v_{ab}$. 
\end{proposition}

\begin{proof}
Our claim is that the corresponding structural maps are as follows:
$$\Delta(\alpha\otimes\beta)=\Delta(\alpha)_{13}\Delta(\beta)_{24}$$
$$\varepsilon(\alpha\otimes\beta)=\varepsilon(\alpha)\varepsilon(\beta)$$
$$S(\alpha\otimes\beta)=S(\alpha)S(\beta)$$

Indeed, the verification for the comultiplication goes as follows:
\begin{eqnarray*}
\Delta(w_{ia,jb})
&=&\Delta(u_{ij})_{13}\Delta(v_{ab})_{24}\\
&=&\sum_{kc}u_{ik}\otimes v_{ac}\otimes u_{kj}\otimes v_{cb}\\
&=&\sum_{kc}w_{ia,kc}\otimes w_{kc,jb}
\end{eqnarray*}

For the counit, we have the following computation:
\begin{eqnarray*}
\varepsilon(w_{ia,jb})
&=&\varepsilon(u_{ij})\varepsilon(v_{ab})\\
&=&\delta_{ij}\delta_{ab}\\
&=&\delta_{ia,jb}
\end{eqnarray*}

As for the antipode, here we have the following computation:
\begin{eqnarray*}
S(w_{ia,jb})
&=&S(u_{ij})S(v_{ab})\\
&=&v_{ba}^*u_{ji}^*\\
&=&(u_{ji}v_{ba})^*\\
&=&w_{jb,ia}^*
\end{eqnarray*}

We refer to Wang's paper \cite{wa1} for more details regarding this construction.
\end{proof}

We will need one more ingredient, which is a definition, as follows:

\begin{definition}
We call a closed quantum subgroup $G\subset U_N^+$ self-transpose when we have an automorphism $T:C(G)\to C(G)$ given by $T(u_{ij})=u_{ji}$.
\end{definition}

Observe that in the classical case, this amounts in our closed subgroup $G\subset U_N$ to be closed under the transposition operation $g\to g^t$.

\bigskip

With the above notions and general theory in hand, let us go back to the affine homogeneous spaces. As a first result here, any closed subgroup $G\subset U_N^+$ appears as an affine homogeneous space over an appropriate quantum group, as follows:

\begin{theorem}
Given a closed subgroup $G\subset U_N^+$, we have an identification 
$$X_{\mathcal G,I}^{min}\simeq G$$
given at the level of standard coordinates by $x_{ij}=\frac{1}{\sqrt{N}}u_{ij}$, where:
\begin{enumerate}
\item $\mathcal G=G^t\times G\subset U_{N^2}^+$, with coordinates $w_{ia,jb}=u_{ji}\otimes u_{ab}$.

\item $I\subset\{1,\ldots,N\}^2$ is the diagonal set, $I=\{(k,k)|k=1,\ldots,N\}$.
\end{enumerate}
In the self-transpose case we can choose as well $\mathcal G=G\times G$, with $w_{ia,jb}=u_{ij}\otimes u_{ab}$.
\end{theorem}

\begin{proof}
As a first observation, our closed subgroup $G\subset U_N^+$ appears as an algebraic submanifold of the free complex sphere on $N^2$ variables, as follows:
$$G\subset S^{N^2-1}_{\mathbb C,+}\quad,\quad 
x_{ij}=\frac{1}{\sqrt{N}}u_{ij}$$

Let us construct now the affine homogeneous space structure. Our claim is that, with $\mathcal G=G^t\times G$ and $I=\{(k,k)\}$ as in the statement, the structural maps are:
$$\alpha=\Delta$$
$$\Phi=(\Sigma\otimes id)\Delta^{(2)}$$

Indeed, in what regards $\alpha=\Delta$, this is given by the following formula:
\begin{eqnarray*}
\alpha(u_{ij})
&=&\sum_ku_{ik}\otimes u_{kj}\\
&=&\sum_kw_{kk,ij}
\end{eqnarray*}

Thus, by dividing by $\sqrt{N}$, we obtain the usual affine homogeneous space formula:
$$\alpha(x_{ij})=\frac{1}{\sqrt{|I|}}\sum_kw_{kk,ij}$$

Regarding now $\Phi=(\Sigma\otimes id)\Delta^{(2)}$, the formula here is as follows:
\begin{eqnarray*}
\Phi(u_{ij})
&=&(\Sigma\otimes id)\sum_{kl}u_{ik}\otimes u_{kl}\otimes u_{lj}\\
&=&\sum_{kl}u_{kl}\otimes u_{ik}\otimes u_{lj}\\
&=&\sum_{kl}u_{kl}\otimes w_{kl,ij}
\end{eqnarray*}

Thus, by dividing by $\sqrt{N}$, we obtain the usual affine homogeneous space formula:
$$\Phi(x_{ij})=\sum_{kl}x_{kl}\otimes w_{kl,ij}$$

The ergodicity condition being clear as well, this gives the first assertion. Regarding now the second assertion, assume that we are in the self-transpose case, and so that we have an automorphism $T:C(G)\to C(G)$ given by $T(u_{ij})=u_{ji}$. With the notation $w_{ia,jb}=u_{ij}\otimes u_{ab}$, the modified map $\alpha=(T\otimes id)\Delta$ is then given by:
\begin{eqnarray*}
\alpha(u_{ij})
&=&(T\otimes id)\sum_ku_{ik}\otimes u_{kj}\\
&=&\sum_ku_{ki}\otimes u_{kj}\\
&=&\sum_kw_{kk,ij}
\end{eqnarray*}

As for the modified map $\Phi=(id\otimes T\otimes id)(\Sigma\otimes id)\Delta^{(2)}$, this is given by:
\begin{eqnarray*}
\Phi(u_{ij})
&=&(id\otimes T\otimes id)\sum_{kl}u_{kl}\otimes u_{ik}\otimes u_{lj}\\
&=&\sum_{kl}u_{kl}\otimes u_{ki}\otimes u_{lj}\\
&=&\sum_{kl}u_{kl}\otimes w_{kl,ij}
\end{eqnarray*}

Thus we have the correct affine homogeneous space formulae, and once again the ergodicity condition being clear as well, this gives the result.
\end{proof}

Let us discuss now the generalization of the above result, to the context of the spaces introduced in \cite{bss}. We recall from there that we have the following construction:

\index{row space}
\index{row algebra}

\begin{definition}
Given a closed subgroup $G\subset U_N^+$ and an integer $M\leq N$ we set
$$C(G_{MN})=\left<u_{ij}\Big|i\in\{1,\ldots,M\},j\in\{1,\ldots,N\}\right>\subset C(G)$$
and we call row space of $G$ the underlying quotient space $G\to G_{MN}$.
\end{definition}

As a basic example here, at $M=N$ we obtain $G$ itself. Also, at $M=1$ we obtain the space whose coordinates are those on the first row of coordinates on $G$. Finally, in the case of the basic quantum unitary and reflection groups, these are particular cases of the partial isometry spaces discussed in chapter 6. See \cite{bss}.

\bigskip

Given $G_N\subset U_N^+$ and an integer $M\leq N$, we can consider the quantum group $G_M=G_N\cap U_M^+$, with the intersection taken inside $U_N^+$, and with $U_M^+\subset U_N^+$ given by:
$$u=diag(v,1_{N-M})$$

Observe that we have a quotient map $C(G_N)\to C(G_M)$, given by $u_{ij}\to v_{ij}$. With these conventions, we have the following extension of Theorem 7.27:

\begin{theorem}
Given a closed subgroup $G_N\subset U_N^+$, we have an identification 
$$X_{\mathcal G,I}^{min}\simeq G_{MN}$$
given at the level of standard coordinates by $x_{ij}=\frac{1}{\sqrt{M}}u_{ij}$, where:
\begin{enumerate}
\item $\mathcal G=G_M^t\times G_N\subset U_{NM}^+$, where $G_M=G_N\cap U_M^+$, with coordinates as follows:
$$w_{ia,jb}=u_{ji}\otimes v_{ab}$$

\item $I\subset\{1,\ldots,M\}\times\{1,\ldots,N\}$ is the diagonal set, namely:
$$I=\left\{(k,k)\Big|k=1,\ldots,M\right\}$$
\end{enumerate}
In the self-transpose case we can choose as well $\mathcal G=G_M\times G_N$, with $w_{ia,jb}=u_{ij}\otimes v_{ab}$.
\end{theorem}

\begin{proof}
Consider the row space $X=G_{MN}$ constructed in Definition 7.28, with its standard row space coordinates, namely:
$$x_{ij}=\frac{1}{\sqrt{M}}u_{ij}$$

In order to prove the result, we have to show that this space coincides with the space $X_{\mathcal G,I}^{min}$ constructed in the statement, with its standard coordinates. For this purpose, consider the following composition of morphisms, where in the middle we have the comultiplication, and at left and right we have the canonical maps:
$$C(X)\subset C(G_N)\to C(G_N)\otimes C(G_N)\to C(G_M)\otimes C(G_N)$$

The standard coordinates are then mapped as follows:
\begin{eqnarray*}
x_{ij}
&=&\frac{1}{\sqrt{M}}u_{ij}\\
&\to&\frac{1}{\sqrt{M}}\sum_ku_{ik}\otimes u_{kj}\\
&\to&\frac{1}{\sqrt{M}}\sum_{k\leq M}u_{ik}\otimes v_{kj}\\
&=&\frac{1}{\sqrt{M}}\sum_{k\leq M}w_{kk,ij}
\end{eqnarray*}

Thus we obtain the standard coordinates on the space $X_{\mathcal G,I}^{min}$, as claimed. Finally, the last assertion is standard as well, by suitably modifying the above morphism.
\end{proof}

Summarizing, our notion of affine homogeneous space covers, but in a somewhat tricky and technical way, all the examples of homogeneous spaces discussed so far in this book. Again, and as mentioned on several occasions, the point here is that this theory comes from a long series of papers, namely \cite{bgo}, followed by \cite{bss}, then by \cite{ba5}, then by \cite{ba6}, and a number of secondary papers as well, which each paper considerably building on the previous ones. And so this theory is quite abstract and advanced. We will keep studying such spaces in the next chapter, with a number of algebraic and analytic results.

\section*{7e. Exercises} 

The material in this chapter has been quite abstract and technical, and so will be our exercises here. To start with, in relation with the axioms, we have:

\begin{exercise}
Study, with some general theory, examples and counterexamples, the validity of the well-known statement ``the quotient by a normal subgroup is a group'' from group theory, in the compact quantum group setting.
\end{exercise}

This is a very good exercise, the answer to it being folklore. Enjoy.

\begin{exercise}
Try developing a theory of noncommutative homogeneous spaces, both algebraic and geometric, paralleling what is known classically.
\end{exercise}

As already mentioned on several occasions, such theories are not very useful in the free case, where things are quite wild and special. However, we will be back to more ``tame'' geometries, such as the half-classical ones, or the twisted ones, which are quite close to the classical world, in chapters 9-12 below, and such prior knowledge can only help.

\chapter{Liberation theory}

\section*{8a. Integration results}

We discuss in this chapter a number of further topics, in relation with what was said in chapters 5-7, namely liberation theory, Bercovici-Pata bijection and Tannakian duality for the affine homogeneous spaces, along with the question of axiomatizing the free manifolds, following \cite{ba6} and related papers, and then the formalism of row spaces from \cite{bss} and related papers, which goes in a rather opposite direction, namely particularization.

\bigskip

Let us also mention that things will be basically about open problems that we don't know how to solve, with the whole material being quite recent, and research-grade. Many questions here are waiting for enthusiastic young people. Like you.

\bigskip

Let us first discuss the liberation operation, in the context of the affine homogeneous spaces, following \cite{ba6}. In the easy case, we have the following result:

\begin{proposition}
When $G\subset U_N^+$ is easy, coming from a category of partitions $D$, the space $X_{G,I}\subset S^{N-1}_{\mathbb C,+}$ appears by imposing the relations
$$\sum_{i_1\ldots i_k}\delta_\pi(i_1\ldots i_k)x_{i_1}^{e_1}\ldots x_{i_k}^{e_k}=|I|^{|\pi|-k/2},\quad\forall k,\forall\pi\in D(k)$$
where $D(k)=D(0,k)$, and where $|.|$ denotes the number of blocks.
\end{proposition}

\begin{proof}
We know by easiness that $Fix(u^{\otimes k})$ is spanned by the vectors $\xi_\pi=T_\pi$, with $\pi\in D(k)$. But these latter vectors are given by:
$$\xi_\pi=\sum_{i_1\ldots i_k}\delta_\pi(i_1\ldots i_k)e_{i_1}\otimes\ldots\otimes e_{i_k}$$

We deduce that $X_{G,I}\subset S^{N-1}_{\mathbb C,+}$ appears by imposing the following relations:
$$\sum_{i_1\ldots i_k}\delta_\pi(i_1\ldots i_k)x_{i_1}^{e_1}\ldots x_{i_k}^{e_k}=\frac{1}{\sqrt{|I|^k}}\sum_{j_1\ldots j_k\in I}\delta_\pi(j_1\ldots j_k),\quad\forall k,\forall\pi\in D(k)$$

Now since the sum on the right equals $|I|^{|\pi|}$, this gives the result.
\end{proof}

More generally now, in view of the examples given at the end of chapter 7, making the link with \cite{bss}, it is interesting to work out what happens when $G$ is a product of easy quantum groups, and the index set $I$ above appears as $I=\{(c,\ldots,c)|c\in J\}$, for a certain set $J$. The result here, in its most general form, also from \cite{ba6}, is as follows:

\index{product of easy quantum groups}

\begin{theorem}
For a product of easy quantum groups
$$G=G_{N_1}^{(1)}\times\ldots\times G_{N_s}^{(s)}$$
and with $I=\{(c,\ldots,c)|c\in J\}$, the space $X_{G,I}\subset S^{N-1}_{\mathbb C,+}$ appears via the relations
$$\sum_{i_1\ldots i_k}\delta_\pi(i_1\ldots i_k)x_{i_1}^{e_1}\ldots x_{i_k}^{e_k}=|J|^{|\pi_1\vee\ldots\vee\pi_s|-k/2}$$
for any $k\in\mathbb N$ and any partition of the following type,
$$\pi\in D^{(1)}(k)\times\ldots\times D^{(s)}(k)$$
where $D^{(r)}\subset P$ is the category of partitions associated to $G_{N_r}^{(r)}\subset U_{N_r}^+$, and where 
$$\pi_1\vee\ldots\vee\pi_s\in P(k)$$
is the partition obtained by superposing $\pi_1,\ldots,\pi_s$.
\end{theorem}

\begin{proof}
Since we are in a direct product situation, $G=G_{N_1}^{(1)}\times\ldots\times G_{N_s}^{(s)}$, the general product theory of Wang \cite{wa1} applies, and shows that a basis for $Fix(u^{\otimes k})$ is provided by the vectors $\rho_\pi=\xi_{\pi_1}\otimes\ldots\otimes\xi_{\pi_s}$ associated to the following partitions:
$$\pi=(\pi_1,\ldots,\pi_s)\in D^{(1)}(k)\times\ldots\times D^{(s)}(k)$$

We conclude that the space $X_{G,I}\subset S^{N-1}_{\mathbb C,+}$ appears by imposing the following relations to the standard coordinates:
$$\sum_{i_1\ldots i_k}\delta_\pi(i_1\ldots i_k)x_{i_1}^{e_1}\ldots x_{i_k}^{e_k}=\frac{1}{\sqrt{|I|^k}}\sum_{j_1\ldots j_k\in I}\delta_\pi(j_1\ldots j_k),\ \forall k,\forall\pi\in D^{(1)}(k)\times\ldots\times D^{(s)}(k)$$

Since the conditions $j_1,\ldots,j_k\in I$ read $j_1=(l_1,\ldots,l_1),\ldots,j_k=(l_k,\ldots,l_k)$, for certain elements $l_1,\ldots l_k\in J$, the sums on the right are given by:
\begin{eqnarray*}
\sum_{j_1\ldots j_k\in I}\delta_\pi(j_1\ldots j_k)
&=&\sum_{l_1\ldots l_k\in J}\delta_\pi(l_1,\ldots,l_1,\ldots\ldots,l_k,\ldots,l_k)\\
&=&\sum_{l_1\ldots l_k\in J}\delta_{\pi_1}(l_1\ldots l_k)\ldots\delta_{\pi_s}(l_1\ldots l_k)\\
&=&\sum_{l_1\ldots l_k\in J}\delta_{\pi_1\vee\ldots\vee\pi_s}(l_1\ldots l_k)
\end{eqnarray*}

Now since the sum on the right equals $|J|^{|\pi_1\vee\ldots\vee\pi_s|}$, this gives the result.
\end{proof}

We can now discuss probabilistic aspects. Following \cite{ba6}, we first have:

\begin{proposition}
The moments of the variable 
$$\chi_T=\sum_{i\leq T}x_{i\ldots i}$$
are given by the following formula,
$$\int_X\chi_T^k\simeq\frac{1}{\sqrt{M^k}}\sum_{\pi\in D^{(1)}(k)\cap\ldots\cap D^{(s)}(k)}\left(\frac{TM}{N}\right)^{|\pi|}$$
in the $N_i\to\infty$ limit, $\forall i$, where $M=|I|$, and $N=N_1\ldots N_s$.
\end{proposition}

\begin{proof}
We have the following formula:
$$\pi(x_{i_1\ldots i_s})=\frac{1}{\sqrt{M}}\sum_{c\in J}u_{i_1c}\otimes\ldots\otimes u_{i_sc}$$

For the variable in the statement, we therefore obtain:
$$\pi(\chi_T)=\frac{1}{\sqrt{M}}\sum_{i\leq T}\sum_{c\in J}u_{ic}\otimes\ldots\otimes u_{ic}$$

Now by raising to the power $k$ and integrating, we obtain:
\begin{eqnarray*}
\int_X\chi_T^k
&=&\frac{1}{\sqrt{M^k}}\sum_{i_1\ldots i_k\leq T}\sum_{c_1\ldots c_k\in J}\int_{G^{(1)}}u_{i_1c_1}\ldots u_{i_kc_k}\ldots\ldots\int_{G^{(s)}}u_{i_1c_1}\ldots u_{i_kc_k}\\
&=&\frac{1}{\sqrt{M^k}}\sum_{ic}\sum_{\pi\sigma}\delta_{\pi_1}(i)\delta_{\sigma_1}(c)W_{kN_1}^{(1)}(\pi_1,\sigma_1)\ldots\delta_{\pi_s}(i)\delta_{\sigma_s}(c)W_{kN_s}^{(s)}(\pi_s,\sigma_s)\\
&=&\frac{1}{\sqrt{M^k}}\sum_{\pi\sigma}T^{|\pi_1\vee\ldots\vee\pi_s|}M^{|\sigma_1\vee\ldots\vee\sigma_s|}
W_{kN_1}^{(1)}(\pi_1,\sigma_1)\ldots W_{kN_s}^{(s)}(\pi_s,\sigma_s)
\end{eqnarray*}

We use now the standard fact that the Weingarten functions are concentrated on the diagonal. Thus in the limit we must have $\pi_i=\sigma_i$ for any $i$, and we obtain:
\begin{eqnarray*}
\int_X\chi_T^k
&\simeq&\frac{1}{\sqrt{M^k}}\sum_\pi T^{|\pi_1\vee\ldots\vee\pi_s|}M^{|\pi_1\vee\ldots\vee\pi_s|}N_1^{-|\pi_1|}\ldots N_s^{-|\pi_s|}\\
&\simeq&\frac{1}{\sqrt{M^k}}\sum_{\pi\in D^{(1)}\cap\ldots\cap D^{(s)}}T^{|\pi|}M^{|\pi|}(N_1\ldots N_s)^{-|\pi|}\\
&=&\frac{1}{\sqrt{M^k}}\sum_{\pi\in D^{(1)}\cap\ldots\cap D^{(s)}}\left(\frac{TM}{N}\right)^{|\pi|}
\end{eqnarray*}

But this gives the formula in the statement, and we are done.
\end{proof}

As a consequence, we have the following result, also from \cite{ba6}:

\index{Bercovici-Pata bijection}

\begin{theorem}
In the context of a liberation operation for quantum groups
$$G^{(i)}\to G^{(i)+}$$
the laws of the variables $\sqrt{M}\chi_T$ are in Bercovici-Pata bijection, in the $N_i\to\infty$ limit.
\end{theorem}

\begin{proof}
Assume indeed that we have easy quantum groups $G^{(1)},\ldots,G^{(s)}$, with free versions $G^{(1)+},\ldots,G^{(s)+}$. At the level of the categories of partitions, we have:
$$\bigcap_i\left(D^{(i)}\cap NC\right)=\left(\bigcap_iD^{(i)}\right)\cap NC$$

Since the intersection of Hom-spaces is the Hom-space for the generated quantum group, we deduce that at the quantum group level, we have:
$$<G^{(1)+},\ldots,G^{(s)+}>=<G^{(1)},\ldots,G^{(s)}>^+$$

Thus the result follows from Proposition 8.3, and from the Bercovici-Pata bijection result for truncated characters for this latter liberation operation \cite{bsp}, \cite{twe}.
\end{proof}

The above result is of course not the end of the story, among others because it leads into the question of enlarging the theory of easy quantum groups, as to cover the products of such quantum groups. And the answer to this latter question is not known.

\section*{8b. Tannakian duality}

We recall from the beginning of chapter 5 that one of our main goals is to axiomatize the ``free manifolds''. In this section we discuss this question, not with the idea of solving it, but rather with the idea of explaining why this question is difficult.

\bigskip

To be more precise, we will be interested, as a warm-up to the axiomatization question for the free manifolds, in the question of axiomatizing the affine homogeneous spaces, as submanifolds of the free sphere $S^{N-1}_{\mathbb C,+}$. As a starting point here, we have:

\begin{proposition}
Any affine homogeneous space $X_{G,I}\subset S^{N-1}_{\mathbb C,+}$ is algebraic, with
$$\sum_{i_1\ldots i_k}\xi_{i_1\ldots i_k}x_{i_1}^{e_1}\ldots x_{i_k}^{e_k}=\frac{1}{\sqrt{|I|^k}}\sum_{b_1\ldots b_k\in I}\xi_{b_1\ldots b_k}\quad\forall k,\forall\xi\in Fix(u^{\otimes k})$$
as defining relations. Moreover, we can use vectors $\xi$ belonging to a basis of $Fix(u^{\otimes k})$.
\end{proposition}

\begin{proof}
This follows indeed from the various results from chapter 7.
\end{proof}

In order to reach to a more categorical description of $X_{G,I}$, the idea will be that of using Frobenius duality. We use colored indices, and we denote by $k\to\bar{k}$ the operation on the colored indices which consists in reversing the index, and switching all the colors. Also, we agree to identify the linear maps $T:(\mathbb C^N)^{\otimes k}\to(\mathbb C^N)^{\otimes l}$ with the corresponding rectangular matrices $T\in M_{N^l\times N^k}(\mathbb  C)$, written $T=(T_{i_1\ldots i_l,j_1\ldots j_k})$. With these conventions, the precise formulation of Frobenius duality that we will need is as follows:

\index{Frobenius duality}

\begin{proposition}
We have an isomorphism of complex vector spaces
$$T\in Hom(u^{\otimes k},u^{\otimes l})\ \leftrightarrow\ \xi\in Fix(u^{\otimes l}\otimes u^{\otimes\bar{k}})$$
given by the following formulae,
$$T_{i_1\ldots i_l,j_1\ldots j_k}=\xi_{i_1\ldots i_lj_k\ldots j_1}$$
$$\xi_{i_i\ldots i_lj_1\ldots j_k}=T_{i_1\ldots i_l,j_k\ldots j_1}$$
and called Frobenius duality.
\end{proposition}

\begin{proof}
This is a well-known result, which follows from the general theory in \cite{wo1}. To be more precise, given integers $K,L\in\mathbb N$, consider the following standard isomorphism, which in matrix notation makes $T=(T_{IJ})\in M_{L\times K}(\mathbb C)$ correspond to $\xi=(\xi_{IJ})$:
$$T\in\mathcal L(\mathbb C^{\otimes K},\mathbb C^{\otimes L})\ \leftrightarrow\ \xi\in\mathbb C^{\otimes L+K}$$

Given now two arbitrary corepresentations $v\in M_K(C(G))$ and $w\in M_L(C(G))$, the abstract Frobenius duality result established by Woronowicz in \cite{wo1} states that the above isomorphism restricts into an isomorphism of vector spaces, as follows:
$$T\in Hom(v,w)\ \leftrightarrow\ \xi\in Fix(w\otimes\bar{v})$$

In our case, we can apply this result with $v=u^{\otimes k}$ and $w=u^{\otimes l}$. Since, according to our conventions, we have $\bar{v}=u^{\otimes\bar{k}}$, this gives the isomorphism in the statement.
\end{proof}

With the above result in hand, we can enhance the construction of $X_{G,I}$, as follows:

\begin{theorem}
Any affine homogeneous space $X_{G,I}$ is algebraic, with
$$\sum_{i_1\ldots i_l}\sum_{j_1\ldots j_k}T_{i_1\ldots i_l,j_1\ldots j_k}x_{i_1}^{e_1}\ldots x_{i_l}^{e_l}(x_{j_1}^{f_1}\ldots x_{j_k}^{f_k})^*=\frac{1}{\sqrt{|I|^{k+l}}}\sum_{b_1\ldots b_l\in I}\sum_{c_1\ldots c_k\in I}T_{b_1\ldots b_l,c_1\ldots c_k}$$
for any $k,l$, and any $T\in Hom(u^{\otimes k},u^{\otimes l})$, as defining relations.
\end{theorem}

\begin{proof}
We must prove that the relations in the statement are satisfied, over $X_{G,I}$. We know from Proposition 8.5 that, with $k\to l\bar{k}$, the following relation holds:
$$\sum_{i_1\ldots i_l}\sum_{j_1\ldots j_k}\xi_{i_1\ldots i_lj_k\ldots j_1}x_{i_1}^{e_1}\ldots x_{i_l}^{e_l}x_{j_k}^{\bar{f}_k}\ldots x_{j_1}^{\bar{f}_1}=\frac{1}{\sqrt{|I|^{k+l}}}\sum_{b_1\ldots b_l\in I}\sum_{c_1\ldots c_k\in I}\xi_{b_1\ldots b_lc_k\ldots c_1}$$

In terms of the matrix $T_{i_1\ldots i_l,j_1\ldots j_k}=\xi_{i_1\ldots i_lj_k\ldots j_1}$ from Proposition 8.6, we obtain:
$$\sum_{i_1\ldots i_l}\sum_{j_1\ldots j_k}T_{i_1\ldots i_l,j_1\ldots j_k}x_{i_1}^{e_1}\ldots x_{i_l}^{e_l}x_{j_k}^{\bar{f}_k}\ldots x_{j_1}^{\bar{f}_1}=\frac{1}{\sqrt{|I|^{k+l}}}\sum_{b_1\ldots b_l\in I}\sum_{c_1\ldots c_k\in I}T_{b_1\ldots b_l,c_1\ldots c_k}$$

But this gives the formula in the statement, and we are done.
\end{proof}

The above results suggest the following notion:

\begin{definition}
Given a submanifold $X\subset S^{N-1}_{\mathbb C,+}$ and a subset $I\subset\{1,\ldots,N\}$, we say that $X$ is $I$-affine when $C(X)$ is presented by relations of type
$$\sum_{i_1\ldots i_l}\sum_{j_1\ldots j_k}T_{i_1\ldots i_l,j_1\ldots j_k}x_{i_1}^{e_1}\ldots x_{i_l}^{e_l}(x_{j_1}^{f_1}\ldots x_{j_k}^{f_k})^*=\frac{1}{\sqrt{|I|^{k+l}}}\sum_{b_1\ldots b_l\in I}\sum_{c_1\ldots c_k\in I}T_{b_1\ldots b_l,c_1\ldots c_k}$$
with the operators $T$ belonging to certain linear spaces 
$$F(k,l)\subset M_{N^l\times N^k}(\mathbb C)$$
which altogether form a tensor category $F=(F(k,l))$.
\end{definition}

According to Theorem 8.7, any affine homogeneous space $X_{G,I}$ is an $I$-affine manifold, with the corresponding tensor category being the one associated to the quantum group $G\subset U_N^+$ which produces it, formed by the following linear spaces: 
$$F(k,l)=Hom(u^{\otimes k},u^{\otimes l})$$

Let us study now the quantum isometry groups $G^+(X)$ of the manifolds $X\subset S^{N-1}_{\mathbb C,+}$ which are $I$-affine, in the above sense. We have here the following result:

\begin{proposition}
For an $I$-affine manifold $X\subset S^{N-1}_{\mathbb C,+}$ we have 
$$G\subset G^+(X)$$
where $G\subset U_N^+$ is the Tannakian dual of the associated tensor category $F$.
\end{proposition}

\begin{proof}
We recall from chapter 3 that the relations defining $G^+(X)$ are those expressing the vanishing of the following quantities:
$$P(X_1,\ldots,X_N)=\sum_r\alpha_r\sum_{j_1^r\ldots j_{s(r)}^r}u_{i_1^rj_1^r}\ldots u_{i_{s(r)}^rj_{s(r)}^r}\otimes x_{j_1^r}\ldots x_{j_{s(r)}^r}$$

In the case of an $I$-affine manifold, the defining relations are those from Definition 8.8 above, with the corresponding polynomials $P$ being indexed by the elements of $F$. But the vanishing of the associated relations $P(X_1,\ldots,X_N)=0$ corresponds precisely to the Tannakian relations defining $G\subset U_N^+$, and so we obtain $G\subset G^+(X)$, as claimed.
\end{proof}

We have now all the needed ingredients, and we can prove:

\begin{theorem}
Assuming that an algebraic manifold $X\subset S^{N-1}_{\mathbb C,+}$ is $I$-affine, with associated tensor category $F$, the following happen:
\begin{enumerate}
\item We have an inclusion $G\subset G^+(X)$, where $G$ is the Tannakian dual of $F$.

\item $X$ is an affine homogeneous space, $X=X_{G,I}$, over this quantum group $G$.
\end{enumerate}
\end{theorem}

\begin{proof}
In the context of Definition 8.8, the tensor category $F$ there gives rise, by the Tannakian duality of Woronowicz \cite{wo2}, to a quantum group $G\subset U_N^+$. What is left is to construct the affine space morphisms $\alpha,\Phi$, and the proof here goes as follows:

(1) Construction of $\alpha$. We want to construct a morphism, as follows:
$$\alpha:C(X)\to C(G)\quad,\quad x_i\to X_i=\frac{1}{\sqrt{|I|}}\sum_{j\in I}u_{ij}$$

In view of Definition 8.8, we must therefore prove that we have:
$$\sum_{i_1\ldots i_l}\sum_{j_1\ldots j_k}T_{i_1\ldots i_l,j_1\ldots j_k}X_{i_1}^{e_1}\ldots X_{i_l}^{e_l}(X_{j_1}^{f_1}\ldots X_{j_k}^{f_k})^*=\frac{1}{\sqrt{|I|^{k+l}}}\sum_{b_1\ldots b_l\in I}\sum_{c_1\ldots c_k\in I}T_{b_1\ldots b_l,c_1\ldots c_k}$$

By replacing the variables $X_i$ by their above values, we want to prove that:
$$\sum_{i_1\ldots i_l}\sum_{j_1\ldots j_k}\sum_{r_1\ldots r_l\in I}\sum_{s_1\ldots s_k\in I}T_{i_1\ldots i_l,j_1\ldots j_k}u_{i_1r_1}^{e_1}\ldots u_{i_lr_l}^{e_l}(u_{j_1s_1}^{f_1}\ldots u_{j_ks_k}^{f_k})^*=\sum_{b_1\ldots b_l\in I}\sum_{c_1\ldots c_k\in I}T_{b_1\ldots b_l,c_1\ldots c_k}$$

Now observe that from the relation $T\in Hom(u^{\otimes k},u^{\otimes l})$ we obtain:
$$\sum_{i_1\ldots i_l}\sum_{j_1\ldots j_k}T_{i_1\ldots i_l,j_1\ldots j_k}u_{i_1r_1}^{e_1}\ldots u_{i_lr_l}^{e_l}(u_{j_1s_1}^{f_1}\ldots u_{j_ks_k}^{f_k})^*=T_{r_1\ldots r_l,s_1\ldots s_k}$$

Thus, by summing over indices $r_i\in I$ and $s_i\in I$, we obtain the desired formula.

\medskip

(2) Construction of $\Phi$. We want to construct a morphism, as follows:
$$\Phi:C(X)\to C(G)\otimes C(X)\quad,\quad x_i\to X_i=\sum_ju_{ij}\otimes x_j$$

But this is precisely the coaction map constructed in Proposition 8.9.

\medskip

(3) Proof of the ergodicity. If we go back to the general theory in chapter 7, we see that the ergodicty condition is equivalent to a number of Tannakian conditions, which are automatic in our case. Thus, the ergodicity condition is automatic, and we are done.
\end{proof}

The above result, based on the notion of $I$-affine manifold, remains quite theoretical. The problem is that Definition 8.8 still makes reference to a tensor category, and so the abstract characterization of the affine homogeneous spaces that we obtain in this way is not totally intrinsic. We believe that some deeper results should hold as well. To be more precise, the work on noncommutative spheres in \cite{bme} suggests that the relevant category $F$ should appear in a more direct way from $X$. Let us formulate:

\begin{definition}
Given a submanifold $X\subset S^{N-1}_{\mathbb C,+}$ and a subset $I\subset\{1,\ldots,N\}$, we let $F_{X,I}(k,l)\subset M_{N^l\times N^k}(\mathbb C)$ be the linear space of linear maps $T$ such that
$$\sum_{i_1\ldots i_l}\sum_{j_1\ldots j_k}T_{i_1\ldots i_l,j_1\ldots j_k}x_{i_1}^{e_1}\ldots x_{i_l}^{e_l}(x_{j_1}^{f_1}\ldots x_{j_k}^{f_k})^*=\frac{1}{\sqrt{|I|^{k+l}}}\sum_{b_1\ldots b_l\in I}\sum_{c_1\ldots c_k\in I}T_{b_1\ldots b_l,c_1\ldots c_k}$$
holds over $X$. We say that $X$ is $I$-saturated when 
$$F_{X,I}=(F_{X,l}(k,l))$$
is a tensor category, and the collection of the above relations presents $C(X)$.
\end{definition}

Observe that any $I$-saturated manifold is automatically $I$-affine. The results in \cite{bme} seem to suggest that the converse of this fact should hold. We do not have a proof of this fact, but we would like to present a few observations on this subject. First, we have:

\begin{proposition}
The linear spaces $F_{X,I}(k,l)\subset M_{N^l\times N^k}(\mathbb C)$ are as follows:
\begin{enumerate}
\item They contain the units.

\item They are stable by conjugation.

\item They satisfy the Frobenius duality condition.
\end{enumerate}
\end{proposition}

\begin{proof}
All these assertions are elementary, as follows:

\medskip

(1) Consider indeed the unit map. The associated relation is:
$$\sum_{i_1\ldots i_k}x_{i_1}^{e_1}\ldots x_{i_k}^{e_k}(x_{i_1}^{e_1}\ldots x_{i_k}^{e_k})^*=1$$

But this relation holds indeed, due to the defining relations for $S^{N-1}_{\mathbb C,+}$.

\medskip

(2) We have indeed the following sequence of equivalences:
\begin{eqnarray*}
&&T^*\in F_{X,I}(l,k)\\
&\iff&\sum_{i_1\ldots i_l}\sum_{j_1\ldots j_k}T^*_{j_1\ldots j_k,i_1\ldots i_l}x_{j_1}^{f_1}\ldots x_{j_k}^{f_k}(x_{i_1}^{e_1}\ldots x_{i_l}^{e_l})^*=\frac{1}{\sqrt{|I|^{k+l}}}\sum_{b_1\ldots i_l\in I}\sum_{c_1\ldots c_k\in I}T^*_{c_1\ldots c_k,b_1\ldots b_l}\\
&\iff&\sum_{i_1\ldots i_l}\sum_{j_1\ldots j_k}T_{i_1\ldots i_l,j_1\ldots j_k}x_{i_1}^{e_1}\ldots x_{i_l}^{e_l}(x_{j_1}^{f_1}\ldots x_{j_k}^{f_k})^*=\frac{1}{\sqrt{|I|^{k+l}}}\sum_{b_1\ldots b_l\in I}\sum_{c_1\ldots c_k\in I}T_{b_1\ldots b_l,c_1\ldots c_k}\\
&\iff&T\in F_{X,I}(k,l)
\end{eqnarray*}

(3) We have indeed a correspondence $T\in F_{X,I}(k,l)\ \leftrightarrow\ \xi\in F_{X,I}(\emptyset,l\bar{k})$, given by the usual formulae for the Frobenius isomorphism.
\end{proof}

Based on the above result, we can now formulate our observations, as follows:

\begin{theorem}
Given a closed subgroup $G\subset U_N^+$, and an index set $I\subset\{1,\ldots,N\}$, consider the corresponding affine homogeneous space $X_{G,I}\subset S^{N-1}_{\mathbb C,+}$.
\begin{enumerate}
\item $X_{G,I}$ is $I$-saturated precisely when the collection of spaces $F_{X,I}=(F_{X,I}(k,l))$ is stable under compositions, and under tensor products.

\item We have $F_{X,I}=F$ precisely when we have
\begin{eqnarray*}
&&\!\!\!\sum_{j_1\ldots j_l\in I}\big(\sum_{i_1\ldots i_l}\xi_{i_1\ldots i_l}u_{i_1j_1}^{e_1}\ldots u_{i_lj_l}^{e_l}-\xi_{j_1\ldots j_l}\big)=0\\
&\implies&\sum_{i_1\ldots i_l}\xi_{i_1\ldots i_l}u_{i_1j_1}^{e_1}\ldots u_{i_lj_l}^{e_l}-\xi_{j_1\ldots j_l}=0
\end{eqnarray*}
for any choice of the indices $j_1,\ldots,j_l$.
\end{enumerate}
\end{theorem}

\begin{proof}
We use the fact, from Theorem 8.7, that with $F(k,l)=Hom(u^{\otimes k},u^{\otimes l})$, we have inclusions of vector spaces $F(k,l)\subset F_{X,I}(k,l)$. Moreover, once again by Theorem 8.7, the relations coming from the elements of the category formed by the spaces $F(k,l)$ present $X_{G,I}$. Thus, the relations coming from the elements of $F_{X,I}$ present $X_{G,I}$ as well. With this observation in hand, our assertions follow from Proposition 8.12:

\medskip

(1) According to Proposition 8.12 (1,2) the unit and conjugation axioms are satisfied, so the spaces $F_{X,I}(k,l)$ form a tensor category precisely when the remaining axioms, namely the composition and the tensor product one, are satisfied. Now by assuming that these two axioms are satisfied, $X$ follows to be $I$-saturated, by the above observation.

\medskip

(2) Since we already have inclusions in one sense, the equality $F_{X,I}=F$ from the statement means that we must have inclusions in the other sense, as follows:
$$F_{X,I}(k,l)\subset F(k,l)$$

By using now Proposition 8.12 (3), it is enough to discuss the case $k=0$. And here, assuming that we have $\xi\in F_{X,L}(0,l)$, the following condition must be satisfied:
$$\sum_{i_1\ldots i_l}\xi_{i_1\ldots i_l}x_{i_1}^{e_1}\ldots x_{i_l}^{e_l}=\sum_{j_1\ldots j_l\in I}\xi_{j_1\ldots j_l}$$

By applying now the morphism $\alpha:C(X_{G,I})\to C(G)$, we deduce that we have:
$$\sum_{i_1\ldots i_l}\xi_{i_1\ldots i_l}\sum_{j_1\ldots j_l\in I}u_{i_1j_1}^{e_1}\ldots u_{i_lj_l}^{e_l}=\sum_{j_1\ldots j_l\in I}\xi_{j_1\ldots j_l}$$

Now recall that $F(0,l)=Fix(u^{\otimes l})$ consists of the vectors $\xi$ satisfying:
$$\sum_{i_1\ldots i_l}\xi_{i_1\ldots i_l}u_{i_1j_1}^{e_1}\ldots u_{i_lj_l}^{e_l}=\xi_{j_1\ldots j_l},\forall j_1,\ldots,j_l$$

We are therefore led to the conclusion in the statement.
\end{proof}

It is quite unclear on how to advance on these questions, and a more advanced algebraic trick, in the spirit of those used in \cite{bme}, seems to be needed. Nor is it clear on how to explicitely ``capture'' the relevant subgroup $G\subset G^+(X)$, in terms of our given manifold $X=X_{G,I}$, in a direct, geometric way. Summarizing, further improving Theorem 8.13 above is an interesting question, that we would like to raise here.

\bigskip

We will be back to such questions later on in this book, towards the end, when talking about the work in \cite{bme}. In fact, after a break in chapters 9-12 below, for talking about geometries other than classical and free, which are of interest too, we will be back to free geometry in the whole last part of this book, chapters 13-16 below.

\section*{8c. Row spaces}

We discuss in what follows some constructions from \cite{bss}, which go somehow in a direction which is opposite to what has been said in the above, namely particularization. The idea will be that of looking at the ``minimal'' theory of quantum homogeneous spaces generalizing at the same time the spheres $S$, and the unitary quantum groups $U$ they come form, and this time with very precise results. Such homogeneous spaces are technically covered by the general affine homogeneous space formalism from chapter 7, which is from the paper \cite{ba6}, which came some time after \cite{bss}, but the difference of generality level being notable, there are many things that can be said, sharper than in general.

\bigskip

We first discuss the construction in the classical case. Given a closed subgroup $G\subset U_N$ and a number $k\leq N$, we can consider the compact group $H=G\cap U_k$, computed inside $U_N$, where the embedding $U_k\subset U_N$ that we use is given by the following formula:
$$g\to
\begin{pmatrix}
g&0\\
0&1_{N-k}
\end{pmatrix}$$

We can form the homogeneous space $X=G/H$, and we have the following result:

\begin{proposition}
Let $G\subset U_N$ be a closed subgroup, and construct as above the closed subgroup $H\subset G$ given by the formula
$$H=G\cap U_k$$
with the intersection being computed inside $U_N$. Then the subalgebra
$$C(G/H)\subset C(G)$$
that we obtain is generated by the last $N-k$ rows of coordinates on $G$.
\end{proposition}

\begin{proof}
Let $u_{ij}\in C(G)$ be the standard coordinates on $G$, given as usual by the formula $u_{ij}(g)=g_{ij}$, and consider the following subalgebra of $C(G)$:
$$A=\left<u_{ij}\Big|i>k,j>0\right>$$

Since each coordinate function $u_{ij}$ with $i>k$ is constant on each coset $Hg\in G/H$, we have an inclusion as follows, between subalgebras of $C(G)$:
$$A\subset C(G/H)$$

In order to prove that this inclusion in a isomorphism, as to finish, we use the Stone-Weierstrass theorem. Indeed, in view of this theorem, it is enough to show that the following family of functions separates the cosets $\{Hg|g\in G\}$:
$$\left\{u_{ij}\Big|i>k,j>0\right\}$$

But this is the same as saying that $Hg\neq Hh$ implies $g_{ij}\neq h_{ij}$, for some $i>k,j>0$. Equivalently, we must prove that $g_{ij}=h_{ij}$ for any $i>k,j>0$ implies:
$$Hg=Hh$$

Now since $Hg=Hh$ is equivalent to $gh^{-1}\in H$, the result follows from the usual matrix formula of $gh^{-1}$, and from the fact that $g,h$ are unitary.
\end{proof}

In the quantum case now, we can proceed in a similar way. Let $k\leq N$, and consider the embedding $U_k^+\subset U_N^+$ given by the same formula as before, namely:
$$g\to
\begin{pmatrix}
g&0\\
0&1_{N-k}
\end{pmatrix}$$

That is, at the level of algebras, we use the quotient map $C(U_N^+)\to C(U_k^+)$ given by the following formula, where $v$ is the fundamental corepresentation of $U_k^+$:
$$u\to 
\begin{pmatrix}
v&0\\
0&1_{N-k}
\end{pmatrix}$$

With this convention, we have the following definition, from \cite{bss}:

\index{row space}
\index{row algebra}

\begin{definition}
Associated to any quantum subgroup $G\subset U_N^+$ and any $k\leq N$ are:
\begin{enumerate}
\item The compact quantum group $H=G\cap U_k^+$.

\item The algebra $C(G/H)\subset C(G)$ constructed before.

\item The algebra $C_\times(G/H)\subset C(G/H)$ generated by $\{u_{ij}|i>k,j>0\}$.
\end{enumerate}
\end{definition}

Regarding (3), let $u,v$ be the fundamental corepresentations of $G,H$, so that the quotient map $\pi:C(G)\to C(H)$ is given by $u\to diag(v,1_{N-k})$. We have then:
\begin{eqnarray*}
(\pi\otimes id)\Delta(u_{ij})
&=&\sum_s\pi(u_{is})\otimes u_{sj}\\
&=&\begin{cases}
\sum_{s\leq k}v_{is}\otimes u_{sj}&i\leq k\\
1\otimes u_{ij}&i>k
\end{cases}
\end{eqnarray*}

In particular we see that the equality $(\pi\otimes id)\Delta f=1\otimes f$ defining $C(G/H)$ holds on all the coefficients $f=u_{ij}$ with $i>k$, and this justifies the inclusion appearing in (3).

\bigskip

Let us first try to understand what happens in the group dual case. We will do our study here in two steps, first in the ``diagonal'' case, and then in the general case. We recall that given a discrete group $\Gamma=<g_1,\ldots,g_N>$, the matrix $D=diag(g_i)$ is biunitary, and produces a surjective morphism $C(U_N^+)\to C^*(\Gamma)$. This morphism can be viewed as corresponding to a quantum embedding $\widehat{\Gamma}\subset U_N^+$, that we call ``diagonal''.

\bigskip

We recall also that the normal closure of a subgroup $\Lambda\subset\Gamma$ is the biggest subgroup $\Lambda'\subset\Gamma$ containing $\Lambda$ as a normal subgroup. Note that $\Lambda'$ can be different from the normalizer $N(\Lambda)$. With these conventions, we have the following result, from \cite{bss}:

\begin{proposition}
Assume that we have a group dual $G=\widehat{\Gamma}$, with 
$$\Gamma=<g_1,\ldots,g_N>$$
diagonally embedded, and let $H=G\cap U_k^+$.
\begin{enumerate}
\item $H=\widehat{\Theta}$, where $\Theta=\Gamma/<g_{k+1}=1,\ldots,g_N=1>$.

\item $C_\times(G/H)=C^*(\Lambda)$, where $\Lambda=<g_{k+1},\ldots,g_N>$.

\item $C(G/H)=C^*(\Lambda')$, where ``prime'' is the normal closure.

\item $C_\times(G/H)=C(G/H)$ if and only if $\Lambda\triangleleft\Gamma$.
\end{enumerate}
\end{proposition}

\begin{proof}
We use the standard fact that for any group $\Gamma=<a_i,b_j>$, the kernel of the quotient map $\Gamma\to\Gamma/<a_i=1>$ is the normal closure of the subgroup $<a_i>\subset\Gamma$.

\medskip

(1) Since the map $C(U_N^+)\to C(U_k^+)$ is given on diagonal coordinates by $u_{ii}\to v_{ii}$ for $i\leq k$ and $u_{ii}\to 1$ for $i>k$, the result follows from definitions.

\medskip

(2) Once again, this assertion follows from definitions.

\medskip

(3) From the above and from (1) we get $G/H=\widehat{\Lambda'}$, where $\Lambda'=\ker(\Gamma\to\Theta)$. By the above observation, this kernel is exactly the normal closure of $\Lambda$.

\medskip

(4) This follows from (2) and (3).
\end{proof}

Let us try now to understand the general group dual case. We recall that the group dual subgroups $\widehat{\Gamma}\subset U_N^+$ appear by taking a discrete group $\Gamma=<g_1,\ldots,g_N>$ and a unitary $J\in U_N$, and constructing the morphism $C(U_N^+)\to C^*(\Gamma)$ given by $u\to JDJ^*$, where $D=diag(g_i)$. With this in hand, Proposition 8.16 generalizes as follows:

\begin{theorem}
Assume that we have a group dual $G=\widehat{\Gamma}$, with 
$$\Gamma=<g_1,\ldots,g_N>$$
embedded via $u\to JDJ^*$, and let $H=G\cap U_k^+$.
\begin{enumerate}
\item $H=\widehat{\Theta}$, where $\Theta=\Gamma/<g_r=1|\exists\,i>k,J_{ir}\neq 0>$, embedded $u_{ij}\to (JDJ^*)_{ij}$.

\item $C_\times(G/H)=C^*(\Lambda)$, where $\Lambda=<g_r|\exists\,i>k,J_{ir}\neq 0>$.

\item $C(G/H)=C^*(\Lambda')$, where ``prime'' is the normal closure.

\item $C_\times(G/H)=C(G/H)$ if and only if $\Lambda\triangleleft\Gamma$.
\end{enumerate}
\end{theorem}

\begin{proof}
We basically follow the proof of Proposition 8.16:

\medskip

(1) Let $\Lambda=<g_1,\ldots,g_N>$, let $J\in U_N$, and consider the embedding $\widehat{\Lambda}\subset U_N^+$ corresponding to the following morphism, where $D=diag(g_i)$:
$$C(U_N^+)\to C^*(\Lambda)\quad,\quad u\to JDJ^*$$

Let $G=\widehat{\Lambda}\cap U_k^+$. Since we have $G\subset\widehat{\Lambda}$, the algebra $C(G)$ is cocommutative, so we have $G=\widehat{\Theta}$ for a certain discrete group $\Theta$. Moreover, the inclusion $\widehat{\Theta}\subset\widehat{\Lambda}$ must come from a group morphism $\varphi:\Lambda\to\Theta$. Also, since $\widehat{\Theta}\subset U_k^+$, we have a morphism as follows, where $V$ is a certain $k\times k$ biunitary matrix over the algebra $C^*(\Theta)$:
$$C(U_k^+)\to C^*(\Theta)\quad,\quad v\to V$$

With these observations in hand, let us look now at the intersection operation. We must have a group morphism $\varphi:\Lambda\to\Theta$ such that the following diagram commutes:
$$\xymatrix@R=50pt@C=50pt
{C(U_N^+)\ar[r]\ar[d]&C(U_k^+)\ar[d]\\
C^*(\Lambda)\ar[r]&C^*(\Theta)}$$

Thus we must have the following equality:
$$(id\otimes\varphi)(JDJ^*)=diag(V,1_{N-k})$$

With $f_i=\varphi(g_i)$, we obtain from this:
$$\sum_rJ_{ir}\bar{J}_{jr}f_r
=\begin{cases}
V_{ij}&\mbox{if }i,j\leq k\\
\delta_{ij}&\mbox{otherwise}
\end{cases}$$

Now since $J$ is unitary, the second part of the above condition is equivalent to ``$f_r=1$ whenever there exists $i>k$ such that $J_{ir}\neq 0$''. Indeed, this condition is easily seen to be equivalent to the ``$=1$'' conditions, and implies the ``$=0$'' conditions. We claim that:
$$\Theta=\Lambda\Big/\left<g_r=1\Big|\exists\,i>k,J_{ir}\neq 0\right>$$

Indeed, the above discussion shows that $\Theta$ must be a quotient of the group on the right, say $\Theta_0$. On the other hand, since in $C^*(\Theta_0)$ we have $J_{ir}g_r=J_{ir}1$ for any $i>k$, we obtain that $(JDJ^*)_{ij}=\delta_{ij}$ unless $i,j\leq k$, so we have, for a certain matrix $V$:
$$JDJ^*=diag(V,1_{N-k})$$

But the matrix $V$ must be a biunitary, so we have a morphism $C(U_k^+)\to C^*(\Theta_0)$ mapping $v\to V$, which completes the proof of our claim.

\medskip

(2) Consider the standard generators of the algebra $C_\times(G/H)$ constructed in  Definition 8.15 (3), which are as follows, with indices $i>k,j>0$:
$$A_{ij}=\sum_rJ_{ir}\bar{J}_{jr}g_r$$

We have then the following formula:
$$\sum_jA_{ij}J_{jm}=J_{im}g_m$$

We conclude that $C_\times(G/H)$ contains any $g_r$ such that there exists $i>k$ with $J_{ir}\neq 0$, i.e. contains any $g_r\in\Lambda$. Conversely, if $g_r\in\Gamma-\Lambda$ then $J_{ir}g_r=0$ for any $i>k$, so $g_r$ doesn't appear in the formula of any of the generators $A_{ij}$.

\medskip

(3,4) The proof here is similar to the proof of Proposition 8.16 (3,4).
\end{proof}

Summarizing, we have a good understanding of the row algebras for the compact quantum groups, both in the classical case, and in the group dual case.

\section*{8d. Uniformity}

Following \cite{bss}, we discuss in what follows the structure of the row algebras in the case where the underlying quantum group is easy, which is the case that we are mostly interested in. As in \cite{bss}, which was written some time ago, and based on \cite{bsp} dealing with the orthogonal case, $G\subset O_N^+$, we will restrict the attention to the orthogonal case. With the remark of course that the unitary extension looks quite straightforward. We will need the following key result, coming from \cite{bss}, \cite{bsp}:

\index{uniform quantum group}
\index{removing blocks}
\index{category of partitions}

\begin{theorem}
For an easy subgroup $G_N\subset O_N^+$, the following are equivalent:
\begin{enumerate}
\item $G=(G_N)$ is uniform, in the sense that we have
$$G_N\cap O_k^+=G_k$$
for any $k\leq N$, with respect to the standard embedding $O_k^+\subset O_N^+$.

\item The corresponding category of partitions
$$D=(D(k,l))$$
is stable under the operation which consists in removing blocks.
\end{enumerate}
\end{theorem}

\begin{proof}
This can proved in several steps, as follows:

\medskip

(1) In order to establish the equivalence between the above two conditions, we will prove that $G_N\cap O_k^+=G_k'$, where $G'=(G_N')$ is the easy quantum group associated to the category $D'$ generated by all subpartitions of the partitions in $D$. 

\medskip

(2) We know that the correspondence between categories of partitions and easy quantum groups comes from Woronowicz's Tannakian duality in \cite{wo2}, with the quantum group $G_N\subset O_N^+$ associated to a category of partitions $D=(D(s))$ obtained by imposing to the fundamental representation of $O_N^+$ the fact that its $s$-th tensor power must fix $\xi_\pi$, for any $s\in\mathbb N$ and $\pi\in D(s)$. In other words, we have the following formula:
$$C(G_N)=C(O_N^+)\Big/\Big<\xi_\pi\in Fix(u^{\otimes s}),\forall s,\,\forall\pi\in D(s)\Big>$$

Now since $\xi_\pi\in Fix(u^{\otimes s})$ means $u^{\otimes s}(\xi_\pi\otimes 1)=\xi_\pi\otimes 1$, this condition is equivalent to the following collection of equalities, one for each multi-index $i\in\{1,\ldots,N\}^s$:
$$\sum_{j_1\ldots j_s}\delta_\pi(j)u_{i_1j_1}\ldots u_{i_sj_s}=\delta_\pi(i)1$$

Summarizing, we have the following presentation result:
$$C(G_N)=C(O_N^+)\Big/\left<\sum_{j_1\ldots j_s}\delta_\pi(j)u_{i_1j_1}\ldots u_{i_sj_s}=\delta_\pi(i)1,\forall s,\,\forall\pi\in D(s),\,\forall i\right>$$

(3) Let now $k\leq N$, assume that we have a compact quantum group $K\subset O_k^+$, with fundamental representation denoted $u$, and consider the following $N\times N$ matrix:
$$\tilde{u}=\begin{pmatrix}u&0\\ 0&1_{N-k}\end{pmatrix}$$

Our claim is that for any $s\in\mathbb N$ and any $\pi\in P(s)$, we have:
$$\xi_\pi\in Fix(\tilde{u}^{\otimes s})\iff\xi_{\pi'}\in Fix(u^{\otimes s'}),\,\forall\pi'\subset\pi$$

Here $\pi'\subset\pi$ means that $\pi'\in P(s')$ is obtained from $\pi\in P(s)$ by removing some of its blocks. The proof of this claim is standard. Indeed, when making the replacement $u\to\tilde{u}$ and trying to check the condition $\xi_\pi\in Fix(\tilde{u}^{\otimes s})$, we have two cases:

\medskip

-- $\delta_\pi(i)=1$. Here the $>k$ entries of $i$ must be joined by certain blocks of $\pi$, and we can consider the partition $\pi'\in D(s')$ obtained by removing these blocks. The point now is that the collection of $\delta_\pi(i)=1$ equalities to be checked coincides with the collection of $\delta_\pi(i)=1$ equalities expressing the fact that we have $\xi_\pi\in Fix(u^{\otimes s'})$, for any $\pi'\subset\pi$.

\medskip

-- $\delta_\pi(i)=0$. In this case the situation is quite similar. Indeed, the collection of $\delta_\pi(i)=0$ equalities to be checked coincides, modulo some $0=0$ identities, which hold automatically, with the collection of $\delta_\pi(i)=0$ equalities expressing the fact that we have $\xi_\pi\in Fix(u^{\otimes s'})$, for any $\pi'\subset\pi$.

\medskip

(4) Our second claim is that given a quantum group $K\subset O_N^+$, with fundamental representation denoted $v$, the algebra of functions on $H=K\cap O_k^+$ is given by:
$$C(H)=C(O_k^+)\Big/\Big<\xi\in Fix(\tilde{u}^{\otimes s}),\,\forall\xi\in Fix(v^{\otimes s})\Big>$$

But this follows indeed from Woronowicz's results in \cite{wo2}, because the algebra on the right comes from the Tannakian formulation of the intersection operation.

\medskip

(5) Now with the above two claims in hand,  we can conclude that we have the following formula, where $G'=(G_N')$ is the easy quantum group associated to the category $D'$ generated by all the subpartitions of the partitions in $D$:
$$G_N\cap U_k^+=G_k'$$

In particular we see that the condition $G_N\cap U_k^+=G_k^+$ for any $k\leq N$ is equivalent to $D=D'$, and this gives the result.
\end{proof}

Let us study now the following inclusions of algebras, constructed in Definition 8.15, where $G=(G_n)$ is a uniform easy quantum group:
$$C_\times(G_N/G_k)\subset C(G_N/G_k)$$

For classification purposes the uniformity axiom is something very natural and useful, substantially cutting from complexity, and we have the following result, from \cite{bsp}:

\begin{theorem}
The classical and free uniform orthogonal easy quantum groups, with inclusions between them, are as follows:
$$
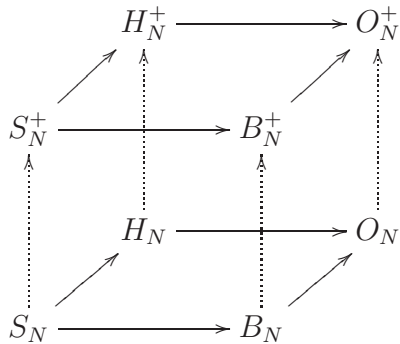
\xymatrix@R=20pt@C=20pt{
&H_N^+\ar[rr]&&O_N^+\\
S_N^+\ar[rr]\ar[ur]&&B_N^+\ar[ur]\\
&H_N\ar[rr]\ar[uu]&&O_N\ar[uu]\\
S_N\ar[uu]\ar[ur]\ar[rr]&&B_N\ar[uu]\ar[ur]
}$$
Moreover, this is an intersection/easy generation diagram, in the sense that for any of its square subdiagrams $P\subset Q,R\subset S$ we have $P=Q\cap R$ and $\{Q,R\}=S$.
\end{theorem}

\begin{proof}
In this statement all the quantum groups are objects that we are familiar with, and that we know to be easy, except for $B_N\subset O_N$ and $B_N^+\subset O_N^+$, which are the bistochastic group and its free analogue, constructed via the relation $\xi\in Fix(u)$, where $\xi$ is the all-one vector. Since this all-one vector corresponds to the singleton partition, the quantum groups $B_N,B_N^+$ follow to be easy too, coming from the categories $P_{12},NC_{12}$ of singletons and pairings. Thus, the quantum groups in the statement are all easy, and clearly uniform too, the corresponding categories of partitions being as follows:
$$
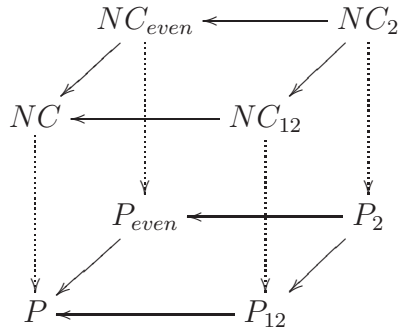
\xymatrix@R=21pt@C=6pt{
&NC_{even}\ar[dl]\ar[dd]&&NC_2\ar[dl]\ar[ll]\ar[dd]\\
NC\ar[dd]&&NC_{12}\ar[dd]\ar[ll]\\
&P_{even}\ar[dl]&&P_2\ar[dl]\ar[ll]\\
P&&P_{12}\ar[ll]
}$$

Since this latter diagram is an intersection and generation diagram, we conclude that we have an intersection and easy generation diagram of quantum groups, as stated.

\medskip

Regarding now the classification, consider an easy quantum group $S_N\subset G_N\subset O_N$. This most come from a category $P_2\subset D\subset P$, and if we assume $G=(G_N)$ to be uniform, then $D$ is uniquely determined by the subset $L\subset\mathbb N$ consisting of the sizes of the blocks of the partitions in $D$. Our claim is that the admissible sets are as follows:

\begin{enumerate}
\item $L=\{2\}$, producing $O_N$.

\smallskip

\item $L=\{1,2\}$, producing $B_N$.

\smallskip

\item $L=\{2,4,6,\ldots\}$, producing $H_N$.

\smallskip

\item $L=\{1,2,3,\ldots\}$, producing $S_N$.
\end{enumerate}

In one sense, this follows indeed from our easiness results for $O_N,B_N,H_N,S_N$. In the other sense now, assume that $L\subset\mathbb N$ is such that the set $P_L$ consisting of partitions whose sizes of the blocks belong to $L$ is a category of partitions. We know from the axioms of the categories of partitions that the semicircle $\cap$ must be in the category, so we have $2\in L$. Our claim is that the following conditions must be satisfied as well:
$$k,l\in L,\,k>l\implies k-l\in L$$
$$k\in L,\,k\geq 2\implies 2k-2\in L$$

Indeed, we will prove that both conditions follow from the axioms of the categories of
partitions. Let us denote by $b_k\in P(0,k)$ the one-block partition, namely:
$$b_k=\left\{\begin{matrix}\sqcap\hskip-0.7mm \sqcap&\ldots&\sqcap\\
1\hskip2mm 2&\ldots&k\end{matrix} \right\}$$

For $k>l$, we can write the partition $b_{k-l}$ in the following way:
$$b_{k-l}=\left\{\begin{matrix}\sqcap\hskip-0.7mm
\sqcap&\ldots&\ldots&\ldots&\ldots&\sqcap\\ 1\hskip2mm 2&\ldots&l&l+1&\ldots&k\\
\sqcup\hskip-0.7mm \sqcup&\ldots&\sqcup&|&\ldots&|\\ &&&1&\ldots&k-l\end{matrix}\right\}$$

In other words, we have the following formula:
$$b_{k-l}=(b_l^*\otimes |^{\otimes k-l})b_k$$

Since all the terms of this composition are in $P_L$, we have $b_{k-l}\in P_L$, and this proves our first claim. As for the second claim, this can be proved in a similar way, by capping two adjacent $k$-blocks with a $2$-block, in the middle. Now, we can conclude as follows:

\medskip

\underline{Case 1}. Assume $1\in L$. By using the first condition with $l=1$ we get:
$$k\in L\implies k-1\in L$$

This condition shows that we must have $L=\{1,2,\ldots,m\}$, for a certain number $m\in\{1,2,\ldots,\infty\}$. On the other hand, by using the second condition we get:
\begin{eqnarray*}
m\in L
&\implies&2m-2\in L\\
&\implies&2m-2\leq m\\
&\implies&m\in\{1,2,\infty\}
\end{eqnarray*}

The case $m=1$ being excluded by the condition $2\in L$, we reach to one of the two sets producing the groups $S_N,B_N$.

\medskip

\underline{Case 2}. Assume $1\notin L$. By using the first condition with $l=2$ we get:
$$k\in L\implies k-2\in L$$

This condition shows that we must have $L=\{2,4,\ldots,2p\}$, for a certain number $p\in\{1,2,\ldots,\infty\}$. On the other hand, by using the second condition we get:
\begin{eqnarray*}
2p\in L
&\implies&4p-2\in L\\
&\implies&4p-2\leq 2p\\
&\implies&p\in\{1,\infty\}
\end{eqnarray*}

Thus $L$ must be one of the two sets producing $O_N,H_N$, and we are done. In the free case, $S_N^+\subset G_N\subset O_N^+$, the situation is quite similar, the admissible sets being once again the above ones, producing this time $O_N^+,B_N^+,H_N^+,S_N^+$. See \cite{bsp}. 
\end{proof}

Let us go back now to the inclusions $C_\times(G_N/G_k)\subset C(G_N/G_k)$. Following \cite{bss}, we first work out a few simple cases, where these inclusions are isomorphisms:

\begin{proposition}
For the basic easy quantum groups, the inclusion of algebras
$$C_\times(G_N/G_k)\subset C(G_N/G_k)$$ 
is an isomorphism at $N=1$, at $k=0$, at $k=N$, as well as in the following cases:
\begin{enumerate}
\item $G=B^+$: at $k=1$.

\item $G=S^+$: at $k=1$, and at $k=2,N=3$.
\end{enumerate}
\end{proposition}

\begin{proof}
First, the results at $N=1$, at $k=0$, and at $k=N$ are clear from definitions. Regarding now the special cases, the situation here is as follows:

\medskip

(1) Since the coordinates of $B_N^+$ sum up to 1 on each column, we have:
$$u_{1j}=1-\sum_{i>1}u_{ij}$$

Thus the following inclusion is an isomorphism:
$$C_\times(B_N^+/B_1^+)\subset C(B_N^+)$$

Thus the inclusion $C_\times(B_N^+/B_1^+)\subset C(B_N^+/B_1^+)$ must be as well an isomorphism.

\medskip

(2) By using the same argument as above we obtain that the following inclusion is as well an isomorphism:
$$C_\times(S_N^+/S_1^+)\subset C(S_N^+/S_1^+)$$

In the remaining case $k=2,N=3$, or more generally at any $k\in\mathbb N$ and $N<4$, it is known from Wang \cite{wa2} that we have $S_N=S_N^+$, so the inclusion in the statement is:
$$C(S_N/S_k)\subset C(S_N/S_k)$$

Thus, in this case we are done again.
\end{proof}

The axiomatization of the algebras $C_\times(G_N/G_k)$ is a quite tricky task. However, following \cite{bss}, we can axiomatize some bigger algebras, as follows:

\begin{definition}
Associated to $k\leq N$ is the universal $C^*$-algebra $C_+(G_N/G_k)$ generated by the entries of a matrix $p=(p_{ij})_{i>k,j>0}$ subject to the following conditions:
\begin{enumerate}
\item $G=O_N^+$: $p$ is a transposed ``orthogonal isometry'', in the sense that its entries $p_{ij}$ are self-adjoint, and $pp^t=1$.

\item $G=S_N^+$: $p$ is a transposed ``magic isometry'', in the sense that $p^t$ is an orthogonal isometry, and $p_{ij}$ are projections, orthogonal on columns.

\item $G=H_N^+$: $p$ is a transposed ``cubic isometry'', in the sense that $p^t$ is an orthogonal isometry, with $xy=0$ for any $x\neq y$ on the same row of $p$

\item $G=B_N^+$: $p$ is a transposed ``stochastic isometry'', in the sense that $p^t$ is an orthogonal isometry, with sum $1$ on rows.
\end{enumerate}
\end{definition}

Observe that we have surjective maps, as follows:
$$C_+(G_N/G_k)\to C_\times(G_N/G_k)$$

Still following \cite{bss}, we have the following result:

\begin{theorem}
The algebras $C_+(G_N/G_k)$ and $C_\times(G_N/G_k)$ are as follows:
\begin{enumerate}
\item They have coactions of $G_N$, given by $\alpha(p_{ij})=\sum_sp_{is}\otimes u_{sj}$.

\item They have unique $G_N$-invariant states, which are tracial.

\item Their reduced algebra versions are isomorphic.

\item Their abelianized versions are isomorphic.
\end{enumerate}
\end{theorem}

\begin{proof}
This is something quite long, the idea being as follows:

\medskip

(1) For the algebra $C_\times(G_N/G_k)$ this is clear, because as explained in \cite{bss}, this algebra is ``embeddable'', and the coaction of $G_N$ is simply the restriction of the comultiplication map. For the algebra $C_+(G_N/G_k)$, consider the following elements:
$$P_{ij}=\sum_{s=1}^Np_{is}\otimes u_{sj}$$

We have to check that these elements satisfy the same relations as those in Definition 8.21, presenting the algebra $C_+(G_n/G_k)$, and the proof here goes as follows:

\medskip

\underline{$O^+$ case}. First, since $p_{ij},u_{ij}$ are self-adjoint, so is $P_{ij}$. Also, we have:
\begin{eqnarray*}
\sum_jP_{ij}P_{rj}
&=&\sum_{jst}p_{is}p_{rt}\otimes u_{sj}u_{tj}\\
&=&\sum_{st}p_{is}p_{rt}\otimes\delta_{st}\\
&=&\sum_sp_{is}p_{rs}\otimes 1\\
&=&\delta_{ir}
\end{eqnarray*}

\underline{$H^+$ case}. The condition $xy=0$ on rows is checked as follows ($j\neq r$):
\begin{eqnarray*}
P_{ij}P_{ir}
&=&\sum_{st}p_{is}p_{it}\otimes u_{sj}u_{tr}\\
&=&\sum_sp_{is}\otimes u_{sj}u_{sr}\\
&=&0
\end{eqnarray*}

\underline{$B^+$ case}. The sum 1 condition on rows is checked as follows:
\begin{eqnarray*}
\sum_jP_{ij}
&=&\sum_{js}p_{is}\otimes u_{sj}\\
&=&\sum_sp_{is}\otimes 1\\
&=&1
\end{eqnarray*}

\underline{$S^+$ case}. Since $P^t$ is cubic and stochastic, we just check the projection condition:
\begin{eqnarray*}
P_{ij}^2
&=&\sum_{st}p_{is}p_{it}\otimes u_{sj}u_{tj}\\
&=&\sum_sp_{is}\otimes u_{sj}\\
&=&P_{ij}
\end{eqnarray*}

Summmarizing, the matrix $P$ satisfies the same conditions as $p$, so we can define a morphism of $C^*$-algebras, as follows:
$$\alpha:C_+(G_N/G_k)\to C_+(G_N/G_k)\otimes C(G_N)\quad,\quad 
\alpha(p_{ij})=P_{ij}$$

We have the following computation:
\begin{eqnarray*}
(\alpha\otimes id)\alpha (p_{ij})
&=&\sum_s\alpha(p_{is})\otimes u_{sj}\\
&=&\sum_{st}p_{it}\otimes u_{ts}\otimes u_{sj}
\end{eqnarray*}

On the other hand, we have as well the following computation:
\begin{eqnarray*}
(id\otimes\Delta)\alpha(p_{ij})
&=&\sum_tp_{it}\otimes\Delta(u_{ij})\\
&=&\sum_{st}p_{it}\otimes u_{ts}\otimes u_{sj}
\end{eqnarray*}

Thus our map $\alpha$ is coassociative. The density conditions can be checked by using dense subalgebras generated by $p_{ij}$ and $u_{st}$, and we are done.

\medskip

(2) For the existence part we can use the following composition, where the first two maps are the canonical ones, and the map on the right is the integration over $G_N$:
$$C_+(G_N/G_k)\to C_\times(G_N/G_k)\subset C(G_N)\to\mathbb C$$

Also, the uniqueness part is clear for the algebra $C_\times(G_N/G_k)$. Regarding now the uniqueness for $C_+(G_N/G_k)$, let $\int$ be the Haar state on $G_N$, and $\varphi$ be the $G_N$-invariant state constructed above. We claim that $\alpha$ is ergodic:
$$\left(id\otimes\int\right)\alpha=\varphi(.)1$$

Indeed, let us recall that the Haar state is given by the following Weingarten formula, where $W_{sN}=G_{sN}^{-1}$, with $G_{sN}(\pi,\sigma)=N^{|\pi\vee\sigma|}$:
$$\int u_{i_1j_1}\ldots u_{i_sj_s}=\sum_{\pi,\sigma\in D(s)}\delta_\pi(i)\delta_\sigma(j)W_{sN}(\pi,\sigma)$$

Now, let us go back now to our claim. By linearity it is enough to check the above equality on a product of basic generators $p_{i_1j_1}\ldots p_{i_sj_s}$. The left term is as follows:
\begin{eqnarray*}
\left(id\otimes\int\right)\alpha(p_{i_1j_1}\ldots p_{i_sj_s})
&=&\sum_{l_1\ldots l_s}p_{i_1l_1}\ldots p_{i_sl_s}\int u_{l_1j_1}\ldots u_{l_sj_s}\\
&=&\sum_{l_1\ldots l_s}p_{i_1l_1}\ldots p_{i_sl_s}\sum_{\pi,\sigma\in D(s)}\delta_\pi(l)\delta_\sigma(j)W_{sN}(\pi,\sigma)\\
&=&\sum_{\pi,\sigma\in D(s)}\delta_\sigma(j)W_{sN}(\pi,\sigma)\sum_{l_1\ldots l_s}\delta_\pi(l)p_{i_1l_1}\ldots p_{i_sl_s}
\end{eqnarray*}

Let us look now at the sum on the right. We have to sum the elements of type $p_{i_1l_1}\ldots p_{i_sl_s}$, over all multi-indices $l=(l_1,\ldots,l_s)$ which fit into our partition $\pi\in D(s)$. In the case of a one-block partition this sum is simply $\sum_lp_{i_1l}\ldots p_{i_sl}$, and we claim that:
$$\sum_lp_{i_1l}\ldots p_{i_sl}=\delta_\pi(i)$$

Indeed, the proof of this formula goes as follows:

\medskip

\underline{$O^+$ case}. Here our one-block partition must be a semicircle, $\pi=\cap$, and the formula to be proved, namely $\sum_lp_{il}p_{jl}=\delta_{ij}$, follows from $pp^t=1$.

\medskip

\underline{$S^+$ case}. Here our one-block partition can be any $s$-block, $1_s\in P(s)$, and the formula to be proved, namely $\sum_lp_{i_1l}\ldots p_{i_sl}=\delta_{i_1,\ldots,i_s}$, follows from orthogonality on columns, and from the fact that the sum is 1 on rows.

\medskip

\underline{$B^+$ case}. Here our one-block partition must be a semicircle or a singleton. We are already done with the semicircle, and for the singleton the formula to be proved, namely $\sum_lp_{il}=1$, follows from the fact that the sum is 1 on rows.

\medskip

\underline{$H^+$ case}. Here our one-block partition must have an even number of legs, $s=2r$, and due to the cubic condition the formula to be proved reduces to $\sum_lp_{il}^{2r}=1$. But since $p_{il}^{2r}=p_{il}^2$, independently on $r$, the result follows from the orthogonality on rows.

\medskip

In the general case now, since $\pi$ noncrossing, the computations over the blocks will not interfere, and we will obtain the same result, namely:
$$\sum_lp_{i_1l}\ldots p_{i_sl}=\delta_\pi(i)$$

Now by plugging this formula into the computation that we have started, we get:
\begin{eqnarray*}
\left(id\otimes\int\right)\alpha(p_{i_1j_1}\ldots p_{i_sj_s})
&=&\sum_{\pi,\sigma\in D(s)}\delta_\pi(i)\delta_\sigma(j)W_{sN}(\pi,\sigma)\\
&=&\int u_{i_1j_1}\ldots u_{i_sj_s}\\
&=&\varphi(p_{i_1j_1}\ldots p_{i_sj_s})
\end{eqnarray*}

This finishes the proof of our claim. So, let us get back now to the original question. Let $\tau:C_+(G_N/G_k)\to\mathbb C$ be a linear form as in the statement. We have:
\begin{eqnarray*}
\tau\left(id\otimes\int\right)\alpha(x)
&=&\left(\tau\otimes\int\right)\alpha(x)\\
&=&\int(\tau\otimes id)\alpha(x)\\
&=&\int(\tau(x)1)\\
&=&\tau(x)
\end{eqnarray*}

On the other hand, according to our above claim, we have as well:
$$\tau\left(id\otimes\int\right)\alpha(x)
=\tau(\varphi(x)1)
=\varphi(x)$$

Thus we get $\tau=\varphi$, which finishes the proof of the uniqueness assertion.

\medskip

(3) This follows from the uniqueness assertions in (2), and from some standard facts regarding the reduced versions with respect to Haar states, from Woronowicz \cite{wo1}.

\medskip

(4) We denote by $G^-$ the classical version of $G$, given by $G^-=O,S,H,B$ in the cases $G=O^+,S^+,H^+,B^+$. We have surjective morphisms of algebras, as follows:
\begin{eqnarray*}
C_+(G_N/G_k)
&\to&C_\times(G_k/G_k)\\
&\to&C_\times(G_N^-/G_k^-)\\
&=&C(G_N^-/G_k^-)
\end{eqnarray*}

Thus at the level of abelianized versions, we have surjective morphisms as follows:
\begin{eqnarray*}
C_+(G_N/G_k)_{comm}
&\to&C_\times(G_N/G_k)_{comm}\\
&\to&C(G_N^-/G_k^-)
\end{eqnarray*}

In order to prove our claim, namely that the first surjective morphism is an isomorphism, it is enough to prove that the above composition is an isomorphism.

\medskip

Let $r=N-k$, and denote by $A_{N,r}$ the algebra on the left. This is by definition the algebra generated by the entries of a transposed $N\times r$ isometry, whose entries commute, and which is respectively orthogonal, magic, cubic, bistochastic. 

\medskip

We have a surjective morphism $A_{N,r}\to C(G_N^-/G_k^-)$, and we must prove that this is an isomorphism.

\medskip

\underline{$S^+$ case}. Since $\#(S_N/S_k)=N!/k!$, it is enough to prove that we have: 
$$\dim(A_{N,r})=\frac{N!}{k!}$$

Let $p_{ij}$ be the standard generators of $A_{N,r}$. By using the Gelfand theorem, we can write $p_{ij}=\chi(X_{ij})$, where $X_{ij}\subset X$ are certain subets of a given set $X$. Now at the level of sets the magic isometry condition on $(p_{ij})$ tells us that the matrix of sets $(X_{ij})$ has the property that its entries are disjoint on columns, and form partitions of $X$ on rows.

\medskip

So, let us try to understand this property for $N$ fixed, and $r=1,2,3,\ldots$

\medskip

-- At $r=1$ we simply have a partition $X=X_1\sqcup\ldots\sqcup X_N$. So, the universal model can be any such partition, with $X_i\neq 0$ for any $i$.

\medskip

-- At $r=2$ the universal model is best described as follows: $X$ is the $N\times N$ square in $\mathbb R^2$, regarded as a union of $N^2$ unit tiles, minus the diagonal, the sets $X_{1i}$ are the disjoint unions on rows, and the sets $X_{2i}$ are the disjoint unions on columns.

\medskip

-- At $r\geq 3$, the universal solution is similar: we can take $X$ to be the $N^r$ cube in $\mathbb R^r$, with all tiles having pairs of equal coordinates removed, and say that the sets $X_{si}$ for $s$ fixed are the various ``slices'' of $X$ in the direction of the $s$-th coordinate of $\mathbb R^r$.

\medskip

Summarizing, the above discussion tells us that $\dim(A_{N,r})$ equals the number of tiles in the above set $X\subset\mathbb R^r$. But these tiles correspond by definition to the various $r$-tuples $(i_1,\ldots,i_r)\in\{1,\ldots,N\}^r$ with all $i_k$ different, and since there are exactly $N!/k!$ such $r$-tuples, we obtain, as desired:
$$\dim(A_{N,r})=\frac{N!}{k!}$$

\underline{$H^+$ case}. We can use here the same method as for $S_N^+$. This time the functions $p_{ij}$ take values in $\{-1,0,1\}$, and the algebra generated by their squares $p_{ij}^2$ coincides with the one computed above for $S_N^+$, having dimension $N!/k!$. Now by taking into account the $N-k$ possible signs we obtain the following estimate, which gives the result:
$$\dim(A_{N,r})\leq\frac{2^{N-k}N!}{k!}=\#(H_N/H_k)$$

\underline{$O^+$ case}. We can use the same method, namely a straightforward application of the Gelfand theorem. However, instead of performing a dimension count, which is no longer possible, we have to complete here any transposed $N\times r$ isometry whose entries commute to a $N\times N$ orthogonal matrix. But this is the same as completing a system of $r$ orthogonal norm 1 vectors in $\mathbb R^N$ into an orthonormal basis of $\mathbb R^N$, which is of course possible.

\medskip

\underline{$B^+$ case}. Since we have a surjective map $C(O_N^+)\to C(B_N^+)$, we obtain a surjective map $C_+(O_N^+/O_k^+)\to A_{N,r}$, and hence surjective maps as follows:
$$C(O_N/O_k)\to A_{N,r}\to C(B_N/B_k)$$

The point now is that this composition is the following canonical map:
$$C(O_N/O_k)\to C(B_N/B_k)$$

Now by looking at the column vector $\xi=(1,\ldots,1)^t$, which is fixed by the stochastic matrices, we conclude that the map on the right is an isomorphism, and we are done.
\end{proof}

We refer to \cite{bss} and related papers for more on the above.

\section*{8e. Exercises} 

Things got fairy complicated in this chapter, basically leading to hot research questions, and as a unique exercise on all this, in the same spirit, we have:

\begin{exercise}
Axiomatize the free manifolds, as a continuation of the above.
\end{exercise}

There is of course no need of completely solving this exercise, and some preliminary study, for some very simple classes of manifolds, of your choice, will more than do.

\part{Easy geometries}

\ \vskip50mm

\begin{center}
{\em Give my love to Rose please, won't you mister

Take her all my money, tell her to buy some pretty clothes

Tell my boy that daddy's so proud of him

And don't forget to give my love to Rose}
\end{center}

\chapter{Half-liberation}

\section*{9a. Spheres and tori}

We have seen in chapter 4 that the quadruplets of type $(S,T,U,K)$ can be axiomatized, and that at the level of basic examples we have 4 such quadruplets, corresponding to the usual real and complex geometries $\mathbb R^N,\mathbb C^N$, and to the free versions of these:
$$\xymatrix@R=50pt@C=50pt{
\mathbb R^N_+\ar[r]&\mathbb C^N_+\\
\mathbb R^N\ar[u]\ar[r]&\mathbb C^N\ar[u]
}$$

Our purpose in what follows will be that of extending the above diagram, with the construction of some supplementary examples. There are two methods here:

\medskip

(1) Look for intermediate geometries $\mathbb R^N\subset\mathcal X\subset\mathbb R^N_+$, and their complex analogues.

\medskip

(2) Look for intermediate geometries $\mathbb R^N\subset\mathcal X\subset \mathbb C^N$, and their free analogues.

\medskip

We will see that, in each case, there is a ``standard'' solution, and that these solutions can be combined. Thus, we will end up with a total of $3\times3=9$ solutions, as follows:
$$\xymatrix@R=40pt@C=40pt{
\mathbb R^N_+\ar[r]&\mathbb T\mathbb R^N_+\ar[r]&\mathbb C^N_+\\
\mathbb R^N_*\ar[u]\ar[r]&\mathbb T\mathbb R^N_*\ar[u]\ar[r]&\mathbb C^N_*\ar[u]\\
\mathbb R^N\ar[u]\ar[r]&\mathbb T\mathbb R^N\ar[u]\ar[r]&\mathbb C^N\ar[u]
}$$

There is quite some work to be done here, and the construction of these 9 geometries will take us the whole present chapter, and most of the next chapter as well. We will see also, at the end of the next chapter, that under certain strong axioms, of combinatorial type, these 9 geometries are conjecturally the only ones. Finally, in chapters 11-12 below we will discuss a number of related topics, such as twisting, and matrix models.

\bigskip

Observe that all this is in direct continuation of what we did in Part I, with no obvious relation with Part II. However, and here comes our point, once these intermediate geometries constructed, we will also have to ``develop'' them, meaning looking at various homogeneous spaces $X=G/H$, and other manifolds $X$, and here the theory developed in Part II, while mainly designed for being of help with free geometry, will be of great use. By the way, let us mention too that the intermediate geometries to be developed here, in Part III, will be quite close to the classical geometries, of $\mathbb R^N,\mathbb C^N$, and so our manifolds $X$ will start having interesting geometric features, that we will explore as well. Finally, for our outline to be complete, later in Part IV we will go back to the free geometries, of $\mathbb R^N_+,\mathbb C^N_+$, and develop more theory there, based on all this knowledge.

\bigskip

A few words on our motivations, too. There are many of them, as follows:

\bigskip

(1) The real half-classical geometry, of $\mathbb R^N_*$, is something very interesting in quantum group theory, due to the fact that the corresponding orthogonal group, $O_N^*$, is conjecturally the unique intermediate subgroup $O_N\subset G\subset O_N^+$. Thus, regardless of our precise axioms here, the geometry of $\mathbb R^N_*$ can only be, at least conjecturally, the only intermediate geometry $\mathbb R^N\subset\mathcal X\subset\mathbb R^N_+$, so is definitely worth a study, mathematically speaking.

\bigskip

(2) Still talking $\mathbb R^N_*$, the geometry here is not that far from the geometries of $\mathbb R^N,\mathbb C^N$, so the study here can potentially lead into many things not available in the free case, and not discussed so far in this book, such as differential geometry, Lie theory, K-theory, and many more. Thus, in a certain sense, $\mathbb R^N_*$ is the ``bridge'' between our free geometry and more traditional visions of noncommutative geometry, such as Connes' \cite{co1}.

\bigskip

(3) And pretty much the same goes for the other geometries to be investigated in this Part III, and particularly for the complex half-classical geometry, of $\mathbb C^N_*$, for reasons similar to those in the real case, and also for the twisted geometries, of $\bar{\mathbb R}^N,\bar{\mathbb C}^N$, making a link between our free quantum groups and free geometry with the more traditional vision of quantum groups and noncommutative geometry of Drinfeld-Jimbo \cite{dri}, \cite{jim}. 

\bigskip

(4) Thus, plenty of good reasons for looking into such things, be them philosophical, or more concrete. And also, talking now physics, an interesting discovery, due to Bhowmick-D'Andrea-Dabrowski \cite{bdd}, and fine-tuned in their later paper with Das \cite{bd+}, is that the computations for the free gauge group of the Standard Model, in its Chamseddine-Connes formulation \cite{cc1}, \cite{cc2}, crucially involve the quantum group $U_N^*$.

\bigskip

And that is all, for the moment, more on this later. Getting to work now, our starting point will be the general axioms found in chapter 4, which are as follows: 

\begin{definition}
An abstract noncommutative geometry is described by a quadruplet $(S,T,U,K)$, formed of intermediate objects as follows,
$$S^{N-1}_\mathbb R\subset S\subset S^{N-1}_{\mathbb C,+}$$
$$T_N\subset T\subset\mathbb T_N^+$$
$$O_N\subset U\subset U_N^+$$
$$H_N\subset K\subset K_N^+$$
subject to a set of connecting formulae between them, as follows, 
$$\begin{matrix}
S&=&S_U\\
S\cap\mathbb T_N^+&=&T&=&K\cap\mathbb T_N^+\\
G^+(S)&=&<O_N,T>&=&U\\
K^+(T)&=&U\cap K_N^+&=&K
\end{matrix}$$
with the usual convention that all this is up to the equivalence relation.
\end{definition}

All this is of course quite tricky, and a bit simplified too, in the above form, and for full details on this, and comments, we refer to chapter 4. Now with this in hand, let us get into our first question, namely finding intermediate geometries as follows:
$$\mathbb R^N\subset\mathcal X\subset\mathbb R^N_+$$

Since such a geometry is given by a quadruplet $(S,T,U,K)$, as above, forgetting about correspondences, we are led to 4 different intermediate object questions, as follows:
$$S^{N-1}_\mathbb R\subset S\subset S^{N-1}_{\mathbb R,+}$$
$$T_N\subset T\subset T_N^+$$
$$O_N\subset U\subset O_N^+$$
$$H_N\subset K\subset H_N^+$$

At the sphere and torus level, there are obviously uncountably many solutions, without supplementary assumptions, and it is hard to get beyond this, with bare hands. Thus, our hopes will basically come from the unitary and reflection quantum groups, where things are more rigid than for spheres and tori. Let us record, however, the following fact regarding the spheres, from \cite{bme}, which will appear to be relevant, later on:

\index{half-classical sphere}

\begin{theorem}
The algebraic manifold $S^{(k)}\subset S^{N-1}_{\mathbb R,+}$ obtained by imposing the relations $a_1\ldots a_k=a_k\ldots a_1$ to the standard coordinates of $S^{N-1}_{\mathbb R,+}$ is as follows:
\begin{enumerate}
\item At $k=1$ we have $S^{(k)}=S^{N-1}_{\mathbb R,+}$.

\item At $k=2,4,6,\ldots$ we have $S^{(k)}=S^{N-1}_\mathbb R$.

\item At $k=3,5,7,\ldots$ we have $S^{(k)}=S^{(3)}$.
\end{enumerate}
\end{theorem}

\begin{proof}
As a first observation, the commutation relations $ab=ba$ imply the following relations, for any $k\geq2$:
$$a_1\ldots a_k=a_k\ldots a_1$$

Thus, for any $k\geq2$, we have an inclusion $S^{(2)}\subset S^{(k)}$. It is also elementary to check that the relations $abc=cba$ imply the following relations, for any $k\geq3$ odd:
$$a_1\ldots a_k=a_k\ldots a_1$$

Thus, for any $k\geq3$ odd, we have an inclusion $S^{(3)}\subset S^{(k)}$. Our claim now is that we have an inclusion as follows, for any $k\geq2$:
$$S^{(k+2)}\subset S^{(k)}$$

In order to prove this, we must show that the relations $a_1\ldots a_{k+2}=a_{k+2}\ldots a_1$ between the coordinates $x_1,\ldots,x_N$ imply the relations $a_1\ldots a_k=a_k\ldots a_1$ between these coordinates $x_1,\ldots,x_N$. But this holds indeed, because of the following implications:
\begin{eqnarray*}
x_{i_1}\ldots x_{i_{k+2}}=x_{i_{k+2}}\ldots x_{i_1}
&\implies&x_{i_1}\ldots x_{i_k}x_j^2=x_j^2x_{i_k}\ldots x_{i_1}\\
&\implies&\sum_jx_{i_1}\ldots x_{i_k}x_j^2=\sum_jx_j^2x_{i_k}\ldots x_{i_1}\\
&\implies&x_{i_1}\ldots x_{i_k}=x_{i_k}\ldots x_{i_1}
\end{eqnarray*}

Summing up, we have proved that we have inclusions as follows:
$$S^{(2)}\subset\ldots\ldots\subset S^{(6)}\subset S^{(4)}\subset S^{(2)}$$
$$S^{(3)}\subset\ldots\ldots\subset S^{(7)}\subset S^{(5)}\subset S^{(3)}$$

Thus, we are led to the conclusions in the statement.
\end{proof}

As a conclusion, the ``privileged'' intermediate sphere $S^{N-1}_\mathbb R\subset S\subset S^{N-1}_{\mathbb R,+}$ that we are looking for can only be the sphere $S^{(3)}$, obtained via the following relations:
$$abc=cba$$

We should mention that, following \cite{bme}, it is possible to go further in this direction, with a study of the spheres given by relations of the following type, with $\sigma\in S_k$: 
$$a_1\ldots a_k=a_{\sigma(1)}\ldots a_{\sigma(k)}$$ 

But this leads to a similar conclusion, namely that the sphere $S^{(3)}$ constructed above is the only new solution. We will discuss this, which is a bit technical, later, in chapter 13 below. All this remains, however, quite ad-hoc. In short, we have constructed so far a new real sphere, $S^{(3)}$, and we some evidence for the fact that this sphere might be the only new one, under some extra combinatorial axioms, which are quite technical.

\section*{9b. Quantum groups}

Let us focus now on the quantum groups. We will see that there is a lot more rigidity here, with regards to what happens for the spheres and tori, which makes things simpler. Our goal will be that of finding the intermediate quantum groups as follows:
$$O_N\subset U\subset O_N^+$$
$$H_N\subset K\subset H_N^+$$

Quite surprisingly, these two questions are of quite different nature. Indeed, regarding $O_N\subset U\subset O_N^+$, there is a solution here, denoted $O_N^*$, coming via the relations $abc=cba$, and conjecturally nothing more. Regarding however $H_N\subset K\subset H_N^+$, here it is possible to use for instance crossed products, for constructing uncountably many solutions.

\medskip

In short, in connection with our intermediate geometry question, we do have in principle our solution, coming via the relations $abc=cba$, and this is compatible with our above $S^{(3)}$ guess for the spheres. In order to get started, let us recall that we have:

\begin{theorem}
The basic quantum unitary and reflection groups, namely
$$
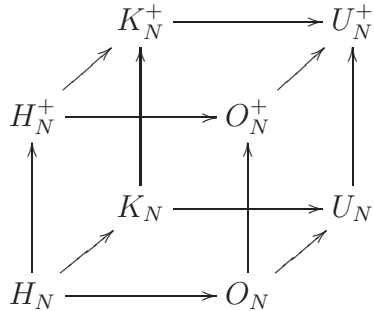
\xymatrix@R=16pt@C=16pt{
&K_N^+\ar[rr]&&U_N^+\\
H_N^+\ar[rr]\ar[ur]&&O_N^+\ar[ur]\\
&K_N\ar[rr]\ar[uu]&&U_N\ar[uu]\\
H_N\ar[uu]\ar[ur]\ar[rr]&&O_N\ar[uu]\ar[ur]
}$$
are all easy, coming from certain categories of partitions.
\end{theorem}

\begin{proof}
This is something that we already discussed, in chapter 2, the corresponding categories of partitions being as follows: 
$$
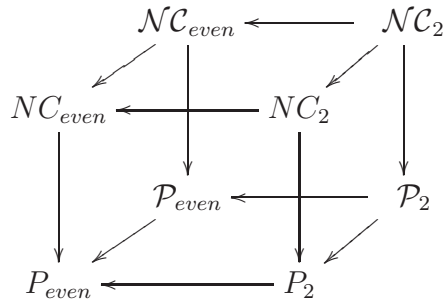
\xymatrix@R=17pt@C=0pt{
&\mathcal{NC}_{even}\ar[dl]\ar[dd]&&\ \ \mathcal {NC}_2\ar[dl]\ar[ll]\ar[dd]\\
NC_{even}\ar[dd]&&NC_2\ar[dd]\ar[ll]\\
&\mathcal P_{even}\ar[dl]&&\ \ \mathcal P_2\ar[dl]\ar[ll]\\
P_{even}&&P_2\ar[ll]
}$$

Thus, we are led to the conclusion in the statement.
\end{proof}

Getting back now to the half-liberation question, let us start by constructing the solutions. The result here, which is well-known as well, is as follows:

\index{half-classical group}
\index{half-classical orthogonal group}
\index{half-classical unitary group}
\index{half-classical hyperoctahedral group}
\index{half-classical reflection group}

\begin{theorem}
We have quantum groups as follows, obtained via the half-commutation relations $abc=cba$, which fit into the diagram of basic quantum groups:
$$\xymatrix@R=55pt@C=55pt{
K_N^*\ar[r]&U_N^*\\
H_N^*\ar[u]\ar[r]&O_N^*\ar[u]
}$$
These quantum groups are all easy, and the corresponding categories of partitions fit into the diagram of categories of partitions for the basic quantum groups.
\end{theorem}

\begin{proof}
This is standard, from \cite{bsp}, the idea being that the half-commutation relations $abc=cba$ come from the map $T_{\slash\hskip-1.6mm\backslash\hskip-1.1mm|\hskip0.5mm}$ associated to the half-classical crossing:
$$\slash\hskip-1.9mm\backslash\hskip-1.7mm|\hskip0.5mm\in P(3,3)$$

Thus, the quantum groups in the statement are indeed easy, obtained by adding the half-classical crossing $\slash\hskip-1.9mm\backslash\hskip-1.7mm|\hskip0.5mm$ to the corresponding categories of noncrossing partitions. We obtain the following categories, with $*$ standing for the fact that, when relabelling clockwise the legs $\circ\bullet\circ\bullet\ldots$, the formula $\#\circ=\#\bullet$ must hold in each block:
$$\xymatrix@R=50pt@C=50pt{
\mathcal P_{even}^*\ar[d]&\mathcal P_2^*\ar[l]\ar[d]\\
P_{even}^*&P_2^*\ar[l]
}$$

Finally, the fact that our new quantum groups and categories fit well into the previous diagrams of quantum groups and categories is clear from this. See \cite{bsp}.
\end{proof}

The point now is that we have the following result:

\index{category of pairings}

\begin{theorem}
There is only one proper intermediate easy quantum group
$$O_N\subset G\subset O_N^+$$
namely the half-classical orthogonal group $O_N^*$.
\end{theorem}

\begin{proof}
According to our definition for the easy quantum groups, we must compute here the intermediate categories of pairings, as follows:
$$NC_2\subset D\subset P_2$$

But this can be done via some standard combinatorics, in three steps, as follows:

\medskip

(1) Let $\pi\in P_2-NC_2$, having $s\geq 4$ strings. Our claim is that:

\medskip

-- If $\pi\in P_2-P_2^*$, there exists a semicircle capping $\pi'\in P_2-P_2^*$.

\medskip

-- If $\pi\in P_2^*-NC_2$, there exists a semicircle capping $\pi'\in P_2^*-NC_2$.

\medskip

Indeed, both these assertions can be easily proved, by drawing pictures.

\medskip

(2) Consider now a partition $\pi\in P_2(k,l)-NC_2(k,l)$. Our claim is that:

\medskip

-- If $\pi\in P_2(k, l)-P_2^*(k,l)$ then $<\pi>=P_2$.

\medskip

-- If $\pi\in P_2^*(k,l)-NC_2(k,l)$ then $<\pi>=P_2^*$.

\medskip

This can be indeed proved by recurrence on the number of strings, $s=(k+l)/2$, by using (1), which provides us with a descent procedure $s\to s-1$, at any $s\geq4$.

\medskip

(3) Finally, assume that we are given an easy quantum group $O_N\subset G\subset O_N^+$, coming from certain sets of pairings $D(k,l)\subset P_2(k,l)$. We have three cases:

\medskip

-- If $D\not\subset P_2^*$, we obtain $G=O_N$.

\medskip

-- If $D\subset P_2,D\not\subset NC_2$, we obtain $G=O_N^*$.

\medskip

-- If $D\subset NC_2$, we obtain $G=O_N^+$.

\medskip

Thus, we are led to the conclusion in the statement.
\end{proof}

It is in fact conjectured that the above result holds without the easiness assumption, and we refer here to \cite{bc+}. Thus, we have now an answer to our questions, with the half-classical real geometry being most likely unique, between classical and free real.

\bigskip

In practice now, what we have to do is to construct this geometry, and its complex analogue as well, and check the axioms from chapter 4. Let us begin by constructing the corresponding quadruplets. We have here the following result:

\index{half-classical quadruplets}
\index{half-classical geometry}

\begin{theorem}
We half-classical real and complex quadruplets, as follows,
$$
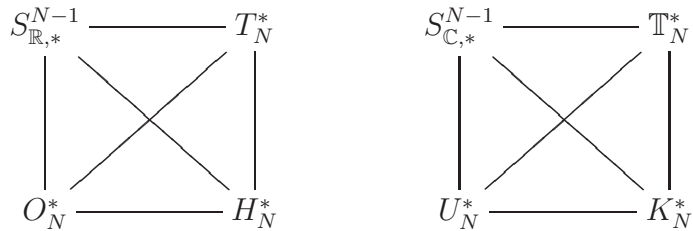
\xymatrix@R=50pt@C=50pt{
S^{N-1}_{\mathbb R,*}\ar@{-}[r]\ar@{-}[d]\ar@{-}[dr]&T_N^*\ar@{-}[l]\ar@{-}[d]\ar@{-}[dl]\\
O_N^*\ar@{-}[u]\ar@{-}[ur]\ar@{-}[r]&H_N^*\ar@{-}[l]\ar@{-}[ul]\ar@{-}[u]
}
\qquad\qquad 
\xymatrix@R=50pt@C=50pt{
S^{N-1}_{\mathbb C,*}\ar@{-}[r]\ar@{-}[d]\ar@{-}[dr]&\mathbb T_N^*\ar@{-}[l]\ar@{-}[d]\ar@{-}[dl]\\
U_N^*\ar@{-}[u]\ar@{-}[ur]\ar@{-}[r]&K_N^*\ar@{-}[l]\ar@{-}[ul]\ar@{-}[u]
}$$
obtained via $abc=cba$, imposed to the standard coordinates and their adjoints.
\end{theorem}

\begin{proof}
This is more of an empty statement, with the real quantum groups being those above, and with the other objects, namely complex quantum groups, and then spheres and tori, being constructed in a similar way, by starting with the free objects, and imposing the relations $abc=cba$ to the standard coordinates, and their adjoints.
\end{proof}

We should mention here that, while the above constructions look trivial, the story with them was not trivial at all. Indeed, while things are certainly clear in the real case, in the complex case there are several possible ways of imposing the half-commutation relations $abc=cba$ to the standard coordinates and their adjoints, as follows:

\bigskip

(1) The above way, imposing $abc=cba$ to everything, both the standard coordinates and their adjoints, is the strongest such way, producing the smallest half-liberations, and in particular the smallest half-classical unitary quantum group, denoted $U_N^*$. 

\bigskip

(2) In an opposite direction, imposing only the relations $ab^*c=cb^*a$ to the standard coordinates is something reasonable too, and this produces the biggest unitary quantum group which can be reasonably called ``half-classical'', denoted $U_N^\times$.

\bigskip

(3) And then, there are all sorts of intermediate objects in between, $U_N^*\subset U_N^\circ\subset U_N^\times$, and notably the quantum group $U_N^{**}$ obtained by stating that the variables $\{ab^*,a^*b\}$ with $a,b$ standard coordinates should all commute, which is something natural too.

\bigskip

All this is quite technical, related to all sorts of advanced quantum group considerations, and there has been fierce debate all over the 10s, often between certain authors and their inner selves, on which relations to use, and more specifically, on which of the quantum groups $U_N^*\subset U_N^{**}\subset U_N^\times$ is the ``correct'' one. And with the literature on the subject, consisting notably of \cite{ba1}, \cite{ba2}, \cite{ba3}, \cite{ba4}, \cite{bb1}, \cite{bb2}, \cite{bdd}, \cite{bd+}, \cite{bic}, \cite{bdu}, \cite{mwe}, \cite{twe} being often confusing, with $U_N^*$ usually denoting the ``correct'' quantum group at the time of the paper, from the viewpoint of the paper, in a somewhat reckless way.

\bigskip

The solution to these questions came quite recently, first from the paper of Mang-Weber \cite{mwe}, who classified all the easy quantum groups $U_N^*\subset U_N^\circ\subset U_N^\times$, which allows one to have a more relaxed, complete perspective on all this, and then with the present noncommutative geometry considerations, coming as a continuation of \cite{bb2}, the point being that by Mang-Weber \cite{mwe} the only ``good'' quantum group among $U_N^*\subset U_N^{**}\subset U_N^\times$, which produces a noncommutative geometry in our sense, is $U_N^*$. We will back to this later, when discussing \cite{mwe}, and classification for our noncommutative geometries.

\section*{9c. Matrix models}

In order to check now our noncommutative geometry axioms, we are in need of a better understanding of the half-liberation operation, via some supplementary results. Let us start with the following simple observation, regarding the real spheres:

\index{matrix model}

\begin{proposition}
We have a morphism of $C^*$-algebras as follows,
$$C(S^{N-1}_{\mathbb R,*})\to M_2(C(S^{N-1}_\mathbb C))\quad,\quad 
x_i\to\begin{pmatrix}0&z_i\\ \bar{z}_i&0\end{pmatrix}$$ 
where $z_i$ are the standard coordinates of $S^{N-1}_\mathbb C$.
\end{proposition}

\begin{proof}
We have to prove that the matrices $X_i$ on the right satisfy the defining relations for $S^{N-1}_{\mathbb R,*}$. But these matrices are self-adjoint, and we have:
\begin{eqnarray*}
\sum_iX_i^2
&=&\sum_i\begin{pmatrix}0&z_i\\ \bar{z}_i&0\end{pmatrix}^2\\
&=&\sum_i\begin{pmatrix}|z_i|^2&0\\0&|z_i|^2\end{pmatrix}\\
&=&\begin{pmatrix}1&0\\0&1\end{pmatrix}
\end{eqnarray*}

As for the half-commutation relations, these follow from the following formula:
\begin{eqnarray*}
X_iX_jX_k
&=&\begin{pmatrix}0&z_i\\ \bar{z}_i&0\end{pmatrix}\begin{pmatrix}0&z_j\\ \bar{z}_j&0\end{pmatrix}\begin{pmatrix}0&z_k\\ \bar{z}_k&0\end{pmatrix}\\
&=&\begin{pmatrix}0&z_i\bar{z}_jz_k\\ \bar{z}_iz_j\bar{z}_k&0\end{pmatrix}
\end{eqnarray*}

Indeed, the quantities on the right being symmetric in $i,k$, this gives the result.
\end{proof}

Regarding the complex spheres, the result here is similar, as follows:

\begin{proposition}
We have a morphism of $C^*$-algebras as follows,
$$C(S^{N-1}_{\mathbb C,*})\to M_2(C(S^{N-1}_\mathbb C\times S^{N-1}_\mathbb C))\quad,\quad 
x_i\to\begin{pmatrix}0&z_i\\ y_i&0\end{pmatrix}$$ 
where $y_i,z_i$ are the standard coordinates of $S^{N-1}_\mathbb C\times S^{N-1}_\mathbb C$.
\end{proposition}

\begin{proof}
We have to prove that the matrices $X_i$ on the right satisfy the defining relations for $S^{N-1}_{\mathbb C,*}$. We have the following computation:
\begin{eqnarray*}
\sum_iX_iX_i^*
&=&\sum_i\begin{pmatrix}0&z_i\\ y_i&0\end{pmatrix}\begin{pmatrix}0&\bar{y}_i\\ \bar{z}_i&0\end{pmatrix}\\
&=&\sum_i\begin{pmatrix}|z_i|^2&0\\0&|y_i|^2\end{pmatrix}\\
&=&\begin{pmatrix}1&0\\0&1\end{pmatrix}
\end{eqnarray*}

We have as well the following computation:
\begin{eqnarray*}
\sum_iX_i^*X_i
&=&\sum_i\begin{pmatrix}0&\bar{y}_i\\ \bar{z}_i&0\end{pmatrix}\begin{pmatrix}0&z_i\\ y_i&0\end{pmatrix}\\
&=&\sum_i\begin{pmatrix}|y_i|^2&0\\0&|z_i|^2\end{pmatrix}\\
&=&\begin{pmatrix}1&0\\0&1\end{pmatrix}
\end{eqnarray*}

As for the half-commutation relations, these follow from the following formula:
\begin{eqnarray*}
X_iX_jX_k
&=&\begin{pmatrix}0&z_i\\ y_i&0\end{pmatrix}\begin{pmatrix}0&z_j\\ y_j&0\end{pmatrix}\begin{pmatrix}0&z_k\\ y_k&0\end{pmatrix}\\
&=&\begin{pmatrix}0&z_iy_jz_k\\ y_iz_jy_k&0\end{pmatrix}
\end{eqnarray*}

Indeed, the quantities on the right being symmetric in $i,k$, this gives the result.
\end{proof}

Our goal now will be that of proving that the morphisms constructed above are faithful, up to the usual equivalence relation for the quantum algebraic manifolds. For this purpose, we will use some projective geometry arguments, the idea being that of proving that the above morphisms are indeed isomorphisms, at the projective version level, and then lifting these isomorphism results, to the affine setting. Let us recall that:

\bigskip

(1) The real projective space $P^{N-1}_\mathbb R$ is the space of lines in $\mathbb R^N$ passing through the origin. We have a quotient map $S^{N-1}_\mathbb R\to P^{N-1}_\mathbb R$, producing an embedding $C(P^{N-1}_\mathbb R)\subset C(S^{N-1}_\mathbb R)$, whose image is the algebra generated by the variables $p_{ij}=x_ix_j$. 

\bigskip

(2) The complex projective space $P^{N-1}_\mathbb C$ has a similar description, namely is the space of complex lines in $\mathbb C^N$ passing through the origin, and we have an embedding $C(P^{N-1}_\mathbb C)\subset C(S^{N-1}_\mathbb C)$, whose image is generated by the variables $p_{ij}=x_i\bar{x}_j$.

\bigskip

The spaces $P^{N-1}_\mathbb R,P^{N-1}_\mathbb C$ have the following functional analytic description:

\index{projective space}
\index{real projective space}
\index{complex projective space}

\begin{theorem}
We have presentation results as follows,
\begin{eqnarray*}
C(P^{N-1}_\mathbb C)&=&C^*_{comm}\left((p_{ij})_{i,j=1,\ldots,N}\Big|p=p^*=p^2,Tr(p)=1\right)\\
C(P^{N-1}_\mathbb R)&=&C^*_{comm}\left((p_{ij})_{i,j=1,\ldots,N}\Big|p=\bar{p}=p^*=p^2,Tr(p)=1\right)
\end{eqnarray*}
where by $C^*_{comm}$ we mean as usual universal commutative $C^*$-algebra.
\end{theorem}

\begin{proof}
We use the elementary fact that the spaces $P^{N-1}_\mathbb C,P^{N-1}_\mathbb R$, as defined above, are respectively the spaces of rank one projections in $M_N(\mathbb C),M_N(\mathbb R)$. With this picture in mind, it is clear that we have arrows $\leftarrow$. In order to construct now arrows $\to$, consider the universal algebras on the right, $A_C,A_R$. These algebras being both commutative, by the Gelfand theorem we can write, with $X_C,X_R$ being certain compact spaces:
$$A_C=C(X_C)\quad,\quad 
A_R=C(X_R)$$

Now by using the coordinate functions $p_{ij}$, we conclude that $X_C,X_R$ are certain spaces of rank one projections in $M_N(\mathbb C),M_N(\mathbb R)$. In other words, we have embeddings:
$$X_C\subset P^{N-1}_\mathbb C\quad,\quad 
X_R\subset P^{N-1}_\mathbb R$$

Bsy transposing we obtain arrows $\to$, as desired.
\end{proof}

The above result suggests constructing free projective spaces $P^{N-1}_{\mathbb R,+},P^{N-1}_{\mathbb C,+}$, simply by lifting the commutativity conditions between the variables $p_{ij}$. However, there is something wrong with this, and more specifically with $P^{N-1}_{\mathbb R,+}$, coming from the fact that if certain noncommutative coordinates $x_1,\ldots,x_N$ are self-adjoint, then the corresponding projective coordinates $p_{ij}=x_ix_j$ are not necessarily self-adjoint:
$$x_i=x_i^*\ \,\not\!\!\!\!\implies x_ix_j=(x_ix_j)^*$$

In short, our attempt to construct free projective spaces $P^{N-1}_{\mathbb R,+},P^{N-1}_{\mathbb C,+}$ as above is not exactly correct, with the space $P^{N-1}_{\mathbb R,+}$ being rather ``irrelevant'', and with the space $P^{N-1}_{\mathbb C,+}$ being probably the good one, but being at the same time ``real and complex''. Observe that there is some similarity here with the following key result, from chapter 4:
$$PO_N^+=PU_N^+$$

To be more precise, we have good evidence here for the fact that, in the free setting, the projective geometry is at the same time real and complex. We will be back to this later, but in the meantime, in view of this, let us formulate the following definition:

\index{free projective space}

\begin{definition}
Associated to any $N\in\mathbb N$ is the following universal algebra,
$$C(P^{N-1}_+)=C^*\left((p_{ij})_{i,j=1,\ldots,N}\Big|p=p^*=p^2,Tr(p)=1\right)$$
whose abstract spectrum is called ``free projective space''.
\end{definition}

Observe that we have embeddings of noncommutative spaces, as follows:
$$P^{N-1}_\mathbb R\subset P^{N-1}_\mathbb C\subset P^{N-1}_+$$

Let us compute now the projective versions of the noncommutative spheres that we have, including the half-classical ones. We use the following formalism here:

\index{projective version}

\begin{definition}
The projective version of a closed subspace $S\subset S^{N-1}_{\mathbb C,+}$ is the quotient space $S\to PS$ determined by the fact that 
$$C(PS)\subset C(S)$$
is the subalgebra generated by $p_{ij}=x_ix_j^*$, called projective coordinates.
\end{definition}

In the classical case, this fits with the usual definition. We will be back with more details in chapter 15 below, which is dedicated to the study of projective geometry. We have the following result, coming from \cite{ba8}, \cite{bgo}, \cite{bme}:

\begin{theorem}
The projective versions of the basic spheres are as follows,
$$\xymatrix@R=10mm@C=20mm{
S^{N-1}_{\mathbb R,+}\ar[r]&S^{N-1}_{\mathbb C,+}\\
S^{N-1}_{\mathbb R,*}\ar[r]\ar[u]&S^{N-1}_{\mathbb C,*}\ar[u]\\
S^{N-1}_\mathbb R\ar[r]\ar[u]&S^{N-1}_\mathbb C\ar[u]}
\qquad
\xymatrix@R=14mm@C=20mm{\\ \to}
\qquad
\xymatrix@R=10mm@C=20mm{
P^{N-1}_+\ar[r]&P^{N-1}_+\\
P^{N-1}_\mathbb C\ar[r]\ar[u]&P^{N-1}_\mathbb C\ar[u]\\
P^{N-1}_\mathbb R\ar[r]\ar[u]&P^{N-1}_\mathbb C\ar[u]}$$
modulo, in the free case, a GNS construction with respect to the uniform integration.
\end{theorem}

\begin{proof}
The formulae on the bottom are true by definition. For the formulae on top, we have to prove first that the variables $p_{ij}=x_ix_j^*$ over the free sphere $S^{N-1}_{\mathbb C,+}$ satisfy the defining relations for $C(P^{N-1}_+)$. In order to check this, we first have:
$$(p^*)_{ij}
=p_{ji}^*
=(x_jx_i^*)^*
=x_ix_j^*
=p_{ij}$$

We have as well the following computation:
$$(p^2)_{ij}
=\sum_kp_{ik}p_{kj}
=\sum_kx_ix_k^*x_kx_j^*
=x_ix_j^*
=p_{ij}$$

Finally, we have as well the following computation:
$$Tr(p)
=\sum_kp_{kk}
=\sum_kx_kx_k^*
=1$$

Thus, we have embeddings of algebraic manifolds, as follows:
$$PS^{N-1}_{\mathbb R,+}\subset PS^{N-1}_{\mathbb C,+}\subset P^{N-1}_+$$

Regarding now the GNS construction assertion, this follows by reasoning as in the case of the free spheres, the idea being that the uniform integration on these projective spaces comes from the uniform integration over the following quantum group:
$$PO_N^+=PU_N^+$$

All this is quite technical, and we will not need this result, in what follows. We refer here to \cite{bme}, and we will back to this in chapter 15 below. Finally, regarding the middle assertions, concerning the projective versions of the half-classical spheres, it is enough to prove here that we have inclusions as follows:
$$P^{N-1}_\mathbb C\subset PS^{N-1}_{\mathbb R,*}\subset PS^{N-1}_{\mathbb C,*}\subset P^{N-1}_\mathbb C$$

But this can be done in 3 steps, as follows:

\medskip

(1) $P^{N-1}_\mathbb C\subset PS^{N-1}_{\mathbb R,*}$. In order to prove this, we recall from Proposition 9.7 that we have a morphism as follows, where $z_i$ are the standard coordinates of $S^{N-1}_\mathbb C$:
$$C(S^{N-1}_{\mathbb R,*})\to M_2(C(S^{N-1}_\mathbb C))\quad,\quad 
x_i\to\begin{pmatrix}0&z_i\\ \bar{z}_i&0\end{pmatrix}$$ 

Now observe that this model maps the projective coordinates as follows:
$$p_{ij}\to P_{ij}=\begin{pmatrix}z_i\bar{z}_j&0\\0&\bar{z}_iz_j\end{pmatrix}$$

Thus, at the level of generated algebras, our model maps:
$$<p_{ij}>\to <P_{ij}>=C(P^{N-1}_\mathbb C)$$

We conclude from this that we have a quotient map as follows:
$$C(PS^{N-1}_{\mathbb R,*})\to C(P^{N-1}_\mathbb C)$$

Thus at the level of corresponding spaces, we have, as desired, an inclusion:
$$P^{N-1}_\mathbb C\subset PS^{N-1}_{\mathbb R,*}$$

(2) $PS^{N-1}_{\mathbb R,*}\subset PS^{N-1}_{\mathbb C,*}$. This is something trivial, coming by functoriality of the operation $S\to PS$, from the inclusion of spheres:
$$S^{N-1}_{\mathbb R,*}\subset S^{N-1}_{\mathbb C,*}$$

(3) $PS^{N-1}_{\mathbb C,*}\subset P^{N-1}_\mathbb C$. This follows from the half-commutation relations, which imply:
$$ab^*cd^*
=cb^*ad^*
=cd^*ab^*$$

Indeed, this shows that the projective version $PS^{N-1}_{\mathbb C,*}$ is classical, and so:
$$PS^{N-1}_{\mathbb C,*}
\subset(P^{N-1}_+)_{class}
=P^{N-1}_\mathbb C$$

Thus, we are led to the conclusion in the statement.
\end{proof}

We can go back now to the spheres, and we have the following result:

\index{matrix model}
\index{antidiagonal matrix}
\index{antidiagonal model}

\begin{theorem}
We have a morphism of $C^*$-algebras as follows,
$$C(S^{N-1}_{\mathbb R,*})\subset M_2(C(S^{N-1}_\mathbb C))\quad,\quad 
x_i\to\begin{pmatrix}0&z_i\\ \bar{z}_i&0\end{pmatrix}$$ 
where $z_i$ are the standard coordinates of $S^{N-1}_\mathbb C$.
\end{theorem}

\begin{proof}
We know from Proposition 9.7 that we have a morphism as above, and the injectivity follows from Theorem 9.12, by using a standard grading trick. See \cite{bic}.
\end{proof}

In the case of the complex spheres we have a similar result, as follows:

\begin{theorem}
We have a morphism of $C^*$-algebras as follows,
$$C(S^{N-1}_{\mathbb C,*})\to M_2(C(S^{N-1}_\mathbb C\times S^{N-1}_\mathbb C))\quad,\quad 
x_i\to\begin{pmatrix}0&z_i\\ y_i&0\end{pmatrix}$$ 
where $y_i,z_i$ are the standard coordinates of $S^{N-1}_\mathbb C\times S^{N-1}_\mathbb C$.
\end{theorem}

\begin{proof}
Again, we know from Proposition 9.8 that we have a morphism as above, and the injectivity follows from Theorem 9.12, via a grading trick, as explained in \cite{bic}.
\end{proof}

Summarizing, we have some interesting affine and projective geometry results regarding the half-classical case. The point now is that the same arguments apply to the tori, and to the quantum groups. We first have the following result:

\begin{proposition}
The real half-classical quadruplet, namely
$$
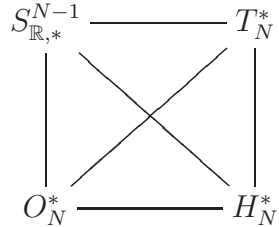
\xymatrix@R=50pt@C=50pt{
S^{N-1}_{\mathbb R,*}\ar@{-}[r]\ar@{-}[d]\ar@{-}[dr]&T_N^*\ar@{-}[l]\ar@{-}[d]\ar@{-}[dl]\\
O_N^*\ar@{-}[u]\ar@{-}[ur]\ar@{-}[r]&H_N^*\ar@{-}[l]\ar@{-}[ul]\ar@{-}[u]
}$$
and the complex real half-classical quadruplet, namely
$$
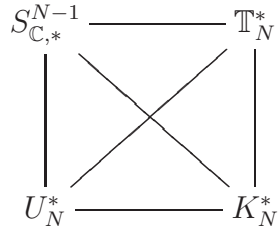
\xymatrix@R=50pt@C=50pt{
S^{N-1}_{\mathbb C,*}\ar@{-}[r]\ar@{-}[d]\ar@{-}[dr]&\mathbb T_N^*\ar@{-}[l]\ar@{-}[d]\ar@{-}[dl]\\
U_N^*\ar@{-}[u]\ar@{-}[ur]\ar@{-}[r]&K_N^*\ar@{-}[l]\ar@{-}[ul]\ar@{-}[u]
}$$
have $2\times2$ matrix models, constructed by using antidiagonal matrices, as for the spheres.
\end{proposition}

\begin{proof}
This is something that we already know from the spheres, from the various results established above. For the other objects which form the quadruplets, this follows by suitably adapting the proof of Proposition 9.7 and Proposition 9.8.
\end{proof}

Next, once again in analogy with the sphere theory, we have the following result:

\begin{theorem}
The real and complex half-classical quadruplets have the same projective version, which is as follows:
$$\xymatrix@R=50pt@C=50pt{
P^{N-1}_\mathbb C\ar@{-}[r]\ar@{-}[d]\ar@{-}[dr]&P\mathbb T_N\ar@{-}[l]\ar@{-}[d]\ar@{-}[dl]\\
PU_N\ar@{-}[u]\ar@{-}[ur]\ar@{-}[r]&PK_N\ar@{-}[l]\ar@{-}[ul]\ar@{-}[u]
}$$
\end{theorem}

\begin{proof}
As before, this is something that we already know from the spheres, from the various results established above. For the other objects which form the quadruplets, this follows from Proposition 9.15, by suitably adapting the proof of Theorem 9.12.
\end{proof}

Finally, completing our study, we have the following result:

\begin{theorem}
The $2\times2$ antidiagonal matrix models for the real and complex half-classical quadruplets, constructed above, are faithful.
\end{theorem}

\begin{proof}
This is something that we already know from the spheres. For the other objects, this follows by suitably adapting the proof of Theorem 9.13 and Theorem 9.14.
\end{proof}

Let us mention that the above results are part of a series of more general results, regarding matrix models for half-liberations. We will be back to this later.

\section*{9d. Axiom check}

Let us check now the axioms, for our real and complex half-classical quadruplets. We first need some quantum isometry group results, which are available from \cite{ba1}, \cite{ba2}, \cite{ba3}, \cite{ba4}, for which we refer for the full details. First, we have the following result:

\begin{theorem}
The quantum isometry groups of the basic spheres are
$$\xymatrix@R=10mm@C=20mm{
S^{N-1}_{\mathbb R,+}\ar[r]&S^{N-1}_{\mathbb C,+}\\
S^{N-1}_{\mathbb R,*}\ar[r]\ar[u]&S^{N-1}_{\mathbb C,*}\ar[u]\\
S^{N-1}_\mathbb R\ar[r]\ar[u]&S^{N-1}_\mathbb C\ar[u]}
\qquad
\xymatrix@R=14.5mm@C=20mm{\\ \to}
\qquad
\xymatrix@R=11.5mm@C=23mm{
O_N^+\ar[r]&U_N^+\\
O_N^*\ar[r]\ar[u]&U_N^*\ar[u]\\
O_N\ar[r]\ar[u]&U_N\ar[u]}$$
modulo identifying, as usual, the various $C^*$-algebraic completions.
\end{theorem}

\begin{proof}
We already know this, from chapter 3, for the spheres on top and bottom, so we just have to prove the results in the middle. So, assume as in chapter 3 that we have an action $G\curvearrowright S^{N-1}_{\mathbb C,*}$. From $\Phi(x_a)=\sum_ix_i\otimes u_{ia}$ we obtain, with $p_{ab}=z_a\bar{z}_b$:
$$\Phi(p_{ab})=\sum_{ij}p_{ij}\otimes u_{ia}u_{jb}^*$$

By multiplying two such arbitrary formulae, we obtain:
\begin{eqnarray*}
\Phi(p_{ab}p_{cd})&=&\sum_{ijkl}p_{ij}p_{kl}\otimes u_{ia}u_{jb}^*u_{kc}u_{ld}^*\\
\Phi(p_{ad}p_{cb})&=&\sum_{ijkl}p_{il}p_{kj}\otimes u_{ia}u_{ld}^*u_{kc}u_{jb}^*
\end{eqnarray*}

The left terms being equal, and the first terms on the right being equal too, we deduce that, with $[a,b,c]=abc-cba$, we must have the following equality:
$$\sum_{ijkl}p_{ij}p_{kl}\otimes u_{ia}[u_{jb}^*,u_{kc},u_{ld}^*]=0$$

Now observe that the products of projective variables $p_{ij}p_{kl}=z_i\bar{z}_jz_k\bar{z}_l$ depend only on the following two cardinalities:
$$|\{i,k\}|,|\{j,l\}|\in\{1,2\}$$

The point now is that this dependence produces the only relations between our variables, we are led in this way to $4$ equations, as follows:

\medskip

(1) $u_{ia}[u_{jb}^*,u_{ka},u_{lb}^*]=0$, $\forall a,b$.

\medskip

(2) $u_{ia}[u_{jb}^*,u_{ka},u_{ld}^*]+u_{ia}[u_{jd}^*,u_{ka},u_{lb}^*]=0$, $\forall a$, $\forall b\neq d$.

\medskip

(3) $u_{ia}[u_{jb}^*,u_{kc},u_{lb}^*]+u_{ic}[u_{jb}^*,u_{ka},u_{lb}^*]=0$, $\forall a\neq c$, $\forall b$.

\medskip

(4) $u_{ia}([u_{jb}^*,u_{kc},u_{ld}^*]+[u_{jd}^*,u_{kc},u_{lb}^*])+u_{ic}([u_{jb}^*,u_{ka},u_{ld}^*]+[u_{jd}^*,u_{ka},u_{lb}^*])=0,\forall a\neq c,b\neq d$.

\bigskip

From (1,2) we conclude that (2) holds with no restriction on the indices. By multiplying now this formula to the left by $u_{ia}^*$, and then summing over $i$, we obtain:
$$[u_{jb}^*,u_{ka},u_{ld}^*]+[u_{jd}^*,u_{ka},u_{lb}^*]=0$$

By applying now the antipode, then the involution, and finally by suitably relabelling all the indices, we successively obtain from this formula:
\begin{eqnarray*}
[u_{dl},u_{ak}^*,u_{bj}]+[u_{bl},u_{ak}^*,u_{dj}]=0
&\implies&[u_{dl}^*,u_{ak},u_{bj}^*]+[u_{bl}^*,u_{ak},u_{dj}^*]=0\\
&\implies&[u_{ld}^*,u_{ka},u_{jb}^*]+[u_{jd}^*,u_{ka},u_{lb}^*]=0
\end{eqnarray*}

Now by comparing with the original relation, above, we conclude that we have:
$$[u_{jb}^*,u_{ka},u_{ld}^*]=[u_{jd}^*,u_{ka},u_{lb}^*]=0$$

Thus we have reached to the formulae defining the quantum group $U_N^*$, and we are done. Finally, in what regards the universality of the action $O_N^*\curvearrowright S^{N-1}_{\mathbb R,*}$, this follows from the universality of the following two actions:
$$U_N^*\curvearrowright S^{N-1}_{\mathbb C,*}\quad,\quad O_N^+\curvearrowright S^{N-1}_{\mathbb R,+}$$

Indeed, we obtain from this that we have $U_N^*\cap O_N^+=O_N^*$, as desired.
\end{proof}

Regarding now the quantum isometry groups of the tori, the computation here, again form \cite{ba1}, \cite{ba2}, \cite{ba3}, \cite{ba4}, and that we partly know from chapter 3, is as follows:

\begin{theorem}
The quantum isometry groups of the basic tori are
$$\xymatrix@R=10mm@C=20mm{
T_N^+\ar[r]&\mathbb T_N^+\\
T_N^*\ar[r]\ar[u]&\mathbb T_N^*\ar[u]\\
T_N\ar[r]\ar[u]&\mathbb T_N\ar[u]}
\qquad
\xymatrix@R=13.5mm@C=20mm{\\ \to}
\qquad
\xymatrix@R=10mm@C=20mm{
H_N^+\ar[r]&K_N^+\\
H_N^*\ar[r]\ar[u]&K_N^*\ar[u]\\
\bar{O}_N\ar[r]\ar@{.}[u]&\bar{U}_N\ar@{.}[u]}$$
with all arrows being inclusions, and with no vertical maps at bottom right.
\end{theorem}

\begin{proof}
As before, we just have to prove the results in the middle. In the real case, we must find the conditions on $G\subset O_N^+$ such that $g_a\to\sum_ig_a\otimes u_{ia}$ defines a coaction. In order for this map to be a coaction, the variables $G_a=\sum_ig_a\otimes u_{ia}$ must satisfy the following relations, which define the groups dual to the tori in the statement:
$$G_a^2=1\quad,\quad 
G_aG_bG_c=G_cG_bG_a$$

In what regards the squares, we have the following formula:
\begin{eqnarray*}
G_a^2
&=&\sum_{ij}g_ig_j\otimes u_{ia}u_{ja}\\
&=&1+\sum_{i\neq j}g_ig_j\otimes u_{ia}u_{ja}
\end{eqnarray*}

As for the products, with the notation $[x,y,z]=xyz-zyx$, we have:
$$\left[G_a,G_b,G_c\right]=\sum_{ijk}g_ig_jg_k\otimes [u_{ia},u_{jb},u_{kc}]$$

From the first relations, $G_a^2=1$, we obtain $G\subset H_N^+$. In order to process now the second relations, $G_aG_bG_c=G_cG_bG_a$, we can split the sum over $i,j,k$, as follows:
\begin{eqnarray*}
\left[G_a,G_b,G_c\right]
&=&\sum_{i,j,k\ distinct}g_ig_jg_k\otimes[u_{ia},u_{jb},u_{kc}]\\
&+&\sum_{i\neq j}g_ig_jg_i\otimes[u_{ia},u_{jb},u_{ic}]\\
&+&\sum_{i\neq j}g_i\otimes[u_{ia},u_{jb},u_{jc}]\\
&+&\sum_{i\neq k}g_k\otimes[u_{ia},u_{ib},u_{kc}]\\
&+&\sum_ig_i\otimes[u_{ia},u_{ib},u_{ic}]
\end{eqnarray*}

Our claim is that the last three sums vanish. Indeed, observe that we have:
$$[u_{ia},u_{ib},u_{ic}]=\delta_{abc}u_{ia}-\delta_{abc}u_{ia}=0$$

Thus the last sum vanishes. Regarding now the fourth sum, we have:
\begin{eqnarray*}
\sum_{i\neq k}[u_{ia},u_{ib},u_{kc}]
&=&\sum_{i\neq k}u_{ia}u_{ib}u_{kc}-u_{kc}u_{ib}u_{ia}\\
&=&\sum_{i\neq k}\delta_{ab}u_{ia}^2u_{kc}-\delta_{ab}u_{kc}u_{ia}^2\\
&=&\delta_{ab}\sum_{i\neq k}[u_{ia}^2,u_{kc}]\\
&=&\delta_{ab}\left[\sum_{i\neq k}u_{ia}^2,u_{kc}\right]\\
&=&\delta_{ab}[1-u_{ka}^2,u_{kc}]\\
&=&0
\end{eqnarray*}

The proof for the third sum is similar. Thus, we are left with the first two sums. By using $g_ig_jg_k=g_kg_jg_i$ for the first sum, the formula becomes:
\begin{eqnarray*}
\left[G_a,G_b,G_c\right]
&=&\sum_{i<k,j\neq i,k}g_ig_jg_k\otimes([u_{ia},u_{jb},u_{kc}]+[u_{ka},u_{jb},u_{ic}])\\
&+&\sum_{i\neq j}g_ig_jg_i\otimes[u_{ia},u_{jb},u_{ic}]
\end{eqnarray*}

In order to have a coaction, the above coefficients must vanish. Now observe that, when setting $i=k$ in the coefficients of the first sum, we obtain twice the coefficients of the second sum. Thus, our vanishing conditions can be formulated as follows:
$$[u_{ia},u_{jb},u_{kc}]+[u_{ka},u_{jb},u_{ic}]=0,\forall j\neq i,k$$

Now observe that at $a=b$ or $b=c$ this condition reads $0+0=0$. Thus, we can formulate our vanishing conditions in a more symmetric way, as follows:
$$[u_{ia},u_{jb},u_{kc}]+[u_{ka},u_{jb},u_{ic}]=0,\forall j\neq i,k,\forall b\neq a,c$$

We use now the trick from \cite{bg1}. We apply the antipode to this formula, and then we relabel the indices $i\leftrightarrow c,j\leftrightarrow b,k\leftrightarrow a$. We succesively obtain in this way:
$$[u_{ck},u_{bj},u_{ai}]+[u_{ci},u_{bj},u_{ak}]=0,\forall j\neq i,k,\forall b\neq a,c$$
$$[u_{ia},u_{jb},u_{kc}]+[u_{ic},u_{jb},u_{ka}]=0,\forall b\neq a,c,\forall j\neq i,k$$

Since we have $[x,y,z]=-[z,y,x]$, by comparing the last formula with the original one, we conclude that our vanishing relations reduce to a single formula, as follows:
$$[u_{ia},u_{jb},u_{kc}]=0,\forall j\neq i,k,\forall b\neq a,c$$

Our first claim is that this formula implies $G\subset H_N^{[\infty]}$, where $H_N^{[\infty]}\subset O_N^+$ is defined via the relations $xyz=0$, for any $x\neq z$ on the same row or column of $u$. In order to prove this, we will just need the $c=a$ particular case of this formula, which reads:
$$u_{ia}u_{jb}u_{ka}=u_{ka}u_{jb}u_{ia},\forall j\neq i,k,\forall a\neq b$$

It is enough to check that the assumptions $j\neq i,k$ and $a\neq b$ can be dropped. But this is what happens indeed, because at $j=i$ we have:
\begin{eqnarray*}
\left[u_{ia},u_{ib},u_{ka}\right]
&=&u_{ia}u_{ib}u_{ka}-u_{ka}u_{ib}u_{ia}\\
&=&\delta_{ab}(u_{ia}^2u_{ka}-u_{ka}u_{ia}^2)\\
&=&0
\end{eqnarray*}

Also, at $j=k$ we have:
\begin{eqnarray*}
\left[u_{ia},u_{kb},u_{ka}\right]
&=&u_{ia}u_{kb}u_{ka}-u_{ka}u_{kb}u_{ia}\\
&=&\delta_{ab}(u_{ia}u_{ka}^2-u_{ka}^2u_{ia})\\
&=&0
\end{eqnarray*}

Finally, at $a=b$ we have:
\begin{eqnarray*}
\left[u_{ia},u_{ja},u_{ka}\right]
&=&u_{ia}u_{ja}u_{ka}-u_{ka}u_{ja}u_{ia}\\
&=&\delta_{ijk}(u_{ia}^3-u_{ia}^3)\\
&=&0
\end{eqnarray*}

Our second claim now is that, due to $G\subset H_N^{[\infty]}$, we can drop the assumptions $j\neq i,k$ and $b\neq a,c$ in the original relations $[u_{ia},u_{jb},u_{kc}]=0$. Indeed, at $j=i$ we have:
\begin{eqnarray*}
[u_{ia},u_{ib},u_{kc}]
&=&u_{ia}u_{ib}u_{kc}-u_{kc}u_{ib}u_{ia}\\
&=&\delta_{ab}(u_{ia}^2u_{kc}-u_{kc}u_{ia}^2)\\
&=&0
\end{eqnarray*}

The proof at $j=k$ and at $b=a$, $b=c$ being similar, this finishes the proof of our claim. We conclude that the half-commutation relations $[u_{ia},u_{jb},u_{kc}]=0$ hold without any assumption on the indices, and so we obtain $G\subset H_N^*$, as claimed. As for the proof in the complex case, this is similar, and we refer here to \cite{ba3} and related papers.
\end{proof}

By intersecting now with $K_N^+$, as required by our $(S,T,U,K)$ axioms, we obtain:

\begin{theorem}
The quantum reflection groups of the basic tori are
$$\xymatrix@R=10mm@C=20mm{
T_N^+\ar[r]&\mathbb T_N^+\\
T_N^*\ar[r]\ar[u]&\mathbb T_N^*\ar[u]\\
T_N\ar[r]\ar[u]&\mathbb T_N\ar[u]}
\qquad
\xymatrix@R=13.5mm@C=20mm{\\ \to}
\qquad
\xymatrix@R=10mm@C=20mm{
H_N^+\ar[r]&K_N^+\\
H_N^*\ar[r]\ar[u]&K_N^*\ar[u]\\
H_N\ar[r]\ar[u]&K_N\ar[u]}$$
with all the arrows being inclusions.
\end{theorem}

\begin{proof}
We already know that the results on the left and on the right hold indeed. As for the results in the middle, these follow from Theorem 9.19.
\end{proof}

We can now formulate our extension result, as follows:

\index{half-classical geometry}

\begin{theorem}
We have basic noncommutative geometries, as follows,
$$\xymatrix@R=30pt@C=80pt{
\mathbb R^N_+\ar[r]&\mathbb C^N_+\\
\mathbb R^N_*\ar[u]\ar[r]&\mathbb C^N_*\ar[u]\\
\mathbb R^N\ar[u]\ar[r]&\mathbb C^N\ar[u]
}$$
with each $\mathbb K^N_\times$ symbol standing for the corresponding $(S,T,U,K)$ quadruplet.
\end{theorem}

\begin{proof}
We have to check the axioms from chapter 4, for the half-classical geometries. The algebraic axioms are all clear, and the quantum isometry axioms follow from the above computations. Next in line, we have to prove the following formulae: 
$$O_N^*=<O_N,T_N^*>$$
$$U_N^*=<U_N,\mathbb T_N^*>$$

By using standard generation results, it is enough to prove the first formula. Moreover, once again by standard generation results, it is enough to check that:
$$H_N^*=<H_N,T_N^*>$$

The inclusion $\supset$ being clear, we are left with proving the inclusion $\subset$. But this follows from the formula $H_N^*=T_N^*\rtimes S_N$, established by Raum-Weber in \cite{rwe}, as follows:
\begin{eqnarray*}
H_N^*
&=&T_N^*\rtimes S_N\\
&=&<S_N,T_N^*>\\
&\subset&<H_N,T_N^*>
\end{eqnarray*}

Alternatively, these formulae can be established by using the technology in \cite{bdu}, or by using categories and easiness. Finally, the axiom $S=S_U$ can be proved as in the classical and free cases, by using the Weingarten formula, and the following ergodicity property:
$$\left(id\otimes\int_U\right)\Phi(x)=\int_Sx$$

Our claim, which will finish the proof, is that this holds as well in the half-classical case. Indeed, in the real case, where $x_i=x_i^*$, it is enough to check the above equality on an arbitrary product of coordinates, $x_{i_1}\ldots x_{i_k}$. The left term is as follows:
\begin{eqnarray*}
\left(id\otimes\int_{O_N^*}\right)\Phi(x_{i_1}\ldots x_{i_k})
&=&\sum_{j_1\ldots j_k}x_{j_1}\ldots x_{j_k}\int_{O_N^*}u_{j_1i_1}\ldots u_{j_ki_k}\\
&=&\sum_{j_1\ldots j_k}\ \sum_{\pi,\sigma\in P_2^*(k)}\delta_\pi(j)\delta_\sigma(i)W_{kN}(\pi,\sigma)x_{j_1}\ldots x_{j_k}\\
&=&\sum_{\pi,\sigma\in P_2^*(k)}\delta_\sigma(i)W_{kN}(\pi,\sigma)\sum_{j_1\ldots j_k}\delta_\pi(j)x_{j_1}\ldots x_{j_k}
\end{eqnarray*}

Let us look now at the last sum on the right. We have to sum there quantities of type $x_{j_1}\ldots x_{j_k}$, over all choices of multi-indices $j=(j_1,\ldots,j_k)$ which fit into our given pairing $\pi\in P_2^*(k)$. But by using the relations $x_ix_jx_k=x_kx_jx_i$,  and then $\sum_ix_i^2=1$ in order to simplify, we conclude that the sum of these quantities is 1. Thus, we obtain:
$$\left(id\otimes\int_{O_N^*}\right)\Phi(x_{i_1}\ldots x_{i_k})
=\sum_{\pi,\sigma\in P_2^*(k)}\delta_\sigma(i)W_{kN}(\pi,\sigma)$$

On the other hand, another application of the Weingarten formula gives:
\begin{eqnarray*}
\int_{S^{N-1}_{\mathbb R,*}}x_{i_1}\ldots x_{i_k}
&=&\int_{O_N^*}u_{1i_1}\ldots u_{1i_k}\\
&=&\sum_{\pi,\sigma\in P_2^*(k)}\delta_\pi(1)\delta_\sigma(i)W_{kN}(\pi,\sigma)\\
&=&\sum_{\pi,\sigma\in P_2^*(k)}\delta_\sigma(i)W_{kN}(\pi,\sigma)
\end{eqnarray*}

Thus, we are done. In the complex case the proof is similar, by adding exponents. For further details, we refer to \cite{bgo} for the real case, and to \cite{ba1} for the complex case.
\end{proof} 

Summarizing, in relation with the plan made in the beginning of this chapter, we have done so far half of our extension program, for the noncommutative geometries that we have. The second half, along with some classification work, is for the next chapter.

\section*{9e. Exercises}

There are many interesting questions regarding the half-classical geometry. First, in relation with what was discussed in the above, we have:

\begin{exercise}
Develop a full theory of the main half-classical groups, 
$$\xymatrix@R=55pt@C=55pt{
K_N^*\ar[r]&U_N^*\\
H_N^*\ar[u]\ar[r]&O_N^*\ar[u]
}$$
in particular by working out in detail their easiness property.
\end{exercise}

This is something which is quite standard, that we already discussed in the above, at least partly. The problem is now that of developing the full theory.

\begin{exercise}
Explain, both at the algebraic and the probabilistic level, how the general theory developed in chapters $5$-$8$ can be applied to the half-classical situation, in order to talk about more general classes of half-classical manifolds, with algebraic and probabilistic results about them, generalizing what we already have.
\end{exercise}

There is quite some work to be done here, with all this being very instructive. In fact, if you love this exercise, write a book on half-classical geometry, afterwards.

\chapter{Hybrid geometries}

\section*{10a. Spheres and tori}

We finish here the extension program outlined in the previous chapter. To be more precise, we have seen so far that have basic noncommutative geometries as follows, with each $F^N_\times$ symbol standing for the corresponding $(S,T,U,K)$ quadruplet:
$$
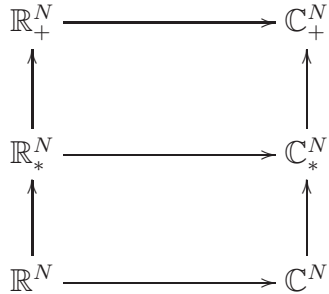
\xymatrix@R=30pt@C=80pt{
\mathbb R^N_+\ar[r]&\mathbb C^N_+\\
\mathbb R^N_*\ar[u]\ar[r]&\mathbb C^N_*\ar[u]\\
\mathbb R^N\ar[u]\ar[r]&\mathbb C^N\ar[u]
}$$

We will see in this chapter that there are some privileged intermediate geometries between the real and the complex ones, completing our diagram as follows:
$$
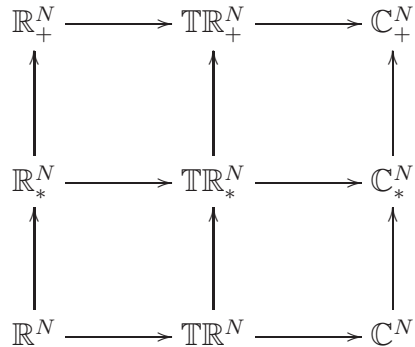
\xymatrix@R=40pt@C=40pt{
\mathbb R^N_+\ar[r]&\mathbb T\mathbb R^N_+\ar[r]&\mathbb C^N_+\\
\mathbb R^N_*\ar[u]\ar[r]&\mathbb T\mathbb R^N_*\ar[u]\ar[r]&\mathbb C^N_*\ar[u]\\
\mathbb R^N\ar[u]\ar[r]&\mathbb T\mathbb R^N\ar[u]\ar[r]&\mathbb C^N\ar[u]
}$$

We will see as well that, that under strong combinatorial axioms, of easiness and uniformity type, these 9 geometries are the only ones. With this being actually the inteersting part of the present chapter, because the general topic of complexification is something quite technical, and not very beautiful, and the new geometries that we will construct in this way have no obvious application. But hey, mathematician is our job, so if you consider that the 9-diagram looks better than the 6-one, for aesthetic reasons, which is something that I do, let's just do the work, without thinking much.

\bigskip

In order to get started, we will solve the classical problem first. An intermediate geometry $\mathbb R^N\subset\mathcal X\subset\mathbb C^N$ is by definition given by a quadruplet $(S,T,U,K)$, whose components are subject to the following conditions, along with a number of axioms:
$$S^{N-1}_\mathbb R\subset S\subset S^{N-1}_\mathbb C$$
$$T_N\subset T\subset\mathbb T_N$$
$$O_N\subset U\subset U_N$$
$$H_N\subset K\subset K_N$$

Our plan will be that of investigating first these intermediate object questions. Then, we will discuss the verification of the geometric axioms, for the solutions that we found. And then, afterwards, we will discuss the half-classical and the free cases as well.

\bigskip

In what regards the intermediate sphere problem, $S^{N-1}_\mathbb R\subset S\subset S^{N-1}_\mathbb C$, there are obviously infinitely many solutions, because there are so many real algebraic manifolds in between. However, we have a ``privileged'' solution, constructed as follows:

\index{intermediate sphere}
\index{hybrid sphere}

\begin{theorem}
We have an intermediate sphere as follows, which consists of the multiples, by scalars in $\mathbb T$, of the points of the real sphere $S^{N-1}_\mathbb R$:
$$S^{N-1}_\mathbb R\subset\mathbb TS^{N-1}_\mathbb R\subset S^{N-1}_\mathbb C$$
Moreover, this sphere appears as the affine lift of $P^{N-1}_\mathbb R$, inside $S^{N-1}_\mathbb C$.
\end{theorem}

\begin{proof}
The first assertion is clear. Regarding now the second assertion, which justified the term ``privileged'' used above, observe first that we have:
$$P\mathbb TS^{N-1}_\mathbb R
=PS^{N-1}_\mathbb R
=P^{N-1}_\mathbb R$$

Conversely, assume that a closed subset $S\subset S^{N-1}_\mathbb C$ satisfies:
$$PS\subset P^{N-1}_\mathbb R$$

For $x\in S$ the projective coordinates $p_{ij}=x_i\bar{x}_j$ must then be real:
$$x_i\bar{x}_j=\bar{x}_ix_j$$

Thus, we must have the following equalities:
$$\frac{x_1}{\bar{x}_1}=\frac{x_2}{\bar{x}_2}=\ldots=\frac{x_N}{\bar{x}_N}$$

Now if we denote by $\lambda\in\mathbb T$ this common number, we succesively have:
\begin{eqnarray*}
\frac{x_i}{\bar{x}_i}=\lambda
&\iff&x_i=\lambda\bar{x}_i\\
&\iff&x_i^2=\lambda |x_i|^2\\
&\iff&x_i=\pm\sqrt{\lambda}|x_i|
\end{eqnarray*}

Thus we obtain $x\in\sqrt{\lambda}S^{N-1}_\mathbb R$, and this gives the result.
\end{proof}

In the case of the tori, we have a similar result, with some new objects added, which are quite natural in the torus setting, as follows:

\index{intermediate torus}
\index{hybrid torus}
\index{Clifford torus}

\begin{theorem}
We have an intermediate torus as follows, which appears as the affine lift of the Clifford torus $PT_N=T_{N-1}$, inside the complex torus $\mathbb T_N$: 
$$T_N\subset\mathbb TT_N\subset\mathbb T_N$$
More generally, we have intermediate tori as follows, with $r\in\mathbb N\cup\{\infty\}$, 
$$T_N\subset\mathbb Z_rT_N\subset\mathbb T_N$$
all whose projective versions equal the Clifford torus $PT_N=T_{N-1}$.
\end{theorem}

\begin{proof}
The first assertion, regarding $\mathbb TT_N$, follows exactly as for the spheres, as in proof of Theorem 10.1. The second assertion is clear as well, because we have:
$$P\mathbb Z_rT_N
=PT_N
=T_{N-1}$$

Thus, we are led to the conclusion in the statement.
\end{proof}

In connection with the above statement, an interesting question is that of classifying the intermediate tori, which in our case are usual compact groups, as follows:
$$T_N\subset T\subset\mathbb T_N$$

At the group dual level, we must classify the following intermediate quotients:
$$\mathbb Z^N\to\Gamma\to\mathbb Z_2^N$$

There are many examples of such groups, and this even when imposing strong supplementary conditions, such as having an action of the symmetric group $S_N$ on the generators. We will not go further in this direction, our main idea being anyway that of basing our study mostly on quantum group theory, and on the related notion of easiness.

\section*{10b. Quantum groups}

At the group level now, the situation is much more rigid, and becomes quite interesting. We have the following result from \cite{bc+}, to start with:

\index{maximality}

\begin{theorem}
The following inclusion of compact groups is maximal,
$$\mathbb TO_N\subset U_N$$
in the sense that there is no intermediate compact group in between.\end{theorem}

\begin{proof}
In order to prove this result, consider as well the following group:
$$\mathbb TSO_N=\left\{wU\Big| w\in\mathbb T,U\in SO_N\right\}$$ 

Observe that we have $\mathbb TSO_N=\mathbb TO_N$ if $N$ is odd. If $N$ is even the group $\mathbb TO_N$ has two connected components, with $\mathbb TSO_N$ being the component containing the identity.

Let us denote by $\mathfrak{so}_N,\mathfrak u_N$ the Lie algebras of $SO_N,U_N$.  It is well-known that $\mathfrak u_N$ consists of the matrices $M\in M_N(\mathbb C)$ satisfying $M^*=-M$, and that:
$$\mathfrak{so}_N=\mathfrak u_N\cap M_N(\mathbb R)$$

Also, it is easy to see that the Lie algebra of $\mathbb TSO_N$ is $\mathfrak{so}_N\oplus i\mathbb R$.

\medskip

\underline{Step 1}. Our first claim is that if $N\geq 2$, the adjoint representation of $SO_N$ on the space of real symmetric matrices of trace zero is irreducible.

\medskip

Let indeed $X \in M_N(\mathbb R)$ be symmetric with trace zero. We must prove that the following space consists of all the real symmetric matrices of trace zero:
$$V=span\left\{UXU^t\Big|U \in SO_N\right\}$$

We first prove that $V$ contains all the diagonal matrices of trace zero. Since we may diagonalize $X$ by conjugating with an element of $SO_N$, our space $V$ contains a nonzero diagonal matrix of trace zero. Consider such a matrix:
$$D=\begin{pmatrix}
d_1\\
&\ddots\\
&&d_N
\end{pmatrix}$$

We can conjugate this matrix by the following matrix:
$$\begin{pmatrix}
0&-1&0\\
1&0&0\\
0&0&I_{N-2}
\end{pmatrix}\in SO_N$$

We conclude that our space $V$ contains as well the following matrix: 
$$D'=\begin{pmatrix}
d_2\\
&d_1\\
&&d_3\\
&&&\ddots\\
&&&&d_N
\end{pmatrix}$$

More generally, we see that for any $1\leq i,j\leq N$ the diagonal matrix obtained from $D$ by interchanging $d_i$ and $d_j$ lies in $V$. Now since $S_N$ is generated by transpositions, it follows that $V$ contains any diagonal matrix obtained by permuting the entries of $D$. 

\medskip

But it is well-known that this representation of $S_N$ on the diagonal matrices of trace zero is irreducible, and hence $V$ contains all such diagonal matrices, as claimed.

\medskip

In order to conclude now, assume that $Y$ is an arbitrary real symmetric matrix of trace zero. We can find then an element $U\in SO_N$ such that $UYU^t$ is a diagonal matrix of trace zero.  But we then have $UYU^t \in V$, and hence also $Y\in V$, as desired.

\medskip

\underline{Step 2}. Our claim is that the inclusion $\mathbb TSO_N\subset U_N$ is maximal in the category of connected compact groups. 

\medskip

Let indeed $G$ be a connected compact group satisfying:
$$\mathbb TSO_N\subset G\subset U_N$$

Then $G$ is a Lie group. Let $\mathfrak g$ denote its Lie algebra, which satisfies:
$$\mathfrak{so}_N\oplus i\mathbb R\subset\mathfrak g\subset\mathfrak u_N$$

Let $ad_{G}$ be the action of $G$ on $\mathfrak g$ obtained by differentiating the adjoint action of $G$ on itself. This action turns $\mathfrak g$ into a $G$-module. Since $SO_N \subset G$, $\mathfrak g$ is also a $SO_N$-module. Now if $G\neq\mathbb TSO_N$, then since $G$ is connected we must have: $$\mathfrak{so}_N\oplus i\mathbb{R}\neq\mathfrak g$$

It follows from the real vector space structure of the Lie algebras $\mathfrak u_N$ and $\mathfrak{so}_N$ that there exists a nonzero symmetric real matrix of trace zero $X$ such that:
$$iX\in\mathfrak g$$

We know that the space of symmetric real matrices of trace zero is an irreducible representation of $SO_N$ under the adjoint action. Thus $\mathfrak g$ must contain all such $X$, and hence $\mathfrak g=\mathfrak u_N$.  But since $U_N$ is connected, it follows that $G=U_N$.  

\medskip

\underline{Step 3}. Let us compute now the commutant of  $SO_N$ in $ M_N(\mathbb C)$. Our first claim is that at $N=2$, this commutant is as follows: 
$$SO_2'
=\left\{\begin{pmatrix}
\alpha&\beta\\
-\beta&\alpha
\end{pmatrix}\Big|\alpha,\beta\in\mathbb C\right\}$$

As for the case $N\geq3$, our claim here is that this commutant is as follows:
$$SO_N'=\left\{\alpha I_N\Big|\alpha\in\mathbb C\right\}$$

Indeed, at $N=2$, the above formula is clear. At $N\geq 3$ now, an element in $X\in SO_N'$ commutes with any diagonal matrix having exactly $N-2$ entries equal to $1$ and two entries equal to $-1$. Hence $X$ is diagonal. Now since $X$ commutes with any even permutation matrix, and we have assumed $N\geq 3$, it commutes in particular with the permutation matrix associated with the cycle $(i,j,k)$ for any $1<i<j<k$, and hence all the entries of $X$ are the same. We conclude that $X$ is a scalar matrix, as claimed.

\medskip 

\underline{Step 4}. Our claim now is that the set of matrices with nonzero trace is dense in $SO_N$.

\medskip  

At $N=2$ this is clear, since the set of elements in $SO_2$ having a given trace is finite.  So assume $N>2$, and consider a matrix as follows:
$$T\in SO_N\simeq SO(\mathbb R^N)\quad,\quad 
Tr(T)=0$$

Let $E\subset\mathbb R^N$ be a 2-dimensional subspace preserved by $T$, such that: 
$$T_{|E} \in SO(E)$$ 
 
Let $\varepsilon>0$ and let $S_\varepsilon \in SO(E)$ satisfying the following condition:
$$||T_{|E}-S_\varepsilon||<\varepsilon$$

Moreover, in the $N=2$ case, we can assume that $T$ satisfies as well:
$$Tr(T_{|E})\neq Tr(S_\varepsilon)$$

Now define $T_\varepsilon\in SO(\mathbb R^N)=SO_N$ by the following formulae:
$$T_{\varepsilon|E}=S_\varepsilon\quad,\quad 
T_{\varepsilon|E^\perp}=T_{|E^\perp}$$

It is clear that we have the following estimate:
$$||T-T_\varepsilon|| \leq ||T_{|E}-S_\varepsilon||<\varepsilon$$

Also, we have the following estimate, which proves our claim:
$$Tr(T_\varepsilon)=Tr(S_\varepsilon)+Tr(T_{|E^\perp})\neq0$$

\underline{Step 5}. Our claim now is that $\mathbb TO_N$ is the normalizer of $\mathbb TSO_N$ in $U_N$, i.e. is the subgroup of $U_N$ consisting of the unitaries $U$ for which, for all $X\in\mathbb TSO_N$:
$$U^{-1}XU \in\mathbb TSO_N$$

Indeed, $\mathbb TO_N$ normalizes $\mathbb TSO_N$, so we must prove that if $U\in U_N$ normalizes $\mathbb TSO_N$ then $U\in\mathbb TO_N$. First note that $U$ normalizes $SO_N$, because if $X \in SO_N$ then:
$$U^{-1}XU \in\mathbb TSO_N$$

Thus we have a formula as follows, for some $\lambda\in\mathbb T$ and $Y\in SO_N$:
$$U^{-1}XU=\lambda Y$$ 

If $Tr(X)\neq0$, we have $\lambda\in\mathbb R$ and hence:
$$\lambda Y=U^{-1}XU \in SO_N$$

The set of matrices having nonzero trace being dense in $SO_N$, we conclude that $U^{-1}XU \in SO_N$ for all $X\in SO_N$. Thus, we have:
\begin{eqnarray*}
X \in SO_N
&\implies&(UXU^{-1})^t(UXU^{-1})=I_N\\
&\implies&X^tU^tUX= U^tU\\
&\implies&U^tU \in SO_N'
\end{eqnarray*}

It follows that at $N\geq 3$ we have $U^tU=\alpha I_N$, with $\alpha \in \mathbb T$, since $U$ is unitary. Hence we have $U=\alpha^{1/2}(\alpha^{-1/2}U)$ with:
$$\alpha^{-1/2}U\in O_N\quad,\quad
U\in\mathbb TO_N$$

If $N=2$, $(U^tU)^t=U^tU$ gives again $U^tU=\alpha I_2$, and we conclude as before.

\medskip

\underline{Step 6}. Our claim is that the inclusion $\mathbb TO_N\subset U_N$ is maximal.

\medskip

Assume indeed that $\mathbb TO_N\subset G\subset U_N$ is a compact group such that $G\neq U_N$. It is a well-known fact that the connected component of the identity in $G$ is a normal subgroup, denoted $G_0$. Since we have $\mathbb TSO_N\subset G_0 \subset U_N$, we must have:
$$G_0=\mathbb TSO_N$$

But since $G_0$ is normal in $G$, the group $G$ normalizes $\mathbb TSO_N$, and hence $G\subset\mathbb TO_N$, which finishes the proof.
\end{proof}

Along the same lines, still following \cite{bc+}, we have as well the following result:

\index{projective unitary group}

\begin{theorem}
The following inclusion of compact groups is maximal,
$$PO_N\subset PU_N$$
in the sense that there is no intermediate compact group in between.
\end{theorem}

\begin{proof}
This follows from Theorem 10.3. Indeed, assuming $PO_N\subset G \subset PU_N$, the preimage of this subgroup under the quotient map $U_N\to PU_N$ would be then a proper intermediate subgroup of $\mathbb TO_N\subset U_N$, which is a contradiction.
\end{proof}

Finally, still following \cite{bc+}, we have as well the following result:

\index{half-classical orthogonal group}

\begin{theorem}
The following inclusion of compact quantum groups is maximal,
$$O_N\subset O_N^*$$
in the sense that there is no intermediate compact quantum group in between.\end{theorem}

\begin{proof}
Consider indeed a sequence of surjective Hopf $*$-algebra maps as follows, whose  composition is the canonical surjection:
$$C(O_N^*)\overset{f}\longrightarrow A\overset{g}\longrightarrow C(O_N)$$

This produces a diagram of Hopf algebra maps with pre-exact rows, as follows:
$$\xymatrix@R=45pt@C=33pt{
\mathbb C\ar[r]&C(PO_N^*)\ar[d]^{f_|}\ar[r]&C(O_N^*)\ar[d]^f\ar[r]&C(\mathbb Z_2)\ar[r]\ar@{=}[d]&\mathbb C\\
\mathbb C\ar[r]&PA\ar[d]^{g_|}\ar[r]&A\ar[d]^g\ar[r]&C(\mathbb Z_2)\ar[r]\ar@{=}[d]&\mathbb C\\
\mathbb C\ar[r]&PC(O_N)\ar[r]&C(O_N)\ar[r]&C(\mathbb Z_2)\ar[r]&\mathbb C}$$
 
Consider now the following composition, with the isomorphism on the left being something well-known, coming from \cite{bdu}, as explained in chapter 9: 
$$C(PU_N)\simeq C(PO_N^*)\overset{f_|}\longrightarrow PA\overset{g_|}\longrightarrow PC(O_N)\simeq C(PO_N)$$

This induces, at the group level, the folowing embedding:
$$PO_N\subset PU_N$$

Thus $f_|$ or $g_|$ is an isomorphism. If $f_|$ is an isomorphism we get a commutative diagram of Hopf algebra morphisms with pre-exact rows, as follows:
$$\xymatrix@R=45pt@C=33pt{
\mathbb C\ar[r]&C(PO_N^*)\ar@{=}[d]\ar[r]&C(O_N^*)\ar[d]^f\ar[r]&C(\mathbb Z_2)\ar[r]\ar@{=}[d]&\mathbb C\\
\mathbb C\ar[r]&C(PO_N^*)\ar[r]&A\ar[r]&C(\mathbb Z_2)\ar[r]&\mathbb C}$$
 
Then $f$ is an isomorphism. Similarly if $g_|$ is an isomorphism, then $g$ is an isomorphism, and this gives the result. See \cite{bc+}.
\end{proof}

In connection now with our question, which is that of classifying the intermediate groups $O_N\subset G\subset U_N$, the above results lead to a dichotomy, coming from:
$$PG\in\{PO_N,PU_N\}$$

In the lack of a classification result here, which is surely well-known, here are some basic examples of such intermediate groups, which are all well-known:

\begin{proposition}
We have compact groups $O_N\subset G\subset U_N$ as follows:
\begin{enumerate}
\item The following groups, depending on a parameter $r\in\mathbb N\cup\{\infty\}$, 
$$\mathbb Z_rO_N=\left\{wU\Big|w\in\mathbb Z_r,U\in O_N\right\}$$ 
whose projective versions equal $PO_N$, and the biggest of which is the group $\mathbb TO_N$, which appears as affine lift of $PO_N$. 

\item The following groups, depending on a parameter $d\in 2\mathbb N\cup\{\infty\}$,
$$U_N^d=\left\{U\in U_N\Big|\det U\in\mathbb Z_d\right\}$$
interpolating between $U_N^2$ and $U_N^\infty=U_N$, whose projective versions equal $PU_N$.
\end{enumerate}
\end{proposition}

\begin{proof}
All the assertions are elementary, the idea being as follows:

\medskip

(1) We have indeed compact groups $\mathbb Z_rO_N$ with $r\in\mathbb N\cup\{\infty\}$ as in the statement, whose projective versions are given by:
$$P\mathbb Z_rO_N=PO_N$$

At $r=\infty$ we obtain the group $\mathbb TO_N$, and the fact that this group appears as the affine lift of $PO_N$ follows exactly as in the sphere case, as in the proof of Theorem 10.1. 

\medskip

(2) As a first observation, the following formula, with $d\in\mathbb N\cup\{\infty\}$, defines indeed a closed subgroup $U_N^d\subset U_N$:
$$U_N^d=\left\{U\in U_N\Big|\det U\in\mathbb Z_d\right\}$$

In the case where $d$ is even, this subgroup contains the orthogonal group $O_N$. As for the last assertion, namely $PU_N^d=PU_N$, this follows either be suitably rescaling the unitary matrices, or by applying the result in Theorem 10.3.
\end{proof}

The above result suggests that the solutions of $O_N\subset G\subset U_N$ should come from $O_N,U_N$, by succesively applying the following constructions:
$$G\to\mathbb Z_rG\quad,\quad 
G\to G\cap U_N^d$$

These operations do not exactly commute, but normally we should be led in this way to a 2-parameter series, unifying the two 1-parameter series from (1,2). However, some other groups like $\mathbb Z_NSO_N$ work too, so all this is probably a bit more complicated. 

\bigskip

In what follows we will be mostly interested in the group $\mathbb TO_N$, which fits with the spheres and tori that we already have, in view of our axiomatization purposes. This particular group $\mathbb TO_N$, and the whole series $\mathbb Z_rO_N$ with $r\in\mathbb N\cup\{\infty\}$ that it is part of, is known to be easy, the precise result, from Tarrago-Weber \cite{twe}, being as follows:

\begin{theorem}
We have the following results:
\begin{enumerate}
\item $\mathbb TO_N$ is easy, the corresponding category $\bar{P}_2\subset P_2$ consisting of the pairings having the property that when flattening, we have $\#\circ=\#\bullet$, globally.

\item $\mathbb Z_rO_N$ is easy, the corresponding category $P_2^r\subset P_2$ consisting of the pairings having the property that when flattening, we have $\#\circ=\#\bullet(r)$, globally.
\end{enumerate}
\end{theorem}

\begin{proof}
These results are standard and well-known, the proof being as follows:

\medskip

(1) If we denote the standard corepresentation by $u=zv$, with $z\in\mathbb T$ and with $v=\bar{v}$, then in order to have $Hom(u^{\otimes k},u^{\otimes l})\neq\emptyset$, the $z$ variabes must cancel, and in the case where they cancel, we obtain the same Hom-space as for $O_N$. 

\medskip

Now since the cancelling property for the $z$ variables corresponds precisely to the fact that $k,l$ must have the same numbers of $\circ$ symbols minus $\bullet$ symbols, the associated Tannakian category must come from the category of pairings $\bar{P}_2\subset P_2$, as claimed.

\medskip

(2) This is something that we already know at $r=1,\infty$, where the group in question is $O_N,\mathbb TO_N$. The proof in general is similar, by writing $u=zv$ as above.
\end{proof}

Quite remarkably, the above result has the following converse, also from \cite{twe}:

\begin{theorem}
The proper intermediate easy compact groups
$$O_N\subset G\subset U_N$$
are precisely the groups $\mathbb Z_rO_N$, with $r\in\{2,3,\ldots,\infty\}$.
\end{theorem}

\begin{proof}
According to our conventions for the easy quantum groups, which apply of course to the classical case, we must compute the following intermediate categories:
$$\mathcal P_2\subset D\subset P_2$$

So, assume that we have such a category, $D\neq\mathcal P_2$, and pick an element $\pi\in D-\mathcal P_2$, assumed to be flat. We can modify $\pi$, by performing the following operations:

\medskip

(1) First, we can compose with the basic crossing, in order to assume that $\pi$ is a partition of type $\cap\ldots\ldots\cap$, consisting of consecutive semicircles. Our assumption $\pi\notin\mathcal P_2$ means that at least one semicircle is colored black, or white.

\medskip

(2) Second, we can use the basic mixed-colored semicircles, and cap with them all the mixed-colored semicircles. Thus, we can assume that $\pi$ is a nonzero partition of type $\cap\ldots\ldots\cap$, consisting of consecutive black or white semicircles.

\medskip

(3) Third, we can rotate, as to assume that $\pi$ is a partition consisting of an upper row of white semicircles, $\cup\ldots\ldots\cup$, and a lower row of white semicircles, $\cap\ldots\ldots\cap$. Our assumption $\pi\notin\mathcal P_2$ means that this latter partition is nonzero.

\medskip

For $a,b\in\mathbb N$ consider the partition consisting of an upper row of $a$ white semicircles, and a lower row of $b$ white semicircles, and set:
$$\mathcal C=\left\{\pi_{ab}\Big|a,b\in\mathbb N\right\}\cap D$$

According to the above, the following happens:
$$\pi\in<\mathcal C>$$ 

The point now is that we have the following results:

\medskip

(1) There exists $r\in\mathbb N\cup\{\infty\}$ such that $\mathcal C$ equals the following set:
$$\mathcal C_r=\left\{\pi_{ab}\Big|a=b(r)\right\}$$

This is indeed standard, by using the categorical axioms.

\medskip

(2) We have the following formula, with $P_2^r$ being as above:
$$<\mathcal C_r>=P_2^r$$

This is standard as well, by doing some diagrammatic work.

\medskip

With these results in hand, the conclusion now follows. Indeed, with $r\in\mathbb N\cup\{\infty\}$ being as above, we know from the beginning of the proof that any $\pi\in D$ satisfies:
$$\pi
\in<\mathcal C>
=<\mathcal C_r>
=P_2^r$$

We conclude from this that we have an inclusion as follows:
$$D\subset P_2^r$$

Conversely, we have as well the following inclusion:
$$P_2^r
=<\mathcal C_r>
=<\mathcal C>\subset<D>
=D$$

Thus we have $D=P_2^r$, and this finishes the proof. See \cite{twe}.
\end{proof}

As a conclusion, $\mathbb TO_N$ is indeed the ``privileged'' unitary group that we were looking for, with the remark that its arithmetic versions $\mathbb Z_rO_N$ are interesting as well.

\bigskip

It remains now to discuss the reflection group case. Here the problem is that of classifying the intermediate compact groups $H_N\subset G\subset K_N$, but the situation is more complicated than in the continuous group case, with the 2-parameter series there being now replaced by a 3-parameter series. Instead of getting into this quite technical subject, let us just formulate a basic result, explaining what the 3 parameters are:

\index{complex reflection group}

\begin{proposition}
We have compact groups $H_N\subset G\subset K_N$ as follows:
\begin{enumerate}
\item The groups $\mathbb Z_rH_N$, with $r\in\mathbb N\cup\{\infty\}$.

\item The groups $H_N^s=\mathbb Z_s\wr S_N$, with $s\in 2\mathbb N$.

\item The groups $H_N^{sd}=H_N^s\cap U_N^d$, with $d|s$ and $s\in 2\mathbb N$.
\end{enumerate}
\end{proposition}

\begin{proof}
The various constructions in the statement produce indeed closed subgroups $G\subset K_N$, and the condition $H_N\subset G$ is clearly satisfied as well.
\end{proof}

The same discussion as in the continuous case applies, the idea being that the constructions $G\to\mathbb Z_rG$ and $G\to G\cap H_N^{sd}$ can be combined, and that all this leads in principle to a 3-parameter series. All this is, however, quite technical. Fortunately, exactly as in the continuous case, a solution to these classification problems comes from the notion of easiness. We have indeed the following result, coming from \cite{bb+}, \cite{twe}:

\begin{theorem}
The following groups are easy:
\begin{enumerate}
\item $\mathbb Z_rH_N$, the corresponding category $P_{even}^r\subset P_{even}$ consisting of the partitions having the property that when flattening, we have the following global formula:
$$\#\circ=\#\bullet(r)$$

\item $H_N^s=\mathbb Z_s\wr S_N$, the corresponding category $P_{even}^{(s)}\subset P_{even}$ consisting of the partitions having the property that we have the following formula, in each block:
$$\#\circ=\#\bullet(s)$$
\end{enumerate}
In addition, the easy solutions of $H_N\subset G\subset K_N$ appear by combining these examples.
\end{theorem}

\begin{proof}
All this is well-known, the idea being as follows:

\medskip

(1) The computation here is similar to the one in the proof of Theorem 10.7, by writing the fundamental representation $u=zv$ as there.

\medskip

(2) This is something very standard and fundamental, known since the paper \cite{bb+}, and which follows from a long, routine computation, perfomed there.

\medskip

As for the last assertion, things here are quite technical, and for the precise statement and proof of the classification result, we refer here to paper \cite{twe}.
\end{proof}

Summarizing, the situation here is more complicated than in the continuous group case. However, in what regards the ``standard'' solution, this is definitely $\mathbb TH_N$.

\section*{10c. Axiom check}

With all this preliminary work done, let us turn now to our main question, namely constructing new geometries. We will be rather brief here, these ``hybrid'' geometries being mostly of theoretical interest. To start with, we have the following result:

\begin{theorem}
We have correspondences as follows,
$$
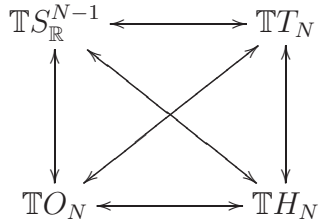
\xymatrix@R=50pt@C=50pt{
\mathbb TS^{N-1}_\mathbb R\ar[r]\ar[d]\ar[dr]&\mathbb TT_N\ar[l]\ar[d]\ar[dl]\\
\mathbb TO_N\ar[u]\ar[ur]\ar[r]&\mathbb TH_N\ar[l]\ar[ul]\ar[u]
}$$
which produce a new geometry.
\end{theorem}

\begin{proof}
We have indeed a quadruplet $(S,T,U,K)$ as in the statement, produced by the various constructions above. Regarding now the verification of the axioms:

\medskip

(1) We have the following computation:
\begin{eqnarray*}
P(\mathbb TS^{N-1}_\mathbb R\cap\mathbb T_N^+)
&=&P(\mathbb TS^{N-1}_\mathbb R\cap\mathbb T_N)\\
&\subset&P\mathbb TS^{N-1}_\mathbb R\cap P\mathbb T_N\\
&=&P^{N-1}_\mathbb R\cap\mathbb T_{N-1}\\
&=&T_{N-1}
\end{eqnarray*}

By lifting, we obtain from this that we have:
$$\mathbb TS^{N-1}_\mathbb R\cap\mathbb T_N^+\subset\mathbb TT_N$$

The inclusion ``$\supset$'' being clear as well, we are done with checking the first axiom.

\medskip

(2) The second axiom states that we must have the following equality:
$$\mathbb TH_N\cap\mathbb T_N^+=\mathbb TT_N$$

But the verification here is similar to the previous verification, for the spheres.

\medskip

(3) The third axiom states that we must have the following equality:
$$\mathbb TO_N\cap K_N^+=\mathbb TH_N$$

But this can be checked either directly, or by proceeding as above, by taking first projective versions, and then lifting.

\medskip

(4) The quantum isometry group axiom states that we must have:
$$G^+(\mathbb TS^{N-1}_\mathbb R)=\mathbb TO_N$$

But the verification here is routine, and this is explained for instance in \cite{ba4}.

\medskip

(5) The quantum reflection group axiom states that we must have:
$$G^+(\mathbb TT_N)\cap K_N^+=\mathbb TH_N$$

But this can be checked in a similar way, by adapting previous computations.

\medskip

(6) Regarding now the hard liberation axiom, this is clear, because we have:
\begin{eqnarray*}
<O_N,\mathbb TT_N>
&=&<O_N,\mathbb T,T_N>\\
&=&<O_N,\mathbb T>\\
&=&\mathbb TO_N
\end{eqnarray*}

(7) Finally, the last axiom, namely $S_{\mathbb TO_N}=\mathbb TS^{N-1}_\mathbb R$, is clear from definitions.
\end{proof}

Let us discuss now the half-classical and free extensions of Theorem 10.11, and of some of the results preceding it. In order to have no redundant discussion and diagrams, later on, we will talk directly about the $\times9$ extension of the theory that we have so far.  We first need to complete our collection of spheres $S$, tori $T$, unitary groups $U$, and reflection groups $K$. In what regards the spheres, the result is as follows:

\begin{proposition}
We have noncommutative spheres as follows,
$$\xymatrix@R=10mm@C=9mm{
S^{N-1}_{\mathbb R,+}\ar[r]&\mathbb TS^{N-1}_{\mathbb R,+}\ar[r]&S^{N-1}_{\mathbb C,+}\\
S^{N-1}_{\mathbb R,*}\ar[r]\ar[u]&\mathbb TS^{N-1}_{\mathbb R,*}\ar[r]\ar[u]&S^{N-1}_{\mathbb C,*}\ar[u]\\
S^{N-1}_\mathbb R\ar[r]\ar[u]&\mathbb TS^{N-1}_\mathbb R\ar[r]\ar[u]&S^{N-1}_\mathbb C\ar[u]}$$
with the middle vertical objects coming via the relations $ab^*=a^*b$.
\end{proposition}

\begin{proof}
We can indeed construct new spheres via the relations $ab^*=a^*b$, and these fit into previous 6-diagram of spheres as indicated. As for the fact that in the classical case we obtain the previously constructed sphere $\mathbb TS^{N-1}_\mathbb R$, this follows from Theorem 10.1 and its proof, because the relations used there are precisely those of type $a\bar{b}=\bar{a}b$.
\end{proof}

There are many things that can be done with the above spheres. As a basic result here, let us record the following fact, regarding the corresponding projective spaces:

\begin{theorem}
The projective spaces associated to the basic spheres are
$$\xymatrix@R=9mm@C=8mm{
P^{N-1}_+\ar[r]&P^{N-1}_+\ar[r]&P^{N-1}_+\\
P^{N-1}_\mathbb C\ar[r]\ar[u]&P^{N-1}_\mathbb C\ar[r]\ar[u]&P^{N-1}_\mathbb C\ar[u]\\
P^{N-1}_\mathbb R\ar[r]\ar[u]&P^{N-1}_\mathbb R\ar[r]\ar[u]&P^{N-1}_\mathbb C\ar[u]}$$
via the standard identifications for noncommutative algebraic manifolds.
\end{theorem}

\begin{proof}
This is something that we already know for the 6 previous spheres. As for the 3 new spheres, this follows from the defining relations $ab^*=a^*b$, which tell us that the coordinates of the corresponding projective spaces must be self-adjoint. See \cite{ba4}.
\end{proof}

At the torus level now, the construction is similar, as follows:

\begin{proposition}
We have noncommutative tori as follows,
$$\xymatrix@R=10mm@C=10mm{
T_N^+\ar[r]&\mathbb TT_N^+\ar[r]&\mathbb T_N^+\\
T_N^*\ar[r]\ar[u]&\mathbb TT_N^*\ar[r]\ar[u]&\mathbb T_N^*\ar[u]\\
T_N\ar[r]\ar[u]&\mathbb TT_N\ar[r]\ar[u]&\mathbb T_N\ar[u]}$$
with the middle vertical objects coming via the relations $ab^*=a^*b$.
\end{proposition}

\begin{proof}
This is clear from Proposition 10.12, by intersecting everything with $\mathbb T_N^+$.
\end{proof}

In what regards the unitary quantum groups, the result is as follows:

\begin{theorem}
We have quantum groups as follows, which are all easy,
$$\xymatrix@R=10mm@C=10mm{
O_N^+\ar[r]&\mathbb TO_N^+\ar[r]&U_N^+\\
O_N^*\ar[r]\ar[u]&\mathbb TO_N^*\ar[r]\ar[u]&U_N^*\ar[u]\\
O_N\ar[r]\ar[u]&\mathbb TO_N\ar[r]\ar[u]&U_N\ar[u]}$$
with the middle vertical objects coming via the relations $ab^*=a^*b$.
\end{theorem}

\begin{proof}
This is standard, indeed, the categories of partitions being as follows:
$$\xymatrix@R=11mm@C=11mm{
NC_2\ar[d]&\bar{NC}_2\ar[l]\ar[d]&\mathcal{NC}_2\ar[d]\ar[l]\\
P_2^*\ar[d]&\bar{P}_2^*\ar[l]\ar[d]&\mathcal P_2^*\ar[l]\ar[d]\\
P_2&\bar{P}_2\ar[l]&\mathcal P_2\ar[l]}$$

Observe that our diagrams are both intersection diagrams. See \cite{ba4}.
\end{proof}

Regarding the quantum reflection groups, we have here:

\begin{theorem}
We have quantum groups as follows, which are all easy,
$$\xymatrix@R=11mm@C=11mm{
H_N^+\ar[r]&\mathbb TH_N^+\ar[r]&K_N^+\\
H_N^*\ar[r]\ar[u]&\mathbb TH_N^*\ar[r]\ar[u]&K_N^*\ar[u]\\
H_N\ar[r]\ar[u]&\mathbb TH_N\ar[r]\ar[u]&K_N\ar[u]}$$
with the middle vertical objects coming via the relations $ab^*=a^*b$.
\end{theorem}

\begin{proof}
This is standard, indeed, the categories of partitions being as follows:
$$\xymatrix@R=11mm@C=7mm{
NC_{even}\ar[d]&\bar{NC}_{even}\ar[l]\ar[d]&\mathcal{NC}_{even}\ar[d]\ar[l]\\
P_{even}^*\ar[d]&\bar{P}_{even}^*\ar[l]\ar[d]&\mathcal P_{even}^*\ar[l]\ar[d]\\
P_{even}&\bar{P}_{even}\ar[l]&\mathcal P_{even}\ar[l]}$$

Observe that our diagrams are both intersection diagrams. See \cite{ba4}.
\end{proof}

Let us point out that we have some interesting questions, regarding the classification of the intermediate compact quantum groups for the following 4 inclusions:
$$\xymatrix@R=20pt@C=20pt{
&K_N^+\ar[rr]&&U_N^+\\
H_N^+\ar[rr]\ar@{-->}[ur]&&O_N^+\ar@{-->}[ur]\\
&K_N^*\ar[rr]\ar[uu]&&U_N^*\ar[uu]\\
H_N^*\ar[uu]\ar@{-->}[ur]\ar[rr]&&O_N^*\ar[uu]\ar@{-->}[ur]
}$$

In what regards the half-classical questions, these can be in principle fully investigated by using the technology in \cite{bdu}, but we do not know what the final answer is. As for the free questions, these are more delicate, but in the easy case, they are solved by \cite{twe}.

\bigskip

Getting back now to the verification of the axioms, we first have:

\begin{theorem}
The quantum isometries of the basic spheres, namely
$$\xymatrix@R=11mm@C=8mm{
S^{N-1}_{\mathbb R,+}\ar[r]&\mathbb TS^{N-1}_{\mathbb R,+}\ar[r]&S^{N-1}_{\mathbb C,+}\\
S^{N-1}_{\mathbb R,*}\ar[r]\ar[u]&\mathbb TS^{N-1}_{\mathbb R,*}\ar[r]\ar[u]&S^{N-1}_{\mathbb C,*}\ar[u]\\
S^{N-1}_\mathbb R\ar[r]\ar[u]&\mathbb TS^{N-1}_\mathbb R\ar[r]\ar[u]&S^{N-1}_\mathbb C\ar[u]}$$
are the basic unitary quantum groups.
\end{theorem}

\begin{proof}
This is routine, by lifting the results that we already have. See \cite{ba4}.
\end{proof}

Regarding now the tori, we first have here:

\begin{proposition}
The quantum isometries of the basic tori are
$$\xymatrix@R=11mm@C=9mm{
H_N^+\ar[r]&\mathbb TH_N^+\ar[r]&K_N^+\\
H_N^*\ar[r]\ar[u]&\mathbb TH_N^*\ar[r]\ar[u]&K_N^*\ar[u]\\
\bar{O}_N\ar[r]\ar[u]&\mathbb T\bar{O}_N\ar[r]\ar[u]&\bar{U}_N\ar[u]}$$
with the bars denoting as usual Schur-Weyl twists.
\end{proposition}

\begin{proof}
This follows again by lifting the results that we already have, with most of the relevant computations here being available from \cite{ba3}, \cite{ba4}.
\end{proof}

By looking now at quantum reflections, we obtain:

\begin{theorem}
The quantum reflections of the tori,
$$\xymatrix@R=11mm@C=11mm{
T_N^+\ar[r]&\mathbb TT_N^+\ar[r]&\mathbb T_N^+\\
T_N^*\ar[r]\ar[u]&\mathbb TT_N^*\ar[r]\ar[u]&\mathbb T_N^*\ar[u]\\
T_N\ar[r]\ar[u]&\mathbb TT_N\ar[r]\ar[u]&\mathbb T_N\ar[u]}$$
are the basic quantum reflection groups.
\end{theorem}

\begin{proof}
This is indeed routine, by intersecting, and with the various technical results regarding the intersections being available from \cite{ba3}, \cite{ba4}.
\end{proof}

Finally, we have hard liberation results, as follows:

\begin{theorem}
We have hard liberation formulae of type
$$U=<O_N,T>$$
for all the basic unitary quantum groups.
\end{theorem}

\begin{proof}
We only need to check this for the ``hybrid'' examples, constructed in this chapter. But for these hybrid examples, $U=\mathbb TO_N^\times$, the results follow from:
\begin{eqnarray*}
\mathbb TO_N^\times
&=&<\mathbb T,O_N^\times>\\
&=&<\mathbb T,<O_N,T_N^\times>>\\
&=&<O_N,<\mathbb T,T_N^\times>>\\
&=&<O_N,\mathbb TT_N^\times>
\end{eqnarray*}

Thus, we have indeed complete hard liberation results, as claimed.
\end{proof}

We can now formulate our main result, as follows:

\begin{theorem}
We have $9$ noncommutative geometries, as follows,
$$\xymatrix@R=35pt@C=35pt{
\mathbb R^N_+\ar[r]&\mathbb T\mathbb R^N_+\ar[r]&\mathbb C^N_+\\
\mathbb R^N_*\ar[u]\ar[r]&\mathbb T\mathbb R^N_*\ar[u]\ar[r]&\mathbb C^N_*\ar[u]\\
\mathbb R^N\ar[u]\ar[r]&\mathbb T\mathbb R^N\ar[u]\ar[r]&\mathbb C^N\ar[u]
}$$
with each of the $\mathbb K^\times$ symbols standing for the corresponding quadruplet.
\end{theorem}

\begin{proof}
This follows indeed by putting everything together, a bit as in the proof of Theorem 10.11, the idea being that the intersection axioms are clear, the quantum isometry axioms follow from the above computations, and the remaining axioms are elementary. Thus, we are led to the conclusion in the statement.
\end{proof}

Summarizing, we are done with the extension program mentioned in chapter 4, and started in the previous chapter, and this with the technical remark that, in what concerns the ``hybrid'' geometries, lying between real and complex, our choice of the group $\mathbb T$ for ``multiplying the real geometries'' might be actually just the ``standard'' one, because the whole family of groups $\mathbb Z_r$ with $r<\infty$ is waiting to be investigated as well. 

\bigskip

As a second comment, it is of course possible to further develop the hybrid geometries that we found here, but the whole subject looks less interesting than, for instance, the subject of further developing the half-classical geometries. Thus, we will stop our study here, and after talking next about classification results, and then in chapter 11 about twists, we will be back in chapter 12 below to the half-classical geometries.

\section*{10d. Classification results}

Getting now into classification results, let us recall from chapter 4 that a geometry coming from a quadruplet $(S,T,U,K)$ is easy when both the quantum groups $U,K$ are easy, and when the following easy generation formula is satisfied:
$$U=\{O_N,K\}$$

Combinatorially, this leads to the following statement:

\index{easy geometry}
\index{category of partitions}
\index{category of pairings}
\index{easy generation}

\begin{proposition}
An easy geometry is uniquely determined by a pair $(D,E)$ of categories of partitions, which must be as follows,
$$\mathcal{NC}_2\subset D\subset P_2$$
$$\mathcal{NC}_{even}\subset E\subset P_{even}$$
and which are subject to the following intersection and generation conditions,
$$D=E\cap P_2$$
$$E=<D,\mathcal{NC}_{even}>$$
and to the usual axioms for the associated quadruplet $(S,T,U,K)$, where $U,K$ are respectively the easy quantum groups associated to the categories $D,E$.
\end{proposition}

\begin{proof}
This comes from the following conditions, with the first one being the one mentioned above, and with the second one being part of our general axioms: 
$$U=\{O_N,K\}$$
$$K=U\cap K_N^+$$

Indeed, $U,K$ must be easy, coming from certain categories of partitions $D,E$. It is clear that $D,E$ must appear as intermediate categories, as in the statement, and the fact that the intersection and generation conditions must be satisfied follows from:
\begin{eqnarray*}
U=\{O_N,K\}&\iff&D=E\cap P_2\\
K=U\cap K_N^+&\iff&E=<D,\mathcal{NC}_{even}>
\end{eqnarray*}

Thus, we are led to the conclusion in the statement.
\end{proof}

In order to discuss now classification results, we will need some technical results regarding the intermediate easy quantum groups as follows:
$$H_N\subset G\subset K_N^+$$
$$O_N\subset G\subset U_N^+$$

Regarding the reflection groups, the complete result known so far, from Raum-Weber \cite{rwe}, concerns only the real case. This result, in a simplified form, is as follows:

\index{Raum-Weber theorem}
\index{quantum reflection group}

\begin{theorem}
The easy quantum groups $H_N\subset G\subset H_N^+$ are as follows,
$$H_N\subset H_N^\Gamma\subset H_N^{[\infty]}\subset H_N^{[r]}\subset H_N^+$$
with the family $H_N^\Gamma$ covering $H_N,H_N^{[\infty]}$, and with the series $H_N^{[r]}$ covering $H_N^+$.
\end{theorem}

\begin{proof}
This is something quite technical, from \cite{rwe}, the idea being as follows:

\medskip

(1) We have a dichotomy concerning the quantum groups $H_N\subset G\subset H_N^+$, which must fall into one of the following two classes:
$$H_N\subset G\subset H_N^{[\infty]}\quad,\quad 
H_N^{[\infty]}\subset G\subset H_N^+$$

This comes indeed from various papers, and more specifically from the final classification paper of Raum-Weber \cite{rwe}, where the quantum groups $S_N\subset G\subset H_N^+$ with $G\not\subset H_N^{[\infty]}$ were classified, and shown to contain $H_N^{[\infty]}$. For details here, we refer to \cite{rwe}.

\medskip

(2) Regarding the first case, namely $H_N\subset G\subset H_N^{[\infty]}$, the result here, from \cite{rwe}, is quite technical. Consider a discrete group generated by real reflections, $g_i^2=1$:
$$\Gamma=<g_1,\ldots,g_N>$$

\index{uniform group}
\index{variety of groups}

We call $\Gamma$ uniform if each $\sigma\in S_N$ produces a group automorphism, as follows:
$$g_i\to g_{\sigma(i)}$$

In this case, we can associate to our group $\Gamma$ a family of subsets $D(k,l)\subset P(k,l)$, which form a category of partitions, as follows:
$$D(k,l)=\left\{\pi\in P(k,l)\Big|\ker\binom{i}{j}\leq\pi\implies g_{i_1}\ldots g_{i_k}=g_{j_1}\ldots g_{j_l}\right\}$$

Observe that we have inclusions of categories as follows, coming respectively from $\eta\in D$, and from the quotient map $\Gamma\to\mathbb Z_2^N$:
$$P_{even}^{[\infty]}\subset D\subset P_{even}$$

Conversely, to any category of partitions $P_{even}^{[\infty]}\subset D\subset P_{even}$ we can associate a uniform reflection group $\mathbb Z_2^{*N}\to\Gamma\to\mathbb Z_2^N$, as follows:
$$\Gamma=\left\langle g_1,\ldots g_N\Big|g_{i_1}\ldots g_{i_k}=g_{j_1}\ldots g_{j_l},\forall i,j,k,l,\ker\binom{i}{j}\in D(k,l)\right\rangle$$

As explained in \cite{rwe}, the correspondences $\Gamma\to D$ and $D\to\Gamma$ constructed above are bijective, and inverse to each other, at $N=\infty$. Thus, we are done with the first case.

\medskip

(3) Regarding now the second case, which is the one left, namely $H_N^{[\infty]}\subset G\subset H_N^+$, the result here, also from \cite{rwe}, is quite technical as well, but has a simple formulation. Let indeed $H_N^{[r]}\subset H_N^+$ be the easy quantum group coming from:
$$\pi_r=\ker\begin{pmatrix}1&\ldots&r&r&\ldots&1\\1&\ldots&r&r&\ldots&1\end{pmatrix}$$

We have then inclusions of quantum groups as follows:
$$H_N^+=H_N^{[1]}\supset H_N^{[2]}\supset H_N^{[3]}\supset\ldots\ldots\supset H_N^{[\infty]}$$

We obtain in this way all the intermediate easy quantum groups $H_N^{[\infty]}\subset G\subset H_N^+$ satisfying the assumption $G\neq H_N^{[\infty]}$, and this finishes the proof. See \cite{rwe}.
\end{proof}

Let us discuss now the rotation groups. Once again, there are only partial results here so far, notably with the results in Mang-Weber \cite{mwe}, concerning the following case:
$$U_N\subset G\subset U_N^+$$ 

A first construction of such quantum groups is as follows:

\index{unitary quantum group}
\index{category of pairings}
\index{category of matching pairings}
\index{cyclic matrix model}

\begin{proposition}
Associated to any $r\in\mathbb N$ is the quantum group $U_N\subset U_N^{(r)}\subset U_N^+$ coming from the category $\mathcal P_2^{(r)}$ of matching pairings having the property that 
$$\#\circ=\#\bullet(r)$$
holds between the legs of each string. These quantum groups have the following properties:
\begin{enumerate}
\item At $r=1$ we obtain the usual unitary group, $U_N^{(1)}=U_N$.

\item At $r=2$ we obtain the half-classical unitary group, $U_N^{(2)}=U_N^*$.

\item For any $r|s$ we have an embedding $U_N^{(r)}\subset U_N^{(s)}$.

\item In general, we have an embedding $U_N^{(r)}\subset U_N^r\rtimes\mathbb Z_r$.

\item We have as well a cyclic matrix model $C(U_N^{(r)})\subset M_r(C(U_N^r))$.

\item In this latter model, $\int_{U_N^{(r)}}$ appears as the restriction of $tr_r\otimes\int_{U_N^r}$.
\end{enumerate}
\end{proposition}

\begin{proof}
This is something quite compact, summarizing various findings from \cite{bb2}, \cite{mwe}. Here are a few brief explanations on all this:

\medskip

(1) This is clear from $\mathcal P_2^{(1)}=\mathcal P_2$, and from a well-known result of Brauer \cite{bra}.

\medskip

(2) This is because $\mathcal P_2^{(2)}$ is generated by the partitions with implement the relations $abc=cba$ between the variables $\{u_{ij},u_{ij}^*\}$, used in \cite{bdu} for constructing $U_N^*$.

\medskip

(3) This simply follows from $\mathcal P_2^{(s)}\subset\mathcal P_2^{(r)}$, by functoriality.

\medskip

(4) This is the original definition of $U_N^{(r)}$, from \cite{bb2}. We refer to \cite{bb2} for the formula of the embedding, and to \cite{mwe} for the compatibility with the Tannakian definition.

\medskip

(5) This is also from \cite{bb2}, more specifically it is an alternative definition for $U_N^{(r)}$.

\medskip

(6) Once again, this is something from \cite{bb2}, and we will be back to it.
\end{proof}

Let us discuss now the second known construction of unitary quantum groups, from \cite{mwe}. This construction uses an additive semigroup $D\subset\mathbb N$, but as pointed out there, using instead the complementary set $C=\mathbb N-D$ leads to several simplifications. So, let us call ``cosemigroup'' any subset $C\subset\mathbb N$ which is complementary to an additive semigroup, $x,y\notin C\implies x+y\notin C$. The construction from \cite{mwe} is then:

\index{unitary quantum group}
\index{category of pairings}
\index{category of matching pairings}
\index{cosemigroup}

\begin{proposition}
Associated to any cosemigroup $C\subset\mathbb N$ is the easy quantum group $U_N\subset U_N^C\subset U_N^+$ coming from the category $\mathcal P_2^C\subset P_2^{(\infty)}$ of pairings having the property 
$$\#\circ-\#\bullet\in C$$
between each two legs colored $\circ,\bullet$ of two strings which cross. We have:
\begin{enumerate}
\item For $C=\emptyset$ we obtain the quantum group $U_N^+$.

\item For $C=\{0\}$ we obtain the quantum group $U_N^\times$.

\item For $C=\{0,1\}$ we obtain the quantum group $U_N^{**}$.

\item For $C=\mathbb N$ we obtain the quantum group $U_N^{(\infty)}$.

\item For $C\subset C'$ we have an inclusion $U_N^{C'}\subset U_N^C$.

\item Each quantum group $U_N^C$ contains each quantum group $U_N^{(r)}$.
\end{enumerate}
\end{proposition}

\begin{proof}
Once again this is something very compact, coming from recent work in \cite{mwe}, with our convention that the semigroup $D\subset\mathbb N$ which is used there is replaced here by its complement $C=\mathbb N-D$. Here are a few explanations on all this:

\medskip

(1) The assumption $C=\emptyset$ means that the condition $\#\circ-\#\bullet\in C$ can never be applied. Thus, the strings cannot cross, we have $\mathcal P_2^\emptyset=\mathcal{NC}_2$, and so $U_N^\emptyset=U_N^+$.

\medskip

(2) As explained in \cite{mwe}, here we obtain indeed the quantum group $U_N^\times$, constructed by using the relations $ab^*c=cb^*a$, with $a,b,c\in\{u_{ij}\}$. 

\medskip

(3) This is also explained in \cite{mwe}, with $U_N^{**}$ being the quantum group from \cite{bb2}, which is the biggest whose full projective version, in the sense there, is classical. 

\medskip

(4) Here the assumption $C=\mathbb N$ simply tells us that the condition $\#\circ-\#\bullet\in C$ in the statement is irrelevant. Thus, we have $\mathcal P_2^\mathbb N=\mathcal P_2^{(\infty)}$, and so $U_N^\mathbb N=U_N^{(\infty)}$.

\medskip

(5) This is clear by functoriality, because $C\subset C'$ implies $\mathcal P_2^{C}\subset\mathcal P_2^{C'}$.

\medskip

(6) This is clear from definitions, and from Proposition 10.24.
\end{proof}

We have the following key result, from Mang-Weber \cite{mwe}:

\index{Mang-Weber theorem}

\begin{theorem}
The easy quantum groups $U_N\subset G\subset U_N^+$ are as follows,
$$U_N\subset\{U_N^{(r)}\}\subset\{U_N^C\}\subset U_N^+$$
with the series covering $U_N$, and the family covering $U_N^+$.
\end{theorem}

\begin{proof}
This is something non-trivial, and we refer here to \cite{mwe}. The general idea is that $U_N^{(\infty)}$ produces a dichotomy for the quantum groups in the statement, as follows:
$$U_N\subset G\subset U_N^{(\infty)}\quad,\quad U_N^{(\infty)}\subset G\subset U_N^+$$

But this leads, via combinatorics, to the series and the family. See \cite{mwe}.
\end{proof}

Observe that there is an obvious similarity here with the dichotomy for the liberations of $H_N$, coming from the work of Raum-Weber \cite{rwe}, explained in the above. To be more precise, the above-mentioned classification results for the liberations of $H_N$ and the liberations of $U_N$ have some obvious similarity between them. We have indeed a family followed by a series, and a series followed by a family. All this suggests the existence of a general ``contravariant duality'' between these quantum groups, as follows:
$$\xymatrix@R=50pt@C=50pt{
U_N\ar[r]\ar@.[d]&U_N^{(r)}\ar[r]\ar@.[d]&U_N^C\ar[r]\ar@.[d]&U_N^+\ar@.[d]\\
H_N^+\ar@.[u]&H_N^{[r]}\ar[l]\ar@.[u]&H_N^\Gamma\ar[l]\ar@.[u]&H_N\ar[l]\ar@.[u]
}$$

At the first glance, this might sound a bit strange. Indeed, we have some natural and well-established correspondences $H_N\leftrightarrow U_N$ and $H_N^+\leftrightarrow U_N^+$, obtained in one sense by taking the real reflection subgroup, $H=U\cap H_N^+$, and in the other sense by setting $U=<H,U_N>$. Thus, our proposal of duality seems to go the wrong way. On the other hand, obvious as well is the fact that these correspondences $H_N\leftrightarrow U_N$ and $H_N^+\leftrightarrow U_N^+$ cannot be extended as to map the series to the series, and the family to the family, because the series/families would have to be ``inverted'', in order to do so. Thus, we are led to the above contravariant duality conjecture, which looks like something quite complicated.

\bigskip

Now back to our abstract noncommutative geometries, as axiomatized here, in the easy case we have the following classification result, based on the above:

\index{uniformity}
\index{purity}
\index{slicing}

\begin{theorem}
There are exactly $4$ geometries which are easy, uniform and pure, with purity meaning that the geometry must be real, classical, complex or free, namely:
$$\xymatrix@R=50pt@C=50pt{
\mathbb R^N_+\ar[r]&\mathbb C^N_+\\
\mathbb R^N\ar[u]\ar[r]&\mathbb C^N\ar[u]
}$$
When lifting the uniformity and purity conditions, and replacing them with a ``slicing'' axiom, we have $9$ such geometries, namely those in Theorem 10.21.
\end{theorem}

\begin{proof}
All this is quite technical, the idea being as follows:

\medskip

(1) Assume first that we have an easy geometry which is pure, in the sense that it lies on one of the 4 edges of the square in the statement. We know from Proposition 10.22 that its unitary group $U$ must come from a category of pairings $D$ satisfying:
$$D=<D,\mathcal{NC}_{even}>\cap  P_2$$

But this equation can be solved by using the results in \cite{mwe}, \cite{rwe}, \cite{twe}, and by using the uniformity axiom, which excludes the half-liberations and the hybrids, we are led to the conclusion that the only solutions are the 4 vertices of the square.

\medskip

(2) Regarding the second assertion, this can be obtained via the same easiness technology, by using the ``slicing'' axiom from \cite{ba8}, which amounts in saying that $U$, or the geometry itself, can be reconstructed from its projections on the edges of the square. All this is quite technical, again, and for details on all this, we refer to \cite{ba8}.
\end{proof}

As a conclusion to all this, we have now a much better understanding of our axioms from chapter 4, and also, generally speaking, of what we have been trying do do, since the beginning of this book. Indeed, our $(S,T,U,K)$ formalism appears to be something quite reasonable, corresponding to the natural thought that there should be 4 main geometries, namely classical/free, real/complex, and that there might be perhaps a few more  geometries, obtained by replacing the commutation relations $ab=ba$ with something ``clever''. With all the above, we have now confirmation for all this. Business doing fine.

\section*{10e. Exercises} 

Things have been quite technical in this chapter, and as unique exercise, we have:

\begin{exercise}
Find a better way of classifying the noncommutative geometries in our sense, say by adding some clever extra axiom, which simplifies the classification.
\end{exercise}

Needless to say, an answer here would be very interesting.

\chapter{Twisted geometry}

\section*{11a. Ad-hoc twists}

We have seen so far that the abstract noncommutative geometries, taken in a ``spherical'' sense, with coordinates bounded by $||x_i||\leq1$, can be axiomatized with the help of quadruplets $(S,T,U,K)$. There are 9 main such geometries, as follows:
$$\xymatrix@R=35pt@C=38pt{
\mathbb R^N_+\ar[r]&\mathbb T\mathbb R^N_+\ar[r]&\mathbb C^N_+\\
\mathbb R^N_*\ar[u]\ar[r]&\mathbb T\mathbb R^N_*\ar[u]\ar[r]&\mathbb C^N_*\ar[u]\\
\mathbb R^N\ar[u]\ar[r]&\mathbb T\mathbb R^N\ar[u]\ar[r]&\mathbb C^N\ar[u]
}$$

An important question, that we would like to investigate now, concerns the twisting  of these geometries, by suitably replacing commutation with anticommutation: 
$$ab=ba\quad\to\quad ab=\pm ba$$
$$abc=cba\quad\to\quad abc=\pm cba$$

We will see that  this is possible, and that we have twisted geometries, as follows:
$$\xymatrix@R=35pt@C=36pt{
\mathbb R^N_+\ar[r]&\mathbb T\mathbb R^N_+\ar[r]&\mathbb C^N_+\\
\bar{\mathbb R}^N_*\ar[u]\ar[r]&\mathbb T\bar{\mathbb R}^N_*\ar[u]\ar[r]&\bar{\mathbb C}^N_*\ar[u]\\
\bar{\mathbb R}^N\ar[u]\ar[r]&\mathbb T\bar{\mathbb R}^N\ar[u]\ar[r]&\bar{\mathbb C}^N\ar[u]
}$$

Here the bars stand, as before in this book, for anticommutation twists. However, all this is quite tricky, and before starting, a few general comments:

\bigskip

(1) First of all, our work is motivated by the general commutation/anticommutation duality from quantum mechanics, in a general sense. If there is one thing to be learned from basic quantum mechanics, and let us recommend here again our favorite books, namely Dirac \cite{dir}, Feynman \cite{fey}, Griffiths \cite{gr2}, von Neumann \cite{von}, Weinberg \cite{we2}, Weyl \cite{wey}, this is the fact that there is no commutation without anticommutation.

\bigskip

(2) Mathematically, and in relation with what we have been doing so far here, we have already met $q=-1$ twists on several occasions, and notably in relation with the computation of the quantum isometry groups $G^+(T)$ of our tori $T$, with one of our 7 noncommutative geometry axioms stating that we must have $K=G^+(T)\cap K_N^+$. And the point is that $G^+(T)$, quite surpringly, often happens to be a $q=-1$ twist.

\bigskip

(3) And there are countless other reasons, both mathematical and physical, to look at anticommutation and $q=-1$ twists. If you are a bit familiar for instance with Drinfeld-Jimbo \cite{dri}, \cite{jim}, you probably know that many geometric objects can be deformed with the help of a parameter $q\in\mathbb C$, the interesting case being $q\in\mathbb T$, and more specifically the case where $q$ is a root of unity. And aren't $q=\pm1$ the simplest roots of unity.

\bigskip

(4) Summarizing, we have motivations. However, when getting to work, several surprises are waiting for us. First if the fact that the Drinfeld-Jimbo deformation procedure \cite{dri}, \cite{jim} produces non-semisimple objects at roots of unity $q\neq1$, and in particular at $q=-1$. Thus, this very popular theory is useless for us, not to say wrong in our opinion, and we must come up with new definitions for everything, at $q=-1$.

\bigskip

(5) Fortunately, this is possible, with the correct objects at $q=-1$ having emerged, in a somewhat discreet way, not to contradict much the popular belief, and get sent to the Inquisition or something, in a number of technical papers, all peer-reviewed and published, on quantum groups and noncommutative geometry, all over the late 00s and 10s, starting with \cite{bbc} which launched everything, with the correct twist of $O_N$.

\bigskip

(6) And so, getting back now to the question of twisting the 9 geometries that we have, this is definitely possible, thanks to all this underground, while ironically public, $q=-1$ knowledge accumulated over the years, and we will explain this, in this chapter. With the technical remark that the twisted geometries do not exactly satisfy our axioms from chapter 4, but are not far from them either, and we will comment on this.

\bigskip

(7) Finally, this chapter will be a modest introduction to all this. The geometries of $\bar{\mathbb R}^N,\bar{\mathbb C}^N$ for instance are potentially as wide as those of $\mathbb R^N,\mathbb C^N$, and with many classical techniques applying well, and there is certainly room for writing a book on this topic, ``twisted geometry''. Let me also mention that, in the lack of such a book, you can always ask my colleague Bichon about such things, he's the one who knows.

\bigskip

In order to get started now, the best is to deform first the simplest objects that we have, namely the quantum spheres. This can be done as follows:

\index{twisting}
\index{twisted sphere}
\index{ad-hoc twisting}
\index{anticommutation}

\begin{theorem}
We have quantum spheres as follows, obtained via the twisted commutation relations $ab=\pm ba$, and twisted half-commutation relations $abc=\pm cba$,
$$\xymatrix@R=13mm@C=10mm{
S^{N-1}_{\mathbb R,+}\ar[r]&\mathbb TS^{N-1}_{\mathbb R,+}\ar[r]&S^{N-1}_{\mathbb C,+}\\
\bar{S}^{N-1}_{\mathbb R,*}\ar[r]\ar[u]&\mathbb T\bar{S}^{N-1}_{\mathbb R,*}\ar[r]\ar[u]&\bar{S}^{N-1}_{\mathbb C,*}\ar[u]\\
\bar{S}^{N-1}_\mathbb R\ar[r]\ar[u]&\mathbb T\bar{S}^{N-1}_\mathbb R\ar[r]\ar[u]&\bar{S}^{N-1}_\mathbb C\ar[u]}$$
with the precise signs being as follows:
\begin{enumerate}
\item The signs on the bottom correspond to anticommutation of distinct coordinates, and their adjoints. That is, with $z_i=x_i,x_i^*$ and $\varepsilon_{ij}=1-\delta_{ij}$, the formula is:
$$z_iz_j=(-1)^{\varepsilon_{ij}}z_jz_i$$

\item The signs in the middle come from functoriality, as for the spheres in the middle to contain those on the bottom. That is, $z_iz_jz_k=(-1)^{\varepsilon_{ij}+\varepsilon_{jk}+\varepsilon_{ik}}z_kz_jz_i$.
\end{enumerate}
\end{theorem}

\begin{proof}
This is something elementary, from \cite{ba1}, the idea being as follows:

\medskip

(1) Here there is nothing to prove, because we can define the spheres on the bottom by the following formulae, with $z_i=x_i,x_i^*$ and $\varepsilon_{ij}=1-\delta_{ij}$ being as above:
$$C(\bar{S}^{N-1}_\mathbb R)=C(S^{N-1}_{\mathbb R,+})\Big/\Big<x_ix_j=(-1)^{\varepsilon_{ij}}x_jx_i\Big>$$
$$C(\bar{S}^{N-1}_\mathbb C)=C(S^{N-1}_{\mathbb C,+})\Big/\Big<z_iz_j=(-1)^{\varepsilon_{ij}}z_jz_i\Big>$$

(2) Here our claim is that, if we want to construct half-classical twisted spheres, via relations of type $abc=\pm cba$ between the coordinates $x_i$ and their adjoints $x_i^*$, as for these spheres to contain the twisted spheres constructed in (1), the only possible choice for these relations is as follows, with $z_i=x_i,x_i^*$ and $\varepsilon_{ij}=1-\delta_{ij}$ being as above:
$$z_iz_jz_k=(-1)^{\varepsilon_{ij}+\varepsilon_{jk}+\varepsilon_{ik}}z_kz_jz_i$$

But this is something clear, coming from the following computation, inside of the quotient algebras corresponding to the twisted spheres constructed in (1) above:
\begin{eqnarray*}
z_iz_jz_k
&=&(-1)^{\varepsilon_{ij}}z_jz_iz_k\\
&=&(-1)^{\varepsilon_{ij}+\varepsilon_{ik}}z_jz_kz_i\\
&=&(-1)^{\varepsilon_{ij}+\varepsilon_{jk}+\varepsilon_{ik}}z_kz_jz_i
\end{eqnarray*}

Thus, we are led to the conclusion in the statement, the spheres being given by:
$$C(\bar{S}^{N-1}_{\mathbb R,*})=C(S^{N-1}_{\mathbb R,+})\Big/\Big<x_ix_jx_k=(-1)^{\varepsilon_{ij}+\varepsilon_{jk}+\varepsilon_{ik}}x_kx_jx_i\Big>$$
$$C(\bar{S}^{N-1}_{\mathbb C,*})=C(S^{N-1}_{\mathbb C,+})\Big/\Big<z_iz_jz_k=(-1)^{\varepsilon_{ij}+\varepsilon_{jk}+\varepsilon_{ik}}z_kz_jz_i\Big>$$

Thus, we have constructed our spheres, and embeddings, as desired.
\end{proof}

Let us twist now the unitary quantum groups $U$. We would like these to act on the corresponding spheres, $U\curvearrowright S$. Thus, we would like to have morphisms, as follows:
$$\Phi(x_i)=\sum_jx_j\otimes u_{ji}$$

Now with $z_i=x_i,x_i^*$ being as before, and with $v_{ij}=u_{ij},u_{ij}^*$ constructed accordingly, the above formula and its adjoint tell us that we must have:
$$\Phi(z_i)=\sum_jz_j\otimes v_{ji}$$

Thus the variables $Z_i=\sum_jz_j\otimes v_{ji}$ on the right must satisfy the twisted commutation or half-commutation relations in Theorem 11.1, and this will lead us to the correct twisted commutation or half-commutation relations to be satisfied by the variables $v_{ij}$. In practice now, let us first discuss the twisting of $O_N,U_N$. Following \cite{bbc} in the orthogonal case, and then \cite{ba1} in the unitary case, the result here is as follows:

\index{twisting}
\index{twisted orthogonal group}
\index{twisted unitary group}
\index{ad-hoc twisting}
\index{anticommutation}

\begin{theorem}
We have twisted orthogonal and unitary groups, as follows,
$$\xymatrix@R=15mm@C=15mm{
O_N^+\ar[r]&U_N^+\\
\bar{O}_N\ar[r]\ar[u]&\bar{U}_N\ar[u]}$$
defined via the following relations, with the convention $\alpha=a,a^*$ and $\beta=b,b^*$:
$$\alpha\beta=\begin{cases}
-\beta\alpha&{\rm for}\ a,b\in\{u_{ij}\}\ {\rm distinct,\ on\ the\ same\ row\ or\ column\ of\ }u\\
\beta\alpha&{\rm otherwise}
\end{cases}$$
These quantum groups act on the corresponding twisted real and complex spheres.
\end{theorem}

\begin{proof}
Let us first discuss the construction of the quantum group $\bar{O}_N$. We must prove that the algebra $C(\bar{O}_N)$ obtained from $C(O_N^+)$ via the relations in the statement has a comultiplication $\Delta$, a counit $\varepsilon$, and an antipode $S$. Regarding $\Delta$, let us set:
$$U_{ij}=\sum_ku_{ik}\otimes u_{kj}$$

For $j\neq k$ we have the following computation:
\begin{eqnarray*}
U_{ij}U_{ik}
&=&\sum_{s\neq t}u_{is}u_{it}\otimes u_{sj}u_{tk}+\sum_su_{is}u_{is}\otimes u_{sj}u_{sk}\\
&=&\sum_{s\neq t}-u_{it}u_{is}\otimes u_{tk}u_{sj}+\sum_su_{is}u_{is}\otimes(-u_{sk}u_{sj})\\
&=&-U_{ik}U_{ij}
\end{eqnarray*}

Also, for $i\neq k,j\neq l$ we have the following computation:
\begin{eqnarray*}
U_{ij}U_{kl}
&=&\sum_{s\neq t}u_{is}u_{kt}\otimes u_{sj}u_{tl}+\sum_su_{is}u_{ks}\otimes u_{sj}u_{sl}\\
&=&\sum_{s\neq t}u_{kt}u_{is}\otimes u_{tl}u_{sj}+\sum_s(-u_{ks}u_{is})\otimes(-u_{sl}u_{sj})\\
&=&U_{kl}U_{ij}
\end{eqnarray*}

Thus, we can define a comultiplication map for $C(\bar{O}_N)$, by setting:
$$\Delta(u_{ij})=U_{ij}$$

Regarding now the counit $\varepsilon$ and the antipode $S$, things are clear here, by using the same method, and with no computations needed, the formulae to be satisfied being trivially satisfied. We conclude that $\bar{O}_N$ is a compact quantum group, and the proof for $\bar{U}_N$ is similar, by adding $*$ exponents everywhere in the above computations. 

\medskip

Finally, the last assertion is clear too, by doing some elementary computations, of the same type as above, and with the remark that the converse holds too, in the sense that if we want a quantum group $U\subset U_N^+$ to be defined by relations of type $ab=\pm ba$, and to have an action $U\curvearrowright S$ on the corresponding twisted sphere, we are led to the relations in the statement. We refer to \cite{ba1} for further details on all this.
\end{proof}

In order to discuss now the half-classical case, given three coordinates $a,b,c\in\{u_{ij}\}$, let us set $span(a,b,c)=(r,c)$, where $r,c\in\{1,2,3\}$ are the number of rows and columns spanned by $a,b,c$. In other words, if we write $a=u_{ij},b=u_{kl},c=u_{pq}$ then $r=\#\{i,k,p\}$ and $l=\#\{j,l,q\}$. With this convention, we have the following result:

\index{twisting}
\index{twisted half-classical orthogonal group}
\index{twisted half-classical unitary group}
\index{ad-hoc twisting}
\index{anticommutation}

\begin{theorem}
We have intermediate quantum groups as follows,
$$\xymatrix@R=12mm@C=12mm{
O_N^+\ar[r]&\mathbb TO_N^+\ar[r]&U_N^+\\
\bar{O}_N^*\ar[r]\ar[u]&\mathbb T\bar{O}_N^*\ar[r]\ar[u]&\bar{U}_N^*\ar[u]\\
\bar{O}_N\ar[r]\ar[u]&\mathbb T\bar{O}_N\ar[r]\ar[u]&\bar{U}_N\ar[u]}$$
defined via the following relations, with $\alpha=a,a^*$, $\beta=b,b^*$ and $\gamma=c,c^*$,
$$\alpha\beta\gamma=\begin{cases}
-\gamma\beta\alpha&{\rm for}\ a,b,c\in\{u_{ij}\}\ {\rm with}\ span(a,b,c)=(\leq 2,3)\ {\rm or}\ (3,\leq 2)\\
\gamma\beta\alpha&{\rm otherwise}
\end{cases}$$
which act on the corresponding twisted half-classical real and complex spheres.
\end{theorem}

\begin{proof}
We use the same method as for Theorem 11.2, but with the combinatorics being now more complicated. Observe first that the rules for the various commutation and anticommutation signs in the statement can be summarized as follows:
$$\begin{matrix}
r\backslash c&1&2&3\\
1&+&+&-\\
2&+&+&-\\
3&-&-&+
\end{matrix}$$

Let us first prove the result for $\bar{O}_N^*$. We must construct here morphisms $\Delta,\varepsilon,S$, and the proof, similar to the proof of Theorem 11.2, goes as follows:

\medskip

(1) We first construct $\Delta$. For this purpose, we must prove that $U_{ij}=\sum_ku_{ik}\otimes u_{kj}$ satisfy the relations in the statement. We have the following computation:
\begin{eqnarray*}
U_{ia}U_{jb}U_{kc}
&=&\sum_{xyz}u_{ix}u_{jy}u_{kz}\otimes u_{xa}u_{yb}u_{zc}\\
&=&\sum_{xyz}\pm u_{kz}u_{jy}u_{ix}\otimes\pm u_{zc}u_{yb}u_{xa}\\
&=&\pm U_{kc}U_{jb}U_{ia}
\end{eqnarray*}

We must show that, when examining the precise two $\pm$ signs in the middle formula, their product produces the correct $\pm$ sign at the end. But the point is that both these signs depend only on $s=span(x,y,z)$, and for $s=1,2,3$ respectively, we have:

\medskip

-- For a $(3,1)$ span we obtain $+-$, $+-$, $-+$, so a product $-$ as needed.

\smallskip

-- For a $(2,1)$ span we obtain $++$, $++$, $--$, so a product $+$ as needed.

\smallskip

-- For a $(3,3)$ span we obtain $--$, $--$, $++$, so a product $+$ as needed.

\smallskip

-- For a $(3,2)$ span we obtain $+-$, $+-$, $-+$, so a product $-$ as needed.

\smallskip

-- For a $(2,2)$ span we obtain $++$, $++$, $--$, so a product $+$ as needed.

\medskip

Together with the fact that our problem is invariant under $(r,c)\to(c,r)$, and with the fact that for a $(1,1)$ span there is nothing to prove, this finishes the proof for $\Delta$.

\medskip

(2) The construction of the counit, via the formula $\varepsilon(u_{ij})=\delta_{ij}$, requires the Kronecker symbols $\delta_{ij}$ to commute/anticommute according to the above table. Equivalently, we must prove that the situation $\delta_{ij}\delta_{kl}\delta_{pq}=1$ can appear only in a case where the above table indicates ``+''. But this is clear, because $\delta_{ij}\delta_{kl}\delta_{pq}=1$ implies $r=c$.

\medskip

(3) Finally, the construction of the antipode, via the formula $S(u_{ij})=u_{ji}$, is clear too, because this requires the choice of our $\pm$ signs to be invariant under transposition, and this is true, the above table being symmetric. 

\medskip

We conclude that $\bar{O}_N^*$ is indeed a compact quantum group, and the proof for $\bar{U}_N^*$ is similar, by adding $*$ exponents everywhere in the above. Finally, the last assertion is clear too, exactly as in the proof of Theorem 11.2. We refer to \cite{ba1} for details.
\end{proof}

The above results can be summarized as follows:

\begin{theorem}
We have quantum groups as follows, obtained via the twisted commutation relations $ab=\pm ba$, and twisted half-commutation relations $abc=\pm cba$,
$$\xymatrix@R=12mm@C=12mm{
O_N^+\ar[r]&\mathbb TO_N^+\ar[r]&U_N^+\\
\bar{O}_N^*\ar[r]\ar[u]&\mathbb T\bar{O}_N^*\ar[r]\ar[u]&\bar{U}_N^*\ar[u]\\
\bar{O}_N\ar[r]\ar[u]&\mathbb T\bar{O}_N\ar[r]\ar[u]&\bar{U}_N\ar[u]}$$
with the various signs coming as follows:
\begin{enumerate}
\item The signs for $\bar{O}_N$ correspond to anticommutation of distinct entries on rows and columns, and commutation otherwise, with this coming from $\bar{O}_N\curvearrowright\bar{S}^{N-1}_\mathbb R$.

\item The signs for $\bar{O}_N^*,\bar{U}_N,\bar{U}_N^*$ come as well from the signs for $\bar{S}^{N-1}_\mathbb R$, either via the requirement $\bar{O}_N\subset U$, or via the requirement $U\curvearrowright S$. 
\end{enumerate}
\end{theorem}

\begin{proof}
This is a summary of Theorem 11.2 and Theorem 11.3, and their proofs.
\end{proof}

Moving ahead now, and back to our geometric program, we have twisted the spheres and unitary groups $S,U$, and we are left with twisting the tori and reflection groups $T,K$. But these are ``discrete'' objects, which can only be rigid, so let us formulate:

\index{twisted torus}
\index{twisted reflection group}
\index{rigidity}

\begin{definition}
The twists of the basic quantum tori and reflection groups,
$$\xymatrix@R=11mm@C=11mm{
T_N^+\ar[r]&\mathbb TT_N^+\ar[r]&\mathbb T_N^+\\
T_N^*\ar[r]\ar[u]&\mathbb TT_N^*\ar[r]\ar[u]&\mathbb T_N^*\ar[u]\\
T_N\ar[r]\ar[u]&\mathbb TT_N\ar[r]\ar[u]&\mathbb T_N\ar[u]}
\qquad\quad\qquad
\xymatrix@R=11mm@C=11mm{
H_N^+\ar[r]&\mathbb TH_N^+\ar[r]&K_N^+\\
H_N^*\ar[r]\ar[u]&\mathbb TH_N^*\ar[r]\ar[u]&K_N^*\ar[u]\\
H_N\ar[r]\ar[u]&\mathbb TH_N\ar[r]\ar[u]&K_N\ar[u]}$$
are by definition these tori and reflection groups themselves.
\end{definition}

With this definition in hand, we are done with our twisting program for the triples $(S,T,U,K)$, and we have now candidates $\bar{\mathbb R}^N$, $\bar{\mathbb C}^N$ and $\bar{\mathbb R}^N_*$, $\bar{\mathbb C}^N_*$ for new noncommutative geometries, to be checked from our axiomatic viewpoint, and then developed.

\section*{11b. Schur-Weyl twists}

In order to discuss these questions, we must first review the above construction of the twists of $S,T,U,K$, which was something quite ad-hoc, and replace this by something more conceptual. We use easiness. Let us start with something that we know, namely:

\index{partitions with even blocks}

\begin{proposition}
The intermediate easy quantum groups 
$$H_N\subset G\subset U_N^+$$
come via Tannakian duality from the intermediate categories of partitions
$$P_{even}\supset D\supset\mathcal {NC}_2$$
with $P_{even}(k,l)\subset P(k,l)$ being the category of partitions whose blocks have even size.
\end{proposition}

\begin{proof}
This is something coming from the general easiness theory for quantum groups, discussed in chapter 2. Indeed, as explained there, the easy quantum groups appear as certain intermediate compact quantum groups, as follows:
$$S_N\subset G\subset U_N^+$$

To be more precise, such a quantum group is easy when the corresponding Tannakian category comes from an intermediate category of partitions, as follows:
$$P\supset D\supset\mathcal {NC}_2$$

Now since this correspondence makes correspond $H_N\leftrightarrow P_{even}$, once again as explained in chapter 2, we are led to the conclusion in the statement.
\end{proof}

The idea now will be that the twisting operation $G\to\bar{G}$, in the easy case, can be implemented, via Tannakian duality as usual, via a signature operation on $P_{even}$. Given a partition $\tau\in P(k,l)$, let us call ``switch'' the operation which consists in switching two neighbors, belonging to different blocks, in the upper row, or in the lower row. Also, we use the standard embedding $S_k\subset P_2(k,k)$, via the pairings having only up-to-down strings. With these conventions, we have the following result, from \cite{ba1}:

\index{signature map}

\begin{theorem}
There is a signature map $\varepsilon:P_{even}\to\{-1,1\}$, given by 
$$\varepsilon(\tau)=(-1)^c$$
where $c$ is the number of switches needed to make $\tau$ noncrossing. In addition:
\begin{enumerate}
\item For $\tau\in S_k$, this is the usual signature.

\item For $\tau\in P_2$ we have $(-1)^c$, where $c$ is the number of crossings.

\item For $\tau\leq\pi\in NC_{even}$, the signature is $1$.
\end{enumerate}
\end{theorem}

\begin{proof}
In order to show that the signature map $\varepsilon:P_{even}\to\{-1,1\}$ in the statement, given by $\varepsilon(\tau)=(-1)^c$, is well-defined, we must prove that the number $c$ in the statement is well-defined modulo 2. It is enough to perform the verification for the noncrossing partitions. More precisely, given $\tau,\tau'\in NC_{even}$ having the same block structure, we must prove that the number of switches $c$ required for the passage $\tau\to\tau'$ is even.

\medskip

In order to do so, observe that any partition $\tau\in P(k,l)$ can be put in ``standard form'', by ordering its blocks according to the appearence of the first leg in each block, counting clockwise from top left, and then by performing the switches as for block 1 to be at left, then for block 2 to be at left, and so on. Here the required switches are also uniquely determined, by the order coming from counting clockwise from top left. 

\medskip

Here is an example of such an algorithmic switching operation, with block 1 being first put at left, by using two switches, then with block 2 left unchanged, and then with block 3 being put at left as well, but at right of blocks 1 and 2, with one switch:
$$
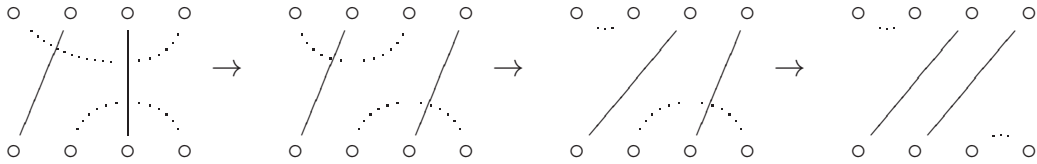
\xymatrix@R=3mm@C=3mm{\circ\ar@/_/@{.}[drr]&\circ\ar@{-}[dddl]&\circ\ar@{-}[ddd]&\circ\\
&&\ar@/_/@{.}[ur]&\\
&&\ar@/^/@{.}[dr]&\\
\circ&\circ\ar@/^/@{.}[ur]&\circ&\circ}
\xymatrix@R=4mm@C=1mm{&\\\to\\&\\& }
\xymatrix@R=3mm@C=3mm{\circ\ar@/_/@{.}[dr]&\circ\ar@{-}[dddl]&\circ&\circ\ar@{-}[dddl]\\
&\ar@/_/@{.}[ur]&&\\
&&\ar@/^/@{.}[dr]&\\
\circ&\circ\ar@/^/@{.}[ur]&\circ&\circ}
\xymatrix@R=4mm@C=1mm{&\\\to\\&\\&}
\xymatrix@R=3mm@C=3mm{\circ\ar@/_/@{.}[r]&\circ&\circ\ar@{-}[dddll]&\circ\ar@{-}[dddl]\\
&&&\\
&&\ar@/^/@{.}[dr]&\\
\circ&\circ\ar@/^/@{.}[ur]&\circ&\circ}
\xymatrix@R=4mm@C=1mm{&\\\to\\&\\& }
\xymatrix@R=3mm@C=3mm{\circ\ar@/_/@{.}[r]&\circ&\circ\ar@{-}[dddll]&\circ\ar@{-}[dddll]\\
&&&\\
&&&\\
\circ&\circ&\circ\ar@/^/@{.}[r]&\circ}$$

The point now is that, under the assumption $\tau\in NC_{even}(k,l)$, each of the moves required for putting a leg at left, and hence for putting a whole block at left, requires an even number of switches. Thus, putting $\tau$ is standard form requires an even number of switches. Now given $\tau,\tau'\in NC_{even}$ having the same block structure, the standard form coincides, so the number of switches $c$ required for the passage $\tau\to\tau'$ is indeed even.

\medskip

Regarding now the remaining assertions, these are all elementary:

\medskip

(1) For $\tau\in S_k$ the standard form is $\tau'=id$, and the passage $\tau\to id$ comes by composing with a number of transpositions, which gives the signature. 

\medskip

(2) For a general $\tau\in P_2$, the standard form is of type $\tau'=|\ldots|^{\cup\ldots\cup}_{\cap\ldots\cap}$, and the passage $\tau\to\tau'$ requires $c$ mod 2 switches, where $c$ is the number of crossings. 

\medskip

(3) Assuming that $\tau\in P_{even}$ comes from $\pi\in NC_{even}$ by merging a certain number of blocks, we can prove that the signature is 1 by proceeding by recurrence.
\end{proof}

With the above result in hand, we can now formulate:

\index{twisted linear map}
\index{twisted Kronecker symbol}

\begin{definition}
Associated to any partition $\pi\in P_{even}(k,l)$ is the linear map
$$\bar{T}_\pi:(\mathbb C^N)^{\otimes k}\to(\mathbb C^N)^{\otimes l}$$
given by the following formula, with $e_1,\ldots,e_N$ being the standard basis of $\mathbb C^N$,
$$\bar{T}_\pi(e_{i_1}\otimes\ldots\otimes e_{i_k})=\sum_{j_1\ldots j_l}\bar{\delta}_\pi\begin{pmatrix}i_1&\ldots&i_k\\ j_1&\ldots&j_l\end{pmatrix}e_{j_1}\otimes\ldots\otimes e_{j_l}$$
and where $\bar{\delta}_\pi\in\{-1,0,1\}$ is $\bar{\delta}_\pi=\varepsilon(\tau)$ if $\tau\geq\pi$, and $\bar{\delta}_\pi=0$ otherwise, with:
$$\tau=\ker\binom{i}{j}$$
\end{definition}

In other words, what we are doing here is to add signatures to the usual formula of $T_\pi$. Indeed, observe that the usual formula for $T_\pi$ can be written as folllows:
$$T_\pi(e_{i_1}\otimes\ldots\otimes e_{i_k})=\sum_{j:\ker(^i_j)\geq\pi}e_{j_1}\otimes\ldots\otimes e_{j_l}$$

Now by inserting signs, coming from the signature map $\varepsilon:P_{even}\to\{\pm1\}$, we are led to the following formula, which coincides with the one given above:
$$\bar{T}_\pi(e_{i_1}\otimes\ldots\otimes e_{i_k})=\sum_{\tau\geq\pi}\varepsilon(\tau)\sum_{j:\ker(^i_j)=\tau}e_{j_1}\otimes\ldots\otimes e_{j_l}$$

We will be back later to this analogy, with more details on what can be done with it. For the moment, we must first prove a key categorical result, as follows:

\begin{proposition}
The assignement $\pi\to\bar{T}_\pi$ is categorical, in the sense that
$$\bar{T}_\pi\otimes\bar{T}_\sigma=\bar{T}_{[\pi\sigma]}\quad,\quad 
\bar{T}_\pi \bar{T}_\sigma=N^{c(\pi,\sigma)}\bar{T}_{[^\sigma_\pi]}\quad,\quad
\bar{T}_\pi^*=\bar{T}_{\pi^*}$$
where $c(\pi,\sigma)$ are certain positive integers.
\end{proposition}

\begin{proof}
We have to go back to the proof from the untwisted case, from chapter 2, and insert signs. We have to check three conditions, as follows:

\medskip

\underline{1. Concatenation}. In the untwisted case, this was based on the following formula:
$$\delta_\pi\begin{pmatrix}i_1\ldots i_p\\ j_1\ldots j_q\end{pmatrix}
\delta_\sigma\begin{pmatrix}k_1\ldots k_r\\ l_1\ldots l_s\end{pmatrix}
=\delta_{[\pi\sigma]}\begin{pmatrix}i_1\ldots i_p&k_1\ldots k_r\\ j_1\ldots j_q&l_1\ldots l_s\end{pmatrix}$$

In the twisted case, it is enough to check the following formula:
$$\varepsilon\left(\ker\begin{pmatrix}i_1\ldots i_p\\ j_1\ldots j_q\end{pmatrix}\right)
\varepsilon\left(\ker\begin{pmatrix}k_1\ldots k_r\\ l_1\ldots l_s\end{pmatrix}\right)=
\varepsilon\left(\ker\begin{pmatrix}i_1\ldots i_p&k_1\ldots k_r\\ j_1\ldots j_q&l_1\ldots l_s\end{pmatrix}\right)$$

Let us denote by $\tau,\nu$ the partitions on the left, so that the partition on the right is of the form $\rho\leq[\tau\nu]$. Now by switching to the noncrossing form, $\tau\to\tau'$ and $\nu\to\nu'$, the partition on the right transforms into $\rho\to\rho'\leq[\tau'\nu']$. Now since the partition $[\tau'\nu']$ is noncrossing, we can use Theorem 11.7 (3), and we obtain the result.

\medskip

\underline{2. Composition}. In the untwisted case, this was based on the following formula:
$$\sum_{j_1\ldots j_q}\delta_\pi\begin{pmatrix}i_1\ldots i_p\\ j_1\ldots j_q\end{pmatrix}
\delta_\sigma\begin{pmatrix}j_1\ldots j_q\\ k_1\ldots k_r\end{pmatrix}
=N^{c(\pi,\sigma)}\delta_{[^\pi_\sigma]}\begin{pmatrix}i_1\ldots i_p\\ k_1\ldots k_r\end{pmatrix}$$

In order to prove now the result in the twisted case, it is enough to check that the signs match. More precisely, we must establish the following formula:
$$\varepsilon\left(\ker\begin{pmatrix}i_1\ldots i_p\\ j_1\ldots j_q\end{pmatrix}\right)
\varepsilon\left(\ker\begin{pmatrix}j_1\ldots j_q\\ k_1\ldots k_r\end{pmatrix}\right)
=\varepsilon\left(\ker\begin{pmatrix}i_1\ldots i_p\\ k_1\ldots k_r\end{pmatrix}\right)$$

Let $\tau,\nu$ be the partitions on the left, so that the partition on the right is of the form $\rho\leq[^\tau_\nu]$. Our claim is that we can jointly switch $\tau,\nu$ to the noncrossing form. Indeed, we can first switch as for $\ker(j_1\ldots j_q)$ to become noncrossing, and then switch the upper legs of $\tau$, and the lower legs of $\nu$, as for both these partitions to become noncrossing. Now observe that when switching in this way to the noncrossing form, $\tau\to\tau'$ and $\nu\to\nu'$, the partition on the right transforms into $\rho\to\rho'\leq[^{\tau'}_{\nu'}]$. Now since the partition $[^{\tau'}_{\nu'}]$ is noncrossing, we can apply Theorem 11.7 (3), and we obtain the result.

\medskip

\underline{3. Involution}. Here we must prove the following formula:
$$\bar{\delta}_\pi\begin{pmatrix}i_1\ldots i_p\\ j_1\ldots j_q\end{pmatrix}=\bar{\delta}_{\pi^*}\begin{pmatrix}j_1\ldots j_q\\ i_1\ldots i_p\end{pmatrix}$$

But this is clear from the definition of $\bar{\delta}_\pi$, and we are done.
\end{proof}

As a conclusion, our twisted construction $\pi\to\bar{T}_\pi$ has all the needed properties for producing quantum groups, via Tannakian duality, and we can now formulate:

\begin{theorem}
Given a category of partitions $D\subset P_{even}$, the construction
$$Hom(u^{\otimes k},u^{\otimes l})=span\left(\bar{T}_\pi\Big|\pi\in D(k,l)\right)$$
produces via Tannakian duality a quantum group $\bar{G}_N\subset U_N^+$, for any $N\in\mathbb N$.
\end{theorem}

\begin{proof}
This follows indeed from the Tannakian results from chapter 2, exactly as in the easy case, by using this time Proposition 11.9 as technical ingredient. To be more precise, Proposition 11.9 shows that the linear spaces on the right form a Tannakian category, and so the results in chapter 2 apply, and give the result.
\end{proof}

We can unify the easy quantum groups, or at least the examples coming from categories $D\subset P_{even}$, with the quantum groups constructed above, as follows:

\index{quizzy quantum group}
\index{q-easy quantum group}
\index{twisting}
\index{Schur-Weyl twisting}

\begin{definition}
A closed subgroup $G\subset U_N^+$ is called $q$-easy, or quizzy, with deformation parameter $q=\pm1$, when its tensor category appears as follows,
$$Hom(u^{\otimes k},u^{\otimes l})=span\left(\dot{T}_\pi\Big|\pi\in D(k,l)\right)$$
for a certain category of partitions $D\subset P_{even}$, where, for $q=-1,1$:
$$\dot{T}=\bar{T},T$$
The Schur-Weyl twist of $G$ is the quizzy quantum group $\bar{G}\subset U_N^+$ obtained via $q\to-q$.
\end{definition}

We will see later on that the easy quantum group associated to $P_{even}$ itself is the hyperochahedral group $H_N$, and so that our assumption $D\subset P_{even}$, replacing $D\subset P$, simply corresponds to $H_N\subset G$, replacing the usual condition $S_N\subset G$.

\bigskip

For the moment, our most pressing task is that of checking that, when applying the Schur-Weyl twisting to the basic unitary quantum groups, we obtain the ad-hoc twists that we previously constructed. This is indeed the case:

\index{twisted orthogonal group}
\index{twisted unitary group}

\begin{theorem}
The twisted unitary quantum groups introduced before,
$$\xymatrix@R=12mm@C=12mm{
O_N^+\ar[r]&\mathbb TO_N^+\ar[r]&U_N^+\\
\bar{O}_N^*\ar[r]\ar[u]&\mathbb T\bar{O}_N^*\ar[r]\ar[u]&\bar{U}_N^*\ar[u]\\
\bar{O}_N\ar[r]\ar[u]&\mathbb T\bar{O}_N\ar[r]\ar[u]&\bar{U}_N\ar[u]}$$
appear as Schur-Weyl twists of the basic unitary quantum groups.
\end{theorem}

\begin{proof}
This is something routine, in several steps, as follows:

\medskip

(1) The basic crossing, $\ker\binom{ij}{ji}$ with $i\neq j$, comes from the transposition $\tau\in S_2$, so its signature is $-1$. As for its degenerated version $\ker\binom{ii}{ii}$, this is noncrossing, so here the signature is $1$. We conclude that the linear map associated to the basic crossing is:
$$\bar{T}_{\slash\!\!\!\backslash}(e_i\otimes e_j)
=\begin{cases}
-e_j\otimes e_i&{\rm for}\ i\neq j\\
e_j\otimes e_i&{\rm otherwise}
\end{cases}$$

For the half-classical crossing, namely $\ker\binom{ijk}{kji}$ with $i,j,k$ distinct, the signature is once again $-1$, and by examining the signatures of the various degenerations of this half-classical crossing, we are led to the following formula:
$$\bar{T}_{\slash\hskip-1.6mm\backslash\hskip-1.1mm|\hskip0.5mm}(e_i\otimes e_j\otimes e_k)
=\begin{cases}
-e_k\otimes e_j\otimes e_i&{\rm for}\ i,j,k\ {\rm distinct}\\
e_k\otimes e_j\otimes e_i&{\rm otherwise}
\end{cases}$$

(2) Our claim now if that for an orthogonal quantum group $G$, the following holds, with the quantum group $\bar{O}_N$ being the one in Theorem 11.2:
$$\bar{T}_{\slash\!\!\!\backslash}\in End(u^{\otimes 2})\iff G\subset\bar{O}_N$$

Indeed, by using the formula of $\bar{T}_{\slash\!\!\!\backslash}$ found in (1) above, we obtain:
\begin{eqnarray*}
(\bar{T}_{\slash\!\!\!\backslash}\otimes1)u^{\otimes 2}(e_i\otimes e_j\otimes1)
&=&\sum_ke_k\otimes e_k\otimes u_{ki}u_{kj}\\
&-&\sum_{k\neq l}e_l\otimes e_k\otimes u_{ki}u_{lj}
\end{eqnarray*}

On the other hand, we have as well the following formula:
\begin{eqnarray*}
u^{\otimes 2}(\bar{T}_{\slash\!\!\!\backslash}\otimes1)(e_i\otimes e_j\otimes1)
&=&\begin{cases}
\sum_{kl}e_l\otimes e_k\otimes u_{li}u_{ki}&{\rm if}\ i=j\\
-\sum_{kl}e_l\otimes e_k\otimes u_{lj}u_{ki}&{\rm if}\ i\neq j
\end{cases}
\end{eqnarray*}

For $i=j$ the conditions are $u_{ki}^2=u_{ki}^2$ for any $k$, and $u_{ki}u_{li}=-u_{li}u_{ki}$ for any $k\neq l$. For $i\neq j$ the conditions are $u_{ki}u_{kj}=-u_{kj}u_{ki}$ for any $k$, and $u_{ki}u_{lj}=u_{lj}u_{ki}$ for any $k\neq l$. Thus we have exactly the relations between the coordinates of $\bar{O}_N$, and we are done.

\medskip

(3) Our claim now if that for an orthogonal quantum group $G$, the following holds, with the quantum group $\bar{O}_N^*$ being the one in Theorem 11.3:
$$\bar{T}_{\slash\hskip-1.6mm\backslash\hskip-1.1mm|\hskip0.5mm}\in End(u^{\otimes 3})\iff G\subset\bar{O}_N^*$$

Indeed, by using the formula of $\bar{T}_{\slash\hskip-1.6mm\backslash\hskip-1.1mm|\hskip0.5mm}$ found in (1) above, we obtain:
\begin{eqnarray*}
(\bar{T}_{\slash\hskip-1.6mm\backslash\hskip-1.1mm|\hskip0.5mm}\otimes1)u^{\otimes 2}(e_i\otimes e_j\otimes e_k\otimes1)
&=&\sum_{abc\ not\ distinct}e_c\otimes e_b\otimes e_a\otimes u_{ai}u_{bj}u_{ck}\\
&-&\sum_{a,b,c\ distinct}e_c\otimes e_b\otimes e_a\otimes u_{ai}u_{bj}u_{ck}
\end{eqnarray*}

On the other hand, we have as well the following formula:
\begin{eqnarray*}
&&u^{\otimes 2}(\bar{T}_{\slash\hskip-1.6mm\backslash\hskip-1.1mm|\hskip0.5mm}\otimes1)(e_i\otimes e_j\otimes e_k\otimes1)\\
&&=\begin{cases}
\sum_{abc}e_c\otimes e_b\otimes e_a\otimes u_{ck}u_{bj}u_{ai}&{\rm for}\ i,j,k\ {\rm not\ distinct}\\
-\sum_{abc}e_c\otimes e_b\otimes e_a\otimes u_{ck}u_{bj}u_{ai}&{\rm for}\ i,j,k\ {\rm distinct}
\end{cases}
\end{eqnarray*}

For $i,j,k$ not distinct the conditions are $u_{ai}u_{bj}u_{ck}=u_{ck}u_{bj}u_{ai}$ for $a,b,c$ not distinct, and $u_{ai}u_{bj}u_{ck}=-u_{ck}u_{bj}u_{ai}$ for $a,b,c$ distinct. For $i,j,k$ distinct the conditions are $u_{ai}u_{bj}u_{ck}=-u_{ck}u_{bj}u_{ai}$ for $a,b,c$ not distinct, and $u_{ai}u_{bj}u_{ck}=u_{ck}u_{bj}u_{ai}$ for $a,b,c$ distinct. Thus we have the relations between the coordinates of $\bar{O}_N^*$, as desired.

\medskip

(4) Now with the above in hand, we obtain that the Schur-Weyl twists of $O_N,O_N^*$ are indeed the quantum groups $\bar{O}_N,\bar{O}_N^*$ from Theorem 11.2 and Theorem 11.3. 

\medskip

(5) The proof in the unitary case is similar, by adding signs in the above computations (2,3), the conclusion being that the Schur-Weyl twists of $U_N,U_N^*$ are $\bar{U}_N,\bar{U}_N^*$. 
\end{proof}

Let us clarify now the relation between the maps $T_\pi,\bar{T}_\pi$. We recall that the M\"obius function of any lattice, and in particular of $P_{even}$, is given by:
$$\mu(\sigma,\pi)=\begin{cases}
1&{\rm if}\ \sigma=\pi\\
-\sum_{\sigma\leq\tau<\pi}\mu(\sigma,\tau)&{\rm if}\ \sigma<\pi\\
0&{\rm if}\ \sigma\not\leq\pi
\end{cases}$$

With this notation, we have the following result:

\index{M\"obius inversion}
\index{twisted linear map}

\begin{proposition}
For any partition $\pi\in P_{even}$ we have the formula
$$\bar{T}_\pi=\sum_{\tau\leq\pi}\alpha_\tau T_\tau$$
where $\alpha_\sigma=\sum_{\sigma\leq\tau\leq\pi}\varepsilon(\tau)\mu(\sigma,\tau)$, with $\mu$ being the M\"obius function of $P_{even}$.
\end{proposition}

\begin{proof}
The linear combinations $T=\sum_{\tau\leq\pi}\alpha_\tau T_\tau$ acts on tensors as follows:
\begin{eqnarray*}
T(e_{i_1}\otimes\ldots\otimes e_{i_k})
&=&\sum_{\tau\leq\pi}\alpha_\tau T_\tau(e_{i_1}\otimes\ldots\otimes e_{i_k})\\
&=&\sum_{\tau\leq\pi}\alpha_\tau\sum_{\sigma\leq\tau}\sum_{j:\ker(^i_j)=\sigma}e_{j_1}\otimes\ldots\otimes e_{j_l}\\
&=&\sum_{\sigma\leq\pi}\left(\sum_{\sigma\leq\tau\leq\pi}\alpha_\tau\right)\sum_{j:\ker(^i_j)=\sigma}e_{j_1}\otimes\ldots\otimes e_{j_l}
\end{eqnarray*}

Thus, in order to have $\bar{T}_\pi=\sum_{\tau\leq\pi}\alpha_\tau T_\tau$, we must have $\varepsilon(\sigma)=\sum_{\sigma\leq\tau\leq\pi}\alpha_\tau$, for any $\sigma\leq\pi$. But this problem can be solved by using the M\"obius inversion formula, and we obtain the numbers $\alpha_\sigma=\sum_{\sigma\leq\tau\leq\pi}\varepsilon(\tau)\mu(\sigma,\tau)$ in the statement.
\end{proof}

With the above results in hand, let us go back now to the question of twisting the quantum reflection groups. It is convenient to include in our discussion two more quantum groups, coming from \cite{rwe} and denoted $H_N^{[\infty]},K_N^{[\infty]}$, constructed as follows:

\begin{proposition}
We have quantum groups $H_N^{[\infty]},K_N^{[\infty]}$ as follows, constructed by using the relations $\alpha\beta\gamma=0$ for any $a\neq c$ on the same row or column of $u$:
$$\xymatrix@R=15mm@C=17mm{
K_N\ar[r]&K_N^*\ar[r]&K_N^{[\infty]}\ar[r]&K_N^+\\
H_N\ar[r]\ar[u]&H_N^*\ar[r]\ar[u]&H_N^{[\infty]}\ar[r]\ar[u]&H_N^+\ar[u]}$$
These quantum groups are both easy, with the corresponding categories of partitions, denoted $P_{even}^{[\infty]}\subset P_{even}$ and $\mathcal P_{even}^{[\infty]}\subset\mathcal P_{even}$, being generated by $\eta=\ker(^{iij}_{jii})$.
\end{proposition}

\begin{proof}
This is routine, by using the fact that the relations $\alpha\beta\gamma=0$ in the statement are equivalent to the condition $\eta\in End(u^{\otimes k})$, with $|k|=3$. For details here, and for more on these two quantum groups, which are very interesting objects, and that we have actually already met in chapter 4 above, we refer to the paper of Raum-Weber \cite{rwe}.
\end{proof}

In order to discuss now the Schur-Weyl twisting of the various quantum reflection groups that we have, we will need the following technical result:

\begin{proposition}
We have the following equalities,
\begin{eqnarray*}
P_{even}^*&=&\left\{\pi\in P_{even}\Big|\varepsilon(\tau)=1,\forall\tau\leq\pi,|\tau|=2\right\}
\\
P_{even}^{[\infty]}&=&\left\{\pi\in P_{even}\Big|\sigma\in P_{even}^*,\forall\sigma\subset\pi\right\}\\
P_{even}^{[\infty]}&=&\left\{\pi\in P_{even}\Big|\varepsilon(\tau)=1,\forall\tau\leq\pi\right\}
\end{eqnarray*}
where $\varepsilon:P_{even}\to\{\pm1\}$ is the signature of even permutations.
\end{proposition}

\begin{proof}
This is routine combinatorics, from \cite{ba3}, \cite{rwe}, the idea being as follows:

\medskip

(1) Given $\pi\in P_{even}$, we have $\tau\leq\pi,|\tau|=2$ precisely when $\tau=\pi^\beta$ is the partition obtained from $\pi$ by merging all the legs of a certain subpartition $\beta\subset\pi$, and by merging as well all the other blocks. Now observe that $\pi^\beta$ does not depend on $\pi$, but only on $\beta$, and that the number of switches required for making $\pi^\beta$ noncrossing is $c=N_\bullet-N_\circ$ modulo 2, where $N_\bullet/N_\circ$ is the number of black/white legs of $\beta$, when labelling the legs of $\pi$ counterclockwise $\circ\bullet\circ\bullet\ldots$ Thus $\varepsilon(\pi^\beta)=1$ holds precisely when $\beta\in\pi$ has the same number of black and white legs, and this gives the result.

\medskip

(2) This simply follows from the equality $P_{even}^{[\infty]}=<\eta>$ coming from Proposition 11.14, by computing $<\eta>$, and for the complete proof here we refer to \cite{rwe}.

\medskip

(3) We use the fact, also from \cite{rwe}, that the relations $g_ig_ig_j=g_jg_ig_i$ are trivially satisfied for real reflections. Thus, we have:
$$P_{even}^{[\infty]}(k,l)=\left\{\ker\begin{pmatrix}i_1&\ldots&i_k\\ j_1&\ldots&j_l\end{pmatrix}\Big|g_{i_1}\ldots g_{i_k}=g_{j_1}\ldots g_{j_l}\ {\rm inside}\ \mathbb Z_2^{*N}\right\}$$

In other words, the partitions in $P_{even}^{[\infty]}$ are those describing the relations between free variables, subject to the conditions $g_i^2=1$. We conclude that $P_{even}^{[\infty]}$ appears from $NC_{even}$ by ``inflating blocks'', in the sense that each $\pi\in P_{even}^{[\infty]}$ can be transformed into a partition $\pi'\in NC_{even}$ by deleting pairs of consecutive legs, belonging to the same block. Now since this inflation operation leaves invariant modulo 2 the number $c\in\mathbb N$ of switches in the definition of the signature, it leaves invariant the signature $\varepsilon=(-1)^c$ itself, and we obtain in this way the inclusion ``$\subset$'' in the statement. 

\medskip

Conversely, given $\pi\in P_{even}$ satisfying $\varepsilon(\tau)=1$, $\forall\tau\leq\pi$, our claim is that:
$$\rho\leq\sigma\subset\pi,|\rho|=2\implies\varepsilon(\rho)=1$$

Indeed, let us denote by $\alpha,\beta$ the two blocks of $\rho$, and by $\gamma$ the remaining blocks of $\pi$, merged altogether. We know that the partitions $\tau_1=(\alpha\wedge\gamma,\beta)$, $\tau_2=(\beta\wedge\gamma,\alpha)$, $\tau_3=(\alpha,\beta,\gamma)$ are all even. On the other hand, putting these partitions in noncrossing form requires respectively $s+t,s'+t,s+s'+t$ switches, where $t$ is the number of switches needed for putting $\rho=(\alpha,\beta)$ in noncrossing form. Thus $t$ is even, and we are done. With the above claim in hand, we conclude, by using the second equality in the statement, that we have $\sigma\in P_{even}^*$. Thus we have $\pi\in P_{even}^{[\infty]}$, which ends the proof of ``$\supset$''.
\end{proof}

With the above result in hand, we can now prove:

\begin{theorem}
The basic quantum reflection groups, namely
$$\xymatrix@R=13mm@C=13mm{
H_N^+\ar[r]&\mathbb TH_N^+\ar[r]&K_N^+\\
H_N^*\ar[r]\ar[u]&\mathbb TH_N^*\ar[r]\ar[u]&K_N^*\ar[u]\\
H_N\ar[r]\ar[u]&\mathbb TH_N\ar[r]\ar[u]&K_N\ar[u]}$$
equal their own Schur-Weyl twists.
\end{theorem}

\begin{proof}
This result, from \cite{ba3}, basically comes from the results that we have:

\medskip

(1) In the real case, the verifications are as follows:

\medskip

-- $H_N^+$. We know from Theorem 11.7 that for $\pi\in NC_{even}$ we have $\bar{T}_\pi=T_\pi$, and since we are in the situation $D\subset NC_{even}$, the definitions of $G,\bar{G}$ coincide.

\medskip

-- $H_N^{[\infty]}$. Here we can use the same argument as in (1), based this time on the description of $P_{even}^{[\infty]}$ involving the signatures found in Proposition 11.15.

\medskip

-- $H_N^*$. We have $H_N^*=H_N^{[\infty]}\cap O_N^*$, so $\bar{H}_N^*\subset H_N^{[\infty]}$ is the subgroup obtained via the defining relations for $\bar{O}_N^*$. But all the $abc=-cba$ relations defining $\bar{H}_N^*$ are automatic, of type $0=0$, and it follows that $\bar{H}_N^*\subset H_N^{[\infty]}$ is the subgroup obtained via the relations $abc=cba$, for any $a,b,c\in\{u_{ij}\}$. Thus we have $\bar{H}_N^*=H_N^{[\infty]}\cap O_N^*=H_N^*$, as claimed.

\medskip

-- $H_N$. We have $H_N=H_N^*\cap O_N$, and by functoriality, $\bar{H}_N=\bar{H}_N^*\cap\bar{O}_N=H_N^*\cap\bar{O}_N$. But this latter intersection is easily seen to be equal to $H_N$, as claimed.

\medskip

(2) In the complex case the proof is similar, and we refer here to \cite{ba3}.
\end{proof}

In relation now with the tori, we have the following result:

\begin{theorem}
The diagonal tori of the twisted quantum groups are
$$\xymatrix@R=13.5mm@C=13.5mm{
T_N^+\ar[r]&\mathbb TT_N^+\ar[r]&\mathbb T_N^+\\
T_N^*\ar[r]\ar[u]&\mathbb TT_N^*\ar[r]\ar[u]&\mathbb T_N^*\ar[u]\\
T_N\ar[r]\ar[u]&\mathbb TT_N\ar[r]\ar[u]&\mathbb T_N\ar[u]}$$
exactly as in the untwisted case.
\end{theorem}

\begin{proof}
This is clear for the quantum reflection groups, which are not twistable, and for the quantum unitary groups this is elementary as well, coming from definitions.
\end{proof}

\section*{11c. Twisted integration}

Before getting into the spheres, let us discuss integration questions. The result here, valid for any Schur-Weyl twist in our sense, is as follows:

\index{twisted integration}
\index{twisted Weingarten formula}

\begin{theorem}
We have the Weingarten type formula
$$\int_{\bar{G}}u_{i_1j_1}^{e_1}\ldots u_{i_kj_k}^{e_k}=\sum_{\pi,\sigma\in D(k)}\bar{\delta}_\pi(i_1\ldots i_k)\bar{\delta}_\sigma(j_1\ldots j_k)W_{kN}(\pi,\sigma)$$
where $W_{kN}=G_{kN}^{-1}$, with $G_{kN}(\pi,\sigma)=N^{|\pi\vee\sigma|}$, for $\pi,\sigma\in D(k)$.
\end{theorem}

\begin{proof}
This follows exactly as in the untwisted case, the idea being that the signs will cancel. Let us recall indeed from Definition 11.8 and the comments afterwards that the twisted vectors $\bar{\xi}_\pi$ associated to the partitions $\pi\in P_{even}(k)$ are as follows: 
$$\bar{\xi}_\pi=\sum_{\tau\geq\pi}\varepsilon(\tau)\sum_{i:\ker(i)=\tau}e_{i_1}\otimes\ldots\otimes e_{i_k}$$

Thus, the Gram matrix of these vectors is given by:
\begin{eqnarray*}
<\xi_\pi,\xi_\sigma>
&=&\sum_{\tau\geq\pi\vee\sigma}\varepsilon(\tau)^2\left|\left\{(i_1,\ldots,i_k)\Big|\ker i=\tau\right\}\right|\\
&=&\sum_{\tau\geq\pi\vee\sigma}\left|\left\{(i_1,\ldots,i_k)\Big|\ker i=\tau\right\}\right|\\
&=&N^{|\pi\vee\sigma|}
\end{eqnarray*}

Thus the Gram matrix is the same as in the untwisted case, and so the Weingarten matrix is the same as well as in the untwisted case, and this gives the result.
\end{proof}

In relation now with the spheres, we have the following result:

\begin{theorem}
The twisted spheres have the following properties:
\begin{enumerate}
\item They have affine actions of the twisted unitary quantum groups.

\item They have unique invariant Haar functionals, which are ergodic.

\item Their Haar functionals are given by Weingarten type formulae.

\item They appear, via the GNS construction, as first row spaces.
\end{enumerate}
\end{theorem}

\begin{proof}
The proofs here are similar to those from the untwisted case, via a routine computation, by adding signs where needed, and with the main technical ingredient, namely the Weingarten formula, being available from Theorem 11.18. See \cite{ba3}.
\end{proof}

As a conclusion now, we have shown that the various quadruplets $(S,T,U,K)$ constructed in chapters 1-10 have twisted counterparts $(\bar{S},T,\bar{U},K)$. The question that we would like to solve now is that of finding correspondences, as follows:
$$
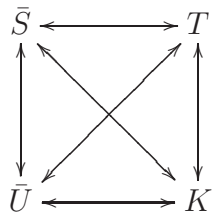
\xymatrix@R=50pt@C=50pt{
\bar{S}\ar[r]\ar[d]\ar[dr]&T\ar[l]\ar[d]\ar[dl]\\
\bar{U}\ar[u]\ar[ur]\ar[r]&K\ar[l]\ar[ul]\ar[u]
}$$

In order to discuss this, let us get back to the axioms from chapter 4. We have seen there that the 12 correspondences come in fact from 7  correspondences, as follows:
$$
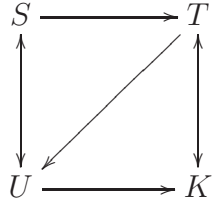
\xymatrix@R=50pt@C=50pt{
S\ar[r]\ar[d]&T\ar[d]\ar[dl]\\
U\ar[u]\ar[r]&K\ar[u]
}$$

In the twisted case, 6 of these correspondences hold as well, but the remaining one, namely $S\to T$, definitely does not hold as stated, and must be modified. Let us begin our discussion with the quantum isometry group results. We have here:

\index{twisted isometry groups}

\begin{theorem}
We have the quantum isometry group formula 
$$\bar{U}=G^+(\bar{S})$$
in all the $9$ main twisted cases.
\end{theorem}

\begin{proof}
The proofs here are similar to those from the untwisted case, via a routine computation, by adding signs where needed, which amounts in replacing the usual commutators $[a,b]=ab-ba$ by twisted commutators, given by:
$$[[a,b]]=ab+ba$$

There is one subtle point, however, coming from the fact that the linear independence of various products of coordinates of length 1,2,3, which was something clear in the untwisted case, is now a non-trivial question. But this can be solved via a technical application of the Weingarten formula, from Theorem 11.18. See \cite{ba1}.
\end{proof}

Regarding now the $K=G^+(T)\cap K_N^+$ axiom, this is something that we already know. However, regarding the correspondence $S\to T$, things here fail in the twisted case. Our ``fix'' for this, or at least the best fix that we could find, is as follows:

\begin{theorem}
Given an algebraic manifold $X\subset S^{N-1}_{\mathbb C,+}$, define its toral isometry group as being the biggest subgroup of $\mathbb T_N^+$ acting affinely on $X$:
$$\mathcal G^+(X)=G^+(X)\cap\mathbb T_N^+$$
With this convention, for the $9$ basic spheres $S$, and for their twists as well, the toral isometry group equals the torus $T$.
\end{theorem}

\begin{proof}
We recall from chapter 3 that the affine quantum isometry group $G^+(X)\subset U_N^+$ of a noncommutative manifold $X\subset S^{N-1}_{\mathbb C,+}$ coming from certain polynomial relations $P$ is constructed according to the following procedure:
$$P(x_i)=0\implies P\left(\sum_jx_j\otimes u_{ji}\right)=0$$

Similarly, the toral isometry group $\mathcal G^+(X)\subset\mathbb T_N^+$ is constructed as follows:
$$P(x_i)=0\implies P\left(x_i\otimes u_i\right)=0$$

In the easy case one can prove that the following formula holds:
$$G^+(\bar{S})=\overline{G^+(S)}$$

By intersecting with $\mathbb T_N^+$, we obtain from this that we have:
$$\mathcal G^+(\bar{S})=\mathcal G^+(S)$$

The result can be of course be proved as well directly. For $\bar{S}^{N-1}_\mathbb R$ we have:
$$\Phi(x_ix_j)=x_ix_j\otimes u_iu_j$$
$$\Phi(x_jx_i)=x_jx_i\otimes u_ju_i$$

Thus we obtain $u_iu_j=-u_ju_i$ for $i\neq j$, and so the quantum group is $T_N$. The proof in the complex, half-liberated and hybrid cases is similar.
\end{proof}

Regarding the hard liberation axiom, this seems to hold indeed in all the cases under consideration, but this is non-trivial, and not known yet. As a conclusion, we conjecturally have an extension of our $(S,T,U,K)$ formalism, with the $S\to T$ axiom needing a modification as above, which covers the twisted objects $(\bar{S},T,\bar{U},K)$ as well.

\section*{11d. Twisted geometry}

There are many things that can be done in the context of the twisted geometry, going beyond what we have so far, namely some theory and computations for the spheres $\bar{S}$, tori $T$, unitary groups $\bar{U}$, and reflection groups $K$. We briefly discuss here, as a main topic, the twisted extension of the various constructions from chapter 6.

\bigskip

So, let us go back to the theory there. As a first observation, we can both liberate the spaces $O_{MN}^L,U_{MN}^L$, and twist them, by proceeding as as follows:

\index{twisted partial isometries}

\begin{definition}
Associated to any integers $L\leq M\leq N$ are the algebras
\begin{eqnarray*}
C(O_{MN}^{L+})&=&C^*\left((u_{ij})_{i=1,\ldots,M,j=1,\ldots,N}\Big|u=\bar{u},uu^t={\rm projection\ of\ trace}\ L\right)\\
C(U_{MN}^{L+})&=&C^*\left((u_{ij})_{i=1,\ldots,M,j=1,\ldots,N}\Big|uu^*,\bar{u}u^t={\rm projections\ of\ trace}\ L\right)
\end{eqnarray*}
and their quotients $C(\bar{O}_{MN}^L),C(\bar{U}_{MN}^L)$, obtained by imposing the twisting relations.
\end{definition}

With this extended formalism, we have inclusions between the various spaces constructed so far, in chapter 6 and then here, as follows:
$$\xymatrix@R=15mm@C=15mm{
U_{MN}^L\ar[r]&U_{MN}^{L+}&\bar{U}_{MN}^L\ar[l]\\
O_{MN}^L\ar[r]\ar[u]&O_{MN}^{L+}\ar[u]&\bar{O}_{MN}^L\ar[l]\ar[u]}$$

More generally, we can perform these constructions for any quizzy quantum group. In order to discuss this, we use the Kronecker symbols $\delta_\pi(i)\in\{-1,0,1\}$, given by:
$$\delta_\sigma(i)=\begin{cases}
\delta_{\ker i\leq\sigma}&{\rm (untwisted\ case)}\\
\varepsilon(\ker i)\delta_{\ker i\leq\sigma}&{\rm (twisted\ case)}
\end{cases}$$

With this convention, we have the following result, from \cite{ba5}:

\begin{proposition}
The various spaces $G_{MN}^L$ constructed so far appear by imposing to the standard coordinates of $U_{MN}^{L+}$ the relations
$$\sum_{i_1\ldots i_s}\sum_{j_1\ldots j_s}\delta_\pi(i)\delta_\sigma(j)u_{i_1j_1}^{e_1}\ldots u_{i_sj_s}^{e_s}=L^{|\pi\vee\sigma|}$$
with $s=(e_1,\ldots,e_s)$ ranging over all the colored integers, and with $\pi,\sigma\in D(0,s)$.
\end{proposition}

\begin{proof}
The relations defining $G_{MN}^L$ are as follows, with $\sigma$ ranging over a family of generators, with no upper legs, of the corresponding category of partitions $D$:
$$\sum_{j_1\ldots j_s}\delta_\sigma(j)u_{i_1j_1}^{e_1}\ldots u_{i_sj_s}^{e_s}=\delta_\sigma(i)$$

We therefore obtain the relations in the statement, as follows:
\begin{eqnarray*}
\sum_{i_1\ldots i_s}\sum_{j_1\ldots j_s}\delta_\pi(i)\delta_\sigma(j)u_{i_1j_1}^{e_1}\ldots u_{i_sj_s}^{e_s}
&=&\sum_{i_1\ldots i_s}\delta_\pi(i)\sum_{j_1\ldots j_s}\delta_\sigma(j)u_{i_1j_1}^{e_1}\ldots u_{i_sj_s}^{e_s}\\
&=&\sum_{i_1\ldots i_s}\delta_\pi(i)\delta_\sigma(i)\\
&=&\sum_{\tau\leq\pi\vee\sigma}\sum_{\ker i=\tau}(\pm1)^2\\
&=&\sum_{\tau\leq\pi\vee\sigma}\sum_{\ker i=\tau}1\\
&=&L^{|\pi\vee\sigma|}
\end{eqnarray*}

As for the converse, this follows by using the relations in the statement, by keeping $\pi$ fixed, and by making $\sigma$ vary over all the partitions in the category.
\end{proof}

Thus, we have unified the twisted and untwisted constructions, in the continuous case. In the general case now, where $G=(G_N)$ is an arbitary uniform quizzy quantum group, we can construct spaces $G_{MN}^L$ by using the above relations, and we have:

\begin{theorem}
The spaces $G_{MN}^L\subset U_{MN}^{L+}$ constructed by imposing the relations 
$$\sum_{i_1\ldots i_s}\sum_{j_1\ldots j_s}\delta_\pi(i)\delta_\sigma(j)u_{i_1j_1}^{e_1}\ldots u_{i_sj_s}^{e_s}=L^{|\pi\vee\sigma|}$$
with $\pi,\sigma$ ranging over all the partitions in the associated category, having no upper legs, are subject to an action map/quotient map diagram, as follows,
$$\xymatrix@R=15mm@C=30mm{
G\times G\ar[r]^m\ar[d]_{p\times id}&G\ar[d]^p\\
X\times G\ar[r]^a&X
}$$
exactly as in the classical case, or the free case.
\end{theorem}

\begin{proof}
We proceed as in chapter 6. We must prove that, if the variables $u_{ij}$ satisfy the relations in the statement, then so do the following variables:
$$U_{ij}=\sum_{kl}a_{ik}\otimes b_{jl}^*\otimes u_{kl}$$
$$V_{ij}=\sum_{l\leq L}a_{il}\otimes b_{jl}^*$$

Regarding the variables $U_{ij}$, the computation here goes as follows:
\begin{eqnarray*}
&&\sum_{i_1\ldots i_s}\sum_{j_1\ldots j_s}\delta_\pi(i)\delta_\sigma(j)U_{i_1j_1}^{e_1}\ldots U_{i_sj_s}^{e_s}\\
&=&\sum_{i_1\ldots i_s}\sum_{j_1\ldots j_s}\sum_{k_1\ldots k_s}\sum_{l_1\ldots l_s}\delta_\pi(i)\delta_\sigma(j)a_{i_1k_1}^{e_1}\ldots a_{i_sk_s}^{e_s}\otimes(b_{j_sl_s}^{e_s}\ldots b_{j_1l_1}^{e_1})^*\otimes u_{k_1l_1}^{e_1}\ldots u_{k_sl_s}^{e_s}\\
&=&\sum_{k_1\ldots k_s}\sum_{l_1\ldots l_s}\delta_\pi(k)\delta_\sigma(l)u_{k_1l_1}^{e_1}\ldots u_{k_sl_s}^{e_s}\\
&=&L^{|\pi\vee\sigma|}
\end{eqnarray*}

For the variables $V_{ij}$ the proof is similar, as follows:
\begin{eqnarray*}
&&\sum_{i_1\ldots i_s}\sum_{j_1\ldots j_s}\delta_\pi(i)\delta_\sigma(j)V_{i_1j_1}^{e_1}\ldots V_{i_sj_s}^{e_s}\\
&=&\sum_{i_1\ldots i_s}\sum_{j_1\ldots j_s}\sum_{l_1,\ldots,l_s\leq L}\delta_\pi(i)\delta_\sigma(j)a_{i_1l_1}^{e_1}\ldots a_{i_sl_s}^{e_s}\otimes(b_{j_sl_s}^{e_s}\ldots b_{j_1l_1}^{e_1})^*\\
&=&\sum_{l_1,\ldots,l_s\leq L}\delta_\pi(l)\delta_\sigma(l)\\
&=&L^{|\pi\vee\sigma|}
\end{eqnarray*}

Thus we have constructed an action map, and a quotient map, and the commutation of the diagram in the statement is then trivial.
\end{proof}

Still by following the material in chapter 6, we can now construct a Haar integration for the above spaces, and we have the following result, also from \cite{ba5}:

\index{twisted Weingarten formula}

\begin{theorem}
We have the Weingarten type formula
$$\int_{G_{MN}^L}u_{i_1j_1}\ldots u_{i_sj_s}=\sum_{\pi\sigma\tau\nu}L^{|\sigma\vee\nu|}\delta_\pi(i)\delta_\tau(j)W_{sM}(\pi,\sigma)W_{sN}(\tau,\nu)$$
where $W_{sM}=G_{sM}^{-1}$, with $G_{sM}(\pi,\sigma)=M^{|\pi\vee\sigma|}$.
\end{theorem}

\begin{proof}
We make use of the usual quantum group Weingarten formula, explained in the above in the twisted case. By using this formula for $G_M,G_N$, we obtain:
\begin{eqnarray*}
\int_{G_{MN}^L}u_{i_1j_1}\ldots u_{i_sj_s}
&=&\sum_{l_1\ldots l_s\leq L}\int_{G_M}a_{i_1l_1}\ldots a_{i_sl_s}\int_{G_N}b_{j_1l_1}^*\ldots b_{j_sl_s}^*\\
&=&\sum_{l_1\ldots l_s\leq L}\sum_{\pi\sigma}\delta_\pi(i)\delta_\sigma(l)W_{sM}(\pi,\sigma)\sum_{\tau\nu}\delta_\tau(j)\delta_\nu(l)W_{sN}(\tau,\nu)\\
&=&\sum_{\pi\sigma\tau\nu}\left(\sum_{l_1\ldots l_s\leq L}\delta_\sigma(l)\delta_\nu(l)\right)\delta_\pi(i)\delta_\tau(j)W_{sM}(\pi,\sigma)W_{sN}(\tau,\nu)
\end{eqnarray*}

Let us compute now the coefficient appearing in the last formula. Since the signature map takes $\pm1$ values, for any multi-index $l=(l_1,\ldots,l_s)$ we have:
\begin{eqnarray*}
\delta_\sigma(l)\delta_\nu(l)
&=&\delta_{\ker l\leq\sigma}\varepsilon(\ker l)\cdot\delta_{\ker l\leq\nu}\varepsilon(\ker l)\\
&=&\delta_{\ker l\leq\sigma\vee\nu}
\end{eqnarray*}

Thus the coefficient is $L^{|\sigma\vee\nu|}$, and we obtain the formula in the statement.
\end{proof}

With this formula in hand, we can derive explicit integration results for the sums of non-overlapping coordinates, exactly as in chapter 6. To be more precise, the laws and their asymptotics are identical in the classical and twisted cases. See \cite{ba5}.

\section*{11e. Exercises} 

As already mentioned in the above, this chapter was just a modest introduction to the twisted geometry, and what we did here, namely a study of the quadruplets $(\bar{S},T,\bar{U},K)$, and of some related homogeneous spaces, is just an epsilon of what can be done. As an initiation to all the unexplored land which is left, we have:

\begin{exercise}
Make a list, based on the existing $q=-1$ literature, carefully checked and doublechecked, of the compact quantum groups $G\subset U_N^+$, not necessarily easy, which can be twisted in a reasonable sense, and develop some theory for them.
\end{exercise}

The key words here are ``careful'' and ``reasonable'', and this due to the fact that the standard Drinfeld-Jimbo twisting procedure, widely used in the literature, produces non-semisimple quantum groups at $q=-1$, and so is not useful. However, there are many interesting semisimple examples, that is, quantum groups $G\subset U_N^+$ in our sense, waiting to be discovered, such as twists $\bar{SU}_2,\bar{SO}_3$ of the much loved groups $SU_2,SO_3$, and with a mysterious isomorphism $S_4^+\simeq \bar{SO}_3$ waiting to be discovered as well.

\begin{exercise}
Develop some systematic geometric theory, based on K-theory, differential geometry, and other techniques of Connes, in the twisted setting.
\end{exercise}

This is another interesting exercise, the point being that all the techniques of Connes apply well to the twisted case, and by the way to the half-classical case too, and generally speaking, to everything which is not wild enough to the point of being free. As an example here, K-theory fails already for the free group algebras, or rather for the free tori, in our parlance, because it yields different groups, depending on whether the full or reduced algebra is considered, and so is not an invariant of the corresponding torus, viewed as a quantum space, but is rather some kind of functional analysis complication, of not much use. However, as said above, when getting away from freeness, and more specifically from non-amenability, K-theory works perfectly, and so do the other geometric techniques of Connes, with potentially very interesting results at stake.

\chapter{Matrix models}

\section*{12a. Matrix models}

You can model everything with random matrices, the saying in analysis goes. We have already seen an instance of this phenomenon in chapter 9, when talking about half-liberation. To be more precise, for certain manifolds $X\subset S^{N-1}_{\mathbb C,*}$, we have constructed embeddings of algebras of the following type, with $Y$ being a certain classical manifold, and $T_1,\ldots,T_N\in M_2(C(Y))$ being certain suitable antidiagonal $2\times 2$ matrices:
$$\pi:C(X)\subset M_2(C(Y))\quad,\quad 
x_i\to T_i$$

These models, which are quite powerful, were used afterwards in order to establish several non-trivial results on the original half-classical manifolds $X\subset S^{N-1}_{\mathbb C,*}$. Indeed, some knowledge and patience helping, any computation inside the target algebra $M_2(C(Y))$ can only be fun and doable, and produce results about $X\subset S^{N-1}_{\mathbb C,*}$ itself.

\bigskip

We discuss here, following \cite{bb1}, modeling questions for general manifolds $X\subset S^{N-1}_{\mathbb C,+}$, by using the same idea, suitably modified and generalized, as to cover most of the manifolds that we are interested in. Let us start with a broad definition, as follows:

\index{model for a manifold}

\begin{definition}
A model for a real algebraic manifold $X\subset S^{N-1}_{\mathbb C,+}$ is a morphism of $C^*$-algebras of the following type,
$$\pi:C(X)\to B$$
with $B$ being a $C^*$-algebra, called target of the model. We say that the model is faithful if $\pi$ is faithful, in the usual sense.
\end{definition}

Obviously, this is something too broad, because we can simply take $B=C(X)$, and we have in this way our faithful model, which is of course something unuseful:
$$id:C(X)\to C(X)$$

Thus, we must suitably restrict the class of target algebras $B$ that we use, to algebras that we ``know well''. However, this is something quite tricky, because if we want our model to be faithful, we cannot use simple algebras like the algebras $M_2(C(Y))$ used in the half-classical setting. In short, we are running into some difficulties here, of functional analytic nature, and a systematic discussion of all this is needed. 

\bigskip

As a first objective, let us try to understand if an arbitrary manifold $X\subset S^{N-1}_{\mathbb C,+}$ can be modeled by using familiar variables such as usual matrices, or operators. The answer here is yes, when using operators on a separable Hilbert space, with this coming from the GNS representation theorem, that we know from chapter 1, which is as follows:

\index{GNS theorem}
\index{GNS representation}

\begin{theorem}
Any $C^*$-algebra $A$ appears as closed $*$-algebra of operators on a Hilbert space, $A\subset B(H)$, in the following way:
\begin{enumerate}
\item In the commutative case, where $A=C(X)$, we can set $H=L^2(X)$, with respect to some probability measure on $X$, and use the embedding $g\to(g\to fg)$.

\item In general, we can set $H=L^2(A)$, with respect to some faithful positive trace $tr:A\to\mathbb C$, and then use a similar embedding, $a\to(b\to ab)$.
\end{enumerate}
\end{theorem}

\begin{proof}
This is something that we already know, from chapter 1, coming from basic measure theory and functional analysis, the idea being as follows:

\medskip

(1) In the commutative case, where $A=C(X)$ by the Gelfand theorem, we can pick a probability measure on $X$, and then we have an embedding as follows:
$$C(X)\subset B(L^2(X))\quad,\quad 
f\to(g\to fg)$$

(2) In general, assuming that a linear form $\varphi:A\to\mathbb C$ has suitable positivity properties, we can define a scalar product on $A$, by the following formula:
$$<a,b>=\varphi(b^*a)$$

By completing we obtain a Hilbert space $H$, and we have a representation as follows, called GNS representation of our algebra with respect to the linear form $\varphi$:
$$A\to B(H)\quad,\quad 
a\to(b\to ab)$$

Moreover, when $\varphi:A\to\mathbb C$ has suitable faithfulness properties, making it analogous to the integration functionals $\int_X:A\to\mathbb C$ from the commutative case, with respect to faithful probability measures on $X$, this representation is faithful, as desired. 
\end{proof}

Now back to our questions, the above result tells us that we have:

\index{algebraic manifold}
\index{faithful model}

\begin{theorem}
Given an algebraic manifold $X\subset S^{N-1}_{\mathbb C,+}$, coming via
$$C(X)=C(S^{N-1}_{\mathbb C,+})\Big/\Big<f_\alpha(x_1,\ldots,x_N)=0\Big>$$
we have a morphism of $C^*$-algebras as follows, 
$$\pi:C(X)\to B(H)\quad,\quad x_i\to T_i$$
whenever the operators $T_i\in B(H)$ satisfy the following relations: 
$$\sum_iT_iT_i^*=\sum_iT_i^*T_i=1\quad,\quad 
f_\alpha(T_1,\ldots,T_N)=0$$
Moreover, we can always find a Hilbert space $H$ and operators $\{T_i\}$ such that $\pi$ is faithful.
\end{theorem}

\begin{proof}
Here the first assertion is more of an empty statement, explaining the definition of the algebra $C(X)$, via generators and relations, and the second assertion is something non-trivial, coming as a consequence of the GNS theorem.
\end{proof}

In practice now, all this is a bit too general, and not very useful. We need a good family of target algebras $B$, that we understand well. And here, we can use:

\index{random matrix}
\index{random matrix algebra}

\begin{definition}
A random matrix $C^*$-algebra is an algebra of type
$$B=M_K(C(T))$$
with $T$ being a compact space, and $K\in\mathbb N$ being an integer.
\end{definition}

The terminology here comes from the fact that, in practice, the space $T$ usually comes with a probability measure on it, which makes the elements of $B$ ``random matrices''. Observe that we can write our random matrix algebra as follows:
$$B=M_K(\mathbb C)\otimes C(T)$$

Thus, the random matrix algebras appear by definition as tensor products of the simplest types of $C^*$-algebras that we know, namely the full matrix algebras, $M_K(\mathbb C)$ with $K\in\mathbb N$, and the commutative algebras, $C(T)$, with $T$ being a compact space. Getting back now to our modelling questions for manifolds, we can formulate:

\index{matrix model}

\begin{definition}
A matrix model for a noncommutative algebraic manifold 
$$X\subset S^{N-1}_{\mathbb C,+}$$
is a morphism of $C^*$-algebras of the following type,
$$\pi:C(X)\to M_K(C(T))$$
with $T$ being a compact space, and $K\in\mathbb N$ being an integer.
\end{definition}

As a first observation, when $X$ happens to be classical, we can take $K=1$ and $T=X$, and we have a faithful model for our manifold, namely:
$$id:C(X)\to M_1(C(X))$$

In general, we cannot use $K=1$, and the smallest value $K\in\mathbb N$ doing the job, if any, will correspond somehow to the ``degree of noncommutativity'' of our manifold. 

\bigskip

Before getting into this, we would like to clarify a few abstract issues. As mentioned above, the algebras of type $B=M_K(C(T))$ are called random matrix $C^*$-algebras. The reason for this is the fact that most of the interesting compact spaces $T$ come by definition with a natural probability measure of them. Thus, $B$ is a subalgebra of the bigger algebra $B''=M_K(L^\infty(T))$, usually known as a ``random matrix algebra''.

\bigskip

This perspective is quite interesting for us, because most of our examples of manifolds $X\subset X^{N-1}_{\mathbb C,+}$ appear as homogeneous spaces, and so are measured spaces too. Thus, we can further ask for our models $C(X)\to M_K(C(T))$ to extend into models of the following type, which can be of help in connection with integration problems:
$$L^\infty(X)\to M_K(L^\infty(T))$$

In short, time now to talk about $L^\infty$-functions, in the noncommutative setting. 

\section*{12b. Von Neumann algebras}

In order to discuss all this, we will need some basic von Neumann algebra theory, coming as a complement to the $C^*$-algebra theory from chapter 1. Let us start with a key result in functional analysis, as follows:

\index{topologies on operators}
\index{weak operator topology}
\index{weak topology}
\index{strong operator topology}

\begin{proposition}
For an operator algebra $A\subset B(H)$, the following are equivalent:
\begin{enumerate}
\item $A$ is closed under the weak operator topology, making each of the linear maps $T\to<Tx,y>$ continuous.

\item $A$ is closed under the strong operator topology, making each of the linear maps $T\to Tx$ continuous.
\end{enumerate}
In the case where these conditions are satisfied, $A$ is closed under the norm topology.
\end{proposition}

\begin{proof}
There are several statements here, the proof being as follows:

\medskip

(1) It is clear that the norm topology is stronger than the strong operator topology, which is in turn stronger than the weak operator topology. At the level of the subsets $S\subset B(H)$ which are closed things get reversed, in the sense that weakly closed implies strongly closed, which in turn implies norm closed. Thus, we are left with proving that for any algebra $A\subset B(H)$, strongly closed implies weakly closed.

\medskip

(2) But this latter fact is something standard, which can be proved via an amplification trick. Consider the Hilbert space obtained by summing $n$ times $H$ with itself:
$$K=H\oplus\ldots\oplus H$$

The operators over $K$ can be regarded as being square matrices with entries in $B(H)$, and in particular, we have a representation $\pi:B(H)\to B(K)$, as follows:
$$\pi(T)=\begin{pmatrix}
T\\
&\ddots\\
&&T
\end{pmatrix}$$

Assume now that we are given an operator $T\in\bar{A}$, with the bar denoting the weak closure. We have then, by using the Hahn-Banach theorem, for any $x\in K$:
\begin{eqnarray*}
T\in\bar{A}
&\implies&\pi(T)\in\overline{\pi(A)}\\
&\implies&\pi(T)x\in\overline{\pi(A)x}\\
&\implies&\pi(T)x\in\overline{\pi(A)x}^{\,||.||}
\end{eqnarray*}

Now observe that the last formula tells us that for any $x=(x_1,\ldots,x_n)$, and any $\varepsilon>0$, we can find $S\in A$ such that the following holds, for any $i$:
$$||Sx_i-Tx_i||<\varepsilon$$

Thus $T$ belongs to the strong operator closure of $A$, as desired.
\end{proof}

In the above the terminology, while standard, is a bit confusing, because the norm topology is stronger than the strong operator topology. As a solution, we agree in what follows to call the norm topology ``strong'', and the weak and strong operator topologies ``weak'', whenever these two topologies coincide. With this convention, the algebras from Proposition 12.6 are those which are weakly closed, and we can formulate:

\index{weak topology}
\index{von Neumann algebra}

\begin{definition}
A von Neumann algebra is a $*$-algebra of operators
$$A\subset B(H)$$
which is closed under the weak topology.
\end{definition}

As basic examples, we have the algebra $B(H)$ itself, then the singly generated von Neumann algebras, $A=<T>$, with $T\in B(H)$, and then the multiply generated von Neumann algebras, namely $A=<T_i>$, with $T_i\in B(H)$. At the level of the general results, we first have the bicommutant theorem of von Neumann, as follows:

\index{bicommutant}

\begin{theorem}
For a $*$-algebra $A\subset B(H)$, the following are equivalent:
\begin{enumerate}
\item $A$ is weakly closed, so it is a von Neumann algebra.

\item $A$ equals its algebraic bicommutant $A''$, taken inside $B(H)$.
\end{enumerate}
\end{theorem}

\begin{proof}
Since the commutants are automatically weakly closed, it is enough to show that weakly closed implies $A=A''$. For this purpose, we will prove something a bit more general, stating that given a $*$-algebra of operators $A\subset B(H)$, the following holds, with $A''$ being the bicommutant inside $B(H)$, and with $\bar{A}$ being the weak closure:
$$A''=\bar{A}$$

We prove this equality by double inclusion, as follows:

\medskip

``$\supset$'' Since any operator commutes with the operators that it commutes with, we have a trivial inclusion $S\subset S''$, valid for any set $S\subset B(H)$. In particular, we have:
$$A\subset A''$$

Our claim now is that the algebra $A''$ is closed, with respect to the strong operator topology. Indeed, assuming that we have $T_i\to T$ in this topology, we have:
\begin{eqnarray*}
T_i\in A''
&\implies&ST_i=T_iS,\ \forall S\in A'\\
&\implies&ST=TS,\ \forall S\in A'\\
&\implies&T\in A
\end{eqnarray*}

Thus our claim is proved, and together with Proposition 12.6, which allows us to pass from the strong to the weak operator topology, this gives the desired inclusion:
$$\bar{A}\subset A''$$

``$\subset$'' Here we must prove that we have the following implication, valid for any $T\in B(H)$, with the bar denoting as usual the weak operator closure:
$$T\in A''\implies T\in\bar{A}$$

For this purpose, we use the same amplification trick as in the proof of Proposition 12.5. Consider the Hilbert space obtained by summing $n$ times $H$ with itself:
$$K=H\oplus\ldots\oplus H$$

The operators over $K$ can be regarded as being square matrices with entries in $B(H)$, and in particular, we have a representation $\pi:B(H)\to B(K)$, as follows:
$$\pi(T)=\begin{pmatrix}
T\\
&\ddots\\
&&T
\end{pmatrix}$$

The idea will be that of doing the computations in this representation. First, in this representation, the image of our algebra $A\subset B(H)$ is given by:
$$\pi(A)=\left\{\begin{pmatrix}
T\\
&\ddots\\
&&T
\end{pmatrix}\Big|T\in A\right\}$$

We can compute the commutant of this image, exactly as in the usual scalar matrix case, and we obtain the following formula:
$$\pi(A)'=\left\{\begin{pmatrix}
S_{11}&\ldots&S_{1n}\\
\vdots&&\vdots\\
S_{n1}&\ldots&S_{nn}
\end{pmatrix}\Big|S_{ij}\in A'\right\}$$

We conclude from this that, given an operator $T\in A''$ as above, we have:
$$\begin{pmatrix}
T\\
&\ddots\\
&&T
\end{pmatrix}\in\pi(A)''$$

In other words, the conclusion of all this is that we have:
$$T\in A''\implies \pi(T)\in\pi(A)''$$

Now given a vector $x\in K$, consider the orthogonal projection $P\in B(K)$ on the norm closure of the vector space $\pi(A)x\subset K$. Since the subspace $\pi(A)x\subset K$ is invariant under the action of $\pi(A)$, so is its norm closure inside $K$, and we obtain from this:
$$P\in\pi(A)'$$

By combining this with what we found above, we conclude that we have:
$$T\in A''\implies \pi(T)P=P\pi(T)$$

Now since this holds for any $x\in K$, we conclude that any $T\in A''$ belongs to the strong operator closure of $A$. By using now Proposition 12.5, which allows us to pass from the strong to the weak operator closure, we conclude that we have $A''\subset\bar{A}$, as desired. 
\end{proof}

In order to develop now some general theory, let us start by investigating the finite dimensional case. Here the ambient operator algebra is $B(H)=M_N(\mathbb C)$, and any subspace $A\subset B(H)$ is automatically closed, for all 3 topologies from Proposition 12.6. Thus, we are left with the question of investigating the $*$-algebras of usual matrices $A\subset M_N(\mathbb C)$. But this is a purely algebraic question, whose answer is as follows:

\index{finite dimensional algebra}

\begin{theorem}
The $*$-algebras $A\subset M_N(\mathbb C)$ are exactly the algebras of the form
$$A=M_{r_1}(\mathbb C)\oplus\ldots\oplus M_{r_k}(\mathbb C)$$
depending on parameters $k\in\mathbb N$ and $r_1,\ldots,r_k\in\mathbb N$ satisfying
$$r_1+\ldots+r_k=N$$
embedded into $M_N(\mathbb C)$ via the obvious block embedding, twisted by a unitary $U\in U_N$.
\end{theorem}

\begin{proof}
We have two assertions to be proved, the idea being as follows:

\medskip

(1) Given numbers $r_1,\ldots,r_k\in\mathbb N$ satisfying $r_1+\ldots+r_k=N$, we have an obvious embedding of $*$-algebras, via matrix blocks, as follows:
$$M_{r_1}(\mathbb C)\oplus\ldots\oplus M_{r_k}(\mathbb C)\subset M_N(\mathbb C)$$

In addition, we can twist this embedding by a unitary $U\in U_N$, as follows:
$$M\to UMU^*$$

(2) In the other sense now, consider an arbitrary $*$-algebra of the $N\times N$ matrices, $A\subset M_N(\mathbb C)$. Let us first look at the center of this algebra, which given by:
$$Z(A)=A\cap A'$$

It is elementary to prove that this center, as an algebra, is of the following form:
$$Z(A)\simeq\mathbb C^k$$

Consider now the standard basis $e_1,\ldots,e_k\in\mathbb C^k$, and let  $p_1,\ldots,p_k\in Z(A)$ be the images of these vectors via the above identification. In other words, these elements $p_1,\ldots,p_k\in A$ are central minimal projections, summing up to 1:
$$p_1+\ldots+p_k=1$$

The idea is then that this partition of the unity will eventually lead to the block decomposition of $A$, as in the statement. We prove this in 4 steps, as follows:

\medskip

\underline{Step 1}. We first construct the matrix blocks, our claim here being that each of the following linear subspaces of $A$ are non-unital $*$-subalgebras of $A$:
$$A_i=p_iAp_i$$

But this is clear, with the fact that each $A_i$ is closed under the various non-unital $*$-subalgebra operations coming from the projection equations $p_i^2=p_i^*=p_i$.

\medskip

\underline{Step 2}. We prove now that the above algebras $A_i\subset A$ are in a direct sum position, in the sense that we have a non-unital $*$-algebra sum decomposition, as follows:
$$A=A_1\oplus\ldots\oplus A_k$$

As with any direct sum question, we have two things to be proved here. First, by using the formula $p_1+\ldots+p_k=1$ and the projection equations $p_i^2=p_i^*=p_i$, we conclude that we have the needed generation property, namely:
$$A_1+\ldots+ A_k=A$$

As for the fact that the sum is indeed direct, this follows as well from the formula $p_1+\ldots+p_k=1$, and from the projection equations $p_i^2=p_i^*=p_i$.

\medskip

\underline{Step 3}. Our claim now, which will finish the proof, is that each of the $*$-subalgebras $A_i=p_iAp_i$ constructed above is a full matrix algebra. To be more precise here, with $r_i=rank(p_i)$, our claim is that we have isomorphisms, as follows:
$$A_i\simeq M_{r_i}(\mathbb C)$$

In order to prove this claim, recall that the projections $p_i\in A$ were chosen central and minimal. Thus, the center of each of the algebras $A_i$ reduces to the scalars:
$$Z(A_i)=\mathbb C$$

But this shows, either via a direct computation, or via the bicommutant theorem, that the each of the algebras $A_i$ is a full matrix algebra, as claimed.

\medskip

\underline{Step 4}. We can now obtain the result, by putting together what we have. Indeed, by using the results from Step 2 and Step 3, we obtain an isomorphism as follows:
$$A
=A_1\oplus\ldots\oplus A_k
\simeq M_{r_1}(\mathbb C)\oplus\ldots\oplus M_{r_k}(\mathbb C)$$

Moreover, a careful look at the isomorphisms established in Step 3 shows that at the global level, of the algebra $A$ itself, the above isomorphism comes by twisting the standard multimatrix embedding $M_{r_1}(\mathbb C)\oplus\ldots\oplus M_{r_k}(\mathbb C)\subset M_N(\mathbb C)$, discussed in the beginning of the proof, (1) above, by a certain unitary $U\in U_N$. Thus, we obtain the result.
\end{proof}

As an application of Theorem 12.9, clarifying the relation with linear algebra, or operator theory in finite dimensions, we have the following result:

\begin{proposition}
Given an operator $T\in B(H)$ in finite dimensions, $H=\mathbb C^N$, the von Neumann algebra $A=<T>$ that it generates inside $B(H)=M_N(\mathbb C)$ is
$$A=M_{r_1}(\mathbb C)\oplus\ldots\oplus M_{r_k}(\mathbb C)$$
with the sizes of the blocks $r_1,\ldots,r_k\in\mathbb N$ coming from the spectral theory of the associated matrix $M\in M_N(\mathbb C)$. In the normal case $TT^*=T^*T$, this decomposition comes from
$$T=UDU^*$$
with $D\in M_N(\mathbb C)$ diagonal, and with $U\in U_N$ unitary.
\end{proposition}

\begin{proof}
This is something standard, by using the basic linear algebra theory and spectral theory for the usual matrices $M\in M_N(\mathbb C)$.
\end{proof}

Let us get now to infinite dimensions, with Proposition 12.10 as our main source of inspiration. We have here the following result:

\index{normal operator}
\index{spectral theorem}
\index{commutative von Neumann algebra}

\begin{theorem}
Given an operator $T\in B(H)$ which is normal, 
$$TT^*=T^*T$$
the von Neumann algebra $A=<T>$ that it generates inside $B(H)$ is
$$<T>=L^\infty(\sigma(T))$$
with $\sigma(T)$ being its spectrum, formed of numbers $\lambda\in\mathbb C$ such that $T-\lambda$ is not invertible.
\end{theorem}

\begin{proof}
This is something standard as well, by using the spectral theory for the normal operators $T\in B(H)$, coming from chapter 1.
\end{proof}

More generally, along the same lines, we have the following result, dealing this time with commuting families of normal operators:

\begin{theorem}
Given operators $T_i\in B(H)$ which are normal, and which commute, the von Neumann algebra $A=<T_i>$ that these operators generates inside $B(H)$ is
$$<T_i>=L^\infty(X)$$
with $X$ being a certain measured space, associated to the family $\{T_i\}$.
\end{theorem}

\begin{proof}
This is again routine, by using this time the spectral theory for the families of commuting normal operators $T_i\in B(H)$. See for instance Blackadar \cite{bla}.
\end{proof}

As an interesting abstract consequence of this, we have:

\index{commutative von Neumann algebra}

\begin{theorem}
The commutative von Neumann algebras are the algebras of type
$$A=L^\infty(X)$$
with $X$ being a measured space.
\end{theorem}

\begin{proof}
We have two assertions to be proved, the idea being as follows:

\medskip

(1) In one sense, we must prove that given a measured space $X$, we can realize the commutative algebra $A=L^\infty(X)$ as a von Neumann algebra, on a certain Hilbert space $H$. But this is something that we already know, coming from the multiplicity operators $T_f(g)=fg$ from the proof of the GNS theorem, the representation being as follows: 
$$L^\infty(X)\subset B(L^2(X))$$

(2) In the other sense, given a commutative von Neumann algebra $A\subset B(H)$, we must construct a certain measured space $X$, and an identification $A=L^\infty(X)$. But this follows from Theorem 12.12, because we can write our algebra as follows:
$$A=<T_i>$$

To be more precise, $A$ being commutative, any element $T\in A$ is normal. Thus, we can pick a basis $\{T_i\}\subset A$, and then we have $A=<T_i>$ as above, with $T_i\in B(H)$ being commuting normal operators. Thus Theorem 12.12 applies, and gives the result.

\medskip

(3) Alternatively, and more explicitly, we can deduce this from Theorem 12.11, applied with $T=T^*$. Indeed, by using $T=Re(T)+iIm(T)$, we conclude that any von Neumann algebra $A\subset B(H)$ is generated by its self-adjoint elements $T\in A$. Moreover, by using measurable functional calculus, we conclude that $A$ is linearly generated by its projections. But then, assuming $A=\overline{span}\{p_i\}$, with $p_i$ being projections, we can set:
$$T=\sum_{i=0}^\infty\frac{p_i}{3^i}$$

Then $T=T^*$, and by functional calculus we have $p_0\in<T>$, then $p_1\in <T>$, and so on. Thus $A=<T>$, and $A=L^\infty(X)$ comes now via Theorem 12.11, as claimed.
\end{proof}

Moving ahead now, we can combine Proposition 12.8 with Theorem 12.13, and by building along the lines of Theorem 12.9, but this time in infinite dimensions, we are led to the following statement, due to Murray-von Neumann and Connes:

\index{reduction theory}
\index{factor}
\index{Connes classification}

\begin{theorem}
Given a von Neumann algebra $A\subset B(H)$, if we write its center as 
$$Z(A)=L^\infty(X)$$
then we have a decomposition as follows, with the fibers $A_x$ having trivial center:
$$A=\int_XA_x\,dx$$
Moreover, the factors, $Z(A)=\mathbb C$, can be basically classified in terms of the ${\rm II}_1$ factors, which are those satisfying $\dim A=\infty$, and having a faithful trace $tr:A\to\mathbb C$.
\end{theorem}

\begin{proof}
This is something that we know to hold in finite dimensions, as a consequence of Theorem 12.9. In general, this is something heavy, the idea being as follows:

\medskip

(1) This is von Neumann's reduction theory main result, whose statement is already quite hard to understand, and whose proof uses advanced functional analysis.

\medskip

(2) This is heavy, due to Murray-von Neumann and Connes, the idea being that the other factors can be basically obtained via crossed product constructions.
\end{proof}

All this is certainly quite advanced, taking substantial time to be fully understood. For general reading on von Neumann algebras we recommend the book of Blackadar \cite{bla}, but be aware tough that, while being at the same time well-written, condensed and reasonably thick, that book is only an introduction to Theorem 12.14. So, if we want to learn the full theory, with the complete proof of Theorem 12.14, you will have to go, every now and then, through the original papers of Murray-von Neumann and Connes.

\bigskip

By the way, talking von Neumann and Connes, we can only warmly recommend their books, \cite{von} and \cite{co1}, as a complement to anything that you might want to learn on operator algebras, from Blackadar \cite{bla} or from somewhere else. The point indeed is that von Neumann algebras come from quantum mechanics, meaning they are designed to help with quantum mechanics, and also happen to have 0 serious applications to pure mathematics, which by the way is not really science anyway, and if there are 2 people who understood all this, what von Neumann algebras are potentially good for, these are von Neumann himself, and Connes. So, read their books, \cite{von} and \cite{co1}.

\bigskip

For the discussion to be complete, yet another 2 people who understood what von Neumann algebras are good for, in more modern times, are Jones and Voiculescu. So, have their main writings, \cite{jo3} and \cite{vdn}, available nearby, ready for some reading. And finally, talking quantum mechanics, always a pleasure to recommend, as usual, the books of Feynman \cite{fey}, Griffiths \cite{gr2} and Weinberg \cite{we2}. And with a preference for Griffiths \cite{gr2}, that's got all the quantum mechanics that you need to know, clearly explained, and comes with a cat on the cover too. But Feynman \cite{fey} and Weinberg \cite{we2} are excellent too, and as a general rule, as long as you stay away from quantum mechanics books which claim to be ``rigorous'', ``axiomatic'', ``mathematical'', and so on, things fine. 

\bigskip

Always remember here, as per Feynman saying, that ``no one understands quantum mechanics''. With this being of course a euphemism for something of type ``quantum mechanics as we know it is wrong, sorry for that, and we're working on the fix''. 

\bigskip

Now back to work, and our noncommutative geometry questions, as a first application of the above, we can extend our noncommutative space setting, as follows:

\index{noncommutative measured space}
\index{quantum measured space}
\index{GNS construction}

\begin{theorem}
Consider the category of ``noncommutative measure spaces'', having as objects the pairs $(A,tr)$ consisting of a von Neumann algebra with a faithful trace, and with the arrows reversed, which amounts in writing $A=L^\infty(X)$ and $tr=\int_X$.
\begin{enumerate}
\item The category of usual measured spaces embeds into this category, and we obtain in this way the objects whose associated von Neumann algebra is commutative.

\item Each $C^*$-algebra given with a trace produces as well a noncommutative measure space, by performing the GNS construction, and taking the weak closure.

\item In what regards the finitely generated group duals, or more generally the compact matrix quantum groups, the corresponding identification is injective.

\item Even more generally, for noncommutative algebraic manifolds having an integratiuon functional, like the spheres, the identification is injective.
\end{enumerate}
\end{theorem}

\begin{proof}
This is clear indeed from the basic properties of the GNS construction, from Theorem 12.2, and from the general theory from Theorem 12.14.
\end{proof}

Before getting back to matrix models, we would like to formulate the following result, in relation with our axiomatization questions from chapters 1-4:

\begin{theorem}
In the context of noncommutative geometries coming from quadruplets $(S,T,U,K)$, we have von Neumann algebras, with traces, as follows,
$$\xymatrix@R=50pt@C=50pt{
L^\infty(S)\ar[r]\ar[d]\ar[dr]&L^\infty(T)\ar[l]\ar[d]\ar[dl]\\
L^\infty(U)\ar[u]\ar[ur]\ar[r]&L^\infty(K)\ar[l]\ar[ul]\ar[u]
}$$
with $L^\infty(S)\subset L^\infty(U)$ being obtained by taking the first row algebra.
\end{theorem}

\begin{proof}
This follows indeed from the results that we already have, from chapters 1-4 above, by using the general formalism from Theorem 12.15.
\end{proof}

This statement, which is quite interesting, philosophically speaking, raises the question of axiomatizing, or rather re-axiomatizing, the quadruplets $(S,T,U,K)$ that we are interested in directly in terms of the associated von Neumann algebras, as above. Indeed, in view of our general quantum mechanics motivations, we are after all mostly interested in integrating over our quantum manifolds, and so with this is mind, the von Neumann algebra formalism seems to be the one which is best adapted to our questions.

\bigskip

However, this is wrong. The above result is something theoretical, because it assumes the existence of Haar measures on our spaces $S,T,U,K$, which itself is something coming as a theorem. Thus, while all this is nice, the good way of doing things is with $C^*$-algebras, as we did in chapters 1-4. And the von Neumann algebras from Theorem 12.16 remain something more advanced and specialized, coming afterwards.

\bigskip

\index{Hilbert space}
\index{operator algebra}

As a side comment here, and for ending with some physics, the question ``does the algebra or the Hilbert space come first'' is a well-known one in quantum mechanics, basically leading to 2 different schools of thought. We obviously adhere here to the ``algebra comes first'' school. But let us not get here into this, perhaps enough controversies discussed, so far in this book. For more on this, get to know about the Bohr vs Einstein debate, which is the mother of all debates, in quantum mechanics. 

\bigskip

And then, once this learned, as an instructive exercise: what do you think, from your reading so far of this book, do we rather side with Bohr, or with Einstein? 

\section*{12c. Matrix truncations}

In relation now with the modelling questions that we are interested in here, with all the above operator algebra material digested, we can now go ahead with our program, and discuss von Neumann algebraic extensions. We have the following result:

\index{stationary model}
\index{faithful model}

\begin{theorem}
Given a matrix model $\pi:C(X)\to M_K(C(T))$, with both $X,T$ being assumed to have integration functionals, the following are equivalent:
\begin{enumerate}
\item $\pi$ is stationary, in the sense that $\int_X=(tr\otimes\smallint_T)\pi$.

\item $\pi$ produces an inclusion $\pi':C_{red}(X)\subset M_K(X(T))$.

\item $\pi$ produces an inclusion $\pi'':L^\infty(X)\subset M_K(L^\infty(T))$.
\end{enumerate}
Moreover, in the quantum group case, these conditions imply that $\pi$ is faithful.
\end{theorem}

\begin{proof}
This is standard functional analysis, as follows:

\medskip

(1) Consider the following diagram, with all the solid arrows being by definition the canonical maps between the algebras concerned:
$$
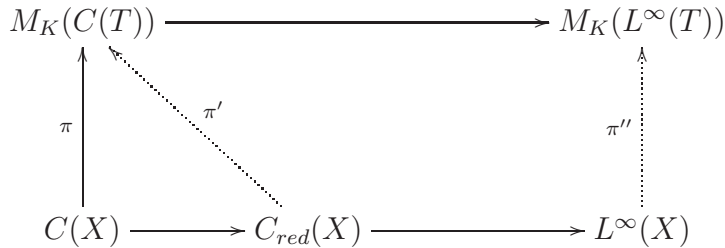
\xymatrix@R=60pt@C=30pt{
M_K(C(T))\ar[rrr]&&&M_K(L^\infty(T))\\
C(X)\ar[u]^\pi\ar[r]&C_{red}(X)\ar[rr]\ar@.[ul]_{\pi'}&&L^\infty(X)\ar@.[u]^{\pi''}
}$$

(2) With this picture in hand, the implications $(1)\iff(2)\iff(3)$ between the conditions (1,2,3) in the statement are all clear, coming from the basic properties of the GNS construction, and of the von Neumann algebras, explained in the above.

\medskip

(3) As for the last assertion, this is something more subtle, coming from the fact that if $L^\infty(G)$ is of type I, as required by (3), then $G$ must be coamenable.
\end{proof}

The above result raises a number of interesting questions, notably in what regards the extension of the last assertion, to the case of more general homogeneous spaces. 

\bigskip

Before going further, we would like to record as well the following key result regarding the matrix models, valid so far in the quantum group case only:

\index{Hopf image}
\index{inner faithfulness}

\begin{theorem}
Consider a matrix model $\pi:C(G)\to M_K(C(T))$ for a closed subgroup $G\subset U_N^+$, with $T$ being assumed to be a compact probability space.
\begin{enumerate}
\item There exists a smallest subgroup $G'\subset G$, producing a factorization of type: 
$$\pi:C(G)\to C(G')\to M_K(C(T))$$
The algebra $C(G')$ is called Hopf image of $\pi$.

\item When $\pi$ is inner faithful, in the sense that $G=G'$, we have the formula 
$$\int_G=\lim_{k\to\infty}\sum_{r=1}^k\varphi^{*r}$$
where $\varphi=(tr\otimes\smallint_T)\pi$, and $\phi*\psi=(\phi\otimes\psi)\Delta$.
\end{enumerate}
\end{theorem}

\begin{proof}
All this is well-known, but quite specialized, the idea being as follows:

\medskip

(1) This follows by dividing the algebra $C(G)$ by a suitable ideal, namely the Hopf ideal generated by the kernel of the matrix model map $\pi:C(G)\to M_K(C(T))$. 

\medskip

(2) This follows by suitably adapting Woronowicz's proof for the existence and formula of the Haar integration functional from \cite{wo1}, to the matrix model situation. 
\end{proof}

The above result is quite important, for a number of reasons. Indeed, as a main application of it, while the existence of a faithful matrix model $\pi:C(G)\subset M_K(C(T))$ forces the $C^*$-algebra $C(G)$ to be of type I, and so $G$ to be coamenable, as already mentioned in the proof of Theorem 12.17, there is no known restriction coming from the existence of an inner faithful model $\pi:C(G)\to M_K(C(T))$. See \cite{ba8}.

\bigskip

In the general manifold setting, talking about such things is in general not possible, unless our manifold $X$ has some extra special structure, as for instance being an homogeneous space, in the spirit of the various such spaces discussed in chapters 5-8. However, in practice, such a theory has not been developed yet.

\bigskip

Let us go back now to our basic notion of a matrix model, from Definition 12.5, and develop some more general theory, in that setting. We first have:

\begin{proposition}
A $1\times1$ model for a manifold $X\subset S^{N-1}_{\mathbb C,+}$ must come from a map
$$p:T\to X_{class} \subset X$$
and $\pi$ is faithful precisely when $X=X_{class}$, and when $p$ is surjective.
\end{proposition}

\begin{proof}
According to our conventions, a $1\times1$ model for a manifold $X\subset S^{N-1}_{\mathbb C,+}$ is simply a morphism of algebras $\pi:C(X)\to C(T)$. Now since $C(T)$ is commutative, this morphism must factorize through the abelianization of $C(X)$, as follows:
$$\pi:C(X)\to C(X_{class})\to C(T)$$

Thus, our morphism $\pi$ must come by transposition from a map $p$, as claimed.
\end{proof}

Following \cite{bb1}, in order to generalize the above trivial fact, we can use:

\begin{definition}
Let  $X\subset S^{N-1}_{\mathbb C,+}$. We define a closed subspace $X^{(K)}\subset X$ by
$$C(X^{(K)})=C(X)/J_K$$
where $J_K$ is the common null space of matrix representations of $C(X)$, of size $L\leq K$,
$$J_K=\bigcap_{L\leq K}\ \bigcap_{\pi:C(X)\to M_L(\mathbb C)}\ker(\pi)$$
and we call $X^{(K)}$ the ``part of $X$ which is realizable with $K\times K$ models''.
\end{definition}

As a basic example here, the first such space, at $K=1$, is the classical version:
$$X^{(1)}=X_{class}$$

Observe that we have embeddings of quantum spaces, as follows:
$$X^{(1)}\subset X^{(2)}\subset X^{(3)}\ldots\ldots\subset X$$

As a first result now on these spaces, we have the following well-known fact:

\begin{theorem}
The increasing union of compact quantum spaces
$$X^{(\infty)}=\bigcup_{K\geq1}X^{(K)}$$ 
equals $X$ precisely when the algebra $C(X)$ is residually finite dimensional.
\end{theorem}

\begin{proof}
This is something well-known. We refer to Chirvasitu \cite{chi} for a discussion on this topic, in the context of the quantum groups, and to \cite{bb1} for more.
\end{proof}

Getting back now to the case $K<\infty$, we first have, following \cite{bb1}:

\begin{proposition}
Consider an algebraic manifold $X\subset S^{N-1}_{\mathbb C,+}$.
\begin{enumerate}
\item Given a closed subspace $Y\subset X\subset S^{N-1}_{\mathbb C,+}$, we have $Y \subset X^{(K)}$ precisely when any irreducible representation of $C(Y)$ has dimension $\leq K$. 
 
\item In particular, we have $X^{(K)} = X$ precisely when any irreducible representation of $C(X)$ has dimension $\leq K$.
\end{enumerate}
\end{proposition}

\begin{proof}
This follows from general $C^*$-algebra theory, as follows:

\medskip

(1) If any irreducible representation of $C(Y)$ has dimension $\leq K$, then we have $Y\subset X^{(K)}$, because the irreducible representations of a $C^*$-algebra separate its points. Conversely, assuming $Y\subset X^{(K)}$, it is enough to show that  any irreducible representation of the algebra $C(X^{(K)})$ has dimension $\leq K$. But this is once again well-known.

\medskip

(2) This follows indeed from (1).
\end{proof}

The connection with the previous considerations comes from:

\begin{theorem}
If $X\subset S^{N-1}_{\mathbb C,+}$ has a faithful matrix model
$$C(X)\to M_K(C(T))$$
then we have $X=X^{(K)}$.
\end{theorem}

\begin{proof}
This follows from the above and from the standard representation theory for the $C^*$-algebras. For full details on all this, we refer as before to \cite{bb1}.
\end{proof}

We can now discuss the universal $K\times K$-matrix model, constructed as follows:

\begin{theorem}
Given  $X\subset S^{N-1}_{\mathbb C,+}$ algebraic, the category of its $K\times K$ matrix models, with $K\geq 1$ being fixed, has a universal object as follows:
$$\pi_K:C(X)\to M_K(C(T_K))$$
That is, given a model $\rho:C(X)\to M_K(C(T))$, we have a diagram of type
$$\xymatrix{C(X)\ar[rr]^{\pi}\ar[rd]_{\rho}&&M_K(C(T_K))\ar[dl]\\&M_K(C(T))&}$$
where the map on the right is unique, and arises from a continuous map $T\to T_K$.
\end{theorem}

\begin{proof}
Consider the universal commutative $C^*$-algebra generated by elements $x_{ij}(a)$, with $1 \leq i,j \leq K$ and $a \in \mathcal O(X)$, subject to the following relations:
$$x_{ij}(a+\lambda b) = x_{ij}(a)+\lambda x_{ij}(b)$$
$$x_{ij}(ab)= \sum_kx_{ik}(a)x_{kj}(b)$$
$$x_{ij}(1)=\delta_{ij}$$
$$x_{ij}(a)^*= x_{ji}(a^*)$$  

This algebra is indeed well-defined because of the following relations:
$$\sum_l \sum_k x_{ik}(z_l^*)x_{ki}(z_l)=1$$

Now let $T_K$ be the spectrum of this algebra. Since $X$ is algebraic, we have:
$$\pi : C(X) \rightarrow M_K(C(T_K))\quad,\quad 
\pi(z_k) = (x_{ij}(z_k))$$

By construction of $T_K$ and $\pi$, we have the universal matrix model. See \cite{bb1}.
\end{proof}

Still following \cite{bb1}, as an illustration for the above, we have:

\begin{proposition}
Let $X\subset S^{N-1}_{\mathbb C,+}$ with $X$ algebraic and $X_{class} \not= \emptyset$, and let \
$$\pi:C(X)\to M_K(C(T_K))$$
be the universal matrix model. Then we have 
$$C(X^{(K)})= C(X)/Ker(\pi)$$
and hence $X=X^{(K)}$ if and only if $X$ has a faithful $K \times K$-matrix model.
\end{proposition}

\begin{proof}
We have to prove that $Ker(\pi)=J_K$, the latter ideal being the intersection of the kernels of all matrix representations as follows, with $ L \leq K$:
$$C(X) \rightarrow M_L(\mathbb C)$$

For $a \not \in Ker(\pi)$, we see that $a \not \in J_K$ by evaluating at an appropriate element of $T_K$. Conversely, assume that we are given $a \in Ker(\pi)$. Let $\rho : C(X) \rightarrow M_L(\mathbb C)$ be a representation with $L \leq K$, and let $\varepsilon : C(X) \rightarrow \mathbb C$ be a representation. We can extend $\rho$ to a representation $\rho' : C(X) \rightarrow M_K(\mathbb C)$ by letting, for any $b \in C(X)$:
$$\rho'(b)=
\begin{pmatrix} 
\rho(b)&0\\
0&\varepsilon(b)I_{K-L}  
\end{pmatrix}$$
   
The universal property of the universal matrix model yields that $\rho'(a)=0$, since $\pi(a)=0$. Thus $\rho(a)=0$. We therefore have $a \in J_K$, and $Ker(\pi)\subset J_K$, and the first statement is proved. The last statement follows from the first one. See \cite{bb1}.
\end{proof}

Next, we have the following result, also from \cite{bb1}:

\begin{proposition}
Let $X\subset S^{N-1}_{\mathbb C,+}$ be algebraic, and satisfying:
$$X_{class} \not=\emptyset$$
Then $X^{(K)}$ is algebraic as well.
\end{proposition}

\begin{proof}
We keep the notations above, and consider the following map:
$$\pi_0 : \mathcal O(X) \rightarrow M_K(C(T_K))\quad,\quad 
z_l \to (x_{ij}(z_l))$$

This induces a $*$-algebra map, as follows: 
$$\tilde{\pi_0} : C^*(\mathcal O(X)/Ker(\pi_0)) \rightarrow M_K(C(T_K))$$

We need to show that $\tilde{\pi}_0$ is injective. For this purpose, observe that the universal model factorizes as follows, where $p$ is canonical surjection:
$$\pi : C(X) \overset{p}\to C^*(\mathcal O(X)/Ker(\pi_0))\overset{\tilde{\pi}_0}\to M_K(C(T_K))$$

We therefore obtain $Ker(\pi)= Ker(p)$, and we conclude that:
$$C(X^{(K)})
=C(X)/Ker(p)
=C^*(\mathcal O(X)/Ker(\pi_0))$$

Thus $X^{(K)}$ is indeed algebraic. Since $\mathcal O(X)/Ker(\pi_0)$ is isomorphic to a $*$-subalgebra of $M_K(C(T_K))$, it satisfies the standard Amitsur-Levitski polynomial identity:
$$S_{2K}(x_1,\ldots,x_{2K})=0$$

By density, so does $C^*(\mathcal O(X)/Ker(\pi_0))$. Thus any irreducible representation of the algebra $C^*(\mathcal O(X)/Ker(\pi_0))$ has dimension $\leq K$. Consider now an element as follows:
$$a \in  C^*(\mathcal O(X)/Ker(\pi_0))$$

Assuming $a\neq0$ we can, by the same reasoning as in the previous proof, find a representation as follows, such that $\rho(a)\not=0$:
$$\rho:C^*(\mathcal O(X)/Ker(\pi_0))\to M_K(\mathbb C)$$

Indeed, a given algebra map $\varepsilon : C(X) \rightarrow \mathbb C$ induces an algebra map as follows:
$$C(T_K)\to\mathbb C\quad,\quad 
x_{ij}(a)\to\delta_{ij}\varepsilon(a)$$

But this map enables us to extend representations, as before. By construction the universal model space yields an algebra  map as follows:
$$M_K(C(T_K)) \rightarrow M_K(\mathbb C)$$

The composition with $\tilde{\pi_0}p=\pi$ is then $\rho p$, so $\tilde{\pi_0}(a)\not=0$, and $\tilde{\pi}_0$ is injective.
\end{proof}

Summarizing, we have proved the following result:

\index{truncations}
\index{matrix model truncations}

\begin{theorem}
Let $X\subset S^{N-1}_{\mathbb C,+}$ be algebraic, satisfying $X_{class}\not=\emptyset$. 
Then we have an increasing sequence of algebraic submanifolds
$$X_{class}= X^{(1)}\subset X^{(2)}\subset X^{(3)}\subset\ldots\ldots\subset X$$
where $X^{(K)}$ is given by the fact that
$$C(X^{(K)}) \subset M_K(C(T_K))$$
is obtained by factorizing the universal matrix model. 
\end{theorem}

\begin{proof}
This follows indeed from the above results. See \cite{bb1}.
\end{proof}

There are many other things that can be said about the above matrix truncations $X^{(K)}$, and we refer here to \cite{bb1} and related papers. However, the main problem remains that of suitably fine-tuning this theory, as to make it compatible with the theory of matrix models for the Woronowicz algebras, which itself is something quite advanced, and rather satisfactory. To be more precise here, the situation is as follows:

\bigskip

(1) As a first observation, when taking as input a quantum group, $X=G$, the above truncation procedure does not produce a quantum group at $K\geq2$, because the compultiplication $\Delta$ does not factorize. Thus, Theorem 12.27 as stated remains something a bit orthogonal to what is known about the matrix models for quantum groups.

\bigskip

(2) Conversely, as already said before, the main results on the matrix models for quantum groups regard the notion of inner faithfulness from Theorem 12.18. And such results cannot extend to general manifolds $X\subset S^{N-1}_{\mathbb C,+}$, unless we are dealing with special classes of homogeneous spaces, in the spirit of those discussed in chapters 5-8.

\bigskip

Summarizing, many things to be done. The main problem is probably that of talking about inner faithful models for affine homogeneous spaces, but the general theory here is unknown, at least so far. Finally, let us mention that, in the quantum group setting, the known theory of matrix models was heavily inspired by the work of Jones \cite{jo1}, \cite{jo2}, \cite{jo3}, in connection with general problems in statistical mechanics, and in what regards the extension of this to the case of more general homogeneous spaces, or other algebraic manifolds, the motivations remain a bit too advanced to be fully understood.

\bigskip

In short, the bet would be that the matrix models for affine homogeneous spaces can be axiomatized and understood mathematically, notably with a notion of inner faithfulness for them, and then can be useful in connection with certain questions at the interface between quantum mechanics and statistical mechanics. And that is all we can say, for the moment. Sometimes authors do not really understand what they are talking about, in their own books, and not that this ever happened to me, but the present discussion starts to be a bit too complicated for unexperienced readers, so time to stop here.

\section*{12d. Half-liberation}

As a nice illustration for the above modelling theory, let us discuss now the half-liberation operation, which is connected to $X^{(2)}$, as a continuation of the material from chapter 9. We first restrict the attention to the real case. Let us start with:

\index{half-classical manifold}
\index{half-classical version}

\begin{definition}
The half-classical version of a manifold $X\subset S^{N-1}_{\mathbb R,+}$ is given by:
$$C(X^*)=C(X)\Big/\left<abc=cba\Big|\forall a,b,c\in\{x_i\}\right>$$
We say that $X$ is half-classical when $X=X^*$.
\end{definition}

Observe the obvious similarity with the construction of the classical version. In fact, philosophically, this definition is some sort of ``next level'' definition for the classical version, assuming that you managed, via some sort of yoga, to be as familiar with half-commutation, $abc=cba$, as you are with usual commutation, $ab=ba$.

\bigskip

In order to understand now the structure of $X^*$, we can use an old matrix model method, which goes back to Bichon-Dubois-Violette \cite{bdu}, and then to Bichon \cite{bic}. This is based on the following observation, that we already met in chapter 9:

\begin{proposition}
For any $z\in\mathbb C^N$, the matrices
$$X_i=\begin{pmatrix}0&z_i\\ \bar{z}_i&0\end{pmatrix}$$
are self-adjoint, and half-commute.
\end{proposition}

\begin{proof}
The matrices $X_i$ are indeed self-adjoint, and their products are given by:
$$X_iX_j
=\begin{pmatrix}0&z_i\\ \bar{z}_i&0\end{pmatrix}\begin{pmatrix}0&z_j\\ \bar{z}_j&0\end{pmatrix}
=\begin{pmatrix}z_i\bar{z}_j&0\\ 0&\bar{z}_iz_j\end{pmatrix}$$

Also, we have as well the following formula:
$$X_iX_jX_k
=\begin{pmatrix}z_i\bar{z}_j&0\\ 0&\bar{z}_iz_j\end{pmatrix}\begin{pmatrix}0&z_k\\ \bar{z}_k&0\end{pmatrix}
=\begin{pmatrix}0&z_i\bar{z}_jz_k\\ \bar{z}_iz_j\bar{z}_k&0\end{pmatrix}$$

Now since this latter quantity is symmetric in $i,k$, we obtain from this that we have the half-commutation formula $X_iX_jX_k=X_kX_jX_i$, as desired. 
\end{proof}

The idea now, following Bichon-Dubois-Violette \cite{bdu} and Bichon \cite{bic}, will be that of using the matrices in Proposition 12.29 in order to model the coordinates of the arbitrary half-classical manifolds. In order to connect the algebra of the classical coordinates $z_i$ to that of the noncommutative coordinates $X_i$, we will need an abstract definition:  

\begin{definition}
Given a noncommutative polynomial $f\in\mathbb R<x_1,\ldots,x_N>$ in $N$ variables, we define a usual polynomial in $2N$ variables
$$f^\circ\in\mathbb R[z_1,\ldots,z_N,\bar{z}_1,\ldots,\bar{z}_N]$$
according to the formula
$$f=x_{i_1}x_{i_2}x_{i_3}x_{i_4}\ldots\implies f^\circ=z_{i_1}\bar{z}_{i_2}z_{i_3}\bar{z}_{i_4}\ldots$$
in the monomial case, and then by extending this correspondence, by linearity.
\end{definition}

As a basic example here, the polynomial defining the free real sphere $S^{N-1}_{\mathbb R,+}$ produces in this way the polynomial defining the complex sphere $S^{N-1}_\mathbb C$:
$$f=x_1^2+\ldots+x_N^2\implies f^\circ=|z_1|^2+\ldots+|z_N|^2$$

Also, given a polynomial $f\in\mathbb R<x_1,\ldots,x_N>$, we can decompose it into its even and odd parts, $f=g+h$, by putting into $g/h$ the monomials of even/odd length. Observe that with $z=(z_1,\ldots,z_N)$, these odd and even parts are given by:
$$g(z)=\frac{f(z)+f(-z)}{2}\quad,\quad 
h(z)=\frac{f(z)-f(-z)}{2}$$

With these conventions, we have the following result:

\begin{proposition}
Given a manifold $X$, coming from a family of noncommutative polynomials $\{f_\alpha\}\subset\mathbb R<x_1,\ldots,x_N>$, we have a morphism algebras
$$\pi:C(X)\to M_2(\mathbb C)\quad,\quad \pi(x_i)=\begin{pmatrix}0&z_i\\ \bar{z}_i&0\end{pmatrix}$$
precisely when $z=(z_1,\ldots,z_N)\in\mathbb C^N$ belongs to the real algebraic manifold 
$$Y=\left\{z\in\mathbb C^N\Big|g^\circ_\alpha(z_1,\ldots,z_N)=h^\circ_\alpha(z_1,\ldots,z_N)=0,\forall\alpha\right\}$$
where $f_\alpha=g_\alpha+h_\alpha$ is the even/odd decomposition of $f_\alpha$.
\end{proposition}

\begin{proof}
Let $X_i$ be the matrices in the statement. In order for $x_i\to X_i$ to define a morphism of algebras, these matrices must satisfy the equations defining $X$. Thus, the space $Y$ in the statement consists of the points $z=(z_1,\ldots,z_N)\in\mathbb C^N$ satisfying:
$$f_\alpha(X_1,\ldots,X_N)=0\quad,\quad\forall\alpha$$

Now observe that the matrices $X_i$ in the statement multiply as follows:
$$X_{i_1}X_{j_1}\ldots X_{i_k}X_{j_k}=\begin{pmatrix}z_{i_1}\bar{z}_{j_1}\ldots z_{i_k}\bar{z}_{j_k}&0\\ 0&\bar{z}_{i_1}z_{j_1}\ldots\bar{z}_{i_k}z_{j_k}\end{pmatrix}$$
$$X_{i_1}X_{j_1}\ldots X_{i_k}X_{j_k}X_{i_{k+1}}=\begin{pmatrix}0&z_{i_1}\bar{z}_{j_1}\ldots z_{i_k}\bar{z}_{j_k}z_{i_{k+1}}\\ \bar{z}_{i_1}z_{j_1}\ldots\bar{z}_{i_k}z_{j_k}\bar{z}_{i_{k+1}}&0\end{pmatrix}$$

We therefore obtain, in terms of the even/odd decomposition $f_\alpha=g_\alpha+h_\alpha$:
$$f_\alpha(X_1,\ldots,X_N)=\begin{pmatrix}g^\circ_\alpha(z_1,\ldots,z_N)&h^\circ_\alpha(z_1,\ldots,z_N)\\ 
\\
\overline{h^\circ_\alpha(z_1,\ldots,z_N)}&\overline{g^\circ_\alpha(z_1,\ldots,z_N)}\end{pmatrix}$$

Thus, we obtain the equations for $Y$ from the statement.
\end{proof}

As a first consequence, of theoretical interest, a necessary condition for $X$ to exist is that the manifold $Y\subset\mathbb C^N$ constructed above must be compact, and we will be back to this later. In order to discuss now modeling questions, we will need as well:

\index{projective version}
\index{projective manifold}

\begin{definition}
Assuming that we are given a manifold $Z$, appearing via
$$C(Z)=C^*\left(z_1,\ldots,z_N\Big|f_\alpha(z_1,\ldots,z_N)=0\right)$$
we define the projective version of $Z$ to be the quotient space $Z\to PZ$ corresponding to the subalgebra $C(PZ)\subset C(Z)$ generated by the variables $x_{ij}=z_iz_j^*$.
\end{definition}

The relation with the half-classical manifolds comes from the fact that the projective version of a half-classical manifold is classical. Indeed, from $abc=cba$ we obtain:
\begin{eqnarray*}
ab\cdot cd
&=&(abc)d\\
&=&(cba)d\\
&=&c(bad)\\
&=&c(dab)\\
&=&cd\cdot ab
\end{eqnarray*}

Finally, let us call as before ``matrix model'' any morphism of unital $C^*$-algebras $f:A\to B$, with target algebra $B=M_K(C(Y))$, with $K\in\mathbb N$, and $Y$ being a compact space. With these conventions, following Bichon \cite{bic}, we have the following result:

\begin{theorem}
Given a half-classical manifold $X$ which is symmetric, in the sense that all its defining polynomials $f_\alpha$ are even, its universal $2\times2$ antidiagonal model,
$$\pi:C(X)\to M_2(C(Y))$$
where $Y$ is the manifold constructed in Proposition 12.31, is faithful. In addition, the construction $X\to Y$ is such that $X$ exists precisely when $Y$ is compact.
\end{theorem}

\begin{proof}
We can proceed as in \cite{bic}. Indeed, the universal model $\pi$ in the statement induces, at the level of projective versions, a certain representation:
$$C(PX)\to M_2(C(PY))$$

By using the multiplication formulae from the proof of Proposition 12.31, the image of this representation consists of diagonal matrices, and the upper left components of these matrices are the standard coordinates of $PY$. Thus, we have an isomorphism: 
$$PX\simeq PY$$

We can conclude then by using a grading trick. See \cite{bic}.
\end{proof}

As a first observation, this result shows that when $X$ is symmetric, we have $X^*\subset X^{(2)}$. Going beyond this observation is an interesting problem.

\bigskip

In what follows, we will rather need a more detailed version of the above result. For this purpose, we can use the following definition:

\begin{definition}
Associated to any compact manifold $Y\subset\mathbb C^N$ is the real compact half-classical manifold $[Y]$, having as coordinates the following variables,
$$X_i=\begin{pmatrix}0&z_i\\ \bar{z}_i&0\end{pmatrix}$$
where $z_1,\ldots,z_N$ are the standard coordinates on $Y$. In other words, $[Y]$ is given by the fact that $C([Y])\subset M_2(C(Y))$ is the algebra generated by these matrices.
\end{definition}

Here the fact that the manifold $[Y]$ is indeed half-classical follows from the results above. As for the fact that $[Y]$ is indeed algebraic, this follows from Theorem 12.33. Now with this notion in hand, we can reformulate Theorem 12.33, as follows:

\begin{theorem}
The symmetric half-classical manifolds $X$ appear as follows:
\begin{enumerate}
\item We have $X=[Y]$, for a certain conjugation-invariant subspace $Y\subset \mathbb C^N$.

\item $PX=P[Y]$, and $X$ is maximal with this property.

\item In addition, we have an embedding $C([X])\subset C(X)\rtimes\mathbb Z_2$.
\end{enumerate}
\end{theorem}

\begin{proof}
This follows from Theorem 12.33, with the embedding in (3) being constructed as in \cite{bic}, by $x_i=z_i\otimes\tau$, where $\tau$ is the standard generator of $\mathbb Z_2$. See \cite{bic}.
\end{proof}

And this is all, on this subject. In the unitary case things are a bit more complicated, and in connection with this, there are also some higher analogues of the above developed, using $K\times K$ matrix models. We refer to \cite{bb1}, \cite{bic}, \cite{bdu} for more on these topics.

\bigskip

As a conclusion now, and by getting back to the real case, for simplifying, the half-classical geometry can be normally developed in a quite efficient way, at a technical level which is close to that of the classical one, by using $2\times2$ matrix models, as indicated above. Of course, all this still remains to be done. In fact, as already mentioned in chapter 9 and afterwards, on several occasions, there are plenty of interesting things to be done here, and there is certainly room for writing a nice book on the subject.

\bigskip

Which reminds a bit the situation with the twisting, from chapter 11, with a nice book to be written there as well. In fact, both the half-classical geometry and the twisted geometry, and their combination the half-classical twisted geometry, which is something which exists as well, are somehow examples of ``tame geometries'', not far from the classical geometry, and with a bewildering array of techniques, including those of Connes \cite{co1}, potentially applying, and with very interesting results at stake.

\bigskip

And a word about physics, to finish with. Although there is nothing much concrete here, at least so far, a quite common belief is that, mathematically speaking somehow, QED is supposed to be something tame, and QCD is supposed to be something wild. And this is why we have mixed tame and wild things in this book, with tame and wild meaning for us something purely mathematical, namely non-free and free, with the belief that things are in correspondence, and that all this can be of help in physics.

\bigskip

We will be back to more speculations in chapters 13-16 below, when discussing more in detail free geometry, in continuation of the material from chapters 5-8, and benefiting too from what we learned from here, chapters 9-12, now coming to an end.

\section*{12e. Exercises} 

The matrix model problematics is quite exciting, making us exit the abstract algebra computations that we have been mainly doing throughout this book, and we have many exercises on the subject, for the most of research level. First, we have: 

\begin{exercise}
Try axiomatizing the quadruplets $(S,T,U,K)$ in terms of the associated von Neumann algebras, and report on what you found.
\end{exercise}

This is something that we discussed in the above, with the comment that this is a ``bad idea'', physically speaking. However, trying to have it done is certainly instructive. Plus hey, maybe I'm wrong with physics, and this is the way to go. Who knows.

\begin{exercise}
Prove that given a compact quantum group $G$, in order for having a faithful model $C(G)\subset M_K(C(T))$, the discrete dual $\Gamma=\widehat{G}$ must be amenable.
\end{exercise}

As a bonus exercise here, try to fully clarify the situation in the case where $\Gamma=\widehat{G}$ is assumed to be a classical discrete group. This is something quite tricky, and in case you do not find the answer, the keyword for a search is ``Thoma theorem''.

\begin{exercise}
Try to come up with a notion of inner faithfulness for the matrix models $C(X)\to M_K(C(T))$, in the case where $X\subset S^{N-1}_{\mathbb C,+}$ is an homogeneous space.
\end{exercise}

This is actually an open question, and any study on it would be very interesting.

\begin{exercise}
In the context of the matrix truncations, comment on what happens when $X_{class}=\emptyset$. Also, comment on the case $X^{(\infty)}=X$. And also, comment on the case where $X=G$ is assumed to be a compact quantum group.
\end{exercise}

Here, in what regards the first 2 questions, the answer normally requires some theory, examples, and counterexamples. As for the last question, the very first problem here is whether $G^{(2)}$ is a quantum group or not, the answer being no in general.

\begin{exercise}
Develop a matrix model theory for the spaces of quantum partial isometries and partial permutations from chapter $6$.
\end{exercise}

Needed here would be especially interesting examples. As a hint, try first finding some interesting models for the quantum permutation groups $S_N^+$, and then suitably modify your construction, as to make it work for the spaces of quantum partial permutations.

\part{Free coordinates}

\ \vskip50mm

\begin{center}
{\em And that seemed the end

But they caught him in vain

Cause a change came for Spain

And El Lute}
\end{center}

\chapter{Free coordinates}

\section*{13a. Easy geometries}

We discuss here and in the next 3 chapters a number of more specialized questions, of algebraic, geometric, analytic and probabilistic nature. We will be interested in the main 9 examples of noncommutative geometries in our sense, which are as follows:
$$\xymatrix@R=40pt@C=40pt{
\mathbb R^N_+\ar[r]&\mathbb T\mathbb R^N_+\ar[r]&\mathbb C^N_+\\
\mathbb R^N_*\ar[u]\ar[r]&\mathbb T\mathbb R^N_*\ar[u]\ar[r]&\mathbb C^N_*\ar[u]\\
\mathbb R^N\ar[u]\ar[r]&\mathbb T\mathbb R^N\ar[u]\ar[r]&\mathbb C^N\ar[u]
}$$

Our purpose will be that of going beyond the basic level, where we are now, with a number of results regarding the coordinates $x_1,\ldots,x_N$ of such spaces:

\smallskip

\begin{enumerate}
\item A first question, which is algebraic, is that of understanding the precise relations satisfied by these coordinates. We will see that this is related to the question of unifying the twisted and untwisted geometries, via intersection.

\smallskip

\item A second question, which is analytic, is that of understanding the fixed $N$ behavior of these coordinates. This can be done via deformation methods. We will see as well that there is an unexpected link with quantum permutations.
\end{enumerate}

\smallskip

Let us begin by discussing algebraic aspects. This is something quite fundamental. Indeed, in the classical case, the algebraic manifolds $X$ can be identified with the corresponding ideals of vanishing polynomials $J$, and the correspondence $X\leftrightarrow J$ is the foundation for all the known algebraic geometric theory, ancient or more modern.

\bigskip

In the free setting, things are in a quite primitive status, and a suitable theory of ``noncommutative algebra'', useful in connection with our present considerations, is so far missing. Computing $J$ for the free spheres, and perhaps for some other spheres as well, is a problem which is difficult enough for us, and that we will investigate here. 

\bigskip

As a starting point, we know that the above 9 geometries are easy, and looking in detail at this easiness property will be our first task. Let us first recall that we have:

\index{easy geometry}
\index{easiness}

\begin{definition}
A geometry $(S,T,U,K)$ is called easy when $U,K$ are easy, and 
$$U=\{O_N,K\}$$
with the operation on the right being the easy generation operation.
\end{definition}

To be more precise, in order for a geometry to be easy, the quantum groups $U,K$ must be of course easy, as stated above. Regarding now the generation condition, the point is that one of our general axioms for the nocommutative geometries, from chapter 4, states that we must have $U=<O_N,K>$, with the operation $<\,,>$ being a usual generation operation. And the above easy generation condition $U=\{O_N,K\}$ is something stronger, and so imposing this condition amounts in saying that we must have:
$$<O_N,K>=\{O_N,K\}$$

The easy geometries in the above sense can be investigated by using:

\begin{proposition}
An easy geometry is uniquely determined by a pair $(D,E)$ of categories of partitions, which must be as follows,
$$\mathcal{NC}_2\subset D\subset P_2$$
$$\mathcal{NC}_{even}\subset E\subset P_{even}$$
and which are subject to the following intersection and generation conditions,
$$D=E\cap P_2$$
$$E=<D,\mathcal{NC}_{even}>$$
and to the usual axioms for the associated quadruplet $(S,T,U,K)$, where $U,K$ are respectively the easy quantum groups associated to the categories $D,E$.
\end{proposition}

\begin{proof}
This statement simply comes from the following conditions: 
$$U=\{O_N,K\}$$
$$K=U\cap K_N^+$$

To be more precise, let us look at Definition 13.1. The main condition there tells us that $U,K$ must be easy, coming from certain categories $D,E$. It is clear that $D,E$ must appear as intermediate categories, as in the statement, and the fact that the intersection and generation conditions must be satisfied follows from:
\begin{eqnarray*}
U=\{O_N,K\}&\iff&D=E\cap P_2\\
K=U\cap K_N^+&\iff&E=<D,\mathcal{NC}_{even}>
\end{eqnarray*}

Thus, we are led to the conclusion in the statement.
\end{proof}

Generally speaking, the idea now is that, in the context of an easy geometry, everything can be reformulated in terms of the categories of partitions $(D,E)$, which must satisfy the conditions in Proposition 13.2. Thus, we have in fact a diagram as follows:
$$
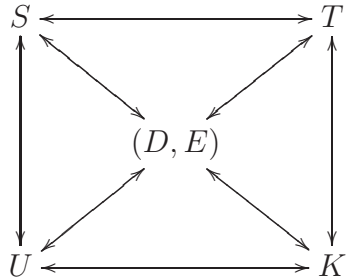
\xymatrix@R=30pt@C=30pt{
S\ar[rr]\ar[dr]\ar[dd]&&T\ar[ll]\ar[dd]\ar[dl]\\
&(D,E)\ar[ul]\ar[ur]\ar[dr]\ar[dl]\\
U\ar[uu]\ar[ur]\ar[rr]&&K\ar[ll]\ar[ul]\ar[uu]
}$$

This is not suprising, because our main examples of geometries are the classical ones, governed by the commutation relations $ab=ba$, then the half-classical ones, coming from the half-commutation relations $abc=cba$, and then the free geometries, coming from no relations at all. Thus, modulo some technical conditions and axioms involving the quadruplets $(S,T,U,K)$, which are there in order for our geometry to really ``work'', everything comes down to the combinatorial structure which replaces the commutation relations $ab=ba$. And the notion of category of partitions is precisely there for that.

\bigskip

This was for the idea. Now instead of discussing the full reformulation of our axions in terms of categories of partitions, which technically speaking will not bring many new things, let us work out at least the construction of the quadruplet $(S,T,U,K)$. In what regards the quantum groups, these come from via Tannakian duality, as follows:

\begin{theorem}
In the context of an easy geometry $(S,T,U,K)$, we have:
$$C(U)=C(U_N^+)\big/\left<T_\pi\in Hom(u^{\otimes k},u^{\otimes l})\Big|\forall k,l,\forall\pi\in D(k,l)\right>$$
Also, we have the following formula:
$$C(K)=C(K_N^+)\big/\left<T_\pi\in Hom(u^{\otimes k},u^{\otimes l})\Big|\forall k,l,\forall\pi\in D(k,l)\right>$$
In fact, these formulae simply follow from the fact that $U$ is easy.
\end{theorem}

\begin{proof}
This follows from general easiness considerations. Indeed, the construction of the easy quantum groups in \cite{bsp}, based on the Tannakian duality of Woronowicz from \cite{wo2}, in its soft form from Malacarne \cite{mal}, amounts in saying that the easy quantum group $G\subset U_N$ associated to a category of partitions $F=(F(k,l))$ is given by:
$$C(G)=C(U_N^+)\big/\left<T_\pi\in Hom(u^{\otimes k},u^{\otimes l})\Big|\forall k,l,\forall\pi\in F(k,l)\right>$$

Thus, for the categories of partitions $D,E$ associated to an easy geometry, as in Proposition 13.2, the corresponding quantum groups are as follows:
$$C(U)=C(U_N^+)\big/\left<T_\pi\in Hom(u^{\otimes k},u^{\otimes l})\Big|\forall k,l,\forall\pi\in D(k,l)\right>$$
$$C(K)=C(U_N^+)\big/\left<T_\pi\in Hom(u^{\otimes k},u^{\otimes l})\Big|\forall k,l,\forall\pi\in E(k,l)\right>$$

But the first formula is the formula for $U$ in the statement. As for the second formula, this can be fine-tuned by using the following formula, again coming from easiness:
$$C(K_N^+)=C(U_N^+)\big/\left<T_\pi\in Hom(u^{\otimes k},u^{\otimes l})\Big|\forall k,l,\forall\pi\in\mathcal{NC}_{even}(k,l)\right>$$

Indeed, by using the formula $E=<D,\mathcal{NC}_{even}>$ from Proposition 13.2, we have:
$$C(K)
=C(U_N^+)\big/\left<T_\pi\in Hom(u^{\otimes k},u^{\otimes l})\Big|\forall k,l,\forall\pi\in <D,\mathcal{NC}_{even}>(k,l)\right>$$

But constructing the algebra on the right amounts in dividing by the ideal coming from the partitions in $\mathcal{NC}_{even}$, which gives the algebra $C(K_N^+)$, and then further dividing by the ideal coming from the partitions in $D$, which gives the algebra in the statement.
\end{proof}

Regarding now the associated torus $T$, which is not exactly covered by the easy quantum group formalism, the result here is a bit different, as follows:

\index{easy torus}

\begin{theorem}
In the context of an easy geometry $(S,T,U,K)$, we have:
$$\Gamma=F_N\Big/\left<g_{i_1}\ldots g_{i_k}=g_{j_1}\ldots g_{j_l}\Big|\forall i,j,k,l,\exists\pi\in D(k,l),\delta_\pi\begin{pmatrix}i\\ j\end{pmatrix}\neq0\right>$$
In fact, this formula simply follows from the fact that $U$ is easy.
\end{theorem}

\begin{proof}
Let us denote by $g_i=u_{ii}$ the standard coordinates on the associated torus $T$, and consider the diagonal matrix formed by these coordinates:
$$g=\begin{pmatrix}g_1\\&\ddots\\&&g_N\end{pmatrix}$$

We have the following computation:
\begin{eqnarray*}
C(T)
&=&\left[C(U_N^+)\Big/\left<T_\pi\in Hom(u^{\otimes k},u^{\otimes l})\Big|\forall\pi\in D\right>\right]\Big/\left<u_{ij}=0\Big|\forall i\neq j\right>\\
&=&\left[C(U_N^+)\Big/\left<u_{ij}=0\Big|\forall i\neq j\right>\right]\Big/\left<T_\pi\in Hom(u^{\otimes k},u^{\otimes l})\Big|\forall\pi\in D\right>\\
&=&C^*(F_N)\Big/\left<T_\pi\in Hom(g^{\otimes k},g^{\otimes l})\Big|\forall\pi\in D\right>
\end{eqnarray*}

Now observe that, with $g=diag(g_1,\ldots,g_N)$ as before, we have:
$$T_\pi g^{\otimes k}(e_{i_1}\otimes\ldots\otimes e_{i_k})=\sum_{j_1\ldots j_l}\delta_\pi\begin{pmatrix}i_1&\ldots&i_k\\ j_1&\ldots&j_l\end{pmatrix}e_{j_1}\otimes\ldots\otimes e_{j_l}\cdot g_{i_1}\ldots g_{i_k}$$

On the other hand, we have as well:
$$g^{\otimes l}T_\pi(e_{i_1}\otimes\ldots\otimes e_{i_k})=\sum_{j_1\ldots j_l}\delta_\pi\begin{pmatrix}i_1&\ldots&i_k\\ j_1&\ldots&j_l\end{pmatrix}e_{j_1}\otimes\ldots\otimes e_{j_l}\cdot g_{j_1}\ldots g_{j_l}$$

Thus, the commutation relation $T_\pi\in Hom(g^{\otimes k},g^{\otimes l})$ reads:
\begin{eqnarray*}
&&\sum_{j_1\ldots j_l}\delta_\pi\begin{pmatrix}i_1&\ldots&i_k\\ j_1&\ldots&j_l\end{pmatrix}e_{j_1}\otimes\ldots\otimes e_{j_l}\cdot g_{i_1}\ldots g_{i_k}\\
&=&\sum_{j_1\ldots j_l}\delta_\pi\begin{pmatrix}i_1&\ldots&i_k\\ j_1&\ldots&j_l\end{pmatrix}e_{j_1}\otimes\ldots\otimes e_{j_l}\cdot g_{j_1}\ldots g_{j_l}
\end{eqnarray*}

Thus we obtain the formula in the statement, and the last assertion is clear.
\end{proof}

Finally, regarding the sphere $S$, which is not a quantum group, but rather an homogeneous space, here the result is a bit more complicated, as follows:

\index{easy sphere}

\begin{theorem}
In the context of an easy geometry $(S,T,U,K)$, we have
$$C(S)=C(S^{N-1}_{\mathbb C,+})\Big/\left<x_{i_1}\ldots x_{i_k}=x_{j_1}\ldots x_{j_k}\Big|\forall i,j,k,l,\exists\pi\in D(k)\cap I_k,\delta_\pi\begin{pmatrix}i\\ j\end{pmatrix}\neq0\right>$$
where the set on the right, $I_k\subset P_2(k,k)$, is the set of colored permutations.
\end{theorem}

\begin{proof}
This follows indeed from Theorem 13.3, by applying the construction $U\to S$, which amounts in taking the first row space.
\end{proof}

Summarizing, in the case of an easy geometry, we can reconstruct $S,T,U,K$ out of $(D,E)$, or simply out of $D$, as done above. It is possible to reformulate everything in terms of $(D,E)$, or just of $D$, by taking our axioms from chapter 4, and plugging in the formulae of $S,T,U,K$ in terms of $(D,E)$, or in terms of $D$, coming from the above results. However, this remains something theoretical, and we will not get into details here.

\section*{13b. Monomial spheres}

Let us discuss now an alternative take on these questions, following \cite{bme}, based on the notion of ``monomiality'', which applies to the spheres, which are not easy. Looking back at the definition of the spheres that we have, and at the precise relations between the coordinates, we are led into the following notion:

\index{monomial sphere}

\begin{definition}
A monomial sphere is a subset $S\subset S^{N-1}_{\mathbb C,+}$ obtained via relations
$$x_{i_1}^{e_1}\ldots x_{i_k}^{e_k}=x_{i_{\sigma(1)}}^{f_1}\ldots x_{i_{\sigma(k)}}^{f_k}\quad,\quad\forall (i_1,\ldots,i_k)\in\{1,\ldots,N\}^k$$
with $\sigma\in S_k$ being certain permutations, and with $e_r,f_r\in\{1,*\}$ being certain exponents.
\end{definition}

This definition is quite broad, and we have for instance as example the sphere $S^{N-1}_{\mathbb C,\times}$ coming from the relations $ab^*c=cb^*a$, corresponding to the following diagram: 
$$
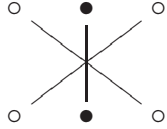
\xymatrix@R=10mm@C=5mm{\circ\ar@{-}[drr]&\bullet\ar@{-}[d]&\circ\ar@{-}[dll]\\\circ&\bullet&\circ}$$ 

This latter sphere is actually a quite interesting object, coming from the projective space considerations in \cite{bdd}, \cite{bd+}. However, while being monomial, this sphere does not exactly fit with our noncommutative geometry considerations here. 

\bigskip

To be more precise, according to the work in \cite{ba4}, \cite{bb2}, this sphere is part of a triple $(S^{N-1}_{\mathbb C,\times},\mathbb T_N^\times,U_N^\times)$, satisfying a simplified set of noncommutative geometry axioms. However, according to the work by Mang-Weber \cite{mwe}, the quantum group $U_N^\times$ has no reflection group counterpart $K_N^\times$. Thus, this sphere does not exactly fit with our axiomatics here.

\bigskip

In view of these difficulties, we will restrict now the attention to the real case. Let us first recall, from the various classification results established in chapter 6, that we have the following fundamental result, dealing with the real case:

\begin{theorem}
There are exactly $3$ real easy geometries, namely
$$\mathbb R^N\subset\mathbb R^N_*\subset\mathbb R^N_+$$
coming from the following categories of pairings $D$,
$$P_2\supset P_2^*\supset NC_2$$
whose associated spheres are as follows,
$$S^{N-1}_\mathbb R\subset S^{N-1}_{\mathbb R,*}\subset S^{N-1}_{\mathbb R,+}$$
and whose tori, unitary and reflection groups are given by similar formulae.
\end{theorem}

\begin{proof}
This is something that we know from chapter 6, coming from the fact that $G=O_N^*$ is the unique intermediate easy quantum group $O_N\subset G\subset O_N^+$.
\end{proof}

Let us focus now on the spheres, and try to better understand their ``easiness'' property, with results getting beyond what has been done above, in the general easy context. That is, our objects of interest in what follows will be the 3 real spheres, namely:
$$S^{N-1}_\mathbb R\subset S^{N-1}_{\mathbb R,*}\subset S^{N-1}_{\mathbb R,+}$$

Our purpose in what follows we will be that of proving that these spheres are the only monomial ones. Following \cite{bme}, in order to best talk about monomiality, in the present real case, it is convenient to introduce the following group:
$$S_\infty=\bigcup_{k\geq0}S_k$$

To be more precise, this group appears by definition as an inductive limit, with the inclusions $S_k\subset S_{k+1}$ that we use being given by: 
$$\sigma\in S_k\implies\sigma(k+1)=k+1$$

In terms of $S_\infty$, the definition of the monomial spheres reformulates as follows:

\index{monomial sphere}

\begin{proposition}
The monomial spheres are the algebraic manifolds $S\subset S^{N-1}_{\mathbb R,+}$ obtained via relations of type
$$x_{i_1}\ldots x_{i_k}=x_{i_{\sigma(1)}}\ldots x_{i_{\sigma(k)}},\ \forall (i_1,\ldots,i_k)\in\{1,\ldots,N\}^k$$
associated to certain elements $\sigma\in S_\infty$, where $k\in\mathbb N$ is such that $\sigma\in S_k$. 
\end{proposition}

\begin{proof}
We must prove that the relations $x_{i_1}\ldots x_{i_k}=x_{i_{\sigma(1)}}\ldots x_{i_{\sigma(k)}}$ are left unchanged when replacing $k\to k+1$. But this follows from $\sum_ix_i^2=1$, because:
\begin{eqnarray*}
&&x_{i_1}\ldots x_{i_k}x_{i_{k+1}}=x_{i_{\sigma(1)}}\ldots x_{i_{\sigma(k)}}x_{i_{k+1}}\\
&\implies&x_{i_1}\ldots x_{i_k}x_{i_{k+1}}^2=x_{i_{\sigma(1)}}\ldots x_{i_{\sigma(k)}}x_{i_{k+1}}^2\\
&\implies&\sum_{i_{k+1}}x_{i_1}\ldots x_{i_k}x_{i_{k+1}}^2=\sum_{i_{k+1}}x_{i_{\sigma(1)}}\ldots x_{i_{\sigma(k)}}x_{i_{k+1}}^2\\
&\implies&x_{i_1}\ldots x_{i_k}=x_{i_{\sigma(1)}}\ldots x_{i_{\sigma(k)}}
\end{eqnarray*}

Thus we can indeed ``simplify at right'', and this gives the result.
\end{proof}

As already mentioned, following \cite{bme}, our goal in what follows will be that of proving that the 3 main spheres are the only monomial ones. In order to prove this result, we will use group theory methods. We call a subgroup $G\subset S_\infty$ filtered when it is stable under concatenation, in the sense that when writing $G=(G_k)$ with $G_k\subset S_k$, we have:
$$\sigma\in G_k,\pi\in G_l\implies \sigma\pi\in G_{k+l}$$

With this convention, we have the following result:

\index{filtered group}
\index{infinite permutation group}

\begin{theorem}
The monomial spheres are the subsets $S_G\subset S^{N-1}_{\mathbb R,+}$ given by
$$C(S_G)=C(S^{N-1}_{\mathbb R,+})\Big/\Big<x_{i_1}\ldots x_{i_k}=x_{i_{\sigma(1)}}\ldots x_{i_{\sigma(k)}},\forall (i_1,\ldots,i_k)\in\{1,\ldots,N\}^k,\forall\sigma\in G_k\Big>$$
where $G=(G_k)$ is a filtered subgroup of $S_\infty=(S_k)$.
\end{theorem}

\begin{proof}
We know from Proposition 13.8 that the construction in the statement produces a monomial sphere. Conversely, given a monomial sphere $S\subset S^{N-1}_{\mathbb R,+}$, let us set:
$$G_k=\left\{\sigma\in S_k\Big|x_{i_1}\ldots x_{i_k}=x_{i_{\sigma(1)}}\ldots x_{i_{\sigma(k)}},\forall (i_1,\ldots,i_k)\in\{1,\ldots,N\}^k\right\}$$

With $G=(G_k)$ we have then $S=S_G$. Thus, it remains to prove that $G$ is a filtered group. But since the relations $x_{i_1}\ldots x_{i_k}=x_{i_{\sigma(1)}}\ldots x_{i_{\sigma(k)}}$ can be composed and reversed, each $G_k$ follows to be stable under composition and inversion, and is therefore a group. Also, since the relations $x_{i_1}\ldots x_{i_k}=x_{i_{\sigma(1)}}\ldots x_{i_{\sigma(k)}}$ can be concatenated as well, our group $G=(G_k)$ is stable under concatenation, and we are done. 
\end{proof}

At the level of examples, according to our definitions, the simplest filtered groups, namely $\{1\}\subset S_\infty$, produce the simplest real spheres, namely:
$$S^{N-1}_{\mathbb R,+}\supset S^{N-1}_\mathbb R$$

In order to discuss now the half-classical case, we need to introduce and study a certain privileged intermediate filtered group $\{1\}\subset S_\infty^*\subset S_\infty$, which will eventually produce the intermediate sphere $S^{N-1}_{\mathbb R,+}\supset S^{N-1}_{\mathbb R,*}\supset S^{N-1}_\mathbb R$. This can be done as follows:

\begin{proposition}
Let $S_\infty^*\subset S_\infty$ be the set of permutations having the property that when labelling cyclically the legs as follows
$$\bullet\circ\bullet\circ\ldots$$
each string joins a black leg to a white leg.
\begin{enumerate}
\item $S_\infty^*$ is a filtered subgroup of $S_\infty$, generated by the half-classical crossing.

\item We have $S_{2k}^*\simeq S_k\times S_k$, and $S^*_{2k+1}\simeq S_k\times S_{k+1}$, for any $k\in\mathbb N$.
\end{enumerate}
\end{proposition}

\begin{proof}
The fact that $S_\infty^*$ is indeed a subgroup of $S_\infty$, which is filtered, is clear. Observe now that the half-classical crossing has the ``black-to-white'' joining property:
$$
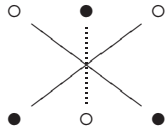
\xymatrix@R=10mm@C=5mm{\circ\ar@{-}[drr]&\bullet\ar@{.}[d]&\circ\ar@{-}[dll]\\\bullet&\circ&\bullet}$$ 

Thus this crossing belongs to $S_3^*$, and it is routine to check that the filtered subgroup of $S_\infty$ generated by it is the whole $S_\infty^*$. Regarding now the last assertion, observe first that the filtered subgroups $S_3^*,S_4^*$ consist of the following permutations:
$$
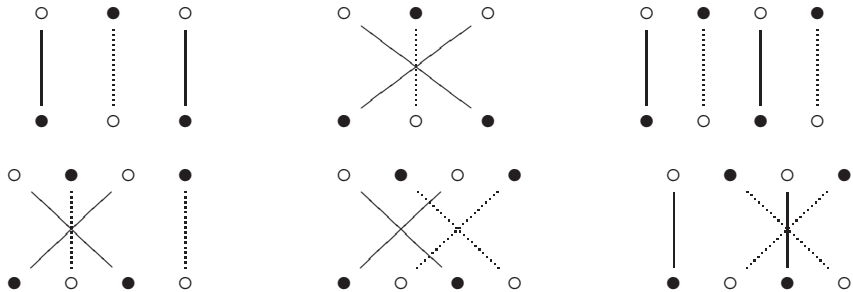
\xymatrix@R=10mm@C=5mm{\circ\ar@{-}[d]&\bullet\ar@{.}[d]&\circ\ar@{-}[d]\\\bullet&\circ&\bullet}\qquad\qquad
\xymatrix@R=10mm@C=5mm{\circ\ar@{-}[drr]&\bullet\ar@{.}[d]&\circ\ar@{-}[dll]\\\bullet&\circ&\bullet}\qquad\qquad 
\xymatrix@R=10mm@C=3mm{\circ\ar@{-}[d]&\bullet\ar@{.}[d]&\circ\ar@{-}[d]&\bullet\ar@{.}[d]\\\bullet&\circ&\bullet&\circ}$$
$$\xymatrix@R=10mm@C=3mm{\circ\ar@{-}[drr]&\bullet\ar@{.}[d]&\circ\ar@{-}[dll]&\bullet\ar@{.}[d]\\\bullet&\circ&\bullet&\circ}\qquad\qquad 
\xymatrix@R=10mm@C=3mm{\circ\ar@{-}[drr]&\bullet\ar@{.}[drr]&\circ\ar@{-}[dll]&\bullet\ar@{.}[dll]\\\bullet&\circ&\bullet&\circ}\qquad\qquad
\xymatrix@R=10mm@C=3mm{\circ\ar@{-}[d]&\bullet\ar@{.}[drr]&\circ\ar@{-}[d]&\bullet\ar@{.}[dll]\\\bullet&\circ&\bullet&\circ}$$

Thus we have $S_3^*=S_1\times S_2$ and $S_4^*=S_2\times S_2$, with the first component coming from dotted permutations, and with the second component coming from the solid line permutations. The same argument works in general, and gives the last assertion.
\end{proof}

Now back to the main 3 real spheres, the result is as follows:

\begin{proposition}
The basic monomial real spheres, namely 
$$S^{N-1}_\mathbb R\subset S^{N-1}_{\mathbb R,*}\subset S^{N-1}_{\mathbb R,+}$$
come respectively from the filtered groups $S_\infty\supset S_\infty^*\supset\{1\}$.
\end{proposition}

\begin{proof}
This is clear by definition in the classical and in the free cases. In the half-liberated case, the result follows from Proposition 13.10 (1).
\end{proof}

Now back to the general case, with the idea in mind of proving the uniqueness of the above spheres, consider a monomial sphere $S_G\subset S^{N-1}_{\mathbb R,+}$, with the filtered group $G\subset S_\infty$ taken to be maximal, as in the proof of Theorem 13.9. We have the following result:

\begin{proposition}
The filtered group $G\subset S_\infty$ associated to a monomial sphere $S\subset S^{N-1}_{\mathbb R,+}$ is stable under the following operations, on the corresponding diagrams:
\begin{enumerate}
\item Removing outer strings.

\item Removing neighboring strings.
\end{enumerate}
\end{proposition}

\begin{proof}
Both these results follow by using the quadratic condition:

\medskip

(1) Regarding the outer strings, by summing over $a$, we have:
\begin{eqnarray*}
Xa=Ya
&\implies&Xa^2=Ya^2\\
&\implies&X=Y
\end{eqnarray*}

We have as well the following computation:
\begin{eqnarray*}
aX=aY
&\implies&a^2X=a^2Y\\
&\implies&X=Y
\end{eqnarray*}

(2) Regarding the neighboring strings, once again by summing over $a$, we have:
\begin{eqnarray*}
XabY=ZabT
&\implies&Xa^2Y=Za^2T\\
&\implies&XY=ZT
\end{eqnarray*}

We have as well the following computation:
\begin{eqnarray*}
XabY=ZbaT
&\implies&Xa^2Y=Za^2T\\
&\implies&XY=ZT
\end{eqnarray*}

Thus $G=(G_k)$ has both the properties in the statement.
\end{proof}

We can now state and prove a main result, from \cite{bme}, as follows:

\begin{theorem}
There is only one intermediate monomial sphere
$$S^{N-1}_\mathbb R\subset S\subset S^{N-1}_{\mathbb R,+}$$
namely the half-classical real sphere $S^{N-1}_{\mathbb R,*}$.
\end{theorem}

\begin{proof}
We will prove that the only filtered groups $G\subset S_\infty$ satisfying the conditions in Proposition 13.12 are those correspoding to our 3 spheres, namely:
$$\{1\}\subset S_\infty^*\subset S_\infty$$

In order to do so, consider such a filtered group $G\subset S_\infty$. We assume this group to be non-trivial, $G\neq\{1\}$, and we want to prove that we have $G=S_\infty^*$ or $G=S_\infty$.

\medskip

\underline{Step 1}. Our first claim is that $G$ contains a 3-cycle. Assume indeed that two permutations $\pi,\sigma\in S_\infty$ have support overlapping on exactly one point, say:
$$supp(\pi)\cap supp(\sigma)=\{i\}$$

The point is then that the commutator $\sigma^{-1}\pi^{-1}\sigma\pi$ is a 3-cycle, namely:
$$(i,\sigma^{-1}(i),\pi^{-1}(i))$$

Indeed the computation of the commutator goes as follows:
$$
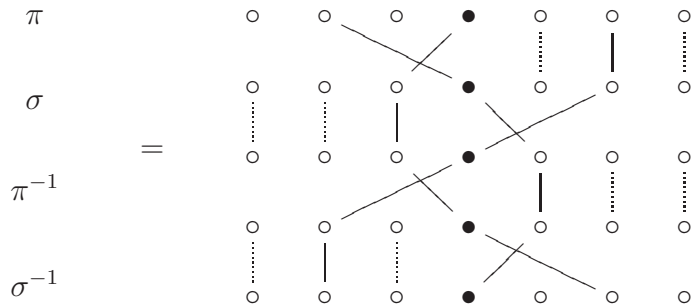
\xymatrix@R=7mm@C=5mm{\pi\\ \sigma\\ \pi^{-1}\\ \sigma^{-1}}\qquad
\xymatrix@R=6mm@C=5mm{\\ \\ =}\qquad
\xymatrix@R=5mm@C=5mm{
\circ&\circ\ar@{-}[drr]&\circ&\bullet\ar@{-}[dl]&\circ\ar@{.}[d]&\circ\ar@{-}[d]&\circ\ar@{.}[d]\\
\circ\ar@{.}[d]&\circ\ar@{.}[d]&\circ\ar@{-}[d]&\bullet\ar@{-}[dr]&\circ&\circ\ar@{-}[dll]&\circ\\
\circ&\circ&\circ\ar@{-}[dr]&\bullet\ar@{-}[dll]&\circ\ar@{-}[d]&\circ\ar@{.}[d]&\circ\ar@{.}[d]\\
\circ\ar@{.}[d]&\circ\ar@{-}[d]&\circ\ar@{.}[d]&\bullet\ar@{-}[drr]&\circ\ar@{-}[dl]&\circ&\circ\\
\circ&\circ&\circ&\bullet&\circ&\circ&\circ
}$$

Now let us pick a non-trivial element $\tau\in G$. By removing outer strings at right and at left we obtain permutations $\tau'\in G_k,\tau''\in G_s$ having a non-trivial action on their right/left leg, and the trick applies, with:
$$\pi=\tau'\otimes id_{s-1}\quad,\quad 
\sigma=id_{k-1}\otimes\tau''$$

Thus, $G$ contains a 3-cycle, as claimed.

\medskip

\underline{Step 2}. Our second claim is $G$ must contain one of the following permutations:
$$
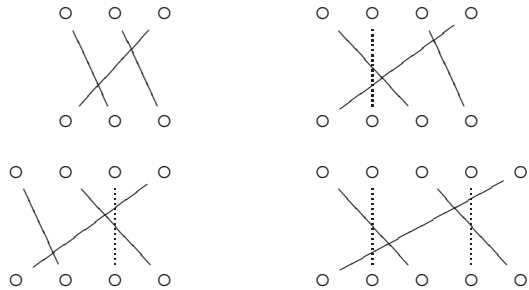
\xymatrix@R=10mm@C=2mm{
\circ\ar@{-}[dr]&\circ\ar@{-}[dr]&\circ\ar@{-}[dll]\\
\circ&\circ&\circ}\qquad\qquad
\xymatrix@R=10mm@C=2mm{
\circ\ar@{-}[drr]&\circ\ar@{.}[d]&\circ\ar@{-}[dr]&\circ\ar@{-}[dlll]\\
\circ&\circ&\circ&\circ}$$
$$\xymatrix@R=10mm@C=2mm{
\circ\ar@{-}[dr]&\circ\ar@{-}[drr]&\circ\ar@{.}[d]&\circ\ar@{-}[dlll]\\
\circ&\circ&\circ&\circ}\qquad\qquad
\xymatrix@R=10mm@C=2mm{
\circ\ar@{-}[drr]&\circ\ar@{.}[d]&\circ\ar@{-}[drr]&\circ\ar@{.}[d]&\circ\ar@{-}[dllll]\\
\circ&\circ&\circ&\circ&\circ}$$

Indeed, consider the 3-cycle that we just constructed. By removing all outer strings, and then all pairs of adjacent vertical strings, we are left with these permutations.

\medskip

\underline{Step 3}. Our claim now is that we must have $S_\infty^*\subset G$. Indeed, let us pick one of the permutations that we just constructed, and apply to it our various diagrammatic rules. From the first permutation we can obtain the basic crossing, as follows:
$$
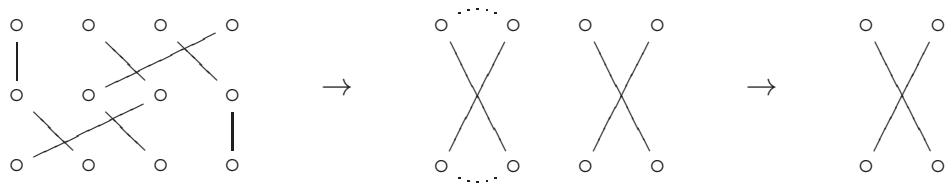
\xymatrix@R=5mm@C=5mm{
\circ\ar@{-}[d]&\circ\ar@{-}[dr]&\circ\ar@{-}[dr]&\circ\ar@{-}[dll]\\
\circ\ar@{-}[dr]&\circ\ar@{-}[dr]&\circ\ar@{-}[dll]&\circ\ar@{-}[d]\\
\circ&\circ&\circ&\circ}
\qquad
\xymatrix@R=5mm@C=5mm{
\\ \to\\}\qquad
\xymatrix@R=6mm@C=5mm{
\circ\ar@{-}[ddr]\ar@/^/@{.}[r]&\circ\ar@{-}[ddl]&\circ\ar@{-}[ddr]&\circ\ar@{-}[ddl]\\
\\
\circ\ar@/_/@{.}[r]&\circ&\circ&\circ}
\qquad
\xymatrix@R=5mm@C=5mm{
\\ \to\\}\qquad
\xymatrix@R=6mm@C=5mm{
\circ\ar@{-}[ddr]&\circ\ar@{-}[ddl]\\
\\
\circ&\circ}$$

Also, by removing a suitable $\slash\hskip-2.1mm\backslash$ shaped configuration, which is represented by dotted lines in the diagrams below, we can obtain the basic crossing from the second and third permutation, and the half-liberated crossing from the fourth permutation:
$$
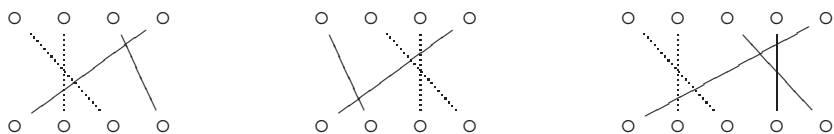
\xymatrix@R=10mm@C=2mm{
\circ\ar@{.}[drr]&\circ\ar@{.}[d]&\circ\ar@{-}[dr]&\circ\ar@{-}[dlll]\\
\circ&\circ&\circ&\circ}\qquad\qquad
\xymatrix@R=10mm@C=2mm{
\circ\ar@{-}[dr]&\circ\ar@{.}[drr]&\circ\ar@{.}[d]&\circ\ar@{-}[dlll]\\
\circ&\circ&\circ&\circ}\qquad\qquad
\xymatrix@R=10mm@C=2mm{
\circ\ar@{.}[drr]&\circ\ar@{.}[d]&\circ\ar@{-}[drr]&\circ\ar@{-}[d]&\circ\ar@{-}[dllll]\\
\circ&\circ&\circ&\circ&\circ}$$

Thus, in all cases we have a basic or half-liberated crossing, and so, as desired: 
$$S_\infty^*\subset G$$

\underline{Step 4}. Our last claim, which will finish the proof, is that there is no proper intermediate subgroup as follows:
$$S_\infty^*\subset G\subset S_\infty$$

In order to prove this, observe that $S_\infty^*\subset S_\infty$ is the subgroup of parity-preserving permutations, in the sense that ``$i$ even $\implies$ $\sigma(i)$ even''. 

\medskip

Now let us pick an element $\sigma\in S_k-S_k^*$, with $k\in\mathbb N$. We must prove that the group $G=<S_\infty^*,\sigma>$ equals the whole $S_\infty$. In order to do so, we use the fact that $\sigma$ is not parity preserving. Thus, we can find $i$ even such that $\sigma(i)$ is odd. In addition, up to passing to $\sigma|$, we can assume that $\sigma(k)=k$, and then, up to passing one more time to $\sigma|$, we can further assume that $k$ is even. Since both $i,k$ are even we have:
$$(i,k)\in S_k^*$$

We conclude that the following element belongs to $G$:
$$\sigma(i,k)\sigma^{-1}=(\sigma(i),k)$$

But, since $\sigma(i)$ is odd, by deleting an appropriate number of vertical strings, $(\sigma(i),k)$ reduces to the basic crossing $(1,2)$. Thus $G=S_\infty$, and we are done.
\end{proof}

As already mentioned in the above, the story is not over with this kind of result, because the complex case still remains to be worked out.

\section*{13c. Twists, intersections}

Our purpose now will be that of going beyond the above results, with a number of more specialized results regarding the coordinates $x_1,\ldots,x_N$ of our real spheres. To be more precise, a first question that we would like to solve, which is of purely algebraic nature, is that of understanding the precise relations satisfied by these coordinates $x_1,\ldots,x_N$ over our real spheres. We will see, in a somewhat unexpected way, that this is related to the question of unifying the twisted and untwisted geometries, via intersection.

\bigskip

Let us begin by recalling the construction of the twisted real spheres, which was discussed in chapter 11. This is something very simple, as follows:

\index{twisting}

\begin{definition}
The subspheres $\bar{S}^{N-1}_\mathbb R,\bar{S}^{N-1}_{\mathbb R,*}\subset S^{N-1}_{\mathbb R,+}$ are constructed by imposing the following conditions on the standard coordinates $x_1,\ldots,x_N$:
\begin{enumerate}
\item $\bar{S}^{N-1}_\mathbb R$: $x_ix_j=-x_jx_i$, for any $i\neq j$.

\item $\bar{S}^{N-1}_{\mathbb R,*}$: $x_ix_jx_k=-x_kx_jx_i$ for any $i,j,k$ distinct, $x_ix_jx_k=x_kx_jx_i$ otherwise.
\end{enumerate}
\end{definition}

Here the fact that we have indeed $\bar{S}^{N-1}_\mathbb R\subset\bar{S}^{N-1}_{\mathbb R,*}$ comes from the following computations, for $a,b,c\in\{x_i\}$ distinct, where $x_1,\ldots,x_N$ are the standard coordinates on $\bar{S}^{N-1}_\mathbb R$:
$$abc=-bac=bca=-cba$$
$$aab=-aba=baa$$

Summarizing, we have a total of 5 real spheres, or rather a total of $3+3=6$ real spheres, with the convention that the free real sphere equals its twist:
$$S^{N-1}_{\mathbb R,+}=\bar{S}^{N-1}_{\mathbb R,+}$$

The point now is that we can intersect these $3+3=6$ spheres, and we end up with a total of $3\times3=9$ real spheres, in a generalized sense, as follows: 

\begin{definition}
Associated to any integer $N\in\mathbb N$ are the generalized spheres
$$\xymatrix@R=13mm@C=3mm{
S^{N-1}_\mathbb R\ar[rr]&&S^{N-1}_{\mathbb R,*}\ar[rrr]&&&S^{N-1}_{\mathbb R,+}\\
S^{N-1}_\mathbb R\cap\bar{S}^{N-1}_{\mathbb R,*}\ar[rr]\ar[u]&&S^{N-1}_{\mathbb R,*}\cap\bar{S}^{N-1}_{\mathbb R,*}\ar[rrr]\ar[u]&&&\bar{S}^{N-1}_{\mathbb R,*}\ar[u]\\
S^{N-1}_\mathbb R\cap\bar{S}^{N-1}_\mathbb R\ar[rr]\ar[u]&&S^{N-1}_{\mathbb R,*}\cap\bar{S}^{N-1}_\mathbb R\ar[rrr]\ar[u]&&&\bar{S}^{N-1}_\mathbb R\ar[u]}$$
obtained by intersecting the $3$ twisted real spheres and the $3$ untwisted real spheres.
\end{definition}

In order to compute the various intersections appearing above, which in general cannot be thought of as being smooth, let us introduce the following objects:

\index{polygonal sphere}

\begin{definition}
The polygonal spheres are real algebraic manifolds, defined as
$$S^{N-1,d-1}_\mathbb R=\left\{x\in S^{N-1}_\mathbb R\Big|x_{i_0}\ldots x_{i_d}=0,\forall i_0,\ldots,i_d\ {\rm distinct}\right\}$$
depending on integers $1\leq d\leq N$.
\end{definition}

These spheres, introduced and studied in \cite{ba2}, are not smooth in general, but recall that we are currently doing algebraic geometry, rather than differential geometry, and with actually the colorful name ``polygonal spheres'', used in \cite{ba2} and that we will use here too, being there for reminding us that. To be more precise, the point is that the problem that we want to solve, namely understanding the precise relations satisfied by the coordinates $x_1,\ldots,x_N$ for the real spheres, naturally leads into polygonal spheres.

\bigskip

More generally now, we have the following construction of ``generalized polygonal spheres'', which applies to the half-classical and twisted cases too:
$$C(\dot{S}^{N-1,d-1}_{\mathbb R,\times})=C\big(\dot{S}^{N-1}_{\mathbb R,\times}\big)\Big/\Big<x_{i_0}\ldots x_{i_d}=0,\forall i_0,\ldots,i_d\ {\rm distinct}\Big>$$

Here the fact that in the classical case we obtain the polygonal spheres from Definition 13.16 comes from a straightforward application of the Gelfand theorem.

\bigskip

With these conventions, we have the following result, dealing with all the spheres that we have so far in real case, namely twisted, untwisted and intersections:

\index{intersections}

\begin{theorem}
The diagram obtained by intersecting the twisted and untwisted real spheres, from Definition 13.15, is given by
$$\xymatrix@R=13mm@C=16mm{
S^{N-1}_\mathbb R\ar[r]&S^{N-1}_{\mathbb R,*}\ar[r]&S^{N-1}_{\mathbb R,+}\\
S^{N-1,1}_\mathbb R\ar[r]\ar[u]&S^{N-1,1}_{\mathbb R,*}\ar[r]\ar[u]&\bar{S}^{N-1}_{\mathbb R,*}\ar[u]\\
S^{N-1,0}_\mathbb R\ar[r]\ar[u]&\bar{S}^{N-1,1}_\mathbb R\ar[r]\ar[u]&\bar{S}^{N-1}_\mathbb R\ar[u]}$$
and so all these spheres are generalized polygonal spheres.
\end{theorem}

\begin{proof}
Consider the 4-diagram obtained by intersecting the 5 main spheres:
$$\xymatrix@R=13mm@C=13mm{
S^{N-1}_\mathbb R\cap\bar{S}^{N-1}_{\mathbb R,*}\ar[r]&S^{N-1}_{\mathbb R,*}\cap\bar{S}^{N-1}_{\mathbb R,*}\\
S^{N-1}_\mathbb R\cap\bar{S}^{N-1}_\mathbb R\ar[r]\ar[u]&S^{N-1}_{\mathbb R,*}\cap\bar{S}^{N-1}_\mathbb R\ar[u]}$$

We must prove that this diagram coincides with the 4-diagram appearing at bottom left in the statement, which is as follows:
$$\xymatrix@R=13mm@C=13mm{
S^{N-1,1}_\mathbb R\ar[r]&S^{N-1,1}_{\mathbb R,*}\\
S^{N-1,0}_\mathbb R\ar[r]\ar[u]&\bar{S}^{N-1,1}_\mathbb R\ar[u]}$$

But this is clear, because combining the commutation and anticommutation relations leads to the vanishing relations defining the spheres of type $\dot{S}^{N-1,d-1}_{\mathbb R,\times}$. More precisely:

\medskip

(1) $S^{N-1}_\mathbb R\cap\bar{S}^{N-1}_\mathbb R$ consists of the points $x\in S^{N-1}_\mathbb R$ such that, for any $i\neq j$:
$$x_ix_j=-x_jx_i$$

Now since we have as well $x_ix_j=x_jx_i$, for any $i,j$, this relation reads $x_ix_j=0$ for $i\neq j$, which means that we have $x\in S^{N-1,0}_\mathbb R$, as desired.

\medskip

(2) $S^{N-1}_\mathbb R\cap\bar{S}^{N-1}_{\mathbb R,*}$ consists of the points $x\in S^{N-1}_\mathbb R$ such that, for $i,j,k$ distinct:
$$x_ix_jx_k=-x_kx_jx_i$$

Once again by commutativity, this relation is equivalent to $x\in S^{N-1,1}_\mathbb R$, as desired.

\medskip

(3) $S^{N-1}_{\mathbb R,*}\cap\bar{S}^{N-1}_\mathbb R$ is obtained from $\bar{S}^{N-1}_\mathbb R$ by imposing to the standard coordinates the half-commutation relations $abc=cba$. On the other hand, we know from $\bar{S}^{N-1}_\mathbb R\subset \bar{S}^{N-1}_{\mathbb R,*}$ that the standard coordinates on $\bar{S}^{N-1}_\mathbb R$ satisfy $abc=-cba$ for $a,b,c$ distinct, and $abc=cba$ otherwise. Thus, the relations brought by intersecting with $S^{N-1}_{\mathbb R,*}$ reduce to the relations $abc=0$ for $a,b,c$ distinct, and so we are led to the sphere $\bar{S}^{N-1,1}_\mathbb R$.

\medskip

(4) $S^{N-1}_{\mathbb R,*}\cap\bar{S}^{N-1}_{\mathbb R,*}$ is obtained from $\bar{S}^{N-1}_{\mathbb R,*}$ by imposing the relations $abc=-cba$ for $a,b,c$ distinct, and $abc=cba$ otherwise. Since we know that $abc=cba$ for any $a,b,c$, the extra relations reduce to $abc=0$ for $a,b,c$ distinct, and so we are led to $S^{N-1,1}_{\mathbb R,*}$.
\end{proof}

Summarizing, whether we want it or not, when talking about intersections between twisted and untwisted geometries, we are led into polygonal spheres, and into non-smooth objects in general. In view of this, and also in connection with general axiomatization questions, let us find now a suitable axiomatic framework for the 9 spheres in Theorem 13.17. We have the following definition, once again from \cite{ba2}, which is based on the signature function $\varepsilon:P_{even}\to\{\pm1\}$ constructed in chapter 11:

\index{twisted relations}
\index{untwisted relations}

\begin{definition}
Given variables $x_1,\ldots,x_N$, any permutation $\sigma\in S_k$ produces two collections of relations between these variables, as follows:
\begin{enumerate}
\item Untwisted relations, namely, for any $i_1,\ldots,i_k$:
$$x_{i_1}\ldots x_{i_k}=x_{i_{\sigma(1)}}\ldots x_{i_{\sigma(k)}}$$

\item Twisted relations, namely, for any $i_1,\ldots,i_k$:
$$x_{i_1}\ldots x_{i_k}=\varepsilon\left(\ker\begin{pmatrix}i_1&\ldots&i_k\\ i_{\sigma(1)}&\ldots&i_{\sigma(k)}\end{pmatrix}\right)x_{i_{\sigma(1)}}\ldots x_{i_{\sigma(k)}}$$
\end{enumerate}
The untwisted relations are denoted $\mathcal R_\sigma$, and the twisted ones are denoted $\bar{\mathcal R}_\sigma$.
\end{definition}

Observe that the untwisted relations $\mathcal R_\sigma$ are trivially satisfied for the standard coordinates on $S^{N-1}_\mathbb R$, for any permutation $\sigma\in S_k$. A twisted analogue of this fact holds, in the sense that the standard coordinates on $\bar{S}^{N-1}_\mathbb R$ satisfy the relations $\bar{\mathcal R}_\sigma$, for any $\sigma\in S_k$. Indeed, by using the anticommutation relations between the distinct coordinates of these latter spheres, we must have a formula of the following type:
$$x_{i_1}\ldots x_{i_k}=\pm x_{i_{\sigma(1)}}\ldots x_{i_{\sigma(k)}}$$

But the sign $\pm$ obtained in this way is precisely the one given above, namely:
$$\pm=\varepsilon\left(\ker\begin{pmatrix}i_1&\ldots&i_k\\ i_{\sigma(1)}&\ldots&i_{\sigma(k)}\end{pmatrix}\right)$$

We have now all the needed ingredients for axiomatizing the various spheres appearing so far, namely the twisted and untwisted ones, and their intersections:

\begin{definition}
We have $3$ types of quantum spheres $S\subset S^{N-1}_{\mathbb R,+}$, as follows,
\begin{enumerate}
\item Monomial, namely $\dot{S}^{N-1}_{\mathbb R,E}$, with $E\subset S_\infty$, obtained via the following relations,
$$\left\{\dot{\mathcal R}_\sigma\Big|\sigma\in E\right\}$$

\item Mixed monomial, which appear as intersections as follows, with $E,F\subset S_\infty$,
$$S^{N-1}_{\mathbb R,E,F}=S^{N-1}_{\mathbb R,E}\cap\bar{S}^{N-1}_{\mathbb R,F}$$

\item Polygonal, which are also intersections, with $E,F\subset S_\infty$, and $d\in\{1,\ldots,N\}$,
$$S^{N-1,d-1}_{\mathbb R,E,F}=S^{N-1}_{\mathbb R,E,F}\cap S^{N-1,d-1}_{\mathbb R,+}$$
\end{enumerate}
and which cover all the quantum spheres that we constructed so far.
\end{definition}

Observe that, as stated, we cover all the spheres appearing so far. More precisely, the 5 basic spheres in are monomial, the 9 spheres in Theorem 13.17 are mixed monomial, and the polygonal sphere formalism covers all the examples constructed so far.

\bigskip

Observe that the set of mixed monomial spheres is closed under intersections. The same holds for the set of polygonal spheres, because we have the following formula:
$$S^{N-1,d-1}_{\mathbb R,E,F}\cap S^{N-1,d'-1}_{\mathbb R,E',F'}=S^{N-1,min(d,d')-1}_{\mathbb R,E\cup E',F\cup F'}$$

Let us try now to understand the structure of the various types of spheres, by using the real sphere technology developed before. We call a group of permutations $G\subset S_\infty$ filtered if, with $G_k=G\cap S_k$, we have $G_k\times G_l\subset G_{k+l}$, for any $k,l$. We have:

\begin{proposition}
The various spheres can be parametrized by groups, as follows:
\begin{enumerate}
\item Monomial case: $\dot{S}^{N-1}_{\mathbb R,G}$, with $G\subset S_\infty$ filtered group.

\item Mixed monomial case: $S^{N-1}_{\mathbb R,G,H}$, with $G,H\subset S_\infty$ filtered groups.

\item Polygonal case: $S^{N-1,d-1}_{\mathbb R,G,H}$, with $G,H\subset S_\infty$ filtered groups, and $d\in\{1,\ldots,N\}$.
\end{enumerate}
\end{proposition}

\begin{proof}
This basically follows from the theory developed before, as follows:

\medskip

(1) As explained before, in order to prove this assertion, for a monomial sphere $S=\dot{S}_{\mathbb R,E}$, we can take $G\subset S_\infty$ to be the set of permutations $\sigma\in S_\infty$ having the property that the relations $\dot{\mathcal R}_\sigma$ hold for the standard coordinates of $S$. We have then $E\subset G$, we have as well $S=\dot{S}^{N-1}_{\mathbb R,G}$, and the fact that $G$ is a filtered group is clear as well.

\medskip

(2) This follows from (1), by taking intersections.

\medskip

(3) Once again this follows from (1), by taking intersections.
\end{proof}

The idea in what follows will be that of writing the 9 main polygonal spheres as in Proposition 13.20 (2), as to reach to a ``standard parametrization'' for our spheres. We recall that the permutations $\sigma\in S_\infty$ having the property that when labelling clockwise their legs $\circ\bullet\circ\bullet\ldots$, and string joins a white leg to a black leg, form a filtered group, denoted $S_\infty^*\subset S_\infty$. This group comes from the general half-liberation considerations from chapter 9, and its algebraic structure is very simple, as follows:
$$S_{2n}^*\simeq S_n\times S_n\quad,\quad 
S_{2n+1}^*\simeq S_n\times S_{n+1}$$

Let us formulate as well the following definition:

\index{standard parametrization}

\begin{definition}
We call a mixed monomial sphere parametrization 
$$S=S^{N-1}_{\mathbb R,G,H}$$
standard when both filtered groups $G,H\subset S_\infty$ are chosen to be maximal.
\end{definition}

In this case, Proposition 13.20 and its proof tell us that $G,H$ encode all the monomial relations which hold in $S$. With these conventions, we have the following result from \cite{ba2}, \cite{ba3}, extending some previous findings from above, regarding the untwisted spheres:

\begin{theorem}
The standard parametrization of the $9$ main spheres is
$$
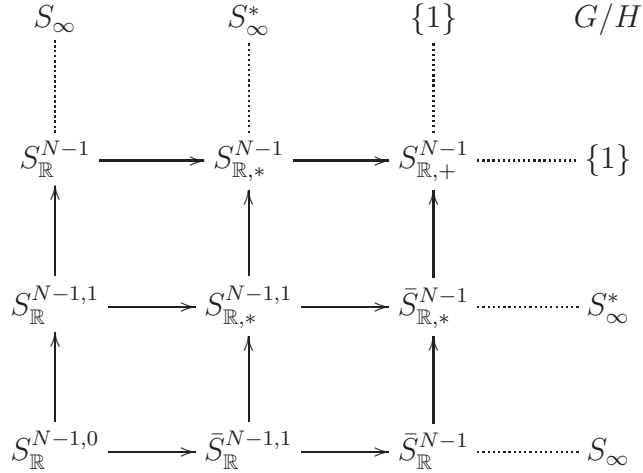
\xymatrix@R=11.5mm@C=11.5mm{
S_\infty\ar@{.}[d]&S_\infty^*\ar@{.}[d]&\{1\}\ar@{.}[d]&G/H\\
S^{N-1}_\mathbb R\ar[r]&S^{N-1}_{\mathbb R,*}\ar[r]&S^{N-1}_{\mathbb R,+}&\{1\}\ar@{.}[l]\\
S^{N-1,1}_\mathbb R\ar[r]\ar[u]&S^{N-1,1}_{\mathbb R,*}\ar[r]\ar[u]&\bar{S}^{N-1}_{\mathbb R,*}\ar[u]&S_\infty^*\ar@{.}[l]\\
S^{N-1,0}_\mathbb R\ar[r]\ar[u]&\bar{S}^{N-1,1}_\mathbb R\ar[r]\ar[u]&\bar{S}^{N-1}_\mathbb R\ar[u]&S_\infty\ar@{.}[l]}$$
so these spheres come from the $3\times 3=9$ pairs of groups among $\{1\}\subset S_\infty^*\subset S_\infty$.
\end{theorem}

\begin{proof}
The fact that we have parametrizations as above is known to hold for the 5 untwisted and twisted spheres. For the remaining 4 spheres the result follows by intersecting, by using the following formula, valid for any $E,F\subset S_\infty$:
$$S^{N-1}_{\mathbb R,E,F}\cap S^{N-1}_{\mathbb R,E',F'}=S^{N-1}_{\mathbb R,E\cup E',F\cup F'}$$

In order to prove now that the parametrizations are standard, we must compute the following two filtered groups, and show that we get the groups in the statement:
$$G=\left\{\sigma\in S_\infty\Big|{\rm the\ relations\ }\mathcal R_\sigma\ {\rm hold\ over\ }S\right\}$$ 
$$H=\left\{\sigma\in S_\infty\Big|{\rm the\ relations\ }\bar{\mathcal R}_\sigma\ {\rm hold\ over\ }S\right\}$$ 

As a first observation, by using the various inclusions between spheres, we just have to compute $G$ for the spheres on the bottom, and $H$ for the spheres on the left:
$$X=S^{N-1,0}_\mathbb R,\bar{S}^{N-1,1}_\mathbb R,\bar{S}^{N-1}_\mathbb R\implies G=S_\infty,S_\infty^*,\{1\}$$
$$X=S^{N-1,0}_\mathbb R,S^{N-1,1}_\mathbb R,S^{N-1}_\mathbb R\implies H=S_\infty,S_\infty^*,\{1\}$$

The results for $S^{N-1,0}_\mathbb R$ being clear, we are left with computing the remaining 4 groups, for the spheres $S^{N-1}_\mathbb R,\bar{S}^{N-1}_\mathbb R,S^{N-1,1}_\mathbb R,\bar{S}^{N-1,1}_\mathbb R$. The proof here goes as follows:

\medskip

(1) $S^{N-1}_\mathbb R$. According to the definition of $H=(H_k)$, we have:
\begin{eqnarray*}
H_k
&=&\left\{\sigma\in S_k\Big|x_{i_1}\ldots x_{i_k}=\varepsilon\left(
\ker\begin{pmatrix}i_1&\ldots&i_k\\
i_{\sigma(1)}&\ldots&i_{\sigma(k)}
\end{pmatrix}\right)
x_{i_{\sigma(1)}}\ldots x_{i_{\sigma(k)}},\forall i_1,\ldots,i_k\right\}\\
&=&\left\{\sigma\in S_k\Big|
\varepsilon\left(
\ker\begin{pmatrix}i_1&\ldots&i_k\\
i_{\sigma(1)}&\ldots&i_{\sigma(k)}
\end{pmatrix}\right)
=1,\forall i_1,\ldots,i_k\right\}\\
&=&\left\{\sigma\in S_k\Big|\varepsilon(\tau)=1,\forall\tau\leq\sigma\right\}\end{eqnarray*}

Now observe that for any permutation $\sigma\in S_k,\sigma\neq1_k$, we can always find a partition $\tau\leq\sigma$ satisfying the following condition:
$$\varepsilon(\tau)=-1$$

We deduce that we have $H_k=\{1_k\}$, and so $H=\{1\}$, as desired.

\medskip

(2) $\bar{S}^{N-1}_\mathbb R$. The proof of $G=\{1\}$ here is similar to the proof of $H=\{1\}$ in (1) above, by using the same combinatorial ingredient at the end.

\medskip

(3) $S^{N-1,1}_\mathbb R$. By definition of $H=(H_k)$, a permutation $\sigma\in S_k$ belongs to $H_k$ when the following condition is satisfied, for any choice of the indices $i_1,\ldots,i_k$:
$$x_{i_1}\ldots x_{i_k}=\varepsilon\left(\ker\begin{pmatrix}i_1&\ldots&i_k\\ i_{\sigma(1)}&\ldots&i_{\sigma(k)}\end{pmatrix}\right)x_{i_{\sigma(1)}}\ldots x_{i_{\sigma(k)}}$$

We have three cases here, as follows:

\medskip

-- When $|\ker i|=1$ this formula reads $x_r^k=x_r^k$, which is true. 

\medskip

-- When $|\ker i|\geq3$ this formula is automatically satisfied as well, because by using the relations $ab=ba$, and $abc=0$ for $a,b,c$ distinct, which both hold over $S^{N-1,1}_\mathbb R$, this formula reduces to $0=0$. 

\medskip

-- Thus, we are left with studying the case $|\ker i|=2$. Here the quantities on the left $x_{i_1}\ldots x_{i_k}$ will not vanish, so the sign on the right must be 1, and we therefore have:
$$H_k=\left\{\sigma\in S_k\Big|\varepsilon(\tau)=1,\forall\tau\leq\sigma,|\tau|=2\right\}$$

Now by coloring the legs of $\sigma$ clockwise $\circ\bullet\circ\bullet\ldots$, the above condition is satisfied when each string of $\sigma$ joins a white leg to a black leg. Thus $H_k=S_k^*$, as desired.

\medskip

(4) $\bar{S}^{N-1,1}_\mathbb R$. The proof of $G=S_\infty^*$ here is similar to the proof of $H=S_\infty^*$ in (3) above, by using the same combinatorial ingredient at the end.
\end{proof}

We will be back to the polygonal spheres in the next chapter, with a better axiomatization, and with a study of the associated quantum groups as well.

\section*{13d. Algebraic geometry}

In relation with the general algebraic geometry questions formulated in the beginning of this chapter, and more precisely with the free algebra needed for developing free geometry, we have now at least one clear result on the subject, namely Theorem 13.22. But, long way to go. Modern algebraic geometry is based on modern commutative algebra, as developed by Hilbert, Noether, Zariski and many others, not to talk about Grothendieck and schemes, and such algebraic knowledge is completely lacking in the free setting, preventing for the moment any serious development of free algebraic geometry. This will most likely take a very long time, needing, to start with, a fresh new generation of mathematicians, finding things like Theorem 13.22 trivial, or even lame.

\bigskip

So, forget about modern times, and let us go back instead to the ancient Greeks, with the idea in mind of having some fun with conics. Let us start with:

\begin{definition}
A real compact hypersurface in $N$ variables, denoted $X_f\subset\mathbb R^N_+$, is the abstract spectrum of a universal $C^*$-algebra of the following type,
$$C(X_f)=C^*\left(x_1,\ldots,x_N\Big|x_i=x_i^*,f(x_1,\ldots,x_N)=0\right)$$
with the noncommutative polynomial $f\in\mathbb R<x_1,\ldots,x_N>$ being such the maximal $C^*$-norm on the complex $*$-algebra $\mathbb C<x_1,\ldots,x_N>/(f)$ is bounded.
\end{definition}

The boundedness condition above is something quite non-trivial, usually related to tricky operator theory, like sums of squares (SOS) theorems, and so on. If this condition is satisfied, we agree to say that ``$X_f$ exists''. As a first result now, we have:

\begin{theorem}
In order for $X_f$ to exist, the real algebraic manifold
$$X_f^\times=\left\{x\in\mathbb R^N\Big|f(x_1,\ldots,x_N)=0\right\}$$
must be compact. In addition, in this case we have $||x_i||_\times\leq||x_i||$, for any $i$.
\end{theorem}

\begin{proof}
Assuming that $X_f$ exists, our claim is that the algebra of continuous functions on the manifold $X_f^\times$ in the statement appears from $C(X_f)$ as follows:
$$C(X_f^\times)=C(X_f)\Big/\Big<[x_i,x_j]=0\Big>$$

But this is clear, by applying the Gelfand theorem, and by using as well the Stone-Weierstrass theorem, in order to have arrows in both directions, mapping $x_i\to x_i$. With this in hand, we have an embedding of compact quantum spaces, as follows:
$$X_f^\times\subset X_f$$

The norm estimate is now clear, because such embeddings increase the norms.
\end{proof}

In practice now, let us first discuss the quadratic case. The existence result here, which is very similar to the one from the classical case, is as follows:

\begin{proposition}
Given a quadratic polynomial $f\in\mathbb R<x_1,\ldots,x_N>$, written as
$$f=\sum_{ij}A_{ij}x_ix_j+\sum_iB_ix_i+C$$
the following conditions are equivalent:
\begin{enumerate}
\item $X_f$ exists.

\item $X_f^\times$ is compact.

\item The symmetric matrix $Q=\frac{A+A^t}{2}$ is positive or negative.
\end{enumerate}
\end{proposition}

\begin{proof}
The implication $(1)\implies(2)$ being known from Theorem 13.24, and the equivalence $(2)\iff(3)$ being well-known, we are left with proving $(3)\implies(1)$. As a first remark, by applying the adjoint, our manifold $X_f$ is defined by:
$$\begin{cases}
\sum_{ij}A_{ij}x_ix_j+\sum_iB_ix_i+C=0\\
\sum_{ij}A_{ij}x_jx_i+\sum_iB_ix_i+C=0
\end{cases}$$

In terms of $P=\frac{A-A^t}{2}$ and $Q=\frac{A+A^t}{2}$, these equations can be written as:
$$\begin{cases}
\sum_{ij}P_{ij}x_ix_j=0\\
\sum_{ij}Q_{ij}x_ix_j+\sum_iB_ix_i+C=0
\end{cases}$$

Let us first examine the second equation. When regarding $x$ as a column vector, and $B$ as a row vector, this equation becomes an equality of $1\times1$ matrices, as follows:
$$x^tQx+Bx+C=0$$

Now let us assume that $Q$ is positive or negative. Up to a sign change, we can assume $Q>0$. We can write $Q=UDU^t$, with $D=diag(d_i)$ and $d_i>0$, and with $U\in O_N$. In terms of the vector $y=U^tx$, and with $E=BU$, our equation becomes:
$$y^tDy+Ey+C=0$$

By reverting back to sums and indices, this equation reads:
$$\sum_id_iy_i^2+\sum_ie_iy_i+C=0$$

Now by making squares, this equation takes the following form:
$$\sum_id_i\left(y_i+\frac{e_i}{2d_i}\right)^2=c$$

By positivity, we deduce that we have the following estimate:
$$\left\|y_i+\frac{e_i}{2d_i}\right\|^2\leq\frac{|c|}{d_i}$$

Thus our hypersurface $X_f$ is well-defined, and we are done. 
\end{proof}

We recall that, up to a linear changes of coordinates, there is only one non-trivial compact quadric in $\mathbb R^N$, namely $S^{N-1}_\mathbb R$.  In the noncommutative setting the situation is more complicated, because the first equation of $X_f$ in the above proof, namely $\sum_{ij}P_{ij}x_ix_j=0$ with $P=\frac{A-A^t}{2}$, that we have neglected so far, and which is trivial in the classical case, is no longer trivial. By taking into account this equation, we are led to:

\begin{theorem}
Up to linear changes of coordinates, the free compact quadrics in $\mathbb R^N_+$ are the empty set, the point, the standard free sphere $S^{N-1}_{\mathbb R,+}$, defined by
$$\sum_ix_i^2=1$$
and some intermediate spheres $S^{N-1}_\mathbb R\subset S\subset S^{N-1}_{\mathbb R,+}$, which can be explicitly characterized. Moreover, for all these free quadrics, we have $||x_i||=||x_i||_\times$, for any $i$.
\end{theorem}

\begin{proof}
We use the computations from the proof of Proposition 13.25. The first equation there, making appear the matrix $P=\frac{A-A^t}{2}$, is as follows:
$$\sum_{ij}P_{ij}x_ix_j=0$$ 

As for the second equation, up to a linear change of the coordinates, this reads:
$$\sum_iz_i^2=c$$

At $c<0$ we obtain the empty set. At $c=0$ we must have $z=0$, and depending on whether the first equation is satisfied or not, we obtain either a point, or the empty set. At $c>0$ now, we can assume by rescaling $c=1$, and our second equation reads:
$$X_f\subset S^{N-1}_{\mathbb R,+}$$

As a conclusion, the solutions here are certain subspaces $S\subset S^{N-1}_{\mathbb R,+}$ which appear via equations of type $\sum_{ij}P_{ij}x_ix_j=0$, with $P\in M_N(\mathbb R)$ being antisymmetric, and with $x_1,\ldots,x_N$ appearing via $z_1,\ldots,z_N$ via a linear change of variables. Since when redoing the above computation with $X_f^\times$ at the place of $X_f$, we obtain $X_f=S^{N-1}_\mathbb R$, we conclude that our subspaces $S\subset S^{N-1}_{\mathbb R,+}$ must satisfy:
$$S^{N-1}_\mathbb R\subset S\subset S^{N-1}_{\mathbb R,+}$$

Thus, we are left with investigating which such subspaces can indeed be solutions. Observe that both the extreme cases can appear as solutions, as shown by:
\begin{eqnarray*}
X_{2x^2+y^2+\frac{3}{2}xy+\frac{1}{2}yx}&=&S^1_\mathbb R\\
X_{2x^2+y^2+xy+yx}&=&S^1_{\mathbb R,+}
\end{eqnarray*}

Finally, the last assertion is clear for the empty set and for the point, and for the remaining hypersurfaces, this follows from $S^{N-1}_\mathbb R\subset S\subset S^{N-1}_{\mathbb R,+}$.
\end{proof}

Here is now yet another version of Proposition 13.25, this time by using an opposite idea, namely using as many linear transformations as possible:

\begin{proposition}
Given $M$ real linear functions $L_1,\ldots,L_M$ in $N$ noncommuting variables $x_1,\ldots,x_N$, the following are equivalent:
\begin{enumerate}
\item $\sum_kL_k(x_1,\ldots,x_N)^2=1$ defines a compact hypersurface in $\mathbb R^N$.

\item $\sum_kL_k(x_1,\ldots,x_N)^2=1$ defines a compact quantum hypersurface.

\item The matrix formed by the coefficients of $L_1,\ldots,L_M$ has rank $N$.
\end{enumerate}
\end{proposition}

\begin{proof}
The equivalence $(1)\iff(2)$ follows from $(1)\iff(2)$ in Proposition 13.25, because the surfaces under investigation are quadrics. As for the equivalence $(2)\iff(3)$, this is well-known. More precisely, our equation can be written as:
\begin{eqnarray*}
1
&=&\sum_kL_k(x_1,\ldots,x_N)^2\\
&=&\sum_k\sum_iL_{ki}x_i\sum_jL_{kj}x_j\\
&=&\sum_{ij}(L^tL)_{ij}x_ix_j
\end{eqnarray*}

Thus, in the context of Proposition 13.25, the underlying square matrix $A\in M_N(\mathbb R)$ is given by $A=L^tL$. It follows that we have $Q=A=L^tL$, and so the condition $Q>0$ is equivalent to $L^tL$ being invertible, and so to $L$ to have rank $N$, as claimed.
\end{proof}

Summarizing, in what concerns the quadrics, the noncommutative theory basically parallels the usual classical theory, with just a few minor twists. In higher degree, however, things look amazingly complicated, because even construcing hypersurfaces via quite trivial sums of squares leads to non-trivial operator theory questions.

\section*{13e. Exercises} 

There has been a lot of non-trivial algebra in this chapter, and our questions here will be on this precise topic, non-trivial algebra. First, we have:

\begin{exercise}
Work out a theory of monomial spheres, in the complex case. Once this done, work out as well a theory of standard parametrization in the complex case.
\end{exercise}

The difficulties in dealing with this question were already explained, in the above.

\begin{exercise}
Extend the theory of standard parametrization that we developed in the above, from the sphere case, to the case of more general manifolds. 
\end{exercise}

For a bonus exercise, try further building on Proposition 13.27, in higher degree. This is difficult, and very interesting, fun guarantee. You will most likely need help from a good algebraic geometer, a good operator theorist, and a computer too.

\chapter{Polygonal spheres}

\section*{14a. Polygonal spheres}

In this chapter we build on the findings from the previous chapter, still following \cite{ba2}, \cite{ba3}, with the idea in mind that all this material belongs to a new and exciting area of noncommutative algebra, which can help in building an algebraic geometry theory for the free manifolds, and which therefore needs to be prioritarily developed.

\bigskip

As in the previous chapter, due to various technical difficulties with the complex case, at least at this stage of the things, we will basically restrict the attention to the real case. The main objects of study here are the 3 real spheres, which are as follows:
$$S^{N-1}_\mathbb R\subset S^{N-1}_{\mathbb R,*}\subset S^{N-1}_{\mathbb R_+}$$

We have seen that the study of the relations between the coordinates $x_1,\ldots,x_N$ of these real spheres naturally leads to the twisted versions of these spheres, namely:
$$\bar{S}^{N-1}_\mathbb R\subset\bar{S}^{N-1}_{\mathbb R,*}\subset S^{N-1}_{\mathbb R_+}$$

More precisely, the study of the algebraic relations between the coordinates $x_1,\ldots,x_N$ of the real spheres leads to the study of the various intersections between the twisted and untwisted spheres. These $3\times3$ intersections form a square diagram, as follows:
$$\xymatrix@R=13mm@C=1mm{
S^{N-1}_\mathbb R\ar[rr]&&S^{N-1}_{\mathbb R,*}\ar[rrr]&&&S^{N-1}_{\mathbb R,+}\\
S^{N-1}_\mathbb R\cap\bar{S}^{N-1}_{\mathbb R,*}\ar[rr]\ar[u]&&S^{N-1}_{\mathbb R,*}\cap\bar{S}^{N-1}_{\mathbb R,*}\ar[rrr]\ar[u]&&&\bar{S}^{N-1}_{\mathbb R,*}\ar[u]\\
S^{N-1}_\mathbb R\cap\bar{S}^{N-1}_\mathbb R\ar[rr]\ar[u]&&S^{N-1}_{\mathbb R,*}\cap\bar{S}^{N-1}_\mathbb R\ar[rrr]\ar[u]&&&\bar{S}^{N-1}_\mathbb R\ar[u]}$$

We have seen as well that these intersections all appear as ``polygonal spheres'', which are certain real algebraic manifolds, according to the following result:

\index{intersections of spheres}

\begin{theorem}
The $5$ main spheres, and the intersections between them, are
$$\xymatrix@R=12mm@C=14mm{
S^{N-1}_\mathbb R\ar[r]&S^{N-1}_{\mathbb R,*}\ar[r]&S^{N-1}_{\mathbb R,+}\\
S^{N-1,1}_\mathbb R\ar[r]\ar[u]&S^{N-1,1}_{\mathbb R,*}\ar[r]\ar[u]&\bar{S}^{N-1}_{\mathbb R,*}\ar[u]\\
S^{N-1,0}_\mathbb R\ar[r]\ar[u]&\bar{S}^{N-1,1}_\mathbb R\ar[r]\ar[u]&\bar{S}^{N-1}_\mathbb R\ar[u]}$$
where $\dot{S}^{N-1,d-1}_{\mathbb R,\times}\subset\dot{S}^{N-1}_{\mathbb R,\times}$ is obtained by assuming $x_{i_0}\ldots x_{i_d}=0$, for $i_0,\ldots,i_d$ distinct.
\end{theorem}

\begin{proof}
This is something that we know from chapter 13, the idea being that commutation and anticommutation produces vanishing relations.
\end{proof}

We refer to chapter 13 for more on these spheres, including their algebraic axiomatization and main properties, and the ``standard parametrization'' result there.

\bigskip

In this chapter we discuss the extension of the axiomatics for abstract noncommutative geometries that we have, from chapter 4, in order to cover both the twisted and untwisted cases, and the intersections. This is a very natural question, in view of our findings from chapter 13. For this purpose, we are in need of some new quantum isometry group computations. In order to deal with the polygonal spheres, we will need:

\index{Bhowmick-Goswami trick}

\begin{proposition}
Assume that $X\subset S^{N-1}_\mathbb R$ is invariant, for any $i$, under: 
$$x_i\to-x_i$$
\begin{enumerate}
\item If the coordinates $x_1,\ldots,x_N$ are linearly independent inside $C(X)$, then the group 
$$G(X)=G^+(X)\cap O_N$$
consists of the usual isometries of $X$.

\item In addition, in the case where the products of coordinates 
$$\left\{x_ix_j\Big|i\leq j\right\}$$
are linearly independent inside $C(X)$, we have $G^+(X)=G(X)$.
\end{enumerate}
\end{proposition}

\begin{proof}
This is a standard trick, that we will heavily use here, which follows from Bhowmick-Goswami \cite{bg1} and Goswami \cite{go2}, the idea being as follows:

\medskip

(1) The assertion here is well-known, $G(X)=G^+(X)\cap O_N$ being by definition the biggest subgroup $G\subset O_N$ acting affinely on $X$. We refer to \cite{go2} for details, and for a number of noncommutative extensions of this fact, with $G(X)$ replaced by $G^+(X)$.

\medskip

(2) Consider an arbitrary coaction map on the algebra $C(X)$, as follows:
$$\Phi:C(X)\to C(X)\otimes C(G)$$
$$\Phi(x_i)=\sum_jx_j\otimes u_{ji}$$

In order to establish the result, we must prove that the variables $u_{ij}$ commute. But this follows by using a strandard trick, from \cite{bg1}, that we will briefly recall now. We can write the action of $\Phi$ on the commutators between the coordinates as follows:
$$\Phi([x_i,x_j])=\sum_{k\leq l}\left(1-\frac{\delta_{kl}}{2}\right)x_kx_l\otimes\left([u_{ki},u_{lj}]-[u_{kj},u_{li}]\right)$$

Now since the variables $\{x_kx_l|k\leq l\}$ were assumed to be linearly independent, we obtain from this that we have the following formula:
$$[u_{ki},u_{lj}]=[u_{kj},u_{li}]$$

Moreover, if we apply now the antipode we further obtain:
$$[u_{jl},u_{ik}]=[u_{il},u_{jk}]$$

By relabelling, this gives the following formula:
$$[u_{ki},u_{lj}]=[u_{li},u_{kj}]$$

Now by comparing with the original equality of commutators, from above, we conclude from this that we have a commutation relation, as follows:
$$[u_{ki},u_{lj}]=0$$

Thus, we are led to the conclusion in the statement. See \cite{bg1}.
\end{proof}

With the above notion in hand, let us investigate the polygonal spheres. We recall that, according to the various computations from the previous chapters, the quantum isometry groups of the 5 main spheres are as follows:
$$\xymatrix@R=18mm@C=18mm{
S^{N-1}_\mathbb R\ar[r]\ar@{~}[d]&S^{N-1}_{\mathbb R,*}\ar[r]\ar@{~}[d]&S^{N-1}_{\mathbb R,+}\ar@{~}[d]&\bar{S}^{N-1}_{\mathbb R,*}\ar[l]\ar@{~}[d]&\bar{S}^{N-1}_\mathbb R\ar[l]\ar@{~}[d]\\
O_N\ar[r]&O_N^*\ar[r]&O_N^+&\bar{O}_N^*\ar[l]&\bar{O}_N\ar[l]}$$

In the polygonal sphere case now, we begin with the computations of the quantum isometry groups in the classical case. We have here the following result, from \cite{ba2}:

\index{polygonal sphere}

\begin{theorem}
The quantum isometry groups of the classical polygonal spheres
$$S^{N-1,d-1}_\mathbb R=\left\{x\in S^{N-1}_\mathbb R\Big|x_{i_0}\ldots x_{i_d}=0,\forall i_0,\ldots,i_d\ {\rm distinct}\right\}$$
are as follows:
\begin{enumerate}
\item At $d=1$ we obtain the free hyperoctahedral group $H_N^+$.

\item At $d=2,\ldots,N-1$ we obtain the hyperoctahedral group $H_N$.

\item At $d=N$ we obtain the orthogonal group $O_N$.
\end{enumerate}
\end{theorem}

\begin{proof}
Observe first that the sphere $S^{N-1,d-1}_\mathbb R$ appears by definition as a union on $\binom{N}{d}$ copies of the sphere $S^{d-1}_\mathbb R$, one for each choice of $d$ coordinate axes, among the coordinate axes of $\mathbb R^N$. We can write this decomposition as follows, with $I_N=\{1,\ldots,N\}$:
$$S^{N-1,d-1}_\mathbb R=\bigcup_{I\subset I_N,|I|=d}(S^{d-1}_\mathbb R)^I$$

With this observation in hand, the proof goes as follows:

\medskip

(1) At $d=1$ our sphere is given by the following formula:
$$S^{N-1,0}_\mathbb R=\mathbb Z_2^{\oplus N}$$

To be more precise, what we have here is the set formed by the endpoints of the $N$ copies of $[-1,1]$ on the coordinate axes of $\mathbb R^N$. Thus by the free wreath product results in \cite{bbc} the corresponding quantum isometry group is $H_N^+$:
\begin{eqnarray*}
G^+(S^{N-1,0}_\mathbb R)
&=&G^+(\mathbb Z_2^{\oplus N})\\
&=&G^+(|\ |\ldots |\ |)\\
&=&\mathbb Z_2\wr_*S_N^+\\
&=&H_N^+
\end{eqnarray*}

(2) In order to discuss now the case $d\geq2$, the idea is to use Proposition 14.2 (2). Our claim is that the following elements are linearly independent:
$$\left\{x_ix_j\Big|i\leq j\right\}$$

Since $S^{N-1,1}_\mathbb R\subset S^{N-1,d}_\mathbb R$, we can restrict attention to the case $d=2$. Here the above decomposition is as follows, where $\mathbb T^{\{i,j\}}$ denote the various copies of $\mathbb T$:
$$S^{N-1,d-1}_\mathbb R=\bigcup_{i<j}\mathbb T^{\{i,j\}}$$

Now observe that the following elements are linearly independent over $\mathbb T\subset\mathbb R^2$:
$$\left\{x^2,y^2,xy\right\}$$

We deduce that the following elements are linearly independent over $S^{N-1,d-1}_\mathbb R$:
$$\left\{x_ix_j\Big|i\leq j\right\}$$

Thus, our claim is proved, and so Proposition 14.2 (2) applies, and gives:
$$G^+(X)=G(X)$$

Thus, we are left with proving the following formula, for any $d\in\{2,\ldots,N-1\}$:
$$G(X)=H_N$$

-- Let us first discuss the case $d=2$. By using the decomposition formula from the beginning of the proof, here any affine isometric action $U\curvearrowright S^{N-1,1}_\mathbb R$ must permute the $\binom{N}{2}$ circles $\mathbb T^I$, so we can write, for a certain permutation of the indices $I\to I'$:
$$U(\mathbb T^I)=\mathbb T^{I'}$$

Now since $U$ is bijective, we deduce that for any $I,J$ we have:
$$U\left(\mathbb T^I\cap\mathbb T^J\right)=\mathbb T^{I'}\cap\mathbb T^{J'}$$

The point now is that for $|I\cap J|=0,1,2$ we have:
$$\mathbb T^I\cap\mathbb T^J\simeq\emptyset,\{-1,1\},\mathbb T$$

By taking now the union over $I,J$ with $|I\cap J|=1$, we deduce that: 
$$U(\mathbb Z_2^{\oplus N})=\mathbb Z_2^{\oplus N}$$

Thus we must have $U\in H_N$, and we are done with the case $d=2$.

\medskip 

-- In the general case now, $d\in\{2,\ldots,N-1\}$, we can proceed similarly, by recurrence. Indeed, for any subsets $I,J\subset I_N$ with $|I|=|J|=d$ we have:
$$(S^{d-1}_\mathbb R)^I\cap(S^{d-1}_\mathbb R)^J=(S^{|I\cap J|-1}_\mathbb R)^{I\cap J}$$

By using $d\leq N-1$, we deduce that we have the following formula:
$$S^{N-1,d-2}_\mathbb R=\bigcup_{|I|=|J|=d,|I\cap J|=d-1}(S^{|I\cap J|-1}_\mathbb R)^{I\cap J}$$

On the other hand, by using exactly the same argument as in the $d=2$ case, we deduce that the space on the right is invariant, under any affine isometric action on $S^{N-1,d-1}_\mathbb R$. Thus by recurrence we obtain, as desired, that we have:
$$G(S^{N-1,d-1}_\mathbb R)
=G(S^{N-1,d-2}_\mathbb R)
=H_N$$

(3) At $d=N$ the result is known since \cite{bgo}, with the proof coming from the equality $G^+(X)=G(X)$, deduced from Proposition 14.2 (2), as explained above.
\end{proof}

The study in the twisted case is considerably more difficult than in the classical case, and we have complete results only at $d=1,2,N$, as follows:

\index{twisted polygonal sphere}

\begin{theorem}
The quantum isometry groups of twisted polygonal spheres, given by
$$C(\bar{S}^{N-1,d-1}_\mathbb R)=C(\bar{S}^{N-1}_\mathbb R)\Big/\Big<x_{i_0}\ldots x_{i_d}=0,\forall i_0,\ldots,i_d\ {\rm distinct}\Big>$$ 
are as follows:
\begin{enumerate}
\item At $d=1$ we obtain the free hyperoctahedral group $H_N^+$.

\item At $d=2$ we obtain the hyperoctahedral group $H_N$.

\item At $d=N$ we obtain the twisted orthogonal group $\bar{O}_N$.
\end{enumerate}
\end{theorem}

\begin{proof}
The idea is to adapt the proof of Theorem 14.3:

\medskip

(1) At $d=1$ the situation is very simple, because we have:
$$\bar{S}^{N-1,0}_\mathbb R=S^{N-1,0}_\mathbb R=\mathbb Z_2^{\oplus N}$$

By Theorem 14.3 (1), coming from the free wreath product computations in \cite{bbc}, the corresponding quantum isometry group is indeed $H_N^+$. 

\medskip

(2) In order to deal now with the case $d=2$, in analogy with what was done before in the classical case, as a first ingredient, we will need the twisted analogue of the trick from \cite{bg1}, explained in the proof of Proposition 14.2 (2). 

\medskip

This twisted trick is known to work for the twisted sphere $\bar{S}^{N-1}_\mathbb R$ itself, and the situation is in fact similar for any closed subset $X\subset\bar{S}^{N-1}_\mathbb R$, having the property that the following variables are linearly indepedent:
$$\left\{x_ix_j\Big|i\leq j\right\}$$

More presisely, our claim is that under this linear independence assumption, if a quantum group $G\subset O_N^+$ acts on $X$, then we must have:
$$G\subset\bar{O}_N$$

Indeed, consider a coaction map, written as follows:
$$\Phi(x_i)=\sum_jx_j\otimes u_{ji}$$

By making products, we have the following formula:
$$\Phi(x_ix_j)=\sum_k x_k^2\otimes u_{ki}u_{kj}+\sum_{k<l}x_kx_l\otimes(u_{ki}u_{lj}-u_{li}u_{kj})$$

We deduce that with $[[a,b]]=ab+ba$ we have the following formula:
$$\Phi([[x_i,x_j]])=\sum_kx_k^2\otimes [[u_{ki},u_{kj}]]+\sum_{k<l}x_kx_l\otimes ([u_{ki},u_{lj}]-[u_{li},u_{kj}])$$

Now assuming $i\neq j$, we have $[[x_i,x_j]]=0$, and we therefore obtain, for any $k$:
$$[[u_{ki},u_{kj}]]=0$$

We also have, for any $k<l$, the following formula:
$$[u_{ki},u_{lj}]=[u_{li},u_{kj}]$$

By applying the antipode and then by relabelling, the latter relation gives:
$$[u_{ki},u_{lj}]=0$$

Thus we have reached to the defining relations for the quantum group $\bar{O}_N$, from chapter 11, and so we have $G\subset\bar{O}_N$, as claimed.

\medskip

Our second claim is that the above trick applies to any $\bar{S}^{N-1,d-1}_\mathbb R$ with $d\geq2$. Consider indeed the following maps, obtained by setting $x_k=0$ for $k\neq i,j$:
$$\pi_{ij}:C(\bar{S}^{N-1,d-1}_\mathbb R)\to C(\bar{S}^1_\mathbb R)$$

By using these maps, we conclude that the following variables are indeed linearly independent over $\bar{S}^{N-1,d-1}_\mathbb R$, as desired:
$$\left\{x_ix_j\Big|i\leq j\right\}$$

Summarizing, we have proved so far that if a compact quantum group $G\subset O_N^+$ acts on a polygonal sphere $\bar{S}^{N-1,d-1}_\mathbb R$ with $d\geq2$, then we must have:
$$G\subset\bar{O}_N$$

In order to finish, we must now adapt the second part of the proof of Proposition 14.2, and since this is quite unobvious at $d\geq3$, due to various technical reasons, we will restrict now attention to the case $d=2$, as in the statement.

\medskip

So, consider a compact quantum group $G\subset\bar{O}_N$. The question is that of understanding when we have a coaction map, as follows:
$$\Phi:C(\bar{S}^{N-1,1}_\mathbb R)\to C(G)\otimes C(\bar{S}^{N-1,1}_\mathbb R)$$
$$\Phi(x_i)=\sum_jx_j\otimes u_{ji}$$

In order for this to happen, the elements $X_i=\sum_jx_j\otimes u_{ji}$ must satisfy the relations $X_iX_jX_k=0$, for any $i,j,k$ distinct.

\medskip

So, let us compute $X_iX_jX_k$ for $i,j,k$ distinct. We have:
\begin{eqnarray*}
X_iX_jX_k
&=&\sum_{abc}x_ax_bx_c\otimes u_{ai}u_{bj}u_{ck}\\
&=&\sum_{a,b,c\ not \ distinct}x_ax_bx_c\otimes u_{ai}u_{bj}u_{ck}\\
&=&\sum_{a\neq b}x_a^2x_b\otimes u_{ai}u_{aj}u_{bk}+\sum_{a\neq b}x_ax_bx_a\otimes u_{ai}u_{bj}u_{ak}\\
&+&\sum_{a\neq b}x_bx_a^2\otimes u_{bi}u_{aj}u_{ak}+\sum_ax_a^3\otimes u_{ai}u_{aj}u_{ak}
\end{eqnarray*}

By using $x_ax_bx_a=-x_a^2x_b$ and $x_bx_a^2=x_a^2x_b$, we deduce that we have:
\begin{eqnarray*}
X_iX_jX_k
&=&\sum_{a\neq b}x_a^2x_b\otimes(u_{ai}u_{aj}u_{bk}-u_{ai}u_{bj}u_{ak}+u_{bi}u_{aj}u_{ak})\\
&+&\sum_ax_a^3\otimes u_{ai}u_{aj}u_{ak}\\
&=&\sum_{ab}x_a^2x_b\otimes(u_{ai}u_{aj}u_{bk}-u_{ai}u_{bj}u_{ak}+u_{bi}u_{aj}u_{ak})
\end{eqnarray*}

By using now the defining relations for $\bar{O}_N$, which apply to the variables $u_{ij}$, this formula can be written in a cyclic way, as follows:
$$X_iX_jX_k=\sum_{ab}x_a^2x_b\otimes(u_{ai}u_{aj}u_{bk}+u_{aj}u_{ak}u_{bi}+u_{ak}u_{ai}u_{bj})$$

We use now the fact that the variables on the left, namely $x_a^2x_b$, are linearly independent. We conclude that, in order for our quantum group $G\subset\bar{O}_N$ to act on $\bar{S}^{N-1,1}_\mathbb R$, its coordinates must satisfy the following relations, for any $i,j,k$ distinct:
$$u_{ai}u_{aj}u_{bk}+u_{aj}u_{ak}u_{bi}+u_{ak}u_{ai}u_{bj}=0$$

By multiplying to the right by $u_{kb}$ and then by summing over $b$, we deduce from this that we have, for any $i,j$:
$$u_{ai}u_{aj}=0$$

Now since the quotient of $C(\bar{O}_N)$ by these latter relations is the algebra $C(H_N)$, we conclude that we have, as claimed:
$$G^+(\bar{S}^{N-1,1}_\mathbb R)=H_N$$

(3) At $d=N$ the result is already known, and its proof follows in fact from the ``twisted trick'' explained in the proof of (2) above, applied to $\bar{S}^{N-1}_\mathbb R$.
\end{proof}

\section*{14b. Quantum groups}

In general now, the idea will be that the quantum isometry groups of the intersections of the spheres will basically appear as intersections of the corresponding quantum isometry groups. To start with, we must compute the intersections between the quantum orthogonal groups and their twists. The result here, which is similar to the one for the corresponding spheres, established in chapter 13, is as follows:

\index{intersection of orthogonal groups}

\begin{proposition}
The $5$ orthogonal groups and their twists, and the intersections between them, are as follows, at any $N\geq3$:
$$\xymatrix@R=12mm@C=16mm{
O_N\ar[r]&O_N^*\ar[r]&O_N^+\\
H_N\ar[r]\ar[u]&H_N^*\ar[r]\ar[u]&\bar{O}_N^*\ar[u]\\
H_N\ar[r]\ar[u]&H_N\ar[r]\ar[u]&\bar{O}_N\ar[u]}$$
At $N=2$ the same holds, with the diagram
$$\xymatrix@R=12mm@C=16mm{
O_2\ar[r]&O_2^+\\
H_2\ar[u]\ar[r]&\bar{O}_2\ar[u]}$$
replacing the lower left square, in the above big square.
\end{proposition}

\begin{proof}
We have to study 4 quantum group intersections, as follows:

\medskip

(1) $O_N\cap\bar{O}_N$. Here an element $U\in O_N$ belongs to the intersection when its entries satisfy $ab=0$ for any $a\neq b$ on the same row or column of $U$. But this means that our matrix $U\in O_N$ must be monomial, and so we get $U\in H_N$, as claimed.

\medskip

(2) $O_N\cap\bar{O}_N^*$. At $N=2$ the defining relations for $\bar{O}_N^*$ dissapear, and so we have the following computation, which leads to the conclusion in the statement:
$$O_2\cap\bar{O}_2^*=O_2\cap O_2^+=O_2$$

At $N\geq3$ now, the following inclusion is clear: 
$$H_N\subset O_N\cap\bar{O}_N^*$$ 

In order to prove the converse inclusion, pick $U\in O_N$ in the intersection, and assume that $U$ is not monomial. By permuting the entries we can further assume:
$$U_{11}\neq0\quad,\quad 
U_{12}\neq0$$

From $U_{11}U_{12}U_{i3}=0$ for any $i$ we deduce that the third column of $U$ is filled with $0$ entries, a contradiction. Thus we must have $U\in H_N$, as claimed.

\medskip

(3) $O_N^*\cap\bar{O}_N$. At $N=2$ we have the following computation, as claimed:
$$O_2^*\cap\bar{O}_2=O_2^+\cap\bar{O}_2=\bar{O}_2$$

At $N\geq3$ now, the best is to use the result in (4) below. Indeed, knowing that we have $O_N^*\cap\bar{O}_N^*=H_N^*$, our intersection is then:
$$G=H_N^*\cap\bar{O}_N$$

Now since the standard coordinates on $H_N^*$ are known to satisfy $ab=0$ for $a\neq b$ on the same row or column of $u$, the commutation/anticommutation relations defining $\bar{O}_N$ reduce to plain commutation relations. Thus $G$ follows to be classical, $G\subset O_N$, and by using (1) above we obtain the following formula, as claimed:
\begin{eqnarray*}
G
&=&H_N^*\cap\bar{O}_N\cap O_N\\
&=&H_N^*\cap H_N\\
&=&H_N
\end{eqnarray*}

(4) $O_N^*\cap\bar{O}_N^*$. The result here is non-trivial, and we must use the half-liberation technology from \cite{bdu}. The quantum group $H_N^\times=O_N^*\cap\bar{O}_N^*$ is indeed half-classical in the sense of \cite{bdu}, and since we have $H_N^*\subset H_N^\times$, this quantum group is not classical. Thus the main result in \cite{bdu} applies, and shows that $H_N^\times\subset O_N^*$ must come, via the crossed product construction there, from an intermediate compact group, as follows:
$$\mathbb T\subset G\subset U_N$$

Now observe that the standard coordinates on $H_N^\times$ are by definition subject to the conditions $abc=0$ when $(r,s)=(\leq2,3),(3,\leq2)$, with the notations and conventions from chapter 11 above. It follows that the standard coordinates on $G$ are subject to the conditions $\alpha\beta\gamma=0$ when $(r,s)=(\leq2,3),(3,\leq2)$, where $r,s=span(a,b,c)$, and $\alpha=a,a^*,\beta=b,b^*,\gamma=c,c^*$. Thus we have an inclusion as follows:
$$G\subset\bar{U}_N^*$$

We deduce that we have an inclusion as follows, with $K_N^\circ=U_N\cap\bar{U}_N^*$:
$$G\subset K_N^\circ$$

But this intersection can be computed exactly as in the real case, in the proof of (2) above, and we obtain $K_2^\circ=U_2$, and $K_N^\circ=\mathbb T\wr S_N$ at $N\geq 3$. 

\medskip

But the half-liberated quantum groups obtained from $U_2$ and $\mathbb T\wr S_N$ via the half-liberation construction in \cite{bdu} are well-known, these being $O_2^*=O_2^+$ and $H_N^*$. Thus by functoriality we have $H_2^\times\subset O_2^+$ and $H_N^\times\subset H_N^*$ at $N\geq 3$, and since the reverse inclusions are clear, we obtain $H_2^\times=O_2^+$ and $H_N^\times=H_N^*$ at $N\geq 3$, as claimed.
\end{proof}

Let us go back now to the sphere left, namely $S^{N-1,1}_{\mathbb R,*}$. Things are quite tricky here, and we will need the following technical result, from Raum-Weber \cite{rwe}:

\index{intermediate hyperoctahedral group}

\begin{proposition}
Let $H_N^{[\infty]}\subset O_N^+$ be the compact quantum group obtained via the relations $abc=0$, whenever $a\neq c$ are on the same row or column of $u$. 
\begin{enumerate}
\item We have inclusions $H_N^*\subset H_N^{[\infty]}\subset H_N^+$.

\item We have $ab_1\ldots b_rc=0$, whenever $a\neq c$ are on the same row or column of $u$.

\item We have $ab^2=b^2a$, for any two entries $a,b$ of $u$.
\end{enumerate}
\end{proposition}

\begin{proof}
We briefly recall the proof in \cite{rwe}, for future use in what follows. Our first claim is that $H_N^{[\infty]}$ comes, as an easy quantum group, from the following diagram:
$$
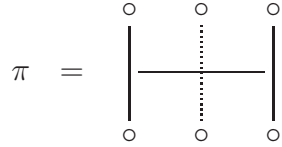
\xymatrix@R=5mm@C=0.1mm{&\\\pi\ \ =\\&}\xymatrix@R=5mm@C=5mm{
\circ\ar@{-}[dd]&\circ\ar@{.}[dd]&\circ\ar@{-}[dd]\\
\ar@{-}[rr]&&\\
\circ&\circ&\circ}$$

Indeed, this diagram acts via the following linear map:
$$T_\pi(e_{ijk})=\delta_{ik}e_{ijk}$$

We therefore have the following formula:
\begin{eqnarray*}
T_\pi u^{\otimes 3}e_{abc}
&=&T_\pi\sum_{ijk}e_{ijk}\otimes u_{ia}u_{jb}u_{kc}\\
&=&\sum_{ijk}e_{ijk}\otimes\delta_{ik}u_{ia}u_{jb}u_{kc}
\end{eqnarray*}

On the other hand, we have as well the following formula:
\begin{eqnarray*}
u^{\otimes 3}T_\pi e_{abc}
&=&u^{\otimes 3}\delta_{ac}e_{abc}\\
&=&\sum_{ijk}e_{ijk}\otimes\delta_{ac}u_{ia}u_{jb}u_{kc}
\end{eqnarray*}

Thus the condition $T_\pi\in End(u^{\otimes 3})$ is equivalent to the following relations:
$$(\delta_{ik}-\delta_{ac})u_{ia}u_{jb}u_{kc}=0$$

The non-trivial cases are $i=k,a\neq c$ and $i\neq k,a=c$, and these produce the relations $u_{ia}u_{jb}u_{ic}=0$ for any $a\neq c$, and $u_{ia}u_{jb}u_{ka}=0$, for any $i\neq k$. Thus, we have reached to the standard relations for the quantum group $H_N^{[\infty]}$.

\medskip

(1) We have the following formula:
$$
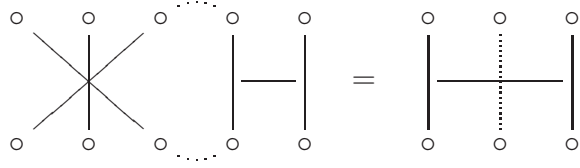
\xymatrix@R=5mm@C=5mm{
\circ\ar@{-}[ddrr]&\circ\ar@{-}[dd]&\circ\ar@{-}[ddll]\ar@/^/@{.}[r]&\circ\ar@{-}[dd]&\circ\ar@{-}[dd]\\
&&&\ar@{-}[r]&\\
\circ&\circ&\circ\ar@/_/@{.}[r]&\circ&\circ}\ \ \xymatrix@R=5mm@C=1mm{&\\=\\&}\xymatrix@R=5mm@C=5mm{
\circ\ar@{-}[dd]&\circ\ar@{.}[dd]&\circ\ar@{-}[dd]\\
\ar@{-}[rr]&&\\
\circ&\circ&\circ}$$

We have as well the following formula:
$$
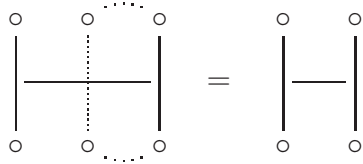
\xymatrix@R=5mm@C=5mm{
\circ\ar@{-}[dd]&\circ\ar@{.}[dd]\ar@/^/@{.}[r]&\circ\ar@{-}[dd]\\
\ar@{-}[rr]&&\\
\circ&\circ\ar@/_/@{.}[r]&\circ}\ \ \xymatrix@R=5mm@C=1mm{&\\=\\&}\xymatrix@R=5mm@C=5mm{
\circ\ar@{-}[dd]&\circ\ar@{-}[dd]\\
\ar@{-}[r]&&\\
\circ&\circ}$$

Thus, we obtain inclusions as desired, namely:
$$H_N^*\subset H_N^{[\infty]}\subset H_N^+$$

(2) At $r=2$, the relations $ab_1b_2c=0$ come indeed from the following diagram:
$$
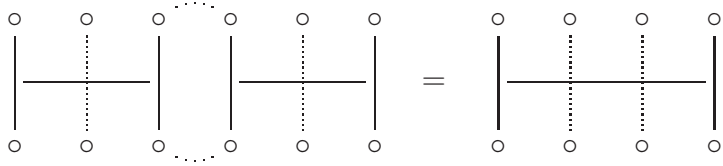
\xymatrix@R=5mm@C=5mm{
\circ\ar@{-}[dd]&\circ\ar@{.}[dd]&\circ\ar@{-}[dd]\ar@/^/@{.}[r]&\circ\ar@{-}[dd]&\circ\ar@{.}[dd]&\circ\ar@{-}[dd]\\
\ar@{-}[rr]&&&\ar@{-}[rr]&&\\
\circ&\circ&\circ\ar@/_/@{.}[r]&\circ&\circ&\circ}
\ \ \xymatrix@R=5mm@C=1mm{&\\=\\&}
\xymatrix@R=5mm@C=5mm{
\circ\ar@{-}[dd]&\circ\ar@{.}[dd]&\circ\ar@{.}[dd]&\circ\ar@{-}[dd]\\
\ar@{-}[rrr]&&&\\
\circ&\circ&\circ&\circ}$$

In the general case $r\geq2$ the proof is similar, see \cite{bcs} for details.

\medskip

(3) We use here an idea from \cite{rwe}. By rotating $\pi$, we obtain:
$$
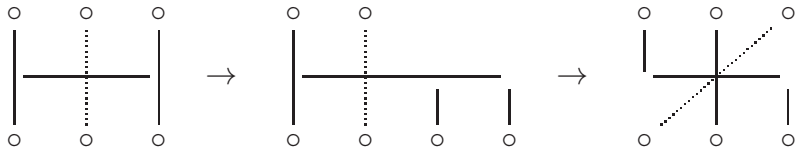
\xymatrix@R=5mm@C=5mm{
\circ\ar@{-}[dd]&\circ\ar@{.}[dd]&\circ\ar@{-}[dd]\\
\ar@{-}[rr]&&\\
\circ&\circ&\circ}
\ \ \xymatrix@R=5mm@C=0.1mm{&\\ \to\\&}\ 
\xymatrix@R=5mm@C=5mm{
\circ\ar@{-}[dd]&\circ\ar@{.}[dd]&\\
\ar@{-}[rrr]&&\ar@{-}[d]&\ar@{-}[d]\\
\circ&\circ&\circ&\circ}
\ \ \xymatrix@R=5mm@C=0.1mm{&\\ \to\\&}\ 
\xymatrix@R=5mm@C=5mm{
\circ\ar@{-}[d]&\circ\ar@{-}[dd]&\circ\ar@{.}[ddll]\\
\ar@{-}[rr]&&\ar@{-}[d]\\
\circ&\circ&\circ}$$

Let us denote by $\sigma$ the partition on the right. Since $T_\sigma(e_{ijk})=\delta_{ij}e_{kji}$, we obtain:
\begin{eqnarray*}
T_\sigma u^{\otimes 3}e_{abc}
&=&T_\sigma\sum_{ijk}e_{ijk}\otimes u_{ia}u_{jb}u_{kc}\\
&=&\sum_{ijk}e_{kji}\otimes\delta_{ij}u_{ia}u_{jb}u_{kc}
\end{eqnarray*}

On the other hand, we obtain as well the following formula:
\begin{eqnarray*}
u^{\otimes 3}T_\sigma e_{abc}
&=&u^{\otimes 3}\delta_{ab}e_{cba}\\
&=&\sum_{ijk}e_{kji}\otimes\delta_{ab}u_{kc}u_{jb}u_{ia}
\end{eqnarray*}

Thus $T_\sigma\in End(u^{\otimes 3})$ is equivalent to the following relations:
$$\delta_{ij}u_{ia}u_{jb}u_{kc}=\delta_{ab}u_{kc}u_{jb}u_{ia}$$

Now by setting $j=i,b=a$ in this formula we obtain the commutation relations in the statement, namely $u_{ia}^2u_{kc}=u_{kc}u_{ia}^2$, and this finishes the proof.
\end{proof}

The relation of $H_N^{[\infty]}$ with the polygonal spheres comes from the following fact:

\begin{proposition}
Let $X\subset S^{N-1}_{\mathbb R,+}$ be closed, let $d\geq2$, and set: $$X^{d-1}=X\cap S^{N-1,d-1}_{\mathbb R,+}$$
Then for a quantum group $G\subset H_N^{[\infty]}$ the following are equivalent:
\begin{enumerate}
\item $x_i\to\sum_jx_j\otimes u_{ji}$ defines a coaction map, as follows:
$$\Phi:C(X^{d-1})\to C(X^{d-1})\otimes C(G)$$

\item $x_i\to\sum_jx_j\otimes u_{ji}$ defines a morphism, as follows:
$$\widetilde{\Phi}:C(X)\to C(X^{d-1})\otimes C(G)$$
\end{enumerate} 
In particular, $G^+(X)\cap H_N^{[\infty]}$ acts on $X^{d-1}$, for any $d\geq2$.
\end{proposition}

\begin{proof}
The idea here is to use the relations in Proposition 14.6 (2):

\medskip

$(1)\implies(2)$ This is clear, by composing $\Phi$ with the following projection map:
$$\pi:C(X)\to C(X^{d-1})$$

$(2)\implies(1)$ In order to prove this implication, we must understand when a coaction map as follows exists:
$$C(X^{d-1})\to C(X^{d-1})\otimes C(G)$$

In order for this to happen, the variables $X_i=\sum_jx_j\otimes u_{ji}$ must satisfy the relations defining $X$, which hold indeed by (2), and must satisfy as well the following relations, with $i_0,\ldots,i_d$ distinct, which define the polygonal spheres $S^{N-1,d-1}_{\mathbb R,+}$:
$$X_{i_0}\ldots X_{i_d}=0$$ 

The point now is that, under the assumption $G\subset H_N^{[\infty]}$, these latter relations are automatic. Indeed, by using Proposition 14.6 (2), for $i_0,\ldots,i_d$ distinct we obtain:
\begin{eqnarray*}
X_{i_0}\ldots X_{i_d}
&=&\sum_{j_0\ldots j_d}x_{j_0}\ldots x_{j_d}\otimes u_{j_0i_0}\ldots u_{j_di_d}\\
&=&\sum_{j_0\ldots j_d\ distinct}0\otimes u_{j_0i_0}\ldots u_{j_di_d}\ \ +\!\!\!\!\sum_{j_0\ldots j_d\ not\ distinct}x_{j_0}\ldots x_{j_d}\otimes0\\
&=&0+0=0
\end{eqnarray*}

Thus the coaction in (1) exists precisely when (2) is satisfied, and we are done.

\medskip

Finally, the last assertion is clear from $(2)\implies(1)$, because the universal coaction of $G=G^+(X)$ gives rise to a map $\widetilde{\Phi}:C(X)\to C(X^{d-1})\otimes C(G)$ as in (2).
\end{proof}

As an illustration, we have the following result:

\begin{theorem}
The compact quantum groups
$$H_N,H_N,H_N^*,H_N^*,H_N^{[\infty]}$$
act respectively on the spheres
$$S^{N-1,d-1}_{\mathbb R},\bar{S}^{N-1,d-1}_{\mathbb R},S^{N-1,d-1}_{\mathbb R,*},\bar{S}^{N-1,d-1}_{\mathbb R,*},S^{N-1,d-1}_{\mathbb R,+}$$ 
at any $d\geq2$.
\end{theorem}

\begin{proof}
We use Proposition 14.7. We know that the quantum isometry groups at $d=N$ are respectively equal to the following quantum groups:
$$O_N,\bar{O}_N,O_N^*,\bar{O}_N^*,O_N^+$$

Our claim is that, by intersecting these quantum groups with $H_N^{[\infty]}$, we obtain the quantum groups in the statement. Indeed:

\medskip

(1) $O_N\cap H_N^{[\infty]}=H_N$ is clear from definitions.

\medskip

(2) $\bar{O}_N\cap H_N^{[\infty]}=H_N$ follows from $\bar{O}_N\cap H_N^+\subset O_N$, which in turn follows from the computation (3) in the proof of Proposition 14.5, with $H_N^*$ replaced by $H_N^+$.

\medskip

(3) $O_N^*\cap H_N^{[\infty]}=H_N^*$ follows from $O_N^*\cap H_N^+=H_N^*$.

\medskip

(4) $\bar{O}_N^*\cap H_N^{[\infty]}\supset H_N^*$ is clear, and the reverse inclusion can be proved by a direct computation, similar to the computation (3) in the proof of Proposition 14.5.

\medskip

(5) $O_N^+\cap H_N^{[\infty]}=H_N^{[\infty]}$ is clear from definitions.
\end{proof}

Observe that the above result is ``sharp'', in the sense that the actions there are the universal ones, in the classical case at any $d\in\{2,\ldots,N-1\}$, as well as in the twisted case at $d=2$. Indeed, this follows from the various results established above. We refer to \cite{ba2} for more comments on all this, and for some related open problems.

\section*{14c. Middle computation}

Let us discuss now the computation of the quantum isometry group for the polygonal sphere which is left, namely $S^{N-1,1}_{\mathbb R,*}$, appearing in the middle of our square diagram of polygonal spheres, from the beginning of this chapter. We have seen in the above that the quantum group $H_N^*$ acts on $S^{N-1,1}_{\mathbb R,*}$. This action, however, is not universal, due to:

\begin{proposition}
The discrete group dual
$$G=\widehat{\mathbb Z_2^{*N}}$$
acts on the polygonal sphere $S^{N-1,1}_{\mathbb R,*}$.
\end{proposition}

\begin{proof}
The standard coordinates on the polygonal sphere $S^{N-1,1}_{\mathbb R,*}$ are by definition subject to the following relations:
$$x_ix_jx_k=\begin{cases}
0&{\rm for}\ i,j,k\ {\rm distinct}\\
x_kx_jx_i&{\rm otherwise}
\end{cases}$$

Let us try now to understand under which assumptions on a compact quantum group $G$, and in particular on a group dual, we have a coaction map, as follows:
$$\Phi:C(S^{N-1,1}_{\mathbb R,*})\to C(S^{N-1,1}_{\mathbb R,*})\otimes C(G)$$
$$\Phi(x_i)=\sum_jx_j\otimes u_{ji}$$

In order for this to happen, the following variables must satisfy the above relations:
$$X_i=\sum_jx_j\otimes u_{ji}$$

For the group dual $G=\widehat{\mathbb Z_2^{*N}}$ we have by definition $u_{ij}=\delta_{ij}g_i$, where $g_1,\ldots,g_N$ are the standard generators of $\mathbb Z_2^{*N}$. We therefore have:
$$X_iX_jX_k=x_ix_jx_k\otimes g_ig_jg_k$$
$$X_kX_jX_i=x_kx_jx_i\otimes g_kg_jg_i$$
 
Thus the formula $X_iX_kX_k=0$ for $i,j,k$ distinct is clear, and the formula $X_iX_jX_k=X_kX_jX_i$ for $i,j,k$ not distinct requires $g_ig_jg_k=g_kg_jg_i$ for $i,j,k$ not distinct, which is clear as well. Indeed, at $i=j$ this latter relation reduces to $g_k=g_k$, at $i=k$ this relation is trivial, $g_ig_jg_i=g_ig_jg_i$, and at $j=k$ this relation reduces to $g_i=g_i$.
\end{proof}

More generally now, we have the following result:

\begin{proposition}
The intermediate liberation of the hyperoctahedral group
$$G=H_N^{[\infty]}$$
acts on the polygonal sphere $S^{N-1,1}_{\mathbb R,*}$.
\end{proposition}

\begin{proof}
We proceed as in the proof of Theorem 14.4. By expanding the formula of $X_iX_jX_k$ and by using the relations for the sphere $S^{N-1,1}_{\mathbb R,*}$, we have:
\begin{eqnarray*}
X_iX_jX_k
&=&\sum_{abc}x_ax_bx_c\otimes u_{ai}u_{bj}u_{ck}\\
&=&\sum_{a,b,c\ not \ distinct}x_ax_bx_c\otimes u_{ai}u_{bj}u_{ck}\\
&=&\sum_{a\neq b}x_a^2x_b\otimes(u_{ai}u_{aj}u_{bk}+u_{bi}u_{aj}u_{ak})\\
&+&\sum_{a\neq b}x_ax_bx_a\otimes u_{ai}u_{bj}u_{ak}\\
&+&\sum_a x_a^3\otimes u_{ai}u_{aj}u_{ak}
\end{eqnarray*}

Now by assuming $G=H_N^{[\infty]}$, as in the statement, and by using the various formulae in Proposition 14.6, we obtain, for any $i,j,k$ distinct:
\begin{eqnarray*}
X_iX_jX_k
&=&\sum_{a\neq b}x_a^2x_b\otimes(0\cdot u_{bk}+u_{bi}\cdot 0)\\
&+&\sum_{a\neq b}x_ax_bx_a\otimes0\\
&+&\sum_ax_a^3\otimes(0\cdot u_{ak})\\
&=&0
\end{eqnarray*}

It remains to prove that we have $X_iX_jX_k=X_kX_jX_i$, for $i,j,k$ not distinct. By replacing $i\leftrightarrow k$ in the above formula of $X_iX_jX_k$, we obtain:
\begin{eqnarray*}
X_kX_jX_i
&=&\sum_{a\neq b}x_a^2x_b\otimes(u_{ak}u_{aj}u_{bi}+u_{bk}u_{aj}u_{ai})\\
&+&\sum_{a\neq b}x_ax_bx_a\otimes u_{ak}u_{bj}u_{ai}\\
&+&\sum_ax_a^3\otimes u_{ak}u_{aj}u_{ai}
\end{eqnarray*}

Let us compare this formula with the above formula of $X_iX_jX_k$. The last sum being 0 in both cases, we must prove that for any $i,j,k$ not distinct and any $a\neq b$ we have:
$$u_{ai}u_{aj}u_{bk}+u_{bi}u_{aj}u_{ak}=u_{ak}u_{aj}u_{bi}+u_{bk}u_{aj}u_{ai}$$
$$u_{ai}u_{bj}u_{ak}=u_{ak}u_{bj}u_{ai}$$

By symmetry the three cases $i=j,i=k,j=k$ reduce to two cases, $i=j$ and $i=k$. The case $i=k$ being clear, we are left with the case $i=j$, where we must prove:
$$u_{ai}u_{ai}u_{bk}+u_{bi}u_{ai}u_{ak}=u_{ak}u_{ai}u_{bi}+u_{bk}u_{ai}u_{ai}$$
$$u_{ai}u_{bi}u_{ak}=u_{ak}u_{bi}u_{ai}$$

By using $a\neq b$, the first equality reads:
$$u_{ai}^2u_{bk}+0\cdot u_{ak}=u_{ak}\cdot 0+u_{bk}u_{ai}^2$$

Since by Proposition 14.6 (3) we have $u_{ai}^2u_{bk}=u_{bk}u_{ai}^2$, this formula is satisfied, and we are done. As for the second equality, this reads:
$$0\cdot u_{ak}=u_{ak}\cdot 0$$

But this is true as well, and this ends the proof.
\end{proof}

We will prove now that the action in Proposition 14.10 is universal. In order to do so, we need to convert the formulae of type $X_iX_jX_k=0$ and $X_iX_jX_k=X_kX_jX_i$ into relations between the quantum group coordinates $u_{ij}$, and this requires a good knowledge of the linear relations between the variables $x_a^2x_b,x_ax_bx_a,x_a^3$ over the sphere $S^{N-1,1}_{\mathbb R,*}$. So, we must first study these variables. The answer here is given by:

\begin{proposition}
The variables 
$$\left\{x_a^2x_b,x_ax_bx_a,x_a^3\Big|a\neq b\right\}$$
are linearly independent over the sphere $S^{N-1,1}_{\mathbb R,*}$.
\end{proposition}

\begin{proof}
We use a trick of Bichon-Dubois-Violette \cite{bdu}. Consider the 1-dimensional polygonal version of the complex sphere $S^{N-1}_\mathbb C$, which is by definition given by:
$$S^{N-1,1}_\mathbb C=\left\{z\in S^{N-1}_\mathbb C\Big|z_iz_jz_k=0,\forall i,j,k\ {\rm distinct}\right\}$$

We have then a $2\times2$ matrix model for the coordinates of $S^{N-1,1}_{\mathbb R,*}$, as follows:
$$x_i\to\gamma_i=\begin{pmatrix}0&z_i\\ \bar{z}_i&0\end{pmatrix}$$

Indeed, the matrices $\gamma_i$ on the right are all self-adjoint, their squares sum up to $1$, they half-commute, and they satisfy, for $i,j,k$ distinct: 
$$\gamma_i\gamma_j\gamma_k=0$$

Thus we have indeed a morphism as follows, as claimed:
$$C(S^{N-1,1}_{\mathbb R,*})\to M_2(C(S^{N-1,1}_\mathbb C))\quad,\quad x_i\to\gamma_i$$

We can use this model in order to prove the linear independence. Consider indeed the variables that we are interested in, namely:
$$x_a^2x_b\quad,\quad 
x_ax_bx_a\quad,\quad 
x_a^3$$

In the model, these variables are mapped to the following variables:
$$\gamma_a^2\gamma_b=\begin{pmatrix}0&|z_a|^2z_b\\ |z_a|^2\bar{z}_b&0\end{pmatrix}$$
$$\gamma_a\gamma_b\gamma_a=\begin{pmatrix}0&z_a^2\bar{z}_b\\ \bar{z}_a^2z_b&0\end{pmatrix}$$
$$\gamma_a^3=\begin{pmatrix}0&|z_a|^2z_a\\ |z_a|^2\bar{z}_a&0\end{pmatrix}$$

Now observe that the following variables are linearly independent over $S^1_\mathbb C$:
$$|z_1|^2z_2,|z_2|^2z_1,z_1^2\bar{z}_2,z_2^2\bar{z}_1,|z_1|^2z_1,|z_2|^2z_2$$

Thus the upper right entries of the above matrices are linearly independent over $S^{N-1,1}_\mathbb C$. Thus the matrices themselves are linearly independent, and this proves our result.
\end{proof}

With the above result in hand, we can now reformulate the coaction problem into a purely quantum group-theoretical problem, as follows:

\begin{proposition}
A quantum group $G\subset O_N^+$ acts on $S^{N-1,1}_{\mathbb R,*}$ precisely when its standard coordinates $u_{ij}$ satisfy the following relations:
\begin{enumerate}
\item $u_{ai}u_{aj}u_{bk}+u_{bi}u_{aj}u_{ak}=0$ for any $i,j,k$ distinct.

\item $u_{ai}u_{bj}u_{ak}=0$ for any $i,j,k$ distinct.

\item $u_{ai}^2u_{bk}=u_{bk}u_{ai}^2$.

\item $u_{ak}u_{ai}u_{bi}=u_{bi}u_{ai}u_{ak}$.

\item $u_{ai}u_{bi}u_{ak}=u_{bk}u_{bi}u_{ai}$.
\end{enumerate}
\end{proposition}

\begin{proof}
We use notations from the beginning of the proof of Proposition 14.10, along with the following formula, also established there:
\begin{eqnarray*}
X_iX_jX_k
&=&\sum_{a\neq b}x_a^2x_b\otimes(u_{ai}u_{aj}u_{bk}+u_{bi}u_{aj}u_{ak})\\
&+&\sum_{a\neq b}x_ax_bx_a\otimes u_{ai}u_{bj}u_{ak}\\
&+&\sum_ax_a^3\otimes u_{ai}u_{aj}u_{ak}
\end{eqnarray*}

In order to have an action as in the statement, these quantities must satisfy $X_iX_kX_k=0$ for $i,j,k$ distinct, and $X_iX_kX_k=X_kX_jX_i$ for $i,j,k$ not distinct. Now by using Proposition 14.11, we conclude that the relations to be satisfied are as follows:

\medskip

(A) For $i,j,k$ distinct, the following must hold:
$$u_{ai}u_{aj}u_{bk}+u_{bi}u_{aj}u_{ak}=0,\forall a\neq b$$
$$u_{ai}u_{bj}u_{ak}=0,\forall a\neq b$$
$$u_{ai}u_{aj}u_{ak}=0,\forall a$$

(B) For $i,j,k$ not distinct, the following must hold:
$$u_{ai}u_{aj}u_{bk}+u_{bi}u_{aj}u_{ak}=u_{ak}u_{aj}u_{bi}+u_{bk}u_{aj}u_{ai},\forall a\neq b$$
$$u_{ai}u_{bj}u_{ak}=u_{ak}u_{bj}u_{ai},\forall a\neq b$$
$$u_{ai}u_{aj}u_{ak}=u_{ak}u_{aj}u_{ai},\forall a$$

The last relations in (A,B) can be merged with the other ones, and we are led to:

\medskip

(A') For $i,j,k$ distinct, the following must hold:
$$u_{ai}u_{aj}u_{bk}+u_{bi}u_{aj}u_{ak}=0,\forall a,b$$
$$u_{ai}u_{bj}u_{ak}=0,\forall a,b$$

(B') For $i,j,k$ not distinct, the following must hold:
$$u_{ai}u_{aj}u_{bk}+u_{bi}u_{aj}u_{ak}=u_{ak}u_{aj}u_{bi}+u_{bk}u_{aj}u_{ai},\forall a,b$$
$$u_{ai}u_{bj}u_{ak}=u_{ak}u_{bj}u_{ai},\forall a,b$$

Observe that the relations (A') are exactly the relations (1,2) in the statement.

\medskip

Let us further process the relations (B'). In the case $i=k$ the relations are automatic, and in the cases $j=i,j=k$ the relations that we obtain coincide, via $i\leftrightarrow k$. Thus (B') reduces to the set of relations obtained by setting $j=i$, which are as follows:
$$u_{ai}u_{ai}u_{bk}+u_{bi}u_{ai}u_{ak}=u_{ak}u_{ai}u_{bi}+u_{bk}u_{ai}u_{ai}$$
$$u_{ai}u_{bi}u_{ak}=u_{ak}u_{bi}u_{ai}$$

Observe that the second relation is the relation (5) in the statement. Regarding now the first relation, with the notation $[x,y,z]=xyz-zyx$, this is as follows:
$$[u_{ai},u_{ai},u_{bk}]=[u_{ak},u_{ai},u_{bi}]$$

By applying the antipode, we obtain from this:
$$[u_{kb},u_{ia},u_{ia}]=[u_{ib},u_{ia},u_{ka}]$$

By relabelling $a\leftrightarrow i$ and $b\leftrightarrow k$, this relation becomes:
$$[u_{bk},u_{ai},u_{ai}]=[u_{ak},u_{ai},u_{bi}]$$

Now since we have $[x,y,z]=-[z,y,x]$, by comparing this latter relation with the original one, a simplification occurs, and the resulting relations are as follows:
$$[u_{ai},u_{ai},u_{bk}]=[u_{ak},u_{ai},u_{bi}]=0$$

But these are exactly the relations (3,4) in the statement, and we are done.
\end{proof}

Now by solving the quantum group problem raised by Proposition 14.12, we obtain:

\begin{proposition}
We have the following formula,
$$G^+(S^{N-1,1}_{\mathbb R,*})=H_N^{[\infty]}$$
with $H_N^{[\infty]}$ being the quantum group constructed above.
\end{proposition}

\begin{proof}
The inclusion $\supset$ is clear from Proposition 14.10. For the converse, we already have the result at $N=2$, so assume $N\geq3$. We will use many times the conditions (1-5) in Proposition 14.12. By using (2), for $i\neq j$ we have:
\begin{eqnarray*}
&&u_{ai}u_{bj}u_{ak}=0,\forall k\neq i,j\\
&\implies&u_{ai}u_{bj}u_{ak}^2=0,\forall k\neq i,j\\
&\implies&u_{ai}u_{bj}\left(\sum_{k\neq i,j}u_{ak}^2\right)=0,\forall i\neq j\\
&\implies&u_{ai}u_{bj}(1-u_{ai}^2-u_{aj}^2)=0,\forall i\neq j
\end{eqnarray*}

Now by using (3), we can move the variable $u_{bj}$ to the right. By further multiplying by $u_{bj}$ to the right, and then summing over $b$, we obtain:
\begin{eqnarray*}
&&u_{ai}u_{bj}(1-u_{ai}^2-u_{aj}^2)=0,\forall i\neq j\\
&\implies&u_{ai}(1-u_{ai}^2-u_{aj}^2)u_{bj}=0,\forall i\neq j\\
&\implies&u_{ai}(1-u_{ai}^2-u_{aj}^2)u_{bj}^2=0,\forall i\neq j\\
&\implies&u_{ai}(1-u_{ai}^2-u_{aj}^2)=0,\forall i\neq j
\end{eqnarray*}

We can proceed now as follows, by summing over $j\neq i$:
\begin{eqnarray*}
&&u_{ia}(1-u_{ai}^2-u_{aj}^2)=0,\forall i\neq j\\
&\implies&u_{ai}u_{aj}^2=u_{ai}-u_{ai}^3,\forall i\neq j\\
&\implies&u_{ai}(1-u_{ai}^2)=(N-1)(u_{ai}-u_{ai}^3)\\
&\implies&u_{ai}=u_{ai}^3
\end{eqnarray*}

Thus the standard coordinates are partial isometries, and so $G\subset H_N^+$. On the other hand, we know from the proof of Proposition 14.6 (3) that the quantum subgroup $G\subset H_N^+$ obtained via the relations $[a,b^2]=0$ is $H_N^{[\infty]}$, and this finishes the proof.
\end{proof}

We have now complete results for the 9 main spheres, as follows:

\begin{theorem}
The quantum isometry groups of the $9$ polygonal spheres are
$$\xymatrix@R=12mm@C=17mm{
O_N\ar[r]&O_N^*\ar[r]&O_N^+\\
H_N\ar[r]\ar[u]&H_N^{[\infty]}\ar[r]\ar[u]&\bar{O}_N^*\ar[u]\\
H_N^+\ar[r]\ar[u]&H_N\ar[r]\ar[u]&\bar{O}_N\ar[u]}$$
where $H_N^+,H_N^{[\infty]}$ and $\bar{O}_N,O_N^*,\bar{O}_N^*,O_N^*$ are quantum versions of $H_N,O_N$.
\end{theorem}

\begin{proof}
This follows indeed by putting together the above results.
\end{proof}

As a conclusion, we have a key computation available, but there are still many questions left, regarding the extension of our $(S,T,U,K)$ formalism, as to cover the intersections between the twisted and untwisted geometries. And here, things are quite complicated, because as explained in chapter 11, even in the plain untwisted case, there is one axiom which needs to be modified, with questions which are currently open.

\bigskip

As explained in the present and in the previous chapter, on several occasions, all this is not exactly something theoretical, but is rather something of practical interest, in connection with the algebraic geometry of the free spheres, and other free manifolds. Finally, let us mention too that none of the questions to be solved is difficult, so to say. It's just that there is so much work to be done, and that hasn't been done yet.

\section*{14d. Complex extension}

Still following \cite{ba2}, let us discuss now a straightforward complex extension of the above results. Our starting point will be the following definition:

\index{complex polygonal sphere}
\index{intersection of complex spheres}

\begin{definition}
The complex polygonal spheres, denoted 
$$S^{N-1,d-1}_{\mathbb C},\bar{S}^{N-1,d-1}_{\mathbb C},\bar{S}^{N-1,d-1}_{\mathbb C,*},S^{N-1,d-1}_{\mathbb C,+}$$
are constructed from $S^{N-1}_{\mathbb C,+}$ in the same way as their real versions, namely 
$$S^{N-1,d-1}_{\mathbb R},\bar{S}^{N-1,d-1}_{\mathbb R},\bar{S}^{N-1,d-1}_{\mathbb R,*},S^{N-1,d-1}_{\mathbb R,+}$$
are constructed from $S^{N-1}_{\mathbb R,+}$, namely by assuming that the corresponding vanishing relations hold between the variables $x_i=z_i,z_i^*$.
\end{definition}

As in the real case, we will restrict now the attention to the 5 main spheres, coming from \cite{ba1}, and to their intersections. We have 9 such spheres here, as follows:
$$\xymatrix@R=13mm@C=3mm{
S^{N-1}_\mathbb C\ar[rr]&&S^{N-1}_{\mathbb C,*}\ar[rrr]&&&S^{N-1}_{\mathbb C,+}\\
S^{N-1}_\mathbb C\cap\bar{S}^{N-1}_{\mathbb C,*}\ar[rr]\ar[u]&&S^{N-1}_{\mathbb C,*}\cap\bar{S}^{N-1}_{\mathbb C,*}\ar[rrr]\ar[u]&&&\bar{S}^{N-1}_{\mathbb C,*}\ar[u]\\
S^{N-1}_\mathbb C\cap\bar{S}^{N-1}_\mathbb C\ar[rr]\ar[u]&&S^{N-1}_{\mathbb C,*}\cap\bar{S}^{N-1}_\mathbb C\ar[rrr]\ar[u]&&&\bar{S}^{N-1}_\mathbb C\ar[u]}$$

The intersections can be computed as in the real case, and we have:

\begin{proposition}
The $5$ main spheres, and the intersections between them, are
$$\xymatrix@R=12mm@C=12mm{
S^{N-1}_\mathbb C\ar[r]&S^{N-1}_{\mathbb C,*}\ar[r]&S^{N-1}_{\mathbb C,+}\\
S^{N-1,1}_\mathbb C\ar[r]\ar[u]&S^{N-1,1}_{\mathbb C,*}\ar[r]\ar[u]&\bar{S}^{N-1}_{\mathbb C,*}\ar[u]\\
S^{N-1,0}_\mathbb C\ar[r]\ar[u]&\bar{S}^{N-1,1}_\mathbb C\ar[r]\ar[u]&\bar{S}^{N-1}_\mathbb C\ar[u]}$$
with all the maps being inclusions.
\end{proposition}

\begin{proof}
This is similar to the proof from the real case, from chapter 13, by replacing in all the computations there the variables $x_i$ by the variables $x_i=z_i,z_i^*$. For full details here, and for more on these spheres, we refer as before to the paper \cite{ba2}.
\end{proof}

As before in the real case, all this raises the question of updating our noncommutative geometry axioms from chapter 4, as to cover such intersections, which are quite interesting objects. And here, we are mostly in need of quantum isometry group results. The computation, from \cite{ba2}, is very similar to the one from the real case, the result being:

\begin{theorem}
The quantum isometry groups of the $9$ main complex spheres are
$$\xymatrix@R=13mm@C=17mm{
U_N\ar[r]&U_N^*\ar[r]&U_N^+\\
K_N\ar[r]\ar[u]&K_N^{[\infty]}\ar[r]\ar[u]&\bar{U}_N^*\ar[u]\\
K_N^+\ar[r]\ar[u]&K_N\ar[r]\ar[u]&\bar{U}_N\ar[u]}$$
where $K_N$ and its versions are the complex analogues of $H_N$ and its versions.
\end{theorem}

\begin{proof}
The idea here is that the proof here is quite similar to the proof in the real case, by replacing $H_N,O_N$ with their complex analogues $K_N,U_N$. For the details here, which are quite straightforward, and for further coments on this, we refer to \cite{ba2}.
\end{proof}

As a conclusion, we have many technical results available, but there are still many questions left, regarding the extension of our $(S,T,U,K)$ formalism, as to cover the intersections between the twisted and untwisted geometries. 

\bigskip

So long for free algebraic geometry. These things are quite recent, going back to the 10s, and there is still an enormous amount of work to be done, in order to have something serious started. And what is quite puzzling here, but not really puzzling if you're a bit familiar with algebraic geometry, and its story, is that everything is somehow ``trivial'', but the massive accumulation of such trivialities eventually leads to difficulties.

\bigskip

Which can only make us think at Grothendieck, the problems that he met when doing algebraic geometry, and his strategy for dealing with them. The story goes that if Grothendieck was to be on top of a mountain, and had to go to the mountain next to it, he would rather fill the valley in between with small rocks, no matter how long it takes, and walk straight, instead of going downhill and then uphill.

\bigskip

But hey, there should be a trick somewhere, in order to avoid such Grothendieck type things, the mathematics that we developed so far is dangerously pushing us into. Up to you to find it, either based on other things discussed in this book, or on some new idea of yours. Plus of course, and obviously, if you feel like doing some Grothendieck type work in free algebraic geometry, just go for that. No fuss, only trivialities, can only work.

\section*{14e. Exercises} 

As with the exercises from the previous chapter, that the present chapter continuates, our exercises here will be quite technical, and algebraic. Let us start with:

\begin{exercise}
Compute the quantum isometry groups of the twisted polygonal spheres, which are by definition given by
$$C(\bar{S}^{N-1,d-1}_\mathbb R)=C(\bar{S}^{N-1}_\mathbb R)\Big/\Big<x_{i_0}\ldots x_{i_d}=0,\forall i_0,\ldots,i_d\ {\rm distinct}\Big>$$ 
at the missing values of the parameter, $d=3,4,\ldots,N-1$.
\end{exercise}

This is something that we discussed in the above, and is a result that will be nice to have. Technically speaking, however, it is not clear how exactly to do this.

\begin{exercise}
Work out the full theory of the quantum group $H_N^{[\infty]}$, notably with a Brauer type result for it, and then with probability computations.
\end{exercise}

This exercise is interesting not only because of the considerations from the present chapter, where the quantum group $H_N^{[\infty]}$ plays a central role, but also in view of the classification discussion in chapter 10, where $H_N^{[\infty]}$ was playing as well a central role.

\begin{exercise}
Work out the full theory of the quantum group $K_N^{[\infty]}$, notably with a Brauer type result for it, and then with probability computations.
\end{exercise}

This exercise is particularly interesting, because unlike $H_N^{[\infty]}$, which is a well-known quantum group, $K_N^{[\infty]}$ is a more recent object, not systematically studied yet.

\begin{exercise}
Find suitable quantum group axioms covering the quantum isometry groups of the main polygonal spheres, namely
$$\xymatrix@R=11mm@C=15mm{
O_N\ar[r]&O_N^*\ar[r]&O_N^+\\
H_N\ar[r]\ar[u]&H_N^{[\infty]}\ar[r]\ar[u]&\bar{O}_N^*\ar[u]\\
H_N^+\ar[r]\ar[u]&H_N\ar[r]\ar[u]&\bar{O}_N\ar[u]}$$
and then do the same in the complex case.
\end{exercise}

The above quantum groups are all quizzy, and the problem is that of looking more carefully, in order to see what are the common features of these quantum groups.

\chapter{Projective geometry}

\section*{15a. Projective spaces}

This chapter will be an introduction to projective geometry, in our sense. The point is that things become considerably simpler in the projective setting. Consider indeed the diagram of 9 main geometries, that we found above:
$$\xymatrix@R=40pt@C=42pt{
\mathbb R^N_+\ar[r]&\mathbb T\mathbb R^N_+\ar[r]&\mathbb C^N_+\\
\mathbb R^N_*\ar[u]\ar[r]&\mathbb T\mathbb R^N_*\ar[u]\ar[r]&\mathbb C^N_*\ar[u]\\
\mathbb R^N\ar[u]\ar[r]&\mathbb T\mathbb R^N\ar[u]\ar[r]&\mathbb C^N\ar[u]
}$$

As explained in chapters 9-10, when looking at the projective versions of these geometries,  the diagram drastically simplifies. To be more precise, the diagram of projective versions of the corresponding spheres are as follows, consisting of 3 objects only:
$$\xymatrix@R=34pt@C=33pt{
P^{N-1}_+\ar[r]&P^{N-1}_+\ar[r]&P^{N-1}_+\\
P^{N-1}_\mathbb C\ar[r]\ar[u]&P^{N-1}_\mathbb C\ar[r]\ar[u]&P^{N-1}_\mathbb C\ar[u]\\
P^{N-1}_\mathbb R\ar[r]\ar[u]&P^{N-1}_\mathbb R\ar[r]\ar[u]&P^{N-1}_\mathbb C\ar[u]}$$

Thus, we are led to the conclusion that, under certain combinatorial axioms, there should be only 3 projective geometries, namely the real, complex and free ones:
$$P^{N-1}_\mathbb R\subset P^{N-1}_\mathbb C\subset P^{N-1}_+$$

We will discuss this in what follows, with analogues and improvements of the affine results. Also, we would like to study the corresponding quadruplets $(P,PT,PU,PK)$, and to axiomatize the projective geometries, with correspondences as follows:
$$
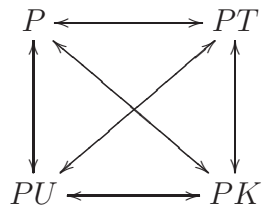
\xymatrix@R=50pt@C=50pt{
P\ar[r]\ar[d]\ar[dr]&PT\ar[l]\ar[d]\ar[dl]\\
PU\ar[u]\ar[ur]\ar[r]&PK\ar[l]\ar[ul]\ar[u]
}$$

Summarizing, there is a lot of work to be done, on one hand in reformulating and improving the results from the affine case, and on the other hand, in starting to develop the projective theory independently from the affine theory.

\bigskip

As a first topic that we would like to discuss, which historically speaking, was at the beginning of everything, we have the following remarkable isomorphism, that we already used in the above, and that we would like to discuss now in detail:
$$PO_N^+=PU_N^+$$

In order to get started, let us first discuss the classical case, and more specifically the precise relation between the orthogonal group $O_N$, and the unitary group $U_N$. Contrary to the passage $\mathbb R^N\to\mathbb C^N$, or to the passage $S^{N-1}_\mathbb R\to S^{N-1}_\mathbb C$, which are both elementary, the passage $O_N\to U_N$ cannot be understood directly. In order to understand this passage we must pass through the corresponding Lie algebras, a follows:

\index{complexification}
\index{Lie algebra complexification}

\begin{theorem}
The passage $O_N\to U_N$ appears via Lie algebra complexification,
$$O_N\to\mathfrak o_N\to\mathfrak u_n\to U_N$$
with the Lie algebra $\mathfrak u_N$ being a complexification of the Lie algebra $\mathfrak o_N$.
\end{theorem}

\begin{proof}
This is something rather philosophical, and advanced as well, that we will not really need here, the idea being as follows:

\medskip

(1) The unitary and orthogonal groups $U_N,O_N$ are both Lie groups, in the sense that they are smooth manifolds. The corresponding Lie algebras $\mathfrak u_N,\mathfrak o_N$, which are by definition the respective tangent spaces at 1, can be computed by differentiating the equations defining $U_N,O_N$, with the conclusion being as follows:
$$\mathfrak u_N=\left\{ A\in M_N(\mathbb C)\Big|A^*=-A\right\}$$
$$\mathfrak o_N=\left\{ B\in M_N(\mathbb R)\Big|B^t=-B\right\}$$

(2) This was for the correspondences $U_N\to\mathfrak u_N$ and $O_N\to\mathfrak o_N$. In the other sense, the correspondences $\mathfrak u_N\to U_N$ and $\mathfrak o_N\to O_N$ appear by exponentiation, the result here stating that, around 1, the unitary matrices can be written as $U=e^A$, with $A\in\mathfrak u_N$, and the orthogonal matrices can be written as $U=e^B$, with $B\in\mathfrak o_N$. 

\bigskip

(3) In view of all this, in order to understand the passage $O_N\to U_N$ it is enough to understand the passage $\mathfrak o_N\to\mathfrak u_N$. But, in view of the above formulae for $\mathfrak o_N,\mathfrak u_N$, this is basically an elementary linear algebra problem. Indeed, let us pick an arbitrary matrix $A\in M_N(\mathbb C)$, and write it as follows, with $B,C\in M_N(\mathbb R)$:
$$A=B+iC$$

In terms of $B,C$, the equation $A^*=-A$ defining the Lie algebra $\mathfrak u_N$ reads:
$$B^t=-B\quad,\quad 
C^t=C$$

(4) As a first observation, we must have $B\in\mathfrak o_N$. Regarding now $C$, let us decompose this matrix as follows, with $D$ being its diagonal, and $C'$ being the reminder:
$$C=D+C'$$

The matrix $C'$ being symmetric with 0 on the diagonal, by swithcing all the signs below the main diagonal we obtain a certain matrix $C'_-\in\mathfrak o_N$. Thus, we have decomposed $A\in\mathfrak u_N$ as follows, with $B,C'\in\mathfrak o_N$, and with $D\in M_N(\mathbb R)$ being diagonal:
$$A=B+iD+iC'_-$$

(5) As a conclusion now, we have shown that we have a direct sum decomposition of real linear spaces as follows, with $\Delta\subset M_N(\mathbb R)$ being the diagonal matrices:
$$\mathfrak u_N\simeq\mathfrak o_N\oplus\Delta\oplus\mathfrak o_N$$

Thus, we can stop our study here, and say that we have reached the conclusion in the statement, namely that $\mathfrak u_N$ appears as a ``complexification'' of $\mathfrak o_N$.
\end{proof}

As before with many other things, that we will not really need in what follows, this was just an introduction to the subject. More can be found in any Lie group book. In the free case now, the situation is much simpler, and we have:

\index{free complexification}
\index{free unitary group}
\index{free rotation}
\index{Haar unitary}

\begin{theorem}
The passage $O_N^+\to U_N^+$ appears via free complexification,
$$U_N^+=\widetilde{O_N^+}$$
where the free complexification of a pair $(G,u)$ is the pair $(\widetilde{G},\widetilde{u})$ with
$$C(\widetilde{G})=<zu_{ij}>\subset C(\mathbb T)*C(G)\quad,\quad 
\widetilde{u}=zu$$
where $z\in C(\mathbb T)$ is the standard generator, given by $x\to x$ for any $x\in\mathbb T$.
\end{theorem}

\begin{proof}
We have embeddings as follows, with the first one coming by using the counit, and with the second one coming from the universality property of $U_N^+$:
$$O_N^+
\subset\widetilde{O_N^+}
\subset U_N^+$$

We must prove that the embedding on the right is an isomorphism, and there are several ways of doing this, all instructive, as follows:

\medskip

(1) If we denote by $v,u$ the fundamental corepresentations of $O_N^+,U_N^+$, we have:
$$Fix(v^{\otimes k})=span\left(\xi_\pi\Big|\pi\in NC_2(k)\right)$$
$$Fix(u^{\otimes k})=span\left(\xi_\pi\Big|\pi\in\mathcal{NC}_2(k)\right)$$

Moreover, the above vectors $\xi_\pi$ are known to be linearly independent at $N\geq2$, and so the above results provide us with bases, and we obtain:
$$\dim(Fix(v^{\otimes k}))=|NC_2(k)|\quad,\quad 
\dim(Fix(u^{\otimes k}))=|\mathcal{NC}_2(k)|$$ 

Now since integrating the character of a corepresentation amounts in counting the fixed points, the above two formulae can be rewritten as follows:
$$\int_{O_N^+}\chi_v^k=|NC_2(k)|\quad,\quad 
\int_{U_N^+}\chi_u^k=|\mathcal{NC}_2(k)|$$ 

But this shows, via standard free probability theory, that $\chi_v$ must follow the Winger semicircle law $\gamma_1$, and that $\chi_u$ must follow the Voiculescu circular law $\Gamma_1$:
$$\chi_v\sim\gamma_1\quad,\quad 
\chi_u\sim\Gamma_1$$

On the other hand, by \cite{vdn}, when freely multiplying a semicircular variable by a Haar unitary we obtain a circular variable. Thus, the main character of $\widetilde{O_N^+}$ is circular:
$$\chi_{zv}\sim\Gamma_1$$

Now by forgetting about circular variables and free probability, the conclusion is that the inclusion $\widetilde{O_N^+}\subset U_N^+$ preserves the law of the main character:
$$law(\chi_{zv})=law(u)$$

Thus by Peter-Weyl we obtain that the inclusion $\widetilde{O_N^+}\subset U_N^+$ must be an isomorphism, modulo the usual equivalence relation for quantum groups.

\medskip

(2) A version of the above proof, not using any prior free probability knowledge, makes use of the easiness property of $O_N^+,U_N^+$ only, namely:
$$Hom(v^{\otimes k},v^{\otimes l})=span\left(\xi_\pi\Big|\pi\in NC_2(k,l)\right)$$
$$Hom(u^{\otimes k},u^{\otimes l})=span\left(\xi_\pi\Big|\pi\in\mathcal{NC}_2(k,l)\right)$$

Indeed, let us look at the following inclusions of quantum groups:
$$O_N^+\subset\widetilde{O_N^+}\subset U_N^+$$

At the level of the associated Hom spaces we obtain reverse inclusions, as follows:
$$Hom(v^{\otimes k},v^{\otimes l})
\supset Hom((zv)^{\otimes k},(zv)^{\otimes l})
\supset Hom(u^{\otimes k},u^{\otimes l})$$

The spaces on the left and on the right are known from easiness, the result being that these spaces are as follows:
$$span\left(T_\pi\Big|\pi\in NC_2(k,l)\right)\supset span\left(T_\pi\Big|\pi\in\mathcal{NC}_2(k,l)\right)$$

Regarding the spaces in the middle, these are obtained from those on the left by ``coloring'', so we obtain the same spaces as those on the right. Thus, by Tannakian duality, our embedding $\widetilde{O_N^+}\subset U_N^+$ is an isomorphism, modulo the usual equivalence relation.
\end{proof}

As an interesting consequence of the above result, we have:

\index{projective quantum group}
\index{projective orthogonal quantum group}
\index{projective unitary quantum group}

\begin{theorem}
We have an identification as follows,
$$PO_N^+=PU_N^+$$
modulo the usual equivalence relation for compact quantum groups.
\end{theorem}

\begin{proof}
As before, we have several proofs for this result, as follows:

\medskip

(1) This follows from Theorem 15.2, because we have:
$$PU_N^+=P\widetilde{O_N^+}=PO_N^+$$

(2) We can deduce this as well directly. With notations as before, we have:
$$Hom\left((v\otimes v)^k,(v\otimes v)^l\right)=span\left(T_\pi\Big|\pi\in NC_2((\circ\bullet)^k,(\circ\bullet)^l)\right)$$
$$Hom\left((u\otimes\bar{u})^k,(u\otimes\bar{u})^l\right)=span\left(T_\pi\Big|\pi\in \mathcal{NC}_2((\circ\bullet)^k,(\circ\bullet)^l)\right)$$

The sets on the right being equal, we conclude that the inclusion $PO_N^+\subset PU_N^+$ preserves the corresponding Tannakian categories, and so must be an isomorphism.
\end{proof}

As a conclusion, the passage $O_N^+\to U_N^+$ is something much simpler than the passage $O_N\to U_N$, with this ultimately coming from the fact that the combinatorics of $O_N^+,U_N^+$ is something much simpler than the combinatorics of $O_N,U_N$. In addition, all this leads as well to the interesting conclusion that the free projective geometry does not fall into real and complex, but is rather unique and ``scalarless''. We will be back to this.

\bigskip

Let us discuss now the projective spaces. We begin with a short summary of the various projective geometry results that we have so far. We will give full details here, with the aim of making the present chapter as independent as possible from the previous chapters, as a beginning of something new. Our starting point is the following functional analytic description of the real and complex projective spaces $P^{N-1}_\mathbb R,P^{N-1}_\mathbb C$:

\index{real projective space}
\index{complex projective space}

\begin{proposition}
We have presentation results as follows,
\begin{eqnarray*}
C(P^{N-1}_\mathbb R)&=&C^*_{comm}\left((p_{ij})_{i,j=1,\ldots,N}\Big|p=\bar{p}=p^t=p^2,Tr(p)=1\right)\\
C(P^{N-1}_\mathbb C)&=&C^*_{comm}\left((p_{ij})_{i,j=1,\ldots,N}\Big|p=p^*=p^2,Tr(p)=1\right)
\end{eqnarray*}
for the algebras of continuous functions on the real and complex projective spaces.
\end{proposition}

\begin{proof}
We use the fact that the projective spaces $P^{N-1}_\mathbb R,P^{N-1}_\mathbb C$ can be respectively identified with the spaces of rank one projections in $M_N(\mathbb R),M_N(\mathbb C)$. With this picture in mind, it is clear that we have arrows $\leftarrow$. In order to construct now arrows $\to$, consider the universal algebras on the right, $A_R,A_C$. These algebras being both commutative, by the Gelfand theorem we can write, with $X_R,X_C$ being certain compact spaces:
$$A_R=C(X_R)\quad,\quad 
A_C=C(X_C)$$

Now by using the coordinate functions $p_{ij}$, we conclude that $X_R,X_C$ are certain spaces of rank one projections in $M_N(\mathbb R),M_N(\mathbb C)$. In other words, we have embeddings:
$$X_R\subset P^{N-1}_\mathbb R\quad,\quad 
X_C\subset P^{N-1}_\mathbb C$$

By transposing we obtain arrows $\to$, as desired.
\end{proof}

The above result suggests the following definition:

\begin{definition}
Associated to any $N\in\mathbb N$ is the following universal algebra,
$$C(P^{N-1}_+)=C^*\left((p_{ij})_{i,j=1,\ldots,N}\Big|p=p^*=p^2,Tr(p)=1\right)$$
whose abstract spectrum is called ``free projective space''.
\end{definition}

Observe that, according to our presentation results for the real and complex projective spaces $P^{N-1}_\mathbb R$ and $P^{N-1}_\mathbb C$, we have embeddings of compact quantum spaces, as follows:
$$P^{N-1}_\mathbb R\subset P^{N-1}_\mathbb C\subset P^{N-1}_+$$

Let us first discuss the relation with the spheres. Given a closed subset $X\subset S^{N-1}_{\mathbb R,+}$, its projective version is by definition the quotient space $X\to PX$ determined by the fact that $C(PX)\subset C(X)$ is the subalgebra generated by the following variables:
$$p_{ij}=x_ix_j$$

In order to discuss the relation with the spheres, it is convenient to neglect the material regarding the complex and hybrid cases, the projective versions of such spheres bringing nothing new. Thus, we are left with the 3 real spheres, and we have:

\index{projective versions of spheres}

\begin{proposition}
The projective versions of the $3$ real spheres are as follows,
$$\xymatrix@R=15mm@C=15mm{
S^{N-1}_\mathbb R\ar[r]\ar[d]&S^{N-1}_{\mathbb R,*}\ar[r]\ar[d]&S^{N-1}_{\mathbb R,+}\ar[d]\\
P^{N-1}_\mathbb R\ar[r]&P^{N-1}_\mathbb C\ar[r]&P^{N-1}_+}$$
modulo the standard equivalence relation for the quantum algebraic manifolds.
\end{proposition}

\begin{proof}
The assertion at left is true by definition. For the assertion at right, we have to prove that the variables $p_{ij}=z_iz_j$ over the free sphere $S^{N-1}_{\mathbb R,+}$ satisfy the defining relations for $C(P^{N-1}_+)$, from Definition 15.5, namely:
$$p=p^*=p^2\quad,\quad 
Tr(p)=1$$

We first have the following computation:
$$(p^*)_{ij}
=p_{ji}^*
=(z_jz_i)^*
=z_iz_j
=p_{ij}$$

We have as well the following computation:
$$(p^2)_{ij}
=\sum_kp_{ik}p_{kj}
=\sum_kz_iz_k^2z_j
=z_iz_j\\
=p_{ij}$$

Finally, we have as well the following computation:
$$Tr(p)
=\sum_kp_{kk}
=\sum_kz_k^2
=1$$

Regarding now $PS^{N-1}_{\mathbb R,*}=P^{N-1}_\mathbb C$, the inclusion ``$\subset$'' follows from $abcd=cbad=cbda$. In the other sense now, the point is that we have a matrix model, as follows:
$$\pi:C(S^{N-1}_{\mathbb R,*})\to M_2(C(S^{N-1}_\mathbb C))\quad,\quad 
x_i\to\begin{pmatrix}0&z_i\\ \bar{z}_i&0\end{pmatrix}$$ 

But this gives the missing inclusion ``$\supset$'', and we are done. See \cite{bgo}.
\end{proof}

In addition to the above result, let us mention that, as already discussed above, passing to the complex case brings nothing new. This is because the projective version of the free complex sphere is equal to the free projective space constructed above:
$$PS^{N-1}_{\mathbb C,+}=P^{N-1}_+$$

And the same goes for the ``hybrid'' spheres. For details on all this, we refer to chapters 9-10. In what regards now the tori, we have here the following result:

\index{projective versions of tori}

\begin{proposition}
The projective versions of the $3$ real tori are as follows,
$$\xymatrix@R=15mm@C=15mm{
T_N\ar[r]\ar[d]&T_N^*\ar[r]\ar[d]&T_N^+\ar[d]\\
PT_N\ar[r]&P\mathbb T_N\ar[r]&PT_N^+}$$
modulo the standard equivalence relation for the quantum algebraic manifolds.
\end{proposition}

\begin{proof}
This follows indeed by using the same arguments as for the spheres.
\end{proof}

In what regards the orthogonal groups, we have here the following result:

\begin{proposition}
The projective versions of the $3$ orthogonal groups are
$$\xymatrix@R=15mm@C=15mm{
O_N\ar[r]\ar[d]&O_N^*\ar[r]\ar[d]&O_N^+\ar[d]\\
PO_N\ar[r]&PU_N\ar[r]&PO_N^+}$$
modulo the standard equivalence relation for the compact quantum groups.
\end{proposition}

\begin{proof}
This follows by using the same arguments as for spheres, or tori.
\end{proof}

Finally, in what regards the reflection groups, we have here the following result:

\begin{proposition}
The projective versions of the $3$ reflection groups are
$$\xymatrix@R=15mm@C=15mm{
H_N\ar[r]\ar[d]&H_N^*\ar[r]\ar[d]&H_N^+\ar[d]\\
PH_N\ar[r]&PK_N\ar[r]&PH_N^+}$$
modulo the standard equivalence relation for the compact quantum groups.
\end{proposition}

\begin{proof}
This follows indeed by using the same arguments as before.
\end{proof}

Let us mention that, as it was the case for the spheres, passing to the complex case brings nothing new. This is indeed because we have isomorphisms as follows, which can be established by using easiness, as explained in the beginning of the present chapter for the isomorphism in the middle, and with the proof of the first and of the last isomorphism being quite similar, based respectively on elementary group theory, and on easiness: 
$$P\mathbb T_N^+=PT_N^+\quad,\quad 
PU_N^+=PO_N^+\quad,\quad
PK_N^+=PH_N^+$$

Thus, passing to the complex case would bring indeed nothing new, and in what follows we will stay with the real formalism. This is of course something quite subtle, which happens only in the free case, and has no classical counterpart. For further details and comments here, we refer to the discussion in the beginning of this chapter.

\bigskip

Getting back now to our general program, we are done with the construction work, for the various projective geometry basic objects. As a conclusion to what we did, with our above constructions, in the projective geometry setting, we have 3 projective quadruplets, whose construction and main properties can be summarized as follows:

\index{projective quadruplet}
\index{real projective geometry}
\index{complex projective geometry}
\index{free projective geometry}
\index{real projective quadruplet}
\index{complex projective quadruplet}
\index{free projective quadruplet}

\begin{theorem}
We have projective quadruplets $(P,PT,PU,PK)$ as follows,
\begin{enumerate}
\item A classical real quadruplet, as follows,
$$\xymatrix@R=50pt@C=50pt{
P^{N-1}_\mathbb R\ar@{-}[r]\ar@{-}[d]\ar@{-}[dr]&PT_N\ar@{-}[l]\ar@{-}[d]\ar@{-}[dl]\\
PO_N\ar@{-}[u]\ar@{-}[ur]\ar@{-}[r]&PH_N\ar@{-}[l]\ar@{-}[ul]\ar@{-}[u]}$$

\item A classical complex quadruplet, as follows,
$$\xymatrix@R=50pt@C=50pt{
P^{N-1}_\mathbb C\ar@{-}[r]\ar@{-}[d]\ar@{-}[dr]&P\mathbb T_N\ar@{-}[l]\ar@{-}[d]\ar@{-}[dl]\\
PU_N\ar@{-}[u]\ar@{-}[ur]\ar@{-}[r]&PK_N\ar@{-}[l]\ar@{-}[ul]\ar@{-}[u]}$$

\item A free quadruplet, as follows,
$$\xymatrix@R=50pt@C=50pt{
P^{N-1}_+\ar@{-}[r]\ar@{-}[d]\ar@{-}[dr]&PT_N^+\ar@{-}[l]\ar@{-}[d]\ar@{-}[dl]\\
PO_N^+\ar@{-}[u]\ar@{-}[ur]\ar@{-}[r]&PH_N^+\ar@{-}[l]\ar@{-}[ul]\ar@{-}[u]}$$
\end{enumerate}
which appear as projective versions of the main $3$ real quadruplets.
\end{theorem}

\begin{proof}
This follows indeed from the results that already have. To be more precise, the details, that we will we need in what comes next, are as follows:

\medskip

(1) Consider the classical affine real quadruplet, which is as follows:
$$
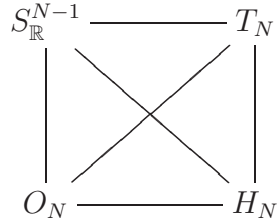
\xymatrix@R=50pt@C=50pt{
S^{N-1}_\mathbb R\ar@{-}[r]\ar@{-}[d]\ar@{-}[dr]&T_N\ar@{-}[l]\ar@{-}[d]\ar@{-}[dl]\\
O_N\ar@{-}[u]\ar@{-}[ur]\ar@{-}[r]&H_N\ar@{-}[l]\ar@{-}[ul]\ar@{-}[u]}$$

The projective version of this quadruplet is then the quadruplet in (1).

\medskip

(2) Consider the half-classical affine real quadruplet, which is as follows:
$$
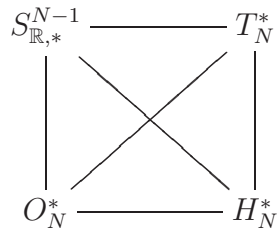
\xymatrix@R=50pt@C=50pt{
S^{N-1}_{\mathbb R,*}\ar@{-}[r]\ar@{-}[d]\ar@{-}[dr]&T_N^*\ar@{-}[l]\ar@{-}[d]\ar@{-}[dl]\\
O_N^*\ar@{-}[u]\ar@{-}[ur]\ar@{-}[r]&H_N^*\ar@{-}[l]\ar@{-}[ul]\ar@{-}[u]}$$

The projective version of this quadruplet is then the quadruplet in (2).

\medskip

(3) Consider the free affine real quadruplet, which is as follows:
$$
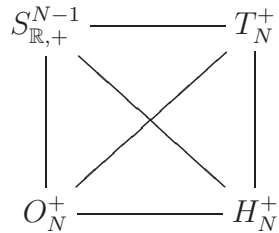
\xymatrix@R=50pt@C=50pt{
S^{N-1}_{\mathbb R,+}\ar@{-}[r]\ar@{-}[d]\ar@{-}[dr]&T_N^+\ar@{-}[l]\ar@{-}[d]\ar@{-}[dl]\\
O_N^+\ar@{-}[u]\ar@{-}[ur]\ar@{-}[r]&H_N^+\ar@{-}[l]\ar@{-}[ul]\ar@{-}[u]}$$

The projective version of this quadruplet is then the quadruplet in (3).
\end{proof}

\section*{15b. The threefold way}

Getting back now to our general projective geometry program, we would like to have axiomatization and classification results for such quadruplets. In order to do this, following \cite{bme}, we can axiomatize our various projective spaces, as follows:

\index{monomial space}
\index{monomial projective space}

\begin{definition}
A monomial projective space is a closed subset $P\subset P^{N-1}_+$ obtained via relations of type
$$p_{i_1i_2}\ldots p_{i_{k-1}i_k}=p_{i_{\sigma(1)}i_{\sigma(2)}}\ldots p_{i_{\sigma(k-1)}i_{\sigma(k)}},\ \forall (i_1,\ldots,i_k)\in\{1,\ldots,N\}^k$$
with $\sigma$ ranging over a certain subset of $\bigcup_{k\in2\mathbb N}S_k$, which is stable under $\sigma\to|\sigma|$.
\end{definition}

Observe the similarity with the corresponding monomiality notion for the spheres, from chapter 13. The only subtlety in the projective case is the stability under the operation $\sigma\to|\sigma|$, which in practice means that if the above relation associated to $\sigma$ holds, then the following relation, associated to $|\sigma|$, must hold as well:
$$p_{i_0i_1}\ldots p_{i_ki_{k+1}}=p_{i_0i_{\sigma(1)}}p_{i_{\sigma(2)}i_{\sigma(3)}}\ldots p_{i_{\sigma(k-2)}i_{\sigma(k-1)}}p_{i_{\sigma(k)}i_{k+1}}$$

As an illustration, the basic projective spaces are all monomial:

\begin{proposition}
The $3$ projective spaces are all monomial, with the permutations
$$
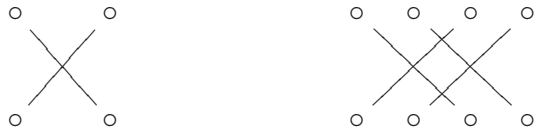
\xymatrix@R=10mm@C=8mm{\circ\ar@{-}[dr]&\circ\ar@{-}[dl]\\\circ&\circ}\qquad\ \qquad\ \qquad 
\xymatrix@R=10mm@C=3mm{\circ\ar@{-}[drr]&\circ\ar@{-}[drr]&\circ\ar@{-}[dll]&\circ\ar@{-}[dll]\\\circ&\circ&\circ&\circ}$$
producing respectively the spaces $P^{N-1}_\mathbb R,P^{N-1}_\mathbb C$, and with no relation needed for $P^{N-1}_+$.
\end{proposition}

\begin{proof}
We must divide the algebra $C(P^{N-1}_+)$ by the relations associated to the diagrams in the statement, as well as those associated to their shifted versions, given by:
$$
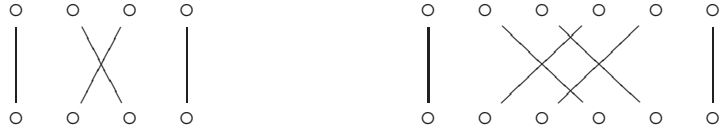
\xymatrix@R=10mm@C=3mm{\circ\ar@{-}[d]&\circ\ar@{-}[dr]&\circ\ar@{-}[dl]&\circ\ar@{-}[d]\\\circ&\circ&\circ&\circ}
\qquad\ \qquad\ \qquad 
\xymatrix@R=10mm@C=3mm{\circ\ar@{-}[d]&\circ\ar@{-}[drr]&\circ\ar@{-}[drr]&\circ\ar@{-}[dll]&\circ\ar@{-}[dll]&\circ\ar@{-}[d]\\\circ&\circ&\circ&\circ&\circ&\circ}$$ 

(1) The basic crossing, and its shifted version, produce the following relations: $$p_{ab}=p_{ba}$$
$$p_{ab}p_{cd}=p_{ac}p_{bd}$$

Now by using these relations several times, we obtain the following formula:
$$p_{ab}p_{cd}
=p_{ac}p_{bd}
=p_{ca}p_{db}
=p_{cd}p_{ab}$$

Thus, the space produced by the basic crossing is classical, $P\subset P^{N-1}_\mathbb C$. By using one more time the relations $p_{ab}=p_{ba}$ we conclude that we have $P=P^{N-1}_\mathbb R$, as claimed.

\medskip

(2) The fattened crossing, and its shifted version, produce the following relations:
$$p_{ab}p_{cd}=p_{cd}p_{ab}$$
$$p_{ab}p_{cd}p_{ef}=p_{ad}p_{eb}p_{cf}$$

The first relations tell us that the projective space must be classical, $P\subset P^{N-1}_\mathbb C$. Now observe that with $p_{ij}=z_i\bar{z}_j$, the second relations read:
$$z_a\bar{z}_bz_c\bar{z}_dz_e\bar{z}_f=z_a\bar{z}_dz_e\bar{z}_bz_c\bar{z}_f$$

Since these relations are automatic, we have $P=P^{N-1}_\mathbb C$, and we are done.
\end{proof}

Following \cite{bme}, we can now formulate our classification result, as follows:

\index{basic projective spaces}

\begin{theorem}
The basic projective spaces, namely 
$$P^{N-1}_\mathbb R\subset P^{N-1}_\mathbb C\subset P^{N-1}_+$$
are the only monomial ones.
\end{theorem}

\begin{proof}
We follow the proof from the affine case. Let $\mathcal R_\sigma$ be the collection of relations associated to a permutation $\sigma\in S_k$ with $k\in 2\mathbb N$, as in Definition 15.11. We fix a monomial projective space $P\subset P^{N-1}_+$, and we associate to it subsets $G_k\subset S_k$, as follows:
$$G_k=\begin{cases}
\{\sigma\in S_k|\mathcal R_\sigma\ {\rm hold\ over\ }P\}&(k\ {\rm even})\\
\{\sigma\in S_k|\mathcal R_{|\sigma}\ {\rm hold\ over\ }P\}&(k\ {\rm odd})
\end{cases}$$

As in the affine case, we obtain in this way a filtered group $G=(G_k)$, which is stable under removing outer strings, and under removing neighboring strings.  Thus the computations in chapter 13 apply, and show that we have only 3 possible situations, corresponding to the 3 projective spaces in Proposition 15.12.
\end{proof}

Let us discuss now similar results for the projective quantum groups. Given a closed subgroup $G\subset O_N^+$, its projective version $G\to PG$ is by definition given by the fact that $C(PG)\subset C(G)$ is the subalgebra generated by the following variables:
$$w_{ij,ab}=u_{ia}u_{jb}$$

In the classical case we recover in this way the usual projective version:
$$PG=G/(G\cap\mathbb Z_2^N)$$

We have the following key result, from \cite{bc+}:

\begin{theorem}
The quantum group $O_N^*$ is the unique intermediate easy quantum group $O_N\subset G\subset O_N^+$. Moreover, in the non-easy case, the following happen:
\begin{enumerate}
\item The group inclusion $\mathbb TO_N\subset U_N$ is maximal.

\item The group inclusion $PO_N\subset PU_N$ is maximal.

\item The quantum group inclusion $O_N\subset O_N^*$ is maximal.
\end{enumerate}
\end{theorem}

\begin{proof}
This is something discussed in chapters 9-10, the idea being that the first assertion comes by classifying the categories of pairings, and then:

\medskip

(1) This can be obtained by using standard Lie group methods.

\medskip

(2) This follows from (1), by taking projective versions. 

\medskip

(3) This follows from (2), via standard algebraic lifting results. 
\end{proof}

Our claim now is that, under suitable assumptions, $PU_N$ is the only intermediate object $PO_N\subset G\subset PO_N^+$. In order to formulate a precise statement here, we first recall the following notion, from \cite{bsp}, that we have already heavily used in this book:

\begin{definition}
A collection of sets $D=\bigsqcup_{k,l}D(k,l)$ with 
$$D(k,l)\subset P(k,l)$$
is called a category of partitions when it has the following properties:
\begin{enumerate}
\item Stability under the horizontal concatenation, $(\pi,\sigma)\to[\pi\sigma]$.

\item Stability under vertical concatenation $(\pi,\sigma)\to[^\sigma_\pi]$, with matching middle symbols.

\item Stability under the upside-down turning $*$, with switching of colors, $\circ\leftrightarrow\bullet$.

\item Each set $P(k,k)$ contains the identity partition $||\ldots||$.

\item The sets $P(\emptyset,\circ\bullet)$ and $P(\emptyset,\bullet\circ)$ both contain the semicircle $\cap$.
\end{enumerate}
\end{definition} 

The above definition is something inspired from the axioms of Tannakian categories, and going hand in hand with it is the following definition, also from \cite{bsp}:

\begin{definition}
An intermediate compact quantum group 
$$O_N\subset G\subset O_N^+$$
is called easy when the corresponding Tannakian category
$$span(NC_2(k,l))\subset Hom(u^{\otimes k},u^{\otimes l})\subset span(P_2(k,l))$$
comes via the following formula, using the standard $\pi\to T_\pi$ construction,
$$Hom(u^{\otimes k},u^{\otimes l})=span(D(k,l))$$
from a certain collection of sets of pairings $D=(D(k,l))$.
\end{definition}

As explained in \cite{bsp}, by ``saturating'' the sets $D(k,l)$, we can assume that the collection $D=(D(k,l))$ is a category of pairings, in the sense that it is stable under vertical and horizontal concatenation, upside-down turning, and contains the semicircle.

\bigskip

In the projective case now, following \cite{bme}, let us formulate:

\index{projective category of pairings}
\index{projective easiness}
\index{projective quantum group}

\begin{definition}
A projective category of pairings is a collection of subsets 
$$NC_2(2k,2l)\subset E(k,l)\subset P_2(2k,2l)$$
stable under the usual categorical operations, and satisfying $\sigma\in E\implies |\sigma|\in E$.
\end{definition}

As basic examples here, we have the following projective categories of pairings, where $P_2^*$ is the category of matching pairings:
$$NC_2\subset P_2^*\subset P_2$$

This follows indeed from definitions. Now with the above notion in hand, we can formulate the following projective analogue of the notion of easiness:

\begin{definition}
An intermediate compact quantum group 
$$PO_N\subset H\subset PO_N^+$$
is called projectively easy when its Tannakian category
$$span(NC_2(2k,2l))\subset Hom(v^{\otimes k},v^{\otimes l})\subset span(P_2(2k,2l))$$
comes via via the following formula, using the standard $\pi\to T_\pi$ construction,
$$Hom(v^{\otimes k},v^{\otimes l})=span(E(k,l))$$
for a certain projective category of pairings $E=(E(k,l))$.
\end{definition}

Thus, we have a projective notion of easiness. Observe that, given an easy quantum group $O_N\subset G\subset O_N^+$, its projective version $PO_N\subset PG\subset PO_N^+$ is projectively easy in our sense. In particular the basic projective quantum groups $PO_N\subset PU_N\subset PO_N^+$ are all projectively easy in our sense, coming from the categories $NC_2\subset P_2^*\subset P_2$. 

\bigskip

We have in fact the following general result, from \cite{bme}:

\begin{theorem}
We have a bijective correspondence between the affine and projective categories of partitions, given by the operation
$$G\to PG$$ 
at the level of the corresponding affine and projective easy quantum groups.
\end{theorem}

\begin{proof}
The construction of correspondence $D\to E$ is clear, simply by setting:
$$E(k,l)=D(2k,2l)$$

Indeed, due to the axioms in Definition 15.15, the conditions in Definition 15.17 are satisfied. Conversely, given $E=(E(k,l))$ as in Definition 15.17, we can set:
$$D(k,l)=\begin{cases}
E(k,l)&(k,l\ {\rm even})\\
\{\sigma:|\sigma\in E(k+1,l+1)\}&(k,l\ {\rm odd})
\end{cases}$$

Our claim is that $D=(D(k,l))$ is a category of partitions. Indeed:

\medskip

(1) The composition action is clear. Indeed, when looking at the numbers of legs involved, in the even case this is clear, and in the odd case, this follows from:
\begin{eqnarray*}
|\sigma,|\sigma'\in E
&\implies&|^\sigma_\tau\in E\\
&\implies&{\ }^\sigma_\tau\in D
\end{eqnarray*}

(2) For the tensor product axiom, we have 4 cases to be investigated, depending on the parity of the number of legs of $\sigma,\tau$, as follows:

\medskip

-- The even/even case is clear. 

\medskip

-- The odd/even case follows from the following computation:
\begin{eqnarray*}
|\sigma,\tau\in E
&\implies&|\sigma\tau\in E\\
&\implies&\sigma\tau\in D
\end{eqnarray*}

-- Regarding now the even/odd case, this can be solved as follows:
\begin{eqnarray*}
\sigma,|\tau\in E
&\implies&|\sigma|,|\tau\in E\\
&\implies&|\sigma||\tau\in E\\
&\implies&|\sigma\tau\in E\\
&\implies&\sigma\tau\in D
\end{eqnarray*}

-- As for the remaining odd/odd case, here the computation is as follows:
\begin{eqnarray*}
|\sigma,|\tau\in E
&\implies&||\sigma|,|\tau\in E\\
&\implies&||\sigma||\tau\in E\\
&\implies&\sigma\tau\in E\\
&\implies&\sigma\tau\in D
\end{eqnarray*}

(3) Finally, the conjugation axiom is clear from definitions. It is also clear that both compositions $D\to E\to D$ and $E\to D\to E$ are the identities, as claimed. As for the quantum group assertion, this is clear as well from definitions.
\end{proof}

Now back to uniqueness issues, we have here the following result, also from \cite{bme}:

\begin{theorem}
We have the following results:
\begin{enumerate}
\item $O_N^*$ is the only intermediate easy quantum group, as follows:
$$O_N\subset G\subset O_N^+$$

\item $PU_N$ is the only intermediate projectively easy quantum group 
$$PO_N\subset G\subset PO_N^+$$
in the projective easiness sense axiomatized above.
\end{enumerate}
\end{theorem}

\begin{proof}
The idea here is as follows:

\medskip

(1) The assertion regarding $O_N\subset O_N^*\subset O_N^+$ is from \cite{bc+}, and this is something that we already know, explained in chapter 9.

\medskip

(2) The assertion regarding $PO_N\subset PU_N\subset PO_N^+$ follows from the classification result in (1), and from the duality in Theorem 15.19.
\end{proof}

Summarizing, we have analogues of the various affine classification results, with the remark that everything becomes simpler in the projective setting. 

\section*{15c. Projective geometry}

We have so far projective analogues of the various affine classification results. In view of this, our next goal will be that of finding projective versions of the quantum isometry group results that we have in the affine setting. We use the following action formalism, which is quite similar to the affine action formalism introduced in chapter 3:

\index{projective action}
\index{projective isometry}
\index{projective affine isometry}
\index{projective quantum isometry}

\begin{definition}
Consider a closed subgroup of the free orthogonal group, $G\subset O_N^+$, and a closed subset of the free real sphere, $X\subset S^{N-1}_{\mathbb R,+}$.
\begin{enumerate}
\item We write $G\curvearrowright X$ when we have a morphism of $C^*$-algebras, as follows:
$$\Phi:C(X)\to C(X)\otimes C(G)$$
$$\Phi(z_i)=\sum_az_a\otimes u_{ai}$$ 

\item We write $PG\curvearrowright PX$ when we have a morphism of $C^*$-algebras, as follows:
$$\Phi:C(PX)\to C(PX)\otimes C(PG)$$
$$\Phi(z_iz_j)=\sum_a z_az_b\otimes u_{ai}u_{bj}$$
\end{enumerate}
\end{definition}

Observe that the above morphisms $\Phi$, if they exist, are automatically coaction maps. Observe also that an affine action $G\curvearrowright X$ produces a projective action $PG\curvearrowright PX$. Let us also mention that given an algebraic subset $X\subset S^{N-1}_{\mathbb R,+}$, it is routine to prove that there exist indeed universal quantum groups $G\subset O_N^+$ acting as (1), and as in (2). We have the following result, from \cite{bgo} and related papers, with respect to the above notions:

\begin{theorem}
The quantum isometry groups of basic spheres and projective spaces,
$$\xymatrix@R=15mm@C=15mm{
S^{N-1}_\mathbb R\ar[r]\ar[d]&S^{N-1}_{\mathbb R,*}\ar[r]\ar[d]&S^{N-1}_{\mathbb R,+}\ar[d]\\
P^{N-1}_\mathbb R\ar[r]&P^{N-1}_\mathbb C\ar[r]&P^{N-1}_+}$$
are the following affine and projective quantum groups,
$$\xymatrix@R=15mm@C=15mm{
O_N\ar[r]\ar[d]&O_N^*\ar[r]\ar[d]&O_N^+\ar[d]\\
PO_N\ar[r]&PU_N\ar[r]&PO_N^+}$$
with respect to the affine and projective action notions introduced above.
\end{theorem}

\begin{proof}
The fact that the 3 quantum groups on top act affinely on the corresponding 3 spheres is known since \cite{bgo}, and is elementary, explained before. By restriction, the 3 quantum groups on the bottom follow to act on the corresponding 3 projective spaces. We must prove now that all these actions are universal. At right there is nothing to prove, so we are left with studying the actions on $S^{N-1}_\mathbb R,S^{N-1}_{\mathbb R,*}$ and on $P^{N-1}_\mathbb R,P^{N-1}_\mathbb C$.

\medskip

\underline{$P^{N-1}_\mathbb R$}. Consider the following projective coordinates:
$$p_{ij}=z_iz_j\quad,\quad 
w_{ij,ab}=u_{ai}u_{bj}$$

In terms of these projective coordinates, the coaction map is given by:
$$\Phi(p_{ij})=\sum_{ab}p_{ab}\otimes w_{ij,ab}$$

Thus, we have the following formulae:
\begin{eqnarray*}
\Phi(p_{ij})&=&\sum_{a<b}p_{ab}\otimes (w_{ij,ab}+w_{ij,ba})+\sum_ap_{aa}\otimes w_{ij,aa}\\
\Phi(p_{ji})&=&\sum_{a<b}p_{ab}\otimes (w_{ji,ab}+w_{ji,ba})+\sum_ap_{aa}\otimes w_{ji,aa}
\end{eqnarray*}

By comparing these two formulae, and then by using the linear independence of the variables $p_{ab}=z_az_b$ for $a\leq b$, we conclude that we must have:
$$w_{ij,ab}+w_{ij,ba}=w_{ji,ab}+w_{ji,ba}$$

Let us apply now the antipode to this formula. For this purpose, observe that:
\begin{eqnarray*}
S(w_{ij,ab})
&=&S(u_{ai}u_{bj})\\
&=&S(u_{bj})S(u_{ai})\\
&=&u_{jb}u_{ia}\\
&=&w_{ba,ji}
\end{eqnarray*}

Thus by applying the antipode we obtain:
$$w_{ba,ji}+w_{ab,ji}=w_{ba,ij}+w_{ab,ij}$$

By relabeling, we obtain the following formula: 
$$w_{ji,ba}+w_{ij,ba}=w_{ji,ab}+w_{ij,ab}$$

Now by comparing with the original relation, we obtain:
$$w_{ij,ab}=w_{ji,ba}$$

But, with $w_{ij,ab}=u_{ai}u_{bj}$, this formula reads:
$$u_{ai}u_{bj}=u_{bj}u_{ai}$$

Thus $G\subset O_N$, and it follows that we have $PG\subset PO_N$, as claimed.

\medskip

\underline{$P^{N-1}_\mathbb C$}. Consider a coaction map, written as follows, with $p_{ab}=z_a\bar{z}_b$:
$$\Phi(p_{ij})=\sum_{ab}p_{ab}\otimes u_{ai}u_{bj}$$

The idea here will be that of using the following formula:
$$p_{ab}p_{cd}=p_{ad}p_{cb}$$

We have the following formulae:
\begin{eqnarray*}
\Phi(p_{ij}p_{kl})&=&\sum_{abcd}p_{ab}p_{cd}\otimes u_{ai}u_{bj}u_{ck}u_{dl}\\
\Phi(p_{il}p_{kj})&=&\sum_{abcd}p_{ad}p_{cb}\otimes u_{ai}u_{dl}u_{ck}u_{bj}
\end{eqnarray*}

The terms at left being equal, and the last terms at right being equal too, we deduce that, with $[a,b,c]=abc-cba$, we must have the following formula:
$$\sum_{abcd}u_{ai}[u_{bj},u_{ck},u_{dl}]\otimes p_{ab}p_{cd}=0$$

Now since the quantities $p_{ab}p_{cd}=z_a\bar{z}_bz_c\bar{z}_d$ at right depend only on the numbers $|\{a,c\}|,|\{b,d\}|\in\{1,2\}$, and this dependence produces the only possible linear relations between the variables $p_{ab}p_{cd}$, we are led to $2\times2=4$ equations, as follows:

\medskip

(1) $u_{ai}[u_{bj},u_{ak},u_{bl}]=0$, $\forall a,b$.

\medskip

(2) $u_{ai}[u_{bj},u_{ak},u_{dl}]+u_{ai}[u_{dj},u_{ak},u_{bl}]=0$, $\forall a$, $\forall b\neq d$.

\medskip

(3) $u_{ai}[u_{bj},u_{ck},u_{bl}]+u_{ci}[u_{bj},u_{ak},u_{bl}]=0$, $\forall a\neq c$, $\forall b$.

\medskip

(4) $u_{ai}[u_{bj},u_{ck},u_{dl}]+u_{ai}[u_{dj},u_{ck},u_{bl}]+u_{ci}[u_{bj},u_{ak},u_{dl}]+u_{ci}[u_{dj},u_{ak},u_{bl}]=0$, $\forall a\neq c,b\neq d$.

\medskip

We will need in fact only the first two formulae. Since (1) corresponds to (2) at $b=d$, we conclude that (1,2) are equivalent to (2), with no restriction on the indices. By multiplying now this formula to the left by $u_{ai}$, and then summing over $i$, we obtain:
$$[u_{bj},u_{ak},u_{dl}]+[u_{dj},u_{ak},u_{bl}]=0$$

We use now the antipode/relabel trick from \cite{bg1}. By applying the antipode we obtain: 
$$[u_{ld},u_{ka},u_{jb}]+[u_{lb},u_{ka},u_{jd}]=0$$

By relabelling we obtain the following formula:
$$[u_{dl},u_{ak},u_{bj}]+[u_{dj},u_{ak},u_{bl}]=0$$

Now by comparing with the original relation, we obtain:
$$[u_{bj},u_{ak},u_{dl}]=[u_{dj},u_{ak},u_{bl}]=0$$

Thus $G\subset O_N^*$, and it follows that we have $PG\subset PU_N$, as desired.
\end{proof}

The above results can be probably improved. As an example, let us say that a closed subgroup $G\subset U_N^+$ acts projectively on $PX$ when we have a coaction map as follows:
$$\Phi(z_iz_j)=\sum_{ab}z_az_b\otimes u_{ai}u_{bj}^*$$

The above proof can be adapted, by putting $*$ signs where needed, and Theorem 15.22 still holds, in this setting. However, establishing general universality results, involving arbitrary subgroups $H\subset PO_N^+$, looks like a quite non-trivial question.

\bigskip

Let us discuss now the axiomatization question for the projective quadruplets of type $(P,PT,PU,PK)$. We recall that we first have a classical real quadruplet, as follows:
$$
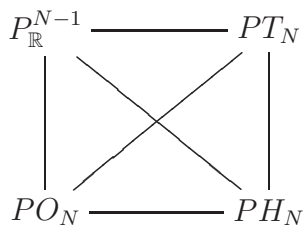
\xymatrix@R=50pt@C=50pt{
P^{N-1}_\mathbb R\ar@{-}[r]\ar@{-}[d]\ar@{-}[dr]&PT_N\ar@{-}[l]\ar@{-}[d]\ar@{-}[dl]\\
PO_N\ar@{-}[u]\ar@{-}[ur]\ar@{-}[r]&PH_N\ar@{-}[l]\ar@{-}[ul]\ar@{-}[u]}$$

We have then a classical complex quadruplet, which can be thought of as well as being a real half-classical quadruplet, which is as follows:
$$
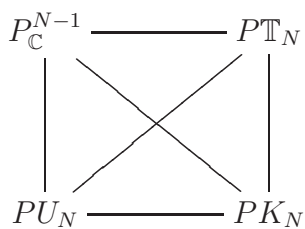
\xymatrix@R=50pt@C=50pt{
P^{N-1}_\mathbb C\ar@{-}[r]\ar@{-}[d]\ar@{-}[dr]&P\mathbb T_N\ar@{-}[l]\ar@{-}[d]\ar@{-}[dl]\\
PU_N\ar@{-}[u]\ar@{-}[ur]\ar@{-}[r]&PK_N\ar@{-}[l]\ar@{-}[ul]\ar@{-}[u]}$$

Finally, we have a free quadruplet, which can be thought of as being the same time real and complex, which is as follows:
$$
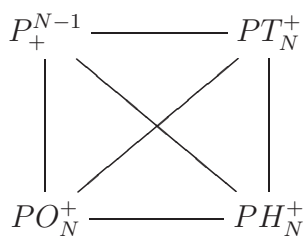
\xymatrix@R=50pt@C=50pt{
P^{N-1}_+\ar@{-}[r]\ar@{-}[d]\ar@{-}[dr]&PT_N^+\ar@{-}[l]\ar@{-}[d]\ar@{-}[dl]\\
PO_N^+\ar@{-}[u]\ar@{-}[ur]\ar@{-}[r]&PH_N^+\ar@{-}[l]\ar@{-}[ul]\ar@{-}[u]}$$

The question is that of axiomatizing these quadruplets.

\bigskip

To be more precise, in analogy with what happens in the affine case, the problem is that of establishing correspondences as follows:
$$
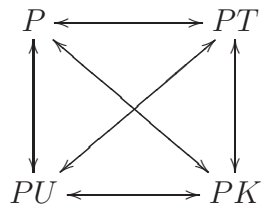
\xymatrix@R=50pt@C=50pt{
P\ar[r]\ar[d]\ar[dr]&PT\ar[l]\ar[d]\ar[dl]\\
PU\ar[u]\ar[ur]\ar[r]&PK\ar[l]\ar[ul]\ar[u]
}$$

Modulo this problem, which is for the moment open, things are potentially quite nice, because we seem to have only 3 geometries, namely real, complex and free. Generally speaking, we are led in this way into several questions:

\bigskip

(1) We first need functoriality results for the operations $<\,,>$ and $\cap$, in relation with taking the projective version, and taking affine lifts, as to deduce most of our 7 axioms, in their obvious projective formulation, from the affine ones.

\bigskip

(2) Then, we need quantum isometry results in the projective setting, for the projective spaces themselves, and for the projective tori, either established ad-hoc, or by using the affine results. For the projective spaces, this was done above.

\bigskip

(3) We need as well some further functoriality results, in order to axiomatize the intermediate objects that we are dealing, the problem here being whether we want to use projective objects, or projective versions, perhaps saturated, of affine objects.

\bigskip

(4) Modulo this, things are quite clear, with the final result being the fact that we have only 3 projective geometries. Technically, the proof should be using the fact that, in the easy setting, $PO_N\subset PU_N\subset PO_N^+$ are the only possible unitary groups.

\bigskip

Let us also mention that, in the noncommutative setting, there are several ways of defining the projective versions, with the one used here being the ``simplest''. As explained in \cite{ba8}, \cite{bb2}, it is possible to construct a left projective version, a right projective version, and a mixed projective version, with all these operations being interesting. Thus, the results and problems presented above are just the ``tip of the iceberg'', with the general projective space and version problematics being much wider then this.

\bigskip

As another remark, our results tend to suggest that the free projective geometry is ``scalarless''. However, things here are quite complicated, because, while there has been some interesting preliminary work on this subject by Bichon and others, it is still not presently known what easiness should mean, over an arbitrary field $F$. 

\bigskip

Finally, and above everything, the free projective geometry remains to be developed. A first piece of homework here is that of developing a theory of free Grassmannians, free flag manifolds, and free Stiefel manifolds, based on the affine theory of the spaces of quantum partial isometries, developed in chapter 6. 

\bigskip

To be more precise, the definition of the free Grassmannians is straightforward, inspired from that of the free projective spaces, as follows: 
$$C(Gr_{LN}^+)=C^*\left((p_{ij})_{i,j=1,\ldots,N}\Big|p=p^*=p^2,Tr(p)=L\right)$$

As for the definition of the free flag manifolds and of the free Stiefel manifolds, this is most likely something very similar, using relations a bit as those above.

\bigskip

Things do not look difficult here, with most of the arguments from the affine case carrying over in the projective setting, and with solid affine results to rely upon being available from chapter 6. But work to be done for sure, which has not been done yet.

\section*{15d. Small dimensions} 

We would like to end this chapter with something refreshing, namely a preliminary study of the free analogue of $P^2_\mathbb R$. We recall that the projective space $P^{N-1}_\mathbb R$ is the space of lines in $\mathbb R^N$ passing through the origin, the basic examples being as follows:

\bigskip

(1) At $N=2$ each such a line, in $\mathbb R^2$ passing through the origin, corresponds to 2 opposite points on the unit circle $\mathbb T\subset\mathbb R^2$, according to the following scheme:
$$\xymatrix@R=6.7pt@C=2pt{
&&&&&&\ar@{-}[d]\\
&&&&&&\bullet\ar@{-}[dddddd]\ar@{-}@/_/[dllll]&\\
&&\ar@{-}@/_/[ddl]&&&&&&&&\ar@.[ddddllllllll]\ar@{-}@/_/[ullll]\\
&&&&&&&&&&\\
\ar@{-}[r]&\bullet\ar@{-}@/_/[ddr]\ar@{-}[rrrrr]&&&&&\ar@{-}[rrrrr]&&&&&\bullet\ar@{-}@/_/[uul]\ar@{-}[r]&\\
&&&&&&&&&&\\
&&\ar@{-}@/_/[drrrr]\ar@.[uuuurrrrrrrr]&&&&&&&&\ar@{-}@/_/[uur]\\
&&&&&&\bullet\ar@{-}@/_/[urrrr]\ar@{-}[d]\\
&&&&&&}$$

Thus, $P^1_\mathbb R$ corresponds to the upper semicircle of $\mathbb T$, with the endpoints identified, so we obtain a circle, $P^1_\mathbb R=\mathbb T$, according to the following scheme: 
$$\xymatrix@R=7pt@C=2pt{
&&&&&&\ar@{-}[d]\\
&&&&&&\bullet\ar@{-}[dddd]\ar@{-}@/_/[dllll]&\\
&&\ar@{-}@/_/[ddl]&&&&&&&&\ar@{-}@/_/[ullll]\\
&&&&&&&&&&\\
\ar@{-}[r]&\bullet\ar@.[rrrrrrrrrr]&&&&&&&&&&\bullet\ar@{-}@/_/[uul]\ar@{-}[r]\ar@.[llllllllll]&\\
&&&&&&&&
}$$

(2) At $N=3$ the situation is similar, with $P^2_\mathbb R$ corresponding to the upper hemisphere of the sphere $S^2_\mathbb R\subset\mathbb R^3$, with the points on the equator identified via $x=-x$. Topologically speaking, we can deform if we want the upper hemisphere into a square, with the equator becoming the boundary of this square, and in this picture, the $x=-x$ identification corresponds to the ``identify opposite edges, with opposite orientations'' folding method for the square, leading to a space $P^2_\mathbb R$ which is obviously not embeddable into $\mathbb R^3$:
$$\xymatrix@R=55pt@C=55pt{
\circ\ar[r]&\circ\ar@{-->}[d]\\
\circ\ar@{-->}[u]&\circ\ar[l]}$$

In what follows we will be interested in the free analogue $P^2_+$ of this projective space $P^2_\mathbb R$. Our main motivation comes from the fact that, according to the work of Bhowmick-D'Andrea-Dabrowski \cite{bdd}, later on continued with Das \cite{bd+},  the quantum isometry group $PO_3^+=PU_3^+$ of the free projective space $P^2_+$ acts on the quark part of the Standard Model spectral triple, in Chamseddine-Connes formulation \cite{cc1}, \cite{cc2}.

\bigskip

We recall that the free projective space is defined by the following formula:
$$C(P^{N-1}_+)=C^*\left((p_{ij})_{i,j=1,\ldots,N}\Big|p=p^*=p^2,Tr(p)=1\right)$$

Let us first discuss, as a warm-up, the 2D case. Here the above matrix of projective coordinates is as follows, with $a=a^*$, $b=b^*$, $a+b=1$:
$$p=\begin{pmatrix}a&c\\ c^*&b\end{pmatrix}$$

We have the following computation:
$$p^2=\begin{pmatrix}a&c\\ c^*&b\end{pmatrix}\begin{pmatrix}a&c\\ c^*&b\end{pmatrix}
=\begin{pmatrix}a^2+cc^*&ac+cb\\ c^*a+bc^*&c^*c+b^2\end{pmatrix}$$

Thus, the equations to be satisfied are as follows:
$$a^2+cc^*=a$$
$$b^2+c^*c=b$$
$$ac+cb=c$$
$$c^*a+bc^*=c^*$$

The 4th equation is the conjugate of the 3rd equation, so we remove it. By using $a+b=1$, the remaining equations can be written as:
$$cc^*=c^*c=ab$$
$$ac+ca=0$$

We have several explicit models for this, using the spheres $S^1_{\mathbb R,+}$ and $S^1_{\mathbb C,+}$, as well as the first row spaces of $O_2^+$ and $U_2^+$, which ultimately lead us to $SU_2$ and $\bar{SU}_2$. These models are known to be all equivalent under Haar, and the question is whether they are identical. Thus, we must do computations as above in all models, and compare. These are all interesting questions, whose precise answers are not known, so far.

\bigskip

In the 3D case now, that of projective space $P^2_+$, that we are mainly interested in here, the matrix of coordinates is as follows, with $r,s,t$ self-adjoint, $r+s+t=1$:
$$p=\begin{pmatrix}
r&a&b\\ a^*&s&c\\ b^*&c^*&t\end{pmatrix}$$

The square of this matrix is given by:
$$p^2=\begin{pmatrix}
r&a&b\\ a^*&s&c\\ b^*&c^*&t\end{pmatrix}
\begin{pmatrix}
r&a&b\\ a^*&s&c\\ b^*&c^*&t\end{pmatrix}$$

We obtain the following formula:
$$p^2=\begin{pmatrix}
r^2+aa^*+bb^*&ra+as+bc^*&rb+ac+bt\\
a^*r+sa^*+cb^*&a^*a+s^2+cc^*&a^*b+sc+ct\\
b^*r+c^*a^*+tb^*&b^*a+c^*s+tc^*&b^*b+c^*c+t^2
\end{pmatrix}$$

On the diagonal, the equations for $p^2=p$ are as follows:
$$aa^*+bb^*=r-r^2$$
$$a^*a+cc^*=s-s^2$$
$$b^*b+c^*c=t-t^2$$

On the off-diagonal upper part, the equations for $p^2=p$ are as follows:
$$ra+as+bc^*=a$$
$$rb+ac+bt=b$$
$$a^*b+sc+ct=c$$

On the off-diagonal lower part, the equations for $p^2=p$ are those above, conjugated. Thus, we have 6 equations. The first problem is that of using $r+s+t=1$, in order to make these equations look better. Again, many interesting questions here.

\bigskip

Observe the analogy with the basic discussion about hypersurfaces, and about basic affine geometry in general, from the end of chapter 13. In both cases indeed we are led to an interesting mix of basic algebraic geometry and operator theory, and with the operator theory component potentially ranging from very basic to very complicated.

\section*{15e. Exercises} 

Projective noncommutative geometry is virtually as big as affine noncommutative geometry, and there are countless questions. We first have the following exercise:

\begin{exercise}
Prove that we have the free complexification formula
$$\widetilde{SU_2}=U_2^+$$
and then look for analogues of this formula, at arbitrary $N\in\mathbb N$.
\end{exercise}

Here the first formula can be established directly, say with character computations, but if you want to solve the second question as well, things become more complicated. As a hint, try first to develop of theory of ``free symplectic groups''.

\begin{exercise}
Prove that we have the formula
$$PO_N^+=PU_N^+$$
directly, without using the free complexification results.
\end{exercise}

This was something that was already discussed in the above, the idea being that of comparing the corresponding Hom spaces, obtained via easiness.

\begin{exercise}
Prove that we have the formula
$$PU_2^+=SO_3$$
and then look for analogues of this, at arbitrary $N\in\mathbb N$.
\end{exercise}

As before, the first assertion is quite elementary, but if you want to solve the second question as well, things become more complicated, and you are led into free symplectic groups, and their easiness type property, which is not exactly the usual easiness.

\begin{exercise}
Find axioms for the projective quadruplets
$$
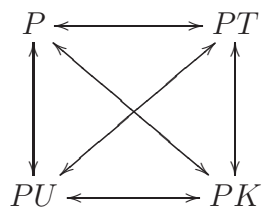
\xymatrix@R=50pt@C=50pt{
P\ar[r]\ar[d]\ar[dr]&PT\ar[l]\ar[d]\ar[dl]\\
PU\ar[u]\ar[ur]\ar[r]&PK\ar[l]\ar[ul]\ar[u]
}$$
covering the real, complex and free quadruplets, constructed above.
\end{exercise}

This is something that we already discussed in the above, and good luck. In addition to all this, there are many interesting questions regarding the development of free projective geometry, that were partly discussed in the above. The whole area is obviously very wide, and interesting, and anything here would be welcome.

\chapter{Hyperspherical laws}

\section*{16a. Calculus}

We have kept the best for this final chapter. Calculus, more calculus, and even more calculus, in relation with the integration over spheres. Our motivations are varied:

\bigskip

(1) First of all, calculus is a good thing, and calculus over spheres, using spherical coordinates, is even better. Mathematicians usually snub spherical coordinates, deemed ``unconceptual'', but physicists just love them. Want to do some electrodynamics? Spherical coordinates. Want to solve the hydrogen atom? Spherical coordinates, too. So, following the physicists, we will love these spherical coordinates too, in this chapter. And let me recommend here again the delightful books of Griffiths \cite{gr1}, \cite{gr2}.

\bigskip

(2) Second, the spheres themselves are a very good thing too, be that in the context of the Connes noncommutative geometry \cite{cdu}, \cite{cla}, \cite{ddl}, or in the context of our noncommutative geometry, following \cite{bgo} and related papers, and as explained so far in this book, or in the context of any other kind of noncommutative geometry theory. Also, in our setting, everything more advanced, as for instance of analysis over free manifolds type, like the work in \cite{cfk}, \cite{dfw}, \cite{dgo}, starts of course with a study in the sphere case.

\bigskip

(3) Finally, all this, calculus over spheres, will naturally lead us into all sorts of advanced considerations. At the core of all this will be a tough computation from \cite{bcz}, as well as a subtle twisting result from \cite{bbs}, relating the free orthogonal/unitary projective quantum group $PO_N^+=PU_N^+$ to the quantum permutation group $S_{N^2}^+$. And with this being virtually related to pretty much everything, mathematics and physics alike, including \cite{bdd}, \cite{bd+}, \cite{cc1}, \cite{cc2}, \cite{dif}, \cite{jo1}, \cite{jo2}, \cite{jo3}, \cite{mpa}, \cite{nas}, \cite{vdn}, \cite{wig}, \cite{wit}.

\bigskip

As a starting point, we have the very natural question, first investigated in \cite{bgo}, of computing the laws of individual coordinates of the main 3 real spheres, namely:
$$S^{N-1}_\mathbb R\subset S^{N-1}_{\mathbb R,*}\subset S^{N-1}_{\mathbb R,+}$$

We already know from chapter 5 the $N\to\infty$ behavior of these laws, called ``hyperspherical''. To be more precise, for $S^{N-1}_\mathbb R$ we obtain the normal law, and for $S^{N-1}_{\mathbb R,+}$ we obtain the semicircle law. As for the sphere $S^{N-1}_{\mathbb R,*}$, this has the same projective version as $S^{N-1}_\mathbb C$, where the corresponding law becomes complex Gaussian with $N\to\infty$, as explained in chapter 5, and so we obtain a symmetrized Rayleigh variable. See \cite{bcs}.

\bigskip

The problem that we want to investigate is that of computing these hyperspherical laws at fixed values of $N\in\mathbb N$. Let us begin with a discussion in the classical case. At $N=2$ the sphere is the unit circle $\mathbb T$, with $z=e^{it}$ the coordinates are
$x=\cos t$, $y=\sin t$, and the integrals of the products of such coordinates can be computed as follows:

\index{trigonometric integral}
\index{partial integration}
\index{double factorial}

\begin{theorem}
We have the following formula,
$$\int_0^{\pi/2}\cos^pt\sin^qt\,dt=\left(\frac{\pi}{2}\right)^{\varepsilon(p)\varepsilon(q)}\frac{p!!q!!}{(p+q+1)!!}$$
where $\varepsilon(p)=1$ if $p$ is even, and $\varepsilon(p)=0$ if $p$ is odd, and where
$$m!!=(m-1)(m-3)(m-5)\ldots$$
with the product ending at $2$ if $m$ is odd, and ending at $1$ if $m$ is even.
\end{theorem}

\begin{proof}
This is standard calculus, with particular cases of this formula being very familiar to everyone loving and teaching calculus, as we all should. Let us set:
$$I_p=\int_0^{\pi/2}\cos^pt\,dt$$

We compute $I_p$ by partial integration. We have the following formula:
\begin{eqnarray*}
(\cos^pt\sin t)'
&=&p\cos^{p-1}t(-\sin t)\sin t+\cos^pt\cos t\\
&=&p\cos^{p+1}t-p\cos^{p-1}t+\cos^{p+1}t\\
&=&(p+1)\cos^{p+1}t-p\cos^{p-1}t
\end{eqnarray*}

By integrating between $0$ and $\pi/2$, we obtain the following formula:
$$(p+1)I_{p+1}=pI_{p-1}$$

Thus we can compute $I_p$ by recurrence, and we obtain:
\begin{eqnarray*}
I_p
&=&\frac{p-1}{p}\,I_{p-2}\\
&=&\frac{p-1}{p}\cdot\frac{p-3}{p-2}\,I_{p-4}\\
&=&\frac{p-1}{p}\cdot\frac{p-3}{p-2}\cdot\frac{p-5}{p-4}\,I_{p-6}\\
&&\vdots\\
&=&\frac{p!!}{(p+1)!!}\,I_{1-\varepsilon(p)}
\end{eqnarray*}

Together with $I_0=\frac{\pi}{2}$ and $I_1=1$, which are both clear, we obtain:
$$I_p=\left(\frac{\pi}{2}\right)^{\varepsilon(p)}\frac{p!!}{(p+1)!!}$$

Summarizing, we have proved the following formula, with one equality coming from the above computation, and with the other equality coming from this, via $t=\frac{\pi}{2}-s$: 
$$\int_0^{\pi/2}\cos^pt\,dt=\int_0^{\pi/2}\sin^pt\,dt=\left(\frac{\pi}{2}\right)^{\varepsilon(p)}\frac{p!!}{(p+1)!!}$$

In relation with the formula in the statement, we are therefore done with the case $p=0$ or $q=0$. Let us investigate now the general case. We must compute:
$$I_{pq}=\int_0^{\pi/2}\cos^pt\sin^qt\,dt$$

In order to do the partial integration, observe that we have:
\begin{eqnarray*}
(\cos^pt\sin^qt)'
&=&p\cos^{p-1}t(-\sin t)\sin^qt\\
&+&\cos^pt\cdot q\sin^{q-1}t\cos t\\
&=&-p\cos^{p-1}t\sin^{q+1}t+q\cos^{p+1}t\sin^{q-1}t
\end{eqnarray*}

By integrating between $0$ and $\pi/2$, we obtain, for $p,q>0$:
$$pI_{p-1,q+1}=qI_{p+1,q-1}$$

Thus, we can compute $I_{pq}$ by recurrence. When $q$ is even we have:
\begin{eqnarray*}
I_{pq}
&=&\frac{q-1}{p+1}\,I_{p+2,q-2}\\
&=&\frac{q-1}{p+1}\cdot\frac{q-3}{p+3}\,I_{p+4,q-4}\\
&=&\frac{q-1}{p+1}\cdot\frac{q-3}{p+3}\cdot\frac{q-5}{p+5}\,I_{p+6,q-6}\\
&=&\vdots\\
&=&\frac{p!!q!!}{(p+q)!!}\,I_{p+q}
\end{eqnarray*}

But the last term was already computed above, and we obtain the result:
\begin{eqnarray*}
I_{pq}
&=&\frac{p!!q!!}{(p+q)!!}\,I_{p+q}\\
&=&\frac{p!!q!!}{(p+q)!!}\left(\frac{\pi}{2}\right)^{\varepsilon(p+q)}\frac{(p+q)!!}{(p+q+1)!!}\\
&=&\left(\frac{\pi}{2}\right)^{\varepsilon(p)\varepsilon(q)}\frac{p!!q!!}{(p+q+1)!!}
\end{eqnarray*}

Observe that this gives the result for $p$ even as well, by symmetry. Indeed, we have $I_{pq}=I_{qp}$, by using the following change of variables:
$$t=\frac{\pi}{2}-s$$

In the remaining case now, where both $p,q$ are odd, we can use once again the formula $pI_{p-1,q+1}=qI_{p+1,q-1}$ established above, and the recurrence goes as follows:
\begin{eqnarray*}
I_{pq}
&=&\frac{q-1}{p+1}\,I_{p+2,q-2}\\
&=&\frac{q-1}{p+1}\cdot\frac{q-3}{p+3}\,I_{p+4,q-4}\\
&=&\frac{q-1}{p+1}\cdot\frac{q-3}{p+3}\cdot\frac{q-5}{p+5}\,I_{p+6,q-6}\\
&=&\vdots\\
&=&\frac{p!!q!!}{(p+q-1)!!}\,I_{p+q-1,1}
\end{eqnarray*}

In order to compute the last term, observe that we have:
\begin{eqnarray*}
I_{p1}
&=&\int_0^{\pi/2}\cos^pt\sin t\,dt\\
&=&-\frac{1}{p+1}\int_0^{\pi/2}(\cos^{p+1}t)'\,dt\\
&=&\frac{1}{p+1}
\end{eqnarray*}

Thus, we can finish our computation in the case $p,q$ odd, as follows:
\begin{eqnarray*}
I_{pq}
&=&\frac{p!!q!!}{(p+q-1)!!}\,I_{p+q-1,1}\\
&=&\frac{p!!q!!}{(p+q-1)!!}\cdot\frac{1}{p+q}\\
&=&\frac{p!!q!!}{(p+q+1)!!}
\end{eqnarray*}

Thus, we obtain the formula in the statement, the exponent of $\pi/2$ appearing there being $\varepsilon(p)\varepsilon(q)=0\cdot 0=0$ in the present case, and this finishes the proof.
\end{proof}

\section*{16b. More calculus}

In order to discuss the higher spheres, we will use spherical coordinates:

\index{spherical coordinates}
\index{Jacobian}

\begin{theorem}
We have spherical coordinates in $N$ dimensions,
$$\begin{cases}
x_1\!\!\!&=\ r\cos t_1\\
x_2\!\!\!&=\ r\sin t_1\cos t_2\\
\vdots\\
x_{N-1}\!\!\!&=\ r\sin t_1\sin t_2\ldots\sin t_{N-2}\cos t_{N-1}\\
x_N\!\!\!&=\ r\sin t_1\sin t_2\ldots\sin t_{N-2}\sin t_{N-1}
\end{cases}$$
the corresponding Jacobian being given by the following formula,
$$J(r,t)=r^{N-1}\sin^{N-2}t_1\sin^{N-3}t_2\,\ldots\,\sin^2t_{N-3}\sin t_{N-2}$$
and with this generalizing the formulae that we know well at $N=2,3$.
\end{theorem}

\begin{proof}
The fact that we have spherical coordinates as above is clear. Regarding now the Jacobian, by developing the determinant over the last column, we have:
\begin{eqnarray*}
J_N
&=&r\sin t_1\ldots\sin t_{N-2}\sin t_{N-1}\times \sin t_{N-1}J_{N-1}\\
&+&r\sin t_1\ldots \sin t_{N-2}\cos t_{N-1}\times\cos t_{N-1}J_{N-1}\\
&=&r\sin t_1\ldots\sin t_{N-2}(\sin^2 t_{N-1}+\cos^2 t_{N-1})J_{N-1}\\
&=&r\sin t_1\ldots\sin t_{N-2}J_{N-1}
\end{eqnarray*}

Thus, we obtain the formula in the statement, by recurrence.
\end{proof}

With the above results in hand, we can now compute arbitrary polynomial integrals, over the spheres of arbitrary dimension, the result being is as follows:

\index{spherical integral}

\begin{theorem}
The spherical integral of $x_{i_1}\ldots x_{i_k}$ vanishes, unless each index $a\in\{1,\ldots,N\}$ appears an even number of times in the sequence $i_1,\ldots,i_k$. We have 
$$\int_{S^{N-1}_\mathbb R}x_{i_1}\ldots x_{i_k}\,dx=\frac{(N-1)!!l_1!!\ldots l_N!!}{(N+\sum l_i-1)!!}$$
with $l_a$ being this number of occurrences.
\end{theorem}

\begin{proof}
First, the result holds indeed at $N=2$, due to the following formula proved above, where $\varepsilon(p)=1$ when $p\in\mathbb N$ is even, and $\varepsilon(p)=0$ when $p$ is odd:
$$\int_0^{\pi/2}\cos^pt\sin^qt\,dt=\left(\frac{\pi}{2}\right)^{\varepsilon(p)\varepsilon(q)}\frac{p!!q!!}{(p+q+1)!!}$$

In general, we can assume $l_a\in 2\mathbb N$, since the other integrals vanish. The integral in the statement can be written in spherical coordinates, as follows:
$$I=\frac{2^N}{A}\int_0^{\pi/2}\ldots\int_0^{\pi/2}x_1^{l_1}\ldots x_N^{l_N}J\,dt_1\ldots dt_{N-1}$$

In this formula $A$ is the area of the sphere, $J$ is the Jacobian, and the $2^N$ factor comes from the restriction to the $1/2^N$ part of the sphere where all the coordinates are positive. The normalization constant in front of the integral is well-known, given by:
$$\frac{2^N}{A}
=\left(\frac{2}{\pi}\right)^{[N/2]}(N-1)!!$$

As for the unnormalized integral, this is given by:
\begin{eqnarray*}
I'=\int_0^{\pi/2}\ldots\int_0^{\pi/2}
&&(\cos t_1)^{l_1}
(\sin t_1\cos t_2)^{l_2}\\
&&\vdots\\
&&(\sin t_1\sin t_2\ldots\sin t_{N-2}\cos t_{N-1})^{l_{N-1}}\\
&&(\sin t_1\sin t_2\ldots\sin t_{N-2}\sin t_{N-1})^{l_N}\\
&&\sin^{N-2}t_1\sin^{N-3}t_2\ldots\sin^2t_{N-3}\sin t_{N-2}\\
&&dt_1\ldots dt_{N-1}
\end{eqnarray*}

By rearranging the terms, we obtain:
\begin{eqnarray*}
I'
&=&\int_0^{\pi/2}\cos^{l_1}t_1\sin^{l_2+\ldots+l_N+N-2}t_1\,dt_1\\
&&\int_0^{\pi/2}\cos^{l_2}t_2\sin^{l_3+\ldots+l_N+N-3}t_2\,dt_2\\
&&\vdots\\
&&\int_0^{\pi/2}\cos^{l_{N-2}}t_{N-2}\sin^{l_{N-1}+l_N+1}t_{N-2}\,dt_{N-2}\\
&&\int_0^{\pi/2}\cos^{l_{N-1}}t_{N-1}\sin^{l_N}t_{N-1}\,dt_{N-1}
\end{eqnarray*}

Now by using the above-mentioned formula at $N=2$, this gives:
\begin{eqnarray*}
I'
&=&\frac{l_1!!(l_2+\ldots+l_N+N-2)!!}{(l_1+\ldots+l_N+N-1)!!}\left(\frac{\pi}{2}\right)^{\varepsilon(N-2)}\\
&&\frac{l_2!!(l_3+\ldots+l_N+N-3)!!}{(l_2+\ldots+l_N+N-2)!!}\left(\frac{\pi}{2}\right)^{\varepsilon(N-3)}\\
&&\vdots\\
&&\frac{l_{N-2}!!(l_{N-1}+l_N+1)!!}{(l_{N-2}+l_{N-1}+l_N+2)!!}\left(\frac{\pi}{2}\right)^{\varepsilon(1)}\\
&&\frac{l_{N-1}!!l_N!!}{(l_{N-1}+l_N+1)!!}\left(\frac{\pi}{2}\right)^{\varepsilon(0)}
\end{eqnarray*}

By simplifying this, and putting everything together, we obtain:
\begin{eqnarray*}
I
&=&\left(\frac{2}{\pi}\right)^{[N/2]}(N-1)!!\cdot\left(\frac{\pi}{2}\right)^{[N/2]}
\frac{l_1!!\ldots l_N!!}{(N+\sum l_i-1)!!}\\
&=&\frac{(N-1)!!l_1!!\ldots l_N!!}{(N+\sum l_i-1)!!}
\end{eqnarray*}

Thus, we have obtained the formula in the statement.
\end{proof}

In connection now with our probabilistic questions, we have:

\index{hyperspherical laws}
\index{hyperspherical variables}
\index{asymptotic independence}

\begin{theorem}
The even moments of the hyperspherical variables are
$$\int_{S^{N-1}_\mathbb R}x_i^kdx=\frac{(N-1)!!k!!}{(N+k-1)!!}$$
and the variables $y_i=\sqrt{N}x_i$ become normal and independent with $N\to\infty$.
\end{theorem}

\begin{proof}
The moment formula in the statement follows from Theorem 16.3. Now observe that with $N\to\infty$ we have the following estimate:
\begin{eqnarray*}
\int_{S^{N-1}_\mathbb R}x_i^kdx
&=&\frac{(N-1)!!}{(N+k-1)!!}\times k!!\\
&\simeq&N^{-k/2}\times k!!\\
&=&N^{-k/2}M_k(g_1)
\end{eqnarray*}

Thus we have, as claimed, the following asymptotic formula:
$$\sqrt{N}x_i\sim g_1$$

Finally, the independence assertion follows as well from the formula in Theorem 16.3, via some standard probability theory.
\end{proof}

In the case of the half-classical sphere, we have the following result:

\index{half-classical integration}

\begin{theorem}
The half-classical integral of $x_{i_1}\ldots x_{i_k}$ vanishes, unless each index $a$ appears the same number of times at odd and even positions in $i_1,\ldots,i_k$. We have
$$\int_{S^{N-1}_{\mathbb R,*}}x_{i_1}\ldots x_{i_k}\,dx=\frac{(N-1)!l_1!\ldots l_n!}{(N+\sum l_i-1)!}$$
where $l_a$ denotes this number of common occurrences.
\end{theorem}

\begin{proof}
As before, we can assume that $k$ is even, $k=2l$. The corresponding integral can be viewed as an integral over $S^{N-1}_\mathbb C$, as follows:
$$I=\int_{S^{N-1}_\mathbb C}z_{i_1}\bar{z}_{i_2}\ldots z_{i_{2l-1}}\bar{z}_{i_{2l}}\,dz$$

In order to get started, and prove the first assertion, let us apply to this integral transformations of the following type, with $|\lambda|=1$:
$$p\to\lambda p$$

We conclude from this that the above integral $I$ vanishes, unless each $z_a$ appears as many times as $\bar{z}_a$ does, and this gives the first assertion.

\medskip

Assume now that we are in the non-vanishing case. Then the $l_a$ copies of $z_a$ and the $l_a$ copies of $\bar{z}_a$ produce by multiplication a factor $|z_a|^{2l_a}$, so we have:
$$I=\int_{S^{N-1}_\mathbb C}|z_1|^{2l_1}\ldots|z_N|^{2l_N}\,dz$$

Now by using the standard identification $S^{N-1}_\mathbb C\simeq S^{2N-1}_\mathbb R$, we obtain:
\begin{eqnarray*}
I
&=&\int_{S^{2N-1}_\mathbb R}(x_1^2+y_1^2)^{l_1}\ldots(x_N^2+y_N^2)^{l_N}\,d(x,y)\\
&=&\sum_{r_1\ldots r_N}\binom{l_1}{r_1}\ldots\binom{l_N}{r_N}\int_{S^{2N-1}_\mathbb R}x_1^{2l_1-2r_1}y_1^{2r_1}\ldots x_N^{2l_N-2r_N}y_N^{2r_N}\,d(x,y)\\
&=&\sum_{r_1\ldots r_N}\binom{l_1}{r_1}\ldots\binom{l_N}{r_N}\frac{(2N-1)!!(2r_1)!!\ldots(2r_N)!!(2l_1-2r_1)!!\ldots (2l_N-2r_N)!!}{(2N+2\sum l_i-1)!!}\\
&=&\sum_{r_1\ldots r_N}\binom{l_1}{r_1}\ldots\binom{l_N}{r_N}
\frac{(N-1)!(2r_1)!\ldots (2r_N)!(2l_1-2r_1)!\ldots (2l_N-2r_N)!}{4^{\sum l_i}(N+\sum l_i-1)!r_1!\ldots r_N!(l_1-r_1)!\ldots (l_N-r_N)!}
\end{eqnarray*}

Observe now that can rewrite this quantity in the following way:
\begin{eqnarray*}
I
&=&\sum_{r_1\ldots r_N}\frac{l_1!\ldots l_N!(N-1)!(2r_1)!\ldots (2r_N)!(2l_1-2r_1)!\ldots (2l_N-2r_N)!}{4^{\sum l_i}(N+\sum l_i-1)!(r_1!\ldots r_N!(l_1-r_1)!\ldots (l_N-r_N)!)^2}\\
&=&\sum_{r_1}\binom{2r_1}{r_1}\binom{2l_1-2r_1}{l_1-r_1}\ldots\sum_{r_N}\binom{2r_N}{r_N}\binom{2l_N-2r_N}{l_N-r_N}\frac{(N-1)!l_1!\ldots l_N!}{4^{\sum l_i}(N+\sum l_i-1)!}\\
&=&4^{l_1}\times\ldots\times 4^{l_N}\times\frac{(N-1)!l_1!\ldots l_N!}{4^{\sum l_i}(N+\sum l_i-1)!}\\
&=&\frac{(N-1)!l_1!\ldots l_N!}{(N+\sum l_i-1)!}
\end{eqnarray*}

Here we used the following well-known formula, whose proof is elementary:
$$\sum_r\binom{2r}{r}\binom{2l-2r}{l-r}=4^l$$

Thus, we are led to the formula in the statement.
\end{proof}

As before, we can deduce from this a probabilistic result, as follows:

\index{half-classical hyperspherical law}
\index{Rayleigh variable}

\begin{theorem}
The even moments of the half-classical hyperspherical variables are
$$\int_{S^{N-1}_{\mathbb R,*}}x_i^{2l}dx=\frac{(N-1)!l!}{(N+l-1)!}$$
and the variables $y_i=\sqrt{N}x_i$ become symmetrized Rayleigh with $N\to\infty$.
\end{theorem}

\begin{proof}
The moment formula in the statement follows from Theorem 16.5. Now observe that with $N\to\infty$ we have the following estimate:
\begin{eqnarray*}
\int_{S^{N-1}_{\mathbb R,*}}x_i^{2l}dx
&=&\frac{(N-1)!}{(N+l-1)!}\times l!\\
&\simeq&N^{-l}\times l!\\
&=&N^{-l}M_l(|c|)
\end{eqnarray*}

Here $c$ is a standard complex Gaussian variable, and this gives the result.
\end{proof}

As a comment here, it is possible to prove, based once again on the integration formula from Theorem 16.5, that the rescaled variables $y_i=\sqrt{N}x_i$ become ``half-independent'' with $N\to\infty$. For a discussion about half-independence, we refer to \cite{bcs}.

\section*{16c. Advanced calculus}

In the case of the free sphere now, the computations are substantially more complicated than those in the classical and half-classical cases. Let us start with the following result, that we basically know from chapter 5, and that we will recall now:

\index{semicircle law}
\index{asymptotic freeness}

\begin{theorem}
For the free sphere $S^{N-1}_{\mathbb R,+}$, the rescaled coordinates 
$$y_i=\sqrt{N}x_i$$
become semicircular and free, in the $N\to\infty$ limit.
\end{theorem}

\begin{proof}
As explained in chapter 5, the Weingarten formula for the free sphere, together with the standard fact that the Gram matrix, and hence the Weingarten matrix too, is asymptotically diagonal, gives the following estimate:
$$\int_{S^{N-1}_{\mathbb R,+}}x_{i_1}\ldots x_{i_k}\,dx\simeq N^{-k/2}\sum_{\sigma\in NC_2(k)}\delta_\sigma(i_1,\ldots,i_k)$$

With this formula in hand, we can compute the asymptotic moments of each coordinate $x_i$. Indeed, by setting $i_1=\ldots=i_k=i$, all Kronecker symbols are 1, and we obtain:
$$\int_{S^{N-1}_{\mathbb R,+}}x_i^k\,dx\simeq N^{-k/2}|NC_2(k)|$$

Thus the rescaled coordinates $y_i=\sqrt{N}x_i$ become semicircular in the $N\to\infty$ limit, as claimed. As for the asymptotic freeness result, this follows as well from the above general joint moment estimate, via standard free probability theory. See \cite{ba8}, \cite{bgo}.
\end{proof}

Summarizing, we have good results for the free sphere, with $N\to\infty$. The problem now is that of computing the moments of the coordinates of the free sphere at fixed values of $N\in\mathbb N$. The answer here, from \cite{bcz}, which is based on advanced quantum group and calculus techniques, that we will briefly explain here, is as follows:

\index{free hyperspherical law}
\index{special functions}
\index{twisting}

\begin{theorem}
The moments of the free hyperspherical law are given by
$$\int_{S^{N-1}_{\mathbb R,+}}x_1^{2l}\,dx=\frac{1}{(N+2)^l}\cdot\frac{q+1}{q-1}\cdot\frac{1}{l+1}\sum_{r=-l-1}^{l+1}(-1)^r\begin{pmatrix}2l+2\cr l+r+1\end{pmatrix}\frac{r}{1+q^r}$$
where $q\in [-1,0)$ is such that $q+q^{-1}=-N$.
\end{theorem}

\begin{proof}
This is something quite tricky. Following \cite{bcz}, this will follow in 4 steps, none of which is something trivial, which are as follows:

\medskip

(1) $x_1\in C(S^{N-1}_{\mathbb R,+})$ has the same law as $u_{11}\in C(O_N^+)$.

\medskip

(2) $u_{11}\in C(O_N^+)$ has the same law as a certain variable $w\in C(SU^q_2)$.

\medskip

(3) $w\in C(SU^q_2)$ can be in turn modeled by an explicit operator $T\in B(l^2(\mathbb N))$.

\medskip

(4) The law of $T\in B(l^2(\mathbb N))$ can be computed by using advanced calculus.

\medskip

Let us first explain the relation between $O_N^+$ and $SU^q_2$. To any matrix $F\in GL_N(\mathbb R)$ satisfying $F^2=1$ we associate the following universal algebra:
$$C(O_F^+)=C^*\left((u_{ij})_{i,j=1,\ldots,N}\Big|u=F\bar{u}F={\rm unitary}\right)$$

Observe that $O_{I_N}^+=O_N^+$. In general, the above algebra satisfies Woronowicz's generalized axioms in \cite{wo1}, which do not include the strong antipode axiom $S^2=id$.

\medskip

At $N=2$, up to a trivial equivalence relation on the matrices $F$, and on the quantum groups $O_F^+$, we can assume that $F$ is as follows, with $q\in [-1,0)$:
$$F=\begin{pmatrix}0&\sqrt{-q}\\
1/\sqrt{-q}&0\end{pmatrix}$$

Our claim is that for this matrix we have an isomorphism as follows:
$$O_F^+=SU^q_2$$

Indeed, the relations $u=F\bar{u}F$ tell us that $u$ must be of the following form:
$$u=\begin{pmatrix}\alpha&-q\gamma^*\cr \gamma&\alpha^*\end{pmatrix}$$

Thus $C(O_F^+)$ is the universal algebra generated by two elements $\alpha,\gamma$, with the relations making the above matrix $u$ unitary. But these unitarity conditions are:
$$\alpha\gamma=q\gamma\alpha$$
$$\alpha\gamma^*=q\gamma^*\alpha$$
$$\gamma\gamma^*=\gamma^*\gamma$$
$$\alpha^*\alpha+\gamma^*\gamma=1$$
$$\alpha\alpha^*+q^2\gamma\gamma^*=1$$

We recognize here the relations in \cite{wo1} defining the algebra $C(SU^q_2)$, and it follows that we have an isomorphism of Hopf $C^*$-algebras, as follows:
$$C(O_F^+)\simeq C(SU^q_2)$$

Now back to the general case, let us try to understand the integration over $O_F^+$. Given $\pi\in NC_2(2k)$ and $i=(i_1,\ldots,i_{2k})$, we set:
$$\delta_\pi^F(i)=\prod_{s\in\pi}F_{i_{s_l}i_{s_r}}$$

Here the product is over all strings of $\pi$, denoted as follows:
$$s=\{s_l\curvearrowright s_r\}$$

Our claim now is that the following family of vectors, with $\pi\in NC_2(2k)$, spans the space of fixed vectors of $u^{\otimes 2k}$:
$$\xi_\pi=\sum_i\delta_\pi^F(i)e_{i_1}\otimes\ldots\otimes e_{i_{2k}}$$ 

Indeed, having $\xi_\cap$ fixed by $u^{\otimes 2}$ is equivalent to assuming that $u=F\bar{u}F$ is unitary. By using now the above vectors, we obtain the following Weingarten formula:
$$\int_{O_F^+}u_{i_1j_1}\ldots u_{i_{2k}j_{2k}}=\sum_{\pi\sigma}\delta_\pi^F(i)\delta_\sigma^F(j)W_{kN}(\pi,\sigma)$$

With these preliminaries in hand, let us start the computation. Let $N\in\mathbb N$, and consider the number $q\in [-1,0)$ satisfying:
$$q+q^{-1}=-N$$

Our claim is that we have the following formula:
$$\int_{O_N^+}\varphi(\sqrt{N+2}\,u_{ij})=\int_{SU^q_2}\varphi(\alpha+\alpha^*+\gamma-q\gamma^*)$$

Indeed, the moments of the variable on the left are given by:
$$\int_{O_N^+}u_{ij}^{2k}=\sum_{\pi\sigma}W_{kN}(\pi,\sigma)$$

On the other hand, the moments of the variable on the right, which in terms of the fundamental corepresentation $v=(v_{ij})$ is given by $w=\sum_{ij}v_{ij}$, are given by:
$$\int_{SU^q_2}w^{2k}=\sum_{ij}\sum_{\pi\sigma}\delta_\pi^F(i)\delta_\sigma^F(j)W_{kN}(\pi,\sigma)$$

We deduce that $w/\sqrt{N+2}$ has the same moments as $u_{ij}$, which proves our claim. In order to do now the computation over $SU^q_2$, we can use a matrix model due to Woronowicz \cite{wo1}, where the standard generators $\alpha,\gamma$ are mapped as follows:
\begin{eqnarray*}
\pi_u(\alpha)e_k&=&\sqrt{1-q^{2k}}e_{k-1}\\
\pi_u(\gamma)e_k&=&uq^k e_k
\end{eqnarray*}

Here $u\in\mathbb T$ is a parameter, and $(e_k)$ is the standard basis of $l^2(\mathbb N)$. The point with this representation is that it allows the computation of the Haar functional. Indeed, if $D$ is the diagonal operator given by $D(e_k)=q^{2k}e_k$, then the formula is as follows:
$$\int _{SU^q_2}x=(1-q^2)\int_{\mathbb T}tr(D\pi_u(x))\frac{du}{2\pi iu}$$

With the above model in hand, the law of the variable that we are interested in is of the following form:
$$\int_{SU^q_2}\varphi(\alpha+\alpha^*+\gamma-q\gamma^*)=(1-q^2)\int_{\mathbb T}tr(D\varphi(M))\frac{du}{2\pi iu}$$

To be more precise, this formula holds indeed, with:
$$M(e_k)=e_{k+1}+q^k(u-qu^{-1})e_k+(1-q^{2k})e_{k-1}$$

The point now is that the integral on the right can be computed, by using advanced calculus methods, and this gives the result. We refer here to \cite{bcz}. 
\end{proof}

The computation of the joint free hyperspherical laws remains an open problem. Open as well is the question of finding a more conceptual proof for the above formula.

\section*{16d. Twisting results}

Following now \cite{bbs}, let us discuss an interesting relation of all this with the quantum permutations, and with the free hypergeometric laws. The idea will be that of working out some abstract algebraic results, regarding twists of quantum automorphism groups, which will particularize into results relating quantum rotations and permutations, having no classical counterpart, both at the algebraic and the probabilistic level.

\bigskip

In order to explain this material, from \cite{bbs}, which is quite technical, requiring good algebraic knowledge, let us begin with some generalities. We first have:

\index{finite quantum space}
\index{counting measure}
\index{canonical trace}

\begin{definition}
A finite quantum space $Z$ is the abstract dual of a finite dimensional $C^*$-algebra $B$, according to the following formula:
$$C(Z)=B$$
The number of elements of such a space is $|Z|=\dim B$. By decomposing the algebra $B$, we have a formula of the following type:
$$C(Z)=M_{n_1}(\mathbb C)\oplus\ldots\oplus M_{n_k}(\mathbb C)$$
With $n_1=\ldots=n_k=1$ we obtain in this way the space $Z=\{1,\ldots,k\}$. Also, when $k=1$ the equation is $C(Z)=M_n(\mathbb C)$, and the solution will be denoted $Z=M_n$.
\end{definition}

Following \cite{ba8}, we endow each finite quantum space $Z$ with its counting measure, corresponding as the algebraic level to the integration functional obtained by applying the regular representation, and then the normalized matrix trace:
$$tr:C(Z)\to B(l^2(Z))\to\mathbb C$$

As basic examples, for both $Z=\{1,\ldots,k\}$ and $Z=M_n$ we obtain the usual trace. In general, we can write the algebra $C(Z)$ as follows:
$$C(Z)=M_{n_1}(\mathbb C)\oplus\ldots\oplus M_{n_k}(\mathbb C)$$

In terms of this writing, the weights of $tr$ are as follows:
$$c_i=\frac{n_i^2}{\sum_in_i^2}$$

Let us study now the quantum group actions $G\curvearrowright Z$. It is convenient here to use, in order to get started, the no basis approach from \cite{ba8}. If we denote by $\mu,\eta$ the multiplication and unit map of the algebra $C(Z)$, we have the following result, from \cite{ba8}:

\begin{proposition}
Consider a linear map $\Phi:C(Z)\to C(Z)\otimes C(G)$, written as
$$\Phi(e_i)=\sum_je_j\otimes u_{ji}$$
with $\{e_i\}$ being a linear space basis of $C(Z)$, orthonormal with respect to $tr$.
\begin{enumerate}
\item $\Phi$ is a linear space coaction $\iff$ $u$ is a corepresentation.

\item $\Phi$ is multiplicative $\iff$ $\mu\in Hom(u^{\otimes 2},u)$.

\item $\Phi$ is unital $\iff$ $\eta\in Hom(1,u)$.

\item $\Phi$ leaves invariant $tr$ $\iff$ $\eta\in Hom(1,u^*)$.

\item If these conditions hold, $\Phi$ is involutive $\iff$ $u$ is unitary.
\end{enumerate}
\end{proposition}

\begin{proof}
This is a bit similar to the proof for $S_N^+$ from chapter 2, as follows:

\medskip

(1) There are two axioms to be processed here. First, we have:
\begin{eqnarray*}
(id\otimes\Delta)\Phi=(\Phi\otimes id)\Phi
&\iff&\sum_je_j\otimes\Delta(u_{ji})=\sum_k\Phi(e_k)\otimes u_{ki}\\
&\iff&\sum_je_j\otimes\Delta(u_{ji})=\sum_{jk}e_j\otimes u_{jk}\otimes u_{ki}\\
&\iff&\Delta(u_{ji})=\sum_ku_{jk}\otimes u_{ki}
\end{eqnarray*}

As for the axiom involving the counit, here we have as well, as desired:
\begin{eqnarray*}
(id\otimes\varepsilon)\Phi=id
&\iff&\sum_j\varepsilon(u_{ji})e_j=e_i\\
&\iff&\varepsilon(u_{ji})=\delta_{ji}
\end{eqnarray*}

(2) We have the following formula:
$$\Phi(e_i)
=\left(\sum_{ij}e_{ji}\otimes u_{ji}\right)(e_i\otimes 1)
=u(e_i\otimes 1)$$

By using this formula, we obtain the following identity:
$$\Phi(e_ie_k)
=u(e_ie_k\otimes 1)
=u(\mu\otimes id)(e_i\otimes e_k\otimes 1)$$

On the other hand, we have as well the following identity, as desired:
\begin{eqnarray*}
\Phi(e_i)\Phi(e_k)
&=&\sum_{jl}e_je_l\otimes u_{ji}u_{lk}\\
&=&(\mu\otimes id)\sum_{jl}e_j\otimes e_l\otimes u_{ji}u_{lk}\\\
&=&(\mu\otimes id)\left(\sum_{ijkl}e_{ji}\otimes e_{lk}\otimes u_{ji}u_{lk}\right)(e_i\otimes e_k\otimes 1)\\
&=&(\mu\otimes id)u^{\otimes 2}(e_i\otimes e_k\otimes 1)
\end{eqnarray*}

(3) The formula $\Phi(e_i)=u(e_i\otimes1)$ found above gives by linearity $\Phi(1)=u(1\otimes1)$. But this shows that $\Phi$ is unital precisely when $u(1\otimes1)=1\otimes1$, as desired.

\medskip

(4) This follows from the following computation, by applying the involution:
\begin{eqnarray*}
(tr\otimes id)\Phi(e_i)=tr(e_i)1
&\iff&\sum_jtr(e_j)u_{ji}=tr(e_i)1\\
&\iff&\sum_ju_{ji}^*1_j=1_i\\
&\iff&(u^*1)_i=1_i\\
&\iff&u^*1=1
\end{eqnarray*}

(5) Assuming that (1-4) are satisfied, and that $\Phi$ is involutive, we have:
\begin{eqnarray*}
(u^*u)_{ik}
&=&\sum_lu_{li}^*u_{lk}\\
&=&\sum_{jl}tr(e_j^*e_l)u_{ji}^*u_{lk}\\
&=&(tr\otimes id)\sum_{jl}e_j^*e_l\otimes u_{ji}^*u_{lk}\\
&=&(tr\otimes id)(\Phi(e_i)^*\Phi(e_k))\\
&=&(tr\otimes id)\Phi(e_i^*e_k)\\
&=&tr(e_i^*e_k)1\\
&=&\delta_{ik}
\end{eqnarray*}

Thus $u^*u=1$, and since we know from (1) that $u$ is a corepresentation, it follows that $u$ is unitary. The proof of the converse is standard too, by using similar tricks.
\end{proof}

Following now \cite{ba8}, \cite{wa2}, we have the following result, extending the basic theory of $S_N^+$ from chapter 2 to the present finite quantum space setting:

\index{quantum automorphism group}
\index{quantum symmetry group}

\begin{theorem}
Given a finite quantum space $Z$, there is a universal compact quantum group $S_Z^+$ acting on $Z$, leaving the counting measure invariant. We have
$$C(S_Z^+)=C(U_N^+)\Big/\Big<\mu\in Hom(u^{\otimes2},u),\eta\in Fix(u)\Big>$$
where $N=|Z|$ and where $\mu,\eta$ are the multiplication and unit maps of $C(Z)$. Also:
\begin{enumerate}
\item For $Z=\{1,\ldots,N\}$ we have $S_Z^+=S_N^+$.

\item For $Z=M_n$ we have $S_Z^+=PO_n^+=PU_n^+$.
\end{enumerate}
\end{theorem}

\begin{proof}
Consider a linear map $\Phi:C(Z)\to C(Z)\otimes C(G)$, written as follows, with $\{e_i\}$ being a linear space basis of $C(Z)$, which is orthonormal with respect to $tr$:
$$\Phi(e_j)=\sum_ie_i\otimes u_{ij}$$

It is routine to check, via standard algebraic computations, that $\Phi$ is a coaction precisely when $u$ is a unitary corepresentation, satisfying the following conditions:
$$\mu\in Hom(u^{\otimes2},u)$$
$$\eta\in Fix(u)$$

But this gives the first assertion. Regarding now the statement about $Z=\{1,\ldots,N\}$ is clear. Finally, regarding $Z=M_2$, here we have embeddings as followss:
$$PO_n^+\subset PU_n^+\subset S_Z^+$$

Now since the fusion rules of all these 3 quantum groups are known to be the same as the fusion rules for $SO_3$, these inclusions are isomorphisms. See \cite{ba8}.
\end{proof}

We have as well the following result, also explained in \cite{ba8}:

\begin{theorem}
The quantum groups $S_Z^+$ have the following properties: 
\begin{enumerate}
\item The associated Tannakian categories are $TL(N)$, with $N=|Z|$.

\item The main character follows the Marchenko-Pastur law $\pi_1$, when $N\geq4$.

\item The fusion rules for $S_Z^+$ with $|Z|\geq4$ are the same as for $SO_3$.
\end{enumerate}
\end{theorem}

\begin{proof}
This result is discussed in detail in \cite{ba8}, the idea being as follows:

\medskip

(1) Our first claim is that the fundamental representation is equivalent to its adjoint, $u\sim\bar{u}$. Indeed, let us go back to the coaction formula from Proposition 16.10:
$$\Phi(e_i)=\sum_je_j\otimes u_{ji}$$

We can pick our orthogonal basis $\{e_i\}$ to be the stadard multimatrix basis of $C(Z)$, so that we have, for a certain involution $i\to i^*$ on the index set:
$$e_i^*=e_{i^*}$$

With this convention made, by conjugating the above formula of $\Phi(e_i)$, we obtain:
$$\Phi(e_{i^*})=\sum_je_{j^*}\otimes u_{ji}^*$$

Now by interchanging $i\leftrightarrow i^*$ and $j\leftrightarrow j^*$, this latter formula reads:
$$\Phi(e_i)=\sum_je_j\otimes u_{j^*i^*}^*$$

We therefore conclude, by comparing with the original formula, that we have:
$$u_{ji}^*=u_{j^*i^*}$$

But this shows that we have an equivalence as follows, as claimed:
$$u\sim\bar{u}$$

Now with this result in hand, the proof goes as for the proof for $S_N^+$, from the previous section. To be more precise, the result follows from the fact that the multiplication and unit of any complex algebra, and in particular of the algebra $C(Z)$ that we are interested in here, can be modelled by the following two diagrams:
$$m=|\cup|\qquad,\qquad u=\cap$$

Indeed, this is certainly true algebrically, and this is something well-known. As in what regards the $*$-structure, things here are fine too, because our choice for the trace leads to the following formula, which must be satisfied as well:
$$\mu\mu^*=N\cdot id$$

But the above diagrams $m,u$ generate the Temperley-Lieb algebra $TL(N)$, as stated.

\medskip

(2) The proof here is exactly as for $S_N^+$, by using moments. To be more precise, according to (1) these moments are the Catalan numbers, which are the moments of $\pi_1$.

\medskip

(3) Once again same proof as for $S_N^+$, by using the fact that the moments of $\chi$ are the Catalan numbers, which naturally leads to the Clebsch-Gordan rules.
\end{proof}

Let us discuss now a number of more advanced twisting aspects, which will eventually lead us into probability, and hypergeometric laws. Following \cite{bbs}, we first have:

\begin{proposition}
Given a finite group $G$, the algebra $C(S_{\widehat{G}}^+)$ is isomorphic to the abstract algebra presented by generators $x_{gh}$ with $g,h\in G$, with the following relations:
$$x_{1g}=x_{g1}=\delta_{1g}$$
$$x_{s,gh}=\sum_{t\in G}x_{st^{-1},g}x_{th}$$
$$x_{gh,s}=\sum_{t\in G}x_{gt^{-1}}x_{h,ts}$$
The comultiplication, counit and antipode are given by the formulae
$$\Delta(x_{gh})=\sum_{s\in G}x_{gs}\otimes x_{sh}$$
$$\varepsilon(x_{gh})=\delta_{gh}$$
$$S(x_{gh})=x_{h^{-1}g^{-1}}$$
on the standard generators $x_{gh}$.
\end{proposition}

\begin{proof}
This follows indeed from a direct verification, based either on Theorem 16.11, or on its equivalent formulation from Wang's paper \cite{wa2}.
\end{proof}

\index{twisting}
\index{cocycle twisting}

Let us discuss now the twisted version of the above result. Consider a 2-cocycle on $G$, which is by definition a map $\sigma:G\times G\to\mathbb C^*$ satisfying:
$$\sigma_{gh,s}\sigma_{gh}=\sigma_{g,hs}\sigma_{hs}$$
$$\sigma_{g1}=\sigma_{1g}=1$$

Given such a cocycle, we can construct the associated twisted group algebra $C(\widehat{G}_\sigma)$, as being the vector space $C(\widehat{G})=C^*(G)$, with product as follows:
$$e_ge_h=\sigma_{gh}e_{gh}$$ 

We have then the following generalization of Proposition 16.13:

\begin{proposition}
The algebra $C(S_{\widehat{G}_\sigma}^+)$ is isomorphic to the abstract algebra presented by generators $x_{gh}$ with $g,h\in G$, with the relations $x_{1g}=x_{g1}=\delta_{1g}$ and:
$$\sigma_{gh}x_{s,gh}=\sum_{t\in G}\sigma_{st^{-1},t}x_{st^{-1},g}x_{th}$$
$$\sigma_{gh}^{-1}x_{gh,s}=\sum_{t\in G}\sigma_{t^{-1},ts}^{-1}x_{gt^{-1}}x_{h,ts}$$
The comultiplication, counit and antipode are given by the formulae
$$\Delta(x_{gh})=\sum_{s\in G}x_{gs}\otimes x_{sh}$$
$$\varepsilon(x_{gh})=\delta_{gh}$$
$$S(x_{gh})=\sigma_{h^{-1}h}\sigma_{g^{-1}g}^{-1}x_{h^{-1}g^{-1}}$$
on the standard generators $x_{gh}$.
\end{proposition}

\begin{proof}
Once again, this follows from a direct verification. See \cite{bbs}.
\end{proof}

We prove now that the quantum groups $S_{\widehat{G}}^+$ and $S_{\widehat{G}_\sigma}^+$ are related by a cocycle twisting operation. Let us begin with some preliminaries. Let $A$ be a Hopf algebra. We recall that a left 2-cocycle is a convolution invertible linear map
$\sigma:A\otimes A\to\mathbb C$ satisfying:
$$\sigma_{x_1y_1}\sigma_{x_2y_2,z}=\sigma_{y_1z_1}\sigma_{x,y_2z_2}$$
$$\sigma_{x1}=\sigma_{1x}=\varepsilon(x)$$

Note that $\sigma$ is a left 2-cocycle if and only if $\sigma^{-1}$, the convolution inverse of $\sigma$, is a right 2-cocycle, in the sense that we have:
$$\sigma^{-1}_{x_1y_1,z}\sigma^{-1}_{x_1y_2}=\sigma^{-1}_{x,y_1z_1}\sigma^{-1}_{y_2z_2}$$
$$\sigma^{-1}_{x1}=\sigma^{-1}_{1x}=\varepsilon(x)$$

Given a left 2-cocycle $\sigma$ on $A$, one can form the 2-cocycle twist $A^\sigma$ as follows. As a coalgebra, $A^\sigma=A$, and an element $x\in A$, when considered in $A^\sigma$, is denoted $[x]$. The product in $A^\sigma$ is defined, in Sweedler notation, by: 
$$[x][y]=\sum\sigma_{x_1y_1}\sigma^{-1}_{x_3y_3}[x_2y_2]$$

Note that the cocycle condition ensures the fact that we have indeed a Hopf algebra. With this convention, still following \cite{bbs}, we have the following result:

\begin{theorem}
If $G$ is a finite group and $\sigma$ is a $2$-cocycle on $G$, the Hopf algebras
$$C(S_{\widehat{G}}^+)\quad,\quad C(S_{\widehat{G}_\sigma}^+)$$
are $2$-cocycle twists of each other, in the above sense.
\end{theorem}

\begin{proof}
In order to prove this result, we use the following Hopf algebra map: 
$$\pi:C(S_{\widehat{G}}^+)\to C(\widehat{G})\quad,\quad 
x_{gh}\to\delta_{gh}e_g$$

Our 2-cocycle $\sigma:G\times G\to\mathbb C^*$ can be extended by linearity into a linear map as follows, which is a left and right 2-cocycle in the above sense:
$$\sigma:C(\widehat{G})\otimes C(\widehat{G})\to\mathbb C$$

Consider now the following composition:
$$\alpha=\sigma(\pi \otimes \pi):C(S_{\widehat{G}}^+)\otimes C(S_{\widehat{G}}^+)\to C(\widehat{G})\otimes C(\widehat{G})\to\mathbb C$$

Then $\alpha$ is a left and right 2-cocycle, because it is induced by a cocycle on a group algebra, and so is its convolution inverse $\alpha^{-1}$. Thus we can construct the twisted algebra $C(S_{\widehat{G}}^+)^{\alpha^{-1}}$, and inside this algebra we have the following computation:
\begin{eqnarray*}
[x_{gh}][x_{rs}]
&=&\alpha^{-1}(x_g,x_r)\alpha(x_h,x_s)[x_{gh}x_{rs}]\\
&=&\sigma_{gr}^{-1}\sigma_{hs}[x_{gh}x_{rs}]
\end{eqnarray*}

By using this, we obtain the following formula:
\begin{eqnarray*}
\sum_{t\in G}\sigma_{st^{-1},t}[x_{st^{-1},g}][x_{th}]
&=&\sum_{t\in G}\sigma_{st^{-1},t}\sigma_{st^{-1},t}^{-1}\sigma_{gh}[x_{st^{-1},g}x_{th}]\\
&=&\sigma_{gh}[x_{s,gh}]
\end{eqnarray*}

Similarly, we have the following formula:
$$\sum_{t\in G}\sigma_{t^{-1},ts}^{-1}[x_{g,t^{-1}}][x_{h,ts}]=\sigma_{gh}^{-1}[x_{gh,s}]$$

We deduce from this that there exists a Hopf algebra map, as follows:
$$\Phi:C(S_{\widehat{G}_\sigma}^+)\to C(S_{\widehat{G}}^+)^{\alpha^{-1}}\quad,\quad 
x_{gh}\to [x_{g,h}]$$

This map is clearly surjective, and is injective as well, by a standard fusion semiring argument, because both Hopf algebras have the same fusion semiring.
\end{proof}

Summarizing, we have proved our main twisting result. Our purpose in what follows will be that of working out versions and particular cases of it. We first have: 

\begin{proposition}
If $G$ is a finite group and $\sigma$ is a $2$-cocycle on $G$, then
$$\Phi(x_{g_1h_1}\ldots x_{g_mh_m})=\Omega(g_1,\ldots,g_m)^{-1}\Omega(h_1,\ldots,h_m)x_{g_1h_1}\ldots x_{g_mh_m}$$
with the coefficients on the right being given by the formula
$$\Omega(g_1,\ldots,g_m)=\prod_{k=1}^{m-1}\sigma_{g_1\ldots g_k,g_{k+1}}$$
is a coalgebra isomorphism $C(S_{\widehat{G}_\sigma}^+)\to C(S_{\widehat{G}}^+)$, commuting with the Haar integrals.
\end{proposition}

\begin{proof}
This is indeed just a technical reformulation of Theorem 16.15.
\end{proof}

Here is another useful result, also from \cite{bbs}, that we will need in what follows:

\begin{theorem}
Let $X\subset G$ be such that $\sigma_{gh}=1$ for any $g,h\in X$, and consider the subalgebra 
$$B_X\subset C(S_{\widehat{G}_\sigma}^+)$$
generated by the elements $x_{gh}$, with $g,h\in X$. Then we have an injective algebra map 
$$\Phi_0:B_X\to C(S_{\widehat{G}}^+)$$
given by $x_{g,h}\to x_{g,h}$.
\end{theorem}

\begin{proof}
With the notations in the proof of Theorem 16.15, we have the following equality in $C(S_{\widehat{G}}^+)^{\alpha^{-1}}$, for any $g_i,h_i,r_i,s_i\in X$:
$$[x_{g_1h_1}\ldots x_{g_ph_p}] \cdot [x_{r_1s_1}\ldots x_{r_qs_q}] 
= [x_{g_1h_1}\ldots x_{g_ph_p}x_{r_1s_1}\ldots x_{r_qs_q}]$$

The point now is that $\Phi_0$ can be defined to be the composition of $\Phi_{|B_X}$ with the following linear isomorphism:
$$C(S_{\widehat{G}}^+)^{\alpha^{-1}}\to C(S_{\widehat{G}}^+)$$
$$[x]\to x$$

This being clearly an injective algebra map, we obtain the result.   
\end{proof}

Let us discuss now some concrete applications of the general results established above. Consider the group $G=\mathbb Z_n^2$, let $w=e^{2\pi i/n}$, and consider the following map: 
$$\sigma:G\times G\to\mathbb C^*$$
$$\sigma_{(ij)(kl)}=w^{jk}$$ 

It is easy to see that $\sigma$ is a bicharacter, and hence a 2-cocycle on $G$. Thus, we can apply our general twisting result, to this situation. In order to understand what is the formula that we obtain, we must do some computations. Following \cite{bbs} as before, let $E_{ij}$ with $i,j \in\mathbb Z_n$ be the standard basis of $M_n(\mathbb C)$. We have the following result:

\begin{proposition}
The linear map given by
$$\psi(e_{(i,j)})=\sum_{k=0}^{n-1}{w}^{ki}E_{k,k+j}$$
defines an isomorphism of algebras $\psi:C(\widehat{G}_\sigma)\simeq M_n(\mathbb C)$. 
\end{proposition}

\begin{proof}
Consider indeed the following linear map:
$$\psi'(E_{ij})=\frac{1}{n}\sum_{k=0}^{n-1}{w}^{-ik}e_{(k,j-i)}$$
 
It is routine then to check that $\psi,\psi'$ are inverse morphisms of algebras.
\end{proof}

As a consequence, we have the following result:

\begin{proposition}
The algebra map given by
$$\varphi(u_{ij}u_{kl}) = \frac{1}{n}\sum_{a,b=0}^{n-1}{w}^{ai-bj}x_{(a,k-i),(b,l-j)}$$
defines a Hopf algebra isomorphism $\varphi:C(S_{M_n}^+)\simeq C(S_{\widehat{G}_\sigma}^+)$.
\end{proposition}

\begin{proof}
We use the identification $C(\widehat{G}_\sigma)\simeq M_n(\mathbb C)$ from Proposition 16.18. This identification produces a coaction map, as follows:
$$\gamma:M_n(\mathbb C)\to M_n(\mathbb C)\otimes C(S_{\widehat{G}_\sigma}^+)$$

Now observe that this map is given by the following formula:
 $$\gamma(E_{ij})=\frac{1}{n}\sum_{ab}E_{ab}\otimes\sum_{kr}w^{ar-ik} x_{(r,b-a),(k,j-i)}$$

Thus, we obtain the isomorphism in the statement.
\end{proof}

We will need one more result of this type, as follows:

\begin{proposition}
The algebra map given by
$$\rho(x_{(a,b),(i,j)})=\frac{1}{n^2}\sum_{klrs}w^{ki+lj-ra-sb}p_{(r,s),(k,l)}$$
defines a Hopf algebra isomorphism $\rho:C(S_{\widehat{G}}^+)\simeq C(S_G^+)$.
\end{proposition}

\begin{proof}
This follows by using the Fourier transform isomorphism over the group $G$, which is a map as follows:
$$C(\widehat{G})\simeq C(G)$$

Indeed, by composing with this isomorphism, we obtain the result.
\end{proof}

We can now formulate a concrete twisting result, from \cite{bbs}, as follows:

\begin{theorem}
Let $n\geq 2$ and $w=e^{2\pi i/n}$. Then
$$\Theta(u_{ij}u_{kl})=\frac{1}{n}\sum_{ab=0}^{n-1}w^{-a(k-i)+b(l-j)}p_{ia,jb}$$
defines a coalgebra isomorphism 
$$C(PO_n^+)\to C(S_{n^2}^+)$$
commuting with the Haar integrals.
\end{theorem}
 
\begin{proof}
We recall from Theorem 16.11 that we have identifications as follows:
$$PO_n^+=PU_n^+=S_{M_n}^+$$

With this in hand, the result follows from Theorem 16.15 and Proposition 16.16, by combining them with the various isomorphisms established above.
\end{proof}

Here is a useful version of the above result:

\begin{theorem}
The following two algebras are isomorphic, via $u_{ij}^2\to X_{ij}$:
\begin{enumerate}
\item The algebra generated by the variables $u_{ij}^2\in C(O_n^+)$.

\item The algebra generated by $X_{ij}=\frac{1}{n}\sum_{a,b=1}^np_{ia,jb}\in C(S_{n^2}^+)$
\end{enumerate}
\end{theorem}

\begin{proof}
We have $\Theta(u_{ij}^2)=X_{ij}$, so it remains to prove that if $B$ is the subalgebra of $C(S_{M_n}^+)$ generated by the variables $u_{ij}^2$, then $\Theta_{|B}$ is an algebra morphism. Let us set:
$$X=\{(i,0)|i\in\mathbb Z_n\}\subset\mathbb Z_n^2$$

Then $X$ satisfies the assumption in Theorem 16.17, and $\varphi(B) \subset B_X$. Thus by Theorem 16.17, the map $\Theta_{|B}=\rho F_0\varphi_{|B}$ is indeed an algebra morphism.
\end{proof}

As a probabilistic consequence now, we have:

\begin{theorem}
The following families of variables have the same joint law,
\begin{enumerate}
\item $\{u_{ij}^2\}\in C(O_n^+)$,

\item $\{X_{ij}=\frac{1}{n}\sum_{ab}p_{ia,jb}\}\in C(S_{n^2}^+)$,
\end{enumerate}
where $u=(u_{ij})$ and $p=(p_{ia,jb})$ are the corresponding fundamental corepresentations.
\end{theorem}

\begin{proof}
As explained in \cite{bbs}, this result follows from Theorem 16.22. An alternative approach, also from \cite{bbs}, which is instructive, and that we will excplain now, is by using the Weingarten formula for our two quantum groups, and the shrinking of partitions $\pi\to\pi'$. Let us recall indeed that we have a standard bijection, as follows:
$$NC(k)\simeq NC_2(2k)$$

To be more precise, the application $NC(k)\to NC_2(2k)$ is the ``fattening'' one, obtained by doubling all the legs, and doubling all the strings as well, and its inverse $NC_2(2k)\to NC(k)$ is the ``shrinking'' application, obtained by collapsing pairs of consecutive neighbors. Now back to our questions, observe that we have:
\begin{eqnarray*}
\int_{O_n^+}u_{ij}^{2k}&=&\sum_{\pi,\sigma\in NC_2(2k)}W_{2k,n}(\pi,\sigma)\\
\int_{S_{n^2}^+}X_{ij}^k&=&\sum_{\pi,\sigma\in NC_2(2k)}n^{|\pi'|+|\sigma'|-k}W_{k,n^2}(\pi',\sigma')
\end{eqnarray*}

The point now is that, in the context of the general fattening and shrinking bijection explained above, it is elementary to see that we have:
$$|\pi\vee\sigma|=k+2|\pi'\vee\sigma'|-|\pi'|-|\sigma'|$$

We therefore have the following formula, valid for any $n\in\mathbb N$:
$$n^{|\pi\vee\sigma|}=n^{k+2|\pi'\vee\sigma'|-|\pi'|-|\sigma'|}$$

Thus in our moment formulae above the summands coincide, and so the moments are equal, as desired. The proof in general, dealing with joint moments, is similar.
\end{proof}

In particular, we have the following result:

\index{free hyperspherical law}
\index{free hypergeometric law}
\index{hypergeometric law}

\begin{theorem}
The free hypergeometric variable
$$X_{ij}=\frac{1}{n}\sum_{a,b=1}^nu_{ia,jb}\in C(S_{n^2}^+)$$
has the same law as the squared free hyperspherical variable,
$$x_i^2\in C(S^{N-1}_{\mathbb R,+})$$
and with this being something having no classical counterpart.
\end{theorem}

\begin{proof}
This follows indeed from Theorem 16.23. See \cite{bbs}.
\end{proof}

We refer as well to  \cite{bbs}, \cite{bcs}, \cite{dif} and related papers for some further computations of this type, which are more advanced, involving this time Gram matrix determinants, and for comments, regarding the relevance of such questions. There is a lot of work to be done here, in relation with physics, virtually for everyone interested. 

\bigskip

In what concerns us, our plan is to explain some of these things, and other applications of the nocommutative geometry theory developed in this book to physics, in a series of forthcoming books, dealing with quantum mechanics, and statistical mechanics.

\bigskip

As a conclusion, there is a lot of interesting mathematics in relation with the free spheres and orthogonal groups, and with the quantum permutations and reflections as well. This tends to confirm our intial thought, from the beginning of this book, that the study and axiomatization of the quadruplets $(S,T,U,K)$ is a good question.

\section*{16e. Exercises}

Congratulations for having read this book, and for having survived our various comments, pieces of advice, and of course exercise sessions. Thus, relax and enjoy. However, talking noncommutative geometry, we would have one last exercise, as follows:

\begin{exercise}
Find the free analogue of the stereographic projection.
\end{exercise}

The point indeed is that modern geometry as we know it comes from Riemann, and in his Habilitation, written old style, there is exactly 1 mathematical formula, in relation with the stereographic projection. We believe that looking for free analogues of such things is an interesting question. To be added to other questions raised in this book.

\baselineskip=14pt

\printindex


\begin{thebibliography}{99}

\baselineskip=13.3pt

\bibitem{ar1}V.I. Arnold, Ordinary differential equations, Springer (1973).

\bibitem{ar2}V.I. Arnold, Mathematical methods of classical mechanics, Springer (1974).

\bibitem{ar3}V.I. Arnold, Lectures on partial differential equations, Springer (1997).

\bibitem{ati}M.F. Atiyah, K-theory, CRC Press (1964).

\bibitem{ama}M.F. Atiyah and I.G. MacDonald, Introduction to commutative algebra, Addison-Wesley (1969).

\bibitem{ba1}T. Banica, Liberations and twists of real and complex spheres, {\em J. Geom. Phys.} {\bf 96} (2015), 1--25.

\bibitem{ba2}T. Banica, Quantum isometries of noncommutative polygonal spheres, {\em M\"unster J. Math.} {\bf 8} (2015), 253--284. 

\bibitem{ba3}T. Banica, A duality principle for noncommutative cubes and spheres, {\em J. Noncommut. Geom.} {\bf 10} (2016), 1043--1081.

\bibitem{ba4}T. Banica, Half-liberated manifolds, and their quantum isometries, {\em Glasg. Math. J.} {\bf 59} (2017), 463--492. 

\bibitem{ba5}T. Banica, Liberation theory for noncommutative homogeneous spaces, {\em Ann. Fac. Sci. Toulouse Math.} {\bf 26} (2017), 127--156. 

\bibitem{ba6}T. Banica, Weingarten integration over noncommutative homogeneous spaces, {\em Ann. Math. Blaise Pascal} {\bf 24} (2017), 195--224.

\bibitem{ba7}T. Banica, Principles of operator algebras (2024).

\bibitem{ba8}T. Banica, Introduction to quantum groups, Springer (2023).

\bibitem{bb+}T. Banica, S.T. Belinschi, M. Capitaine and B. Collins, Free Bessel laws, {\em Canad. J. Math.} {\bf 63} (2011), 3--37.

\bibitem{bb1}T. Banica and J. Bichon, Matrix models for noncommutative algebraic manifolds, {\em J. Lond. Math. Soc.} {\bf 95} (2017), 519--540.

\bibitem{bb2}T. Banica and J. Bichon, Complex analogues of the half-classical geometry, {\em M\"unster J. Math.} {\bf 10} (2017), 457--483. 

\bibitem{bbc}T. Banica, J. Bichon and B. Collins, The hyperoctahedral quantum group, {\em J. Ramanujan Math. Soc.} {\bf 22} (2007), 345--384.

\bibitem{bc+}T. Banica, J. Bichon, B. Collins and S. Curran, A maximality result for orthogonal quantum groups, {\em Comm. Algebra} {\bf 41} (2013), 656--665.

\bibitem{bbs}T. Banica, J. Bichon and S. Curran, Quantum automorphisms of twisted group algebras and free hypergeometric laws, {\em Proc. Amer. Math. Soc.} {\bf 139} (2011), 3961--3971.

\bibitem{bcz}T. Banica, B. Collins and P. Zinn-Justin, Spectral analysis of the free orthogonal matrix, {\em Int. Math. Res. Not.} {\bf 17} (2009), 3286--3309.

\bibitem{bcs}T. Banica, S. Curran and R. Speicher, Classification results for easy quantum groups, {\em Pacific J. Math.} {\bf 247} (2010), 1--26.

\bibitem{bgo}T. Banica and D. Goswami, Quantum isometries and noncommutative spheres, {\em Comm. Math. Phys.} {\bf 298} (2010), 343--356.

\bibitem{bme}T. Banica and S. M\'esz\'aros, Uniqueness results for noncommutative spheres and projective spaces, {\em Illinois J. Math.} {\bf 59} (2015), 219--233.

\bibitem{bsk}T. Banica and A. Skalski, Quantum symmetry groups of C*-algebras equipped with orthogonal filtrations, {\em Proc. Lond. Math. Soc.} {\bf 106} (2013), 980--1004. 

\bibitem{bss}T. Banica, A. Skalski and P.M. So\l tan, Noncommutative homogeneous spaces: the matrix case, {\em J. Geom. Phys.} {\bf 62} (2012), 1451--1466.

\bibitem{bsp}T. Banica and R. Speicher, Liberation of orthogonal Lie groups, {\em Adv. Math.} {\bf 222} (2009), 1461--1501.

\bibitem{bma}E.J. Beggs and S. Majid, Quantum Riemannian geometry, Springer (2020).

\bibitem{bpa}H. Bercovici and V. Pata, Stable laws and domains of attraction in free probability theory, {\em Ann. of Math.} {\bf 149} (1999), 1023--1060.

\bibitem{bdd}J. Bhowmick, F. D'Andrea and L. Dabrowski, Quantum isometries of the finite noncommutative geometry of the standard model, {\em Comm. Math. Phys.} {\bf 307} (2011), 101--131.

\bibitem{bd+}J. Bhowmick, F. D'Andrea, B. Das and L. Dabrowski, Quantum gauge symmetries in noncommutative geometry, {\em J. Noncommut. Geom.} {\bf 8} (2014), 433--471.

\bibitem{bg1}J. Bhowmick and D. Goswami, Quantum isometry groups: examples and computations, {\em Comm. Math. Phys.} {\bf 285} (2009), 421--444.

\bibitem{bg2}J. Bhowmick and D. Goswami, Quantum group of orientation preserving Riemannian isometries, {\em J. Funct. Anal.} {\bf 257} (2009), 2530--2572.

\bibitem{bic}J. Bichon, Half-liberated real spheres and their subspaces, {\em Colloq. Math.} {\bf 144} (2016), 273--287.

\bibitem{bdu}J. Bichon and M. Dubois-Violette, Half-commutative orthogonal Hopf algebras, {\em Pacific J. Math.} {\bf 263} (2013), 13--28.

\bibitem{bla}B. Blackadar, Operator algebras: theory of C$^*$-algebras and von Neumann algebras, Springer (2006).

\bibitem{bra}R. Brauer, On algebras which are connected with the semisimple continuous groups, {\em Ann. of Math.} {\bf 38} (1937), 857--872.

\bibitem{cc1}A.H. Chamseddine and A. Connes, The spectral action principle, {\em Comm. Math. Phys.} {\bf 186} (1997), 731--750.

\bibitem{cc2}A.H. Chamseddine and A. Connes, Why the standard model, {\em J. Geom. Phys.} {\bf 58} (2008), 38--47.

\bibitem{cpr}V. Chari and A. Pressley, A guide to quantum groups, Cambridge Univ. Press (1994).

\bibitem{chi}A. Chirvasitu, Residually finite quantum group algebras, {\em J. Funct. Anal.} {\bf 268} (2015), 3508--3533.

\bibitem{cfk}F. Cipriani, U. Franz and A. Kula, Symmetries of L\'evy processes on compact quantum groups, their Markov semigroups and potential theory, {\em J. Funct. Anal.} {\bf 266} (2014), 2789--2844. 

\bibitem{csn}B. Collins and P. \'Sniady, Integration with respect to the Haar measure on unitary, orthogonal and symplectic groups, {\em Comm. Math. Phys.} {\bf 264} (2006), 773--795.

\bibitem{co1}A. Connes, Noncommutative geometry, Academic Press (1994).

\bibitem{co2}A. Connes, On the spectral characterization of manifolds, {\em J. Noncommut. Geom.} {\bf 7} (2013), 1--82.

\bibitem{cdu}A. Connes and M. Dubois-Violette, Moduli space and structure of noncommutative 3-spheres, {\em Lett. Math. Phys.} {\bf 66} (2003), 91--121.

\bibitem{cla}A. Connes and G. Landi, Noncommutative manifolds, the instanton algebra and isospectral deformations, {\em Comm. Math. Phys.} {\bf 221} (2001), 141--160.

\bibitem{cma}A. Connes and M. Marcolli, Noncommutative geometry, quantum fields and motives, AMS (2008).

\bibitem{ddl}F. D'Andrea, L. Dabrowski and G. Landi, The noncommutative geometry of the quantum projective plane, {\em Rev. Math. Phys.} {\bf 20} (2008), 979--1006.

\bibitem{dfw}B. Das, U. Franz and X. Wang, Invariant Markov semigroups on quantum homogeneous spaces, {\em J. Noncommut. Geom.} {\bf 15} (2021), 531--580.

\bibitem{dgo}B. Das and D. Goswami, Quantum Brownian motion on noncommutative manifolds: construction, deformation and exit times, {\em  Comm. Math. Phys.} {\bf 309} (2012), 193--228.

\bibitem{dif}P. Di Francesco, Meander determinants, {\em Comm. Math. Phys.} {\bf 191} (1998), 543--583.

\bibitem{dir}P.A.M. Dirac, Principles of quantum mechanics, Oxford Univ. Press (1930).

\bibitem{doc}M.P. do Carmo, Riemannian geometry, Birkh\"auser (1992).

\bibitem{dri}V.G. Drinfeld, Quantum groups, Proc. ICM Berkeley (1986), 798--820.

\bibitem{fey}R.P. Feynman, R.B. Leighton and M. Sands, The Feynman lectures on physics, Caltech (1963).

\bibitem{fha}W. Fulton and J. Harris, Representation theory, Springer (1991).

\bibitem{go1}D. Goswami, Quantum group of isometries in classical and  noncommutative geometry, {\em Comm. Math. Phys.} {\bf 285} (2009), 141--160.

\bibitem{go2}D. Goswami, Existence and examples of quantum isometry groups for a class of compact metric spaces, {\em Adv. Math.} {\bf 280} (2015), 340--359.

\bibitem{go3}D. Goswami, Non-existence of genuine quantum symmetries of compact, connected smooth manifolds, {\em Adv. Math.} {\bf 369} (2020), 1--19.

\bibitem{gbh}D. Goswami and J. Bhowmick, Quantum isometry groups, Springer (2016).

\bibitem{gvf}J.M. Gracia-Bond\'ia, J.C. V\'arilly and H. Figueroa, Elements of noncommutative geometry, Birkh\"auser (2001).

\bibitem{gr1}D.J. Griffiths, Introduction to electrodynamics, Cambridge Univ. Press (2017).

\bibitem{gr2}D.J. Griffiths and D.F. Schroeter, Introduction to quantum mechanics, Cambridge Univ. Press (2018).

\bibitem{gha}P. Griffiths and J. Harris, Principles of algebraic geometry, Wiley (1994).

\bibitem{ega}A. Grothendieck and J. Dieudonn\'e, \'El\'ements de g\'eom\'etrie alg\'ebrique, IHES (1967).

\bibitem{har}J. Harris, Algebraic geometry, Springer (1992).

\bibitem{hrt}R. Hartshorne, Algebraic geometry, Springer (1977).

\bibitem{jim}M. Jimbo, A $q$-difference analog of $U(\mathfrak g)$ and the Yang-Baxter equation, {\em Lett. Math. Phys.} {\bf 10} (1985), 63--69.

\bibitem{jo1}V.F.R. Jones, Index for subfactors, {\em Invent. Math.} {\bf 72} (1983), 1--25.

\bibitem{jo2}V.F.R. Jones, On knot invariants related to some statistical mechanical models, {\em Pacific J. Math.} {\bf 137} (1989), 311--334.

\bibitem{jo3}V.F.R. Jones, Planar algebras I (1999).

\bibitem{lan}G. Landi, An introduction to noncommutative spaces and their geometry, Springer (1997).

\bibitem{lng}S. Lang, Algebra, Addison-Wesley (1993).

\bibitem{lax}P. Lax, Functional analysis, Wiley (2002).

\bibitem{lin}B. Lindst\"om, Determinants on semilattices, {\em Proc. Amer. Math. Soc.} {\bf 20} (1969), 207--208.

\bibitem{maj}S. Majid, Foundations of quantum group theory, Cambridge Univ. Press (1995).

\bibitem{mal}S. Malacarne, Woronowicz's Tannaka-Krein duality and free orthogonal quantum groups, {\em Math. Scand.} {\bf 122} (2018), 151--160.

\bibitem{mwe}A. Mang and M. Weber, Categories of two-colored pair partitions: Categories indexed by semigroups, {\em J. Combin. Theory Ser. A} {\bf 180} (2021), 1--37.

\bibitem{man}Y.I. Manin, Quantum groups and noncommutative geometry, Springer (2018). 

\bibitem{mpa}V.A. Marchenko and L.A. Pastur, Distribution of eigenvalues in certain sets of random matrices, {\em Mat. Sb.} {\bf 72} (1967), 507--536.

\bibitem{mar}M. Marcolli, Noncommutative cosmology, World Scientific (2018).

\bibitem{nas}J. Nash, The imbedding problem for Riemannian manifolds, {\em Ann. of Math.} {\bf 63} (1956), 20--63.

\bibitem{rwe}S. Raum and M. Weber, The full classification of orthogonal easy quantum groups, {\em Comm. Math. Phys.} {\bf 341} (2016), 751--779.

\bibitem{rud}W. Rudin, Real and complex analysis, McGraw-Hill (1966).

\bibitem{sha}I.R. Shafarevich, Basic algebraic geometry, Springer (1974).

\bibitem{twa}P. Tarrago and J. Wahl, Free wreath product quantum groups and standard invariants of subfactors, {\em Adv. Math.} {\bf 331} (2018), 1--57.

\bibitem{twe}P. Tarrago and M. Weber, Unitary easy quantum groups: the free case and the group case, {\em Int. Math. Res. Not.} {\bf 18} (2017), 5710--5750.

\bibitem{vdn}D.V. Voiculescu, K.J. Dykema and A. Nica, Free random variables, AMS (1992).

\bibitem{voi}C. Voigt, The Baum-Connes conjecture for free orthogonal quantum groups, {\em Adv. Math.} {\bf 227} (2011), 1873--1913.

\bibitem{von}J. von Neumann, Mathematical foundations of quantum mechanics, Princeton Univ. Press (1955).

\bibitem{wa1}S. Wang, Free products of compact quantum groups, {\em Comm. Math. Phys.} {\bf 167} (1995), 671--692.

\bibitem{wa2}S. Wang, Quantum symmetry groups of finite spaces, {\em Comm. Math. Phys.} {\bf 195} (1998), 195--211.

\bibitem{we1}S. Weinberg, Foundations of modern physics, Cambridge Univ. Press (2011).

\bibitem{we2}S. Weinberg, Lectures on quantum mechanics, Cambridge Univ. Press (2012).

\bibitem{wei}D. Weingarten, Asymptotic behavior of group integrals in the limit of infinite rank, {\em J. Math. Phys.} {\bf 19} (1978), 999--1001.

\bibitem{wey}H. Weyl, The theory of groups and quantum mechanics, Princeton Univ. Press (1931).

\bibitem{wig}E. Wigner, Characteristic vectors of bordered matrices with infinite dimensions, {\em Ann. of Math.} {\bf 62} (1955), 548--564.

\bibitem{wit}E. Witten, Quantum field theory and the Jones polynomial, {\em Comm. Math. Phys.} {\bf 121} (1989), 351--399.

\bibitem{wo1}S.L. Woronowicz, Compact matrix pseudogroups, {\em Comm. Math. Phys.} {\bf 111} (1987), 613--665.

\bibitem{wo2}S.L. Woronowicz, Tannaka-Krein duality for compact matrix pseudogroups. Twisted SU(N) groups, {\em Invent. Math.} {\bf 93} (1988), 35--76.

\end{thebibliography}
\end{document}